\documentclass{lms}

\usepackage{longtable}
\usepackage{textcomp}
\usepackage{amsmath}
\usepackage{amssymb}
\usepackage{graphicx}
\usepackage{setspace}
\usepackage{comment}
\usepackage{placeins}
\usepackage{float}

\newtheorem{theorem}{Theorem}[section] 

\newtheorem{proposition}[theorem]{Proposition}

\newnumbered{assertion}{Assertion}    
\newnumbered{conjecture}{Conjecture}  
\newnumbered{definition}{Definition}
\newnumbered{hypothesis}{Hypothesis}
\newnumbered{remark}{Remark}
\newnumbered{note}{Note}
\newnumbered{observation}{Observation}
\newnumbered{problem}{Problem}
\newnumbered{question}{Question}
\newnumbered{algorithm}{Algorithm}
\newnumbered{example}{Example}
\newunnumbered{notation}{Notation} 

\title[Modular Subgroups, Dessins and K3 Surfaces]
 {Modular Subgroups, Dessins d'Enfants and Elliptic K3 Surfaces} 

\author{Yang-Hui He, John McKay and James Read}

\classno{11G32, 11F06, 14J32, 14J28}

\extraline{The first author would like to thank the Science and Technology Facilities Council, UK, for an Advanced Fellowship and for STFC grant ST/J00037X/1, the Chinese Ministry of Education, for a Chang-Jiang Chair Professorship at NanKai University, the US NSF for grant CCF-1048082, as well as City University, London and Merton College, Oxford, for their enduring support.
The second author is grateful to the NSERC of Canada. The third author is indebted to Oriel College and the Department of Physics, Oxford for support.}

\begin{document}
\maketitle

\begin{abstract}
We consider the 33 conjugacy classes of genus zero, torsion-free modular subgroups, computing
ramification data and Grothendieck's \emph{dessins d'enfants}. In the
particular case of the index 36 subgroups, the corresponding Calabi-Yau
threefolds are identified, in analogy with the index 24 cases being associated with K3 surfaces.
In a parallel vein, we study the 112 semi-stable elliptic
fibrations over $\mathbb{P}^{1}$ as extremal K3 surfaces with six singular
fibres. In each case, a representative of the corresponding class of subgroups is
identified by specifying a generating set for that representative.
\end{abstract}

\section{Introduction} 
\label{intro}

\noindent The congruence subgroups of the modular group $\Gamma = PSL(2;\mathbb{Z})$ play a key role in many branches of modern mathematics.
Whereas the upper half of the complex plane $\mathcal{H}$, when adjoining cusps on the real line, can quotient $\Gamma$ to give a Riemann sphere, it is an interesting question as for which subgroups a similar sphere arises.
Such a subgroup is called genus zero.
Inspired by Moonshine, the authors of \cite{2} looked into the classification problem of congruence, genus zero, torsion free subgroups of the modular group and found these to be rather rare: there are only 33 conjugacy classes of these subgroups, of finite index $\mu\in\left\{ 6,12,24,36,48,60\right\} $.

Therefrom emerge many remarkable connections to geometry, graph theory, number theory and physics \cite{3,1,TopYui}.
First, we naturally extend the modular group action to a so-called modular surface for each of the groups whose Euler number is precisely $\mu$.
These turn out to be elliptically fibred over the genus zero $\mathbb{P}^{1}$ formed by quotienting $\mathcal{H}$ by the group and are, in fact, semi-stable in the sense that all singular fibres are of Kodaira type $I_n$.
Moreover, there are exactly $h = \frac{\mu}{6} + 2$ such singular fibres, and the tuple $\{n_1, \ldots, n_h\}$ is called the cusp widths.

Next, one can construct the Schreier coset graphs for these groups because they have finite index.
They turn out to be finite trivalent graphs of non-trivial topology.
We can interpret them as clean cubic dessin d'enfants in the sense of Grothendieck.
Focus has been made on the index $\mu=24$ cases because they are K3 surfaces \cite{1} and on the $\mu=12$ cases because they have interesting associated modular forms \cite{3}.
One can check that the Belyi rational maps for the dessins are precisely the elliptic $j$-invariants of the Weierstra\ss\ models for the modular surfaces with base $\mathbb{P}^{1}$ and we will complete this story in the present paper.

In a slightly different but obviously related direction, the classification of extremal, semi-elliptic K3 surfaces with 6 singular fibres was wonderfully performed by Miranda and Persson in \cite{6}, and the associated dessins nicely compiled in \cite{7}.
It is therefore a natural question whether they can be affiliated with subgroups of the modular group, not necessarily congruence.
Looking at these dessins, one see them to be trivalent and clean, and since these can be interpreted as the Schreier coset graphs of a subgroup of the modular group,  the answer must be affirmative.
We shall undertake this task, using the dessins to find the generating set of a representative of
the corresponding conjugacy class of subgroups in each case.

The organization of the paper is as follows.
In \S\ref{s:33big}, we introduce the 33 conjugacy classes of congruence, torsion-free, genus zero subgroups of the modular group, especially their $j$-invariants, dessins and Schreier coset graphs.
We will see that as the index 24 cases are attached to K3 surfaces, so too can the index 36 cases be associated to Calabi-Yau threefolds.
Then, in \S\ref{s:K3}, we turn to the 112 extremal semi-stable K3 surfaces, studying their dessins and interpreting them as coset graphs of some conjugacy class of index 24 modular subgroups.
Via a discussion of the elliptic $j$-invariants as Belyi maps and the so-called cartographic groups, we find the generating set for a representative of the class of subgroups in each case.
Finally, we conclude with some outlook and prospects.

\section{Torsion-Free, Genus Zero Congruence Subgroups}
\label{s:33big}
In this section, we first recall some standard notation regarding
the modular group $\Gamma$ and its congruence subgroups. We focus
specifically on the conjugacy classes of the \emph{torsion-free, genus zero, congruence
subgroups}; there are 33 such classes \cite{2}. 
In each case, we will compute the
ramification data, and draw the dessin d'enfant which encodes this
data. 
The cases of index 24 and 36 are particularly interesting: the former can be associated to K3 surfaces and we will demonstrate, using the $j$-invariants, that each of the index 36 subgroups can be identified with a specific Calabi-Yau threefold.

\subsection{The Modular Group and Congruence Subgroups}\label{s:Gamma}

First, the central object of our study is the \emph{modular group}
\begin{equation}
\Gamma\equiv\Gamma\left(1\right)=\mathrm{PSL}(2;\mathbb{Z})=\mathrm{SL}(2;\mathbb{Z})/\left\{ \pm I\right\} \ .
\end{equation}
This is the group of linear fractional transformations $z\rightarrow\frac{az+b}{cz+d}$, with $a,b,c,d\in\mathbb{Z}$ and $ad-bc=1$. It is generated by the
transformations $T$ and $S$ defined by:
\begin{equation}
T\tau=\tau+1\quad,\quad S\tau=-1/\tau \ .
\end{equation}
One presentation of $\Gamma$ is $\left\langle S,T|S^{2}=\left(ST\right)^{3}=I\right\rangle $, and we will later discuss the presentations of certain modular subgroups.
$2\times2$
matrices for $S$ and $T$ are as follows:
\begin{equation}
T=\begin{pmatrix}1 & 1\\
0 & 1
\end{pmatrix}\quad,\quad S=\begin{pmatrix}0 & -1\\
1 & 0
\end{pmatrix} \ .
\end{equation}

Letting $x = S$ and $y = ST$ denote the elements of order 2 and 3 respectively,
we see that $\Gamma$ is the free product of the cyclic groups $C_{2}=\left\langle x|x^{2}=I\right\rangle $
and $C_{3}=\left\langle y|y^{3}=I\right\rangle $. 
It follows that $2\times2$ matrices for $x$ and $y$ are:
\begin{equation}\label{xy}
x=\begin{pmatrix}0 & -1\\
1 & 0
\end{pmatrix}\quad,\quad y=\begin{pmatrix}0 & -1\\
1 & 1
\end{pmatrix} \ .
\end{equation}

\paragraph{Cayley Graphs: }
Given the free product structure of $\Gamma$, we see that its Cayley graph is an infinite free trivalent tree, but with each node replaced by an oriented triangle.
Now, for any \emph{finite index} subgroup $H$ of $\Gamma$, we can quotient the Cayley graph to arrive at a finite cubic graph by associating nodes to right cosets and edges between cosets which are related by action of a group element.
In other words, this graph encodes the permutation representation of $\Gamma$ acting on the right cosets of $H$. 
This is called a {\em Schreier coset graph} and will occupy a central role in our ensuing discussions.

\paragraph{Congruence Subgroups: }
The most important subgroups of $\Gamma$ are the so-called congruence subgroups, defined by having the entries in the generating matrices $S$ and $T$ obeying some modular arithmetic.
Some conjugacy classes of congruence subgroups of particular note are the following:
\begin{itemize}
\item Principal congruence subgroups:
\[
\Gamma\left(m\right):=\left\{ A\in\mathrm{SL}(2;\mathbb{Z})\;,\; A\equiv\pm I\;\mathrm{mod}\; m\right\} /\left\{ \pm I\right\} ;
\]

\item Congruence subgroups of level $m$: subgroups of $\Gamma$ containing
$\Gamma\left(m\right)$ but not any $\Gamma\left(n\right)$ for $n<m$;
\item Unipotent matrices:
\[
\Gamma_{1}\left(m\right):=\left\{ A\in\mathrm{SL}(2;\mathbb{Z})\;,\; A\equiv\pm\begin{pmatrix}1 & b\\
0 & 1
\end{pmatrix}\;\mathrm{mod}\; m\right\} /\left\{ \pm I\right\} ;
\]

\item Upper triangular matrices:
\[
\Gamma_{0}\left(m\right):=\left\{ \begin{pmatrix}a & b\\
c & d
\end{pmatrix}\in\Gamma\;,\; c\equiv0\;\mathrm{mod}\; m\right\} /\left\{ \pm I\right\} .
\]
\end{itemize}

In \cite{1}, our attention is drawn to the conjugacy classes of a particular family of subgroups
of $\Gamma$: the so-called\emph{ genus zero, torsion-free} congruence
subgroups. 
By \emph{torsion-free} we mean that the subgroup
contains no element of finite order other than the identity. 
To explain
\emph{genus zero}, first recall that the modular group acts on the
upper half-plane $\mathcal{H}:=\left\{ \tau\in\mathbb{C}\;,\;\mathrm{Im}\left(\tau\right)>0\right\} $
by linear fractional transformations $z\rightarrow\frac{az+b}{cz+d}$.
$\mathcal{H}$ gives rise to a compactification $\mathcal{H}^{*}$ when adjoining
\emph{cusps}, which are points on $\mathbb{R}\cup\infty$ fixed
under some parabolic element (i.e. an element $A\in\Gamma$ not equal to the identity
 and for which $\mathrm{Tr}\left(A\right)=2$). The quotient $\mathcal{H}^{*}/\Gamma$
is a compact Riemann surface of genus 0, i.e. a sphere. It turns out
that with the addition of appropriate cusp points, the extended upper
half plane $\mathcal{H}^{*}$ factored by various congruence subgroups
will also be compact Riemann surfaces, possibly of higher genus. Such
a Riemann surface, as a complex algebraic variety, 
is called a {\em modular curve}. 
The genus of a subgroup of the modular group is the genus of the modular
curve produced in this way. 

Using the Riemann-Hurwitz
theorem, the genus of the modular curve can be shown to be:
\begin{equation}\label{RH}
g_{C}=1+\frac{\mu}{12}-\frac{\nu_{2}}{4}-\frac{\nu_{3}}{3}-\frac{\nu_{\infty}}{2}
\ ,
\end{equation}
where $\mu$ is the index of the subgroup, $\nu_{2}$, $\nu_{3}$
are the number of inequivalent elliptic points of order 2 and 3 respectively,
and $\nu_{\infty}$ is the number of inequivalent cusp points.

Our focus will be on genus zero, torsion-free subgroups,
thus we have that
\begin{equation}
g_{C}=\nu_{2}=\nu_{3}=0 \ ,
\end{equation}
whence, by \eqref{RH}, we have a further constraint that
$\mu=6\nu_{\infty}-12$, so that the index of all
our subgroups of concern is a multiple of 6. This means that the index
of the subgroup in question also constrains $\nu_{\infty}$, so that
for a subgroup of index $6n$, we have $\nu_{\infty}=n+2$.

\subsubsection{The 33 Conjugacy Classes of Subgroups}
As classified in \cite{2}, the conjugacy classes of the genus zero torsion-free congruence subgroups of the modular group are very rare. 
There are only 33 of them:
\begin{itemize}
\item $\Gamma\left(n\right)$ with $n=2,3,4,5$;
\item $\Gamma_{0}\left(n\right)$ with $n=4,6,8,9,12,16,18$;
\item $\Gamma_{1}\left(n\right)$ with $n=5,7,8,9,10,12$;
\item The intersections $\Gamma_{0}\left(a\right)\cap\Gamma\left(b\right)$
for $\left\{ a,b\right\} =\left\{ 4,2\right\} ,\left\{ 3,2\right\} ,\left\{ 8,2\right\} ,\left\{ 2,3\right\} ,\left\{ 25,5\right\} $;
$\Gamma_{1}\left(8\right)\cap\Gamma\left(2\right)$ and $\Gamma_{0}\left(16\right)\cap\Gamma_{1}\left(8\right)$;
\item The congruence subgroups
\[
\Gamma\left(m;\frac{m}{d},\epsilon,\chi\right):=\left\{ \pm\begin{pmatrix}1+\frac{m}{\epsilon\chi}\alpha & d\beta\\
\frac{m}{\chi}\gamma & 1+\frac{m}{\epsilon\chi}\delta
\end{pmatrix}\;,\;\gamma\equiv\alpha\;\mathrm{mod}\;\chi\right\} 
\]
 with $\left(m,d,\epsilon,\chi\right)=\left(8,2,1,2\right)$, $\left(12,2,1,2\right)$,
$\left(16,1,2,2\right)$, $\left(27,1,3,3\right)$, $\left(8,4,1,2\right)$,
$\left(9,3,1,3\right)$, $\left(16,2,2,2\right)$, $\left(24,1,2,2\right)$,
$\left(32,1,4,2\right)$.
\end{itemize}

For reference, in Table 1 we organize these 33 conjugacy classes of subgroups by index, and we will shortly explain the notation of the ramification data.

\begin{table}[b]\vspace*{-3ex}

{\scriptsize

\begin{tabular}{|c|c|c|c|c|c|c|}
\cline{1-3} \cline{5-7} 
Index & Group & Ramification &  & Index & Group & Ramification\tabularnewline
\cline{1-3} \cline{5-7} 
6 & $\Gamma\left(2\right)$ & $\begin{Bmatrix}3^{2}\\
2^{3}\\
2^{3}
\end{Bmatrix}$ &  & 36 & $\Gamma_{0}\left(2\right)\cap\Gamma\left(3\right)$ & $\begin{Bmatrix}3^{12}\\
2^{18}\\
3^{4},6^{4}
\end{Bmatrix}$\tabularnewline
\cline{1-3} \cline{5-7} 
6 & $\Gamma_{0}\left(4\right)$ & $\begin{Bmatrix}3^{2}\\
2^{3}\\
1^{2},4
\end{Bmatrix}$ &  & 36 & $\Gamma_{1}\left(9\right)$ & $\begin{Bmatrix}3^{12}\\
2^{18}\\
1^{3},3^{2},9^{3}
\end{Bmatrix}$\tabularnewline
\cline{1-3} \cline{5-7} 
12 & $\Gamma\left(3\right)$ & $\begin{Bmatrix}3^{4}\\
2^{6}\\
3^{4}
\end{Bmatrix}$ &  & 36 & $\Gamma\left(9;3,1,3\right)$ & $\begin{Bmatrix}3^{12}\\
2^{18}\\
3^{6},9^{2}
\end{Bmatrix}$\tabularnewline
\cline{1-3} \cline{5-7} 
12 & $\Gamma_{0}\left(4\right)\cap\Gamma\left(2\right)$ & $\begin{Bmatrix}3^{4}\\
2^{6}\\
2^{2},4^{2}
\end{Bmatrix}$ &  & 36 & $\Gamma_{1}\left(10\right)$ & $\begin{Bmatrix}3^{12}\\
2^{18}\\
1^{2},2^{2},5^{2},10^{2}
\end{Bmatrix}$\tabularnewline
\cline{1-3} \cline{5-7} 
12 & $\Gamma_{1}\left(5\right)$ & $\begin{Bmatrix}3^{4}\\
2^{6}\\
1^{2},5^{2}
\end{Bmatrix}$ &  & 36 & $\Gamma_{0}\left(18\right)$ & $\begin{Bmatrix}3^{12}\\
2^{18}\\
1^{3},2^{3},9,18
\end{Bmatrix}$\tabularnewline
\cline{1-3} \cline{5-7} 
12 & $\Gamma_{0}\left(6\right)$ & $\begin{Bmatrix}3^{4}\\
2^{6}\\
1,2,3,6
\end{Bmatrix}$ &  & 36 & $\Gamma\left(27;27,3,3\right)$ & $\begin{Bmatrix}3^{12}\\
2^{18}\\
1^{6},3,27
\end{Bmatrix}$\tabularnewline
\cline{1-3} \cline{5-7} 
12 & $\Gamma_{0}\left(8\right)$ & $\begin{Bmatrix}3^{4}\\
2^{6}\\
1^{2},2,8
\end{Bmatrix}$ &  & 48 & $\Gamma_{1}\left(8\right)\cap\Gamma\left(2\right)$ & $\begin{Bmatrix}3^{16}\\
2^{24}\\
2^{4},4^{2},8^{4}
\end{Bmatrix}$\tabularnewline
\cline{1-3} \cline{5-7} 
12 & $\Gamma_{0}\left(9\right)$ & $\begin{Bmatrix}3^{4}\\
2^{6}\\
1^{3},9
\end{Bmatrix}$ &  & 48 & $\Gamma\left(8;2,1,2\right)$ & $\begin{Bmatrix}3^{16}\\
2^{24}\\
4^{8},8^{2}
\end{Bmatrix}$\tabularnewline
\cline{1-3} \cline{5-7} 
24 & $\Gamma\left(4\right)$ & $\begin{Bmatrix}3^{8}\\
2^{12}\\
4^{6}
\end{Bmatrix}$ &  & 48 & $\Gamma_{1}\left(12\right)$ & $\begin{Bmatrix}3^{16}\\
2^{24}\\
1^{2},2,3^{2},4^{2},6,12^{2}
\end{Bmatrix}$\tabularnewline
\cline{1-3} \cline{5-7} 
24 & $\Gamma_{0}\left(3\right)\cap\Gamma\left(2\right)$ & $\begin{Bmatrix}3^{8}\\
2^{12}\\
2^{3},6^{3}
\end{Bmatrix}$ &  & 48 & $\Gamma\left(12;6,1,2\right)$ & $\begin{Bmatrix}3^{16}\\
2^{24}\\
2^{4},4,6^{4},12
\end{Bmatrix}$\tabularnewline
\cline{1-3} \cline{5-7} 
24 & $\Gamma_{1}\left(7\right)$ & $\begin{Bmatrix}3^{8}\\
2^{12}\\
1^{3},7^{3}
\end{Bmatrix}$ &  & 48 & $\Gamma_{0}\left(16\right)\cap\Gamma_{1}\left(8\right)$ & $\begin{Bmatrix}3^{16}\\
2^{24}\\
1^{4},2^{2},4^{2},16^{2}
\end{Bmatrix}$\tabularnewline
\cline{1-3} \cline{5-7} 
24 & $\Gamma_{1}\left(8\right)$ & $\begin{Bmatrix}3^{8}\\
2^{12}\\
1^{2},2,4,8^{2}
\end{Bmatrix}$ &  & 48 & $\Gamma\left(16;8,2,2\right)$ & $\begin{Bmatrix}3^{16}\\
2^{24}\\
2^{8},16^{2}
\end{Bmatrix}$\tabularnewline
\cline{1-3} \cline{5-7} 
24 & $\Gamma_{0}\left(8\right)\cap\Gamma\left(2\right)$ & $\begin{Bmatrix}3^{8}\\
2^{12}\\
2^{4},8^{2}
\end{Bmatrix}$ &  & 48 & $\Gamma\left(24;24,2,2\right)$ & $\begin{Bmatrix}3^{16}\\
2^{24}\\
1^{4},3^{4},8,24
\end{Bmatrix}$\tabularnewline
\cline{1-3} \cline{5-7} 
24 & $\Gamma\left(8;4,1,2\right)$ & $\begin{Bmatrix}3^{8}\\
2^{12}\\
2^{2},4^{3},8
\end{Bmatrix}$ &  & 48 & $\Gamma\left(32;32,4,2\right)$ & $\begin{Bmatrix}3^{16}\\
2^{24}\\
1^{8},8,32
\end{Bmatrix}$\tabularnewline
\cline{1-3} \cline{5-7} 
24 & $\Gamma_{0}\left(12\right)$ & $\begin{Bmatrix}3^{8}\\
2^{12}\\
1^{2},3^{2},4,12
\end{Bmatrix}$ &  & 60 & $\Gamma\left(5\right)$ & $\begin{Bmatrix}3^{20}\\
2^{30}\\
5^{12}
\end{Bmatrix}$\tabularnewline
\cline{1-3} \cline{5-7} 
24 & $\Gamma_{0}\left(16\right)$ & $\begin{Bmatrix}3^{8}\\
2^{12}\\
1^{4},4,16
\end{Bmatrix}$ &  & 60 & $\Gamma_{0}\left(25\right)\cap\Gamma_{1}\left(5\right)$ & $\begin{Bmatrix}3^{20}\\
2^{30}\\
1^{10},25^{2}
\end{Bmatrix}$\tabularnewline
\cline{1-3} \cline{5-7} 
24 & $\Gamma\left(16;16,2,2\right)$ & $\begin{Bmatrix}3^{8}\\
2^{12}\\
1^{2},2^{3},16
\end{Bmatrix}$ & \multicolumn{1}{c}{} & \multicolumn{1}{c}{} & \multicolumn{1}{c}{} & \multicolumn{1}{c}{}\tabularnewline
\cline{1-3} 
\end{tabular}
}
\label{t:33}
\caption[]{The 33 conjugacy classes of genus zero, torsion free, congruence subgroups of the modular group, organized by index (which can be 6, 12, 24, 36, 48 or 60). The ramification data for each group are defined in \eqref{ram}.}
\end{table}

\subsection{Modular Surfaces, Ramification Data and $j$-Functions}
\label{s:surface}
As detailed in \cite{1}, to each torsion-free, genus zero modular
subgroup corresponds a distinct {\em Shioda modular surface} \cite{shioda}.
We briefly recall how this is constructed \cite{TopYui}.
First, we extend the action of any subgroup $G$ of $\Gamma$ on $\mathcal{H}$ to an action
\begin{equation}
\mathcal{H} \times \mathbb{C} \ni (\tau, z)
\longrightarrow \left( \gamma \tau, \frac{z + m \tau + n}{c \tau + d}
\right) \ ,
\end{equation}
for $\gamma = \begin{pmatrix}a & b\\
c & d
\end{pmatrix} \in G$  and $(m,n) \in \mathbb{Z}^{2}$.
Next, the quotient of $\mathcal{H} \times \mathbb{C}$ by the above automorphism defines a surface equipped with a morphism to the modular curve $C$ arising from the quotient of $\mathcal{H}$ by $\tau \to \gamma \tau$.
The fiber over the image of this morphism to the modular curve is generically an elliptic curve corresponding to the lattice $\mathbb{Z} \oplus \mathbb{Z} \tau$.
What we have therefore is a complex surface which is an elliptic fibration
over the modular curve associated to $G$.
This is the modular surface.

In our present case, because the base modular curves are genus zero (being the Riemann sphere $\mathbb{P}^1_C$), our modular surfaces will be elliptic fibrations over $\mathbb{P}^1_C$.
It turns out the Euler number of this surface is the index $\mu$ of the group $H$.
Hence, the astute reader will recognize that the index 24 groups should give us K3 surfaces, in particular elliptic K3 surfaces (we remark that $\Gamma_1(7)$, in particular, was realized by Tate to be related to a K3 surface). We will return to this point later, especially to the observation that index 36 groups, by analogy, give Calabi-Yau threefolds.

Furthermore, the elliptic modular surface is {\em semi-stable} in the sense that all the singular fibres are of Kodaira type $I_n$; the set of integers $(n_1, n_2, \ldots, n_k)$ denoting the indices of the singular fibres are the {\em cusp widths}.
For simplicity, we arrange this set in increasing order and group repeated entries by the short-hand $n_i^{a_i}$ with $a_i$ being the multiplicity.
Now, we are ready to introduce the concept of {\em ramification data} associated with a subgroup (and its modular surface); this is the following triplet of integer vectors:
\begin{equation}\label{ram}
\begin{Bmatrix}
3^{V}\\
2^{E}\\
n_1^{a_1},  n_2^{a_2}, \ldots, n_k^{a_k}  
\end{Bmatrix}
\ ,
\qquad
\begin{array}{rcl}
V &=& \mbox{ number of nodes in the Schreier coset graph}\\
E &=& \mbox{ number of edges}\\
\{n_i^{a_i}\} &=& \mbox{ cusp widths of the elliptic modular surface }
\end{array}
\end{equation}
Note that the $3$ reminds us that all nodes are trivalent and 2, that each edge links two nodes.
We have marked this ramification data in Table 1.

Indeed, given an elliptic fibration, we can write down a Weierstra\ss\ equation, from which we can then extract a (Klein) $j$-function.
Because of the fibration over $\mathbb{P}_{C}^{1}$, the $j$-function will explicitly depend on the projective coordinate of this base Riemann sphere and can thus be treated as a rational map from the base $\mathbb{P}_{C}^{1}$ to its range $\mathbb{P}^{1}$.
All of these $j$-functions are given in \cite{2} and certain simplifications, in \cite{1,TopYui}.
Generalizing \cite{1,6}, we find the nice fact that
\begin{proposition}
The $j$-function $j\left(t \in\mathbb{P}_{C}^{1} \right) \mapsto \mathbb{P}^{1}$ is ramified at only $\left(0,1,\infty\right)$. 
\end{proposition}
In other words, the only points $t\in\mathbb{P}_{C}^{1}$ such that $\frac{\mathrm{d}}{\mathrm{d}t}j\left(t\right)\left|_{\tilde{t}}\right.=0$ are such that $j(\tilde{t}) = 0, 1$ or $\infty$.
Moreover, the order of vanishing of the Taylor series at $t$ is the ramification index
(we remark that depending on convention, we may need to appropriately normalize to the $j$-function by a factor of 1728).

Specifically, for our 33 conjugacy classes of subgroups, we find that for a group of index $6n$,
there are $2n$ pre-images of 0 and the ramification indices of all
these marked points is $3$; there are $3n$ pre-images of $1$ and
ramification indices of all these marked points is $2$; there are
$n+2$ pre-images of $\infty$, the ramification indices of which
are given as the cusp widths in Table 1.

\subsection{Dessins d'Enfants: $j$-Functions as Belyi Maps}
\label{s:dessin}
Maps to $\mathbb{P}^{1}$ ramified only at $\left(0,1,\infty\right)$
are called {\em Belyi maps} and can be represented graphically as
Grothendieck's \emph{dessins d'enfants}. 
To draw such a dessin is
simple: given ramification data $\begin{Bmatrix}
\{r_{0}\left(i\right)\}\\
\{r_{1}\left(i\right)\}\\
\{r_{\infty}\left(i\right)\}
\end{Bmatrix}$ specifying the ramification indices at the various pre-images of 0, 1 and infinity, one marks one white node for the $i$-th pre-image of 0, with
$r_{0}\left(i\right)$ edges emanating therefrom; similarly, one marks
one black node for the $j$-th pre-image of 1, with $r_{1}\left(j\right)$ edges. Thus we have a {\em bipartite graph} embedded on a Riemann sphere, with $W$ white nodes and $B$ black nodes. Now we connect the
nodes with the edges, joining only black with white, such that we
have $I$ faces, each being a polygon with $2r_{\infty}\left(k\right)$
sides. Dessins for all the 33 torsion-free, genus zero modular subgroups
are presented in \S\ref{s:dessin33}.

The Schreier coset graphs for the 33 conjugacy classes of subgroups are discussed and drawn explicitly in \cite{2}.
In our case, let $\Gamma^{*}$ be the torsion-free, genus zero modular subgroup
in question. The coset graphs display the permutation
action of generators $x$ and $y$ with the relations $x^{2}=y^{3}=I$
on 
 the coset space of $\mathrm{\mathit{PSL}}(2,\mathbb{Z})=\bigcup_{i=1}^{\mu}\Gamma^{*}x_{i}$
where $\mu$ is the index of $\Gamma^{*}$ in $\mathit{\mathrm{\mathit{PSL}}}(2,\mathbb{Z})$.
These graphs have $\mu$ nodes, where $\mu\in\left\{ 6,12,24,36,48,60\right\} $ is the index
of the modular subgroup in question. The graphs are built from simple edges
($x$) and triangles ($y$) which are assumed to be positively oriented. Reversal of orientation
corresponds to applying the automorphism of the modular group that inverts each of the
two generators $x$ and $y$.

Interestingly, there is a direct connection between the Schreier coset
graph and the dessins d'enfant for each modular subgroup \cite{dessincoset}.
This connection can be stated as follows: 
\begin{proposition}
The dessins can be constructed from
the Schreier coset graphs by replacing each positively oriented triangle
with a white node, and inserting a black node into every edge. Equivalently,
the Schreier coset graphs can be constructed from the dessins by replacing
each white node with a positively oriented triangle, and removing
the black node from every edge.
\end{proposition}

The process of inserting into an edge a bivalent node of a different colouring is standard in the study of dessins d'enfants; this gives rise to so-called {\it clean} dessins.
The contraction of triangles into a node is a less obvious one and should correspond to taking a cube root of the polynomial factors in the Belyi map.


\subsubsection{Dessins for the 33 Conjugacy Classes of Subgroups}\label{s:dessin33}

In this section,
we draw the dessins d'enfants for all 33 conjugacy classes of
genus zero, torsion-free congruence subgroups, using the ramification
data in Table 1. Each dessin is labelled by its corresponding conjugacy class of modular subgroups, and the
unique cusp widths for that class.
The reader can visually confirm the aforementioned correspondence between these dessins and the Schreier coset graphs constructed in \cite{2}.

\begin{center}
\begin{minipage}[t]{0.25\textwidth}%
\begin{center}
\includegraphics[scale=0.1]{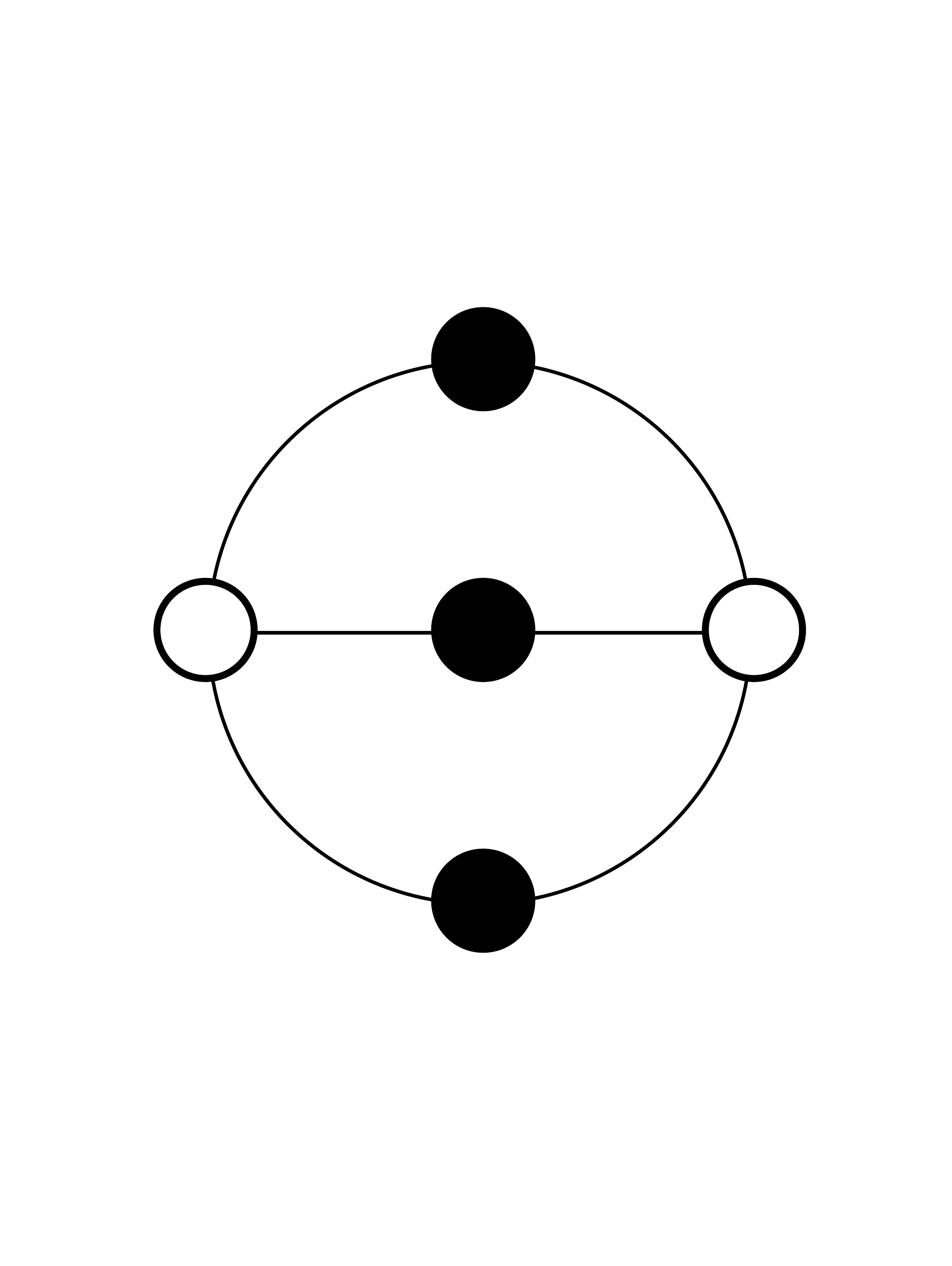}
\par\end{center}

\begin{center}
$\Gamma\left(2\right)$\\
\scriptsize $2,2,2$ \scriptsize
\par\end{center}%
\end{minipage}%
\begin{minipage}[t]{0.25\textwidth}%
\begin{center}
\includegraphics[scale=0.1]{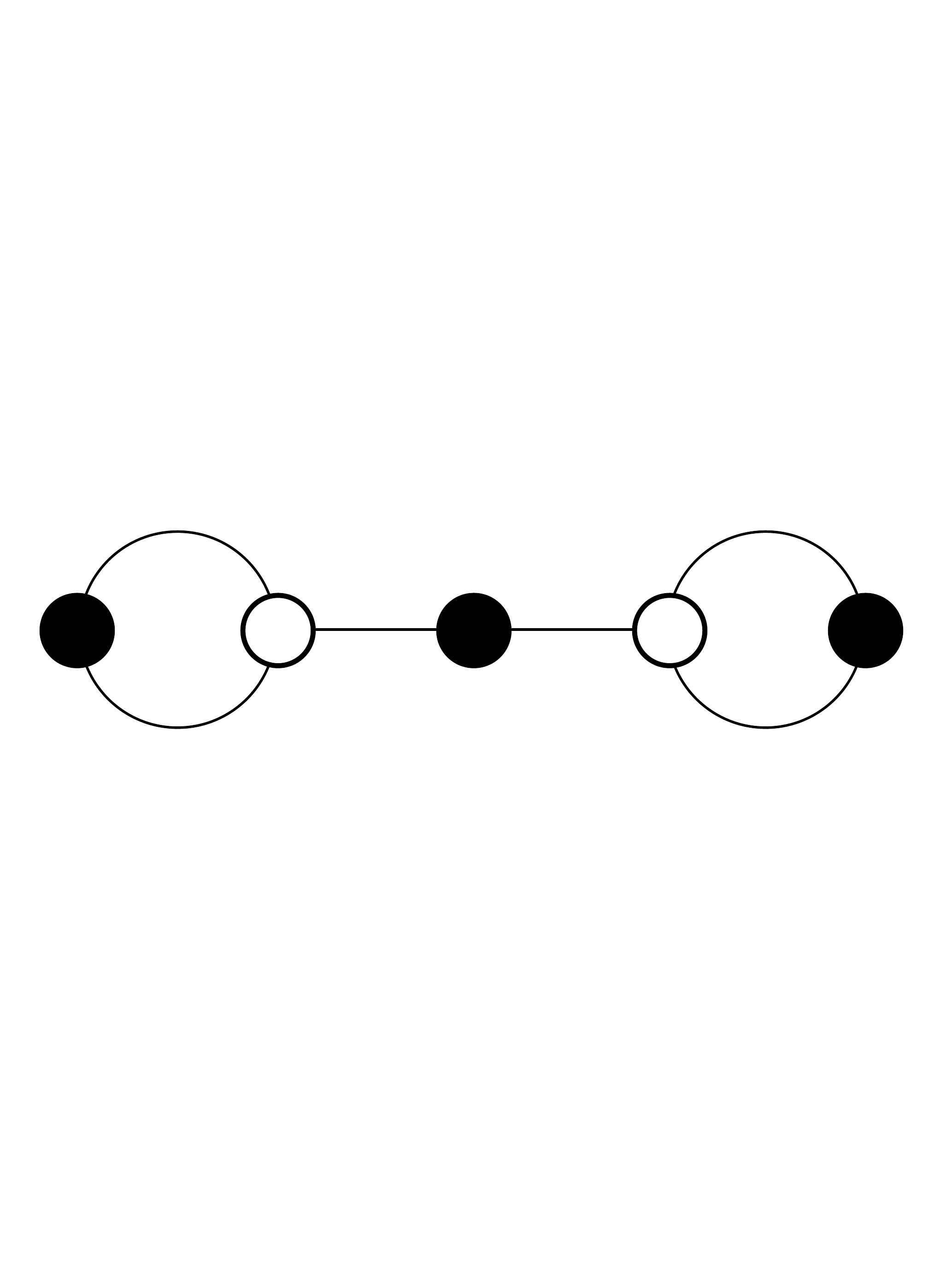}
\par\end{center}

\begin{center}
$\Gamma_{0}\left(4\right)$\\
\scriptsize $4,1,1$ \scriptsize
\par\end{center}%
\end{minipage}%
\begin{minipage}[t]{0.25\textwidth}%
\begin{center}
\includegraphics[scale=0.1]{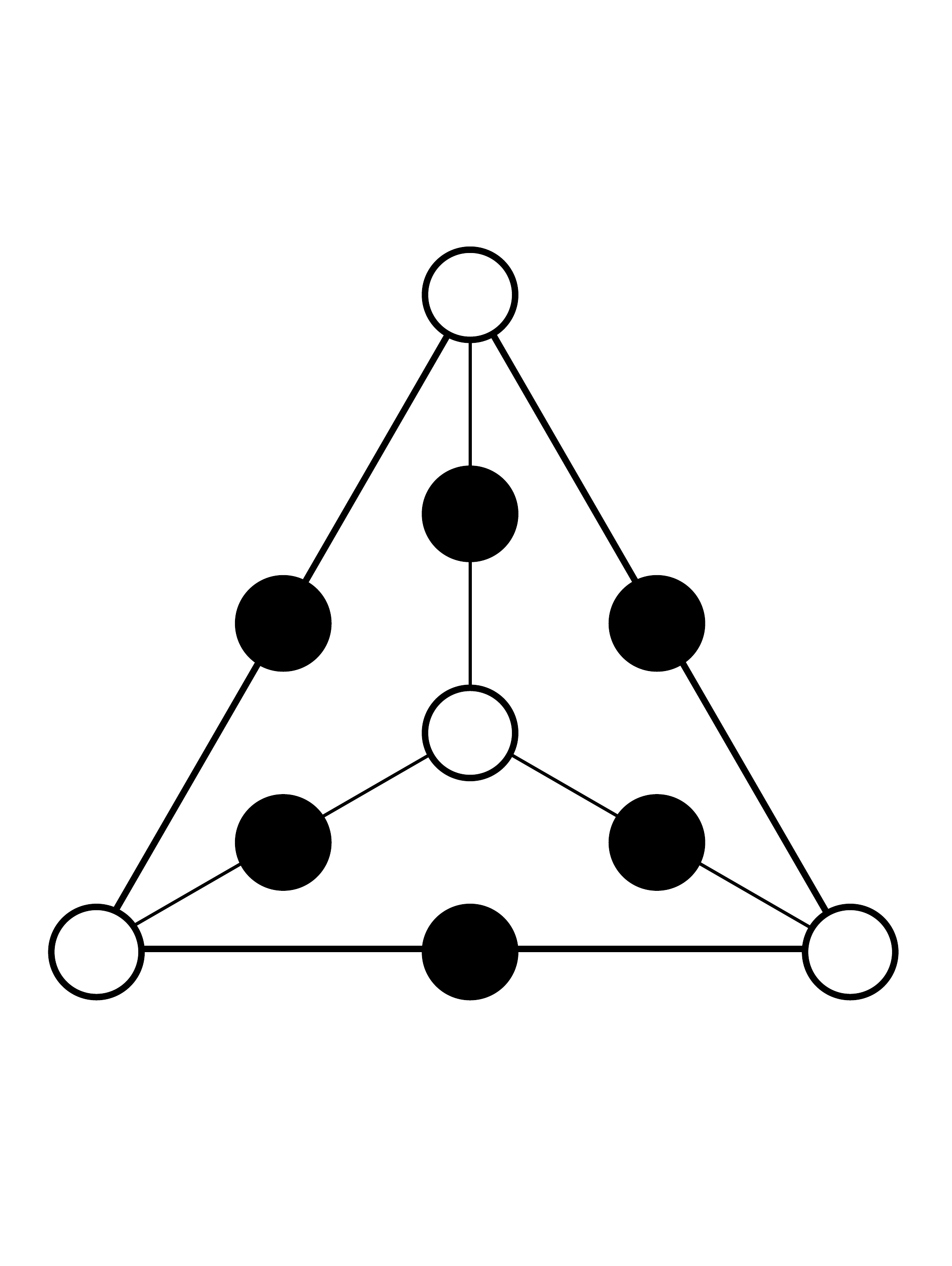}
\par\end{center}

\begin{center}
$\Gamma\left(3\right)$\\
\scriptsize $3,3,3,3$ \scriptsize
\par\end{center}%
\end{minipage}%
\begin{minipage}[t]{0.25\textwidth}%
\begin{center}
\includegraphics[scale=0.1]{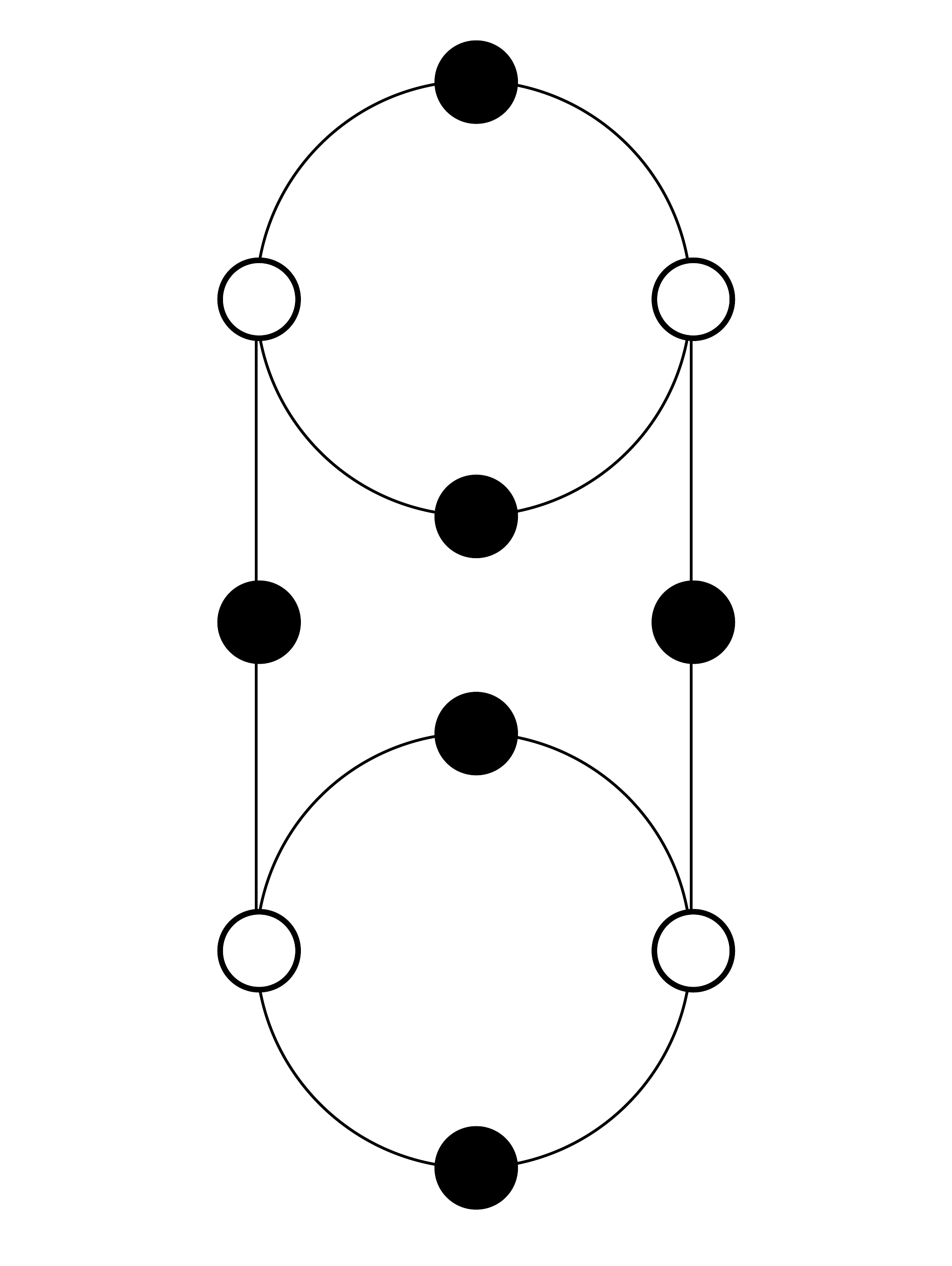}
\par\end{center}

\begin{center}
$\Gamma_{0}\left(4\right)\cap\Gamma\left(2\right)$\\
\scriptsize $4,4,2,2$ \scriptsize
\par\end{center}%
\end{minipage}
\par\end{center}

\begin{center}
\begin{minipage}[t]{0.25\textwidth}%
\begin{center}
\includegraphics[scale=0.15]{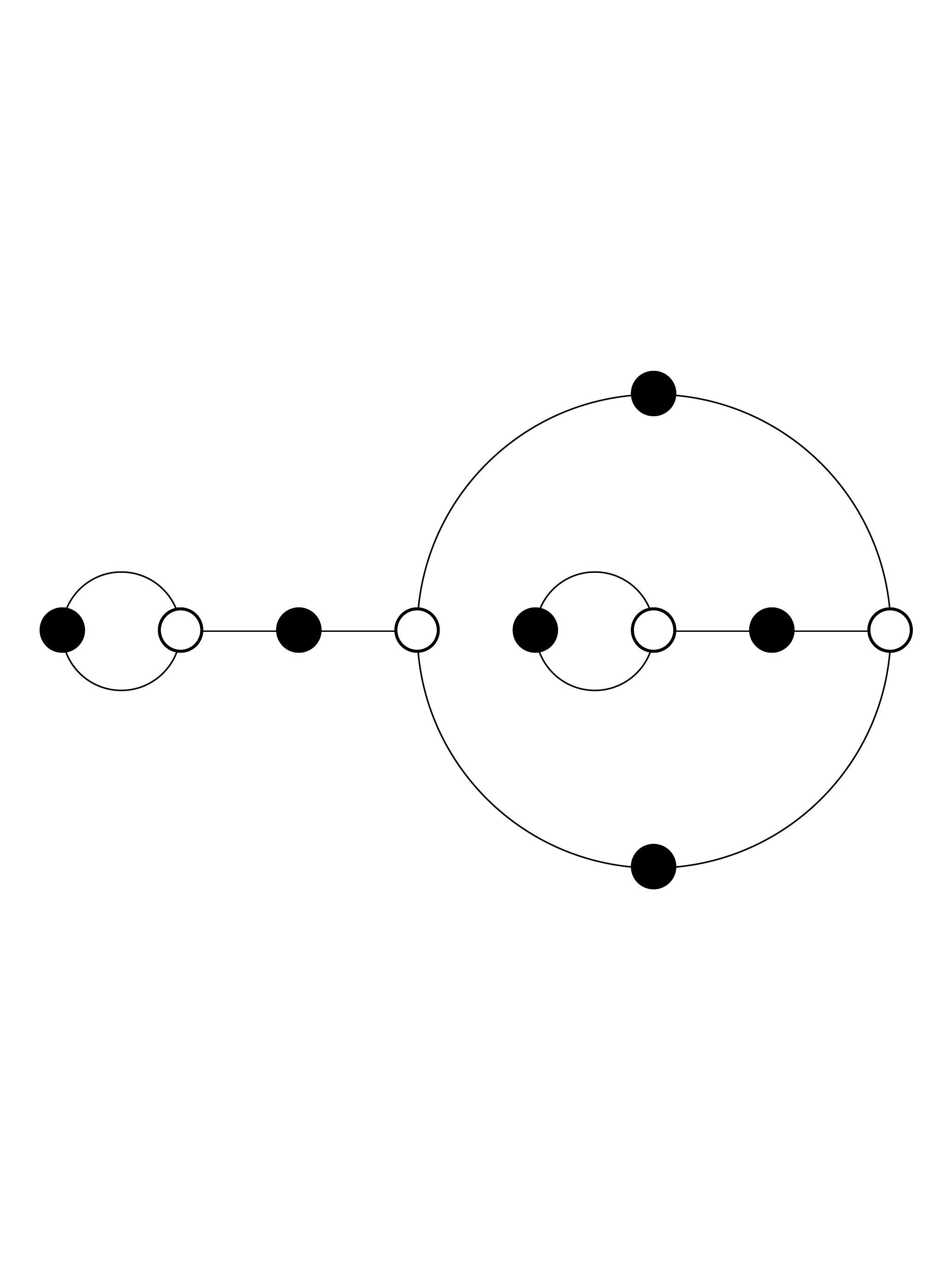}
\par\end{center}

\begin{center}
$\Gamma_{1}\left(5\right)$\\
\scriptsize $5,5,1,1$ \scriptsize
\par\end{center}%
\end{minipage}%
\begin{minipage}[t]{0.25\textwidth}%
\begin{center}
\includegraphics[scale=0.15]{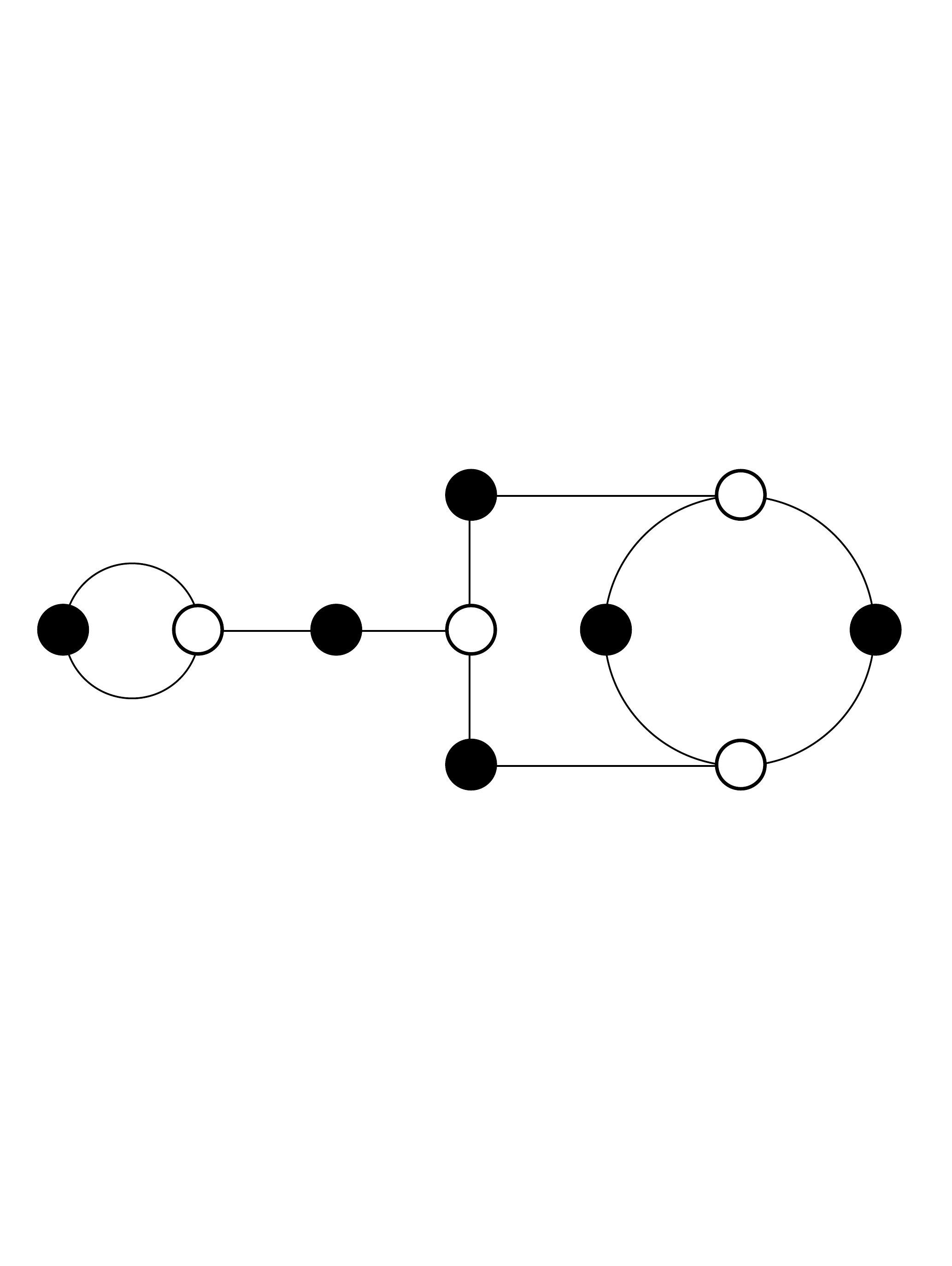}
\par\end{center}

\begin{center}
$\Gamma_{0}\left(6\right)$\\
\scriptsize $6,3,2,1$ \scriptsize
\par\end{center}%
\end{minipage}%
\begin{minipage}[t]{0.25\textwidth}%
\begin{center}
\includegraphics[scale=0.15]{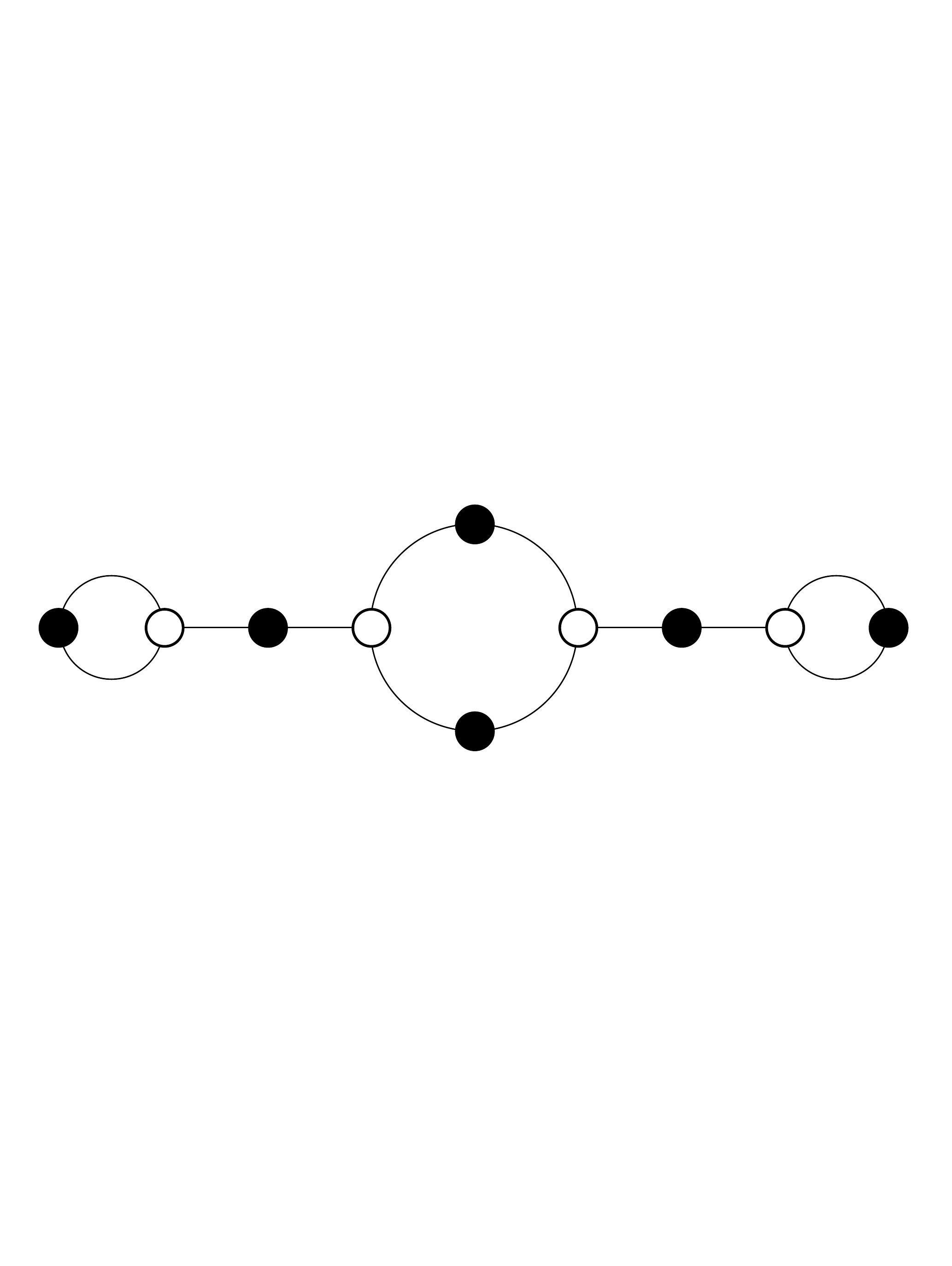}
\par\end{center}

\begin{center}
$\Gamma_{0}\left(8\right)$\\
\scriptsize $8,2,1,1$ \scriptsize
\par\end{center}%
\end{minipage}%
\begin{minipage}[t]{0.25\textwidth}%
\begin{center}
\includegraphics[scale=0.15]{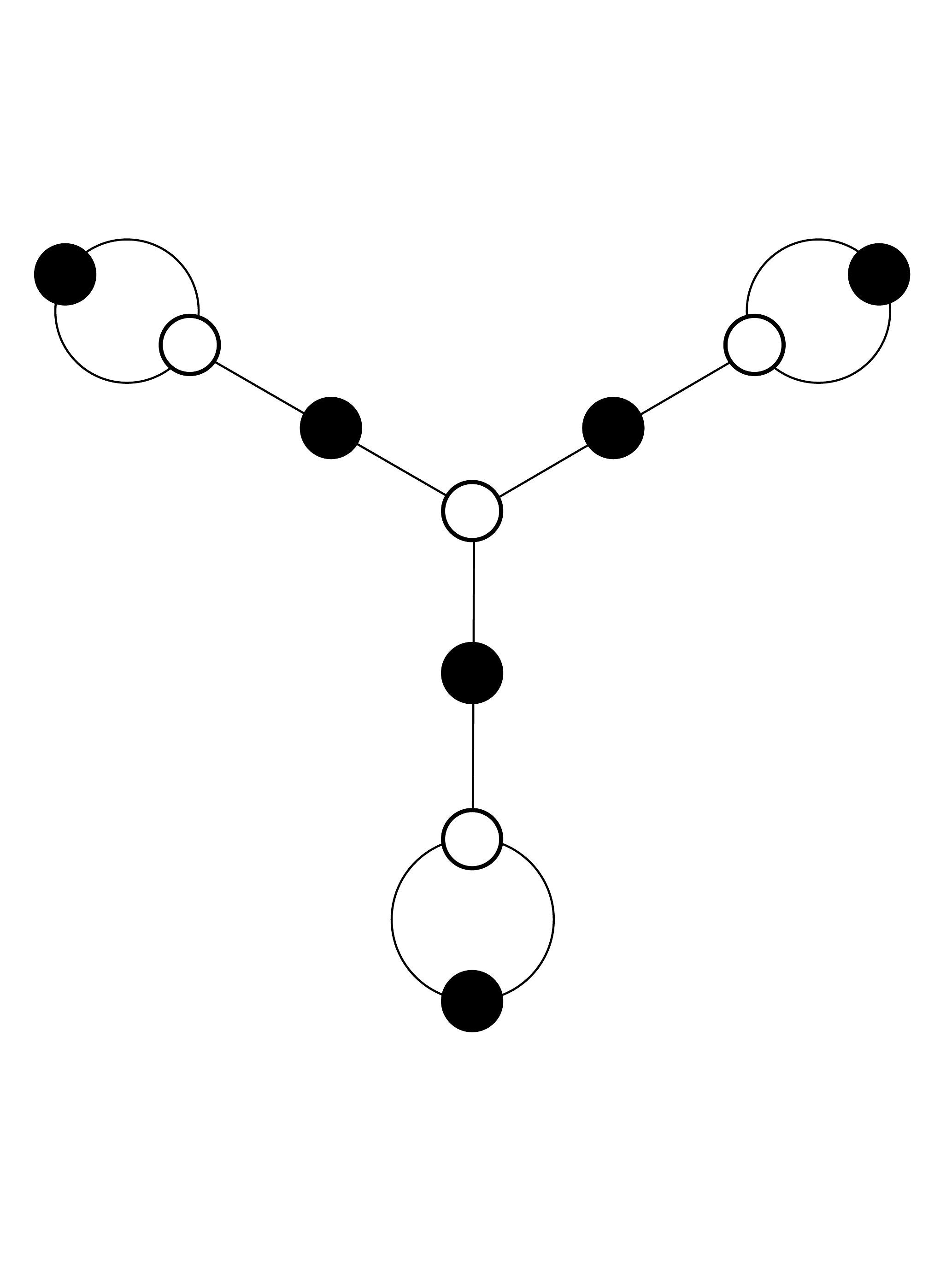}
\par\end{center}

\begin{center}
$\Gamma_{0}\left(9\right)$\\
\scriptsize $9,1,1,1$ \scriptsize
\par\end{center}%
\end{minipage}
\par\end{center}

\begin{center}
\begin{minipage}[t]{0.33\textwidth}%
\begin{center}
\includegraphics[scale=0.15]{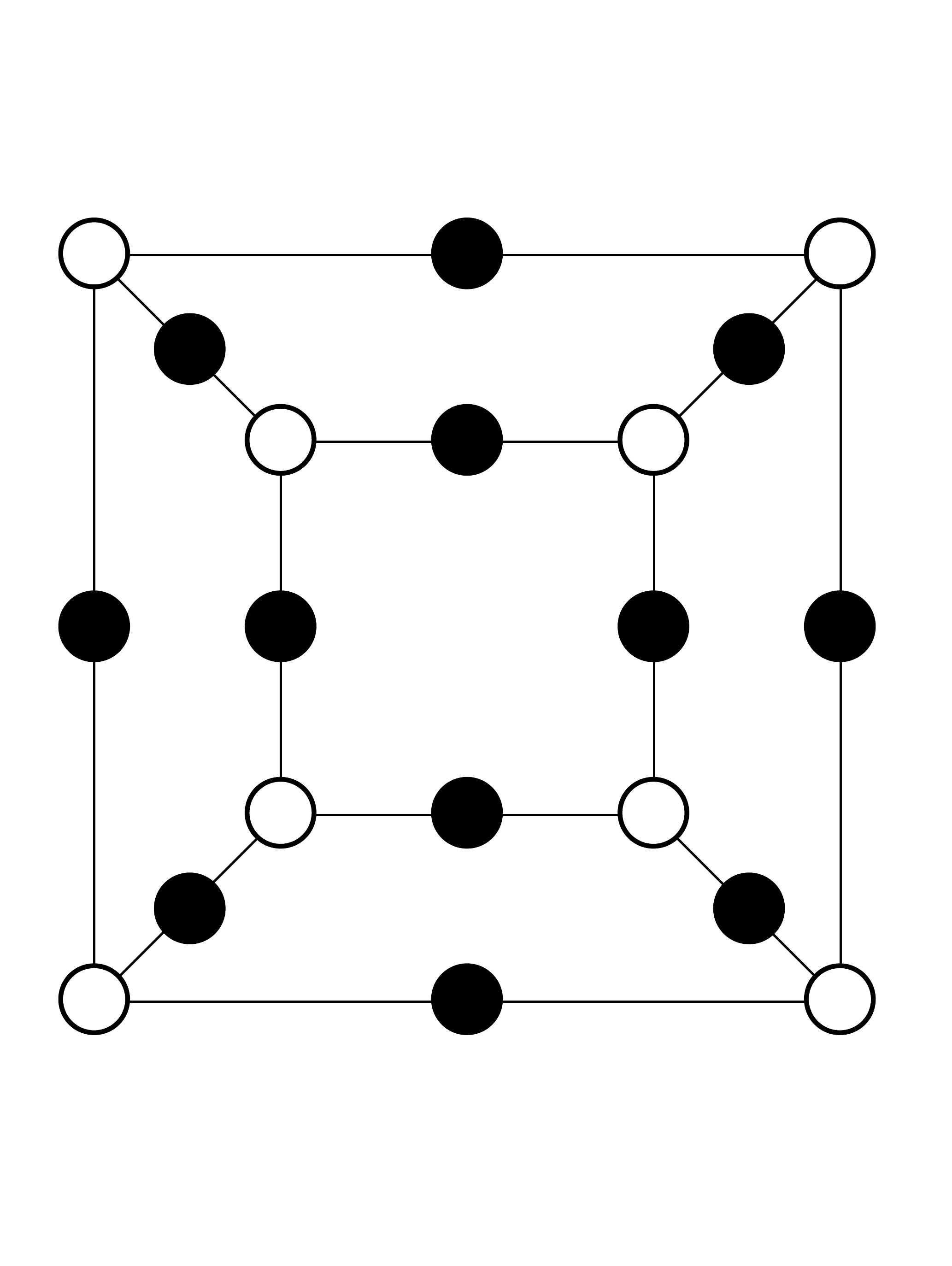}
\par\end{center}

\begin{center}
$\Gamma\left(4\right)$\\
\scriptsize $4,4,4,4,4,4$ \scriptsize
\par\end{center}%
\end{minipage}%
\begin{minipage}[t]{0.33\textwidth}%
\begin{center}
\includegraphics[scale=0.15]{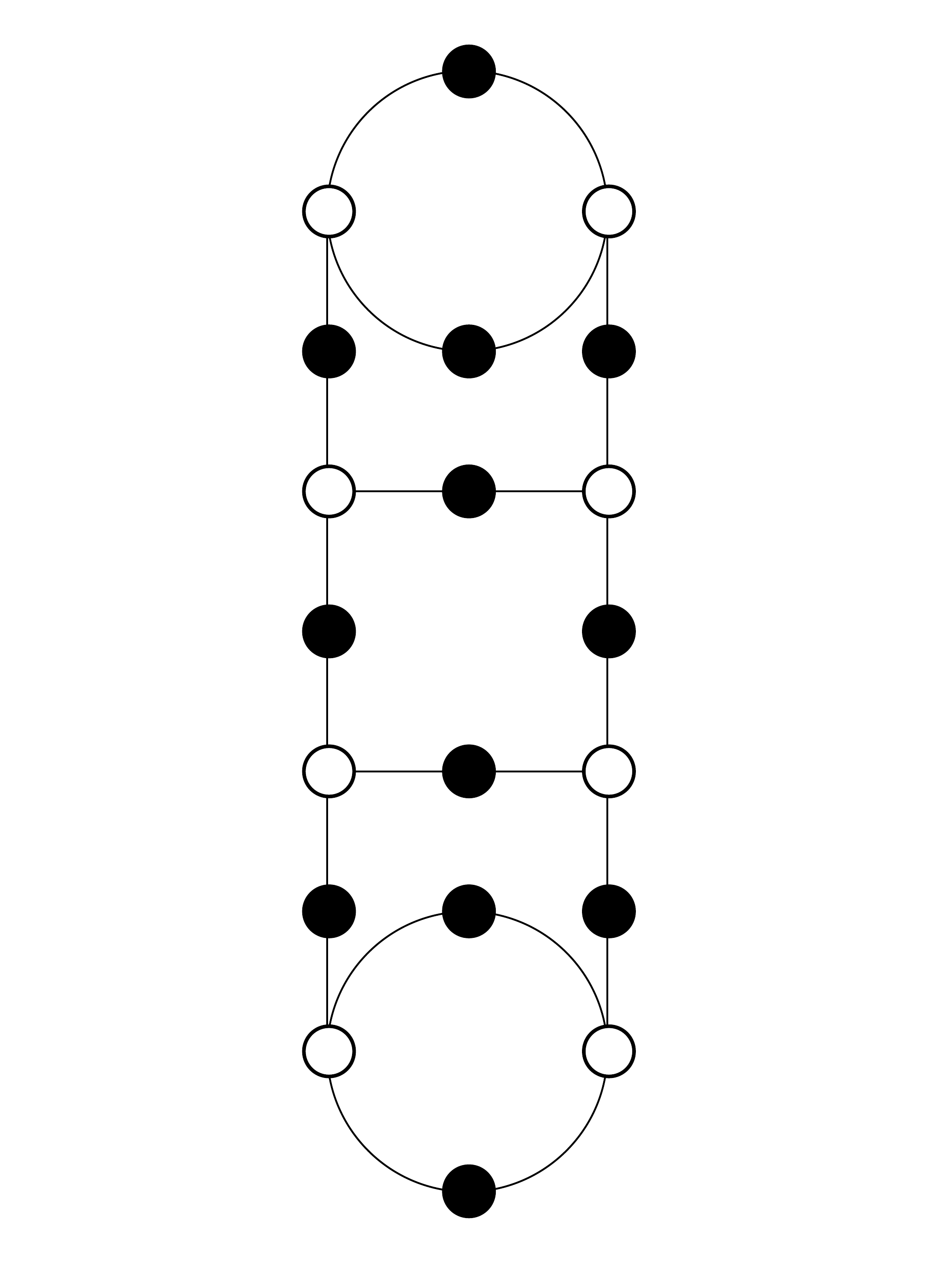}
\par\end{center}

\begin{center}
$\Gamma\left(8;4,1,2\right)$\\
\scriptsize $8,4,4,4,2,2$ \scriptsize
\par\end{center}%
\end{minipage}%
\begin{minipage}[t]{0.33\textwidth}%
\begin{center}
\includegraphics[scale=0.15]{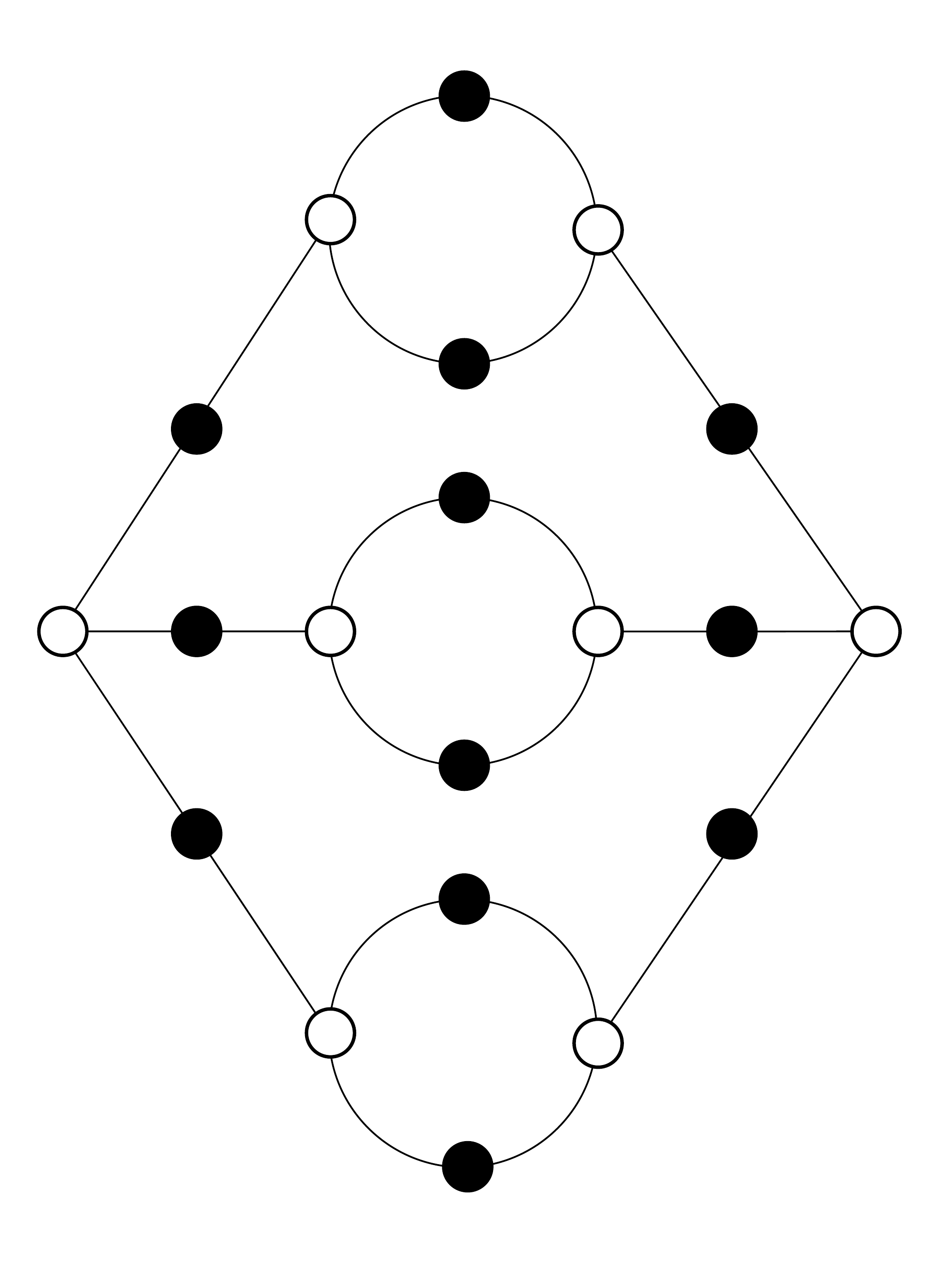}
\par\end{center}

\begin{center}
$\Gamma_{0}\left(3\right)\cap\Gamma\left(2\right)$\\
\scriptsize $6,6,6,2,2,2$ \scriptsize
\par\end{center}%
\end{minipage}
\par\end{center}

\begin{center}
\begin{minipage}[t]{0.33\textwidth}%
\begin{center}
\includegraphics[scale=0.15]{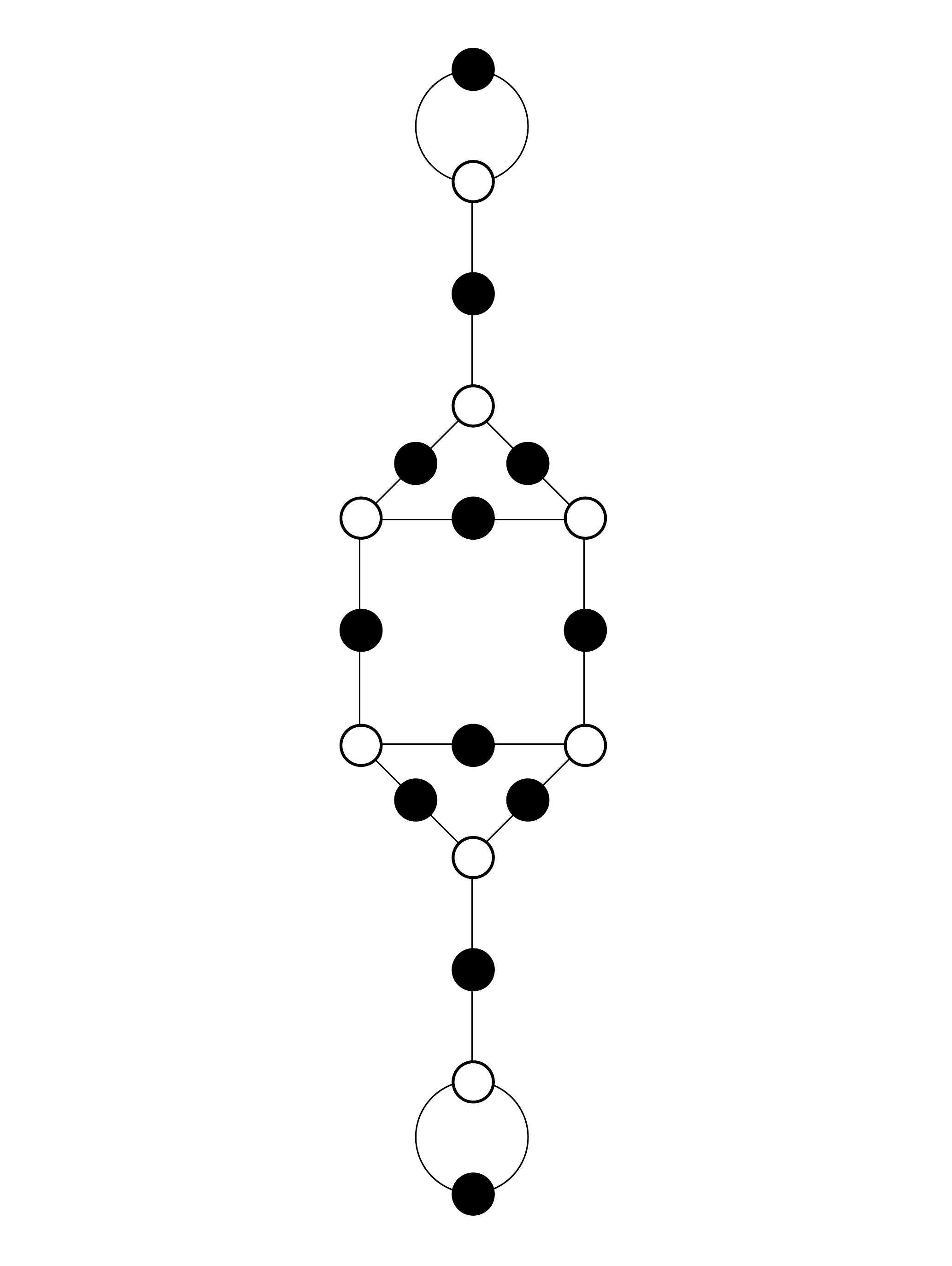}
\par\end{center}

\begin{center}
$\Gamma_{0}\left(12\right)$\\
\scriptsize $12,4,3,3,1,1$ \scriptsize
\par\end{center}%
\end{minipage}%
\begin{minipage}[t]{0.33\textwidth}%
\begin{center}
\includegraphics[scale=0.15]{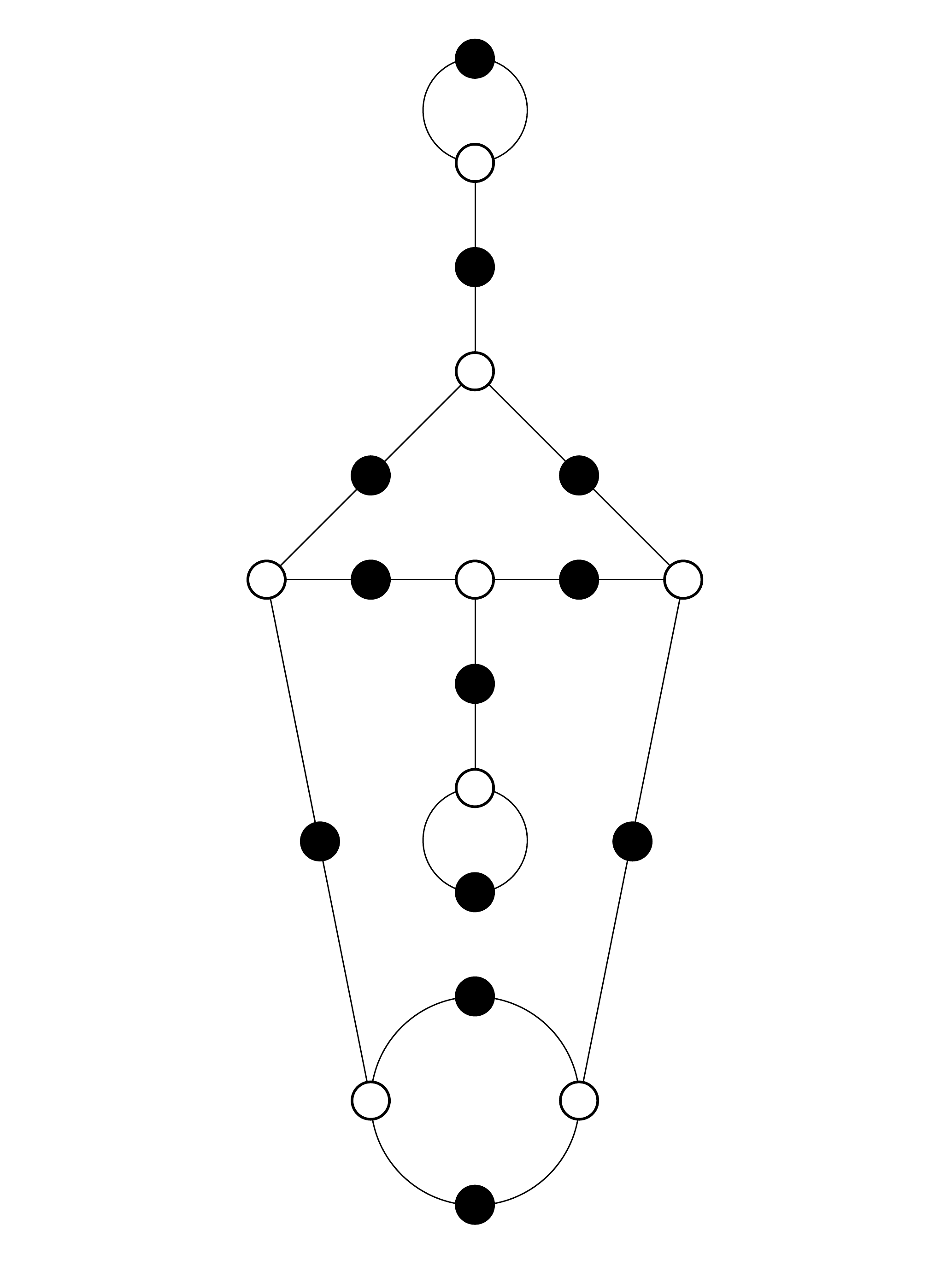}
\par\end{center}

\begin{center}
$\Gamma_{1}\left(8\right)$\\
\scriptsize $8,8,4,2,1,1$ \scriptsize
\par\end{center}%
\end{minipage}%
\begin{minipage}[t]{0.33\textwidth}%
\begin{center}
\includegraphics[scale=0.15]{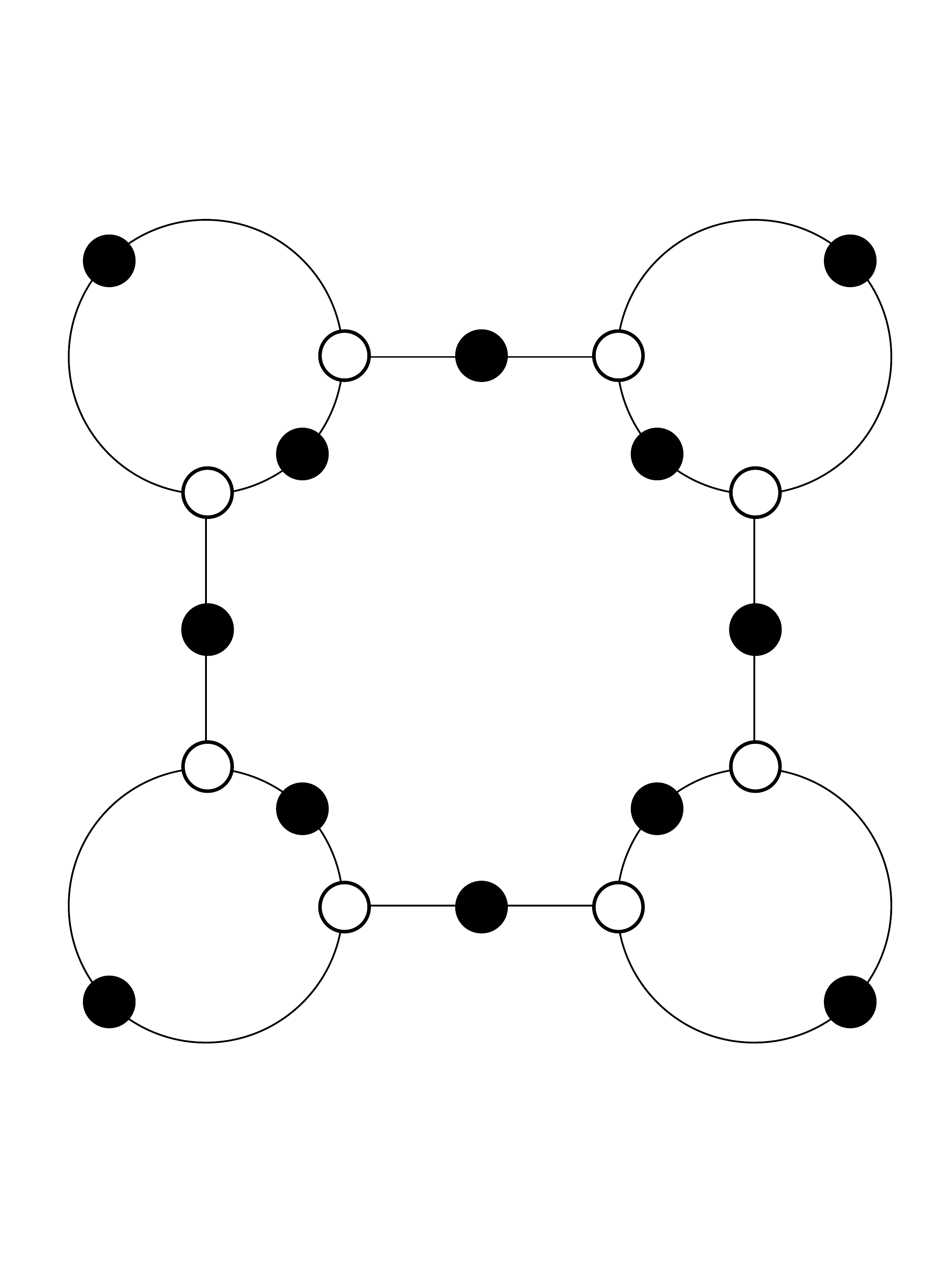}
\par\end{center}

\begin{center}
$\Gamma_{0}\left(8\right)\cap\Gamma\left(2\right)$\\
\scriptsize $8,8,2,2,2,2$ \scriptsize
\par\end{center}%
\end{minipage}
\par\end{center}

\begin{center}
\begin{minipage}[t]{0.33\textwidth}%
\begin{center}
\includegraphics[scale=0.15]{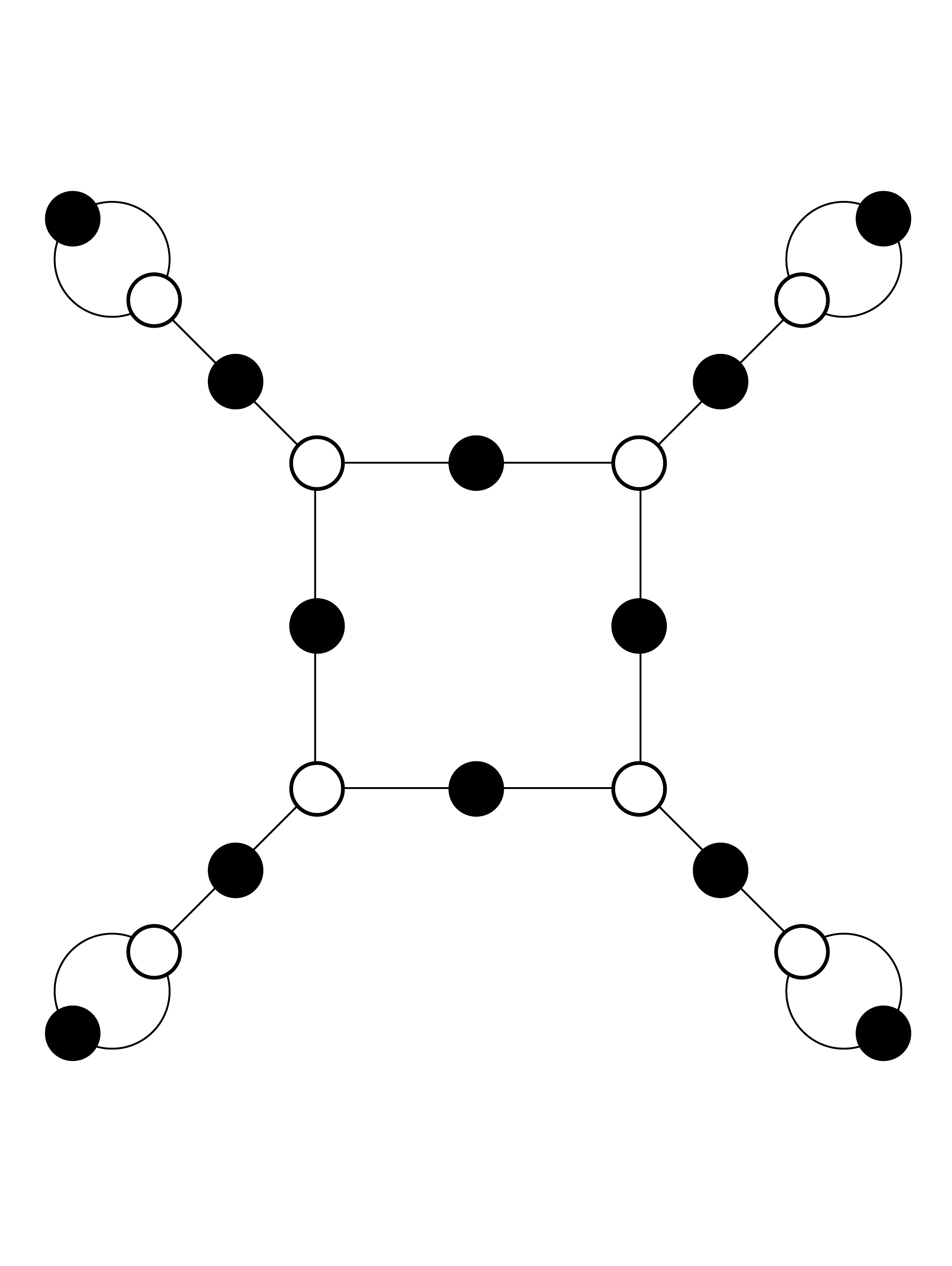}
\par\end{center}

\begin{center}
$\Gamma_{0}\left(16\right)$\\
\scriptsize $16,4,1,1,1,1$ \scriptsize
\par\end{center}%
\end{minipage}%
\begin{minipage}[t]{0.33\textwidth}%
\begin{center}
\includegraphics[scale=0.15]{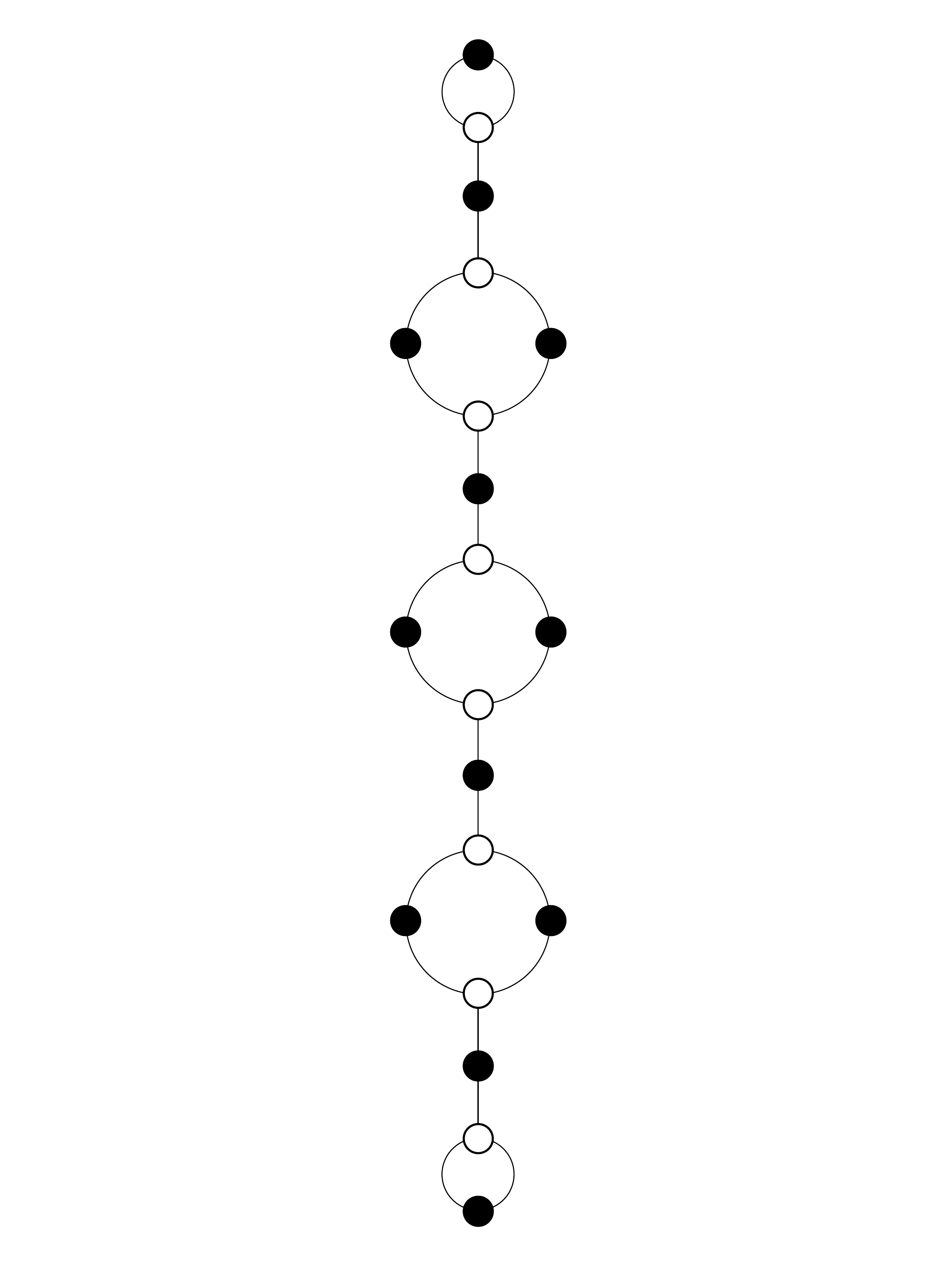}
\par\end{center}

\begin{center}
$\Gamma\left(16;16,2,2\right)$\\
\scriptsize $16,2,2,2,1,1$ \scriptsize
\par\end{center}%
\end{minipage}%
\begin{minipage}[t]{0.33\textwidth}%
\begin{center}
\includegraphics[scale=0.15]{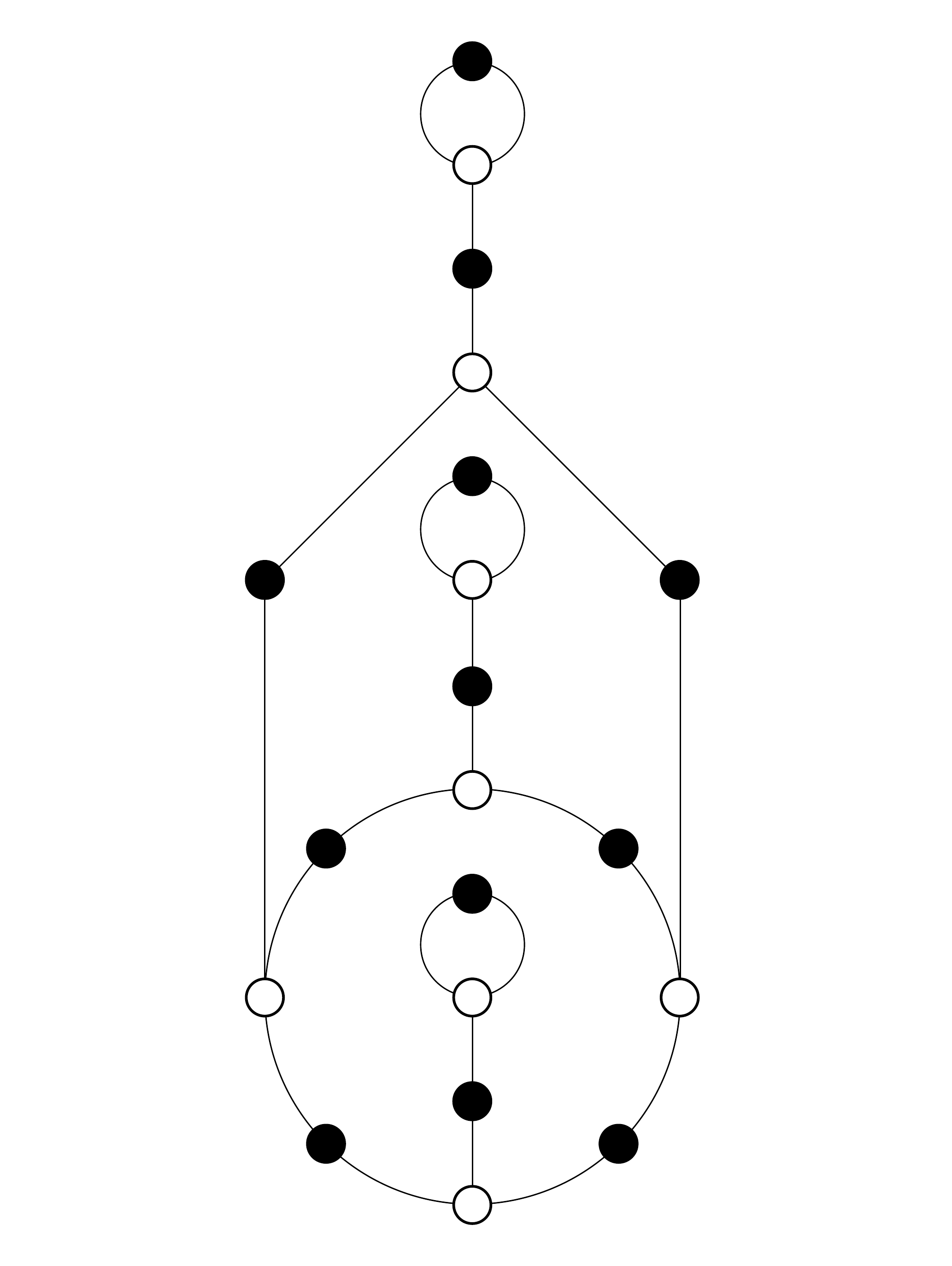}
\par\end{center}

\begin{center}
$\Gamma_{1}\left(7\right)$\\
\scriptsize $7,7,7,1,1,1$ \scriptsize
\par\end{center}%
\end{minipage}
\par\end{center}

\begin{center}
\begin{minipage}[t]{0.33\textwidth}%
\begin{center}
\includegraphics[scale=0.2]{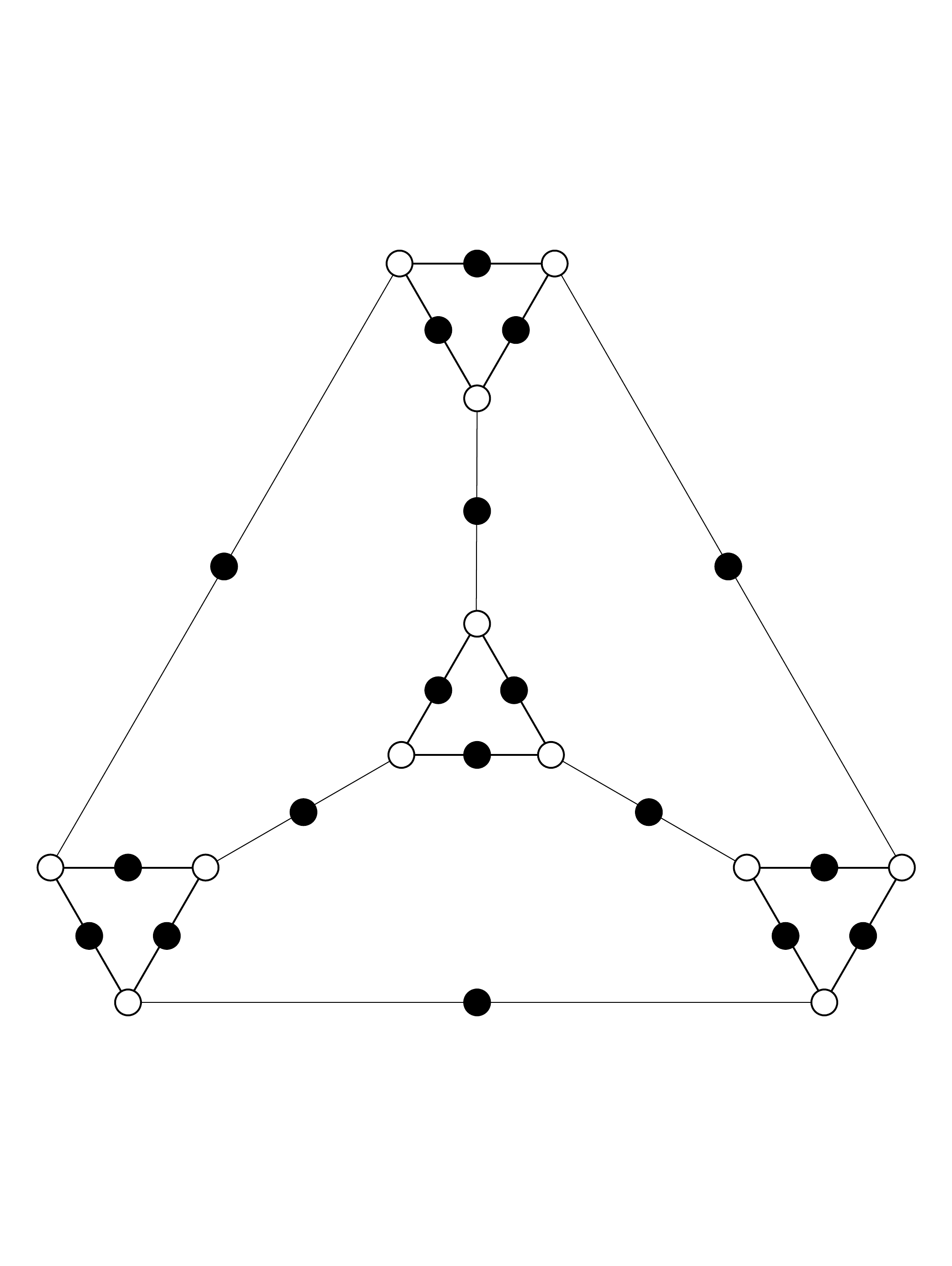}
\par\end{center}

\begin{center}
$\Gamma_{0}\left(2\right)\cap\Gamma\left(3\right)$\\
\scriptsize $6,6,6,6,3,3,3,3$ \scriptsize
\par\end{center}%
\end{minipage}%
\begin{minipage}[t]{0.33\textwidth}%
\begin{center}
\includegraphics[scale=0.2]{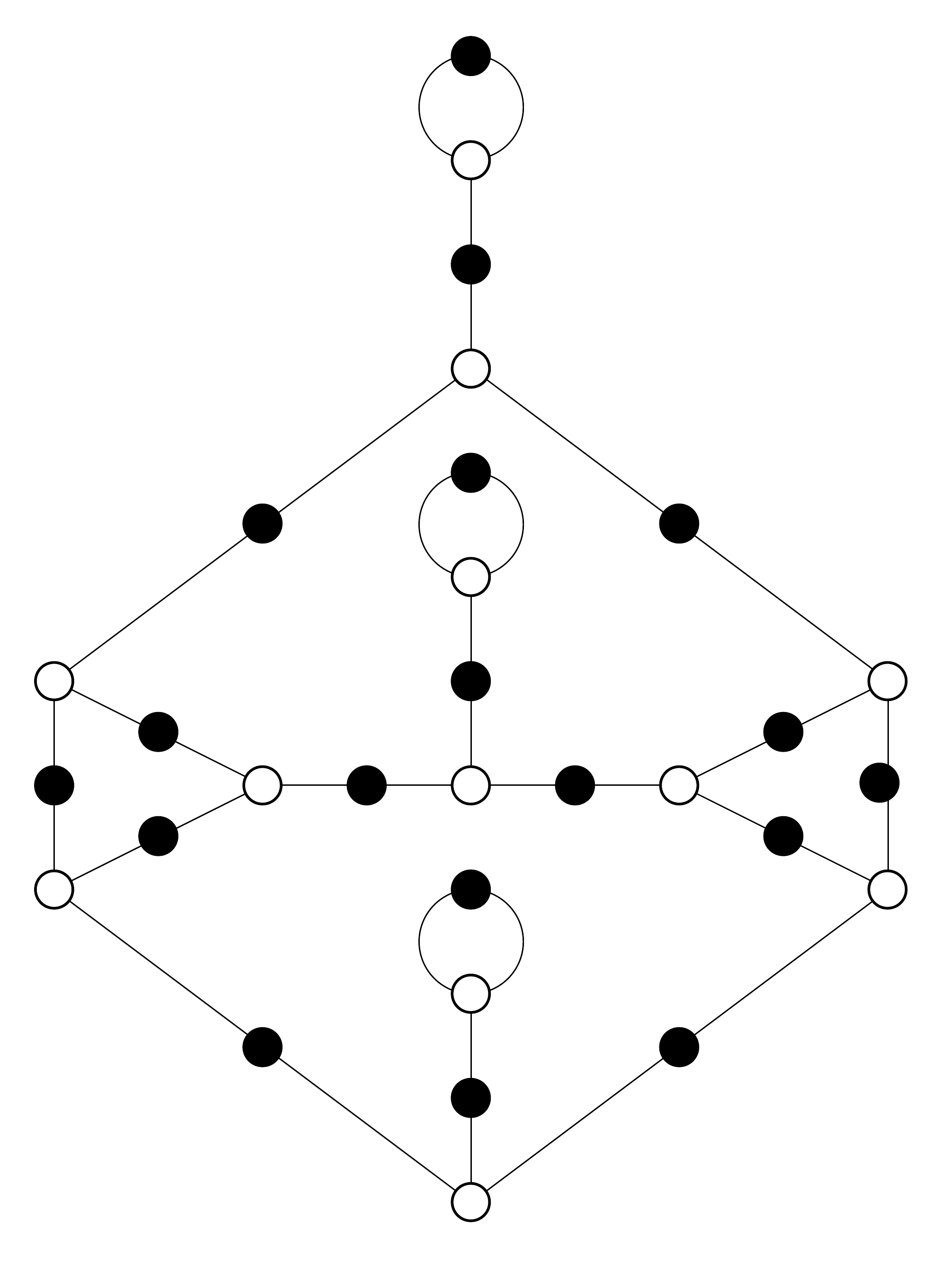}
\par\end{center}

\begin{center}
$\Gamma_{1}\left(9\right)$\\
\scriptsize $9,9,9,3,3,1,1,1$ \scriptsize
\par\end{center}%
\end{minipage}%
\begin{minipage}[t]{0.33\textwidth}%
\begin{center}
\includegraphics[scale=0.2]{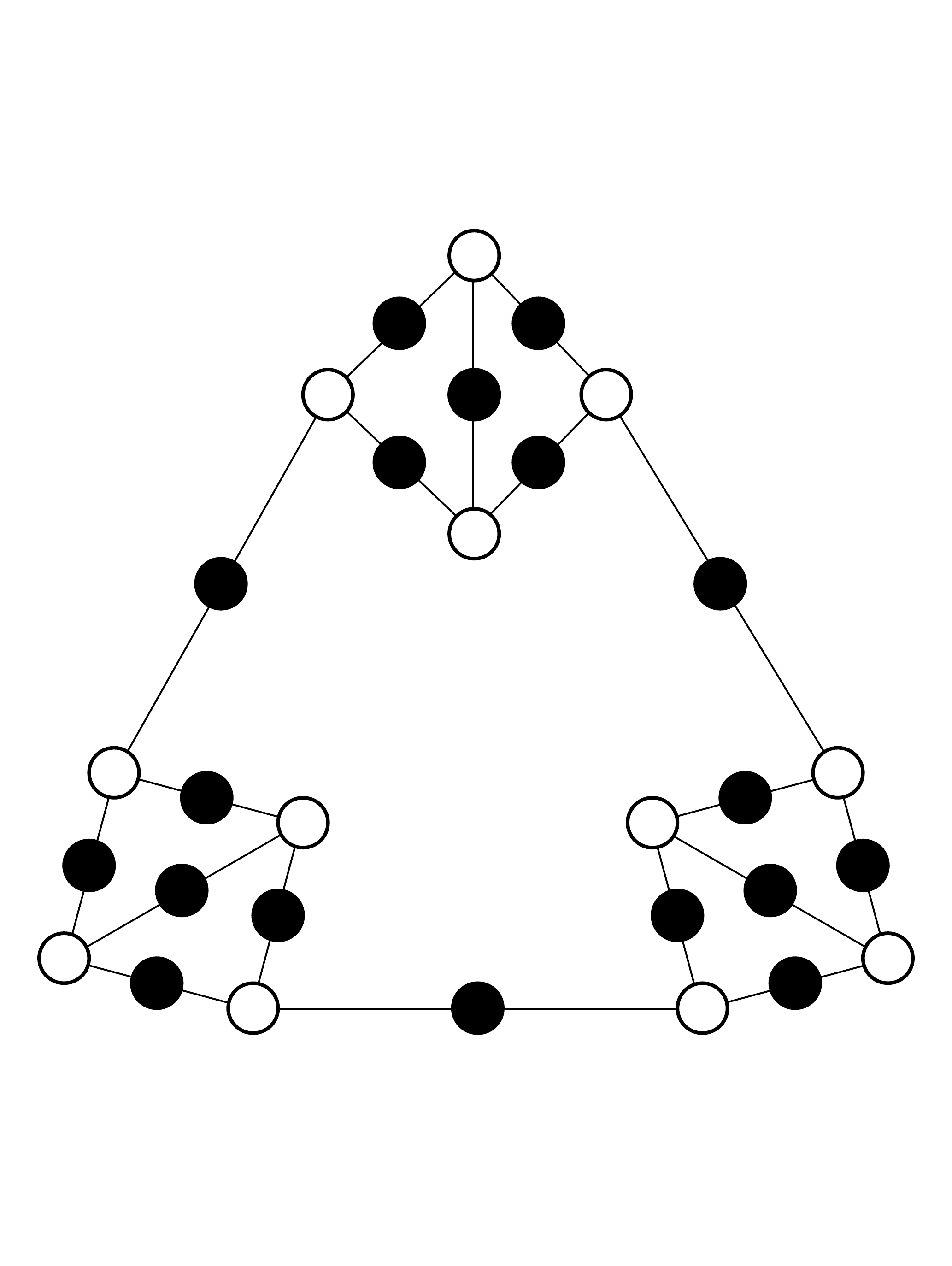}
\par\end{center}

\begin{center}
$\Gamma\left(9;3,1,3\right)$\\
\scriptsize $9,9,3,3,3,3,3,3$ \scriptsize
\par\end{center}%
\end{minipage}
\par\end{center}

\begin{center}
\begin{minipage}[t]{0.33\textwidth}%
\begin{center}
\includegraphics[scale=0.2]{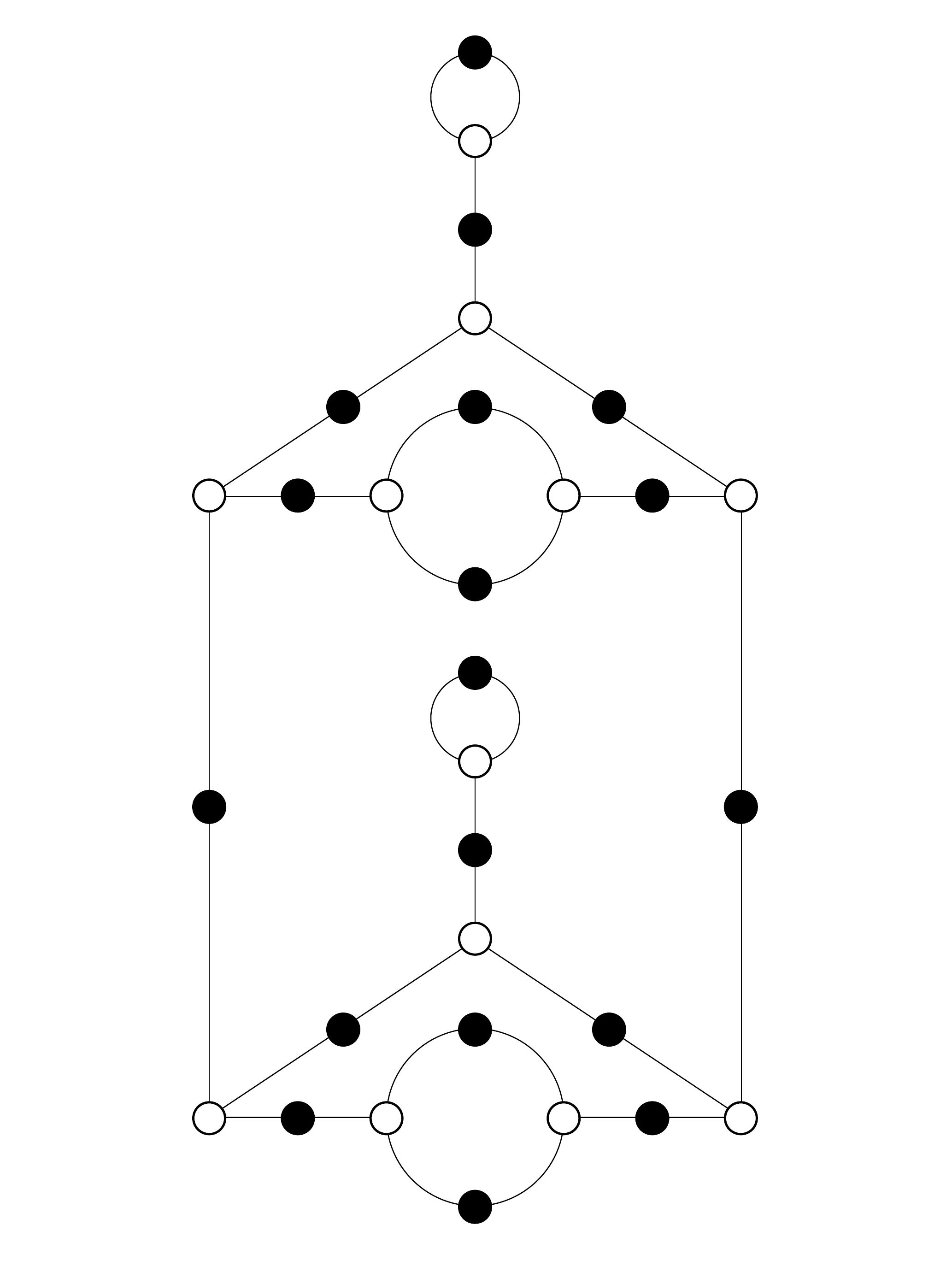}
\par\end{center}

\begin{center}
$\Gamma_{1}\left(10\right)$\\
\scriptsize $10,10,5,5,2,2,1,1$ \scriptsize
\par\end{center}%
\end{minipage}%
\begin{minipage}[t]{0.33\textwidth}%
\begin{center}
\includegraphics[scale=0.2]{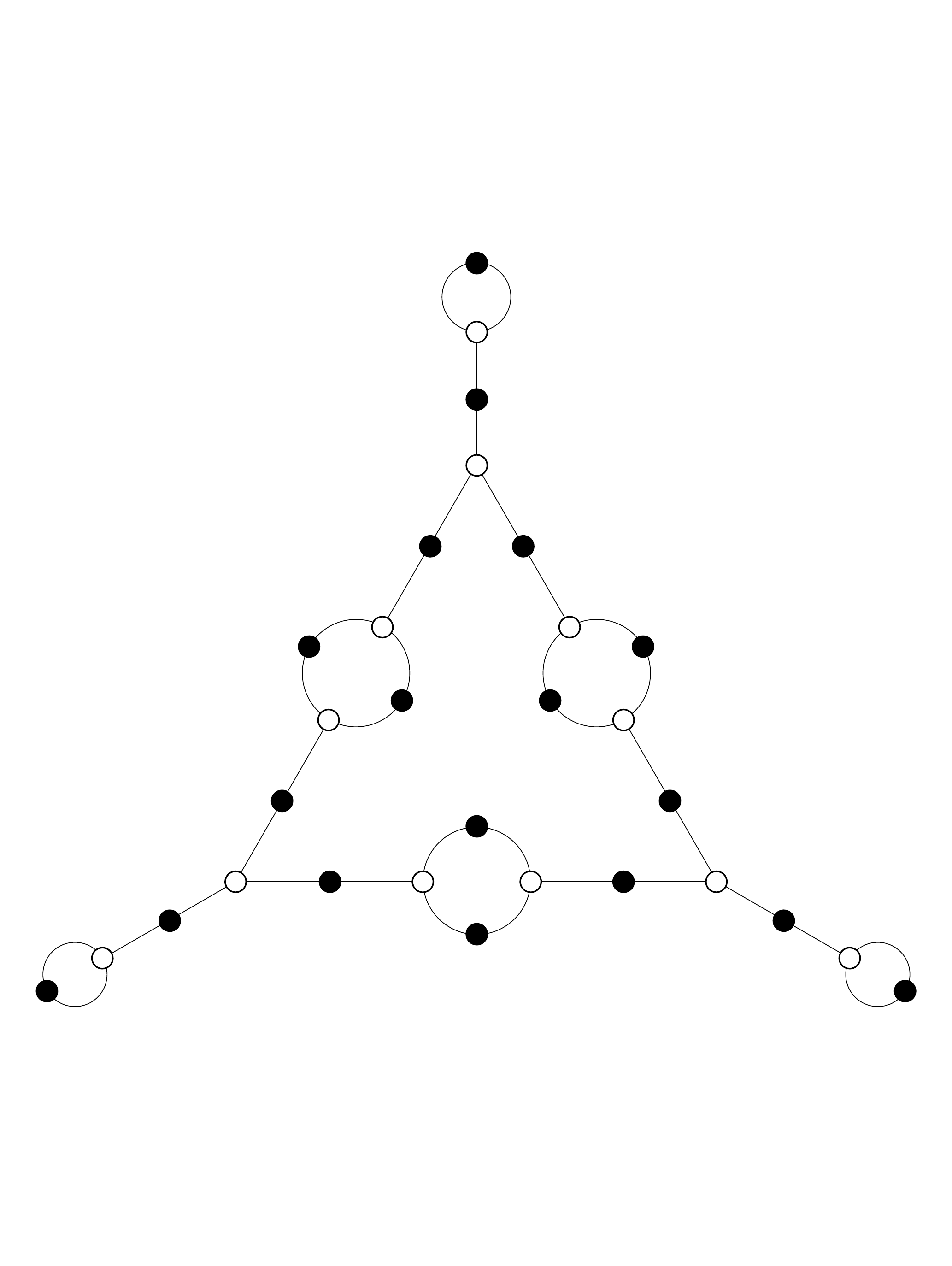}
\par\end{center}

\begin{center}
$\Gamma_{0}\left(18\right)$\\
\scriptsize $18,9,2,2,2,1,1,1$ \scriptsize
\par\end{center}%
\end{minipage}%
\begin{minipage}[t]{0.33\textwidth}%
\begin{center}
\includegraphics[scale=0.2]{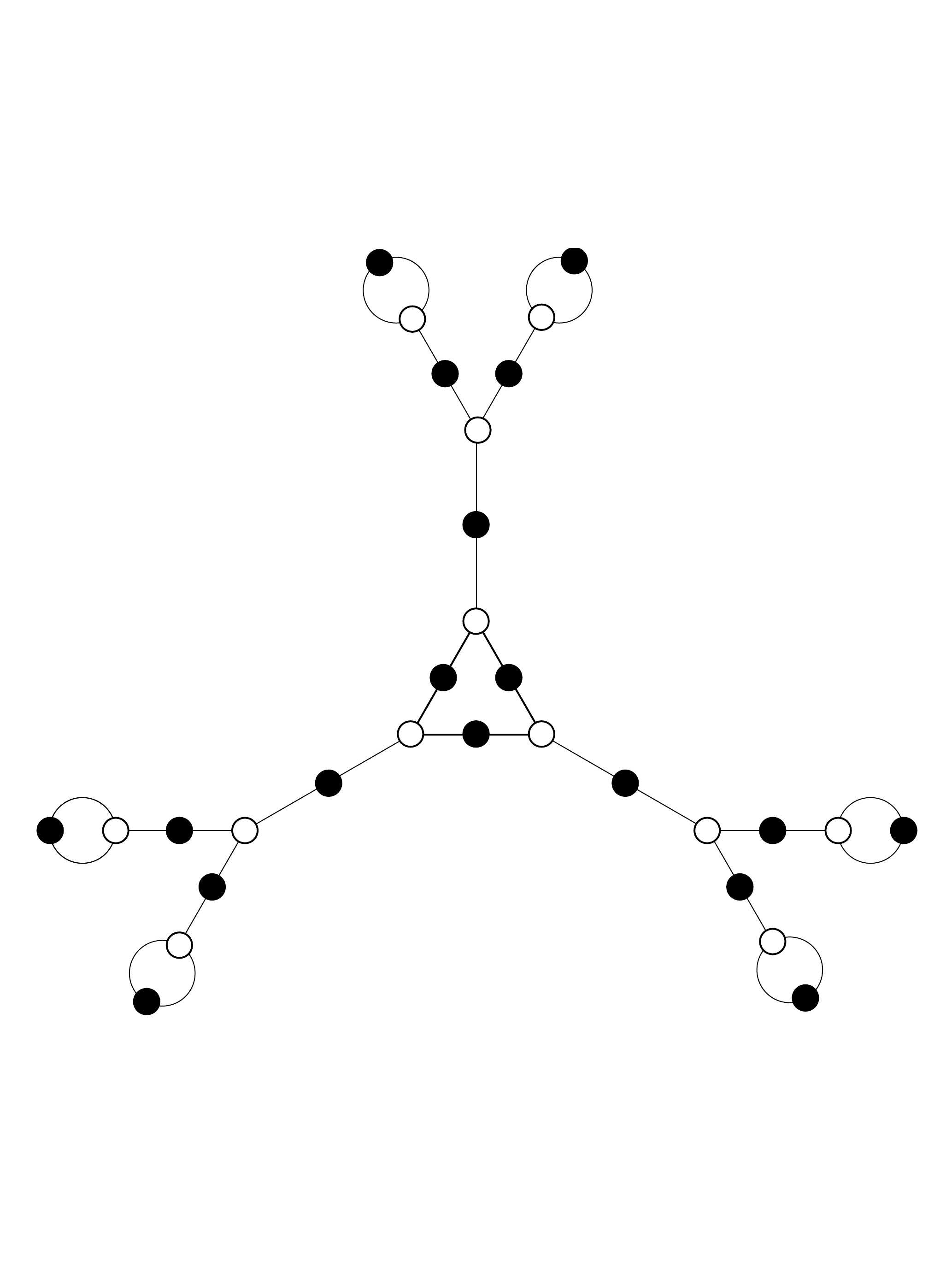}
\par\end{center}

\begin{center}
$\Gamma\left(27;27,3,3\right)$\\
\scriptsize $27,3,1,1,1,1,1,1$ \scriptsize
\par\end{center}%
\end{minipage}
\par\end{center}

\begin{center}
\begin{minipage}[t]{0.33\textwidth}%
\begin{center}
\includegraphics[scale=0.2]{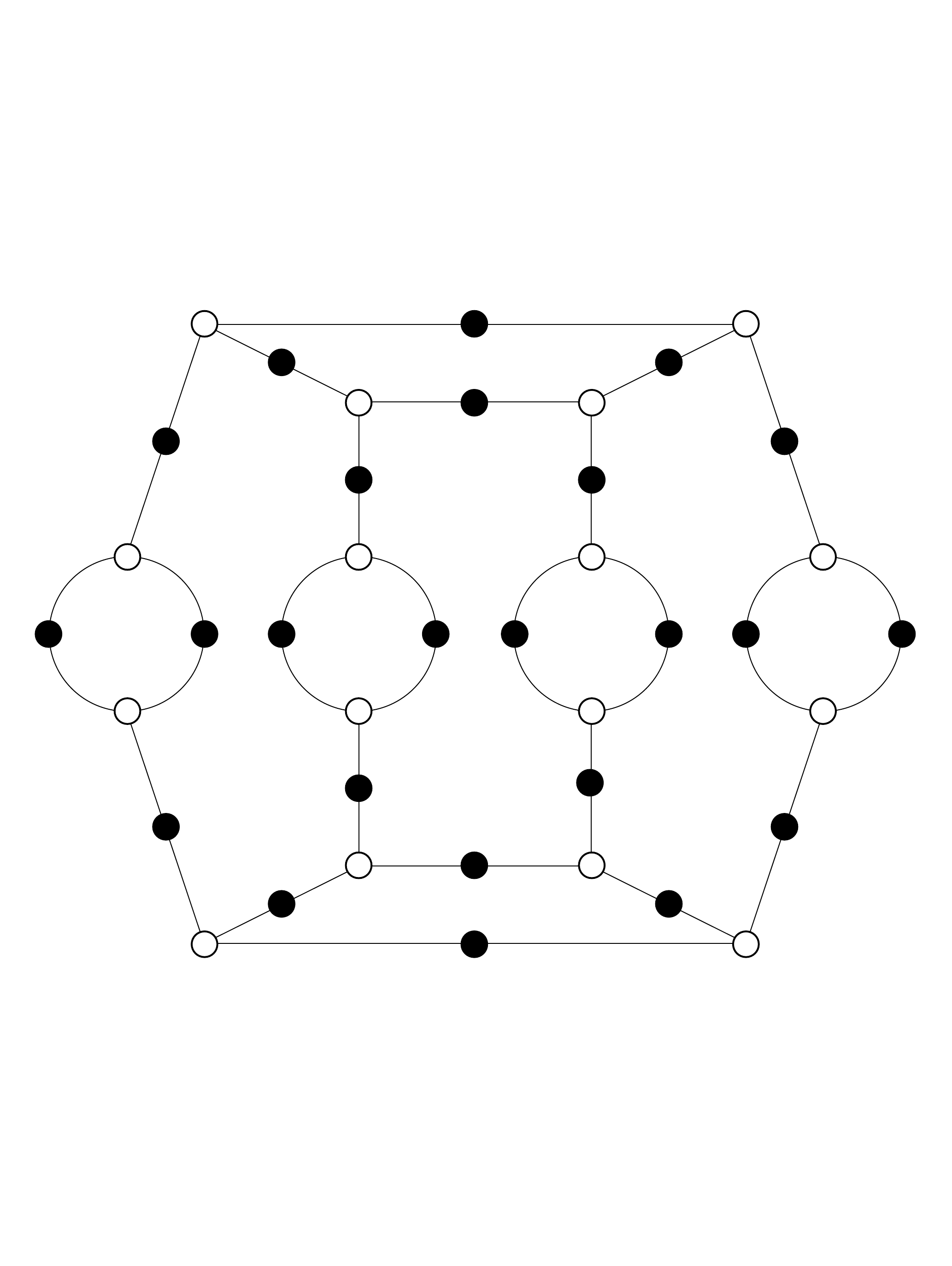}
\par\end{center}

\begin{center}
$\Gamma_{1}\left(8\right)\cap\Gamma\left(2\right)$\\
\scriptsize $8,8,8,8,4,4,2,2,2,2$ \scriptsize
\par\end{center}%
\end{minipage}%
\begin{minipage}[t]{0.33\textwidth}%
\begin{center}
\includegraphics[scale=0.2]{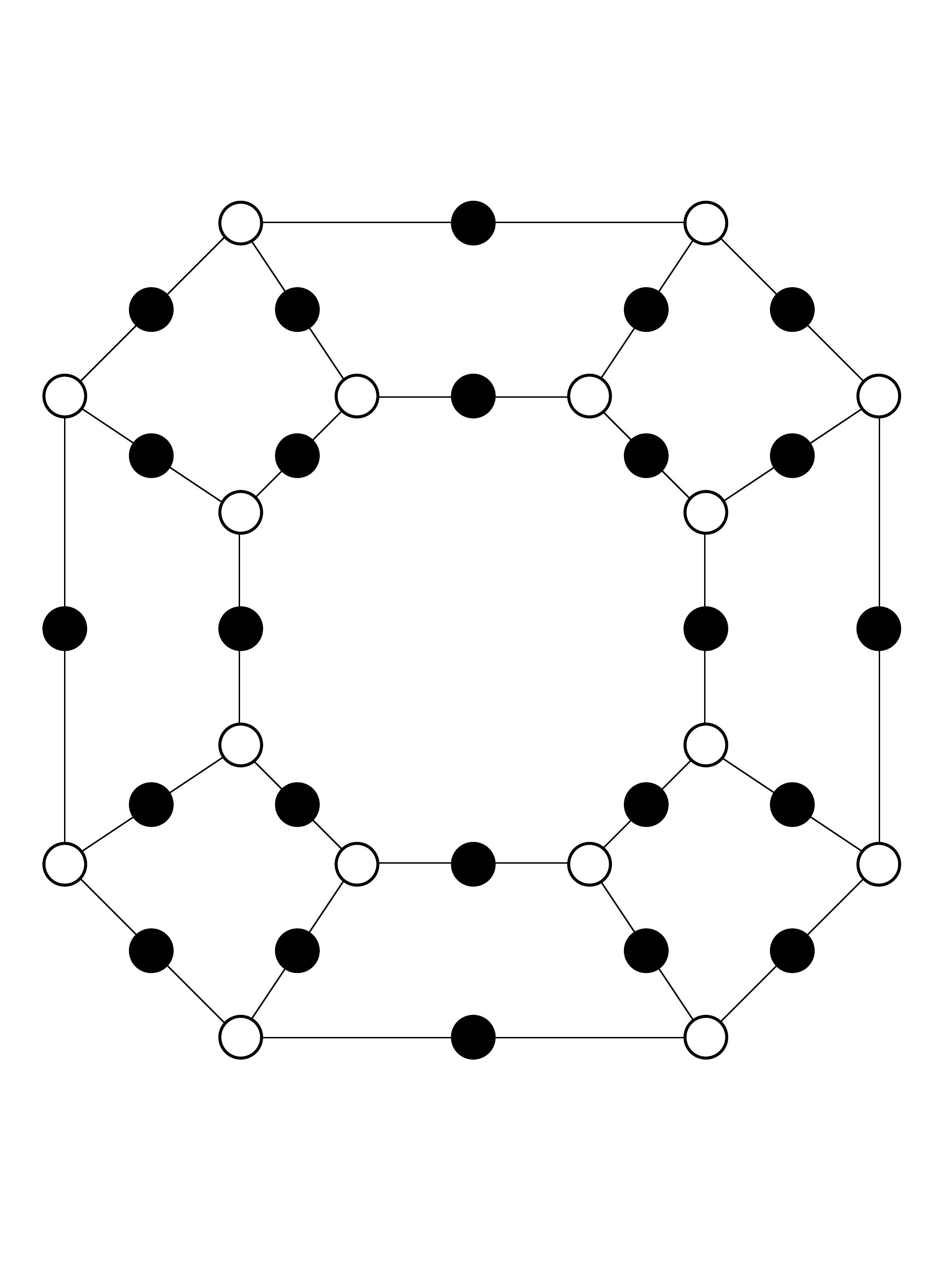}
\par\end{center}

\begin{center}
$\Gamma\left(8;2,1,2\right)$\\
\scriptsize $8,8,4,4,4,4,4,4,4,4$ \scriptsize
\par\end{center}%
\end{minipage}%
\begin{minipage}[t]{0.33\textwidth}%
\begin{center}
\includegraphics[scale=0.2]{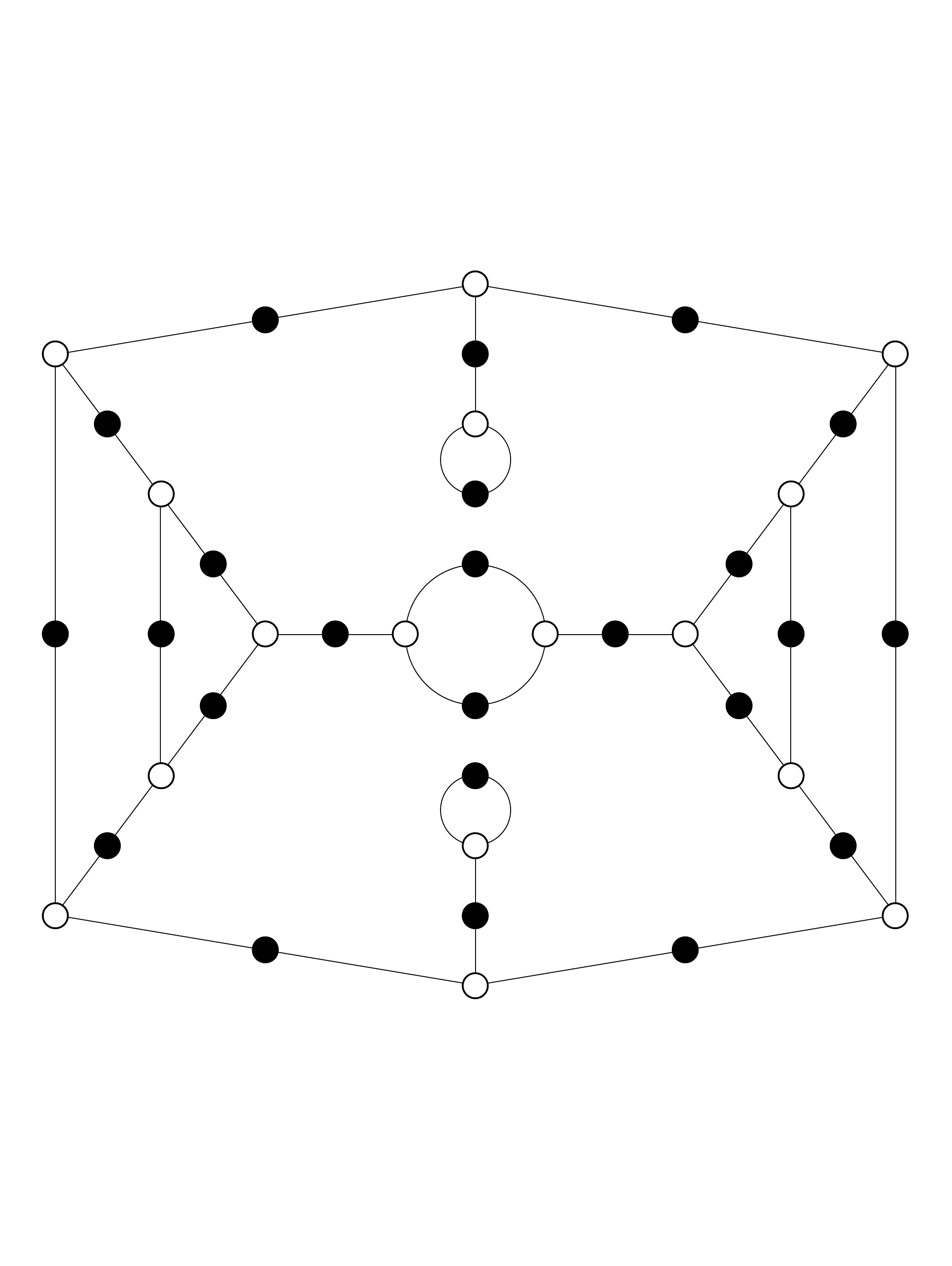}
\par\end{center}

\begin{center}
$\Gamma_{1}\left(12\right)$\\
\scriptsize $12,12,6,4,4,3,3,2,1,1$ \scriptsize
\par\end{center}%
\end{minipage}
\par\end{center}

\begin{center}
\begin{minipage}[t]{0.33\textwidth}%
\begin{center}
\includegraphics[scale=0.2]{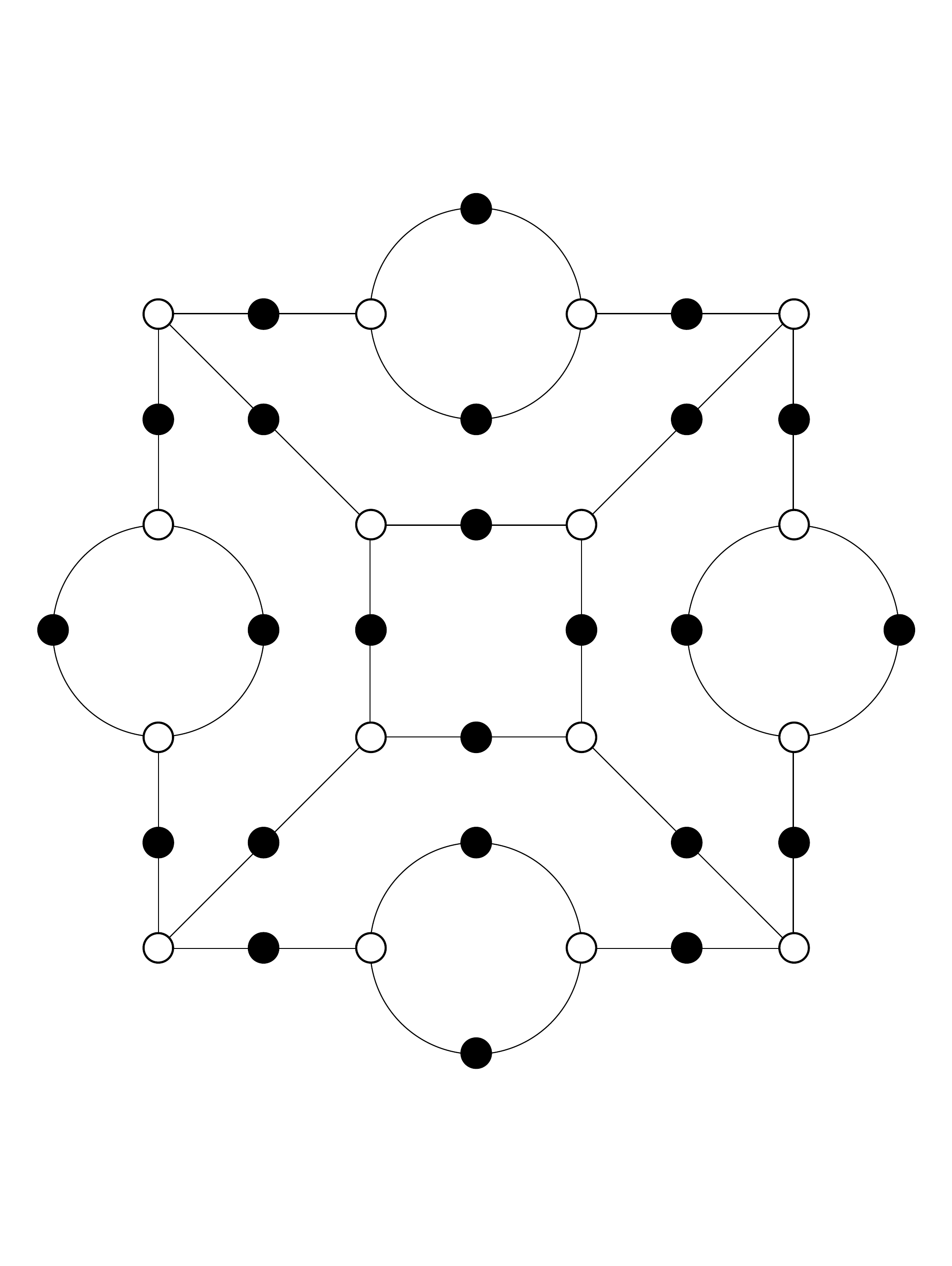}
\par\end{center}

\begin{center}
$\Gamma\left(12;6,1,2\right)$\\
\scriptsize $12,6,6,6,6,4,2,2,2,2$ \scriptsize
\par\end{center}%
\end{minipage}%
\begin{minipage}[t]{0.33\textwidth}%
\begin{center}
\includegraphics[scale=0.2]{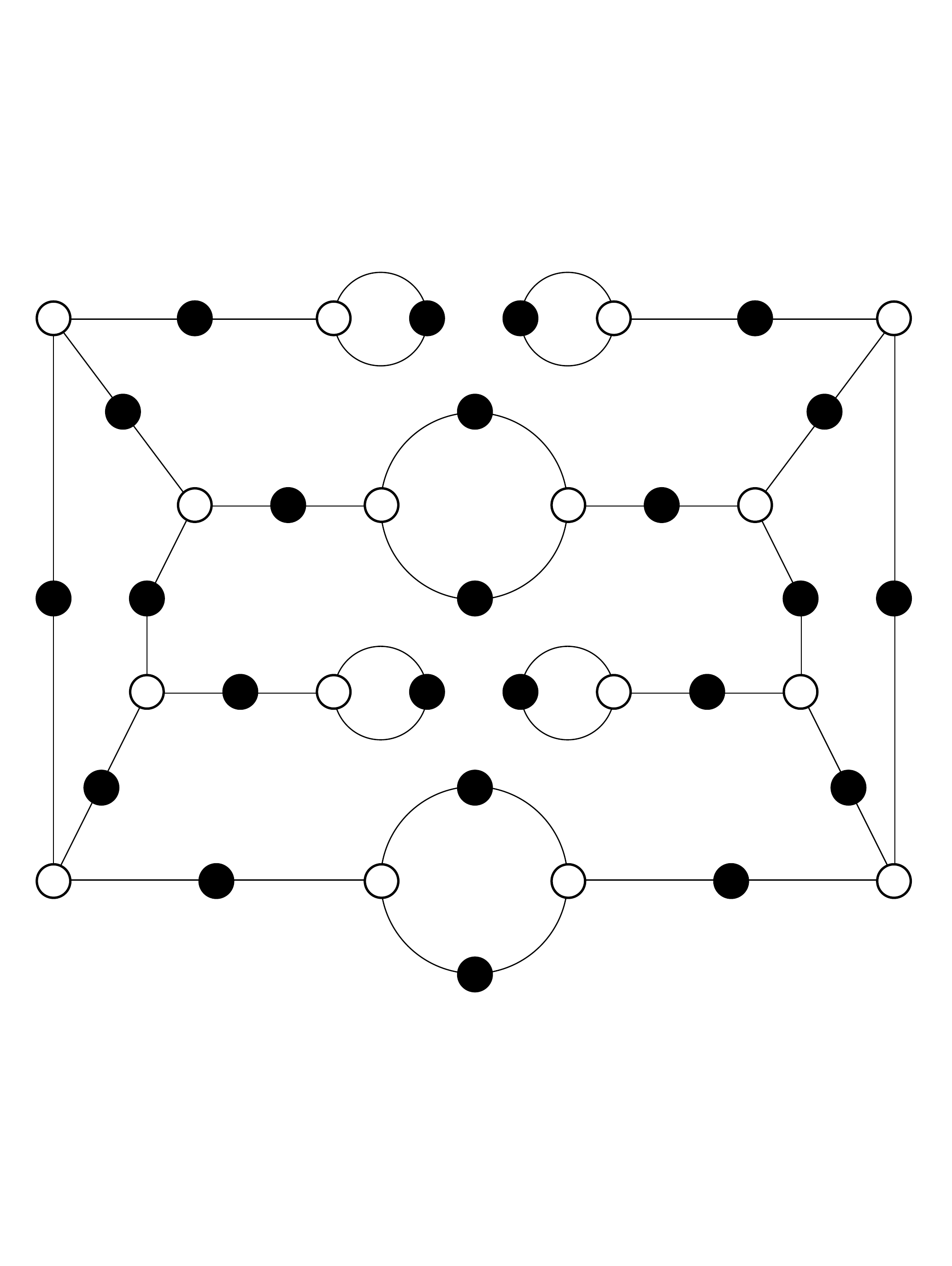}
\par\end{center}

\begin{center}
$\Gamma_{0}\left(16\right)\cap\Gamma_{1}\left(8\right)$\\
\scriptsize $16,16,4,4,2,2,1,1,1,1$ \scriptsize
\par\end{center}%
\end{minipage}%
\begin{minipage}[t]{0.33\textwidth}%
\begin{center}
\includegraphics[scale=0.2]{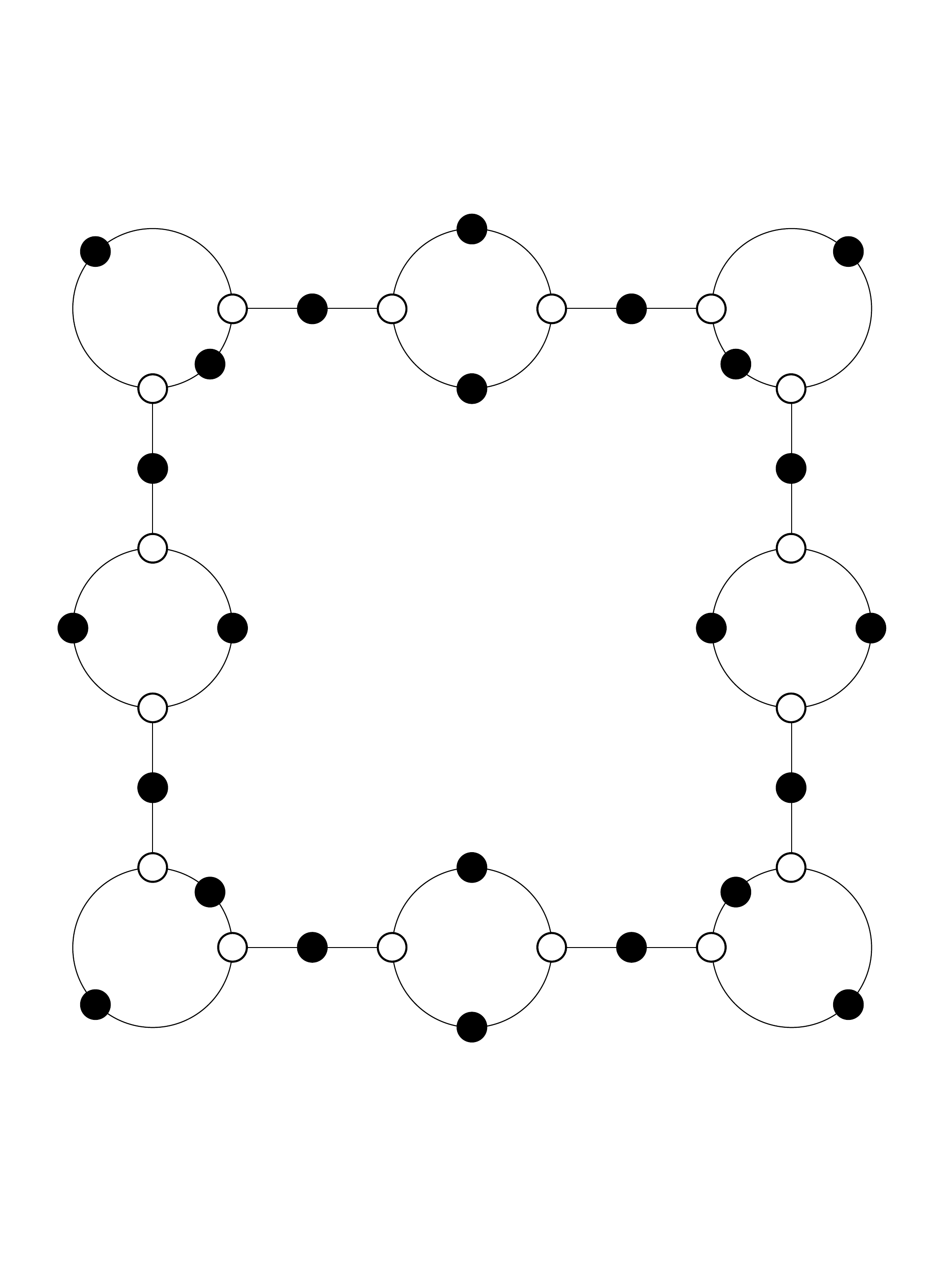}
\par\end{center}

\begin{center}
$\Gamma\left(16;8,2,2\right)$\\
\scriptsize $16,16,2,2,2,2,2,2,2,2$ \scriptsize
\par\end{center}%
\end{minipage}
\par\end{center}

\begin{center}
\begin{minipage}[t]{0.33\textwidth}%
\begin{center}
\includegraphics[scale=0.2]{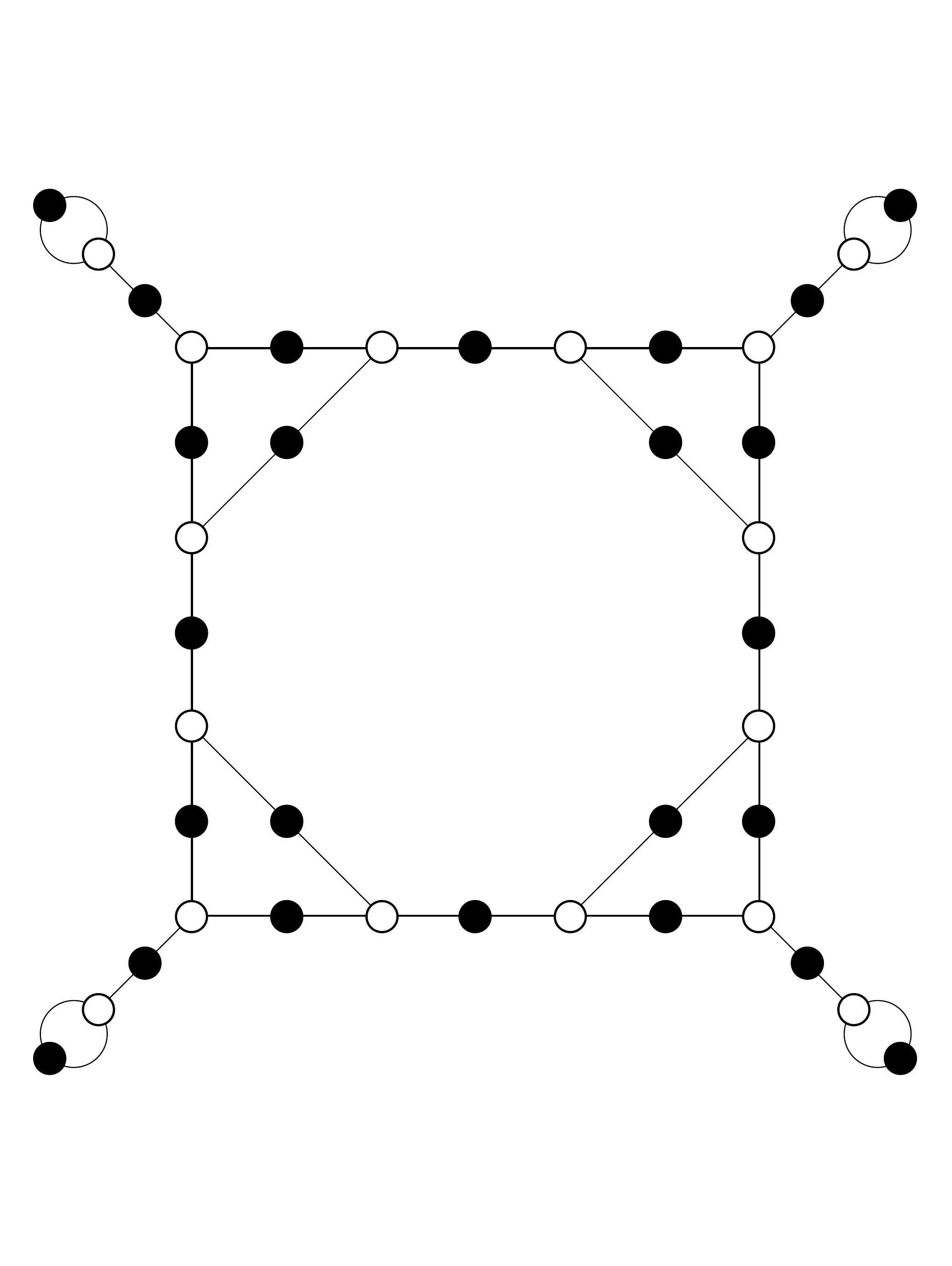}
\par\end{center}

\begin{center}
$\Gamma\left(24;24,2,2\right)$\\
\scriptsize $24,8,3,3,3,3,1,1,1,1$ \scriptsize
\par\end{center}%
\end{minipage}%
\begin{minipage}[t]{0.33\textwidth}%
\begin{center}
\includegraphics[scale=0.2]{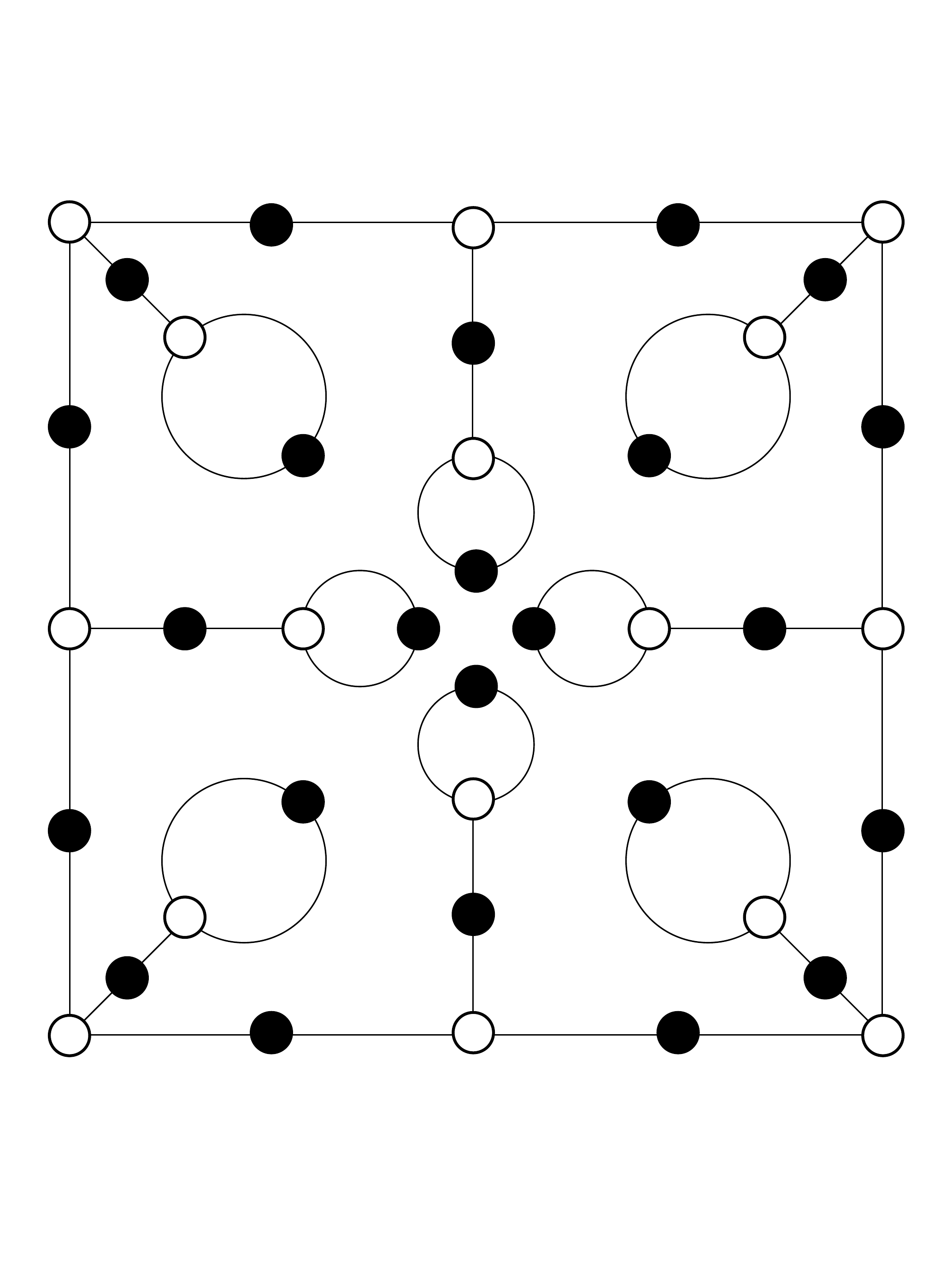}
\par\end{center}

\begin{center}
$\Gamma\left(32;32,4,2\right)$\\
\scriptsize $32,8,1,1,1,1,1,1,1,1$ \scriptsize
\par\end{center}%
\end{minipage}
\par\end{center}

\begin{center}
\begin{minipage}[t]{0.5\textwidth}%
\begin{center}
\includegraphics[scale=0.25]{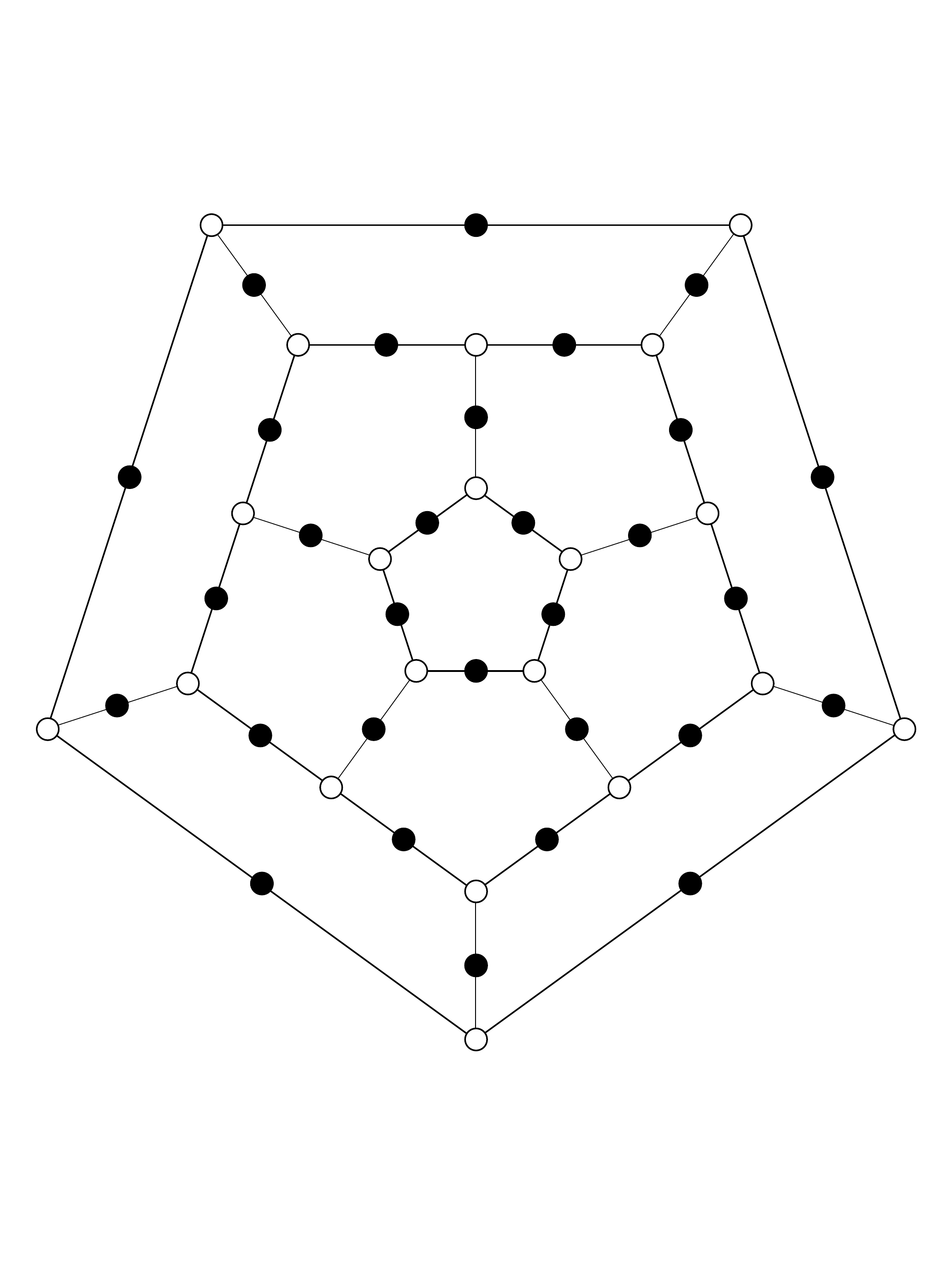}
\par\end{center}

\begin{center}
$\Gamma\left(5\right)$\\
\scriptsize $5,5,5,5,5,5,5,5,5,5,5,5$ \scriptsize
\par\end{center}%
\end{minipage}%
\begin{minipage}[t]{0.5\textwidth}%
\begin{center}
\includegraphics[scale=0.25]{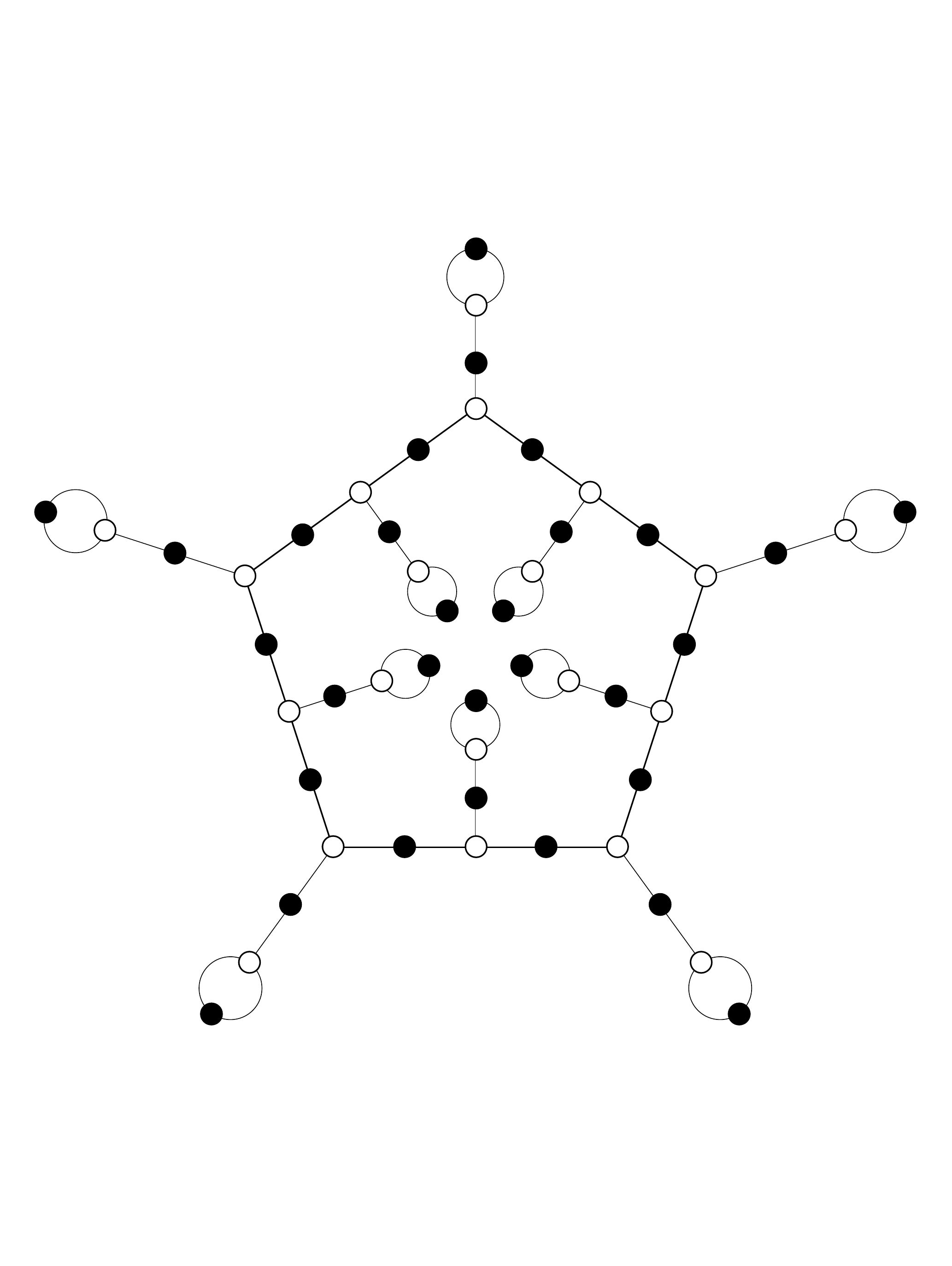}
\par\end{center}

\begin{center}
$\Gamma_{0}\left(25\right)\cap\Gamma_{1}\left(5\right)$\\
\scriptsize $25,25,1,1,1,1,1,1,1,1,1,1$ \scriptsize
\par\end{center}%
\end{minipage}
\par\end{center}


\subsection{Index 36 Subgroups and Calabi-Yau Threefolds}

Since all the index 24 torsion-free, genus zero congruence subgroups
have been associated in \cite{2,TopYui,7} to K3 surfaces 
- in particular elliptically
fibred K3 surfaces over $\mathbb{P}^{1}$ - it is natural to enquire
whether the index 36 groups can be associated with Calabi-Yau threefolds.
We shall see that this is indeed the case, and the threefolds are
a simple generalization of the K3 case, being now elliptic fibrations
over $\mathbb{P}^{2}$.

First, let us recall some standard facts about elliptic Calabi-Yau
threefolds. The classification of base surfaces $B$ over which an
elliptic fibration is a smooth Calabi-Yau threefold was performed
in \cite{11,12}: they can be either Hirzebruch surfaces and their
blowups; del Pezzo surfaces; or the Enriques surface. In particular,
the local Weierstra\ss\ model looks like:
\begin{equation}
y^{2}=x^{3}+f\left(z_{i}\right)x+g\left(z_{i}\right)\subset\mathbb{P}_{\left[x,y\right]}^{2} \ ,
\end{equation}
with
\begin{equation}
f\left(z_{i}\right)\in\Gamma\left(\left(K_{B}^{-1}\right)^{\otimes4}\right)\quad,\quad g\left(z_{i}\right)\in\Gamma\left(\left(K_{B}^{-1}\right)^{\otimes6}\right)
\ .
\end{equation}
This is an elliptic curve in a $\mathbb{P}^{2}$ bundle over $B$,
whose coordinates are $z_{i}$; the $\left(x,y\right)$ are affine
coordinates on the $\mathbb{P}^{2}$. The functions $f$ and $g$
are homogeneous polynomials in the base coordinates. Specifically,
they should be considered as sections $\Gamma$ of appropriate tensor
powers of the anti-canonical bundle $K_{B}^{-1}$ of $B$: for $f$,
it is the fourth power and for $g$, the sixth. 
A nice summary discussion
of this can be found, for example, in \S3 of \cite{11}.

The $j$-invariant of a Weierstra\ss\ model is given by:
\begin{equation}
j\left(z_{i}\right)=1728\cdot\frac{4f\left(z_{i}\right)^{3}}{4f\left(z_{i}\right)^{3}-27g\left(z_{i}\right)^{2}} \ .
\end{equation}
Here, we have marked explicitly the dependence on the base coordinates
$z_{i}$. In particular, we can see that both the numerator and denominator
are sections of $\left(K_{B}^{-1}\right)^{\otimes12}$.

\medskip{}

Let us return to the index 36 torsion-free, genus zero congruence
subgroups. As mentioned, the $j$-invariants for these are computed
and tabulated in \cite{2}; these are rational functions in terms
of the principal modulus $t$ and we interpreted these as Belyi maps in our above discussions to construct the dessins. 
From $j$ we can extract the corresponding
polynomials $f\left(t\right)$ and $g\left(t\right)$ for each subgroup.
We find that the $f\left(t\right)$ will be a degree 12 and $g\left(t\right)$,
degree 18, polynomial in $t$. In general, for an index $\mu$ subgroup
(recall again that $\mu\in\left\{ 6,12,24,36,48,60\right\} $ for genus zero, torsion-free
congruence subgroups), the degrees of $f$ and $g$ are $\mu/3$ and
$\mu/2$.

Finally, let us recall a standard fact about $\mathbb{P}^{2}$ to
complete the story. The anti-canonical bundle of the projective plane
is $K_{\mathbb{P}^{2}}^{-1}=\mathcal{O}_{\mathbb{P}^{2}}\left(3\right)$,
therefore sections of this would be degree 3 homogeneous polynomials
in the coordinates. Therefore, if we choose our base $B$ to be $\mathbb{P}^{2}$
with affine coordinates $z_{1},z_{2}$, then $f\left(z_{i}\right)$
will be a degree $3\times4=12$ and $g\left(z_{i}\right)$ a degree
$3\times6=18$ homogeneous polynomial in $z_{1},z_{2}$. These are
the right degrees to match those of $f\left(t\right)$ and $g\left(t\right)$
for the index 36 case! Therefore, setting a linear change of basis
(normalizing the leading coefficient)
\begin{equation}\label{P1inP2}
z_{1}=t+a\quad,\quad z_{2}=t+b\qquad a,b\in\mathcal{\mathbb{C}}
\end{equation}
and substituting into appropriately chosen $f\left(z_{i}\right)$
and $g\left(z_{i}\right)$ should give the required $f\left(t\right)$
and $g\left(t\right)$. 
Geometrically, this amounts to finding a hyper-plane
$\mathbb{P}^{1}$ defined by \eqref{P1inP2} inside the base $B=\mathbb{P}^{2}$
on which our elliptically fibred Calabi-Yau threefold restricts to
the required modular surface for index 36.

The analogous situation for the index 24 case is clear. 
There, we have K3 surfaces which are elliptic fibrations over a base $\mathbb{P}^{1}$:
we know that the anti-canonical bundle is $K_{\mathbb{P}^{2}}^{-1}=\mathcal{O}_{\mathbb{P}^{2}}\left(2\right)$
whose sections are thus quadratics in the projective coordinate
$t$. Hence, $f\left(t\right)$ and $g\left(t\right)$ should be of
degree $2\times4=8$ and $2\times6=12$ , which are indeed in accord
with the aforementioned fact that the degrees should be respectively
$\mu/3$ and $\mu/2$ with $\mu=24$ here.


\subsubsection{Example: $\Gamma_{0}\left(2\right)\cap\Gamma\left(3\right)$}
Let us see how the above setup works for the specific groups in detail.
The $j$-invariant for the group $\Gamma_{0}\left(2\right)\cap\Gamma\left(3\right)$
is \cite{2}:
\begin{equation}
j\left(t\right)=\frac{\left(\left(t^{3}+4\right)\left(t^{3}+6t^{2}+4\right)\left(t^{6}-6t^{5}+36t^{4}+8t^{3}-24t^{2}+16\right)\right)^{3}}{t^{6}\left(t+1\right)^{3}\left(t^{2}-t+1\right)^{3}\left(t-2\right)^{6}\left(t^{2}+2t+4\right)^{6}}
\end{equation}
Therefore, solving for $\left(\left(t^{3}+4\right)\left(t^{3}+6t^{2}+4\right)\left(t^{6}-6t^{5}+36t^{4}+8t^{3}-24t^{2}+16\right)\right)^{3}=1728\cdot4\cdot f\left(t\right)^{3}$
and $t^{6}\left(t+1\right)^{3}\left(t^{2}-t+1\right)^{3}\left(t-2\right)^{6}\left(t^{2}+2t+4\right)^{6}=4f\left(t\right)^{3}-27g\left(t\right)^{2}$
gives:
\begin{eqnarray}
\nonumber
f\left(t\right) & = & \frac{1}{12\cdot2^{2/3}}\left(t^{12}+232t^{9}+960t^{6}+256t^{3}+256\right)\\
\label{fg-eg1}
g\left(t\right) & = & \frac{1}{216}\left(t^{18}-516t^{15}-12072t^{12}-24640t^{9}-30720t^{6}+6144t^{3}+4096\right)
\end{eqnarray}
which we see to be polynomials of the correct degree. Note here that while $f=0$ is solvable and has relatively
simple roots, the corresponding equation for $g$ is not solvable.

Now, the most general form of degree 12 and 18 polynomials in two
variables are:
\begin{equation}\label{fg-eg}
f\left(z_{1},z_{2}\right) =  \sum_{i=0}^{12}c_{i}z_{1}^{i}z_{2}^{12-i}
\ , \quad
g\left(z_{1},z_{2}\right) =  \sum_{i=0}^{18}c_{i}z_{1}^{i}z_{2}^{18-i} \ ,
\end{equation}
with complex coefficients $c_{i}$ and $d_{i}$. We can readily find
the coefficients $c_{i}$ and $d_{i}$ as well as constants $a,b$
which allow \eqref{fg-eg} to be transformed into the form of \eqref{fg-eg1}.
Setting, for example, $c_1 = d_1 = 0$, we find that $a=1$, $b=-2$, and that:
\begin{equation} 
\begin{array}{rcl} 
f(z_1,z_2) &=& \frac{1}{8748 \cdot 2^{2/3}}
\left(
256 z_1^{12}-256 z_2^3 z_1^9+960 z_2^6 z_1^6-232 z_2^9 z_1^3+z_2^{12}
\right)  \ ;
\\
g(z_1,z_2) &=&
\frac{1}{4251528}
\left(
-4096 z_1^{18}+6144 z_2^3 z_1^{15}+30720 z_2^6 z_1^{12}-24640 z_2^9 z_1^9+12072
   z_2^{12} z_1^6-516 z_2^{15} z_1^3-z_2^{18}
\right)
\end{array}
\end{equation}

We can readily find, using the same technique as above, the explicit Weierstra\ss\ equations, albeit with very complicated {\it algebraic} coefficients, of the remaining index 36 groups, as we had done for the index 24 cases in \cite{1}.
Though the complexity of the numerical expressions makes the presentation unworthwhile, it should be emphasized that all the coefficient we find for the Weierstra\ss\ models must be algebraic numbers.
This is in accord with the deep theorem of Belyi that a rational map from a Riemann surface to $\mathbb{P}^1$ is Belyi (i.e., ramified only on $0,1,\infty$) if and only if there exists an algebraic equation for the Riemann surface over $\overline{\mathbb{Q}}$.


\section{Semi-Stable Elliptic K3 Surfaces}\label{s:K3}
We have mentioned above that the nine conjugacy classes of index 24 subgroups are special in the sense that they correspond to semi-stable elliptic K3 surfaces.
Exploring the physics of these was the subject of \cite{1}.
In this section, we will see how they belong to another important classification scheme which intersects nicely with our above discussion of congruence subgroups.

In \cite{6}, the authors study semi-stable elliptic fibrations over
$\mathbb{P}^{1}$ of K3 surfaces with six singular fibres. 
From the cusp widths (or equivalently the ramification data at $\infty$) presented in Table 1, we see that our subgroups all have six cusps and hence the index 24 cases belong to this class.
In their paper, the authors give a complete list of possible fibre types for such fibrations - there are 112 such cases. 
Our nine conjugacy classes of index 24 subgroups of $\Gamma$ discussed in the previous section correspond to nine of these 112 surfaces. 

Remarkably, the elliptic $j$-invariants for all the 112 surfaces also have the property that, as rational maps from $\mathbb{P}^1$ to $\mathbb{P}^1$, they are Belyi.
Furthermore, all the ramification indices at 0 are 3, and all the ramification indices at 1 are 2; that is, we have trivalent dessins!
In \cite{7}, the dessins d'enfants for
all 112 cases are constructed, along with explicit $j$-functions
with coefficients in $\mathbb{Q}$ (For reference, we include the full list, emphasizing their clean, trivalent nature, in Appendix \ref{ap:112}). Those with coefficients in $\overline{\mathbb{Q}}\setminus\mathbb{Q}$ have yet to be computed, due to the complexity of finding exact roots to polynomials in high degree.
This is thus an interesting open challenge.
Now, because any finite trivalent graph with nodes replaced by oriented triangles is the Schreier coset graph of a finite index subgroup of the modular group, one naturally wonders whether we could attach a conjugacy class of subgroups to each of these dessins.
Of course, only 9 of these will be congruence (they are all genus zero because the K3 surfaces are all fibred over $\mathbb{P}^1$), so what general subgroups can one encounter?

We will see that just as those nine surfaces correspond
to these congruence subgroups of $\Gamma$ (up to modular conjugacy), it turns out that \emph{all}
the 112 surfaces will correspond to some (not necessarily congruence)
subgroup of $\Gamma$. In this section we shall elaborate on the connections
between these surfaces and the subgroups of $\Gamma$.
The structure of this section is as follows: first, we discuss
 the permutation representation of $\Gamma$ on the right cosets
of each modular subgroup, and show how this defines a transitive permutation group known
as a {\em cartographic group}. We then compute
 a generating set for a representative of the conjugacy class of modular subgroups
corresponding to each of the 112 semi-stable elliptic fibrations over $\mathbb{P}^1$ of K3
surfaces with six singular fibres.

\subsection{Permutation Representations and Cartographic Groups}

Returning to our original discussion of the modular group in \S\ref{s:Gamma},
let $\sigma_{0}$ and $\sigma_{1}$ denote the permutations induced
by the respective actions of $S$ and $ST$ on the cosets of each
subgroup.
We can find a third permutation $\sigma_{\infty}$ by imposing
the following condition, thereby constructing a \emph{permutation
triple}:
\begin{equation}\label{triple}
\sigma_{0}\cdot\sigma_{1}\cdot\sigma_{\infty}=1 \ .
\end{equation}

The permutations $\sigma_{0}$, $\sigma_{1}$ and $\sigma_{\infty}$
give the permutation representations of $\Gamma$ on the right cosets
of each modular subgroup in question.
From the point of view of the dessins, there is an equivalent representation in terms of permutation triples \cite{4,leila,LZ}.
Let $d$ be the number of edges, labeled from 1 to $d$, and write elements of the symmetric group on $d$ elements in the standard cycle notation.
Now, at the $i$-th black node write the labels for the incident edges clockwise, forming the cycle $B_i$; likewise, at the $j$-th white node write the labels for the incident edges counter-clockwise, forming the cycle $W_j$.
Then, $\sigma_0$ is simply $\prod_i B_i$ while $\sigma_1 = \prod_j W_j$.
Remarkably, $\sigma_\infty$ as determined by \eqref{triple} will have as many cycle products as there are faces, such that twice the length of each cycle is the number of edges of the corresponding face.
Conversely, a {\it unique} dessin and Belyi map of degree $d$ is thus determined.

As elements of the symmetric group,
$\sigma_{0}$ and $\sigma_{1}$ can be easily computed from
the Schreier coset graphs in \cite{2} by following the procedure
elaborated in \cite{3}, i.e., by noting that the doubly directed edges
represent an element $x$ of order 2, while the positively oriented
triangles represent an element $y$ of order 3. Since the graphs are
connected, the group generated by $x$ and $y$ is transitive on the
vertices. Clearly, $\sigma_{0}$ and $\sigma_{1}$ tell us which vertex of the coset graph is sent to which, i.e. which coset of the modular subgroup in question is sent to which by the action of $\Gamma$ on the right cosets of this subgroup.

To express the cycle decomposition of the permutation $\sigma_{\infty}$,
we choose a vertex and apply $x$, then $y$, i.e., we take an edge
(doubly oriented) followed by a side of a triangle (always in the
positive orientation), until we return back to the original vertex.
The length of this circuit gives the length of a cycle in $\sigma_{\infty}$.
In this way, we can construct the permutation decompositions $\sigma_{0}$,
$\sigma_{1}$ and $\sigma_{\infty}$ for every dessin under
consideration. In turn, any two of these three permutation decompositions
define a transitive permutation group. 
Of course, we only need two of these three permutations since the third will
simply be the inverse of the product of the other two by \eqref{triple}, and we will henceforth focus on $\sigma_0$ and $\sigma_1$.
These permutation groups are known as {\em cartographic groups}, and are discussed in e.g. \cite{4}.

\subsubsection{Cartographic Groups for the 112 K3 Surfaces}
\label{s:carto}
A priori, one would imagine that one could obtain the semi-stable elliptic K3 surfaces with 6 singular (I-type) singularities by partitioning 24 into 6 parts.
There are 199 possibilities and \cite{6} showed that only 112 are allowed. 
In general, one would try to partition 24 into at least 6 parts - the fact that there are at least 6 singular fibres is a restriction coming from the geometry of K3 surfaces. In our case, there are exactly 6 singular fibres,
making these surfaces {\it extremal}. These surfaces are listed in \cite{6}.
 
We will thus distinguish our 112 extremal K3 surfaces by a 6-tuple partition $(I_1,\ldots,I_6)$ of 24; these will correspond to the I-type singular fibres as well as the ramification data (i.e., cusp widths) at infinity.
A table of the cartographic groups for the partitions in \cite{6}
is given below.
Now, dessins have intriguing number theoretic properties and the field of definition of (the coefficients in) the corresponding $j$-function as a Belyi map is an important quantity.
For reference, we include the number field in brackets next to the partition.
For simplicity, where we write $(\sqrt{-3})$, for example, we mean the field is $\mathbb{Q}(\sqrt{-3})$. 
Where it says e.g. ``cubic'', we mean it is $\mathbb{Q}$ extended by the root of some cubic polynomial.

We note that the partitions do not necessarily correspond uniquely to a dessin - when they do, the $j$-function is always defined over $\mathbb{Q}$.
Where a single partition has multiple dessins defined over different
number fields, we have for completeness computed the cartographic
groups for all such cases, indicating explicitly the field as in \cite{7}.
We also keep the order of the entries the same as the order of the dessins
in \cite{7} and in Appendix \ref{ap:112}.

To identify these cartographic groups, we use the
corresponding dessins given in \cite{7} to construct the Schreier
coset graph in each case, following the procedure detailed in \S\ref{s:dessin}.
We then read the permutations
$\sigma_{0}$ and $\sigma_{1}$ directly off these graphs. As mentioned,
together these permutations define a transitive permutation group,
which is the cartographic group in question. 
Extensive use will be made of the computer software 
\textsf{GAP} \cite{GAP}.
The decomposition of each group in terms of better known groups can be identified using commands such as \texttt{StructureDescription()} and \texttt{TransitiveIdentification()} and
it is these decompositions which are given in the ensuing table\footnote{We thank Alexander Hulpke for helping to identify many of the decompositions.}.  We follow \cite{5} and use the following group notation: A (alternating), AL (affine linear), C (cyclic), E (elementary),
S (symmetric). Direct products are denoted by `$\times$', split
extensions by `$:$', and wreath products by `$\wr$'. In
addition, we use $L\left(n\right)$ to denote groups derived from
linear groups as defined in Table 2 \cite{5}. A group name followed by a number
$n$ in round brackets `$\left(n\right)$' denotes that the group operates transitively on $n$
points, or respectively that a point stabilizer has index $n$ in the group \cite{5}.

\begin{table}[H]
\begin{centering}
\small{
\begin{tabular}{|c|c|c|}
\hline 
\textbf{Name} & \textbf{Size} & \textbf{Type}\tabularnewline
\hline 
\hline 
$L\left(6\right)$  & 60 & $PSL\left(2,5\right)=A_{5}\left(6\right)$\tabularnewline
\hline 
$L\left(7\right)$  & 168 & $PSL\left(3,2\right)$\tabularnewline
\hline 
$L\left(8\right)$  & 168 & $PSL\left(2,7\right)$\tabularnewline
\hline 
$L\left(9\right)$  & 504 & $PSL\left(2,8\right)$\tabularnewline
\hline 
$L\left(10\right)$  & 360 & $PSL\left(2,9\right)=A_{6}\left(10\right)$\tabularnewline
\hline 
$L\left(11\right)$  & 660 & $PSL\left(2,11\right)\left(11\right)$\tabularnewline
\hline 
$L\left(12\right)$  & 660 & $PSL\left(2,11\right)$\tabularnewline
\hline 
$L\left(13\right)$  & 5616 & $PSL\left(3,3\right)$\tabularnewline
\hline 
$L\left(14\right)$  & 1092 & $PSL\left(2,13\right)$\tabularnewline
\hline 
$L\left(15\right)$  & 20160 & $PSL\left(4,2\right)=A_{8}\left(15\right)$\tabularnewline
\hline 
\end{tabular}
}
\par\end{centering}

\caption{\em The permutation groups $L\left(n\right)$}
\end{table}

Where no easy decomposition is available, we have simply
given the cartographic group's index in the \textsf{GAP} library of
transitive permutation groups (here a number of the form t24n...). The number after `t' in this index gives the degree of the transitive permutation group, while the number after the `n' gives the location that group in the \textsf{GAP} library of transitive permutation groups for that particular group degree.

{\tiny }%
\begin{longtable}{|c||c||c|}
\hline 
\textbf{\tiny Partition (K3 Surface)} & \textbf{\tiny Cartographic Group} & \textbf{\tiny Cartographic Group Order}\tabularnewline
\hline 
\hline 
{\tiny $19,1,1,1,1,1\;\left(\mathbb{Q}\right)$} & {\tiny $A_{24}$} & {\tiny 24!/2}\tabularnewline
\hline 
\hline 
{\tiny $18,2,1,1,1,1\;\left(\mathbb{Q}\right)$} & {\tiny $\mathrm{t24n18368}=\mathrm{t24n16559}\times C_{2}$} & {\tiny $2^{11}\cdot3^{4}$}\tabularnewline
\hline 
\hline 
{\tiny $18,2,1,1,1,1\;\left(\sqrt{-3}\right)$} & {\tiny $A_{24}$} & {\tiny 24!/2}\tabularnewline
\hline 
\hline 
{\tiny $18,2,1,1,1,1\;\left(\sqrt{-3}\right)$} & {\tiny $A_{24}$} & {\tiny 24!/2}\tabularnewline
\hline 
\hline 
{\tiny $17,3,1,1,1,1\;\left(\mathbb{Q}\right)$} & {\tiny $A_{24}$} & {\tiny 24!/2}\tabularnewline
\hline 
\hline 
{\tiny $17,2,2,1,1,1\;\left(\sqrt{17}\right)$} & {\tiny $A_{24}$} & {\tiny 24!/2}\tabularnewline
\hline 
\hline 
{\tiny $17,2,2,1,1,1\;\left(\sqrt{17}\right)$} & {\tiny $A_{24}$} & {\tiny 24!/2}\tabularnewline
\hline 
\end{longtable}

\subsection{Generators for the Modular Subgroups}\label{s:gens}

There is a so-called {\em Reidemeister-Schreier process} 
(discussed in \cite{9}) which
enables one to compute presentations of subgroups $H$ of finite index
in a group $G$ defined by finite presentations. 
We can use this process to find the presentations of all the modular
subgroups of interest. In our case, $G$ is $\Gamma=\mathit{\mathrm{\mathit{PSL}}}(2,\mathbb{Z})$,
the presentation of which is $\left\langle x,y|x^{2},y^{3}\right\rangle $.
For all our conjugacy classes of modular subgroups, the
Reidemeister-Schreier process yields the following group presentation, when the permutation
data (here $\sigma_{0}$, $\sigma_{1}$) are input in the manner detailed in \cite{9}:
\begin{equation}\label{free5}
\left\langle h,i,j,k,l|\mbox{-----}\right\rangle 
\end{equation}
Thus in every case we have a free group on 5 generators.
Considering the specific example of $\Gamma\left(4\right)$ for illustration,
this should not be surprising, since the modular curve
$X\left(4\right)$ is $\mathbb{P}^{1}$ with six cusps and no elliptic
points, so its fundamental group is the free group on five generators.

We can find the generators  for a representative of all the conjugacy classes of subgroups
of interest, which we now know to be free groups on five generators, using
 \textsf{GAP}\footnote{We thank Alexander Hulpke for guidance here.
}. First, we use the permutation data $\sigma_{0}$, $\sigma_{1}$ obtained
from each of the Schreier coset graphs (in turn obtained from each
of the dessins) to find the group homomorphism by images between $G$
and $H$. We then use this to define the $H$ in question. Finally,
we use the \textsf{GAP} command \texttt{IsomorphismFpGroup(G)}, which
returns an isomorphism from the given $H$ to a finitely presented
group isomorphic to $H$. This function first chooses a set of generators
of $H$ and then computes a presentation in terms of these generators.

To give an example, consider the partition $4,4,4,4,4,4$, which
corresponds to the conjugacy class of subgroups $\Gamma\left(4\right)$. 
We can find a set of generators as $2\times2$ matrices for
a representative of this class of subgroups by implementing the following code in \textsf{GAP}:

\medskip{}

\texttt{\small gap> f:=FreeGroup(``x'',``y'');}{\small \par}

\texttt{\small <free group on the generators {[} x, y {]}>}{\small \par}

\texttt{\small gap>gamma:=f/{[}f.1\textasciicircum{}2,f.2\textasciicircum{}3{]};}{\small \par}

\texttt{\small <fp group on the generators {[} x, y {]}>}{\small \par}

\texttt{\small gap> hom:=GroupHomomorphismByImages(gamma,Group(sigma\_0,sigma\_1),}{\small \par}

\texttt{\small >GeneratorsOfGroup(gamma),{[}sigma\_0,sigma\_1{]});}{\small \par}

\texttt{\small {[}x, y{]} ->}{\small \par}

\texttt{\small {[}(1,10)(2,4)(3,24)(5,7)(6,21)(8,12)(9,18)(11,14)(13,16)(15,23)}{\small \par}

\texttt{\small (17,19)(20,22), (1,2,3)(4,5,6)(7,8,9)(10,11,12)(13,14,15)}{\small \par}

\texttt{\small (16,17,18)(19,20,21)(22,23,24){]} }{\small \par}

\texttt{\small gap> gamma4:=PreImage(hom,Stabilizer(Image(hom),1));}{\small \par}

\texttt{\small Group(<fp, no generators known>)}{\small \par}

\texttt{\small gap> iso:=IsomorphismFpGroup(gamma4);}{\small \par}

\texttt{\small {[} <{[} {[} 1, 1 {]} {]}|y{*}x{*}y{*}x{*}y{*}x\textasciicircum{}-1{*}y{*}x\textasciicircum{}-1>,}{\small \par}

\texttt{\small <{[} {[} 2, 1 {]} {]}|y\textasciicircum{}-1{*}x{*}y\textasciicircum{}-1{*}x{*}y\textasciicircum{}-1{*}x\textasciicircum{}
-1{*}y\textasciicircum{}-1{*}x\textasciicircum{}-1>,}{\small \par}

\texttt{\small <{[} {[} 3, 1 {]} {]}|y{*}x{*}y\textasciicircum{}-1{*}x{*}y\textasciicircum{}-1{*}x\textasciicircum{}-1{*}y\textasciicircum{}-1{*}x\textasciicircum{}-1{*}y>,}{\small \par}

\texttt{\small <{[} {[} 4, 1 {]} {]}|x{*}y{*}x{*}y\textasciicircum{}-1{*}x{*}y\textasciicircum{}-1{*}x\textasciicircum{}-1{*}y\textasciicircum{}-1{*}x\textasciicircum{}-1{*}y{*}x\textasciicircum{}-1>,}{\small \par}

\texttt{\small <{[} {[} 5, 1 {]} {]}|y{*}x{*}y\textasciicircum{}-1{*}x{*}y{*}x{*}y\textasciicircum{}-1{*}x\textasciicircum{}-1{*}y{*}x\textasciicircum{}-1{*}y\textasciicircum{}-1{*}x\textasciicircum{}-1>
{]}}{\small \par}

\texttt{\small -> {[} F1, F2, F3, F4, F5 {]}}{\small \par}

\medskip{}

Here, the only input in each case are the permutation data $\sigma_{0}$,
$\sigma_{1}$. Once the output has been obtained, we will see the five generators $[F1, F2,F3, F4, F5]$, as indicated in the last line of the output, as functions of $x$ and $y$, i.e. the generators of $\Gamma$.
Now, returning to the matrices \eqref{xy}, the only thing
left to do is to multiply together the matrices and their inverses
as indicated. This will produce the generators of each $H$, as desired.
Note that one can achieve the same results here in \textsf{Magma} though
use of the \texttt{Rewrite} command\footnote{We thank an anonymous referee for this observation.} \cite{magma}.

\subsubsection{Catalogue of Generators as Explicit $2\times2$ Matrices}

In this section we give the generators for a representative of each of
 the conjugacy classes of modular subgroups
corresponding to the entries in the table at the end of \S\ref{s:carto}. 
For each of the 112 K3 surfaces, denoted as a 6-tuple partition of 24 (the notation is as before, we also include the number field of definition for the associated $j$-function as a Belyi map), we write the five matrix generators explicitly.

\begin{center}
{\scriptsize }%
\begin{longtable}{|c|c|}
\hline 
\textbf{\scriptsize Partition} & \textbf{\scriptsize Generators}\tabularnewline
\hline 
\hline 
{\scriptsize $19,1,1,1,1,1\;\left(\mathbb{Q}\right)$} & {\scriptsize $\left(\begin{array}{rr}
-1 & 1\\
0 & -1
\end{array}\right)\left(\begin{array}{rr}
5 & 4\\
-9 & -7
\end{array}\right)\left(\begin{array}{rr}
3 & 1\\
-16 & -5
\end{array}\right)\left(\begin{array}{rr}
18 & -5\\
25 & 0
\end{array}\right)\left(\begin{array}{rr}
20 & 9\\
-49 & -22
\end{array}\right)$}\tabularnewline
\hline 
\hline 
{\scriptsize $18,2,1,1,1,1\;\left(\mathbb{Q}\right)$} & {\scriptsize $\left(\begin{array}{rr}
-1 & 1\\
0 & -1
\end{array}\right)\left(\begin{array}{rr}
5 & 4\\
-9 & -7
\end{array}\right)\left(\begin{array}{rr}
-5 & -2\\
18 & 7
\end{array}\right)\left(\begin{array}{rr}
5 & 1\\
-36 & -7
\end{array}\right)\left(\begin{array}{rr}
17 & 4\\
-81 & -19
\end{array}\right)$}\tabularnewline
\hline 
\hline 
{\scriptsize $18,2,1,1,1,1\;\left(\sqrt{-3}\right)$} & {\scriptsize $\left(\begin{array}{rr}
-1 & 1\\
0 & -1
\end{array}\right)\left(\begin{array}{rr}
5 & 4\\
-9 & -7
\end{array}\right)\left(\begin{array}{rr}
3 & 1\\
-16 & -5
\end{array}\right)\left(\begin{array}{rr}
-19 & -8\\
50 & 21
\end{array}\right)\left(\begin{array}{rr}
43 & 16\\
-121 & -45
\end{array}\right)$}\tabularnewline
\hline 
\hline 
{\scriptsize $18,2,1,1,1,1\;\left(\sqrt{-3}\right)$} & {\scriptsize $\left(\begin{array}{rr}
-1 & 1\\
0 & -1
\end{array}\right)\left(\begin{array}{rr}
5 & 4\\
-9 & -7
\end{array}\right)\left(\begin{array}{rr}
9 & 4\\
-25 & -11
\end{array}\right)\left(\begin{array}{rr}
-7 & -2\\
32 & 9
\end{array}\right)\left(\begin{array}{rr}
5 & 1\\
-36 & -7
\end{array}\right)$}\tabularnewline
\hline 
\end{longtable}
\par\end{center}

\section{Conclusions and Outlook}

In this paper, we have further developed the connections between
the torsion-free, genus zero modular subgroups and the theory of dessins
d'enfants. This was accomplished through construction of
the ramification data and the dessins for each conjugacy class of subgroups. 
In addition,
we have shown that the index 24 and 36 subgroups are of special interest,
since they are naturally associated to K3 surfaces and Calabi-Yau
threefolds, respectively.

In a parallel vein, we have found the generating set of a representative of 
the conjugacy classes of modular subgroups associated with 
all the 112
semi-stable elliptic fibrations over $\mathbb{P}^{1}$ which are extremal 
K3 surfaces with six singular fibres. 
Indeed, at the intersection of the 112 and 33 lie the nine K3 surfaces.
This result is important not only for developing
our understanding of these surfaces and their connection to the theory
of modular subgroups, but also because these connections should give
us further tools for categorizing these surfaces, and the relations
between them.

To our intricate web of correspondences we can add physics.
In \cite{13}, Seiberg-Witten curves for four-dimensional
supersymmetric field theories with unitary gauge groups were given an interpretation as clean dessins.
In \cite{1}, the nine K3 surfaces were singled out to relate to the recently advanced class of $\mathcal{N}=2$ Gaiotto theories.
Now, we have an extensive catalogue of dessins which have underlying Calabi-Yau structure, to construct their associated gauge theory is clearly an interesting story which the authors are currently pursuing.

\parskip
\parskip

\appendix


\section{Dessins from the Miranda-Persson Table}\label{ap:112}

In this appendix we give the dessins for all the entires in \cite{6}.
These dessins match those already given in \cite{7}, hence
their being relegated to an appendix. Nevertheless, there are two
advantages to our drawings over those in \cite{7}: first, they are
drawn electronically; second, we have marked the black/white nodes as an emphasis of the bipartite nature of the graphs.

\begin{center}
{\scriptsize }%
\begin{minipage}[t]{0.33\textwidth}%
\begin{center}
{\scriptsize \includegraphics[scale=0.15]{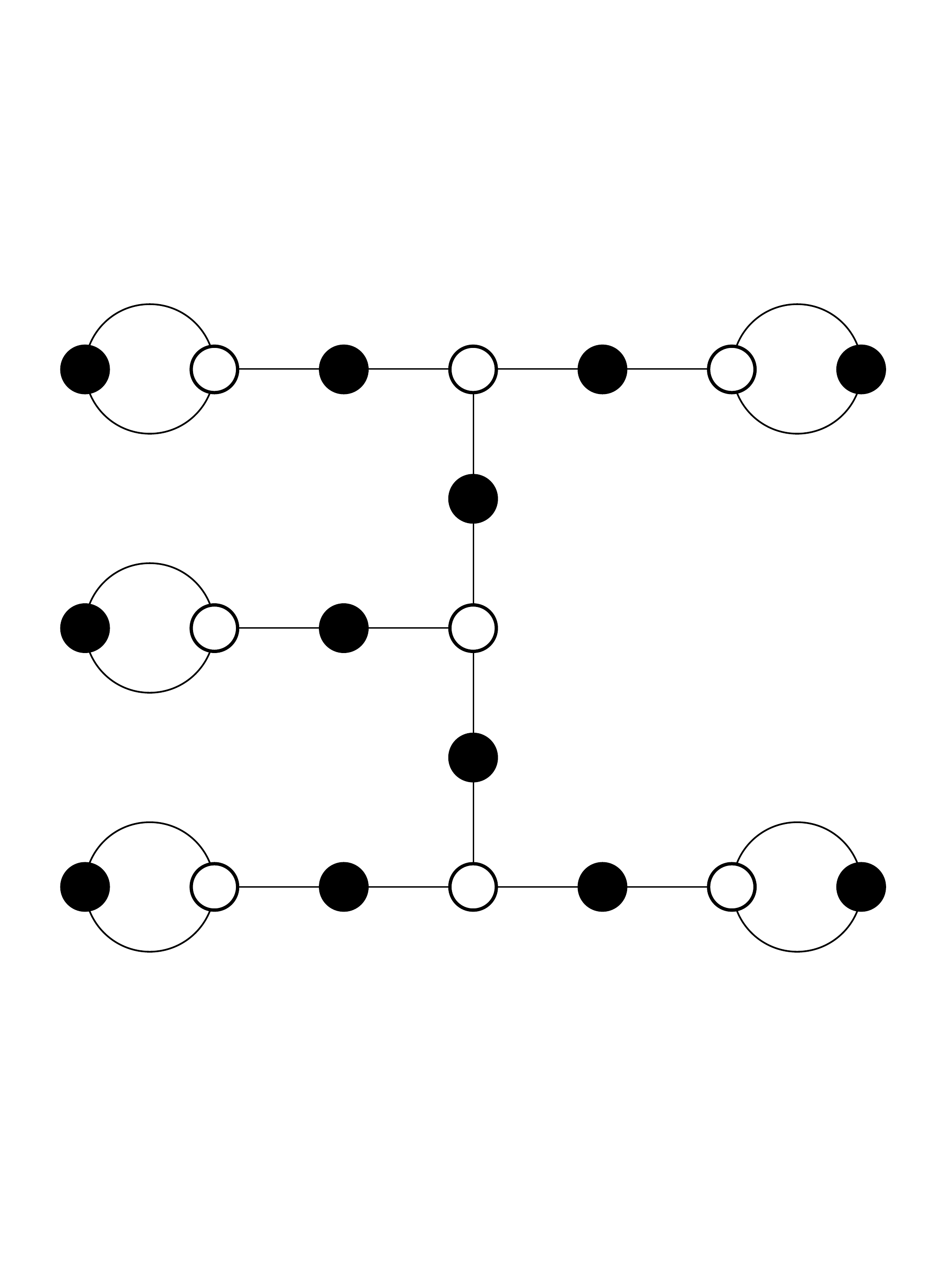}}
\par\end{center}{\scriptsize \par}

\begin{center}
{\scriptsize $19,1,1,1,1,1\;\left(\mathbb{Q}\right)$}
\par\end{center}%
\end{minipage}{\scriptsize }%
\begin{minipage}[t]{0.33\textwidth}%
\begin{center}
{\scriptsize \includegraphics[scale=0.15]{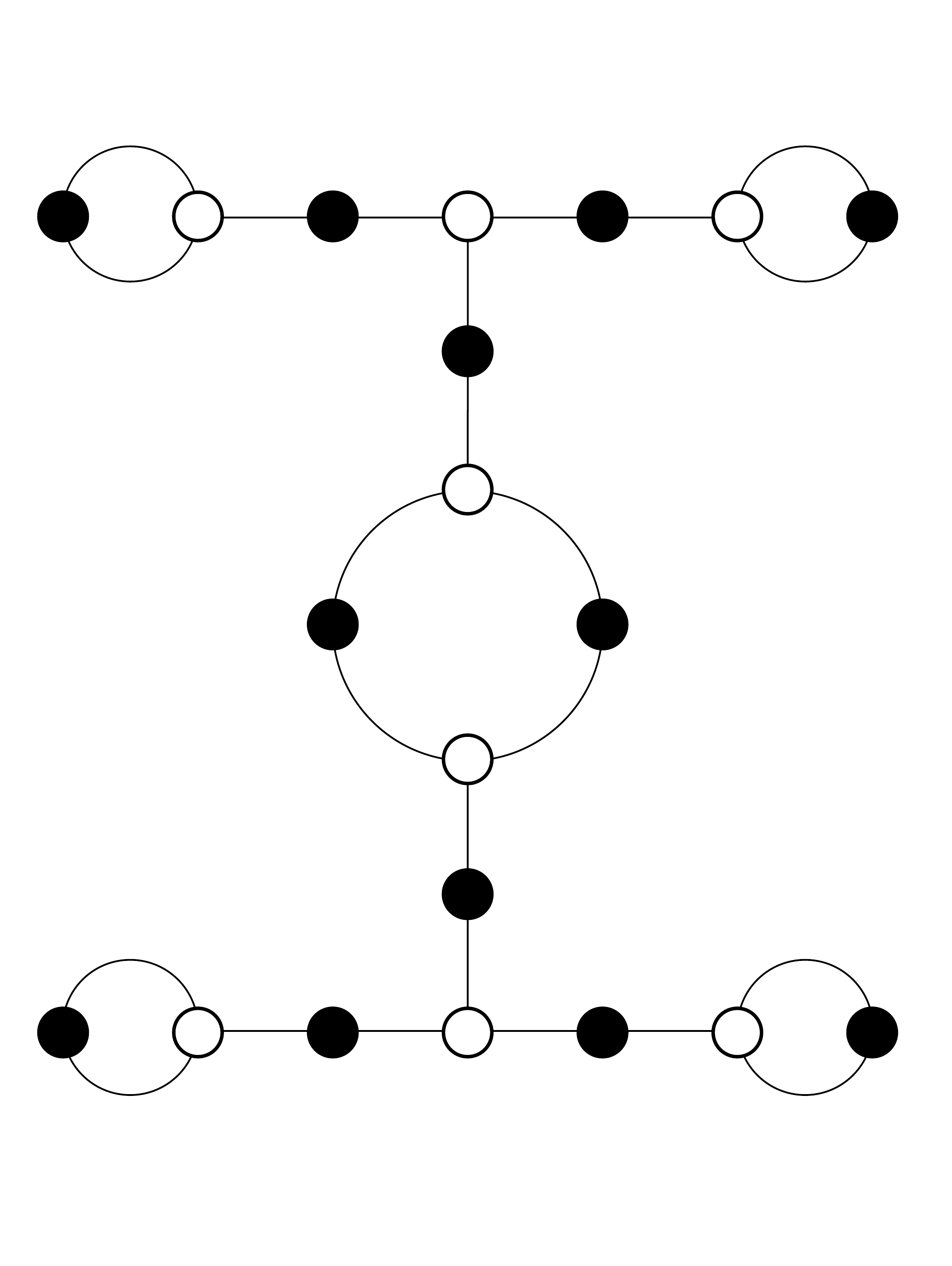}}
\par\end{center}{\scriptsize \par}

\begin{center}
{\scriptsize $18,2,1,1,1,1\;\left(\mathbb{Q}\right)$}
\par\end{center}%
\end{minipage}{\scriptsize }%
\begin{minipage}[t]{0.33\textwidth}%
\begin{center}
{\scriptsize \includegraphics[scale=0.15]{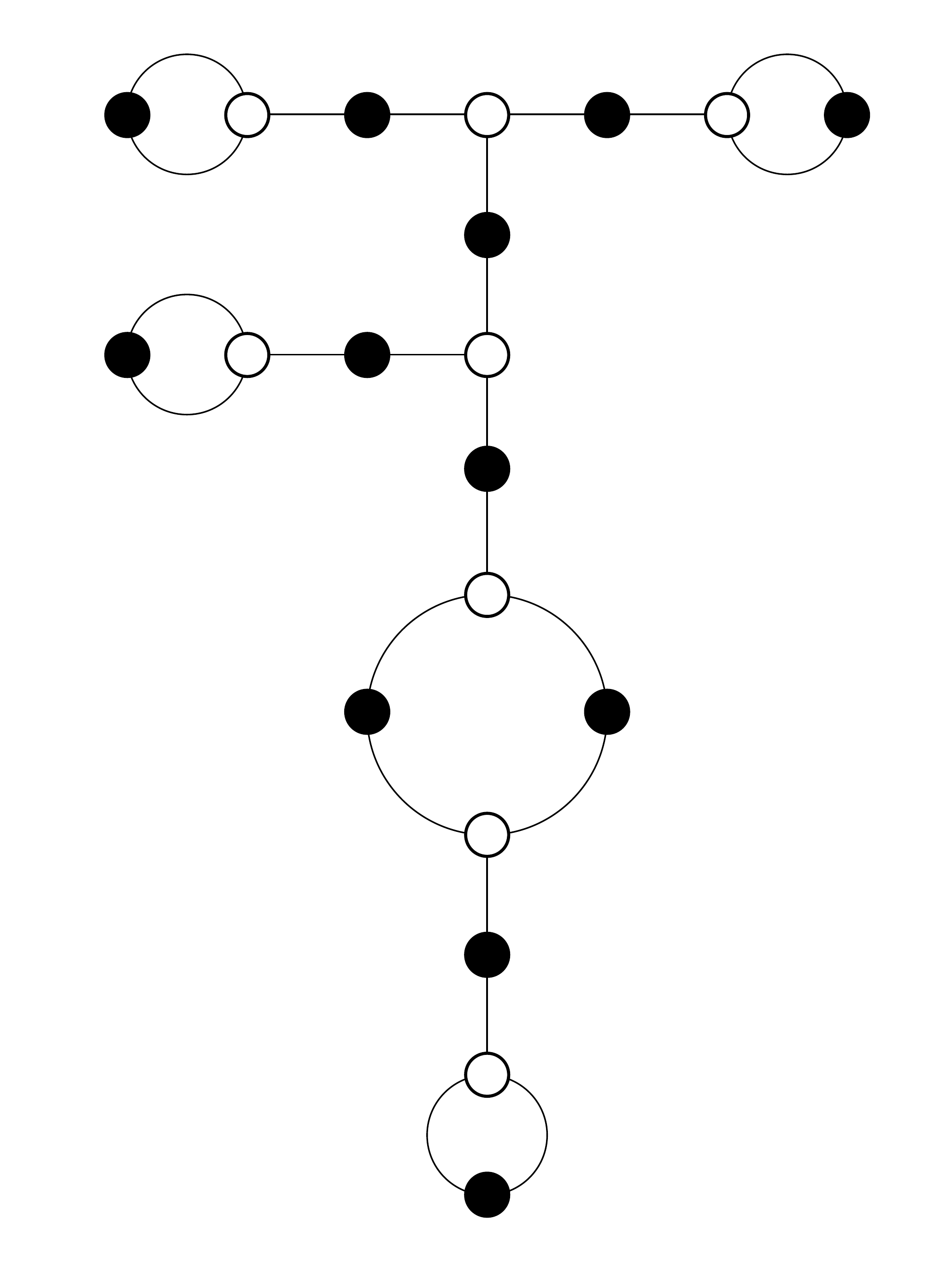}}
\par\end{center}{\scriptsize \par}

\begin{center}
{\scriptsize $18,2,1,1,1,1\;\left(\sqrt{-3}\right)$}
\par\end{center}%
\end{minipage}
\par\end{center}{\scriptsize \par}

\begin{center}
{\scriptsize }%
\begin{minipage}[t]{0.33\textwidth}%
\begin{center}
{\scriptsize \includegraphics[scale=0.15]{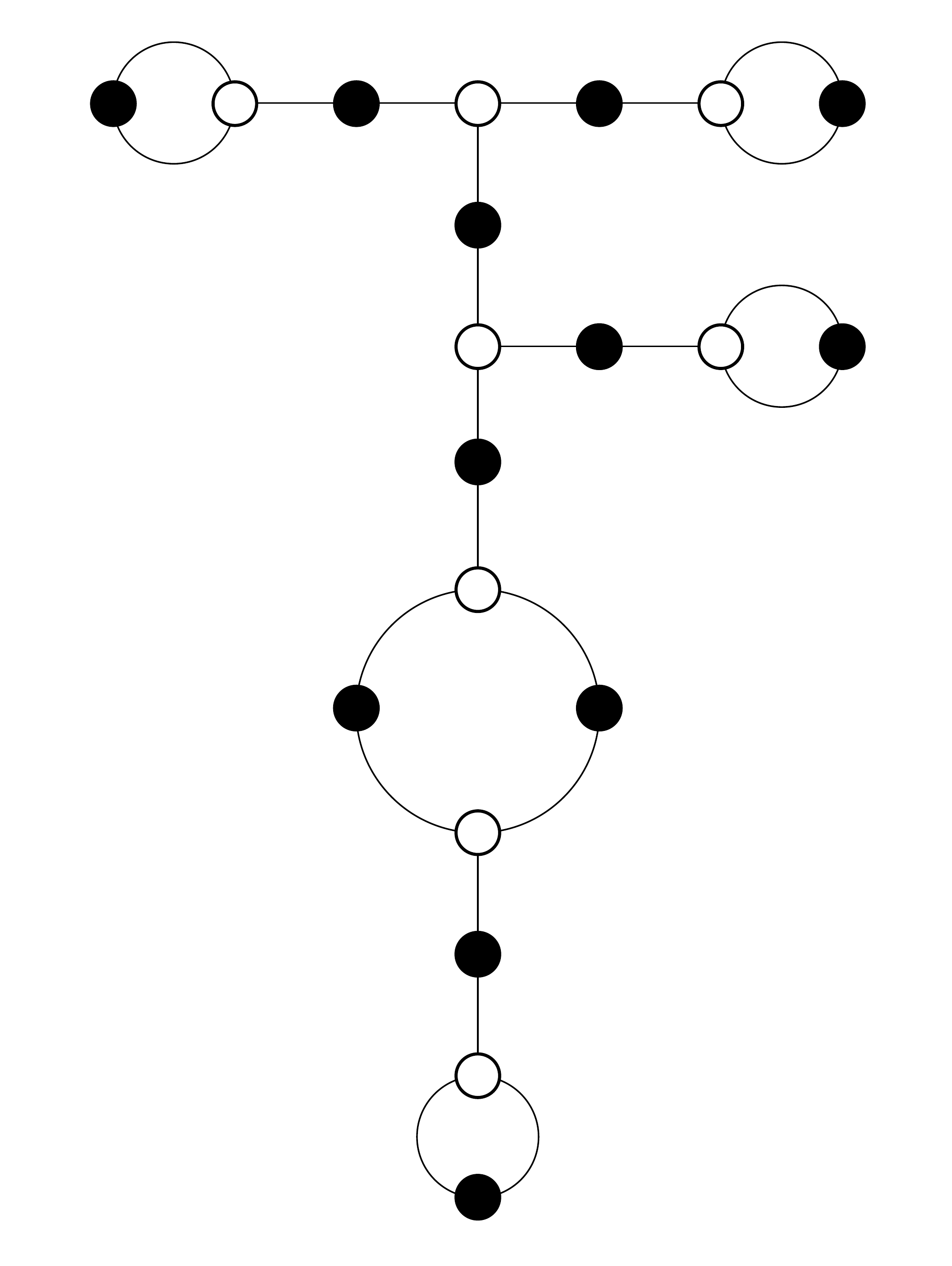}}
\par\end{center}{\scriptsize \par}

\begin{center}
{\scriptsize $18,2,1,1,1,1\;\left(\sqrt{-3}\right)$}
\par\end{center}%
\end{minipage}{\scriptsize }%
\begin{minipage}[t]{0.33\textwidth}%
\begin{center}
{\scriptsize \includegraphics[scale=0.15]{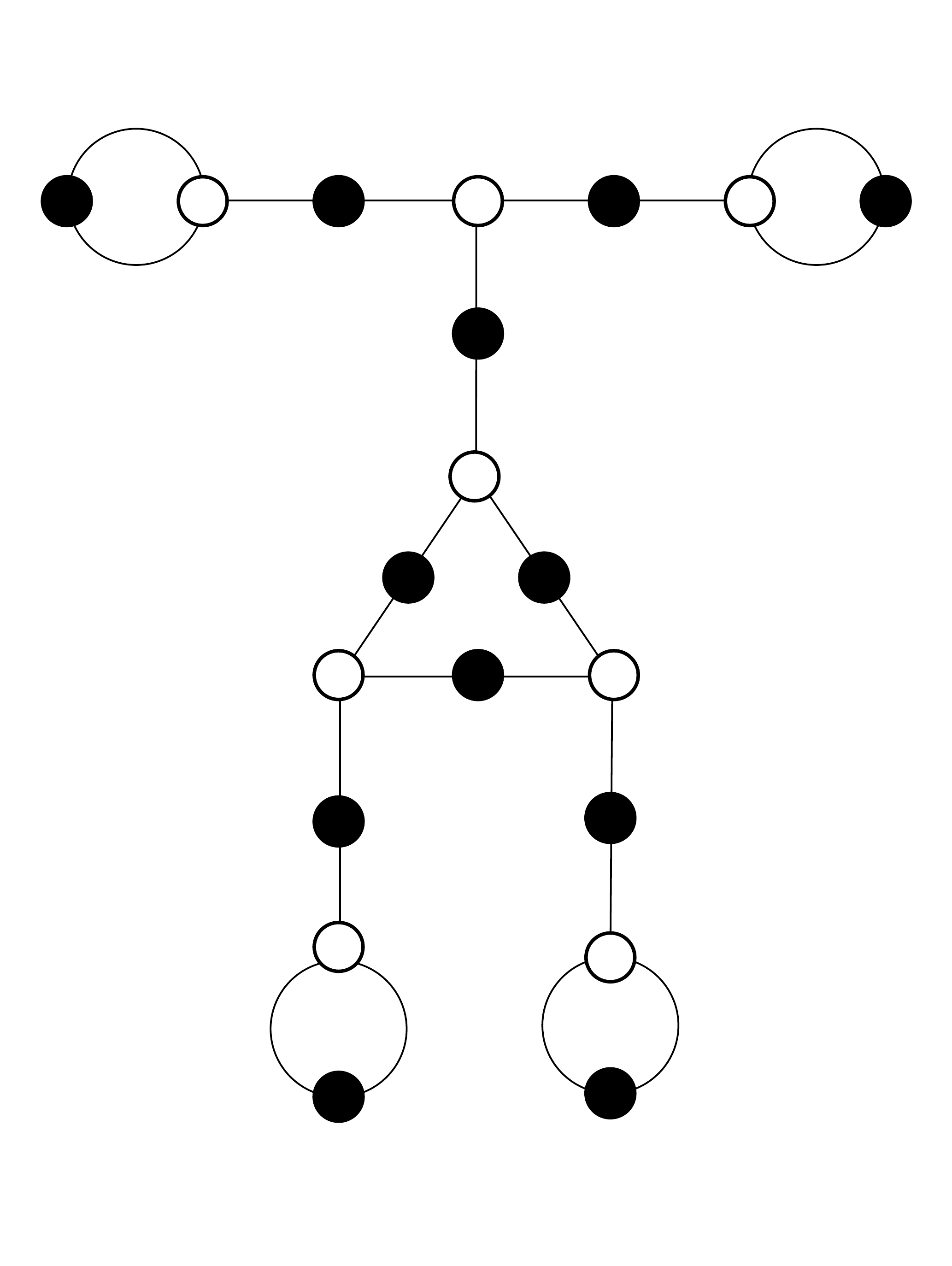}}
\par\end{center}{\scriptsize \par}

\begin{center}
{\scriptsize $17,3,1,1,1,1\;\left(\mathbb{Q}\right)$}
\par\end{center}%
\end{minipage}{\scriptsize }%
\begin{minipage}[t]{0.33\textwidth}%
\begin{center}
{\scriptsize \includegraphics[scale=0.15]{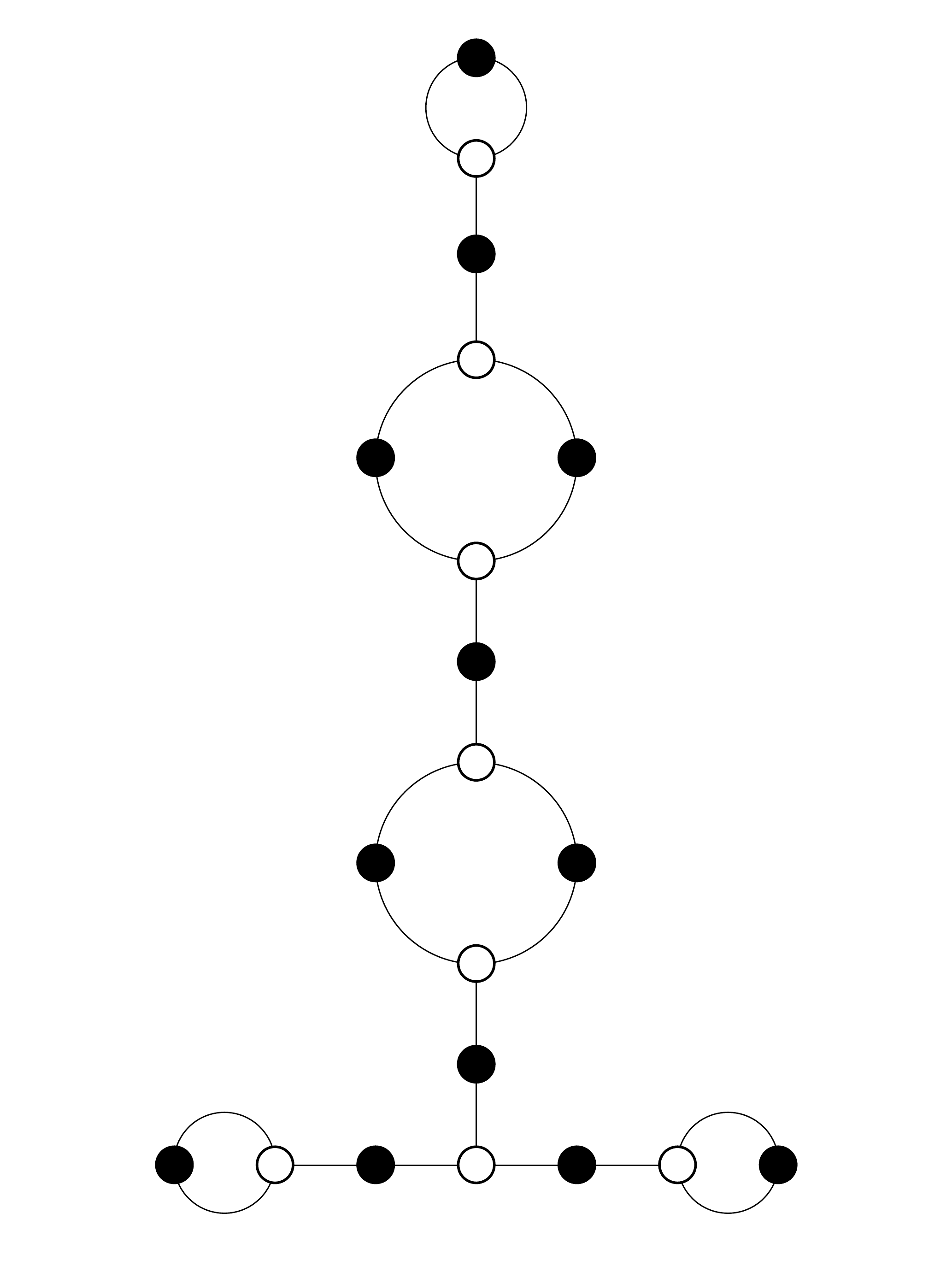}}
\par\end{center}{\scriptsize \par}

\begin{center}
{\scriptsize $17,2,2,1,1,1\;\left(\sqrt{17}\right)$}
\par\end{center}%
\end{minipage}
\par\end{center}{\scriptsize \par}

\begin{center}
{\scriptsize }%
\begin{minipage}[t]{0.33\textwidth}%
\begin{center}
{\scriptsize \includegraphics[scale=0.15]{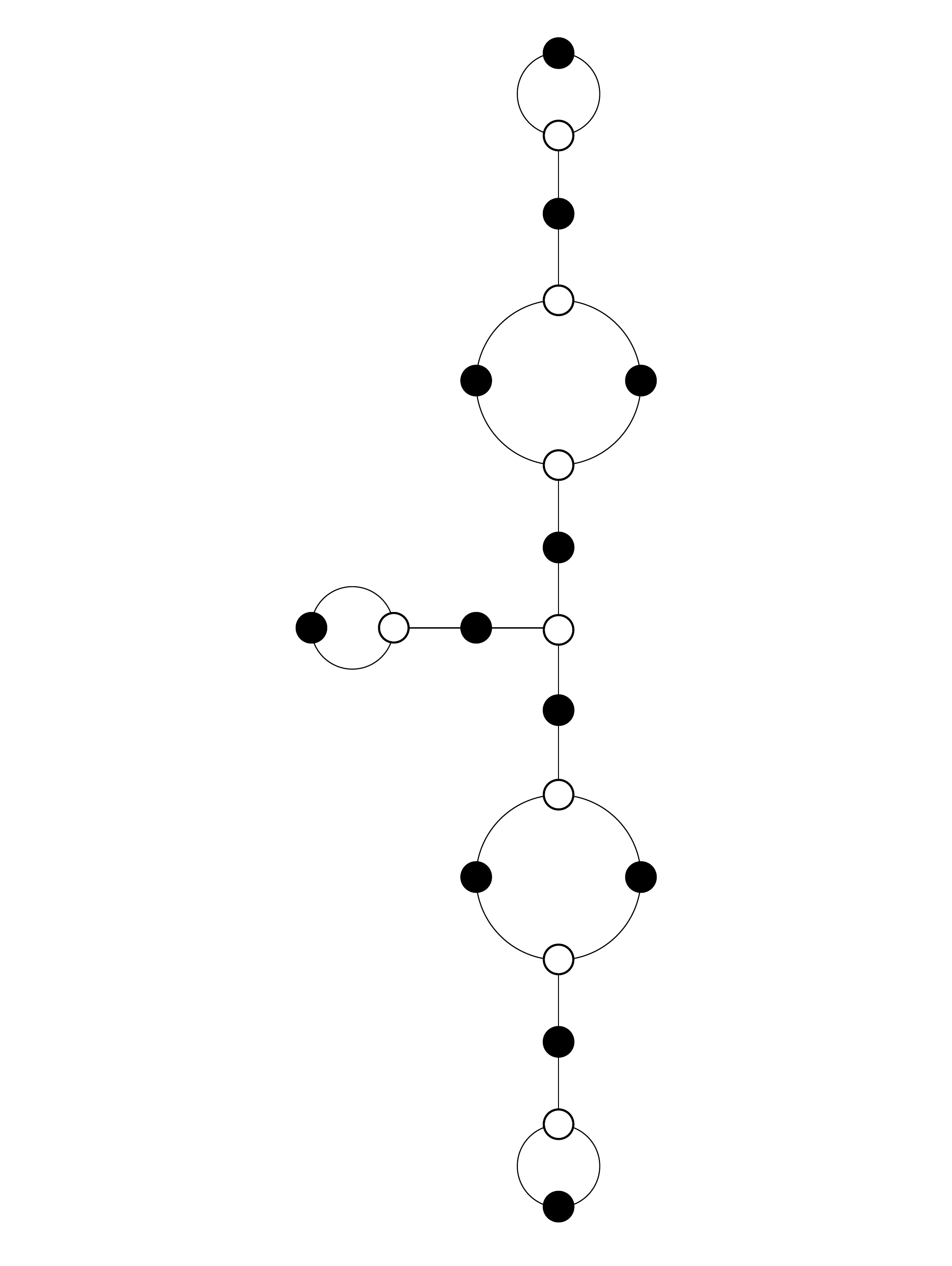}}
\par\end{center}{\scriptsize \par}

\begin{center}
{\scriptsize $17,2,2,1,1,1\;\left(\sqrt{17}\right)$}
\par\end{center}%
\end{minipage}{\scriptsize }%
\begin{minipage}[t]{0.33\textwidth}%
\begin{center}
{\scriptsize \includegraphics[scale=0.15]{\string"PICT/New_16-4-1-1-1-1\string".pdf}}
\par\end{center}{\scriptsize \par}

\begin{center}
{\scriptsize $16,4,1,1,1,1\;\left(\mathbb{Q}\right)$}
\par\end{center}%
\end{minipage}{\scriptsize }%
\begin{minipage}[t]{0.33\textwidth}%
\begin{center}
{\scriptsize \includegraphics[scale=0.15]{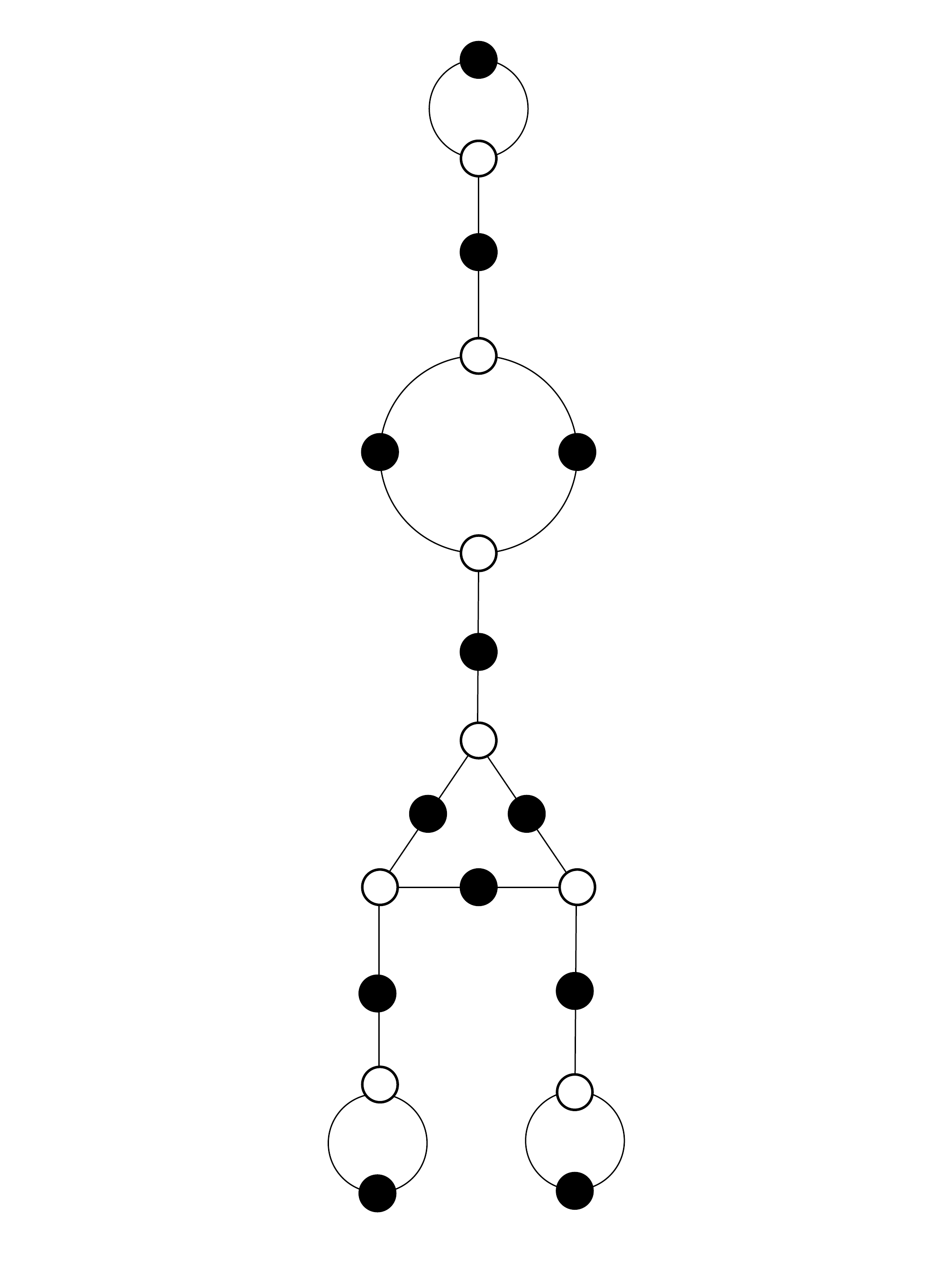}}
\par\end{center}{\scriptsize \par}

\begin{center}
{\scriptsize $16,3,2,1,1,1\;\left(\mathbb{Q}\right)$}
\par\end{center}%
\end{minipage}
\par\end{center}{\scriptsize \par}

\begin{center}
{\scriptsize }%
\begin{minipage}[t]{0.33\textwidth}%
\begin{center}
{\scriptsize \includegraphics[scale=0.15]{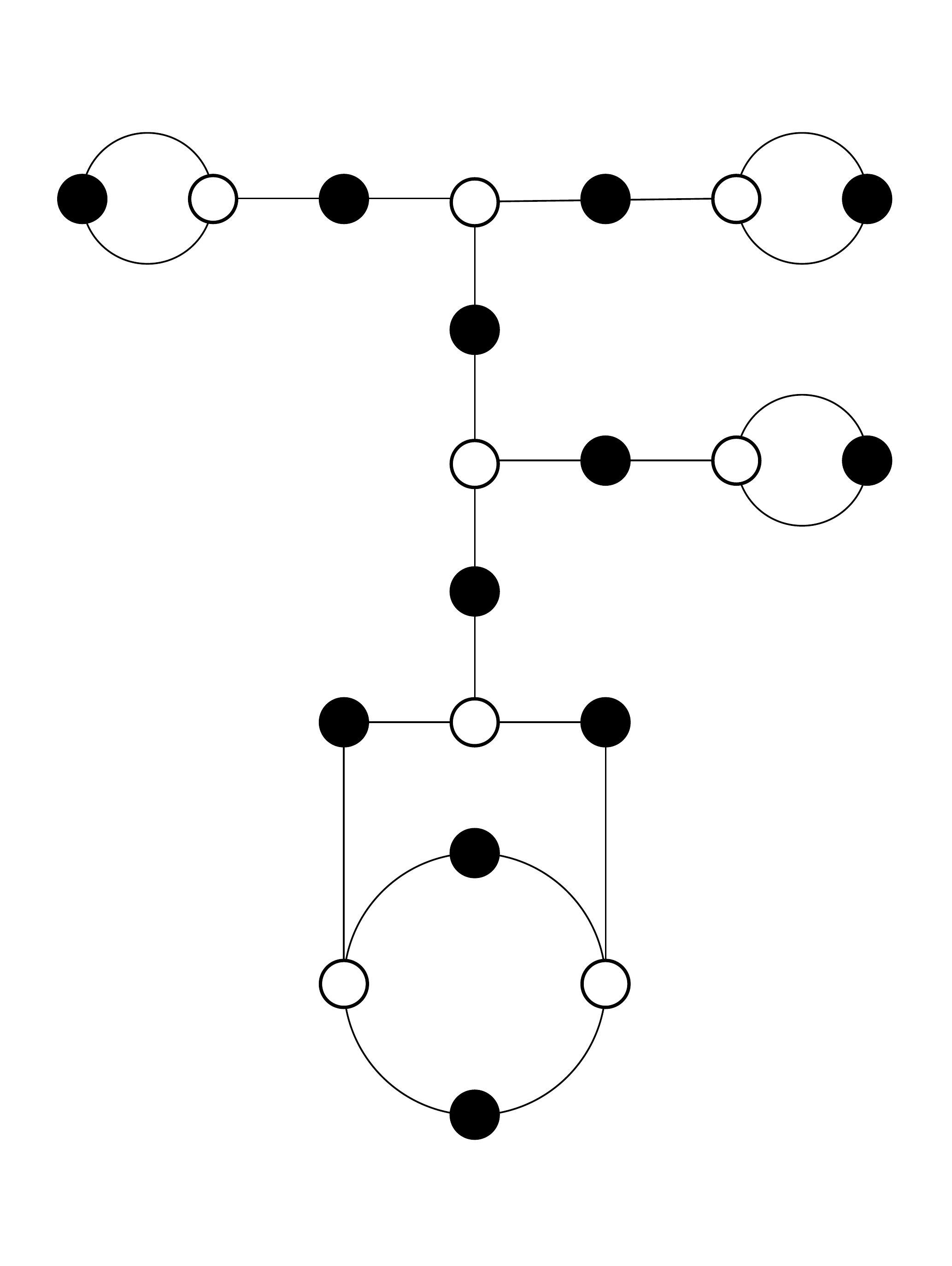}}
\par\end{center}{\scriptsize \par}

\begin{center}
{\scriptsize $16,3,2,1,1,1\;\left(\sqrt{-2}\right)$}
\par\end{center}%
\end{minipage}{\scriptsize }%
\begin{minipage}[t]{0.33\textwidth}%
\begin{center}
{\scriptsize \includegraphics[scale=0.15]{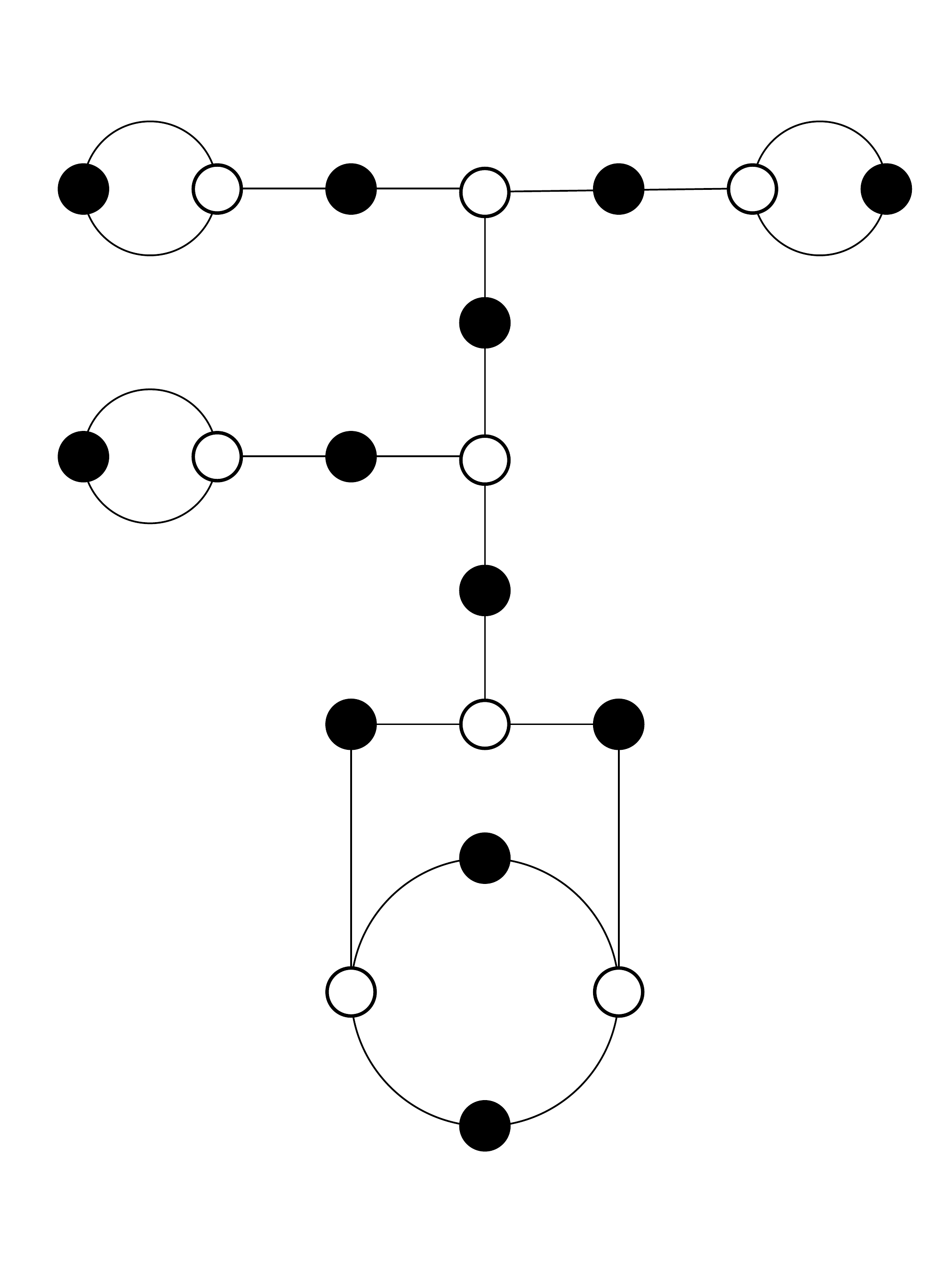}}
\par\end{center}{\scriptsize \par}

\begin{center}
{\scriptsize $16,3,2,1,1,1\;\left(\sqrt{-2}\right)$}
\par\end{center}%
\end{minipage}{\scriptsize }%
\begin{minipage}[t]{0.33\textwidth}%
\begin{center}
{\scriptsize \includegraphics[scale=0.15]{\string"PICT/16-2-2-2-1-1\string".pdf}}
\par\end{center}{\scriptsize \par}

\begin{center}
{\scriptsize $16,2,2,2,1,1\;\left(\mathbb{Q}\right)$}
\par\end{center}%
\end{minipage}
\par\end{center}{\scriptsize \par}

\begin{center}
{\scriptsize }%
\begin{minipage}[t]{0.33\textwidth}%
\begin{center}
{\scriptsize \includegraphics[scale=0.15]{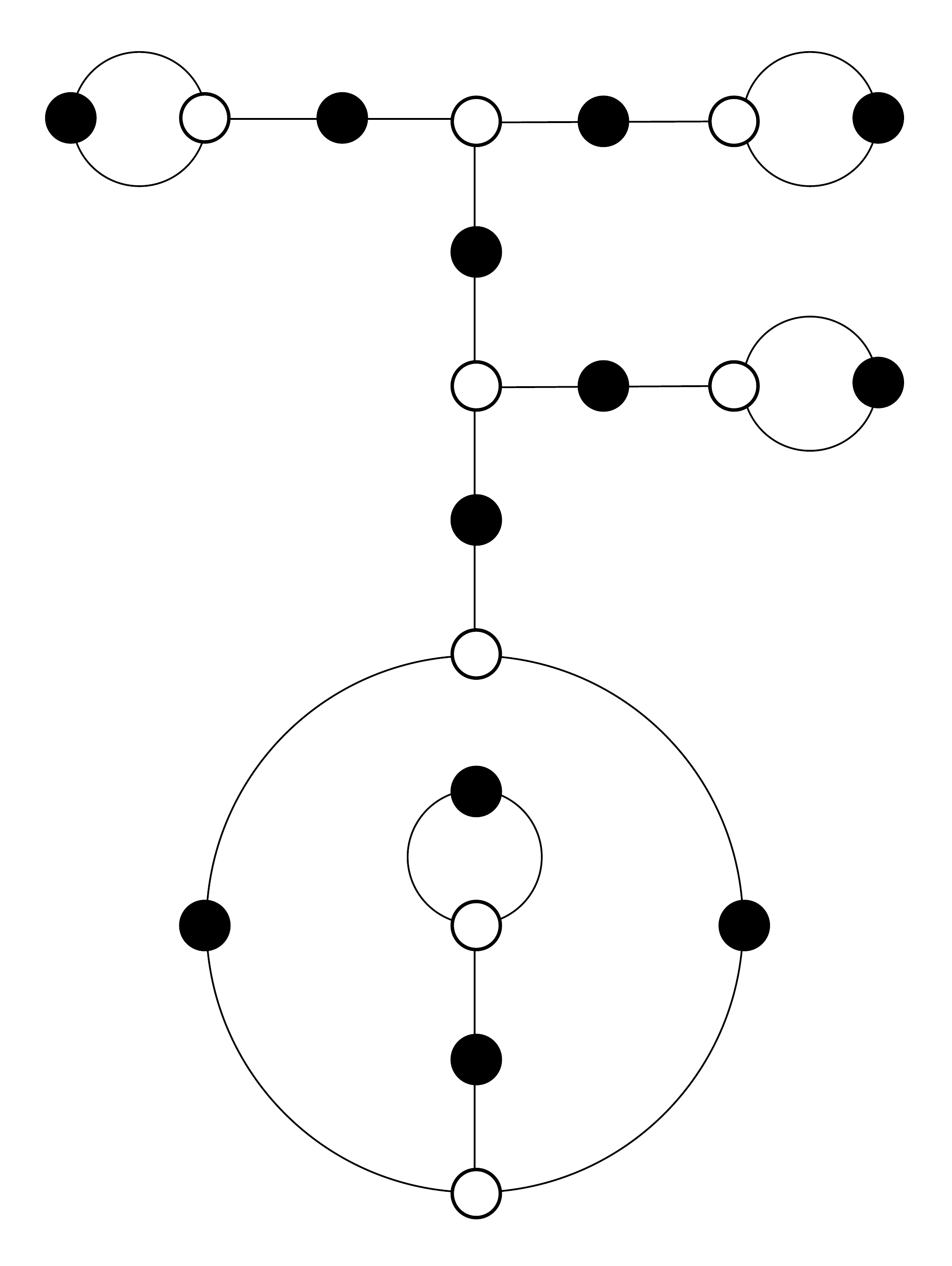}}
\par\end{center}{\scriptsize \par}

\begin{center}
{\scriptsize $15,5,1,1,1,1\;\left(\sqrt{-15}\right)$}
\par\end{center}%
\end{minipage}{\scriptsize }%
\begin{minipage}[t]{0.33\textwidth}%
\begin{center}
{\scriptsize \includegraphics[scale=0.15]{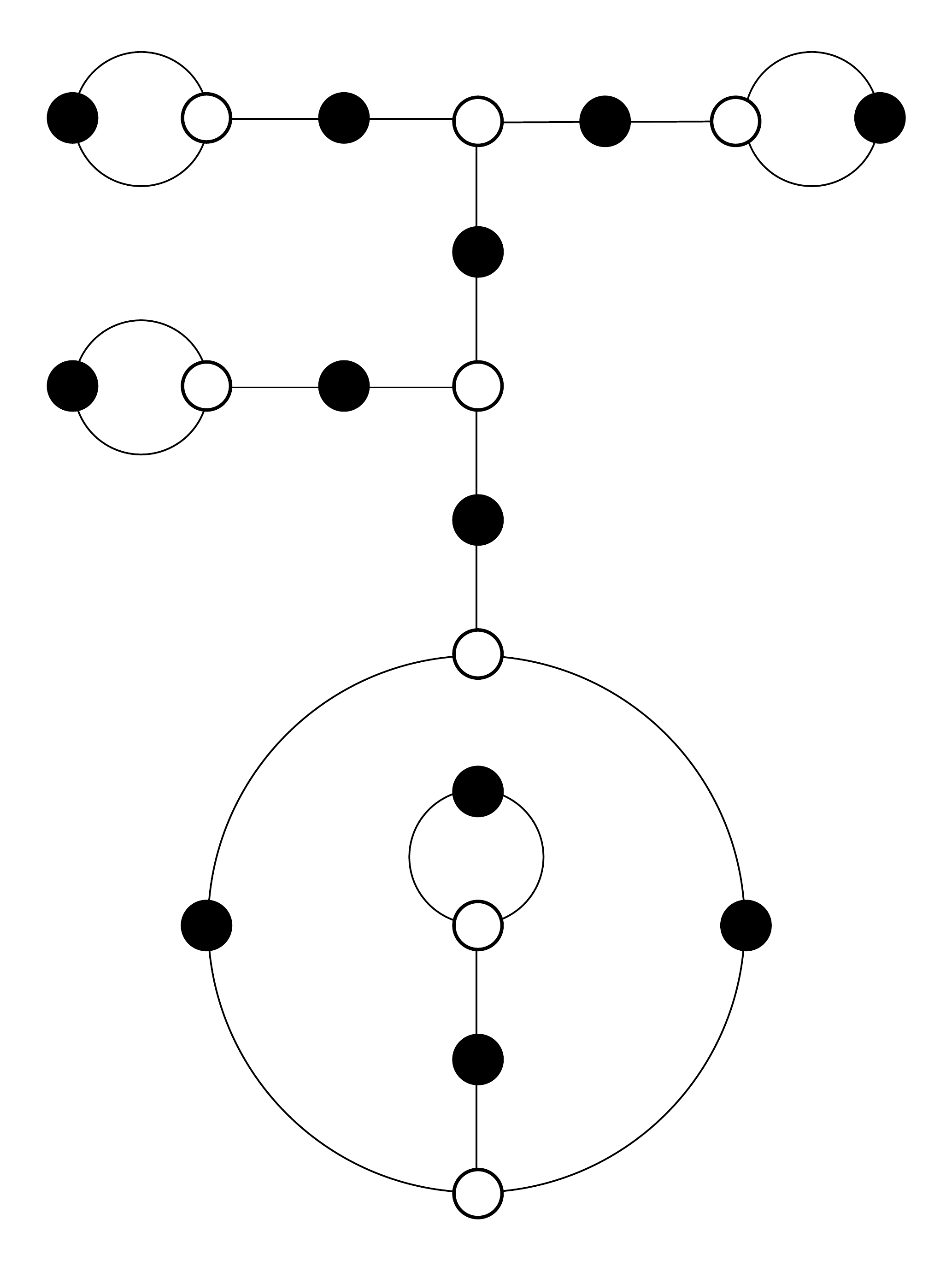}}
\par\end{center}{\scriptsize \par}

\begin{center}
{\scriptsize $15,5,1,1,1,1\;\left(\sqrt{-15}\right)$}
\par\end{center}%
\end{minipage}{\scriptsize }%
\begin{minipage}[t]{0.33\textwidth}%
\begin{center}
{\scriptsize \includegraphics[scale=0.15]{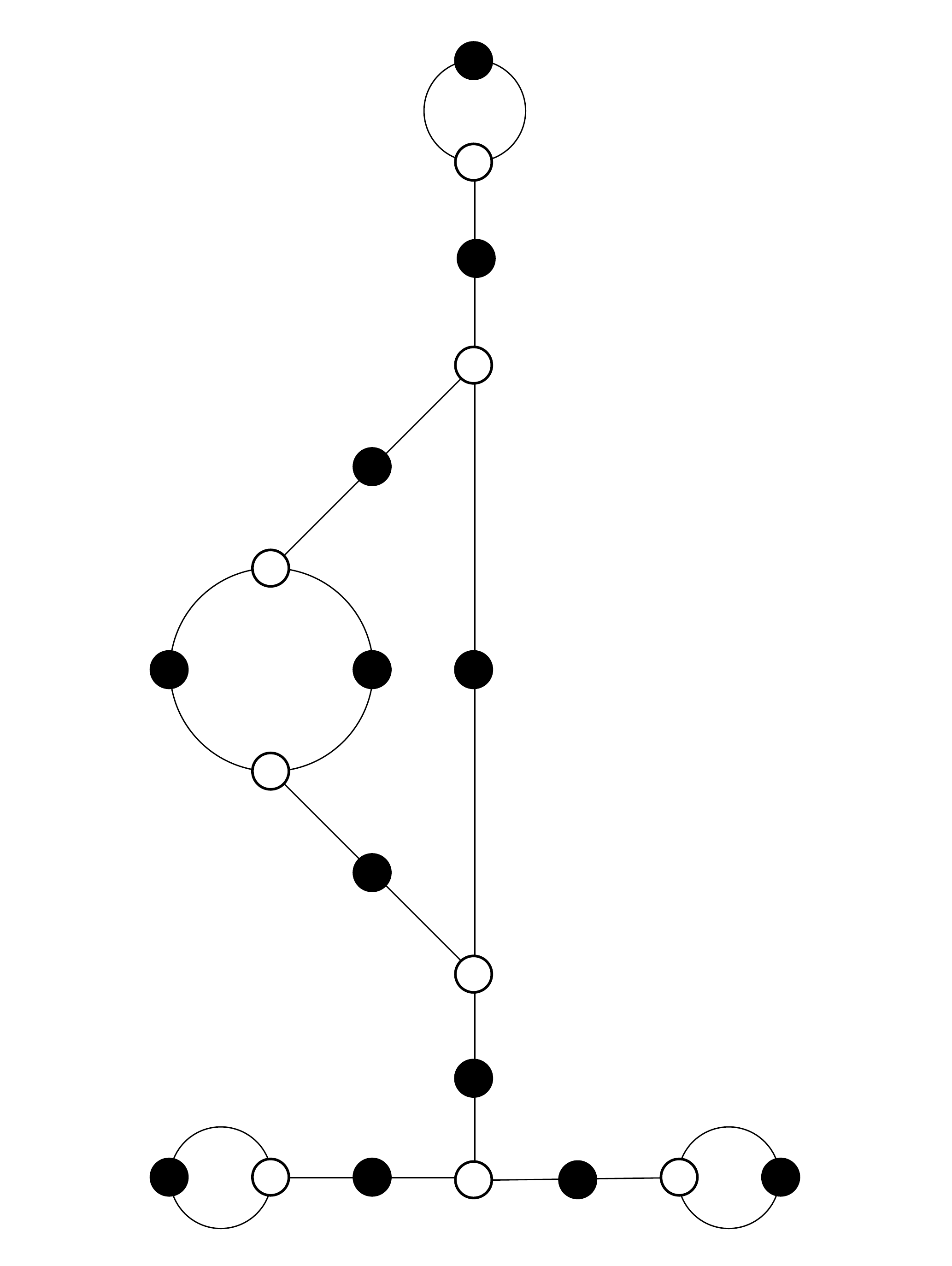}}
\par\end{center}{\scriptsize \par}

\begin{center}
{\scriptsize $15,4,2,1,1,1\;\left(\sqrt{-15}\right)$}
\par\end{center}%
\end{minipage}
\par\end{center}{\scriptsize \par}

\begin{center}
{\scriptsize }%
\begin{minipage}[t]{0.33\textwidth}%
\begin{center}
{\scriptsize \includegraphics[scale=0.15]{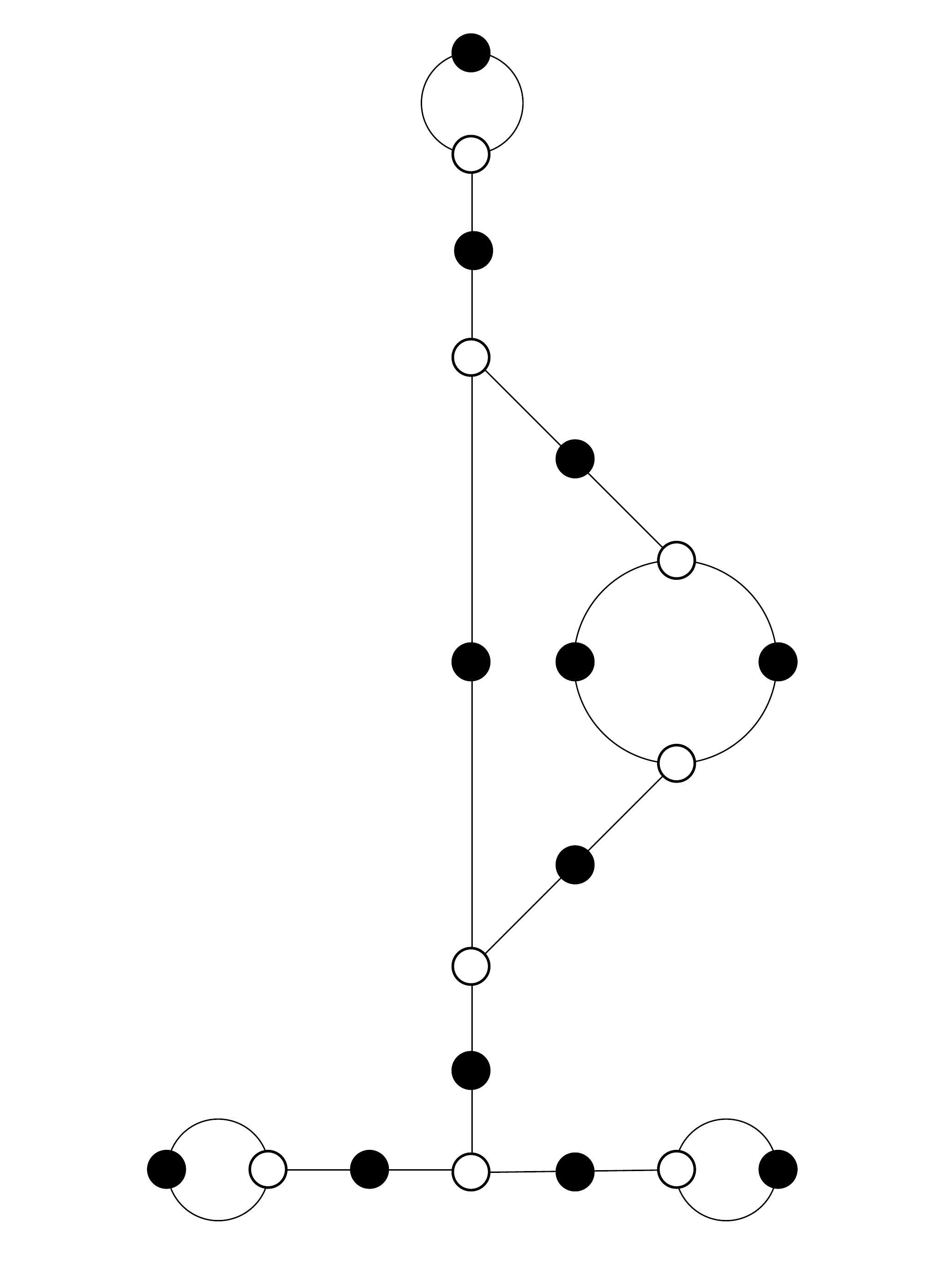}}
\par\end{center}{\scriptsize \par}

\begin{center}
{\scriptsize $15,4,2,1,1,1\;\left(\sqrt{-15}\right)$}
\par\end{center}%
\end{minipage}{\scriptsize }%
\begin{minipage}[t]{0.33\textwidth}%
\begin{center}
{\scriptsize \includegraphics[scale=0.15]{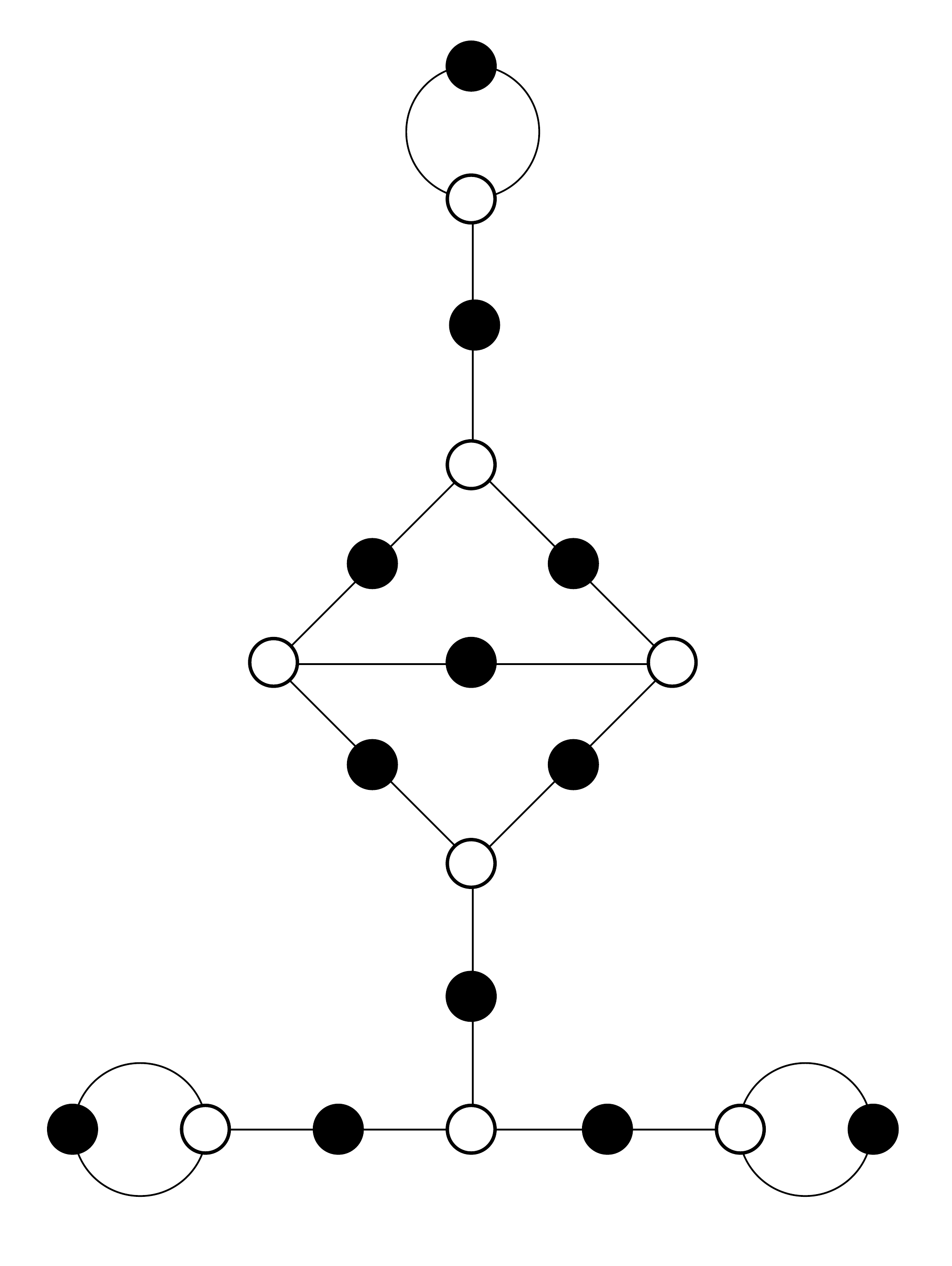}}
\par\end{center}{\scriptsize \par}

\begin{center}
{\scriptsize $15,3,3,1,1,1\;\left(\mathbb{Q}\right)$}
\par\end{center}%
\end{minipage}{\scriptsize }%
\begin{minipage}[t]{0.33\textwidth}%
\begin{center}
{\scriptsize \includegraphics[scale=0.15]{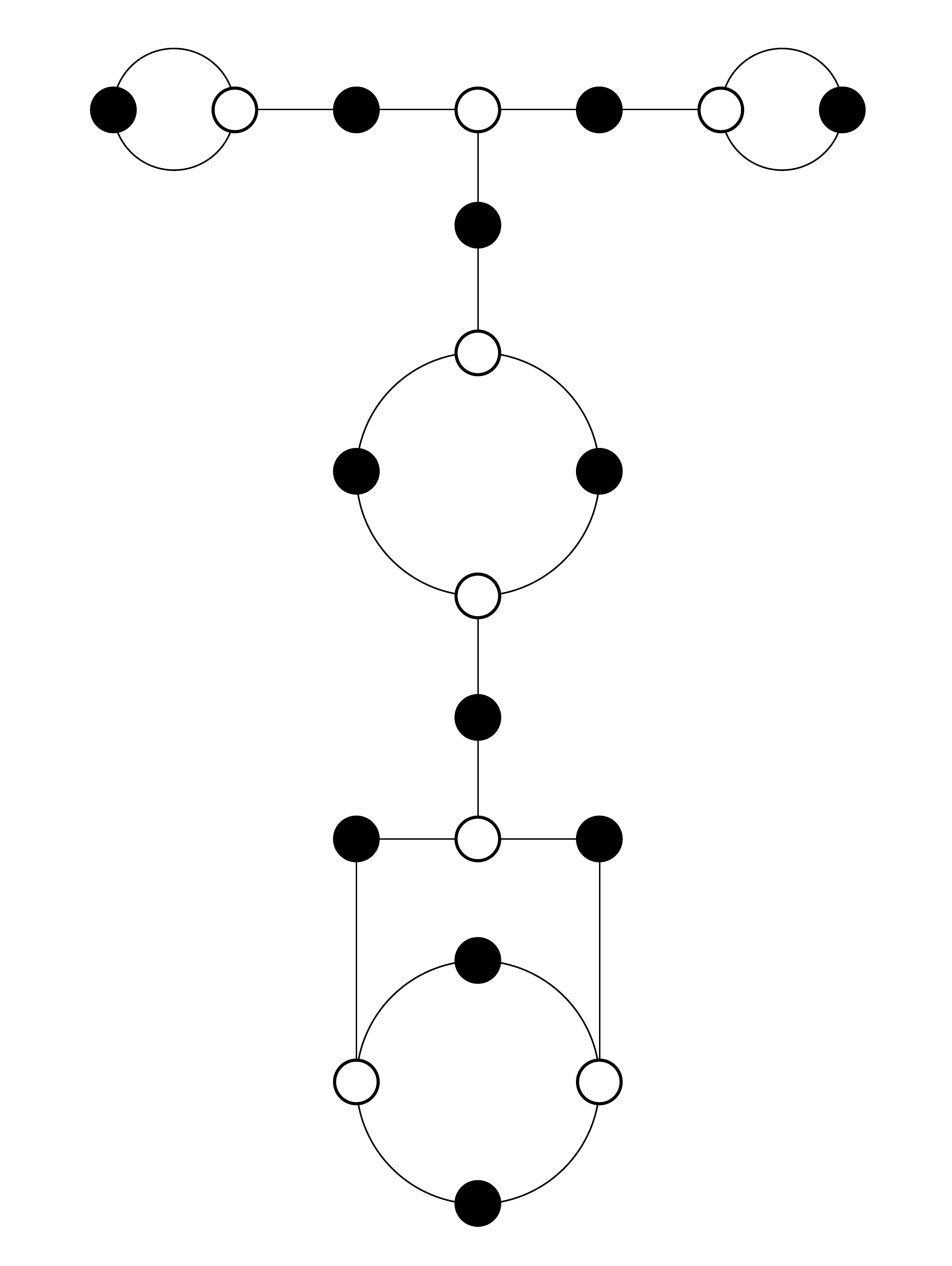}}
\par\end{center}{\scriptsize \par}

\begin{center}
{\scriptsize $15,3,2,2,1,1\;\left(\mathbb{Q}\right)$}
\par\end{center}%
\end{minipage}
\par\end{center}{\scriptsize \par}

\begin{center}
{\scriptsize }%
\begin{minipage}[t]{0.33\textwidth}%
\begin{center}
{\scriptsize \includegraphics[scale=0.15]{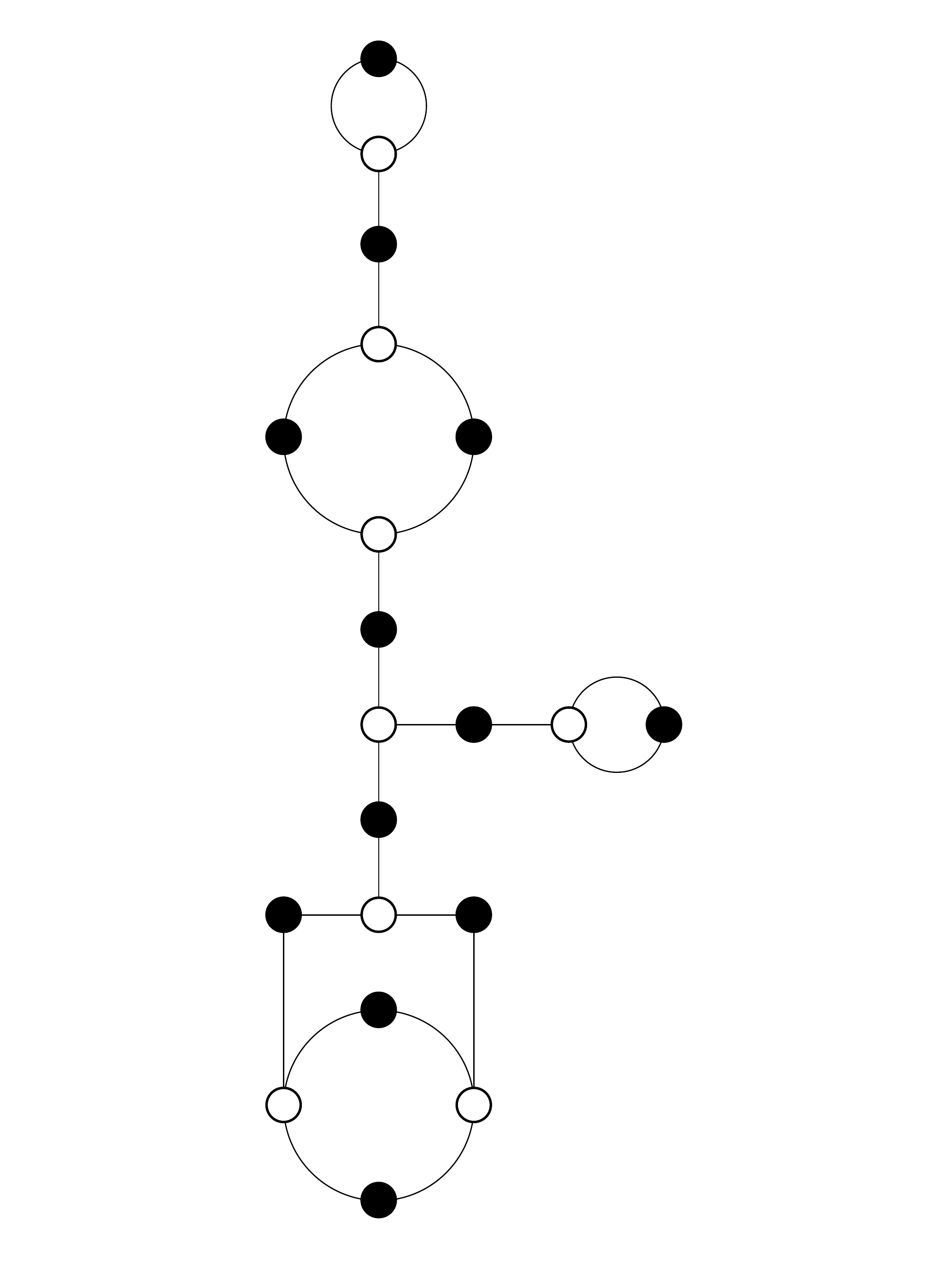}}
\par\end{center}{\scriptsize \par}

\begin{center}
{\scriptsize $15,3,2,2,1,1\;\left(\sqrt{-15}\right)$}
\par\end{center}%
\end{minipage}{\scriptsize }%
\begin{minipage}[t]{0.33\textwidth}%
\begin{center}
{\scriptsize \includegraphics[scale=0.15]{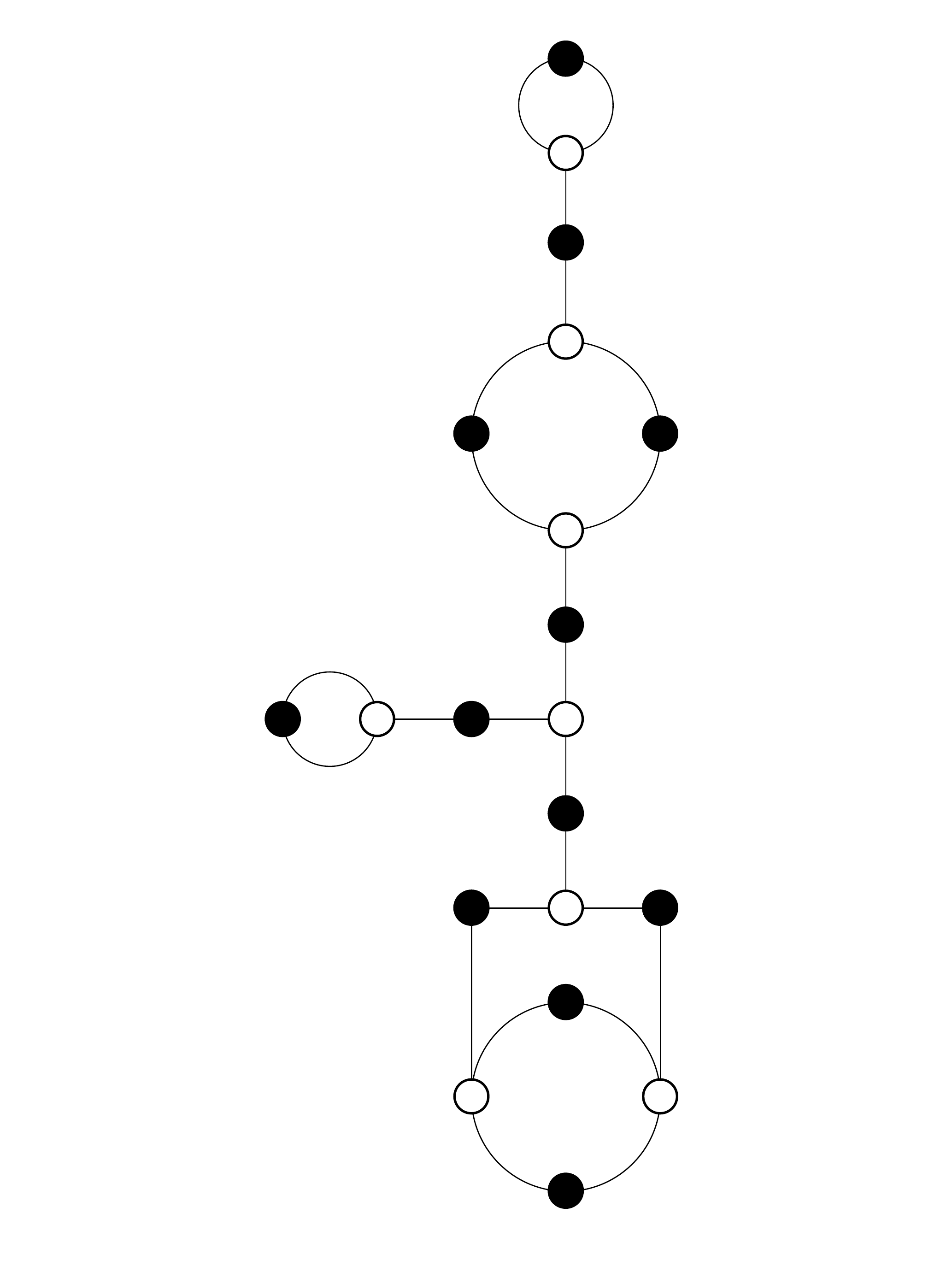}}
\par\end{center}{\scriptsize \par}

\begin{center}
{\scriptsize $15,3,2,2,1,1\;\left(\sqrt{-15}\right)$}
\par\end{center}%
\end{minipage}{\scriptsize }%
\begin{minipage}[t]{0.33\textwidth}%
\begin{center}
{\scriptsize \includegraphics[scale=0.15]{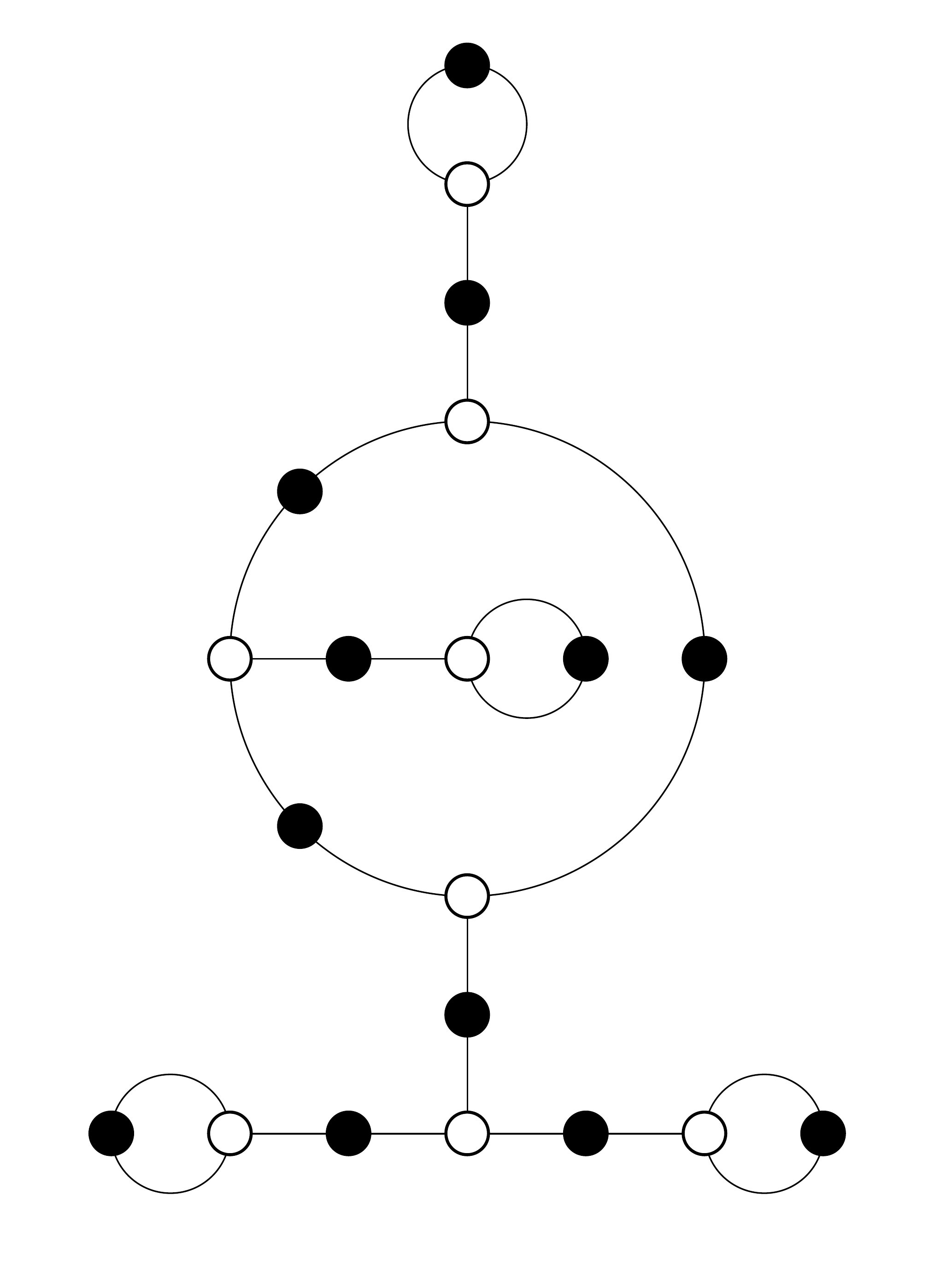}}
\par\end{center}{\scriptsize \par}

\begin{center}
{\scriptsize $14,6,1,1,1,1\;\left(\sqrt{-3}\right)$}
\par\end{center}%
\end{minipage}
\par\end{center}{\scriptsize \par}

\begin{center}
{\scriptsize }%
\begin{minipage}[t]{0.33\textwidth}%
\begin{center}
{\scriptsize \includegraphics[scale=0.15]{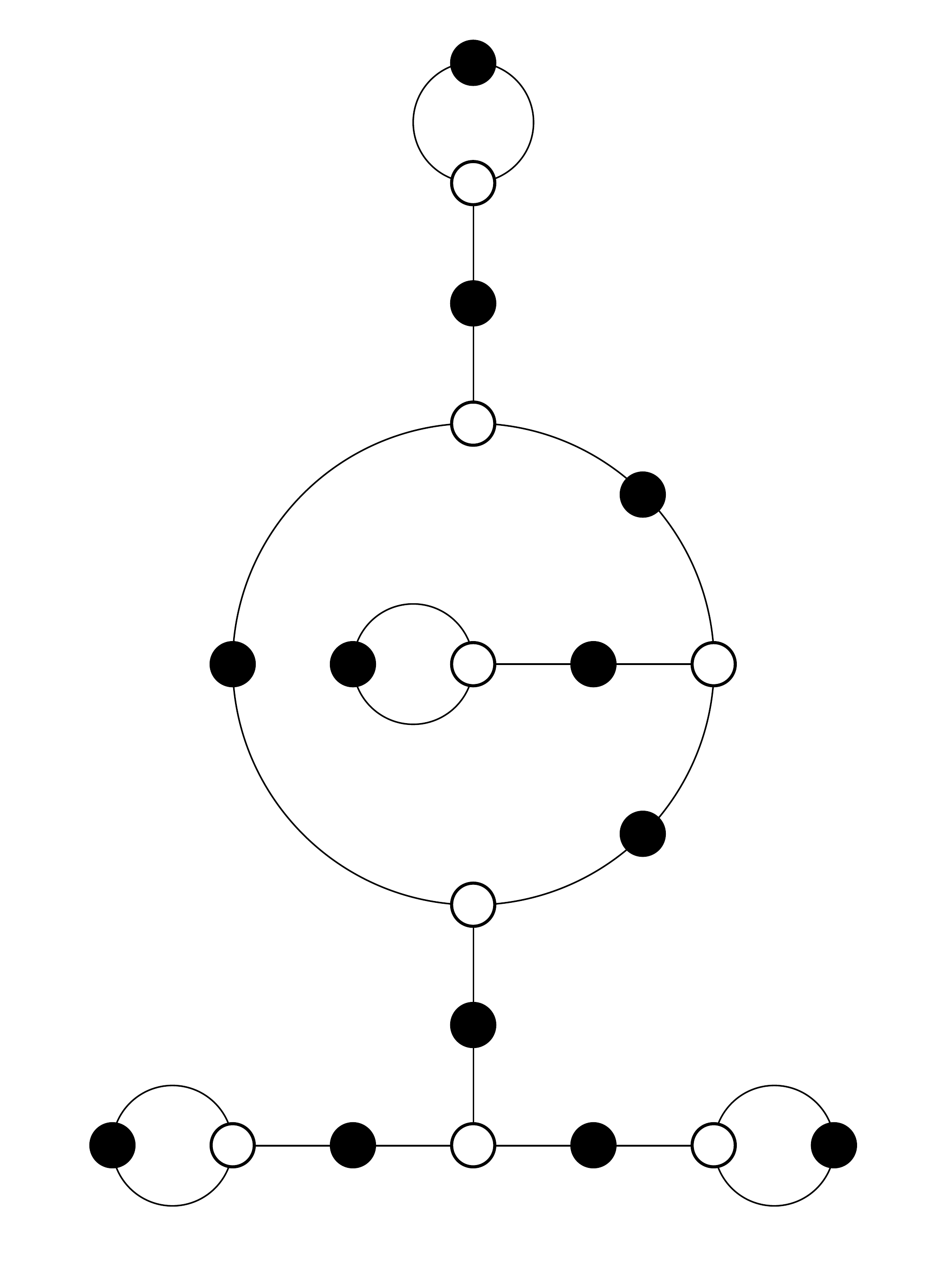}}
\par\end{center}{\scriptsize \par}

\begin{center}
{\scriptsize $14,6,1,1,1,1\;\left(\sqrt{-3}\right)$}
\par\end{center}%
\end{minipage}{\scriptsize }%
\begin{minipage}[t]{0.33\textwidth}%
\begin{center}
{\scriptsize \includegraphics[scale=0.15]{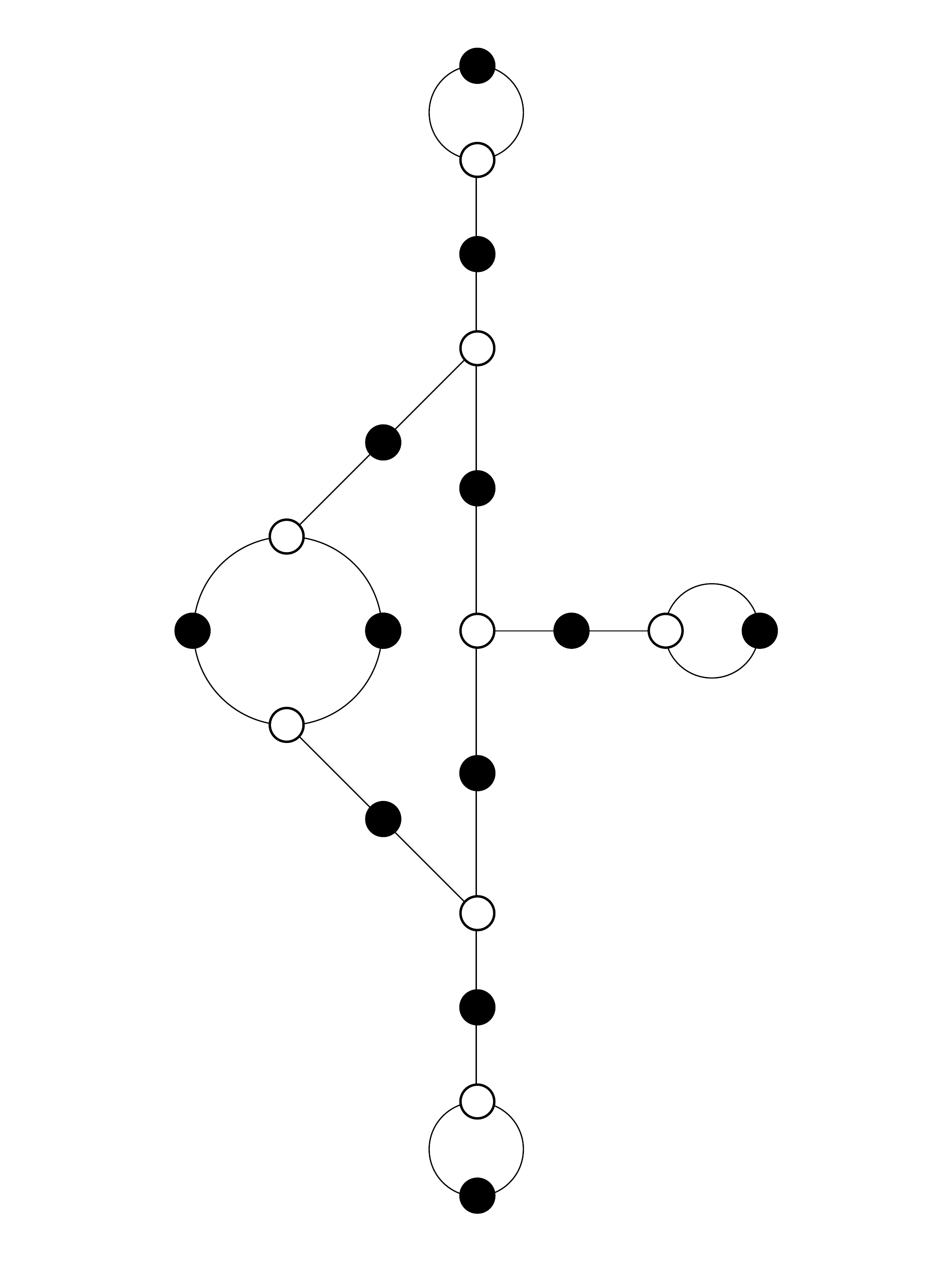}}
\par\end{center}{\scriptsize \par}

\begin{center}
{\scriptsize $14,5,2,1,1,1\;\left(\mathbb{Q}\right)$}
\par\end{center}%
\end{minipage}{\scriptsize }%
\begin{minipage}[t]{0.33\textwidth}%
\begin{center}
{\scriptsize \includegraphics[scale=0.15]{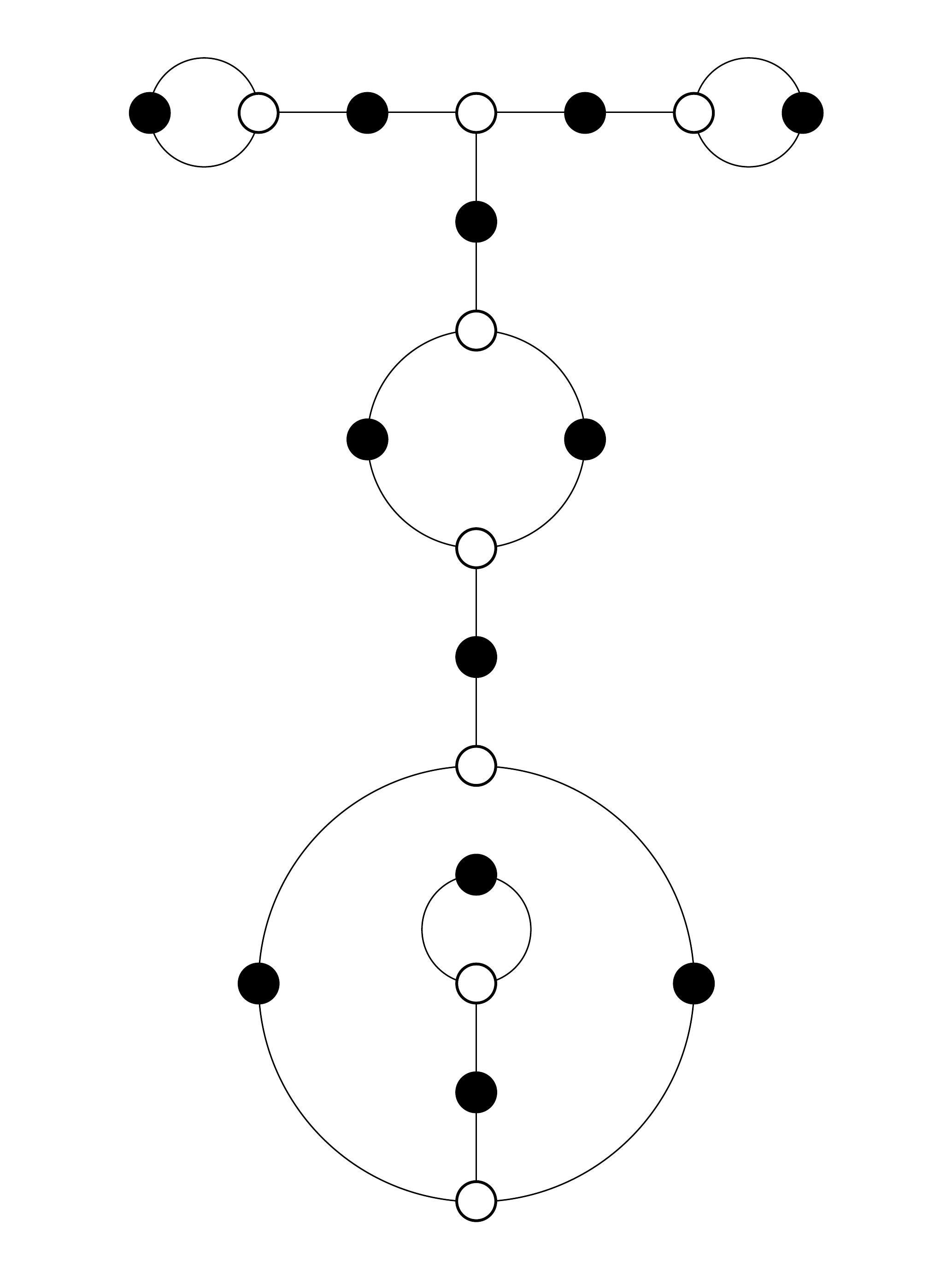}}
\par\end{center}{\scriptsize \par}

\begin{center}
{\scriptsize $14,5,2,1,1,1\;\left(\mathrm{cubic}\right)$}
\par\end{center}%
\end{minipage}
\par\end{center}{\scriptsize \par}

\begin{center}
{\scriptsize }%
\begin{minipage}[t]{0.33\textwidth}%
\begin{center}
{\scriptsize \includegraphics[scale=0.15]{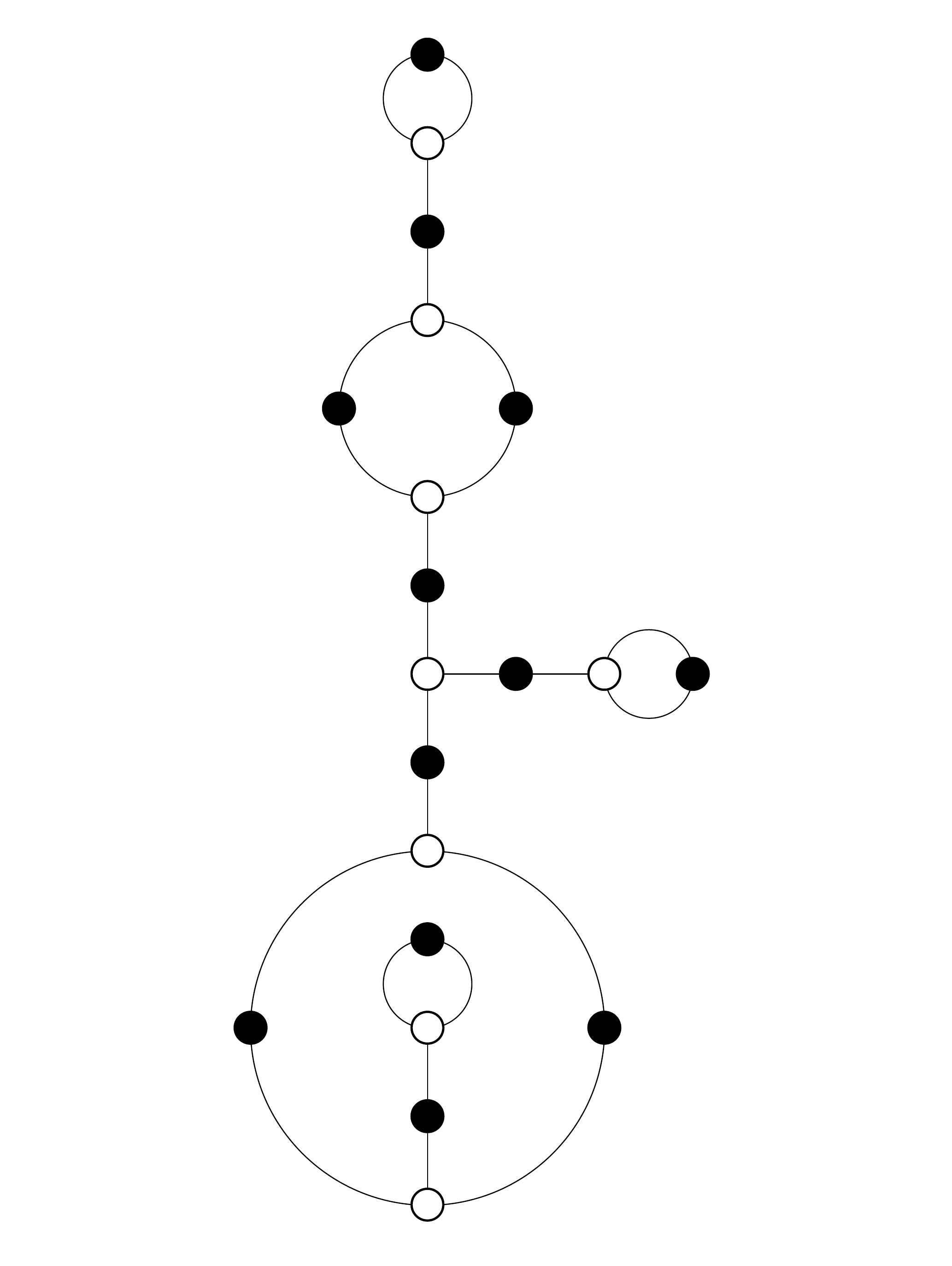}}
\par\end{center}{\scriptsize \par}

\begin{center}
{\scriptsize $14,5,2,1,1,1\;\left(\mathrm{cubic}\right)$}
\par\end{center}%
\end{minipage}{\scriptsize }%
\begin{minipage}[t]{0.33\textwidth}%
\begin{center}
{\scriptsize \includegraphics[scale=0.15]{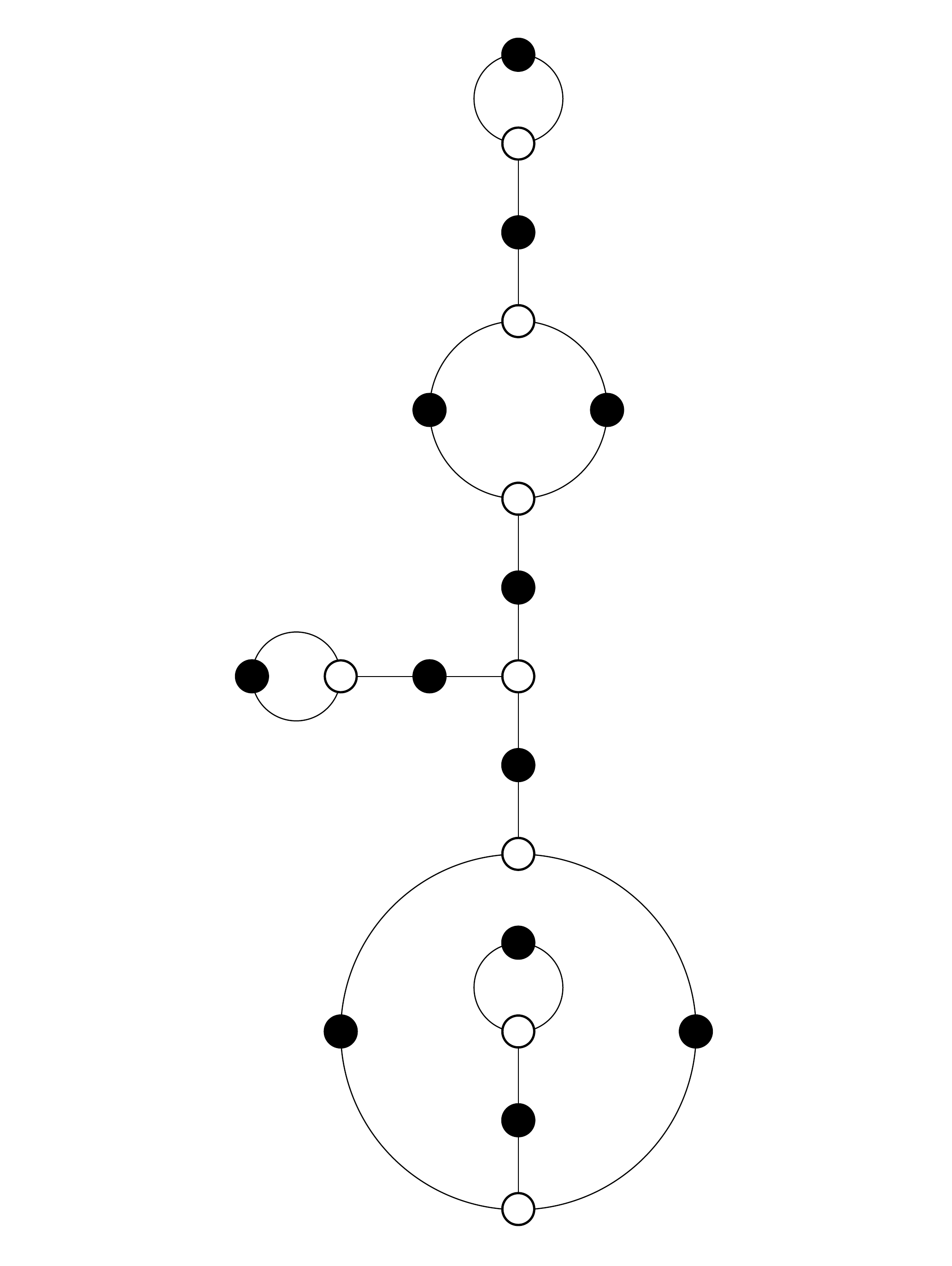}}
\par\end{center}{\scriptsize \par}

\begin{center}
{\scriptsize $14,5,2,1,1,1\;\left(\mathrm{cubic}\right)$}
\par\end{center}%
\end{minipage}{\scriptsize }%
\begin{minipage}[t]{0.33\textwidth}%
\begin{center}
{\scriptsize \includegraphics[scale=0.15]{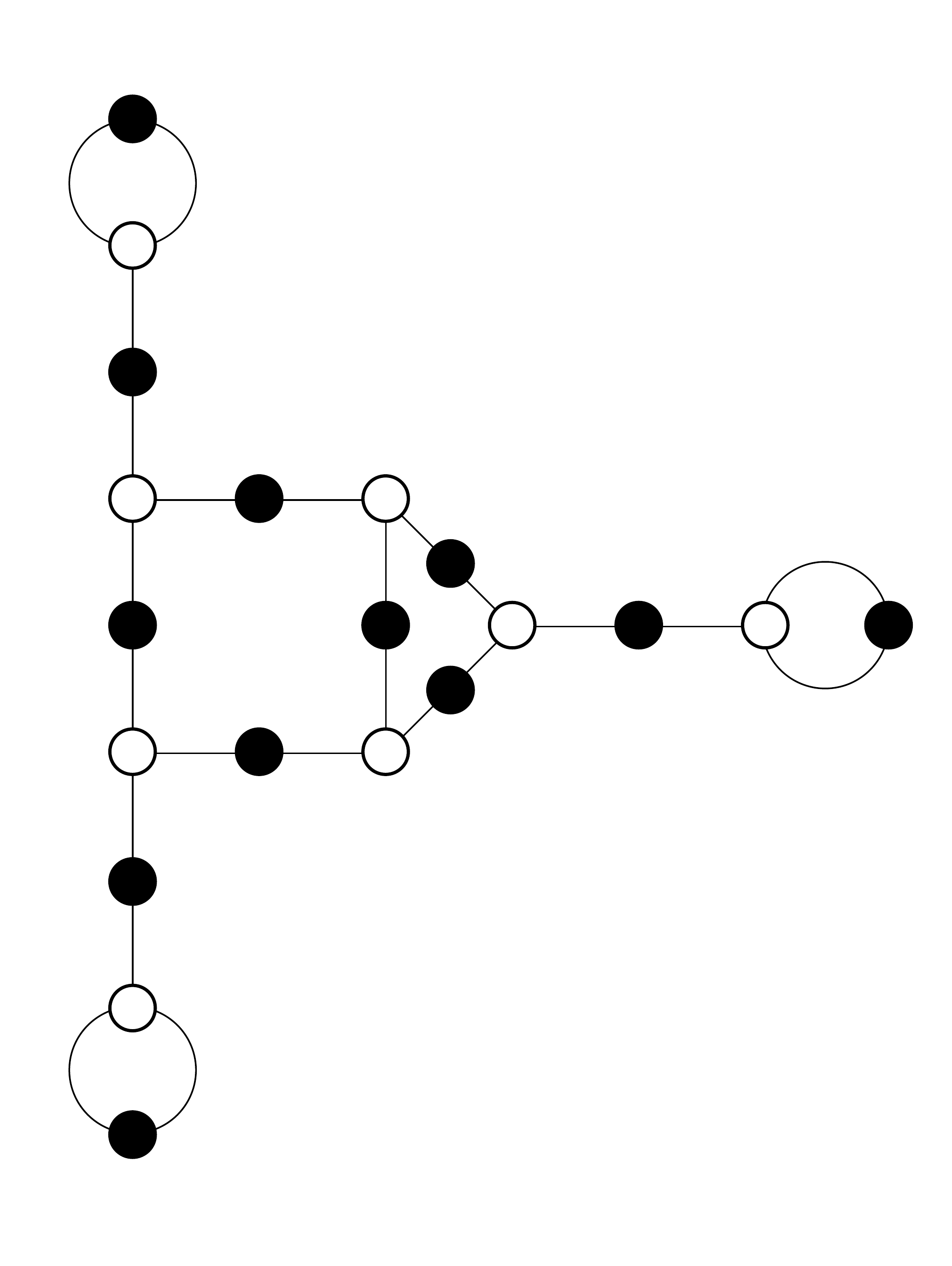}}
\par\end{center}{\scriptsize \par}

\begin{center}
{\scriptsize $14,4,3,1,1,1\;\left(\mathbb{Q}\right)$}
\par\end{center}%
\end{minipage}
\par\end{center}{\scriptsize \par}

\begin{center}
{\scriptsize }%
\begin{minipage}[t]{0.33\textwidth}%
\begin{center}
{\scriptsize \includegraphics[scale=0.15]{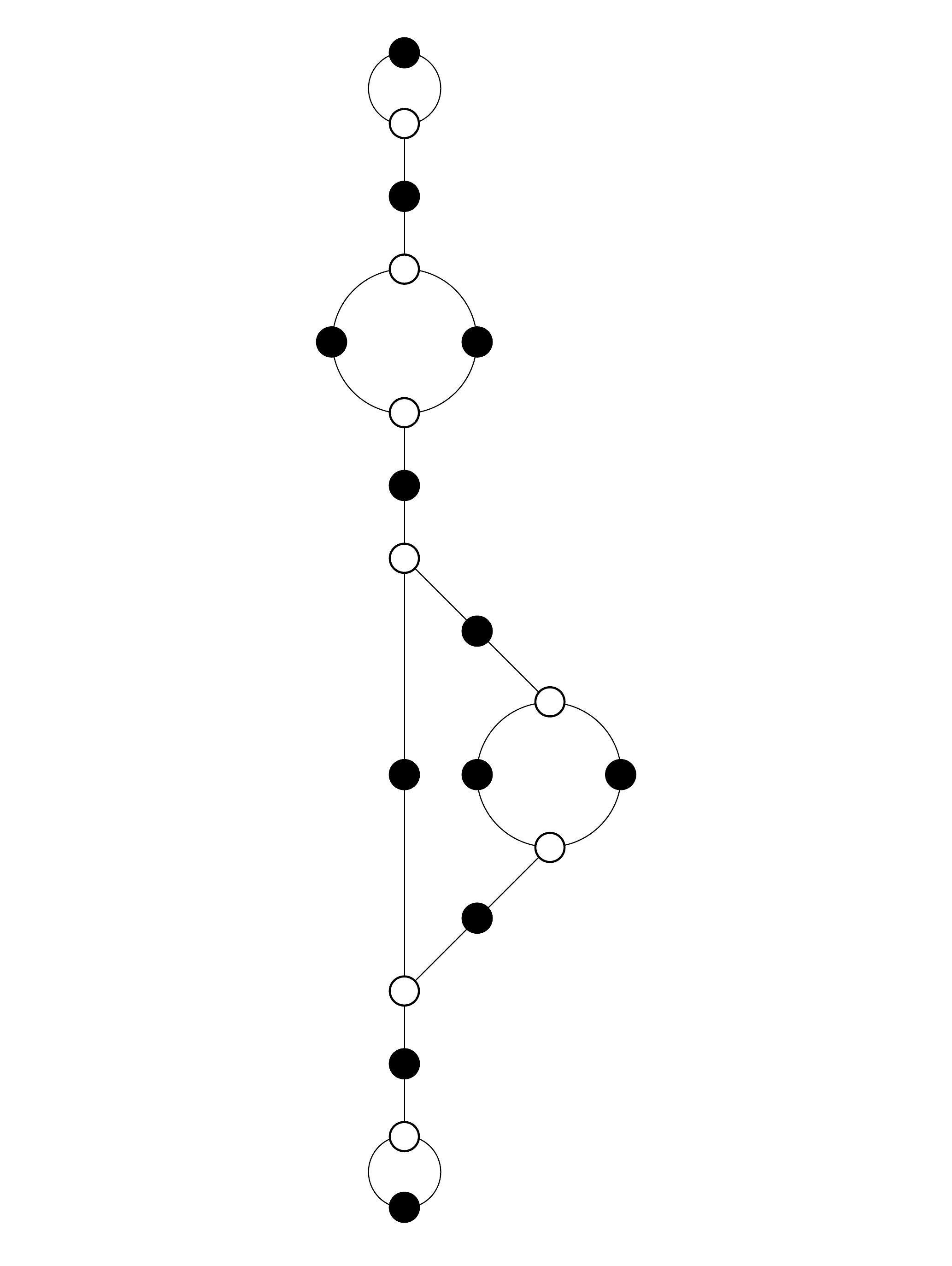}}
\par\end{center}{\scriptsize \par}

\begin{center}
{\scriptsize $14,4,2,2,1,1\;\left(\sqrt{-7}\right)$}
\par\end{center}%
\end{minipage}{\scriptsize }%
\begin{minipage}[t]{0.33\textwidth}%
\begin{center}
{\scriptsize \includegraphics[scale=0.15]{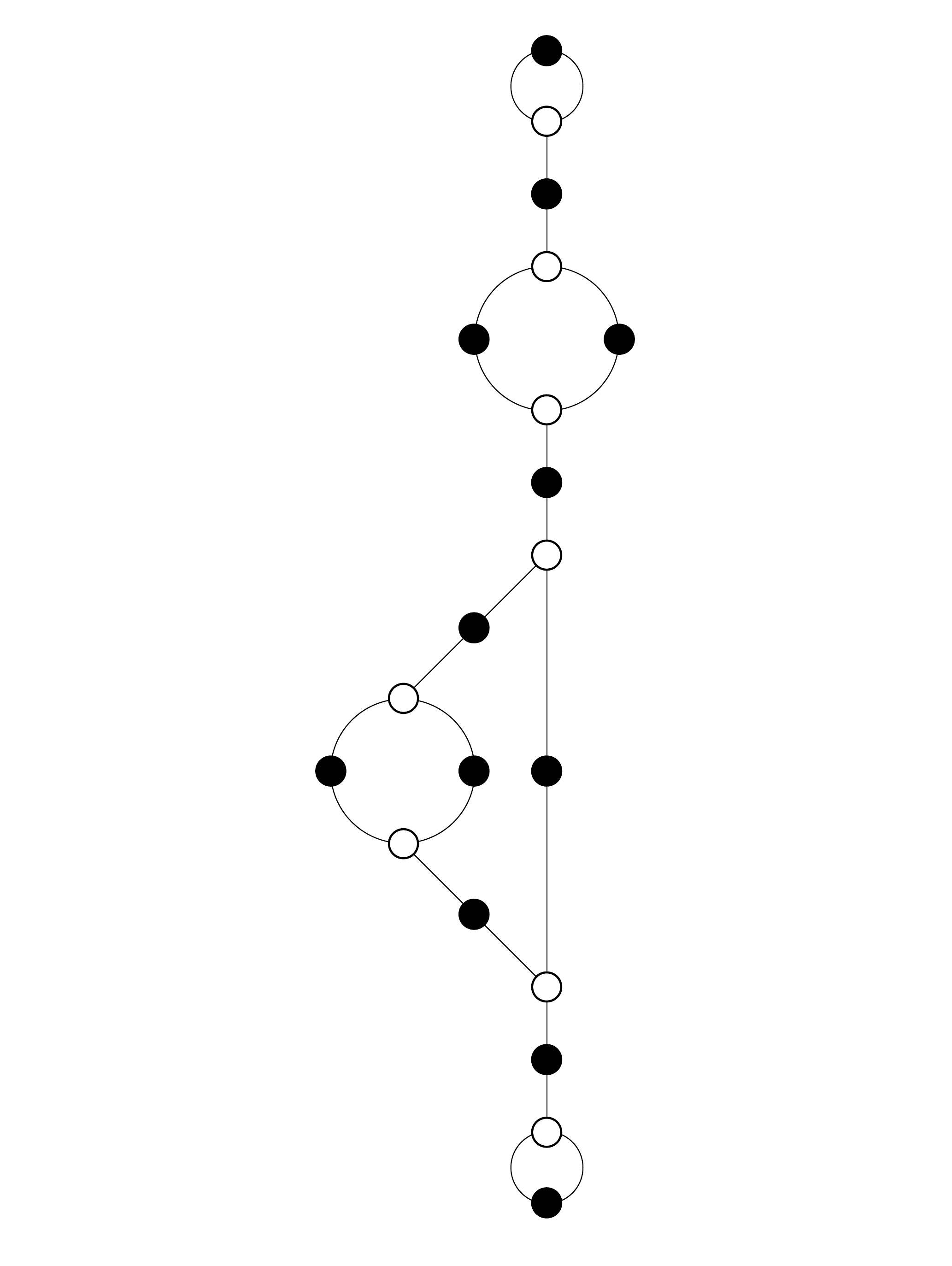}}
\par\end{center}{\scriptsize \par}

\begin{center}
{\scriptsize $14,4,2,2,1,1\;\left(\sqrt{-7}\right)$}
\par\end{center}%
\end{minipage}{\scriptsize }%
\begin{minipage}[t]{0.33\textwidth}%
\begin{center}
{\scriptsize \includegraphics[scale=0.15]{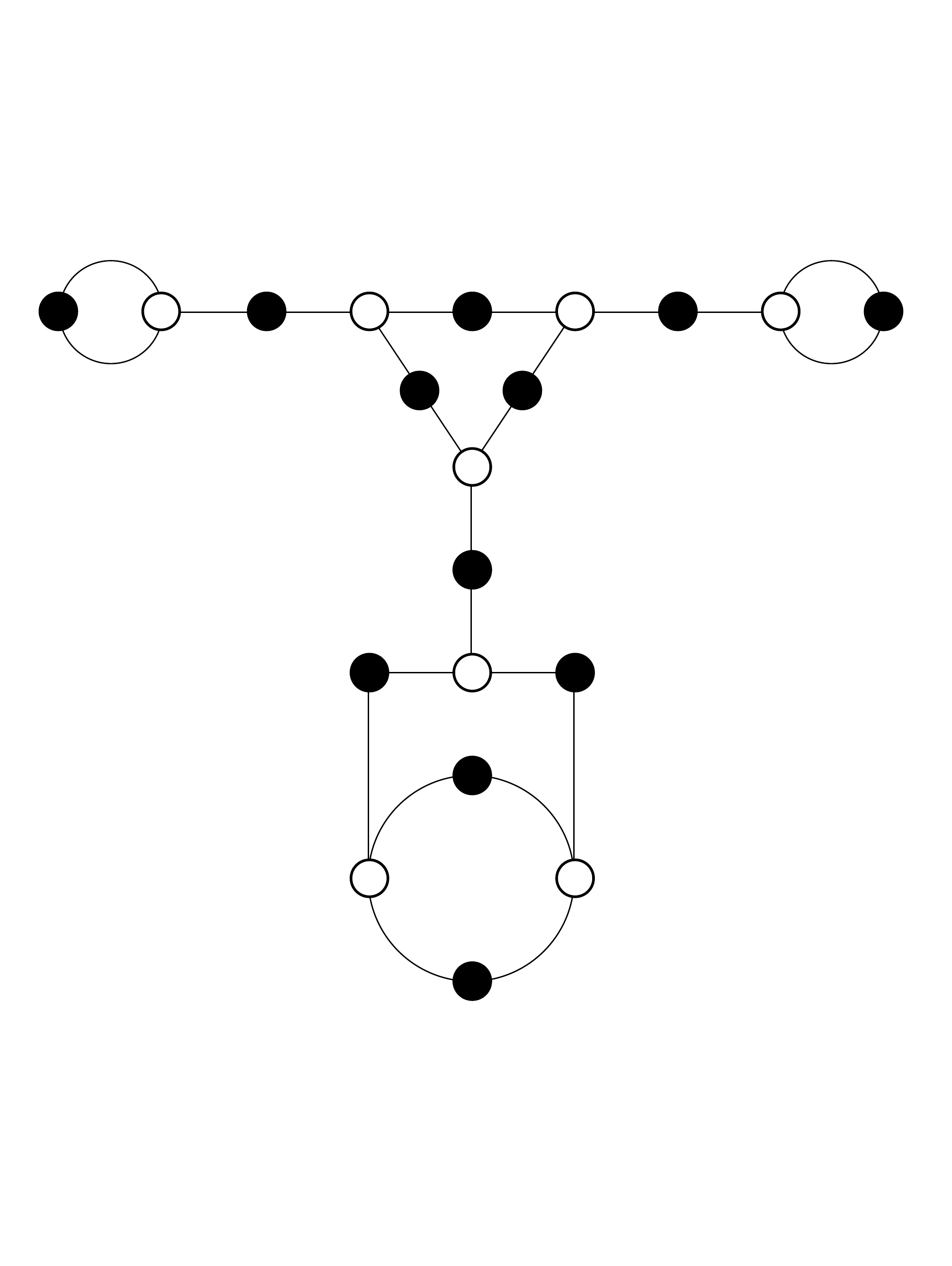}}
\par\end{center}{\scriptsize \par}

\begin{center}
{\scriptsize $14,3,3,2,1,1\;\left(\mathbb{Q}\right)$}
\par\end{center}%
\end{minipage}
\par\end{center}{\scriptsize \par}

\begin{center}
{\scriptsize }%
\begin{minipage}[t]{0.33\textwidth}%
\begin{center}
{\scriptsize \includegraphics[scale=0.15]{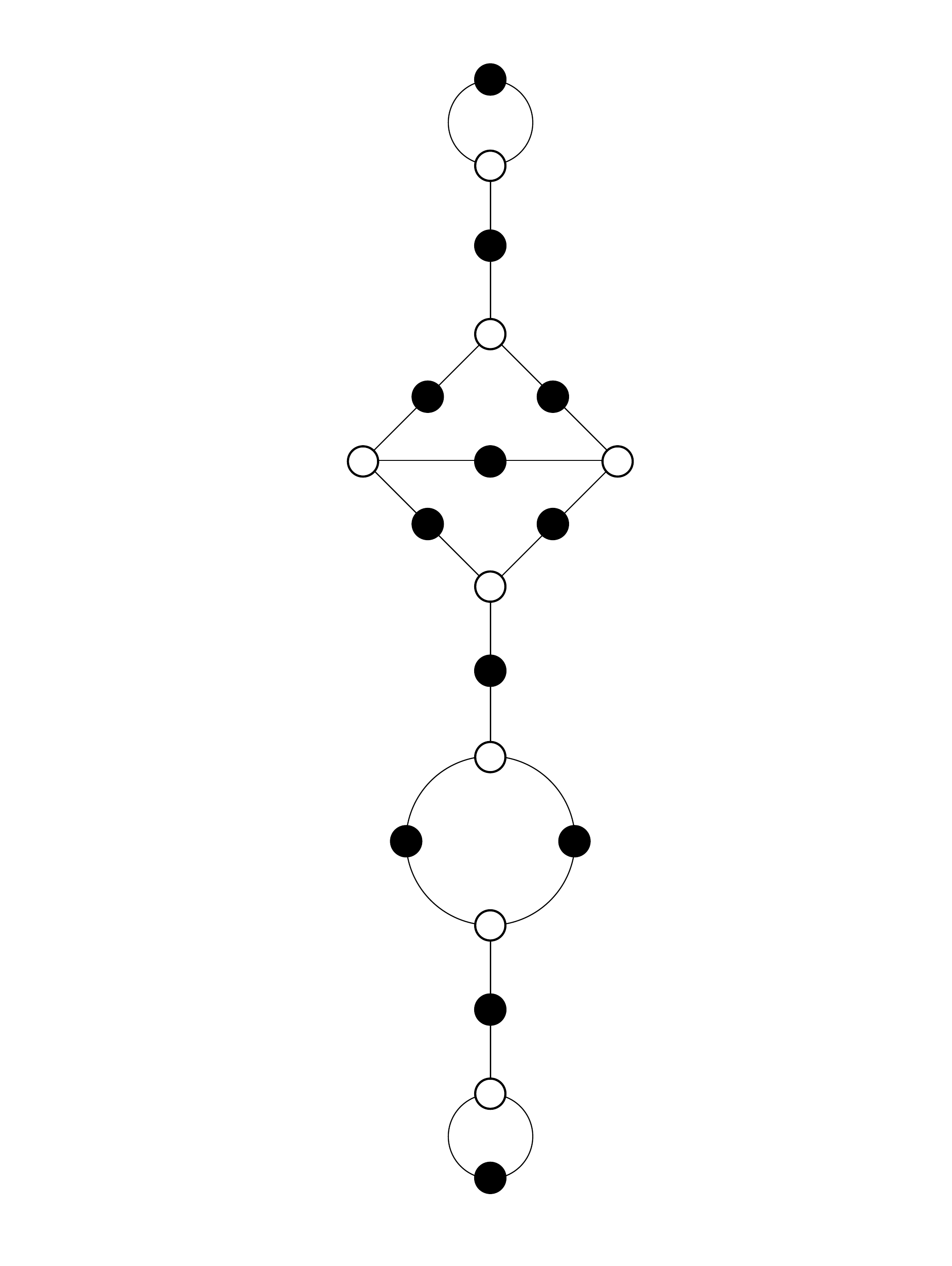}}
\par\end{center}{\scriptsize \par}

\begin{center}
{\scriptsize $14,3,3,2,1,1\;\left(\mathbb{Q}\right)$}
\par\end{center}%
\end{minipage}{\scriptsize }%
\begin{minipage}[t]{0.33\textwidth}%
\begin{center}
{\scriptsize \includegraphics[scale=0.15]{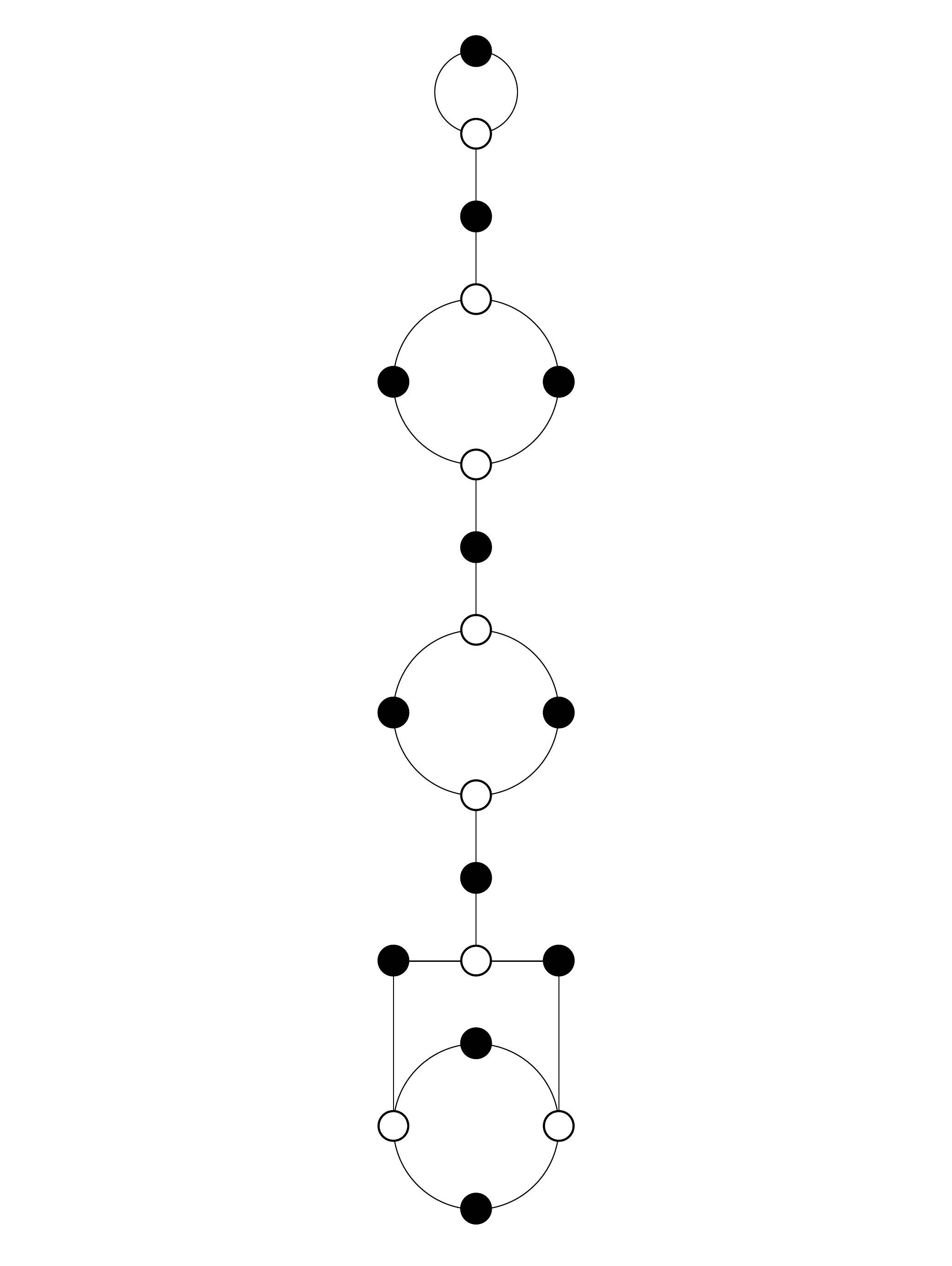}}
\par\end{center}{\scriptsize \par}

\begin{center}
{\scriptsize $14,3,2,2,2,1\;\left(\mathbb{Q}\right)$}
\par\end{center}%
\end{minipage}{\scriptsize }%
\begin{minipage}[t]{0.33\textwidth}%
\begin{center}
{\scriptsize \includegraphics[scale=0.15]{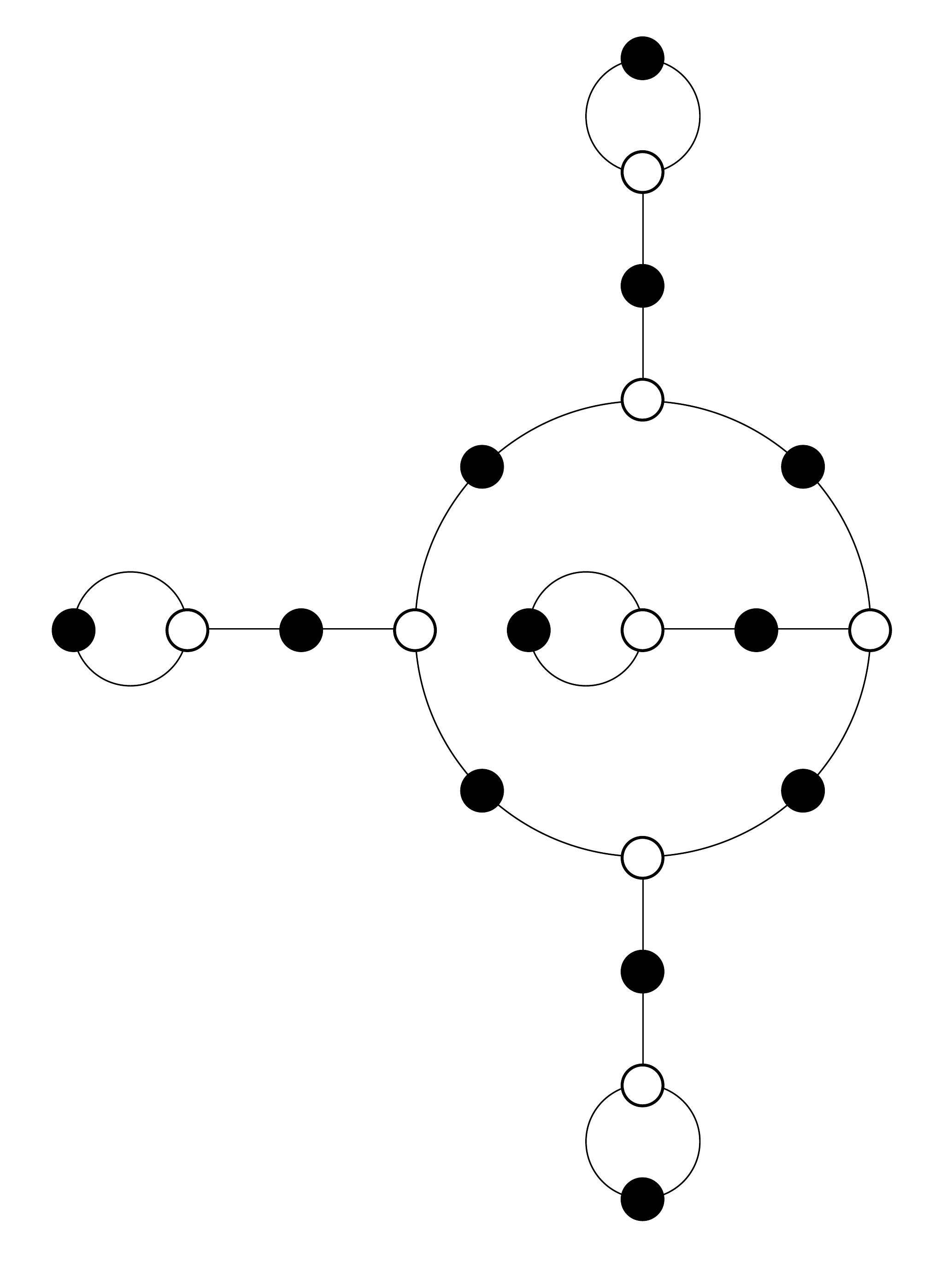}}
\par\end{center}{\scriptsize \par}

\begin{center}
{\scriptsize $13,7,1,1,1,1\;\left(\mathbb{Q}\right)$}
\par\end{center}%
\end{minipage}
\par\end{center}{\scriptsize \par}

\begin{center}
{\scriptsize }%
\begin{minipage}[t]{0.33\textwidth}%
\begin{center}
{\scriptsize \includegraphics[scale=0.15]{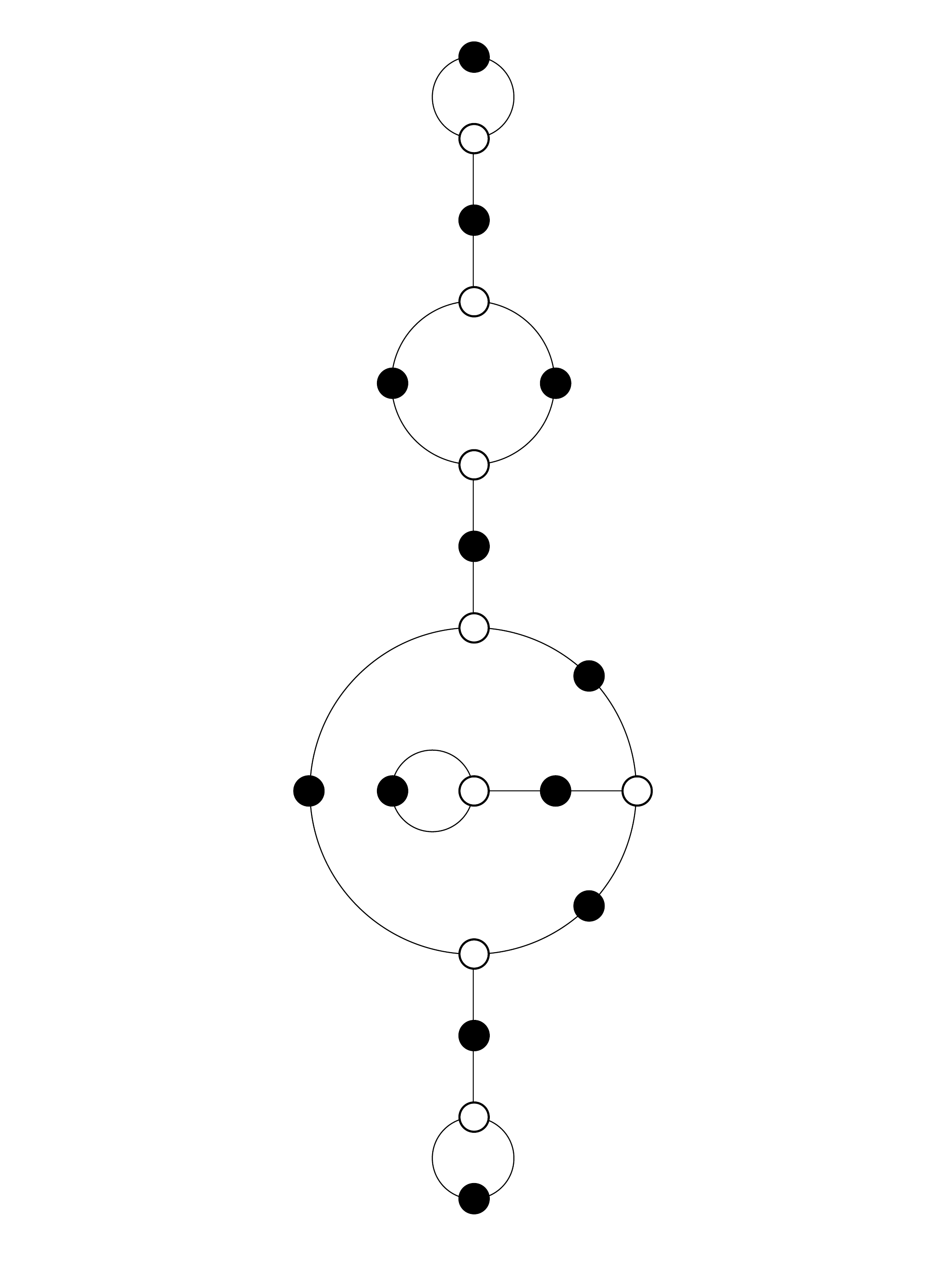}}
\par\end{center}{\scriptsize \par}

\begin{center}
{\scriptsize $13,6,2,1,1,1\;\left(\sqrt{-3}\right)$}
\par\end{center}%
\end{minipage}{\scriptsize }%
\begin{minipage}[t]{0.33\textwidth}%
\begin{center}
{\scriptsize \includegraphics[scale=0.15]{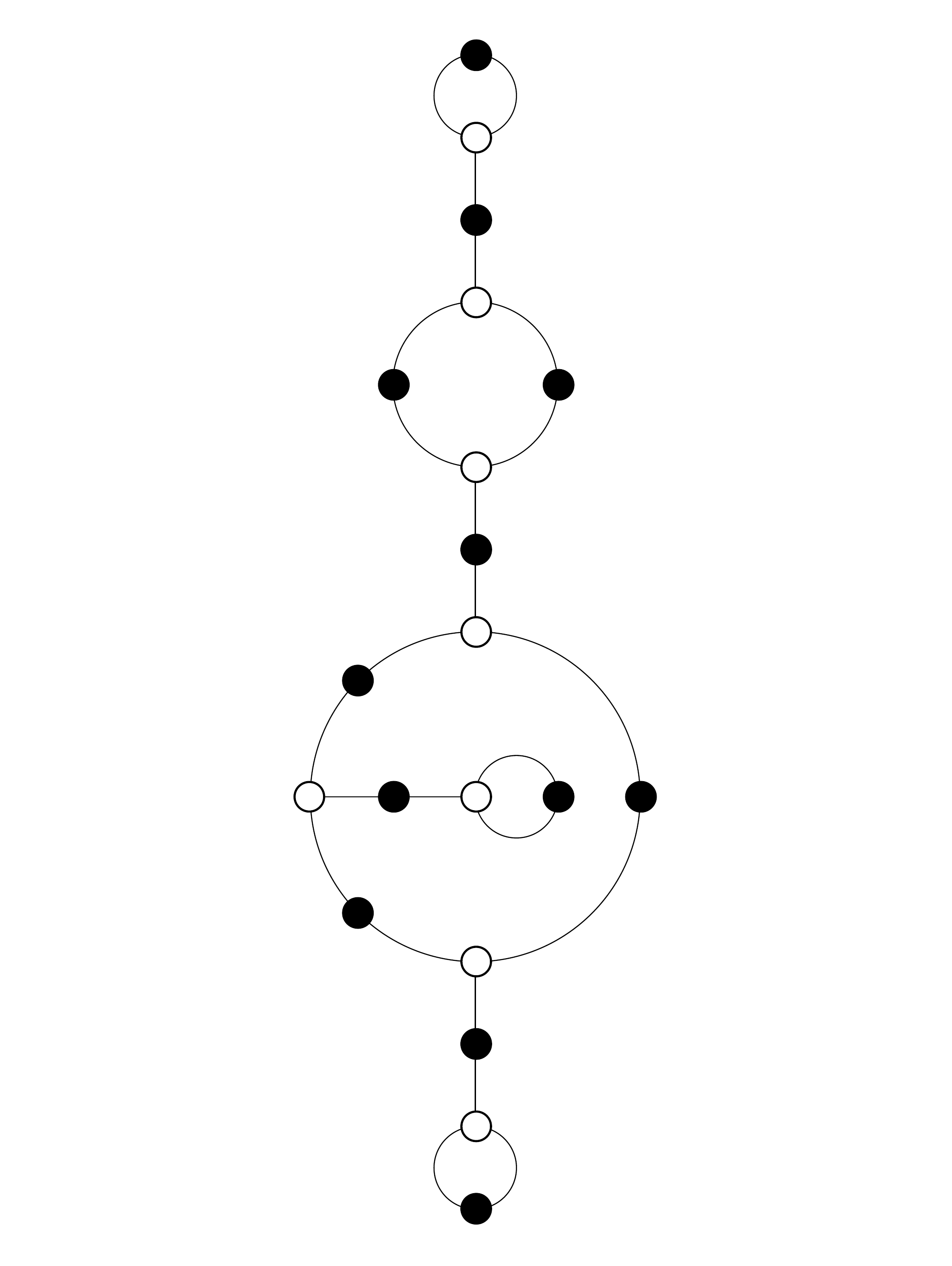}}
\par\end{center}{\scriptsize \par}

\begin{center}
{\scriptsize $13,6,2,1,1,1\;\left(\sqrt{-3}\right)$}
\par\end{center}%
\end{minipage}{\scriptsize }%
\begin{minipage}[t]{0.33\textwidth}%
\begin{center}
{\scriptsize \includegraphics[scale=0.15]{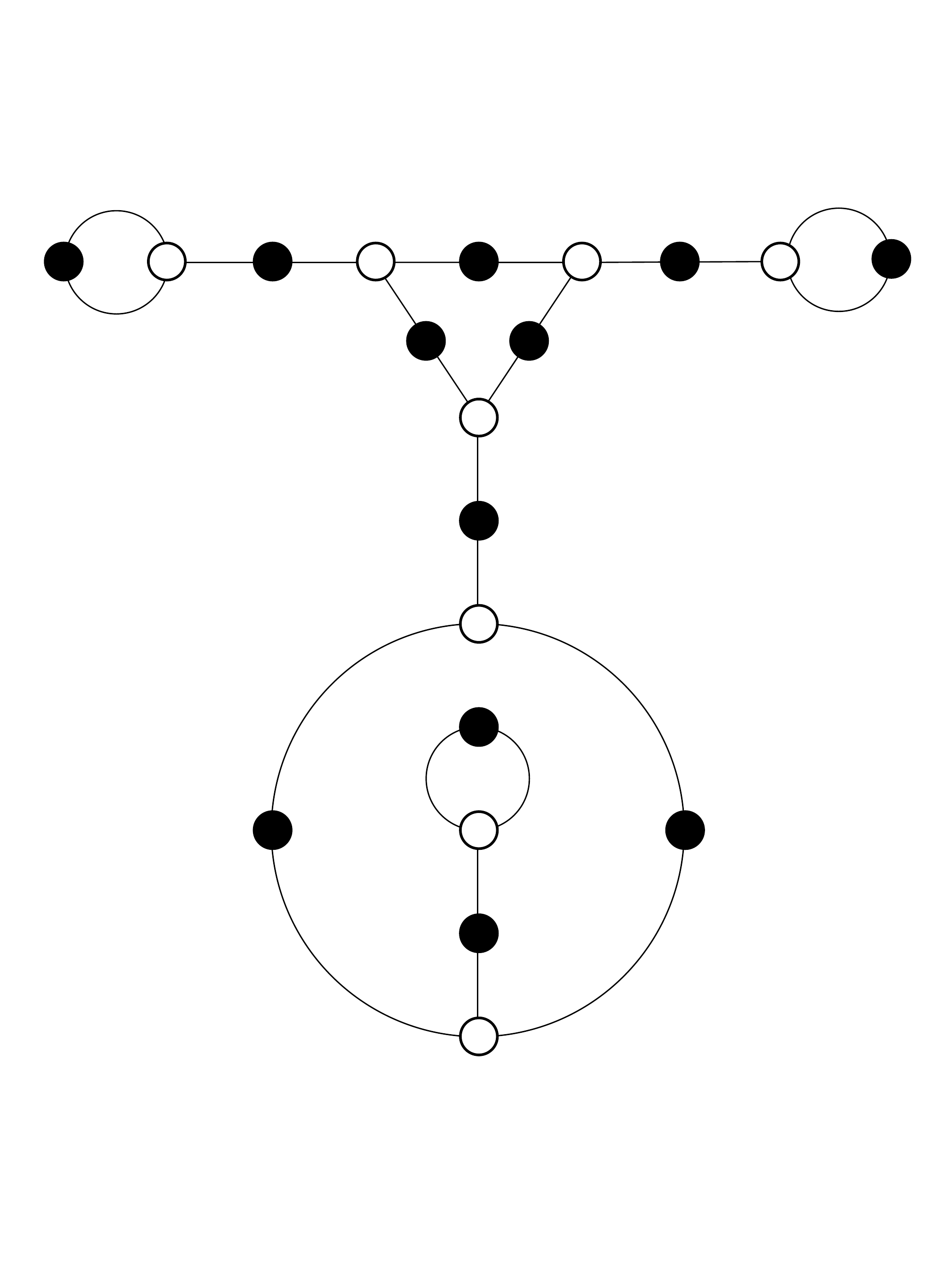}}
\par\end{center}{\scriptsize \par}

\begin{center}
{\scriptsize $13,5,3,1,1,1\;\left(\mathbb{Q}\right)$}
\par\end{center}%
\end{minipage}
\par\end{center}{\scriptsize \par}

\begin{center}
{\scriptsize }%
\begin{minipage}[t]{0.33\textwidth}%
\begin{center}
{\scriptsize \includegraphics[scale=0.15]{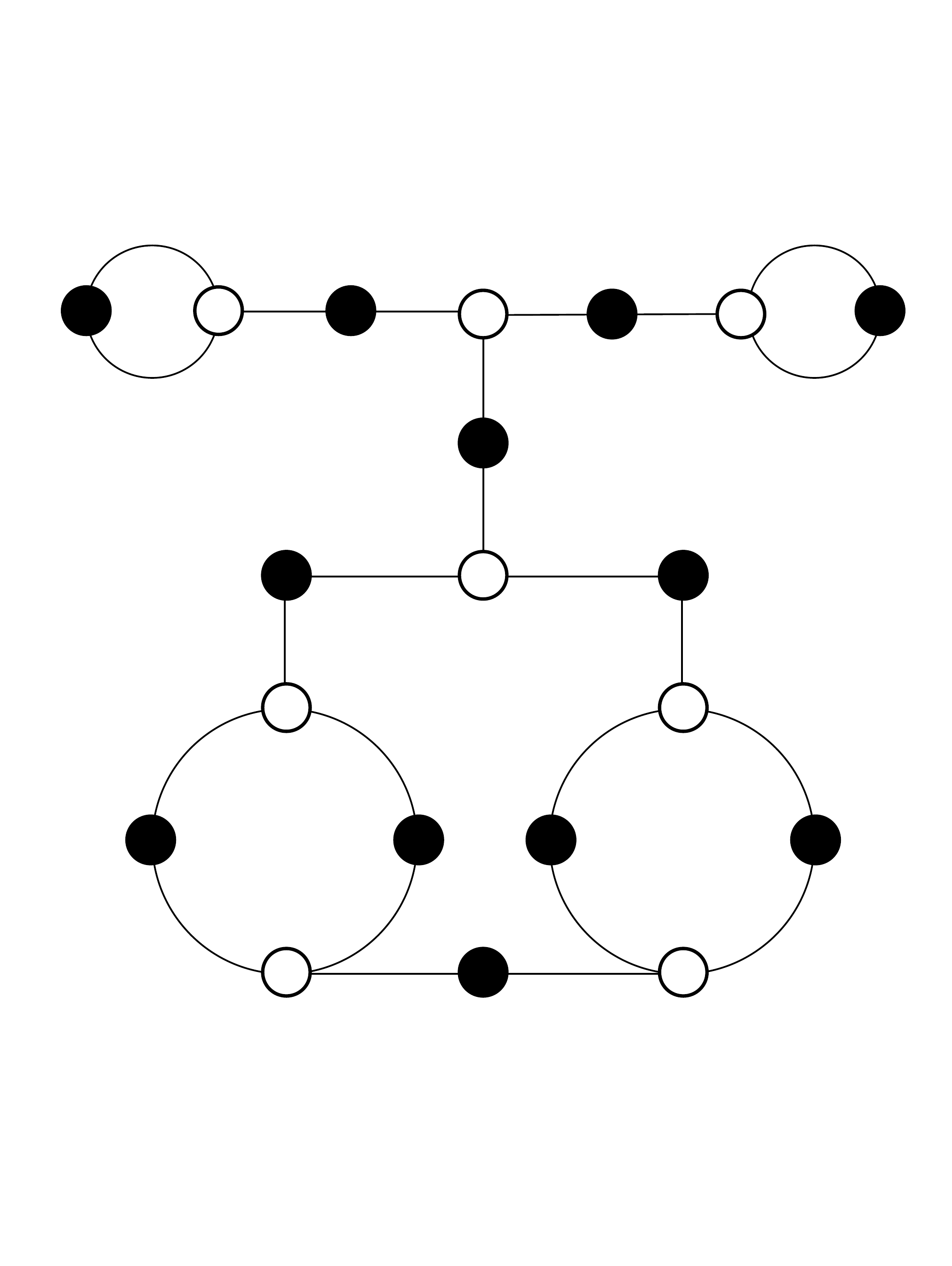}}
\par\end{center}{\scriptsize \par}

\begin{center}
{\scriptsize $13,5,2,2,1,1\;\left(\sqrt{65}\right)$}
\par\end{center}%
\end{minipage}{\scriptsize }%
\begin{minipage}[t]{0.33\textwidth}%
\begin{center}
{\scriptsize \includegraphics[scale=0.15]{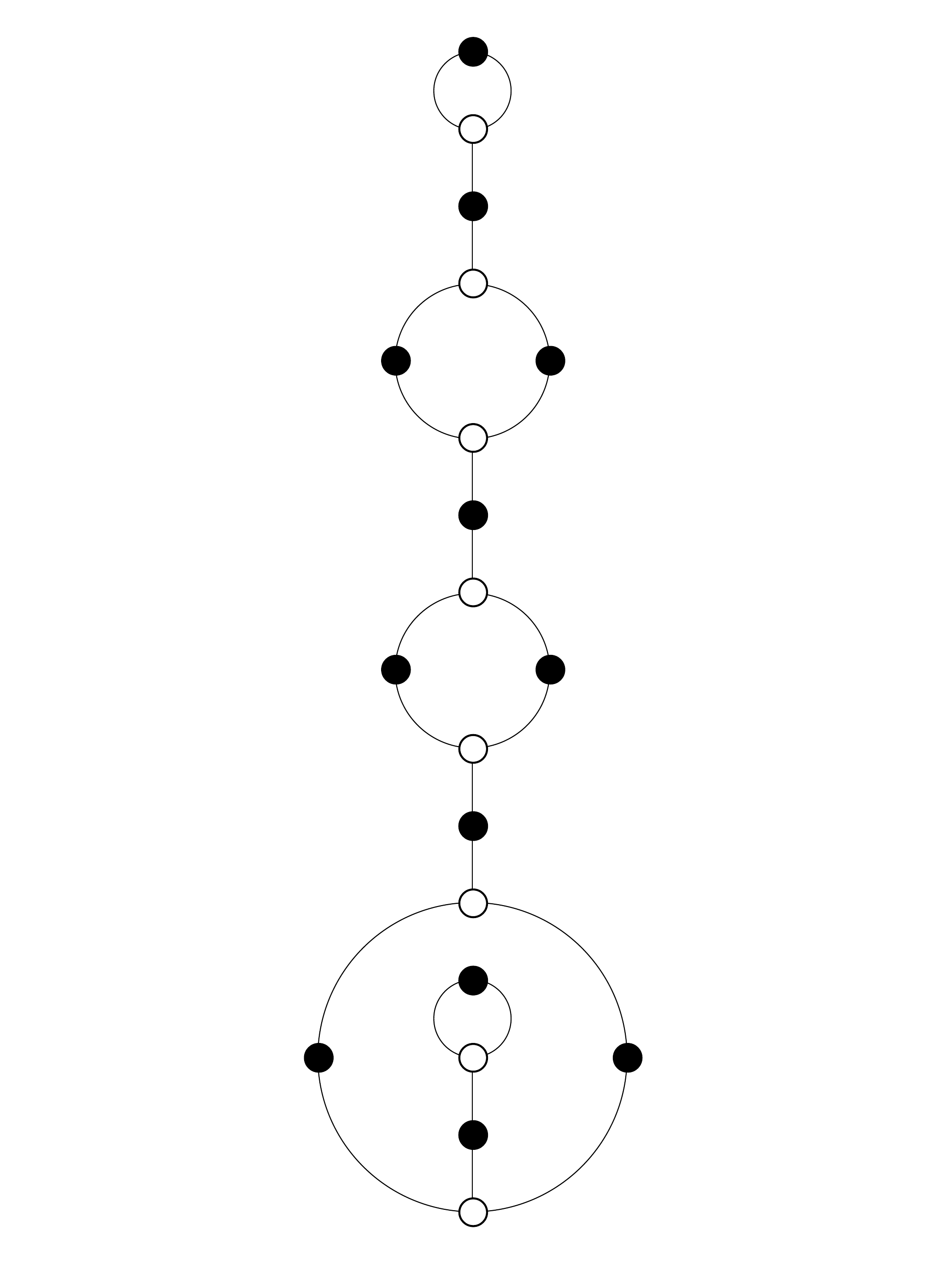}}
\par\end{center}{\scriptsize \par}

\begin{center}
{\scriptsize $13,5,2,2,1,1\;\left(\sqrt{65}\right)$}
\par\end{center}%
\end{minipage}{\scriptsize }%
\begin{minipage}[t]{0.33\textwidth}%
\begin{center}
{\scriptsize \includegraphics[scale=0.15]{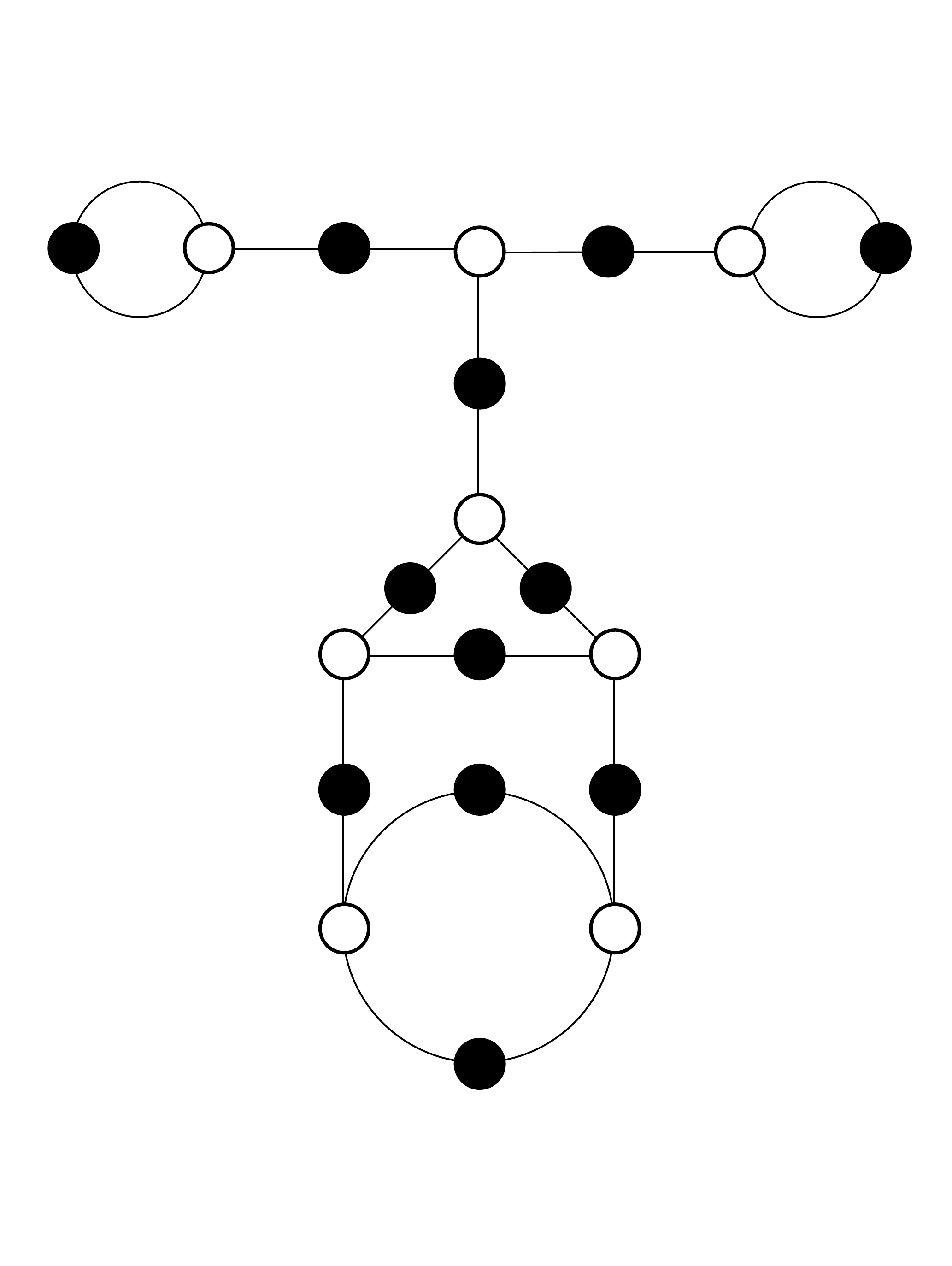}}
\par\end{center}{\scriptsize \par}

\begin{center}
{\scriptsize $13,4,3,2,1,1\;\left(\mathbb{Q}\right)$}
\par\end{center}%
\end{minipage}
\par\end{center}{\scriptsize \par}

\begin{center}
{\scriptsize }%
\begin{minipage}[t]{0.33\textwidth}%
\begin{center}
{\scriptsize \includegraphics[scale=0.15]{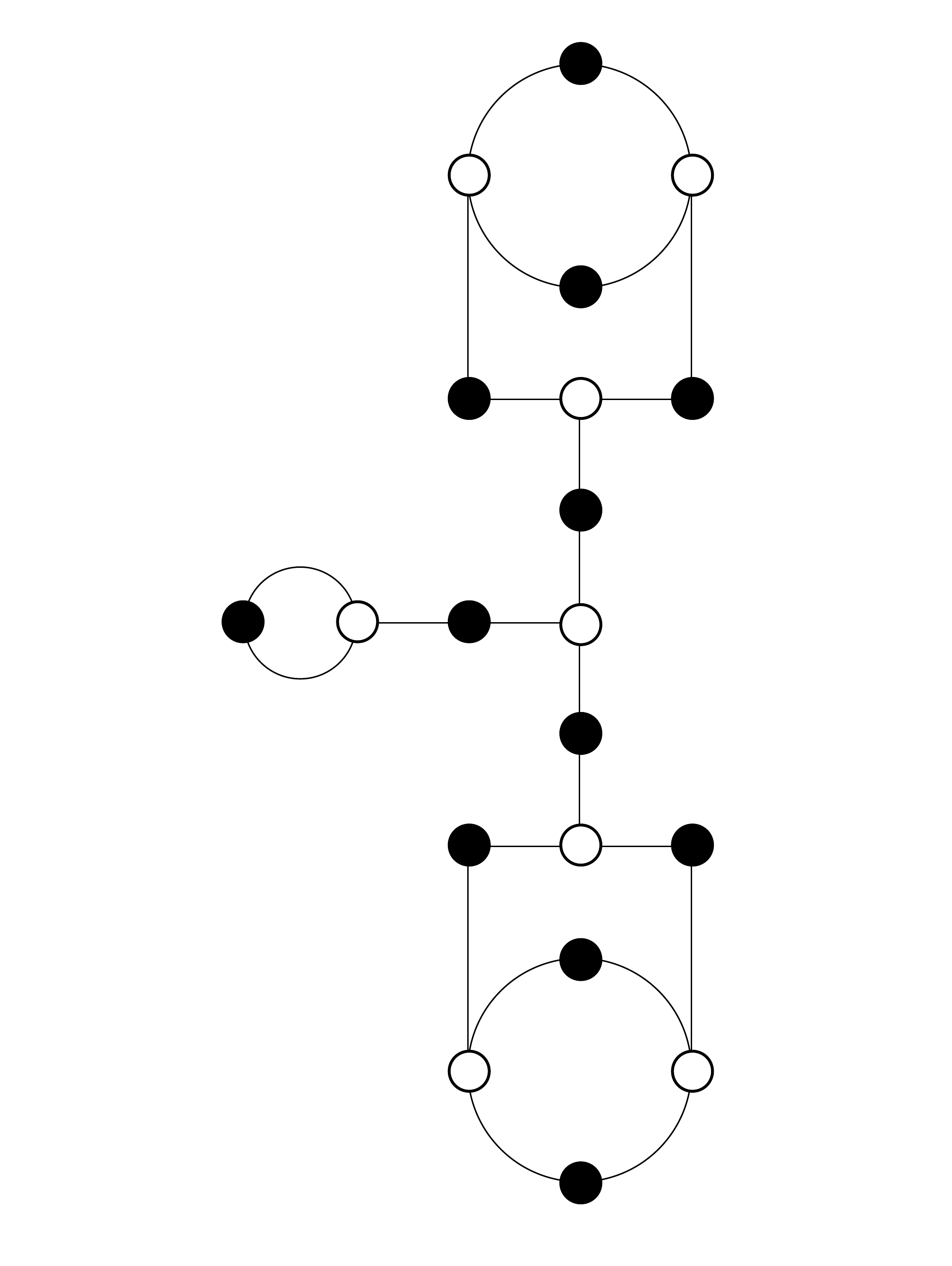}}
\par\end{center}{\scriptsize \par}

\begin{center}
{\scriptsize $13,3,3,2,2,1\;\left(\mathbb{Q}\right)$}
\par\end{center}%
\end{minipage}{\scriptsize }%
\begin{minipage}[t]{0.33\textwidth}%
\begin{center}
{\scriptsize \includegraphics[scale=0.15]{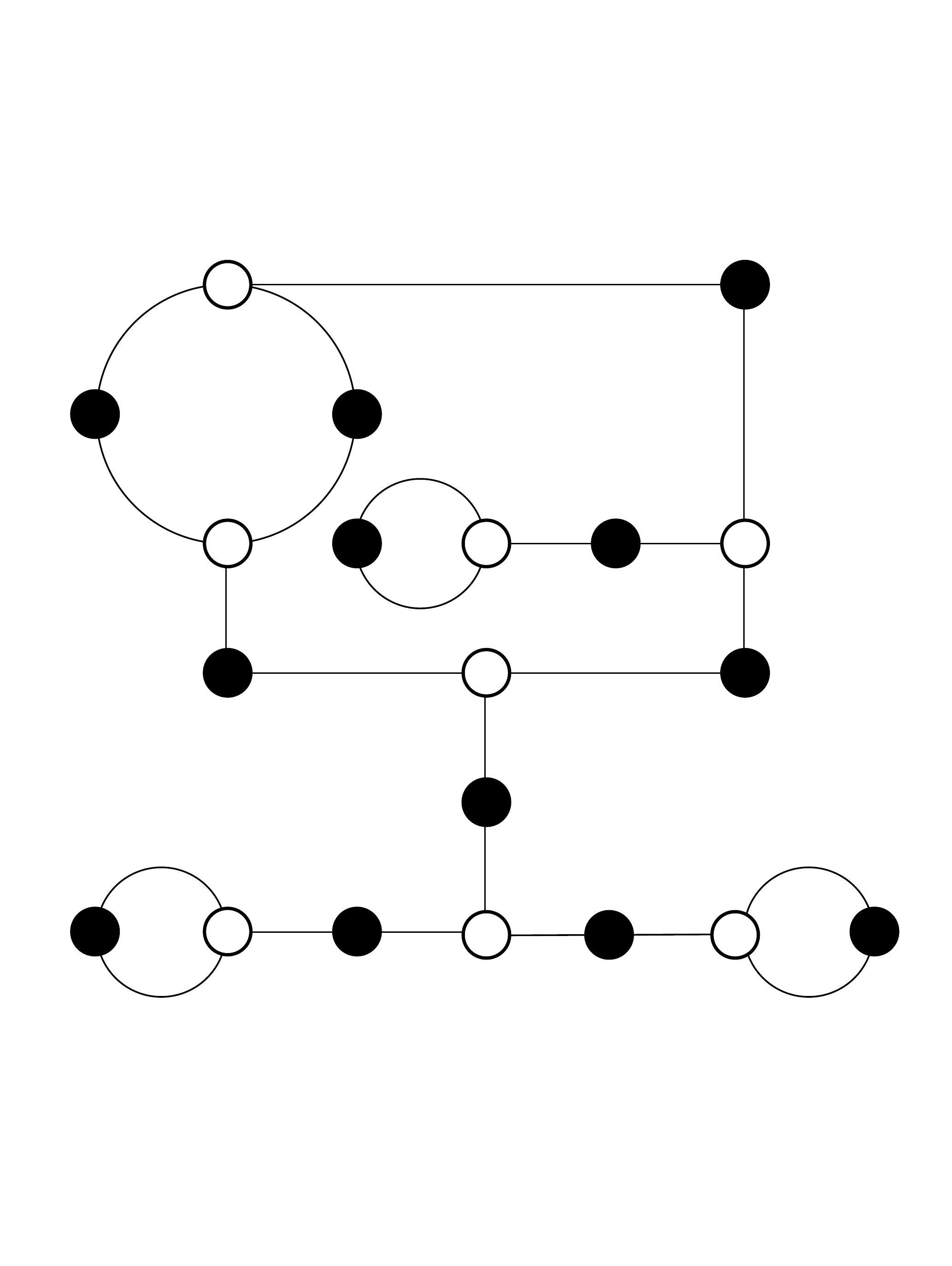}}
\par\end{center}{\scriptsize \par}

\begin{center}
{\scriptsize $12,7,2,1,1,1\;\left(\sqrt{-3}\right)$}
\par\end{center}%
\end{minipage}{\scriptsize }%
\begin{minipage}[t]{0.33\textwidth}%
\begin{center}
{\scriptsize \includegraphics[scale=0.15]{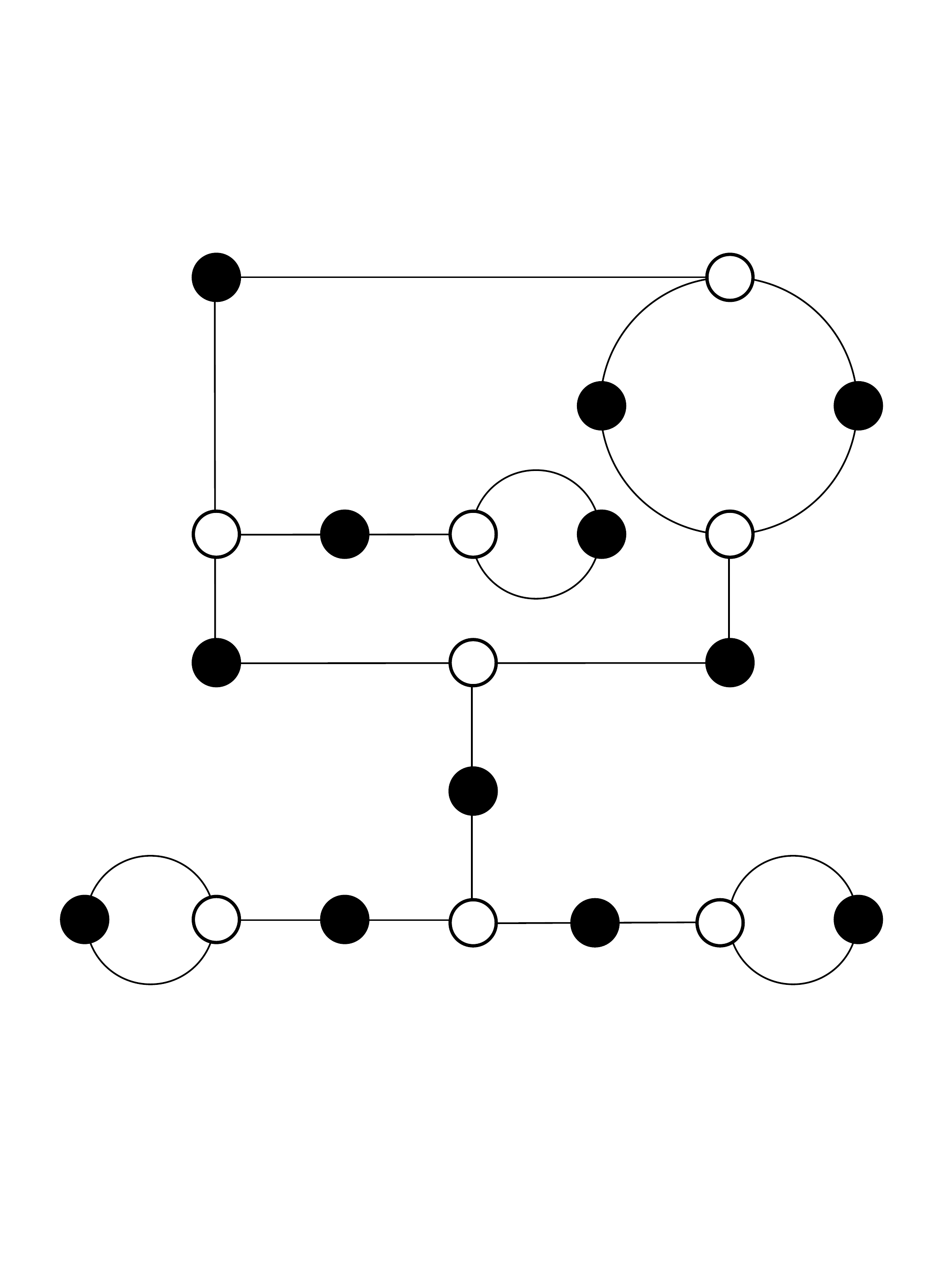}}
\par\end{center}{\scriptsize \par}

\begin{center}
{\scriptsize $12,7,2,1,1,1\;\left(\sqrt{-3}\right)$}
\par\end{center}%
\end{minipage}
\par\end{center}{\scriptsize \par}

\begin{center}
{\scriptsize }%
\begin{minipage}[t]{0.33\textwidth}%
\begin{center}
{\scriptsize \includegraphics[scale=0.15]{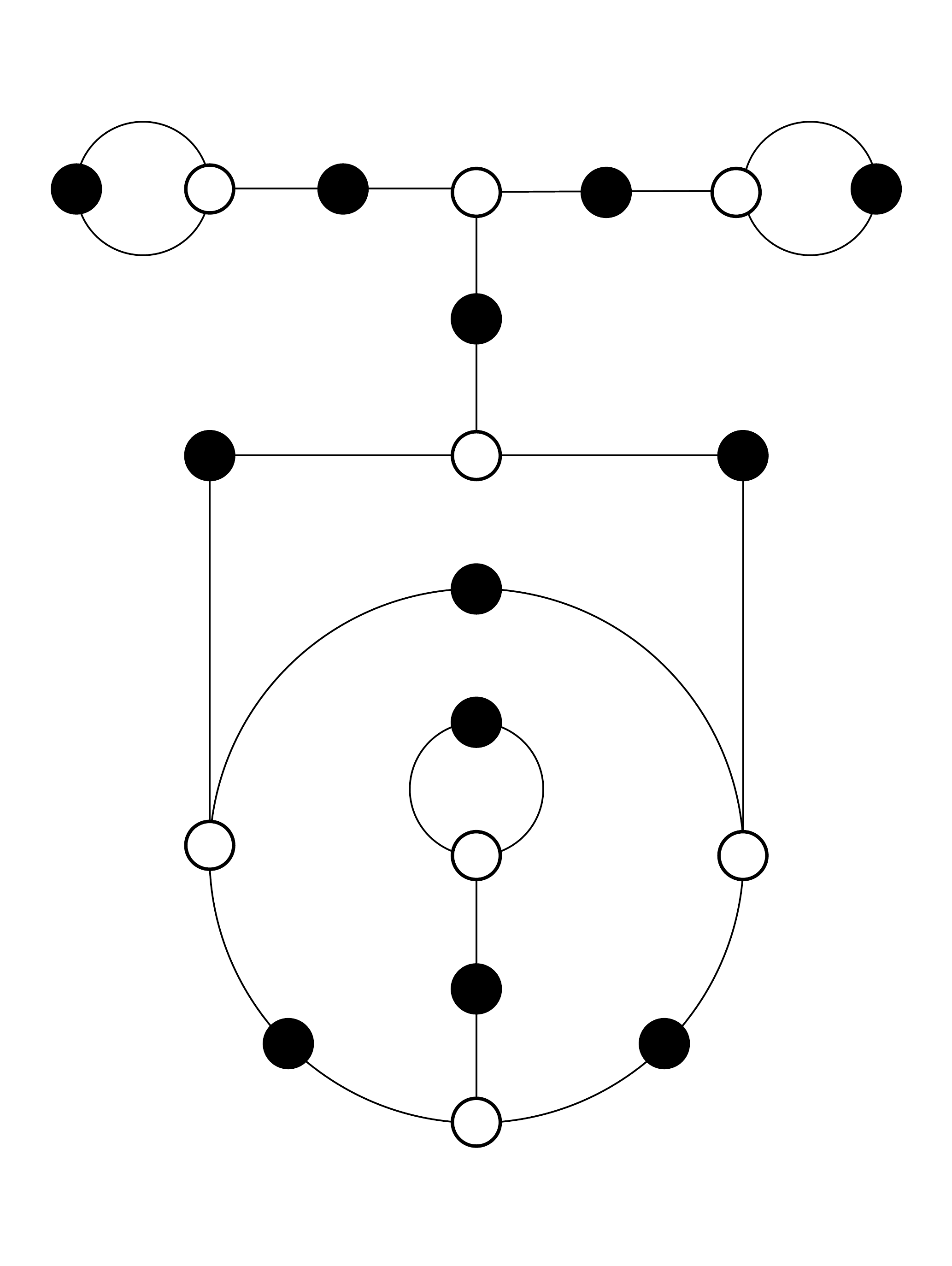}}
\par\end{center}{\scriptsize \par}

\begin{center}
{\scriptsize $12,6,3,1,1,1\;\left(\mathbb{Q}\right)$}
\par\end{center}%
\end{minipage}{\scriptsize }%
\begin{minipage}[t]{0.33\textwidth}%
\begin{center}
{\scriptsize \includegraphics[scale=0.15]{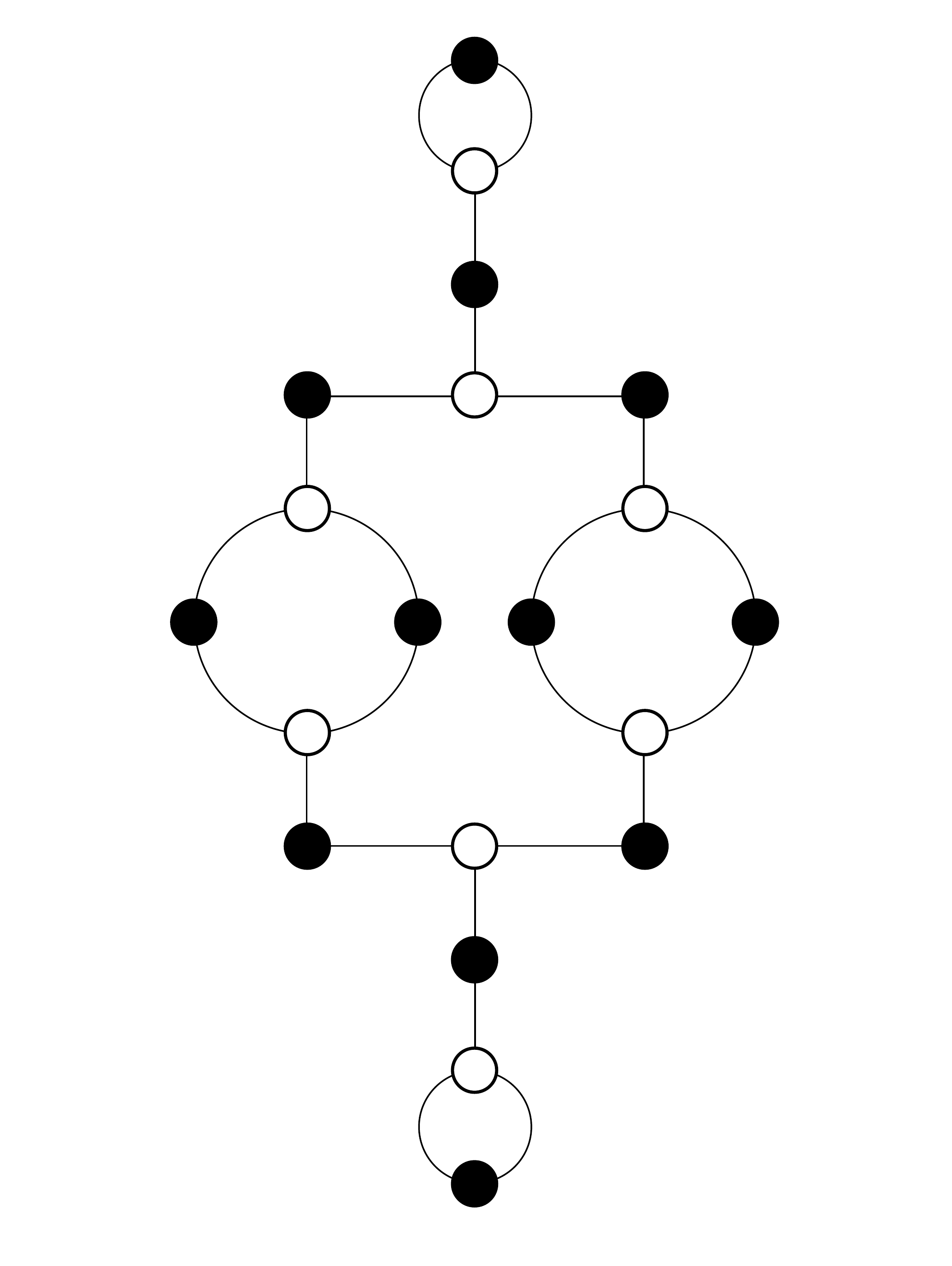}}
\par\end{center}{\scriptsize \par}

\begin{center}
{\scriptsize $12,6,2,2,1,1\;\left(\mathbb{Q}\right)$}
\par\end{center}%
\end{minipage}{\scriptsize }%
\begin{minipage}[t]{0.33\textwidth}%
\begin{center}
{\scriptsize \includegraphics[scale=0.15]{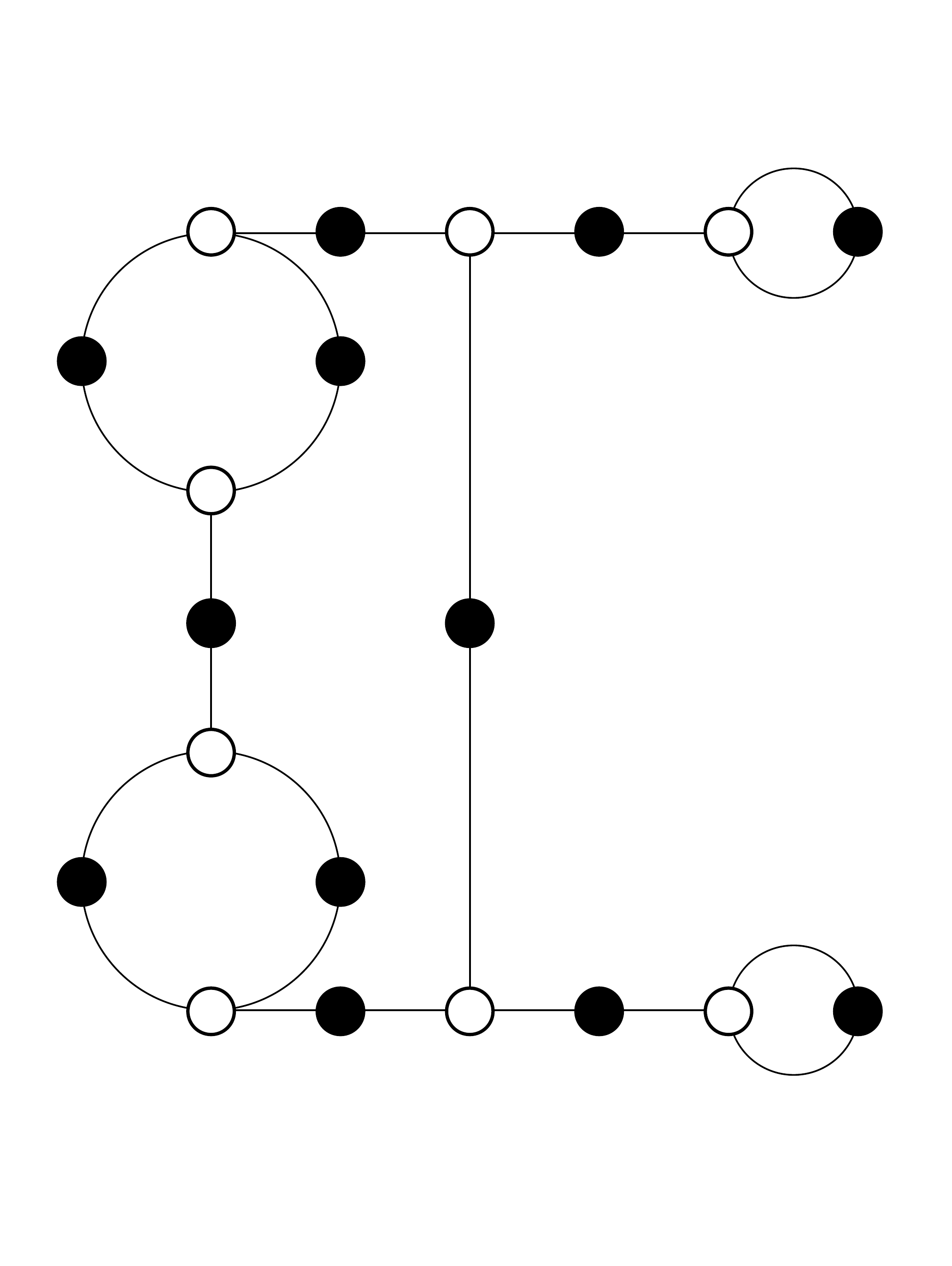}}
\par\end{center}{\scriptsize \par}

\begin{center}
{\scriptsize $12,6,2,2,1,1\;\left(\mathbb{Q}\right)$}
\par\end{center}%
\end{minipage}
\par\end{center}{\scriptsize \par}

\begin{center}
{\scriptsize }%
\begin{minipage}[t]{0.33\textwidth}%
\begin{center}
{\scriptsize \includegraphics[scale=0.15]{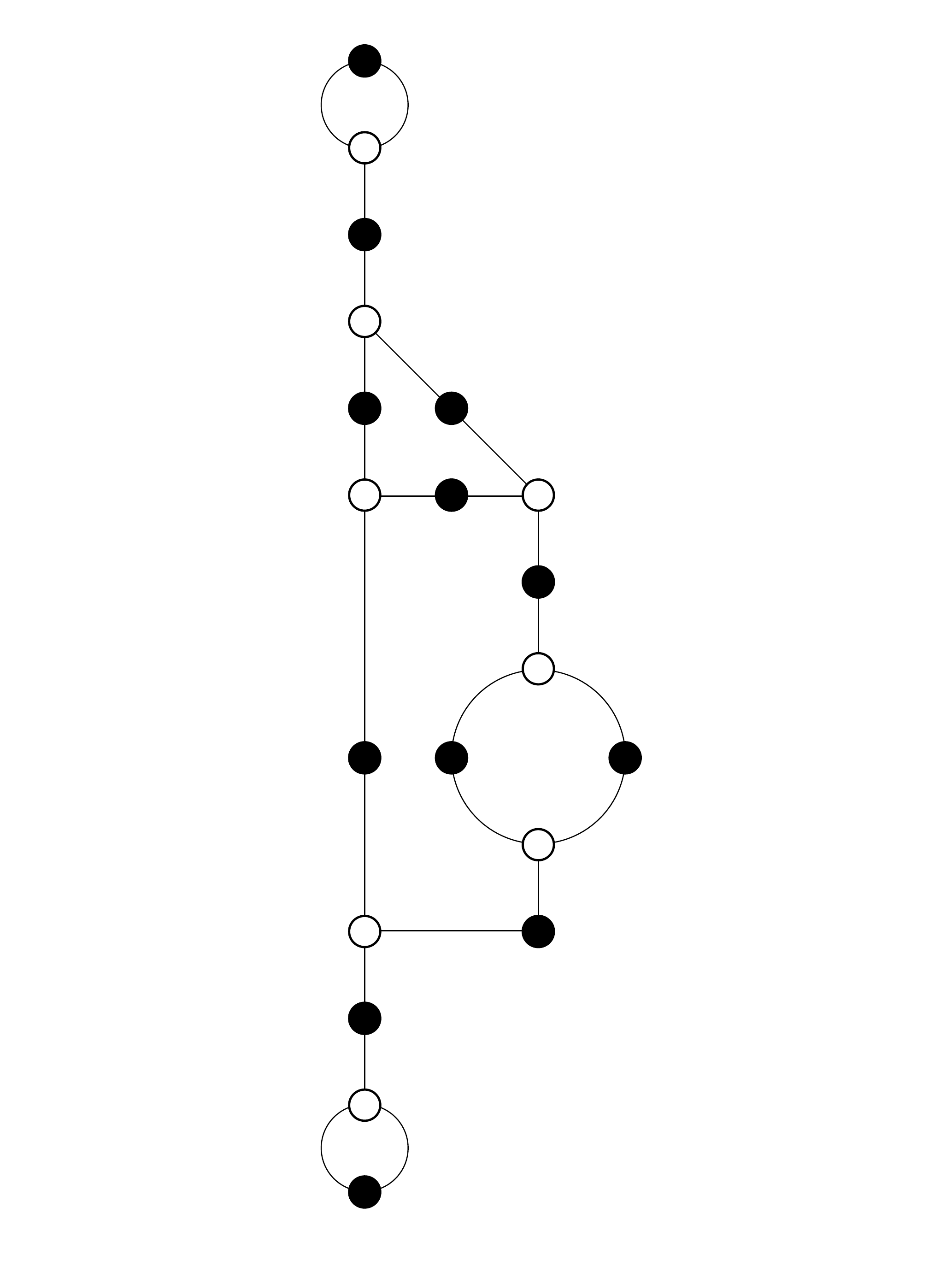}}
\par\end{center}{\scriptsize \par}

\begin{center}
{\scriptsize $12,5,3,2,1,1\;\left(\mathrm{quartic}\right)$}
\par\end{center}%
\end{minipage}{\scriptsize }%
\begin{minipage}[t]{0.33\textwidth}%
\begin{center}
{\scriptsize \includegraphics[scale=0.15]{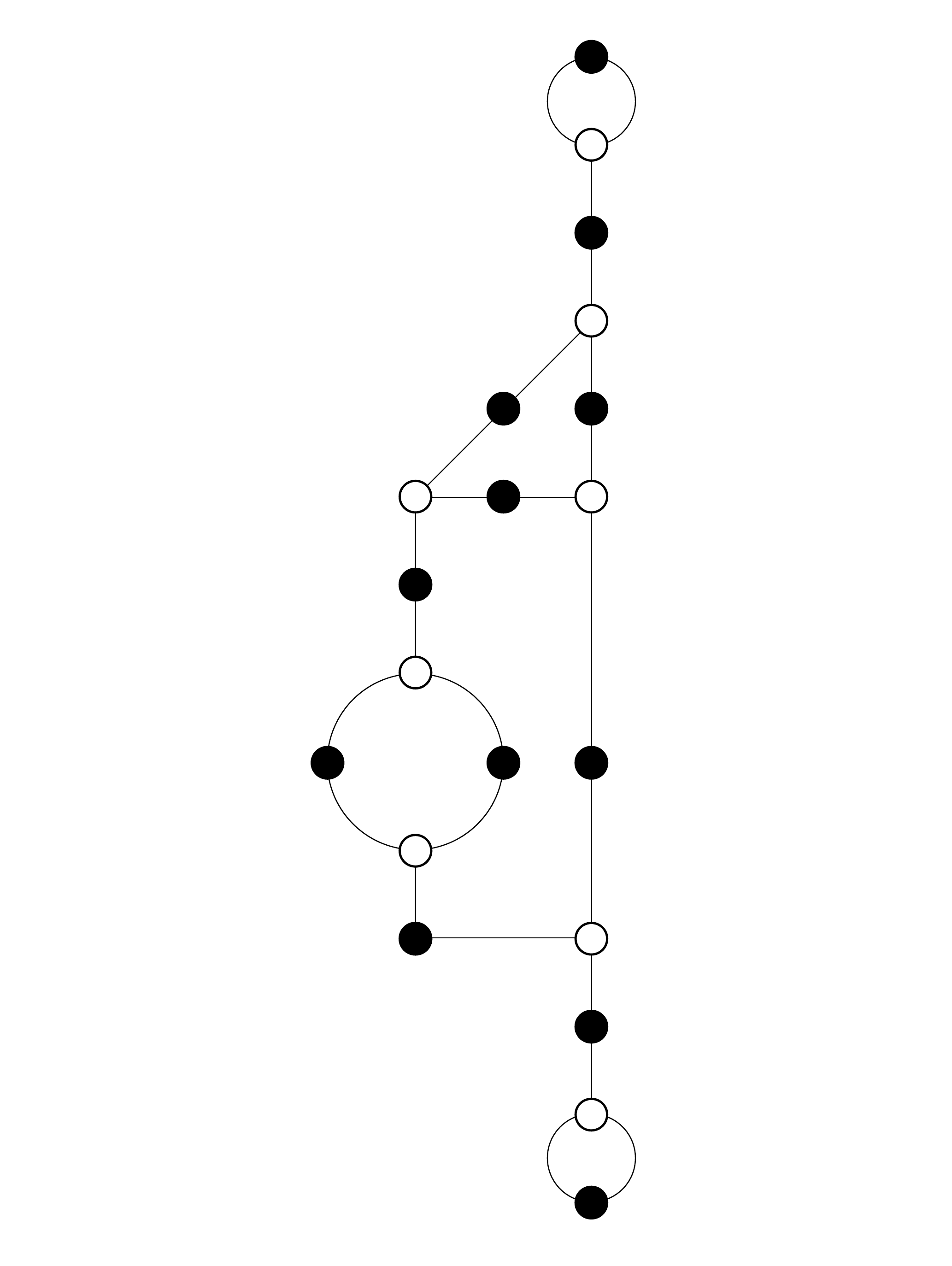}}
\par\end{center}{\scriptsize \par}

\begin{center}
{\scriptsize $12,5,3,2,1,1\;\left(\mathrm{quartic}\right)$}
\par\end{center}%
\end{minipage}{\scriptsize }%
\begin{minipage}[t]{0.33\textwidth}%
\begin{center}
{\scriptsize \includegraphics[scale=0.15]{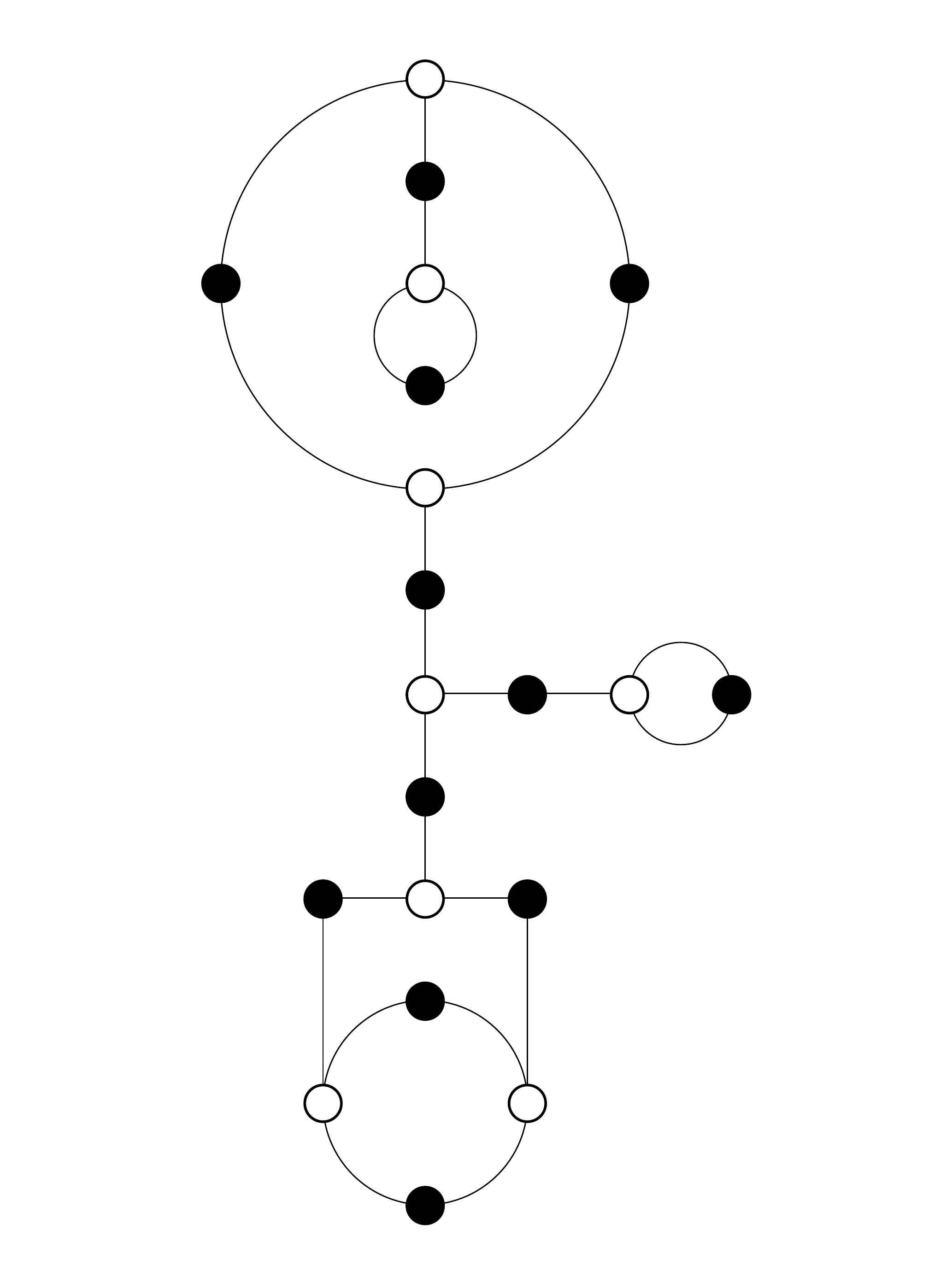}}
\par\end{center}{\scriptsize \par}

\begin{center}
{\scriptsize $12,5,3,2,1,1\;\left(\mathrm{quartic}\right)$}
\par\end{center}%
\end{minipage}
\par\end{center}{\scriptsize \par}

\begin{center}
{\scriptsize }%
\begin{minipage}[t]{0.33\textwidth}%
\begin{center}
{\scriptsize \includegraphics[scale=0.15]{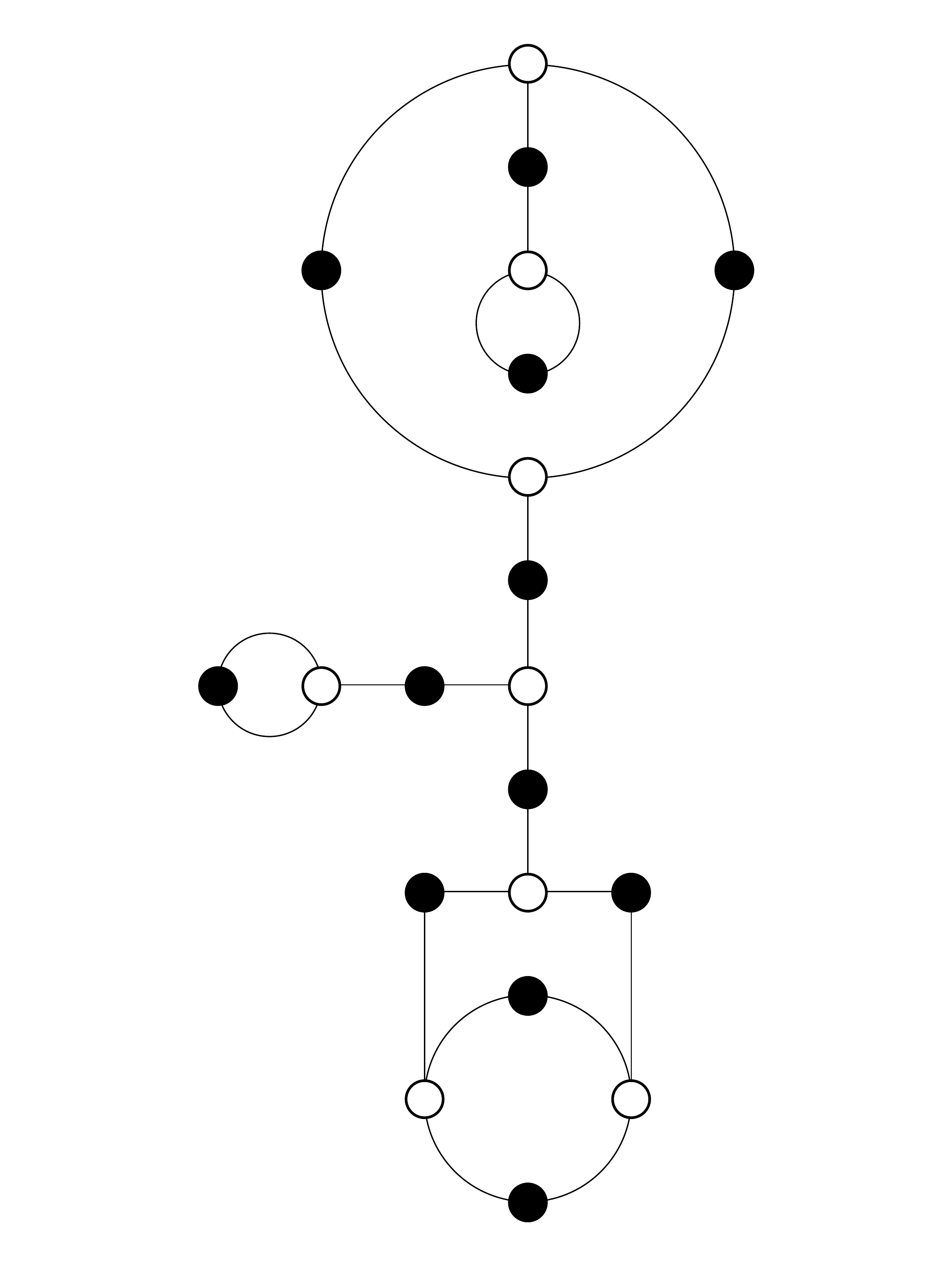}}
\par\end{center}{\scriptsize \par}

\begin{center}
{\scriptsize $12,5,3,2,1,1\;\left(\mathrm{quartic}\right)$}
\par\end{center}%
\end{minipage}{\scriptsize }%
\begin{minipage}[t]{0.33\textwidth}%
\begin{center}
{\scriptsize \includegraphics[scale=0.15]{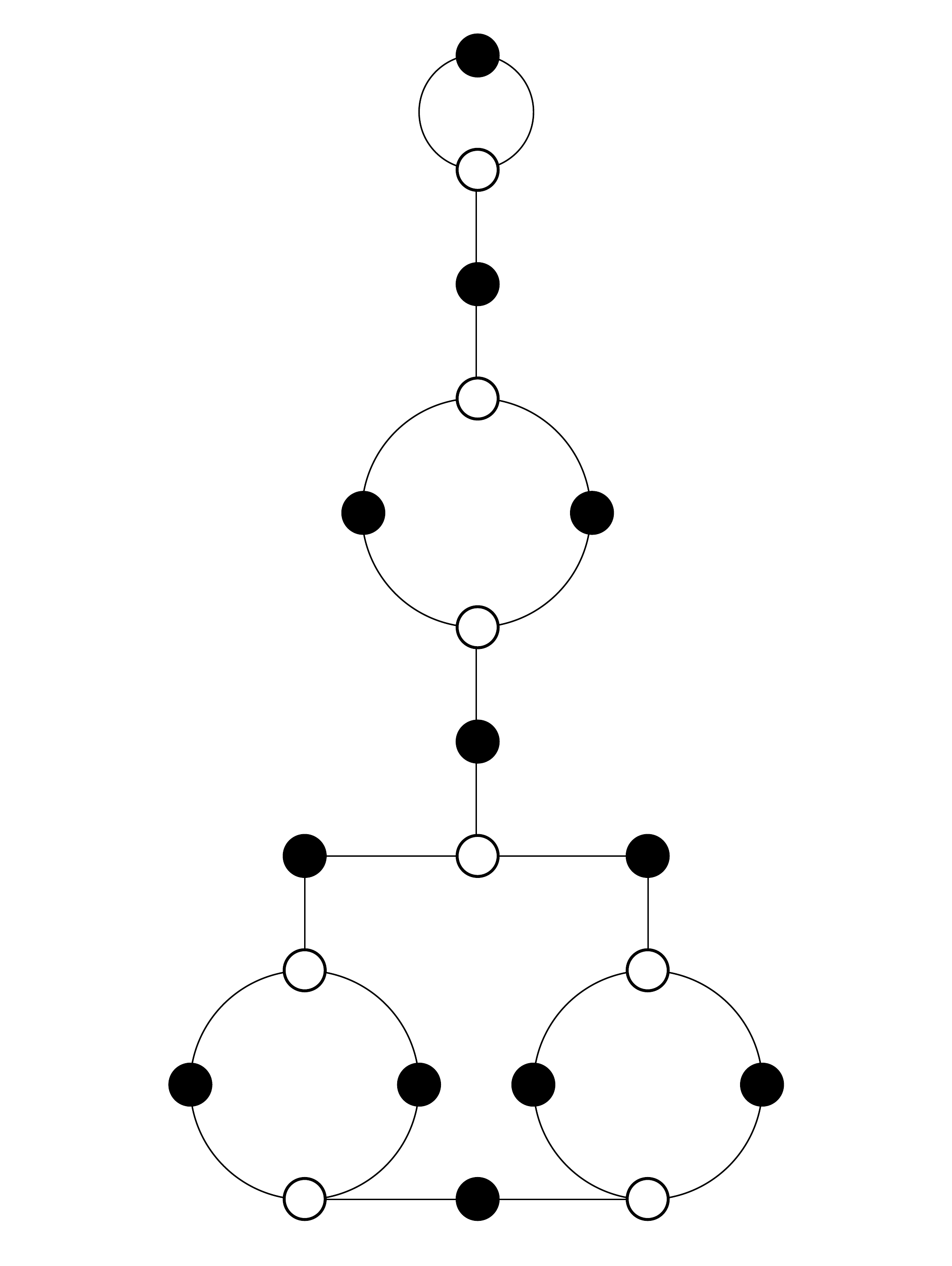}}
\par\end{center}{\scriptsize \par}

\begin{center}
{\scriptsize $12,5,2,2,2,1\;\left(\mathbb{Q}\right)$}
\par\end{center}%
\end{minipage}{\scriptsize }%
\begin{minipage}[t]{0.33\textwidth}%
\begin{center}
{\scriptsize \includegraphics[scale=0.15]{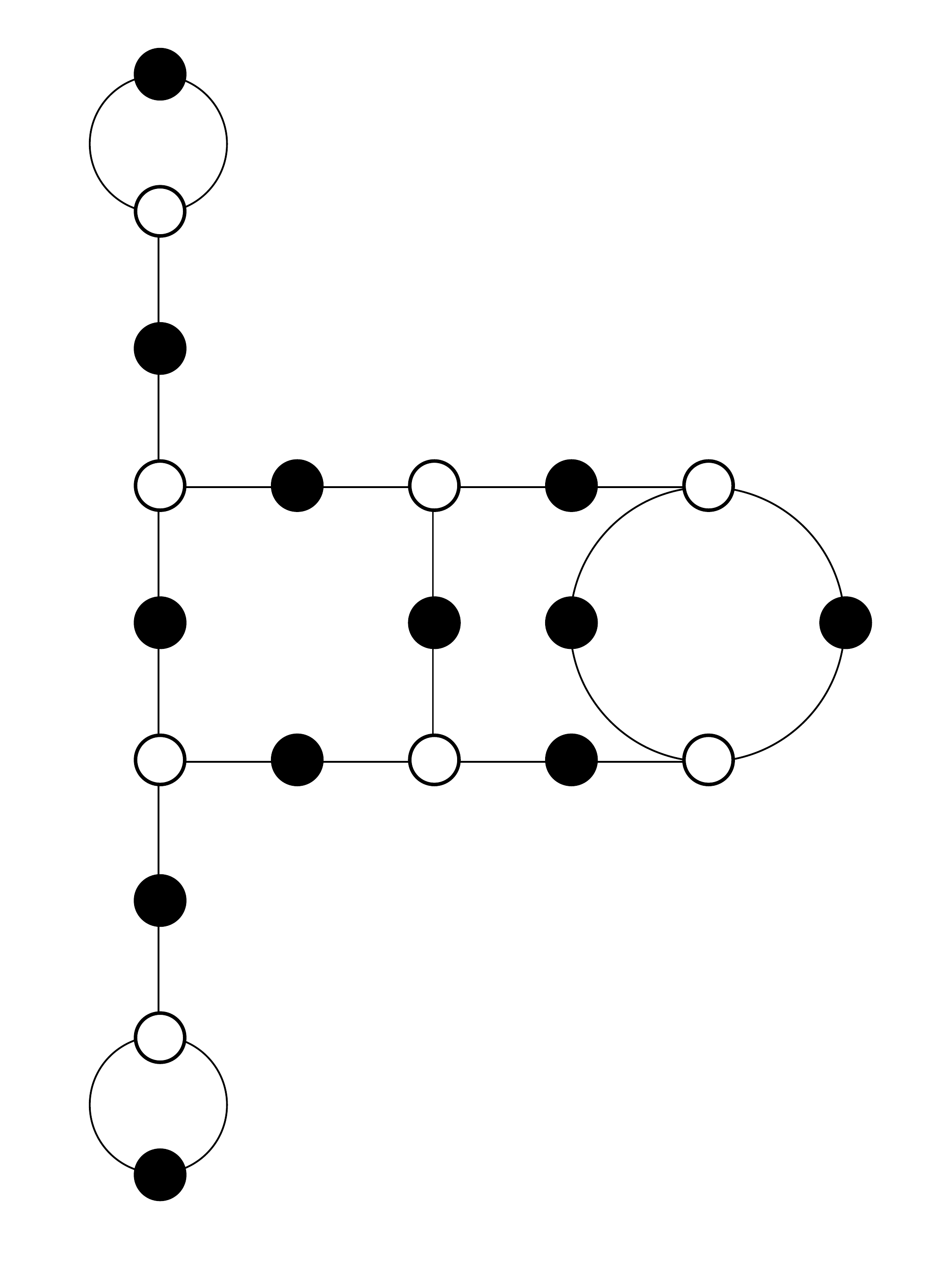}}
\par\end{center}{\scriptsize \par}

\begin{center}
{\scriptsize $12,4,4,2,1,1\;\left(\mathbb{Q}\right)$}
\par\end{center}%
\end{minipage}
\par\end{center}{\scriptsize \par}

\begin{center}
{\scriptsize }%
\begin{minipage}[t]{0.33\textwidth}%
\begin{center}
{\scriptsize \includegraphics[scale=0.15]{\string"PICT/12-4-3-3-1-1_A\string".pdf}}
\par\end{center}{\scriptsize \par}

\begin{center}
{\scriptsize $12,4,3,3,1,1\;\left(\mathbb{Q}\right)$}
\par\end{center}%
\end{minipage}{\scriptsize }%
\begin{minipage}[t]{0.33\textwidth}%
\begin{center}
{\scriptsize \includegraphics[scale=0.15]{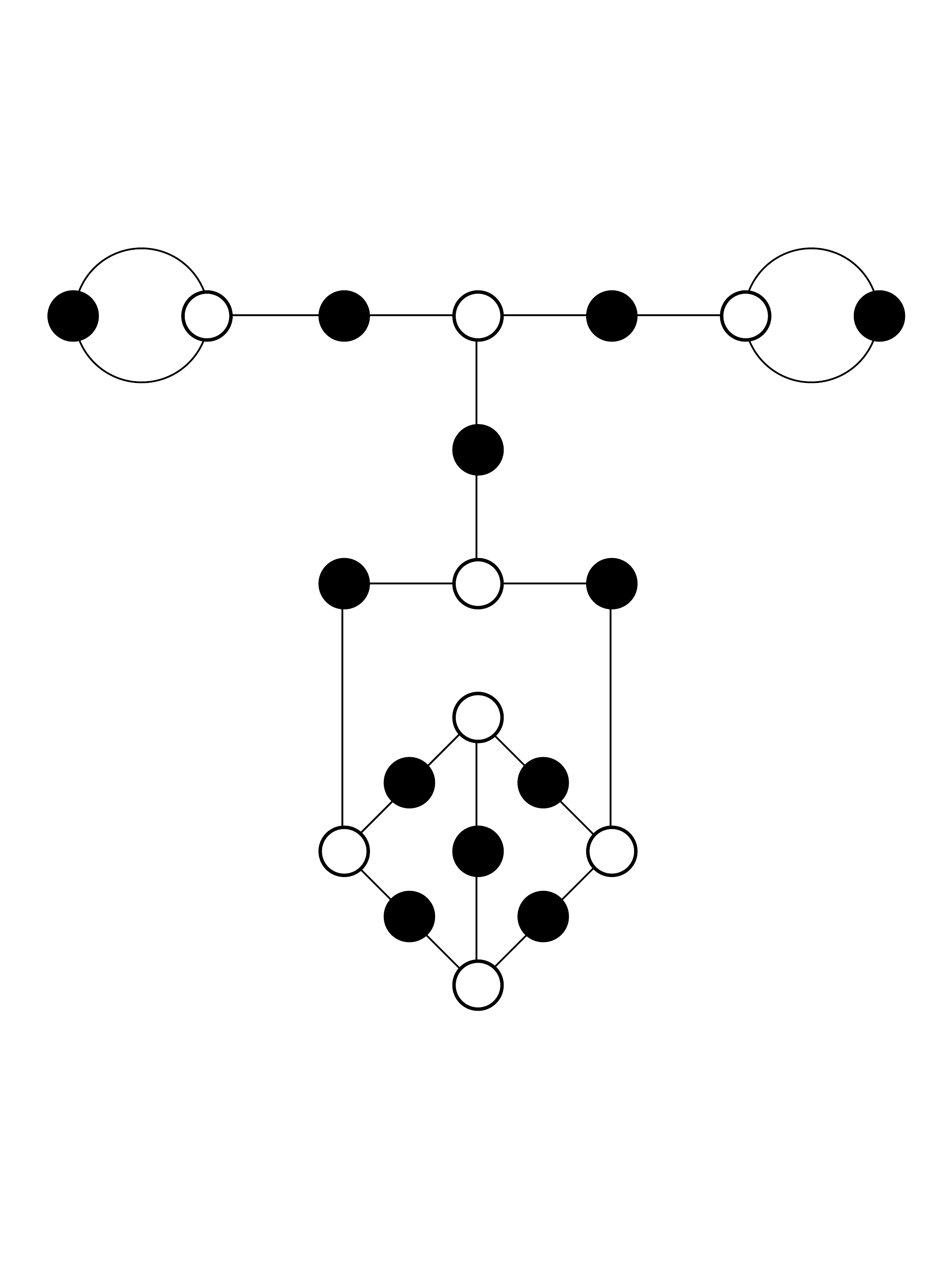}}
\par\end{center}{\scriptsize \par}

\begin{center}
{\scriptsize $12,4,3,3,1,1\;\left(\mathbb{Q}\right)$}
\par\end{center}%
\end{minipage}{\scriptsize }%
\begin{minipage}[t]{0.33\textwidth}%
\begin{center}
{\scriptsize \includegraphics[scale=0.15]{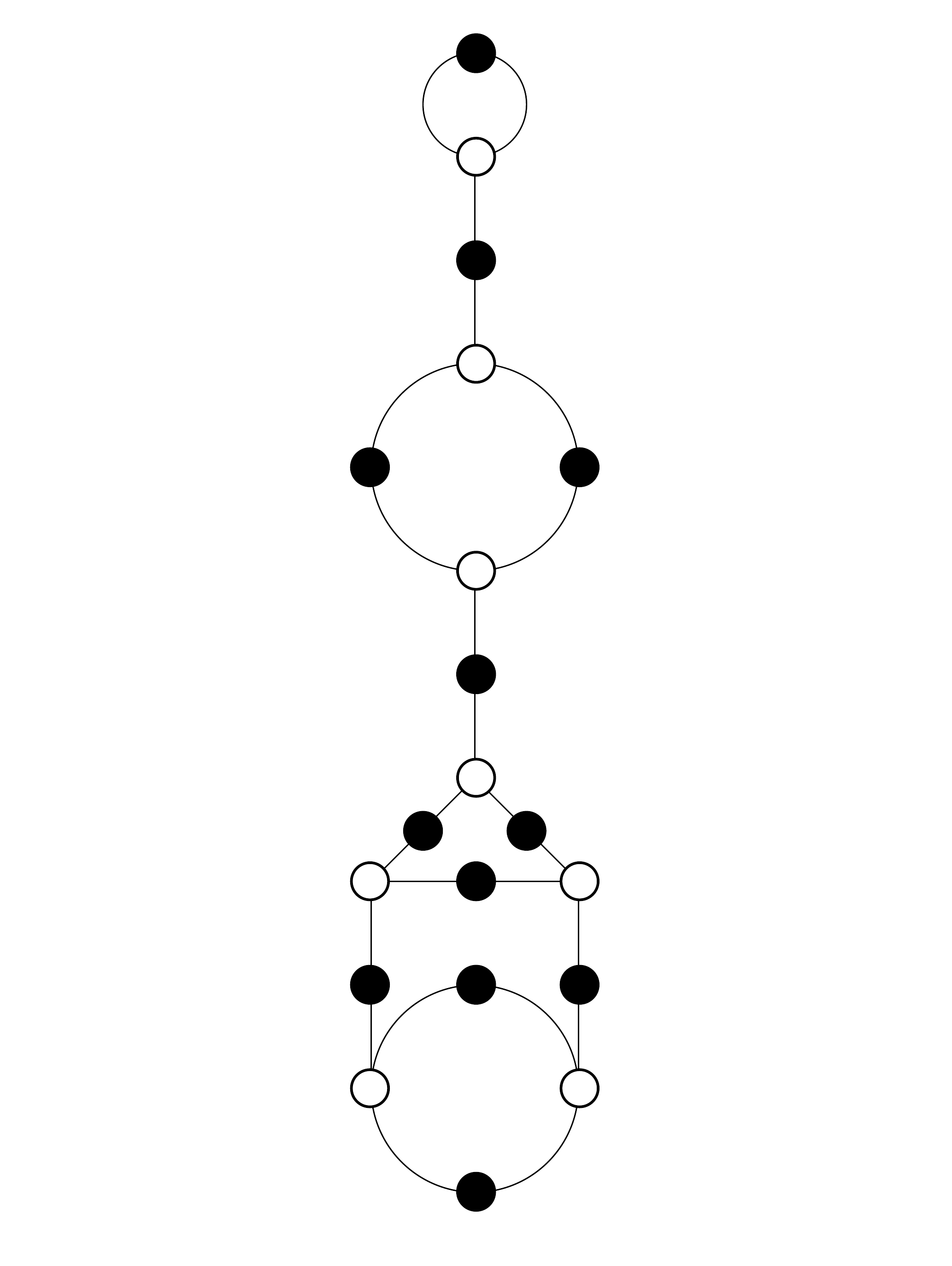}}
\par\end{center}{\scriptsize \par}

\begin{center}
{\scriptsize $12,4,3,2,2,1\;\left(\mathbb{Q}\right)$}
\par\end{center}%
\end{minipage}
\par\end{center}{\scriptsize \par}

\begin{center}
{\scriptsize }%
\begin{minipage}[t]{0.33\textwidth}%
\begin{center}
{\scriptsize \includegraphics[scale=0.15]{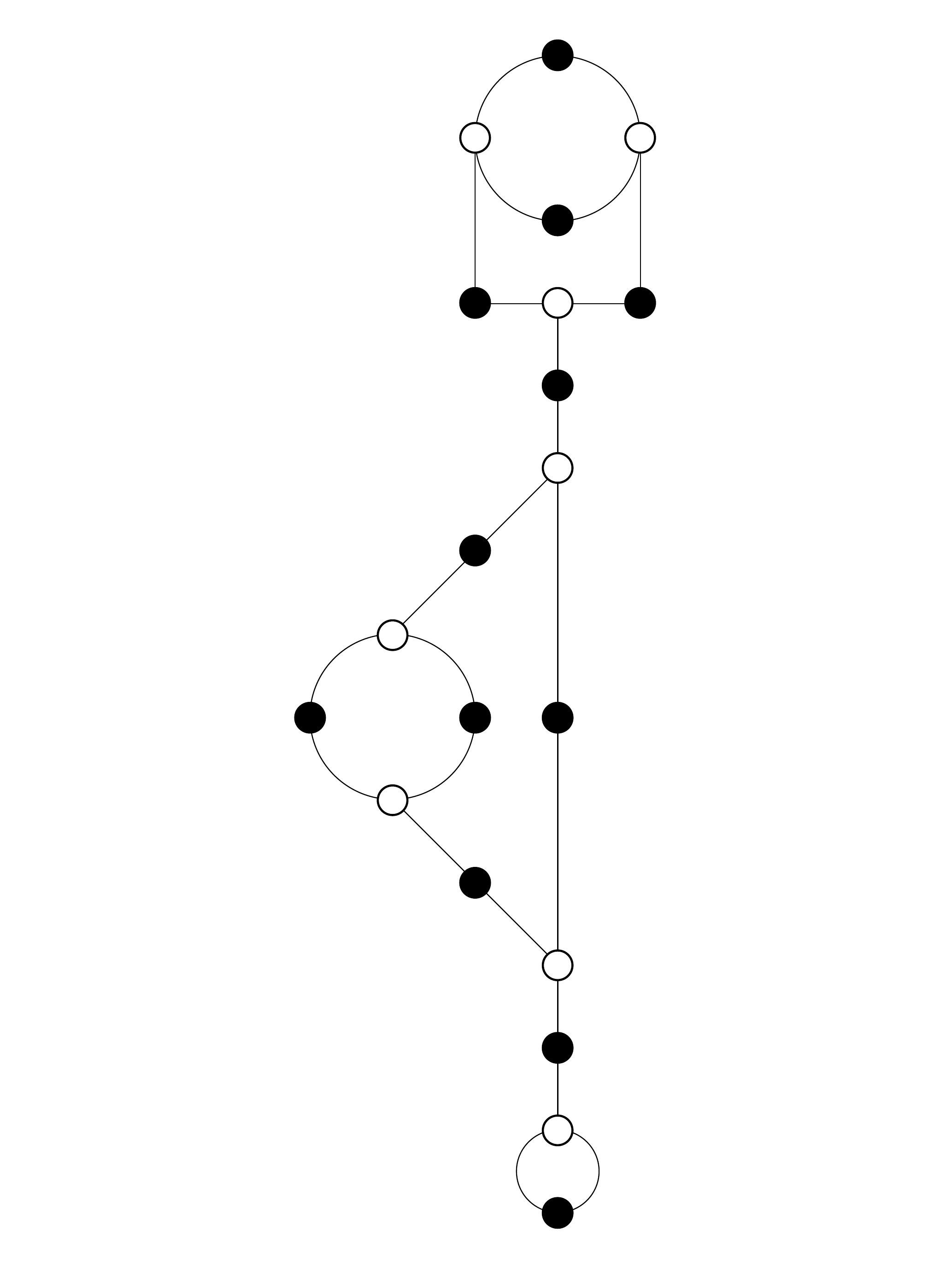}}
\par\end{center}{\scriptsize \par}

\begin{center}
{\scriptsize $12,4,3,2,2,1\;\left(\sqrt{-3}\right)$}
\par\end{center}%
\end{minipage}{\scriptsize }%
\begin{minipage}[t]{0.33\textwidth}%
\begin{center}
{\scriptsize \includegraphics[scale=0.15]{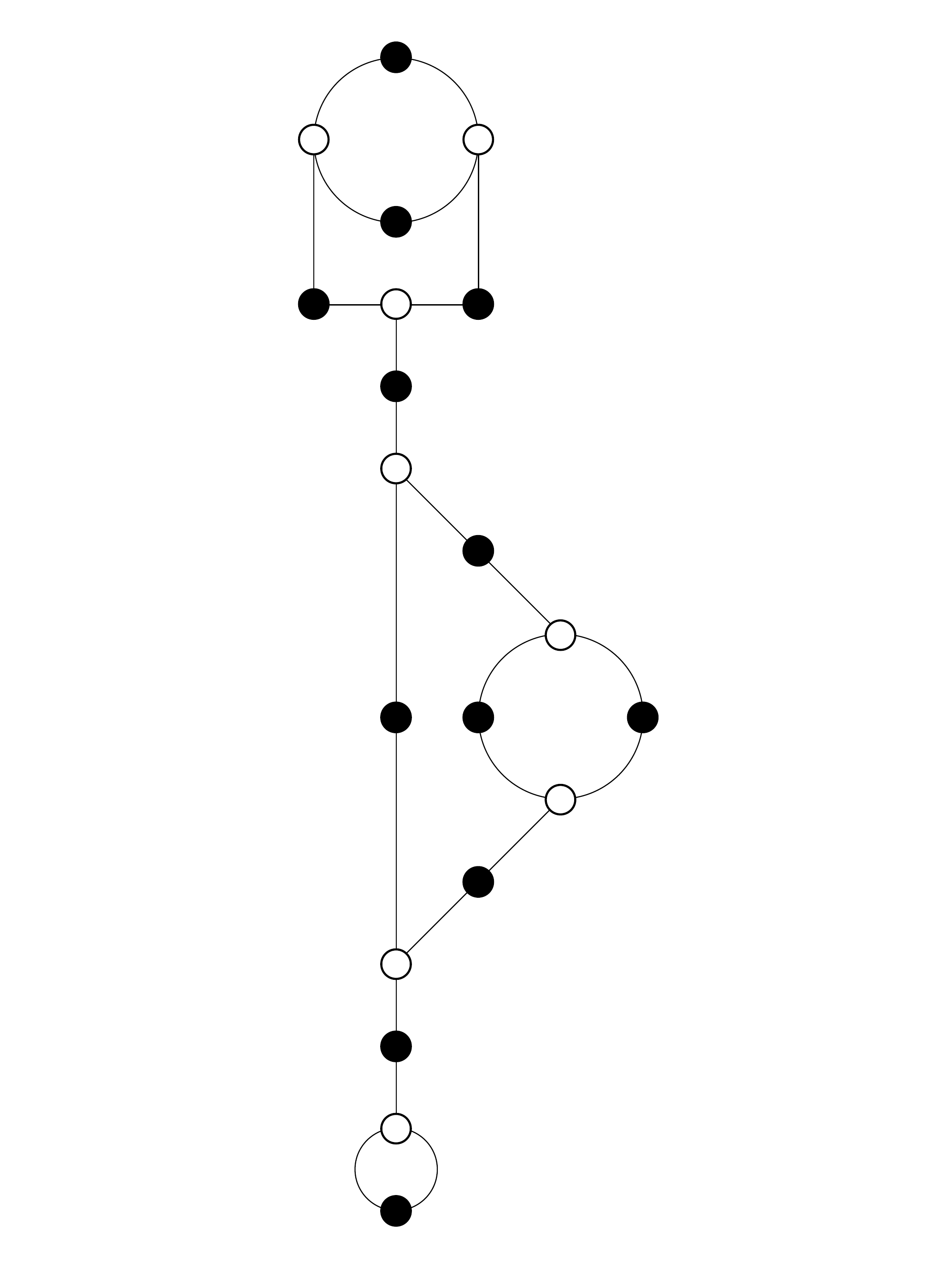}}
\par\end{center}{\scriptsize \par}

\begin{center}
{\scriptsize $12,4,3,2,2,1\;\left(\sqrt{-3}\right)$}
\par\end{center}%
\end{minipage}{\scriptsize }%
\begin{minipage}[t]{0.33\textwidth}%
\begin{center}
{\scriptsize \includegraphics[scale=0.15]{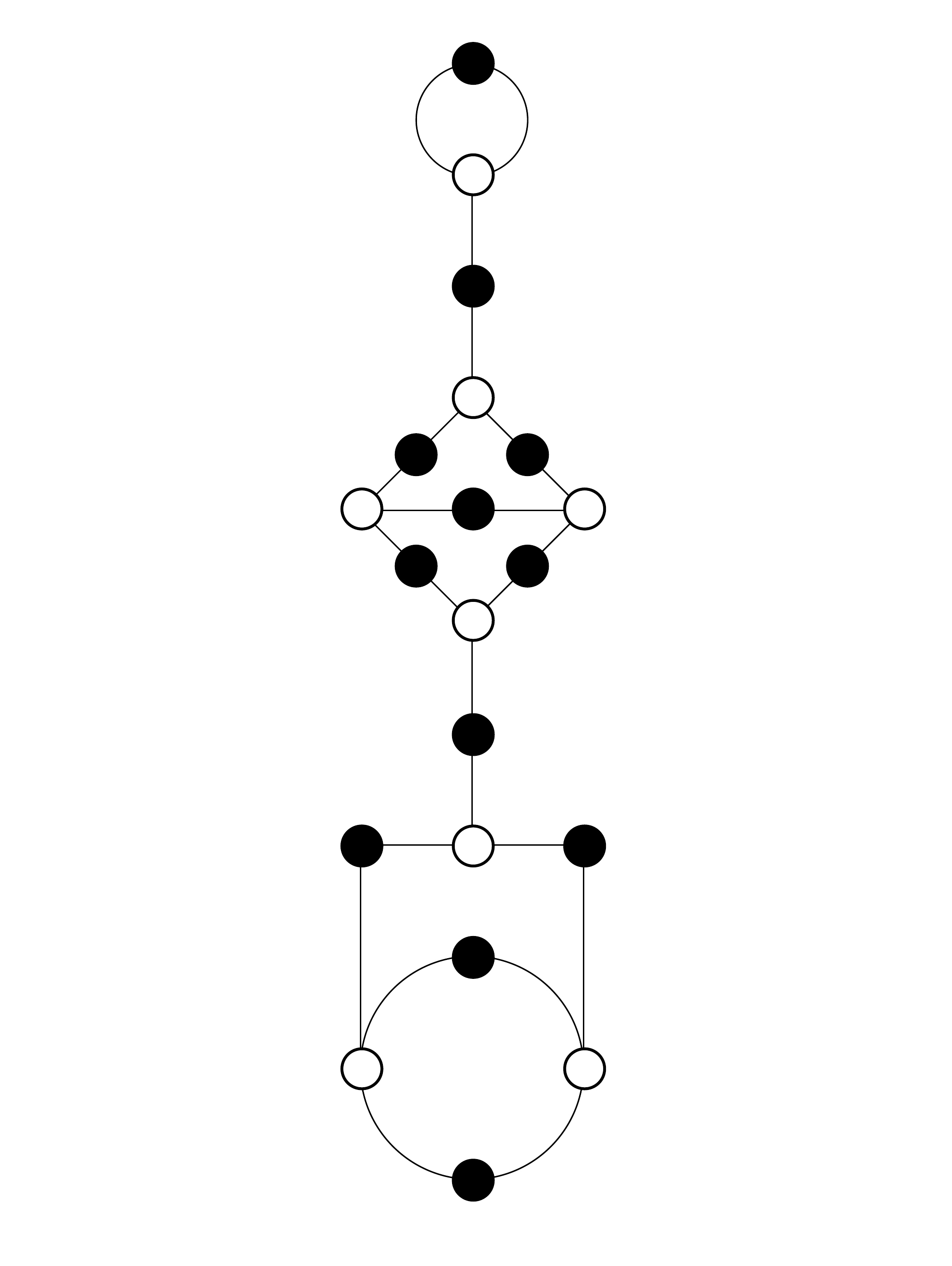}}
\par\end{center}{\scriptsize \par}

\begin{center}
{\scriptsize $12,3,3,3,2,1\;\left(\mathbb{Q}\right)$}
\par\end{center}%
\end{minipage}
\par\end{center}{\scriptsize \par}

\begin{center}
{\scriptsize }%
\begin{minipage}[t]{0.33\textwidth}%
\begin{center}
{\scriptsize \includegraphics[scale=0.15]{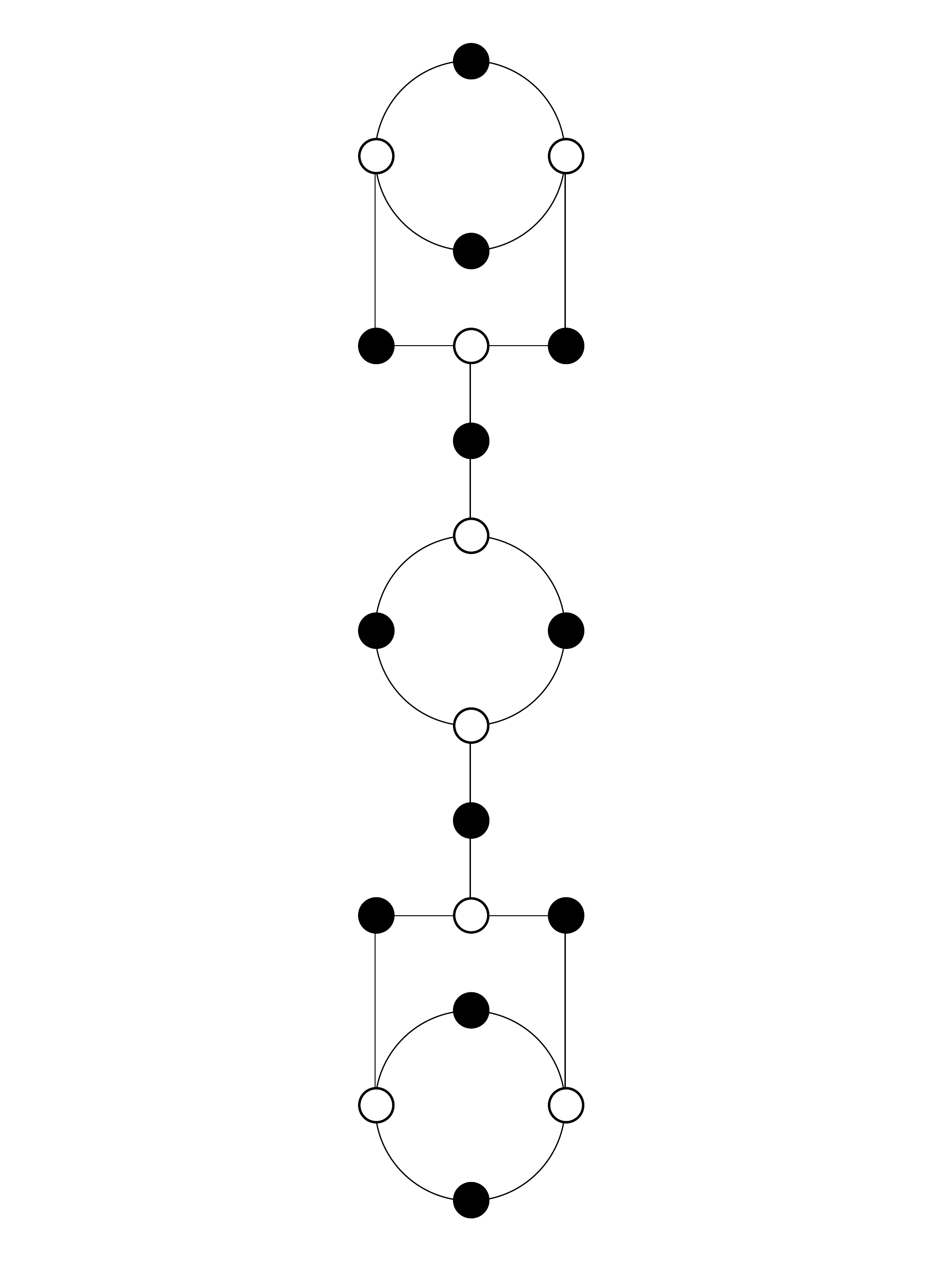}}
\par\end{center}{\scriptsize \par}

\begin{center}
{\scriptsize $12,3,3,2,2,2\;\left(\mathbb{Q}\right)$}
\par\end{center}%
\end{minipage}{\scriptsize }%
\begin{minipage}[t]{0.33\textwidth}%
\begin{center}
{\scriptsize \includegraphics[scale=0.15]{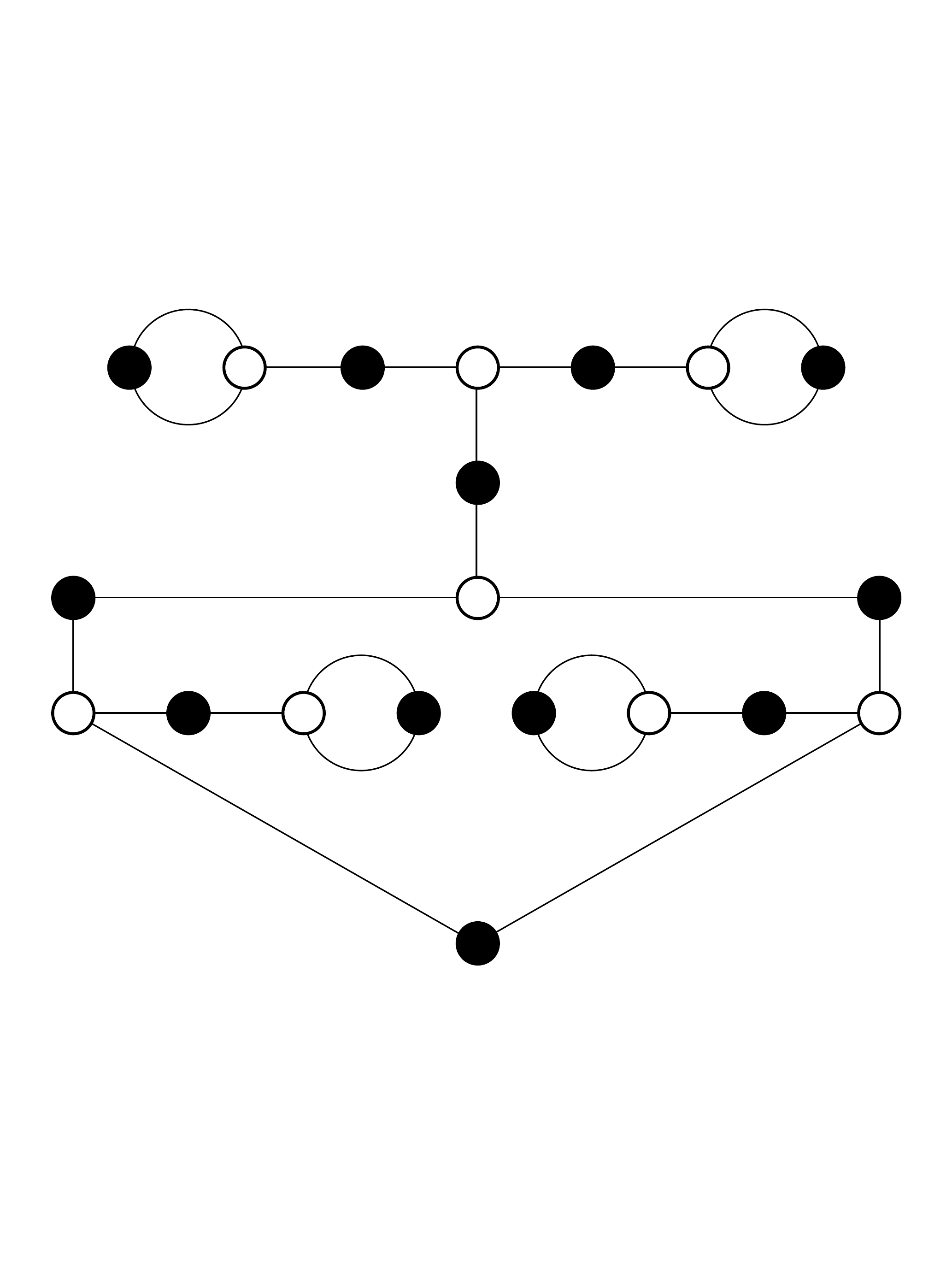}}
\par\end{center}{\scriptsize \par}

\begin{center}
{\scriptsize $11,9,1,1,1,1\;\left(\mathbb{Q}\right)$}
\par\end{center}%
\end{minipage}{\scriptsize }%
\begin{minipage}[t]{0.33\textwidth}%
\begin{center}
{\scriptsize \includegraphics[scale=0.15]{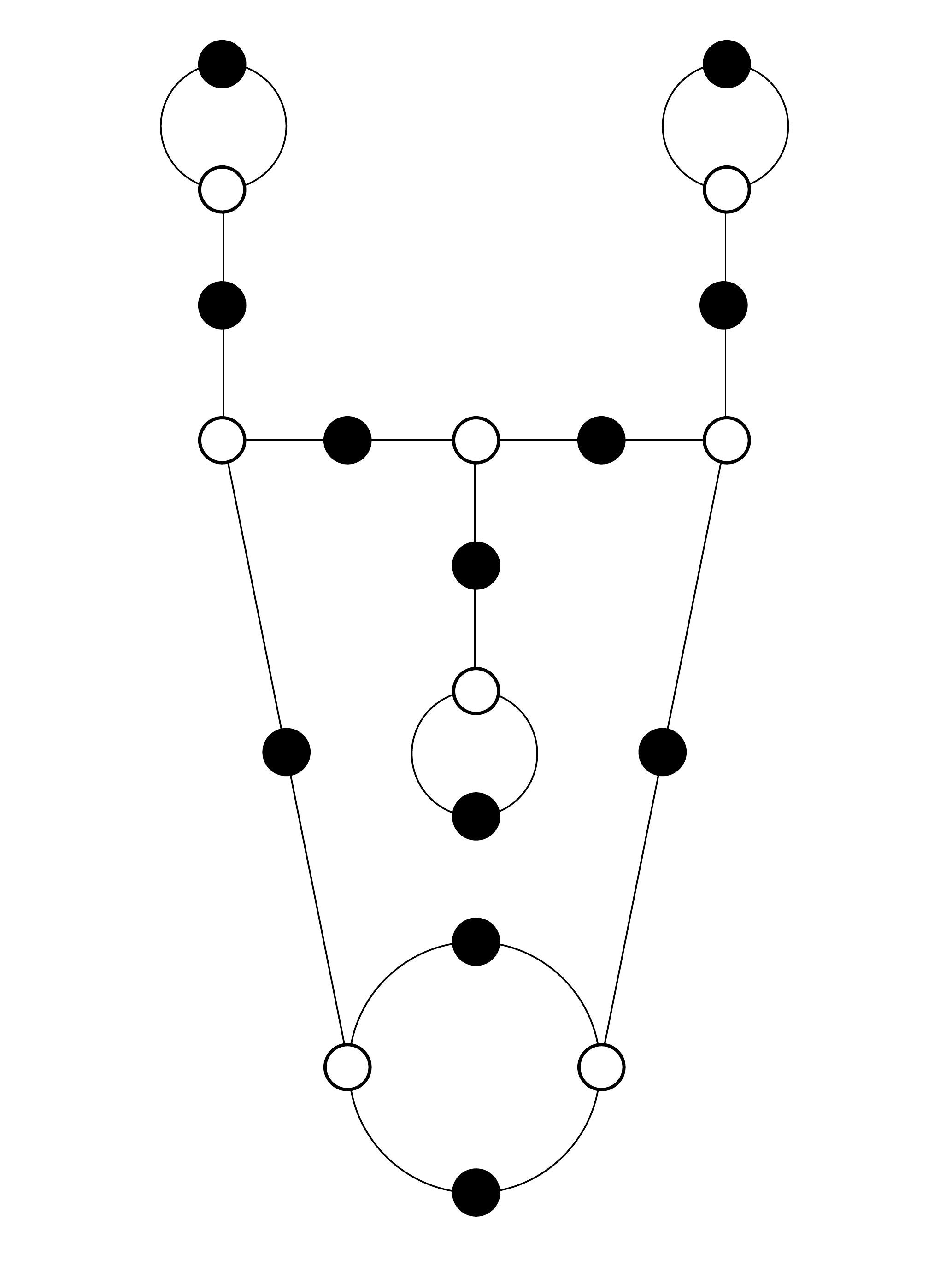}}
\par\end{center}{\scriptsize \par}

\begin{center}
{\scriptsize $11,8,2,1,1,1\;\left(\mathrm{cubic}\right)$}
\par\end{center}%
\end{minipage}
\par\end{center}{\scriptsize \par}

\begin{center}
{\scriptsize }%
\begin{minipage}[t]{0.33\textwidth}%
\begin{center}
{\scriptsize \includegraphics[scale=0.15]{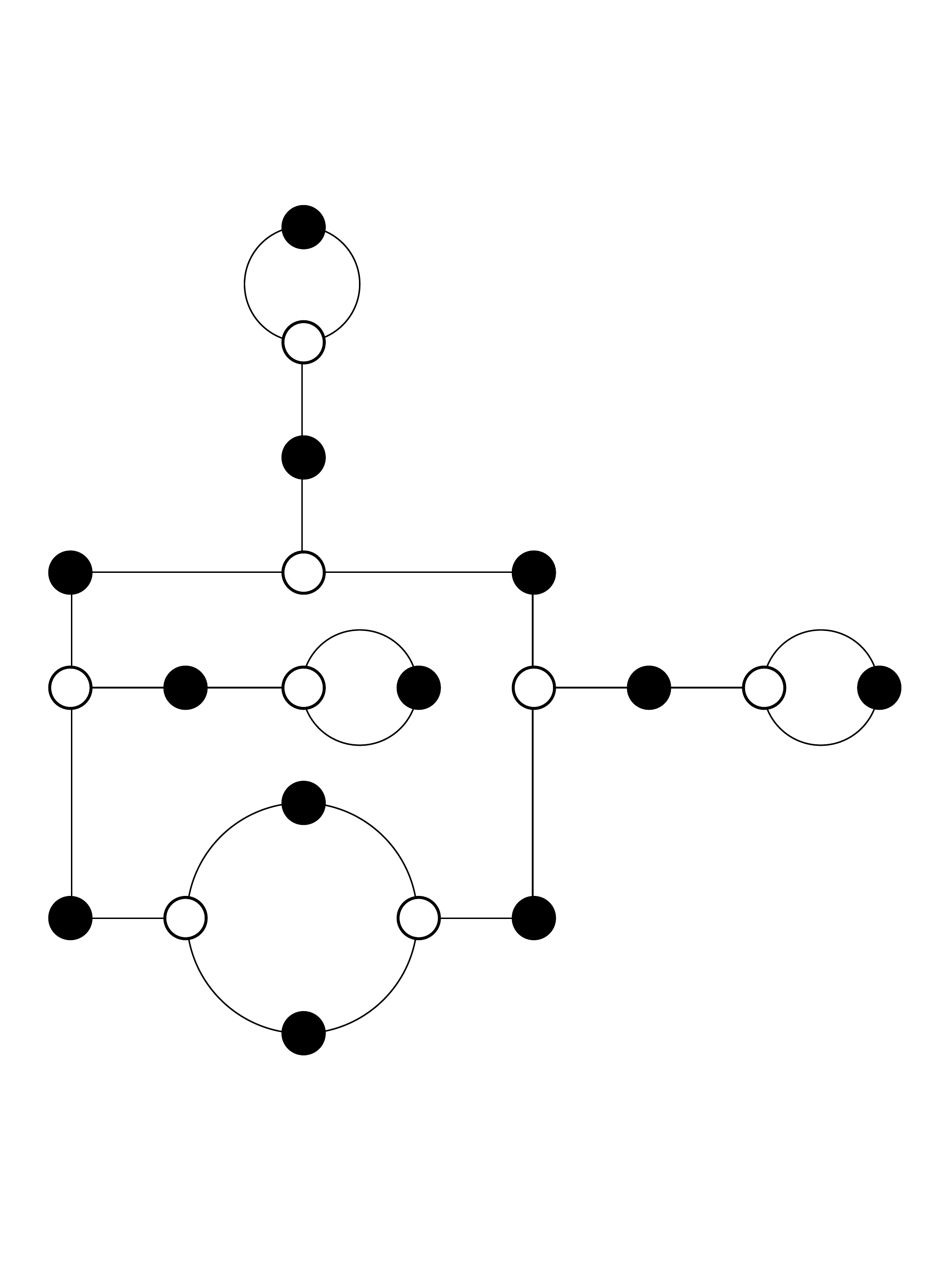}}
\par\end{center}{\scriptsize \par}

\begin{center}
{\scriptsize $11,8,2,1,1,1\;\left(\mathrm{cubic}\right)$}
\par\end{center}%
\end{minipage}{\scriptsize }%
\begin{minipage}[t]{0.33\textwidth}%
\begin{center}
{\scriptsize \includegraphics[scale=0.15]{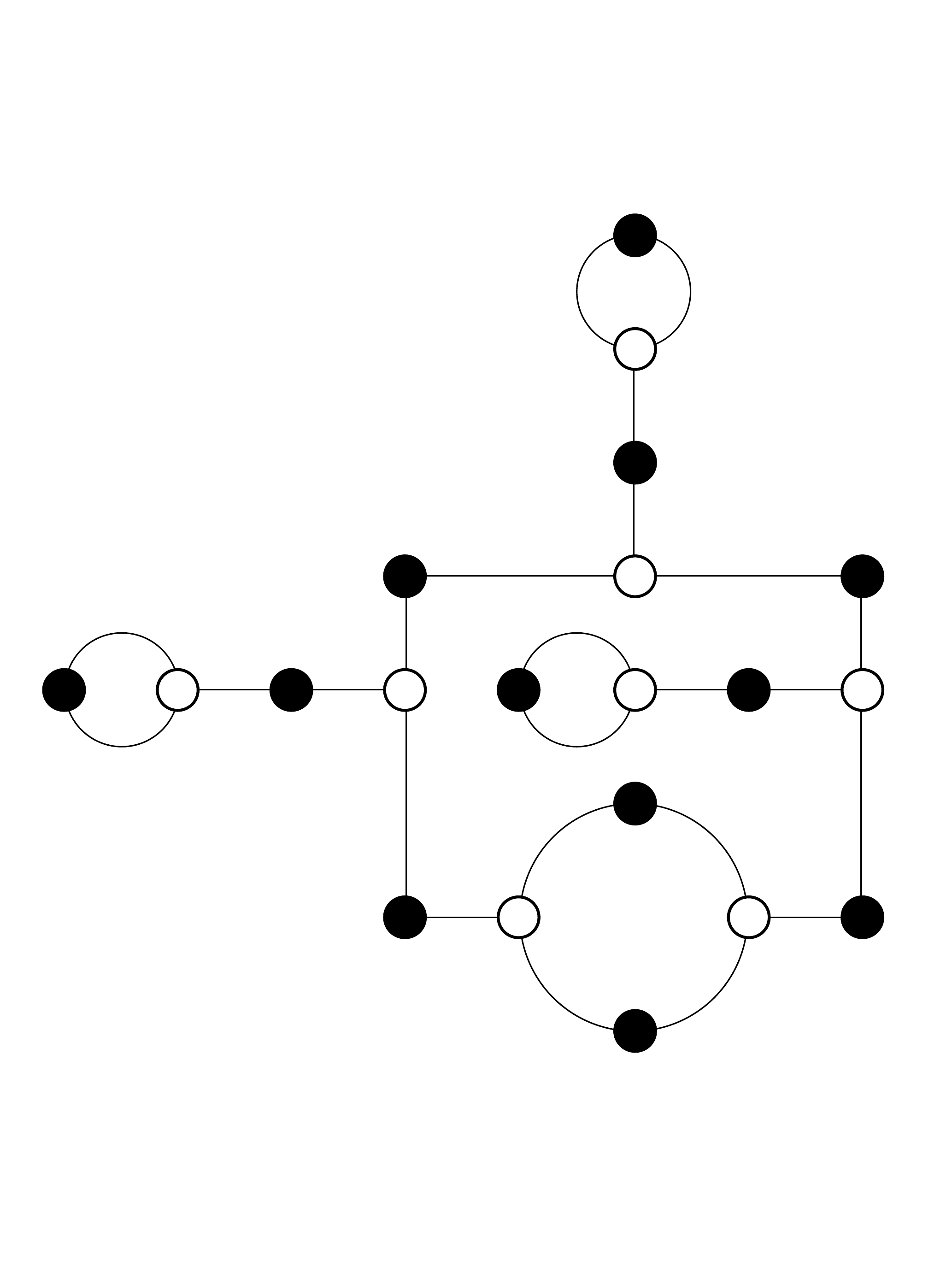}}
\par\end{center}{\scriptsize \par}

\begin{center}
{\scriptsize $11,8,2,1,1,1\;\left(\mathrm{cubic}\right)$}
\par\end{center}%
\end{minipage}{\scriptsize }%
\begin{minipage}[t]{0.33\textwidth}%
\begin{center}
{\scriptsize \includegraphics[scale=0.15]{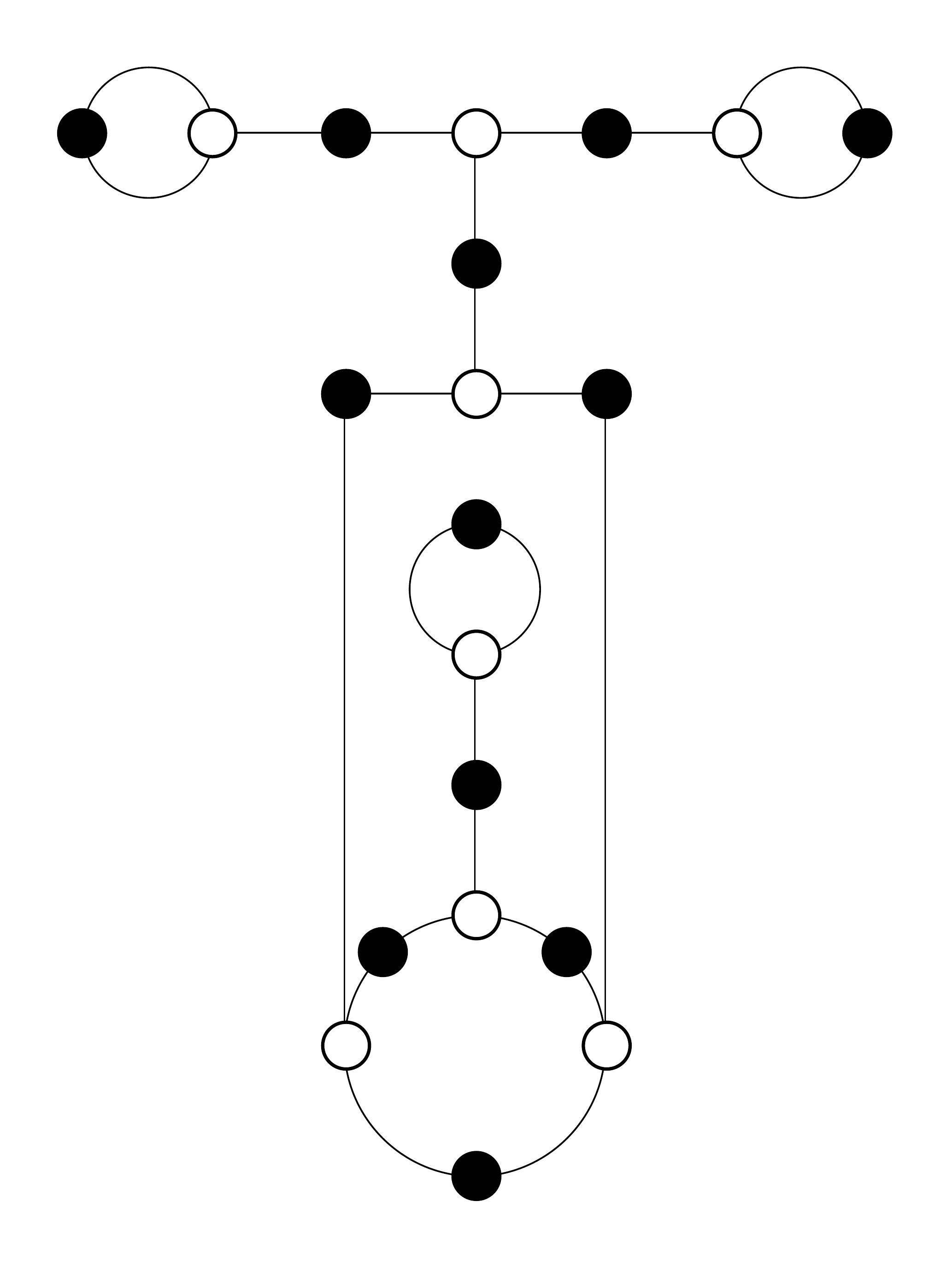}}
\par\end{center}{\scriptsize \par}

\begin{center}
{\scriptsize $11,7,3,1,1,1\;\left(\mathrm{cubic}\right)$}
\par\end{center}%
\end{minipage}
\par\end{center}{\scriptsize \par}

\begin{center}
{\scriptsize }%
\begin{minipage}[t]{0.33\textwidth}%
\begin{center}
{\scriptsize \includegraphics[scale=0.15]{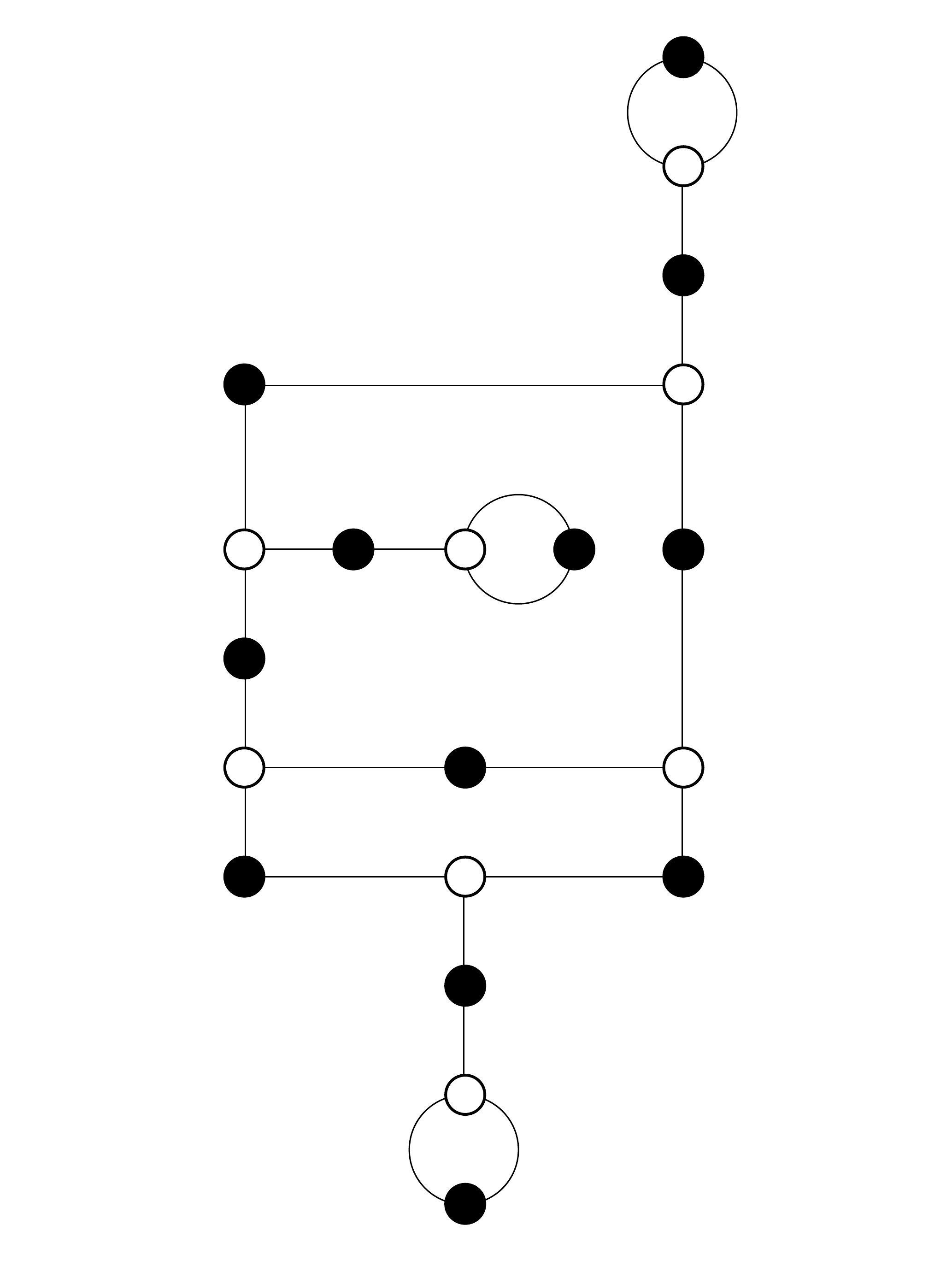}}
\par\end{center}{\scriptsize \par}

\begin{center}
{\scriptsize $11,7,3,1,1,1\;\left(\mathrm{cubic}\right)$}
\par\end{center}%
\end{minipage}{\scriptsize }%
\begin{minipage}[t]{0.33\textwidth}%
\begin{center}
{\scriptsize \includegraphics[scale=0.15]{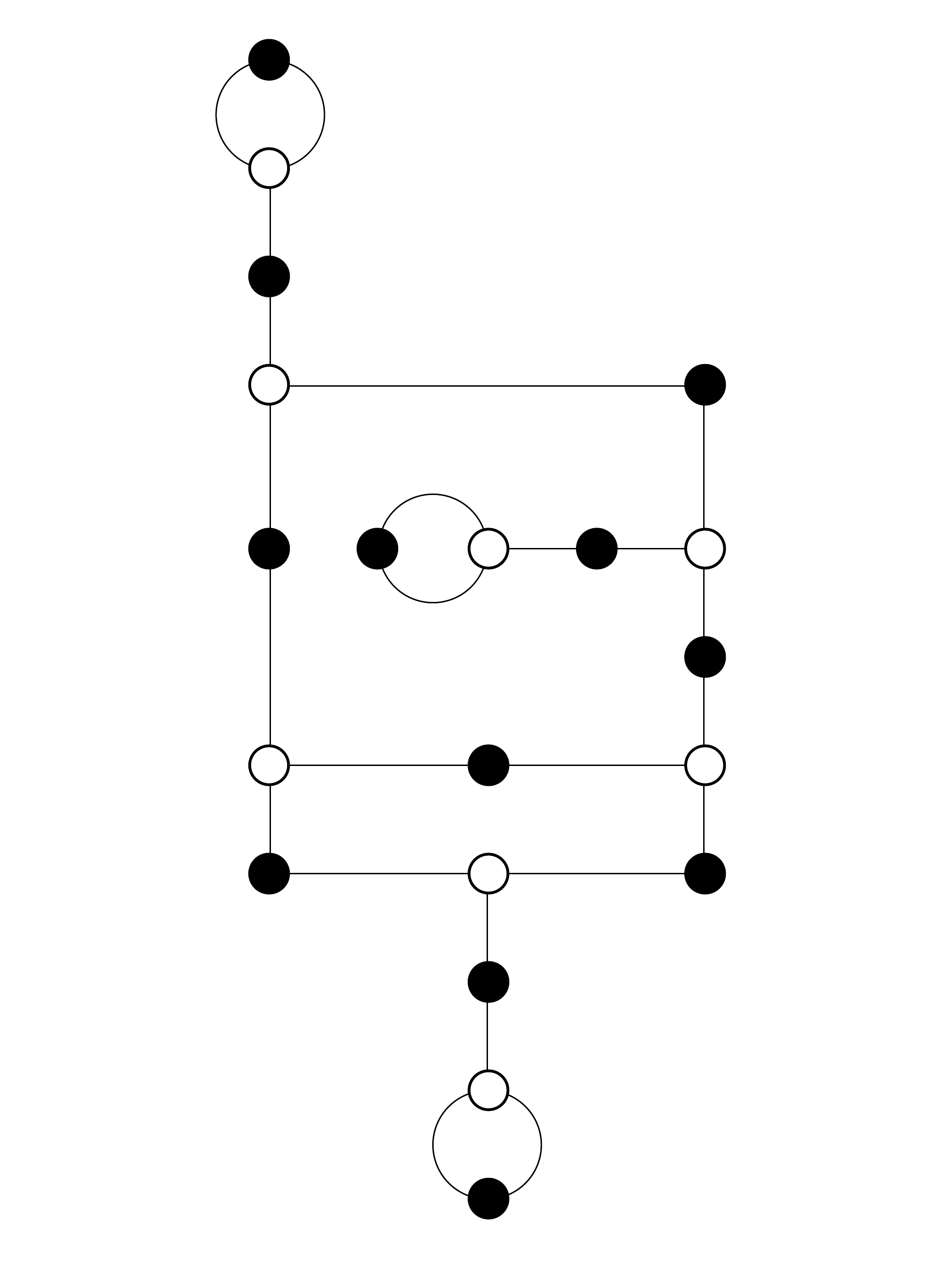}}
\par\end{center}{\scriptsize \par}

\begin{center}
{\scriptsize $11,7,3,1,1,1\;\left(\mathrm{cubic}\right)$}
\par\end{center}%
\end{minipage}{\scriptsize }%
\begin{minipage}[t]{0.33\textwidth}%
\begin{center}
{\scriptsize \includegraphics[scale=0.15]{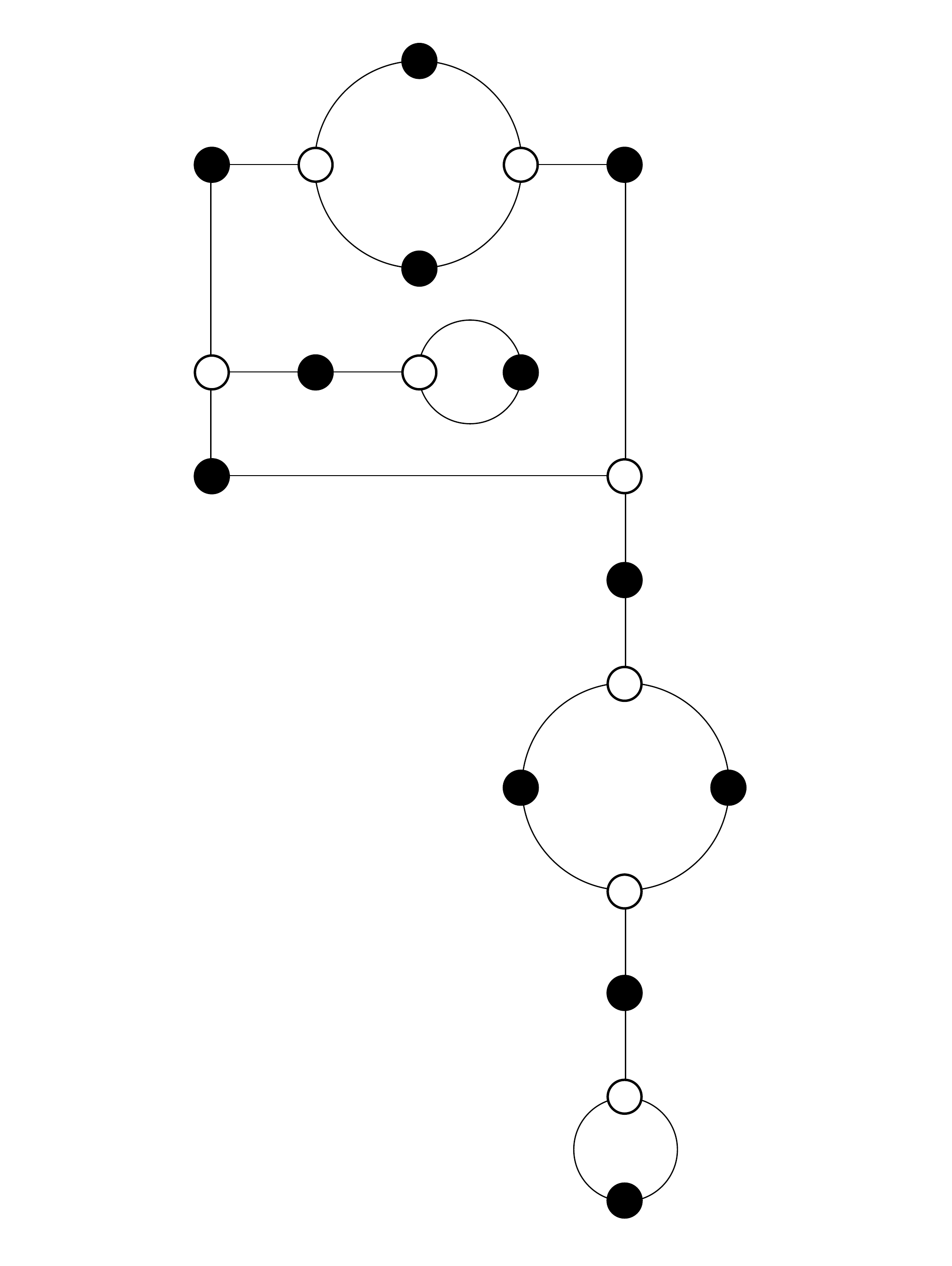}}
\par\end{center}{\scriptsize \par}

\begin{center}
{\scriptsize $11,7,2,2,1,1\;\left(\sqrt{-7}\right)$}
\par\end{center}%
\end{minipage}
\par\end{center}{\scriptsize \par}

\begin{center}
{\scriptsize }%
\begin{minipage}[t]{0.33\textwidth}%
\begin{center}
{\scriptsize \includegraphics[scale=0.15]{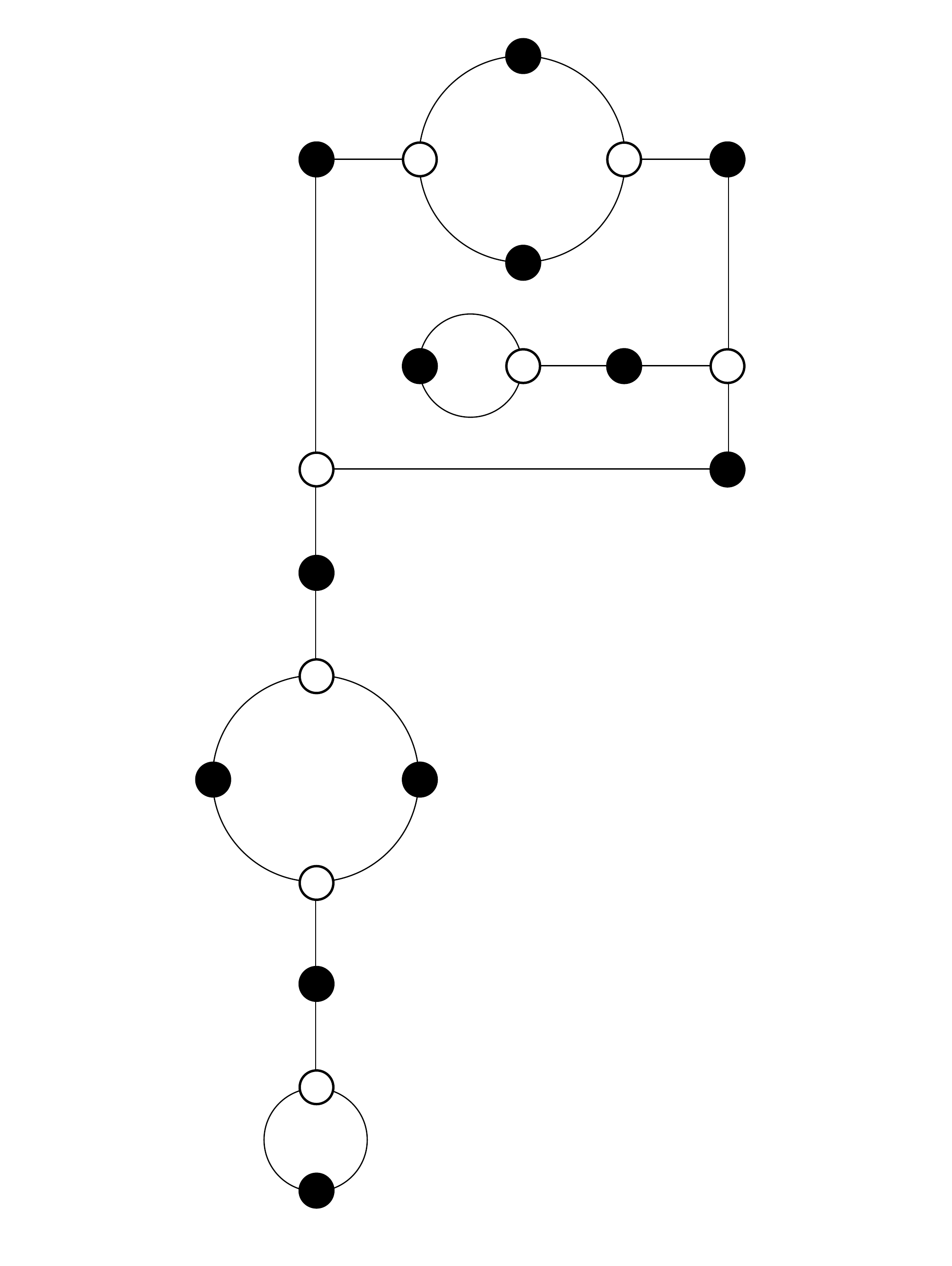}}
\par\end{center}{\scriptsize \par}

\begin{center}
{\scriptsize $11,7,2,2,1,1\;\left(\sqrt{-7}\right)$}
\par\end{center}%
\end{minipage}{\scriptsize }%
\begin{minipage}[t]{0.33\textwidth}%
\begin{center}
{\scriptsize \includegraphics[scale=0.15]{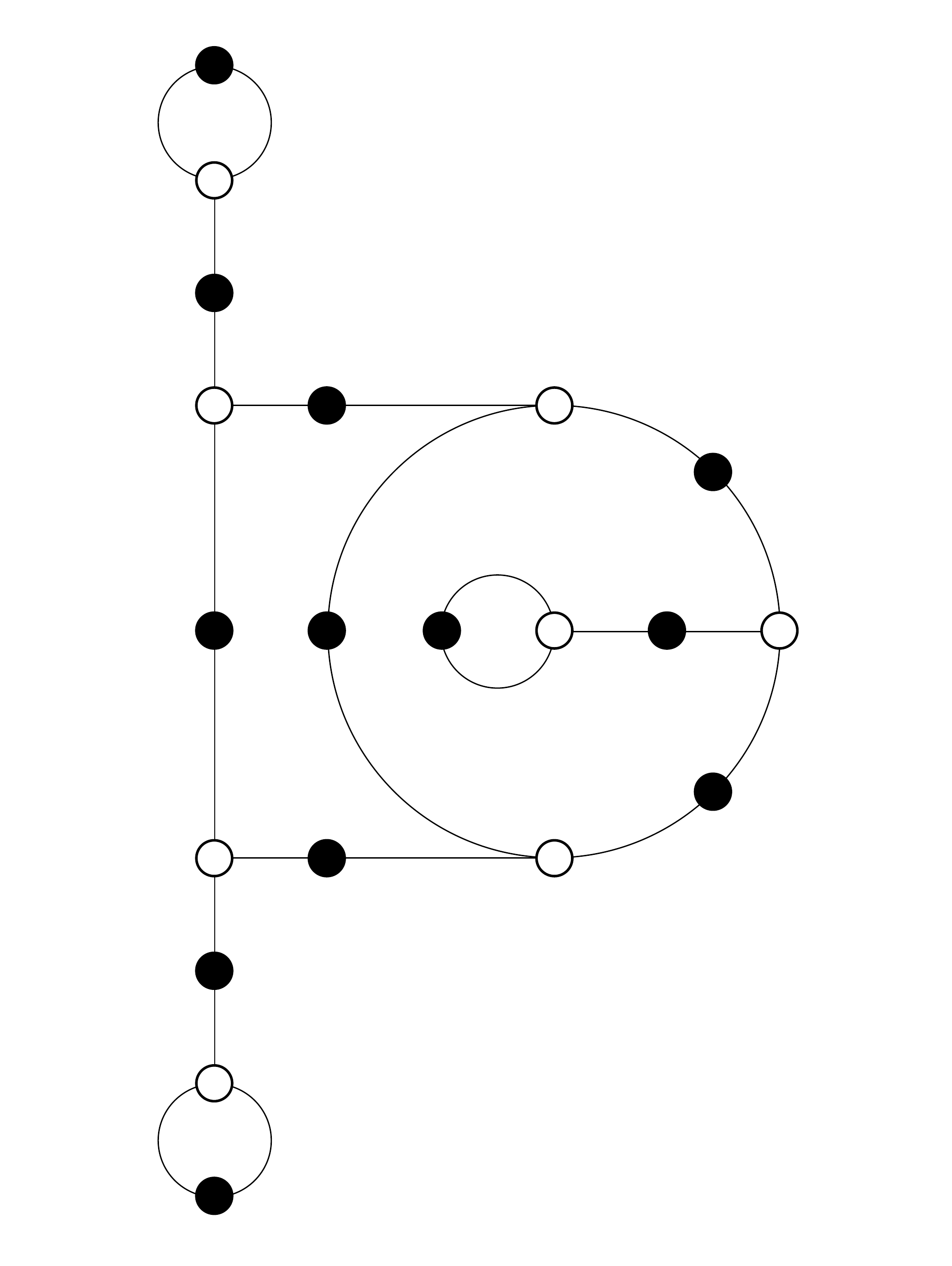}}
\par\end{center}{\scriptsize \par}

\begin{center}
{\scriptsize $11,6,4,1,1,1\;\left(\sqrt{33}\right)$}
\par\end{center}%
\end{minipage}{\scriptsize }%
\begin{minipage}[t]{0.33\textwidth}%
\begin{center}
{\scriptsize \includegraphics[scale=0.15]{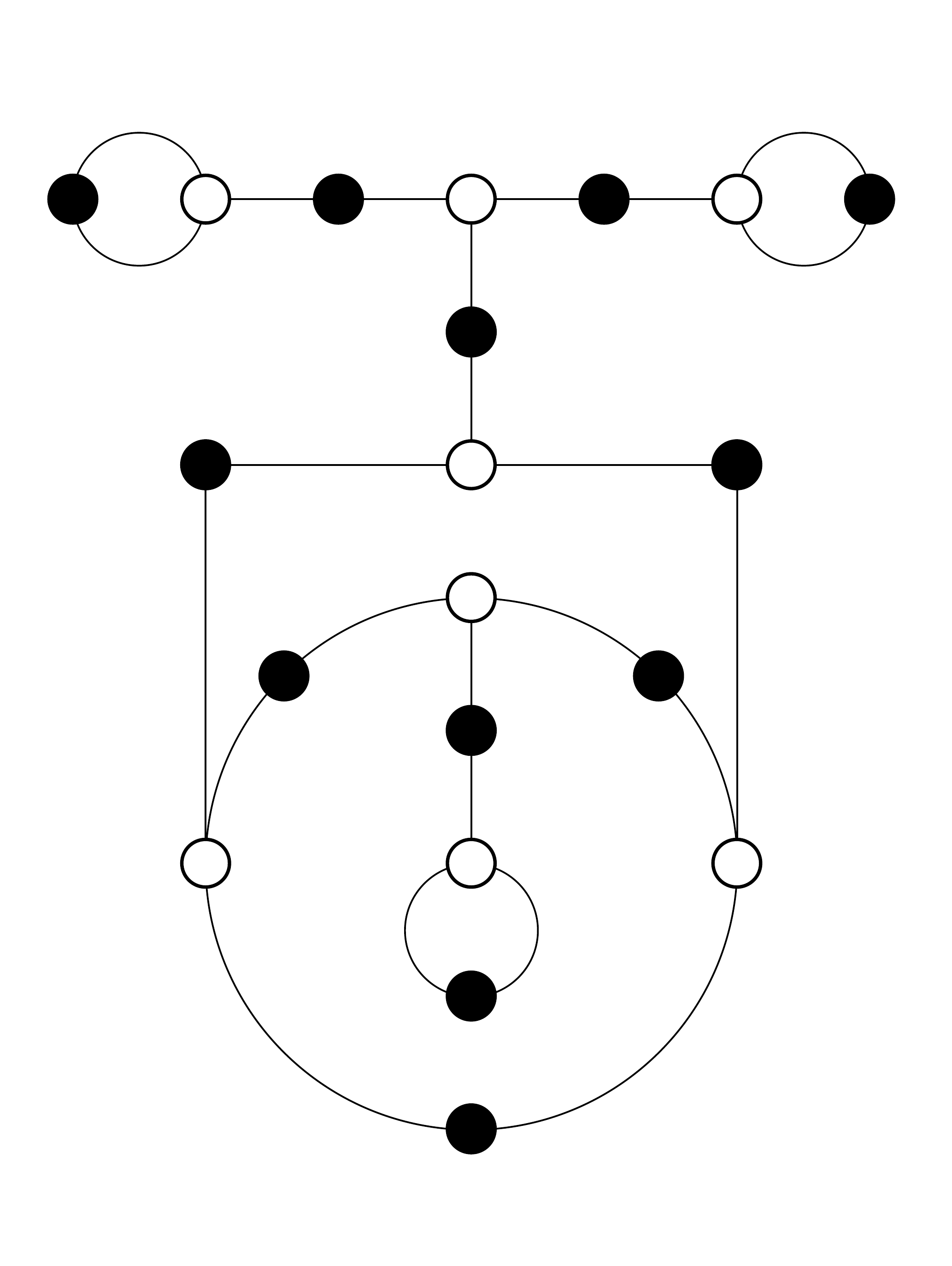}}
\par\end{center}{\scriptsize \par}

\begin{center}
{\scriptsize $11,6,4,1,1,1\;\left(\sqrt{33}\right)$}
\par\end{center}%
\end{minipage}
\par\end{center}{\scriptsize \par}

\begin{center}
{\scriptsize }%
\begin{minipage}[t]{0.33\textwidth}%
\begin{center}
{\scriptsize \includegraphics[scale=0.15]{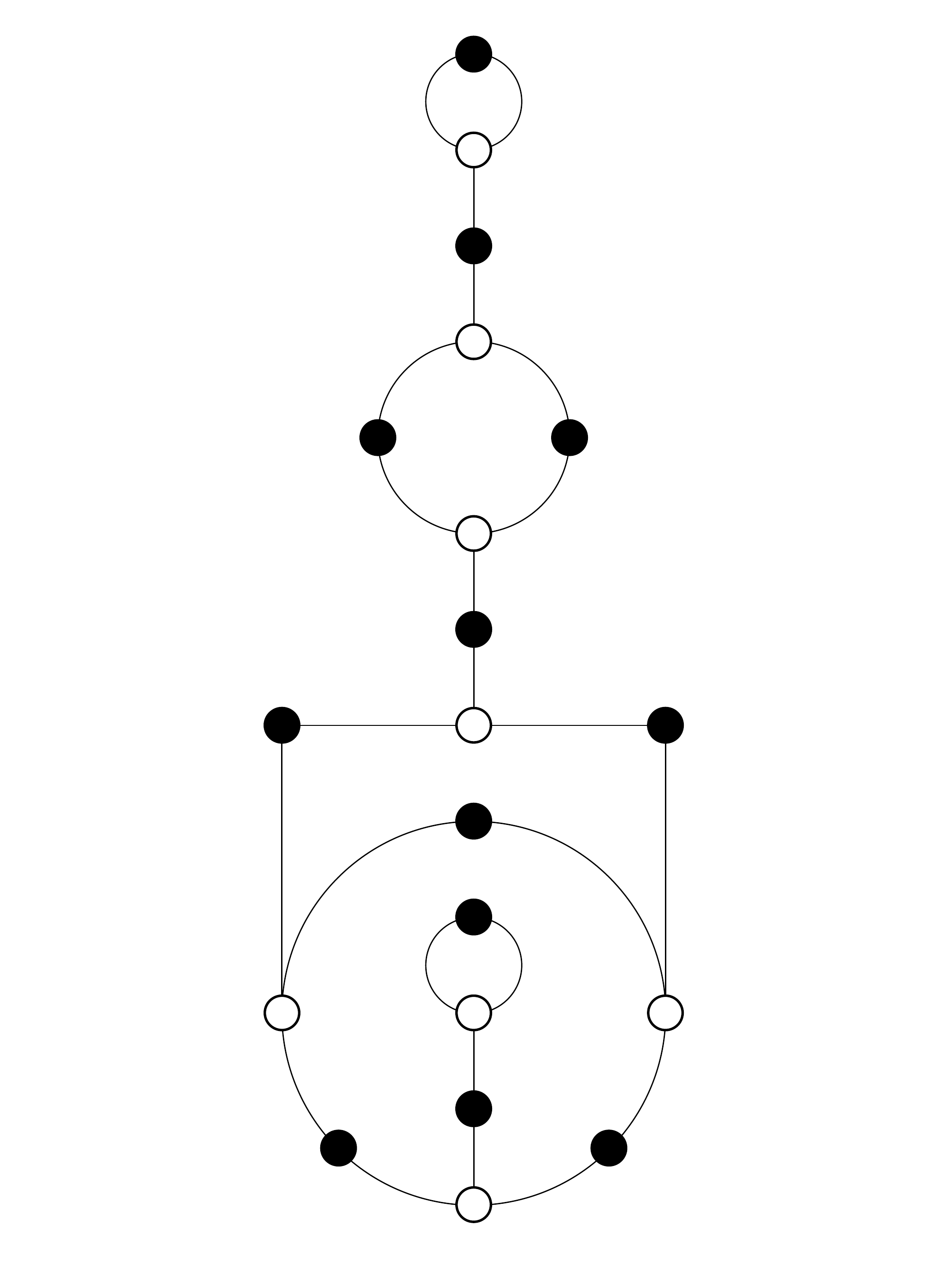}}
\par\end{center}{\scriptsize \par}

\begin{center}
{\scriptsize $11,6,3,2,1,1\;\left(\mathrm{cubic}\right)$}
\par\end{center}%
\end{minipage}{\scriptsize }%
\begin{minipage}[t]{0.33\textwidth}%
\begin{center}
{\scriptsize \includegraphics[scale=0.15]{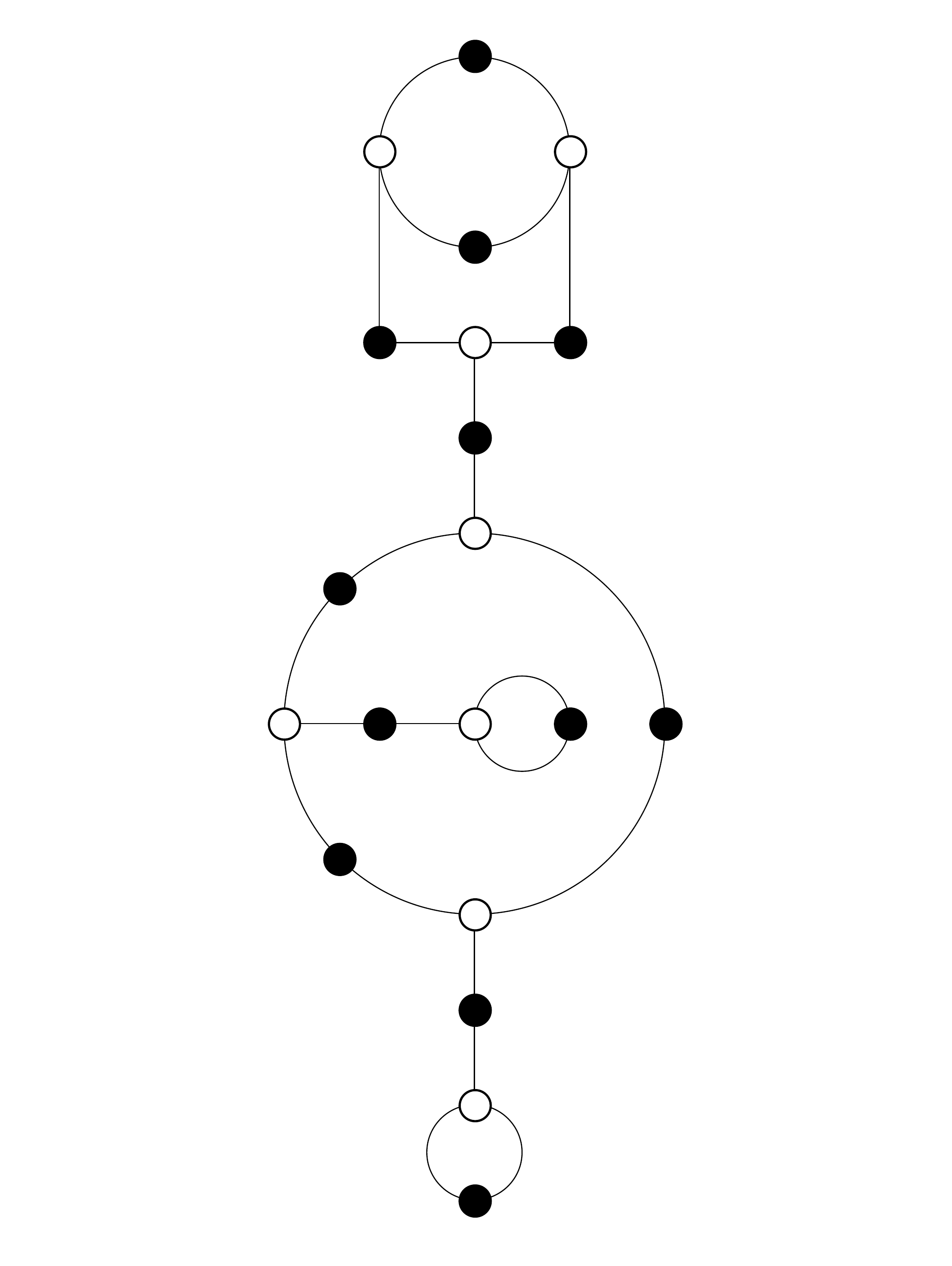}}
\par\end{center}{\scriptsize \par}

\begin{center}
{\scriptsize $11,6,3,2,1,1\;\left(\mathrm{cubic}\right)$}
\par\end{center}%
\end{minipage}{\scriptsize }%
\begin{minipage}[t]{0.33\textwidth}%
\begin{center}
{\scriptsize \includegraphics[scale=0.15]{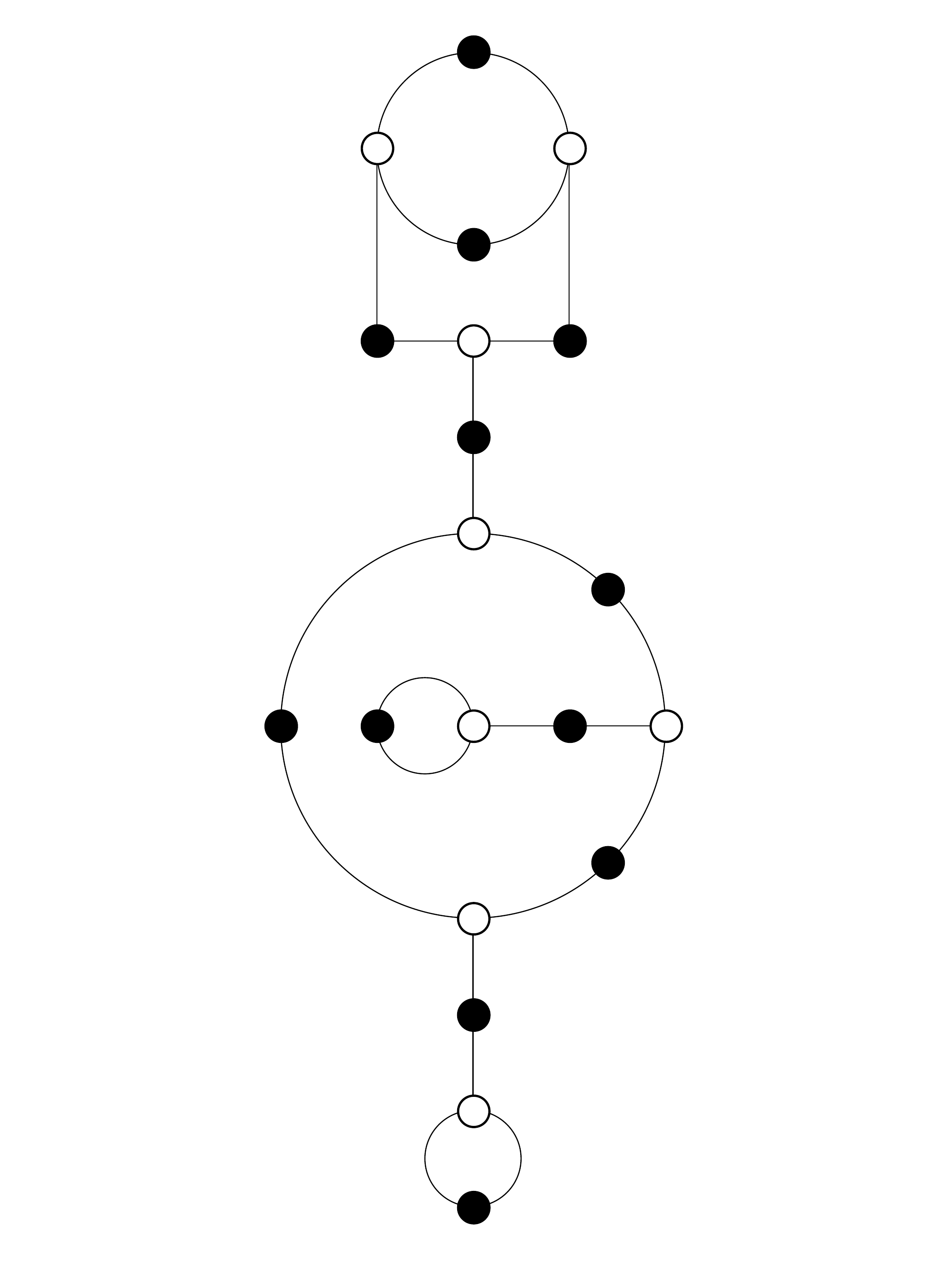}}
\par\end{center}{\scriptsize \par}

\begin{center}
{\scriptsize $11,6,3,2,1,1\;\left(\mathrm{cubic}\right)$}
\par\end{center}%
\end{minipage}
\par\end{center}{\scriptsize \par}

\begin{center}
{\scriptsize }%
\begin{minipage}[t]{0.33\textwidth}%
\begin{center}
{\scriptsize \includegraphics[scale=0.15]{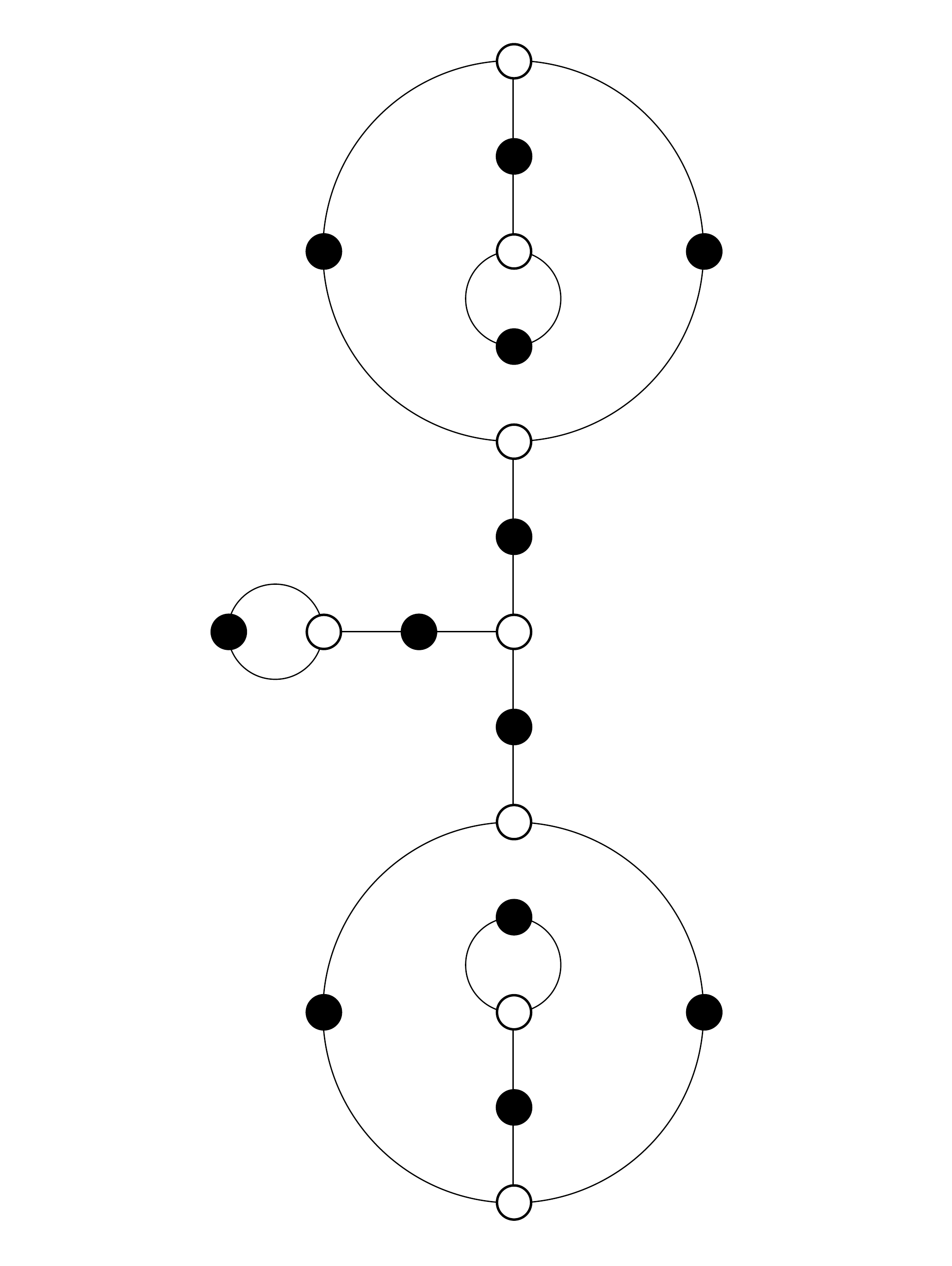}}
\par\end{center}{\scriptsize \par}

\begin{center}
{\scriptsize $11,5,5,1,1,1\;\left(\mathbb{Q}\right)$}
\par\end{center}%
\end{minipage}{\scriptsize }%
\begin{minipage}[t]{0.33\textwidth}%
\begin{center}
{\scriptsize \includegraphics[scale=0.15]{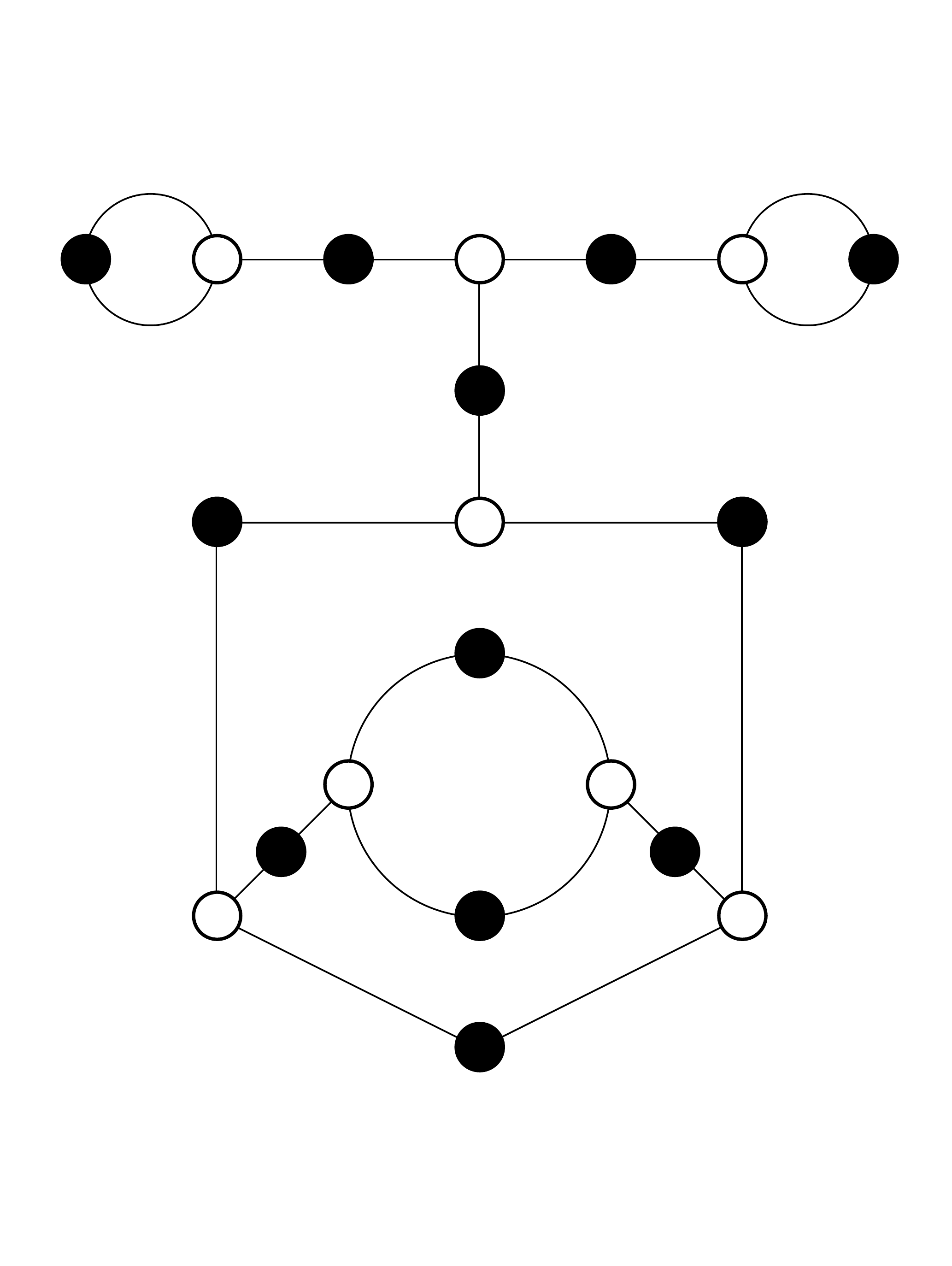}}
\par\end{center}{\scriptsize \par}

\begin{center}
{\scriptsize $11,5,4,2,1,1\;\left(\mathrm{cubic}\right)$}
\par\end{center}%
\end{minipage}{\scriptsize }%
\begin{minipage}[t]{0.33\textwidth}%
\begin{center}
{\scriptsize \includegraphics[scale=0.15]{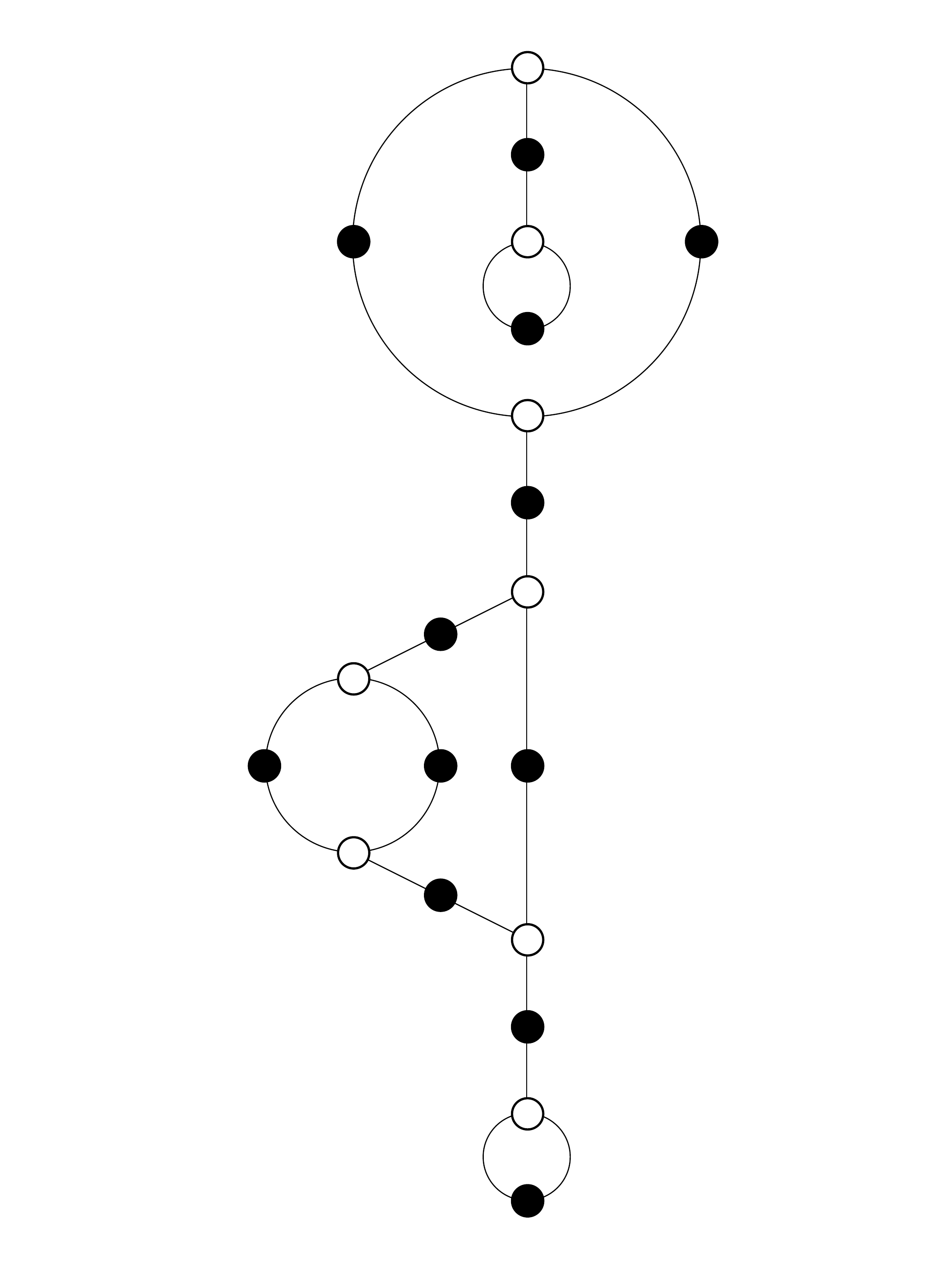}}
\par\end{center}{\scriptsize \par}

\begin{center}
{\scriptsize $11,5,4,2,1,1\;\left(\mathrm{cubic}\right)$}
\par\end{center}%
\end{minipage}
\par\end{center}{\scriptsize \par}

\begin{center}
{\scriptsize }%
\begin{minipage}[t]{0.33\textwidth}%
\begin{center}
{\scriptsize \includegraphics[scale=0.15]{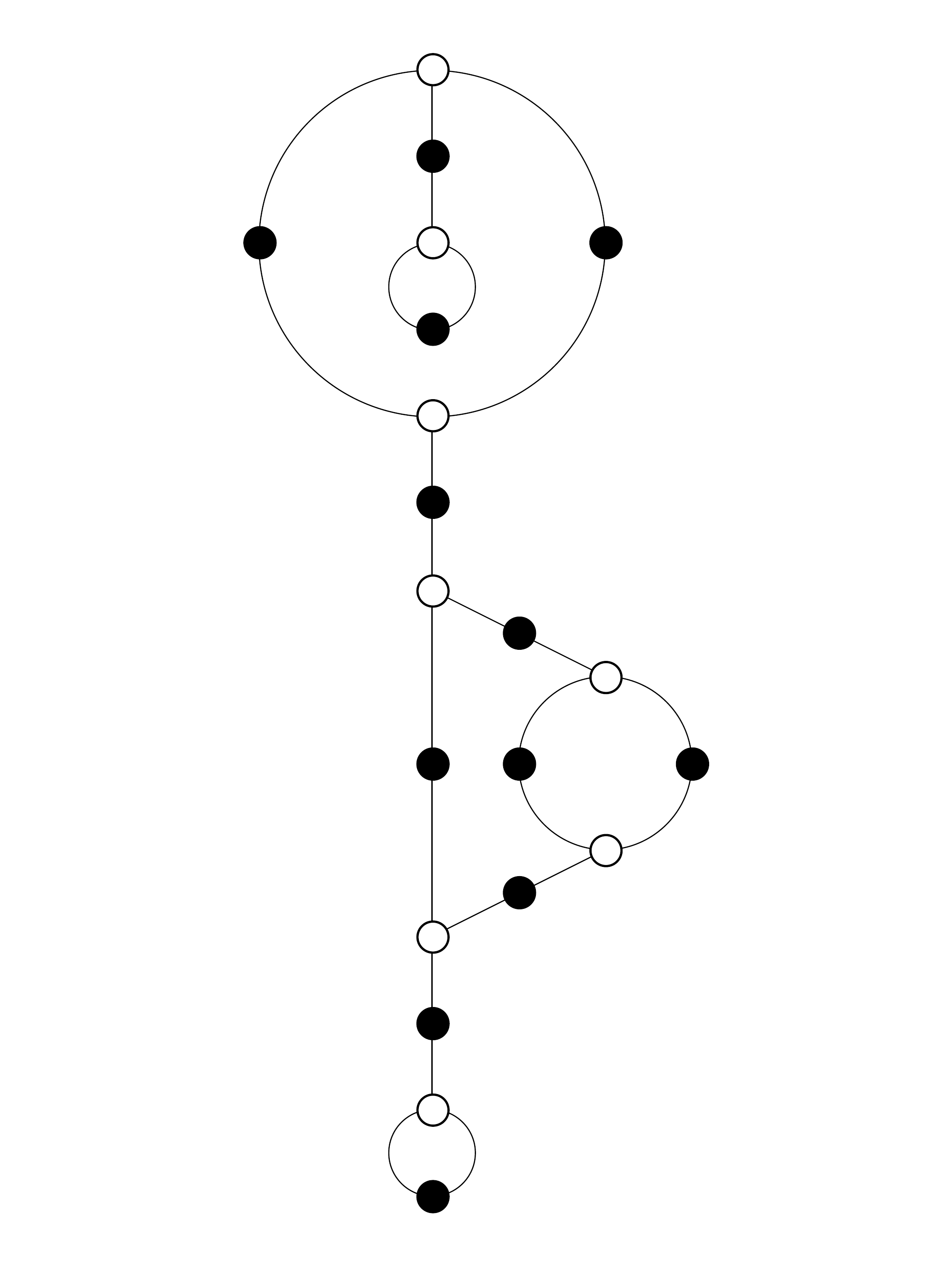}}
\par\end{center}{\scriptsize \par}

\begin{center}
{\scriptsize $11,5,4,2,1,1\;\left(\mathrm{cubic}\right)$}
\par\end{center}%
\end{minipage}{\scriptsize }%
\begin{minipage}[t]{0.33\textwidth}%
\begin{center}
{\scriptsize \includegraphics[scale=0.15]{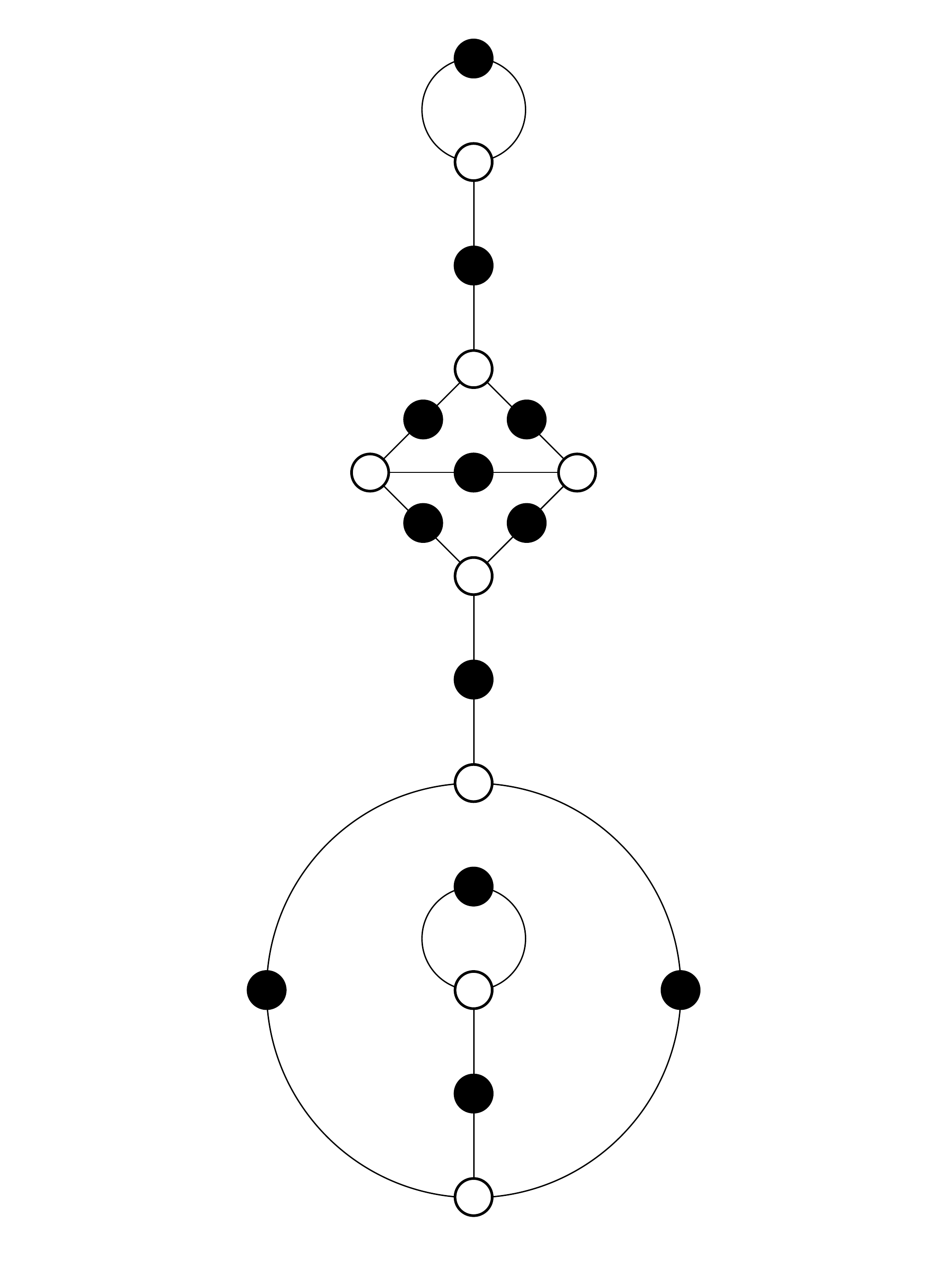}}
\par\end{center}{\scriptsize \par}

\begin{center}
{\scriptsize $11,5,3,3,1,1\;\left(\sqrt{5}\right)$}
\par\end{center}%
\end{minipage}{\scriptsize }%
\begin{minipage}[t]{0.33\textwidth}%
\begin{center}
{\scriptsize \includegraphics[scale=0.15]{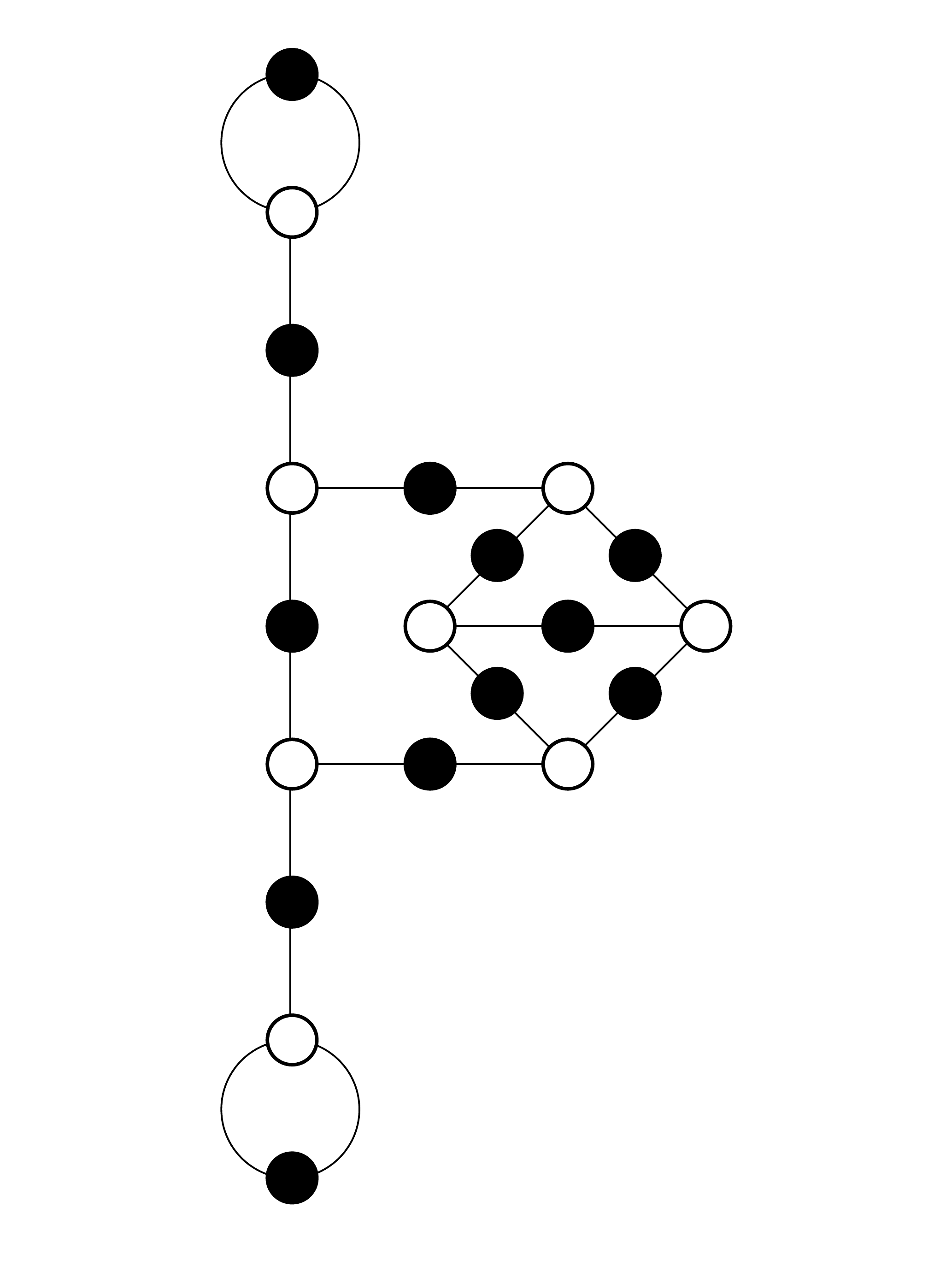}}
\par\end{center}{\scriptsize \par}

\begin{center}
{\scriptsize $11,5,3,3,1,1\;\left(\sqrt{5}\right)$}
\par\end{center}%
\end{minipage}
\par\end{center}{\scriptsize \par}

\begin{center}
{\scriptsize }%
\begin{minipage}[t]{0.33\textwidth}%
\begin{center}
{\scriptsize \includegraphics[scale=0.15]{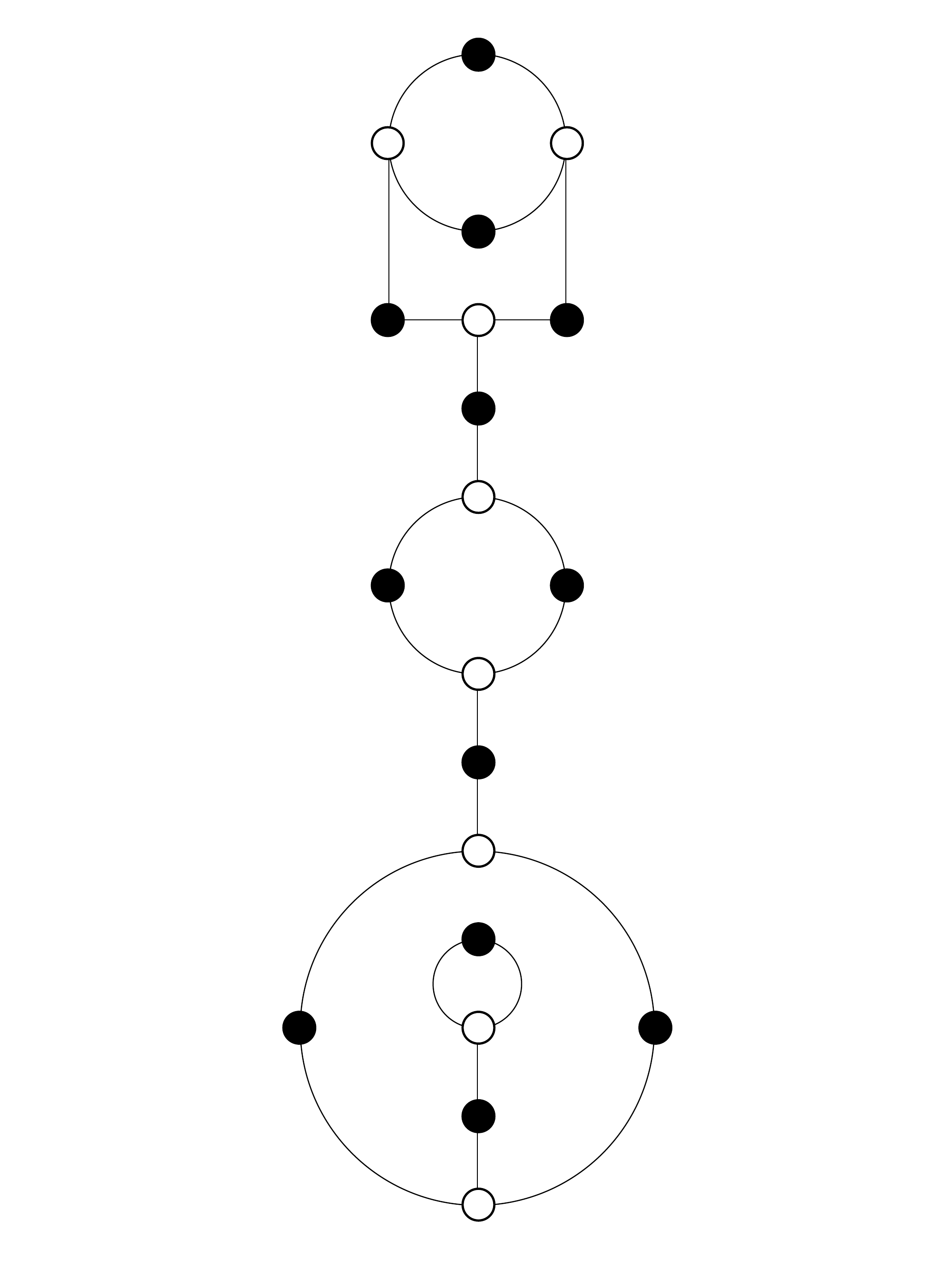}}
\par\end{center}{\scriptsize \par}

\begin{center}
{\scriptsize $11,5,3,2,2,1\;\left(\mathbb{Q}\right)$}
\par\end{center}%
\end{minipage}{\scriptsize }%
\begin{minipage}[t]{0.33\textwidth}%
\begin{center}
{\scriptsize \includegraphics[scale=0.15]{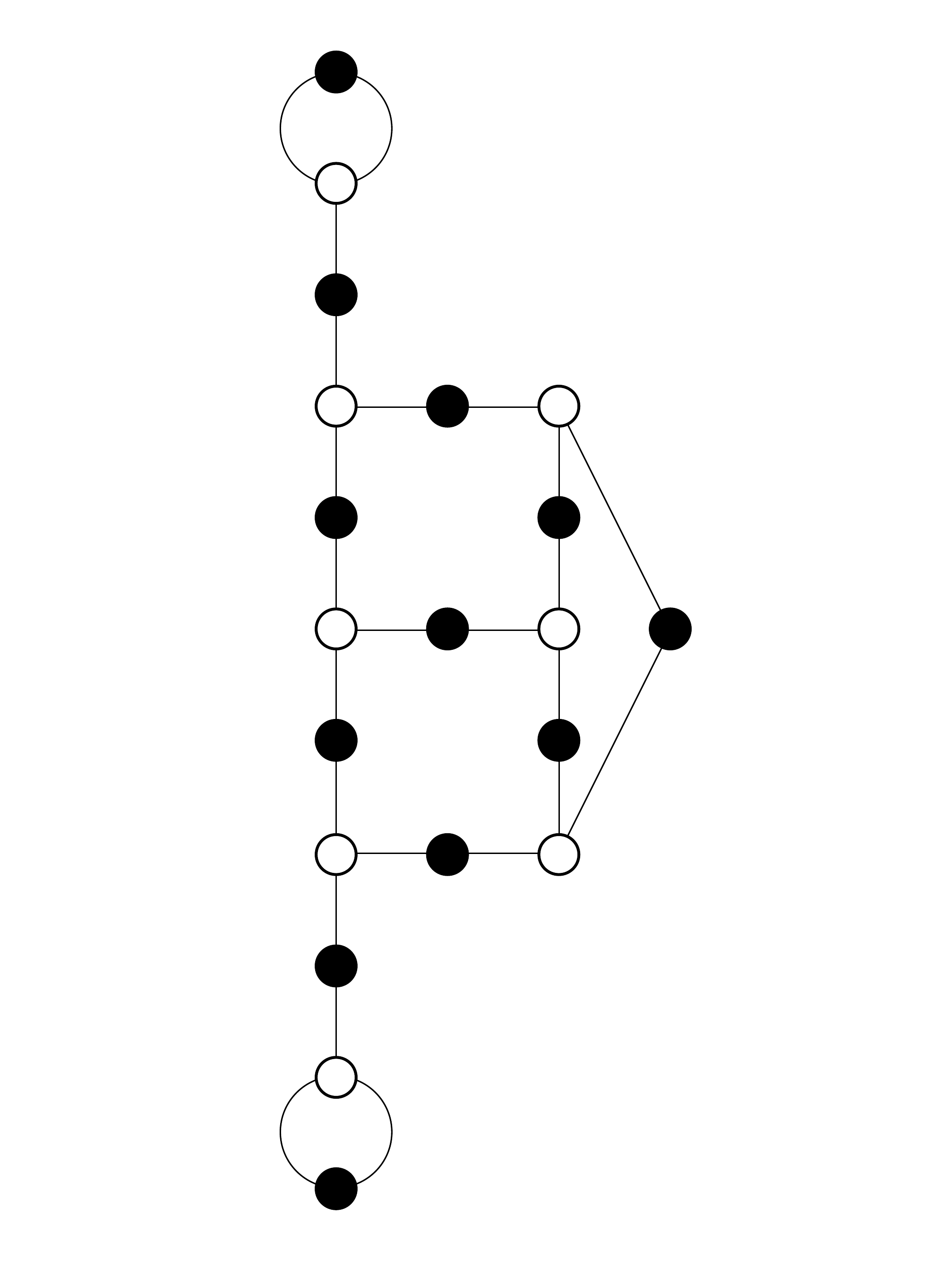}}
\par\end{center}{\scriptsize \par}

\begin{center}
{\scriptsize $11,4,4,3,1,1\;\left(\mathbb{Q}\right)$}
\par\end{center}%
\end{minipage}{\scriptsize }%
\begin{minipage}[t]{0.33\textwidth}%
\begin{center}
{\scriptsize \includegraphics[scale=0.15]{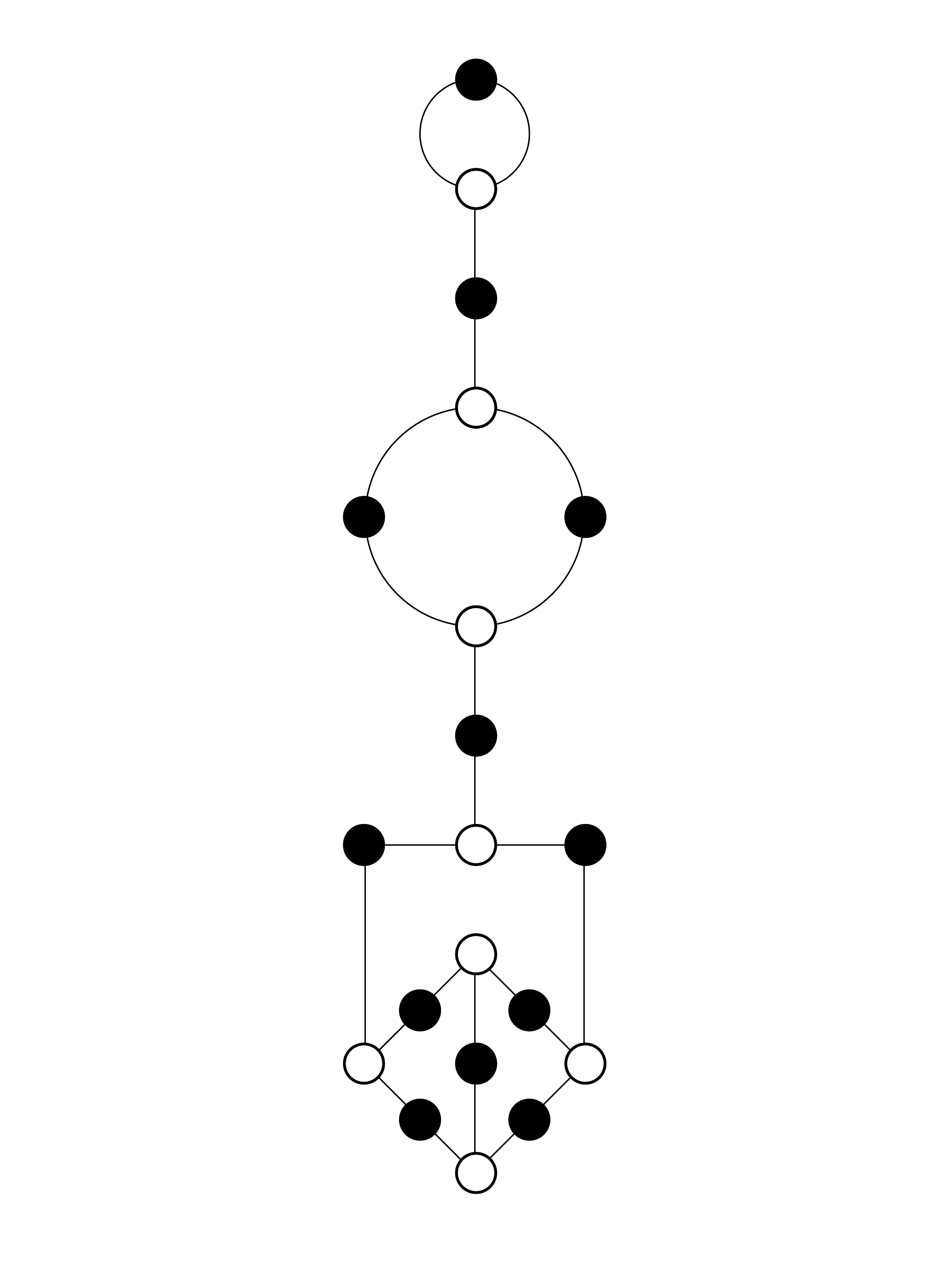}}
\par\end{center}{\scriptsize \par}

\begin{center}
{\scriptsize $11,4,3,3,2,1\;\left(\mathbb{Q}\right)$}
\par\end{center}%
\end{minipage}
\par\end{center}{\scriptsize \par}

\begin{center}
{\scriptsize }%
\begin{minipage}[t]{0.33\textwidth}%
\begin{center}
{\scriptsize \includegraphics[scale=0.15]{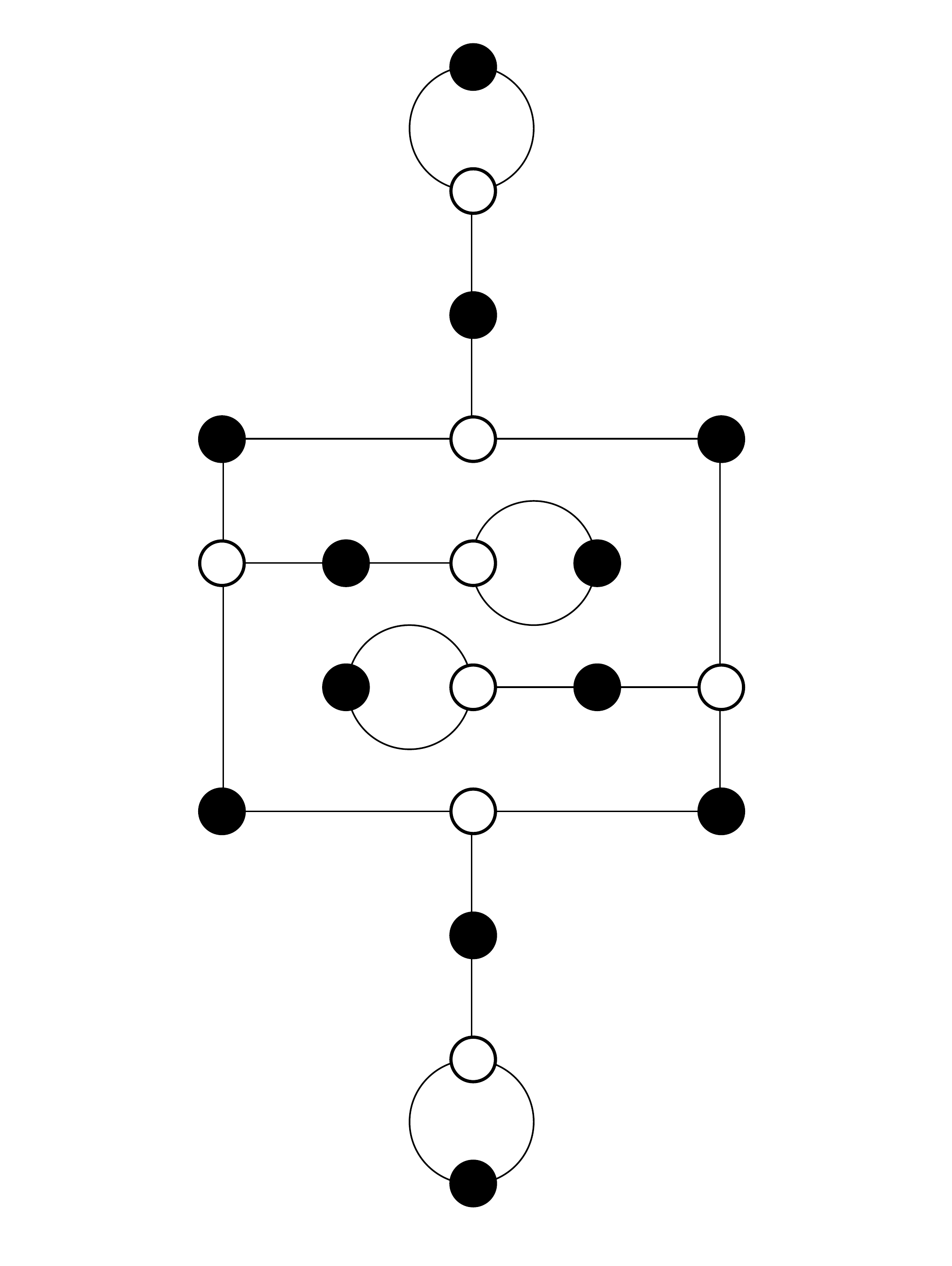}}
\par\end{center}{\scriptsize \par}

\begin{center}
{\scriptsize $10,10,1,1,1,1\;\left(\mathbb{Q}\right)$}
\par\end{center}%
\end{minipage}{\scriptsize }%
\begin{minipage}[t]{0.33\textwidth}%
\begin{center}
{\scriptsize \includegraphics[scale=0.15]{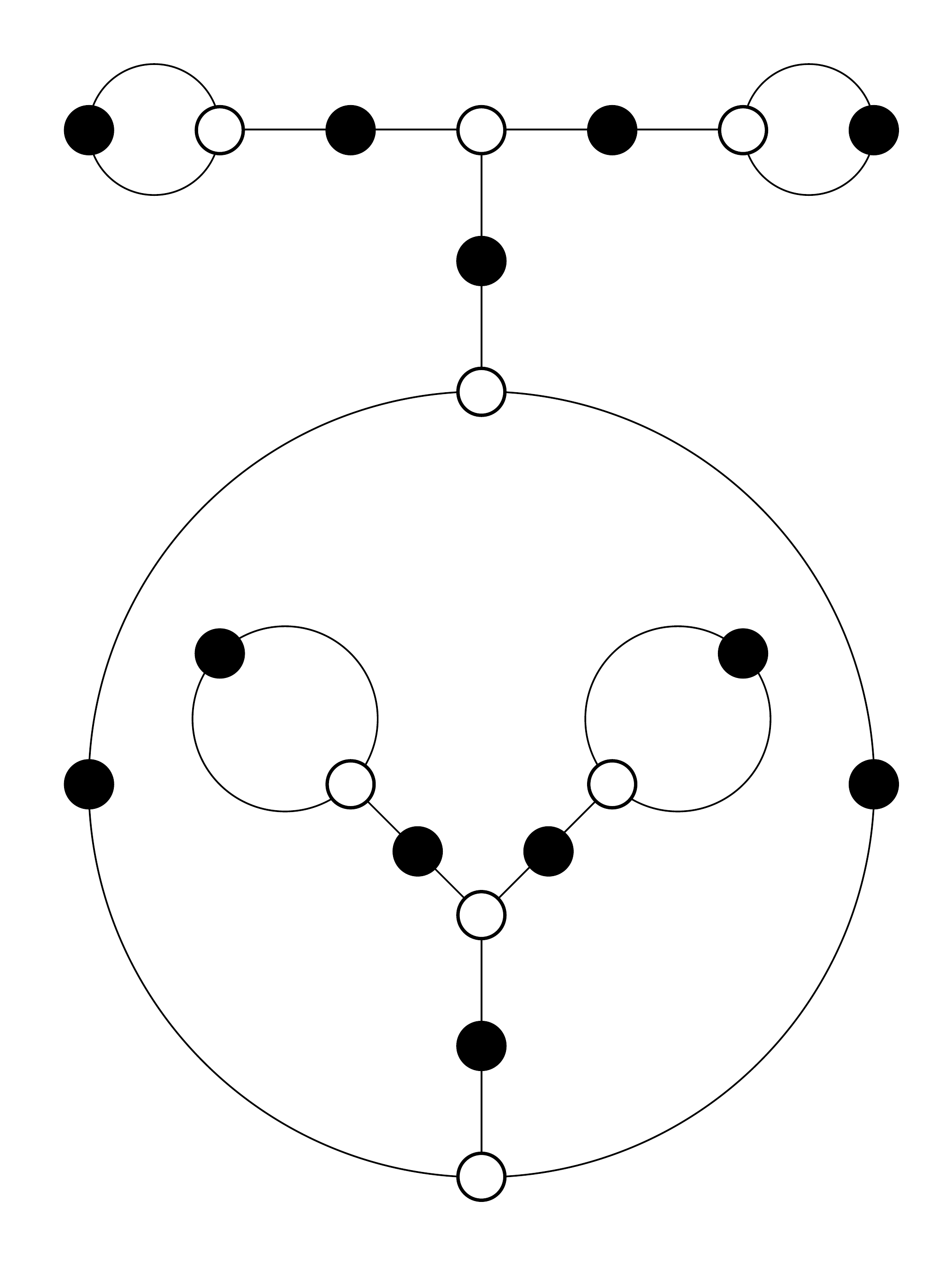}}
\par\end{center}{\scriptsize \par}

\begin{center}
{\scriptsize $10,10,1,1,1,1\;\left(\sqrt{5}\right)$}
\par\end{center}%
\end{minipage}{\scriptsize }%
\begin{minipage}[t]{0.33\textwidth}%
\begin{center}
{\scriptsize \includegraphics[scale=0.15]{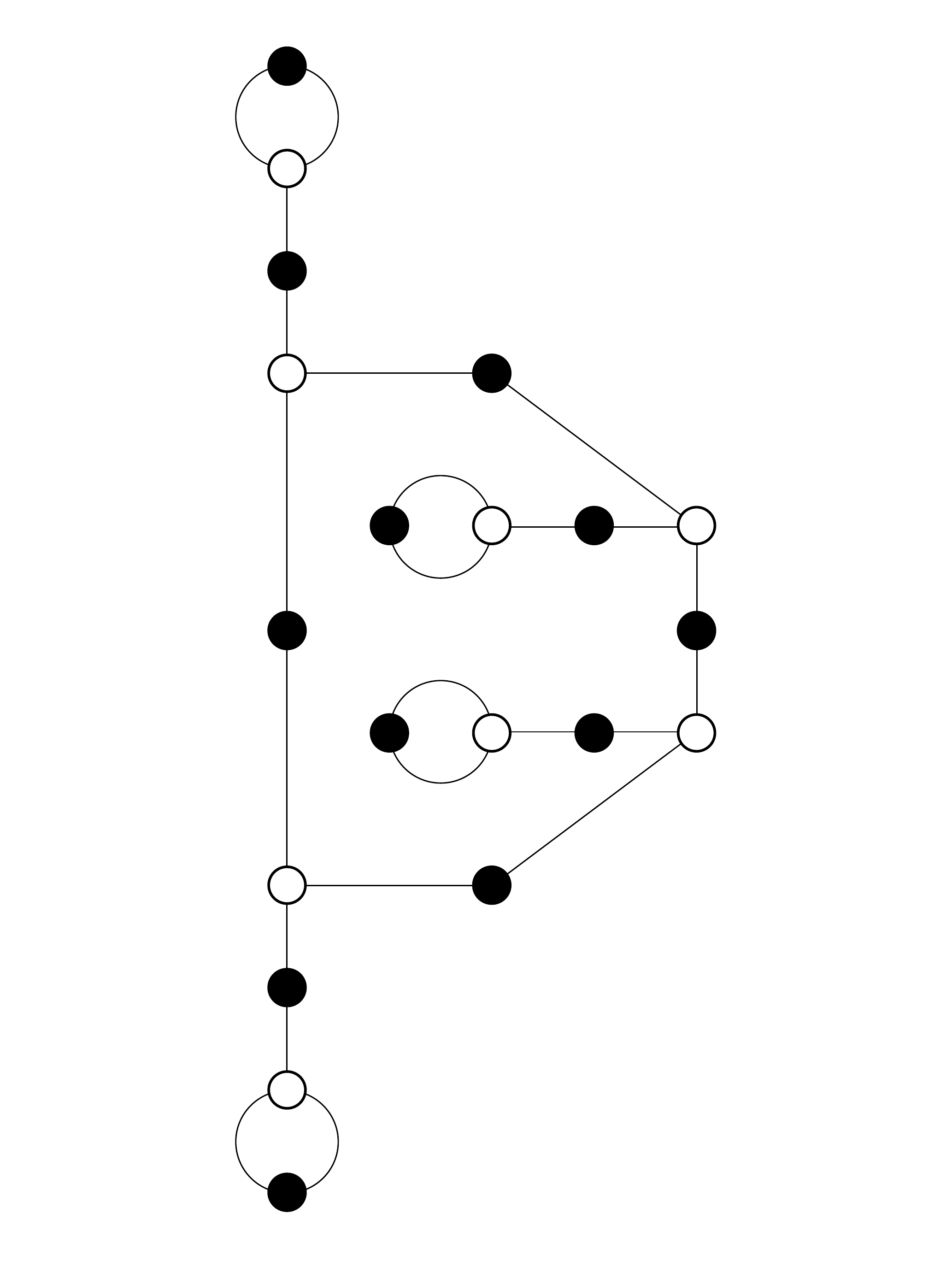}}
\par\end{center}{\scriptsize \par}

\begin{center}
{\scriptsize $10,10,1,1,1,1\;\left(\sqrt{5}\right)$}
\par\end{center}%
\end{minipage}
\par\end{center}{\scriptsize \par}

\begin{center}
{\scriptsize }%
\begin{minipage}[t]{0.33\textwidth}%
\begin{center}
{\scriptsize \includegraphics[scale=0.15]{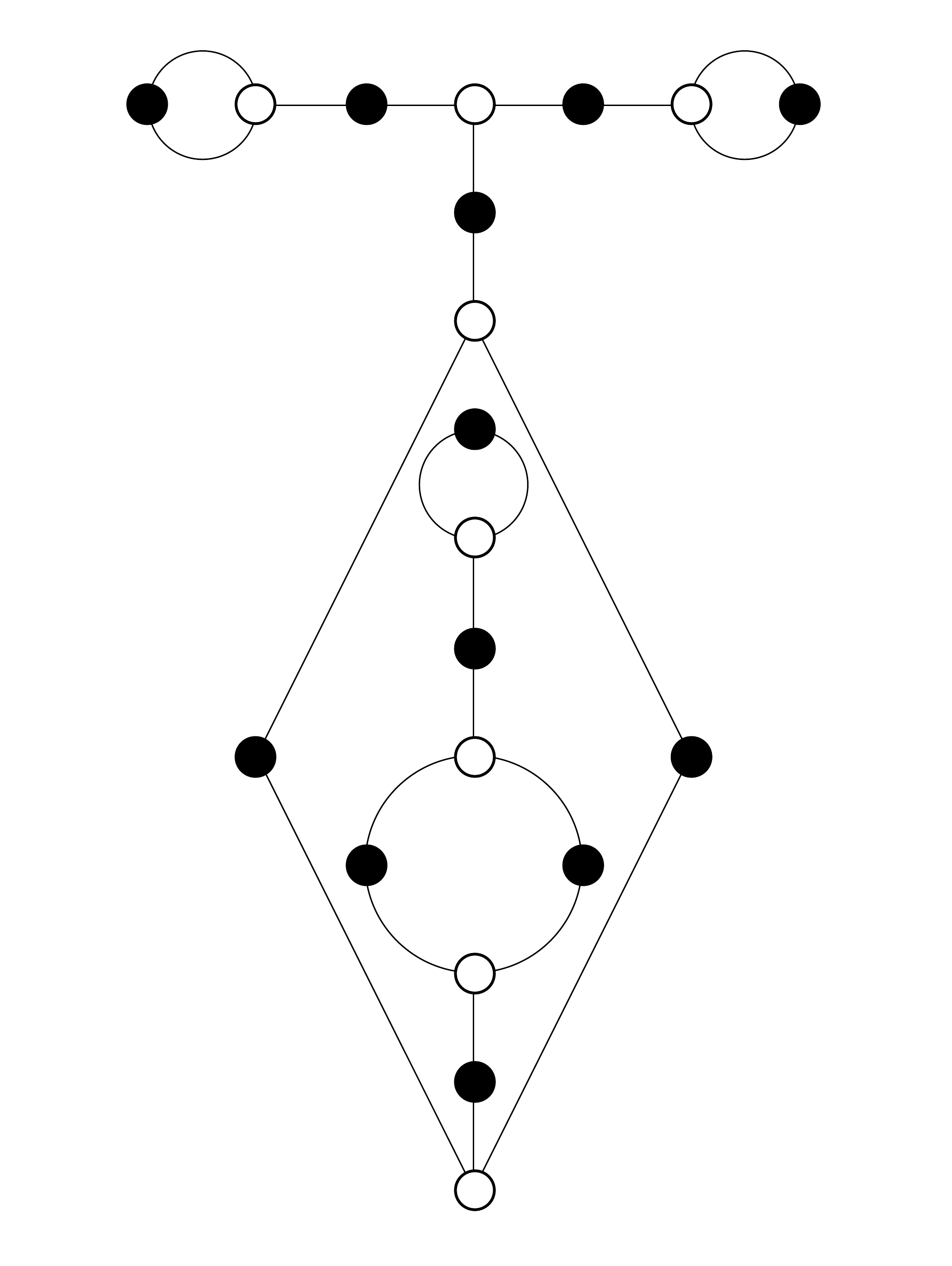}}
\par\end{center}{\scriptsize \par}

\begin{center}
{\scriptsize $10,9,2,1,1,1\;\left(\sqrt{5}\right)$}
\par\end{center}%
\end{minipage}{\scriptsize }%
\begin{minipage}[t]{0.33\textwidth}%
\begin{center}
{\scriptsize \includegraphics[scale=0.15]{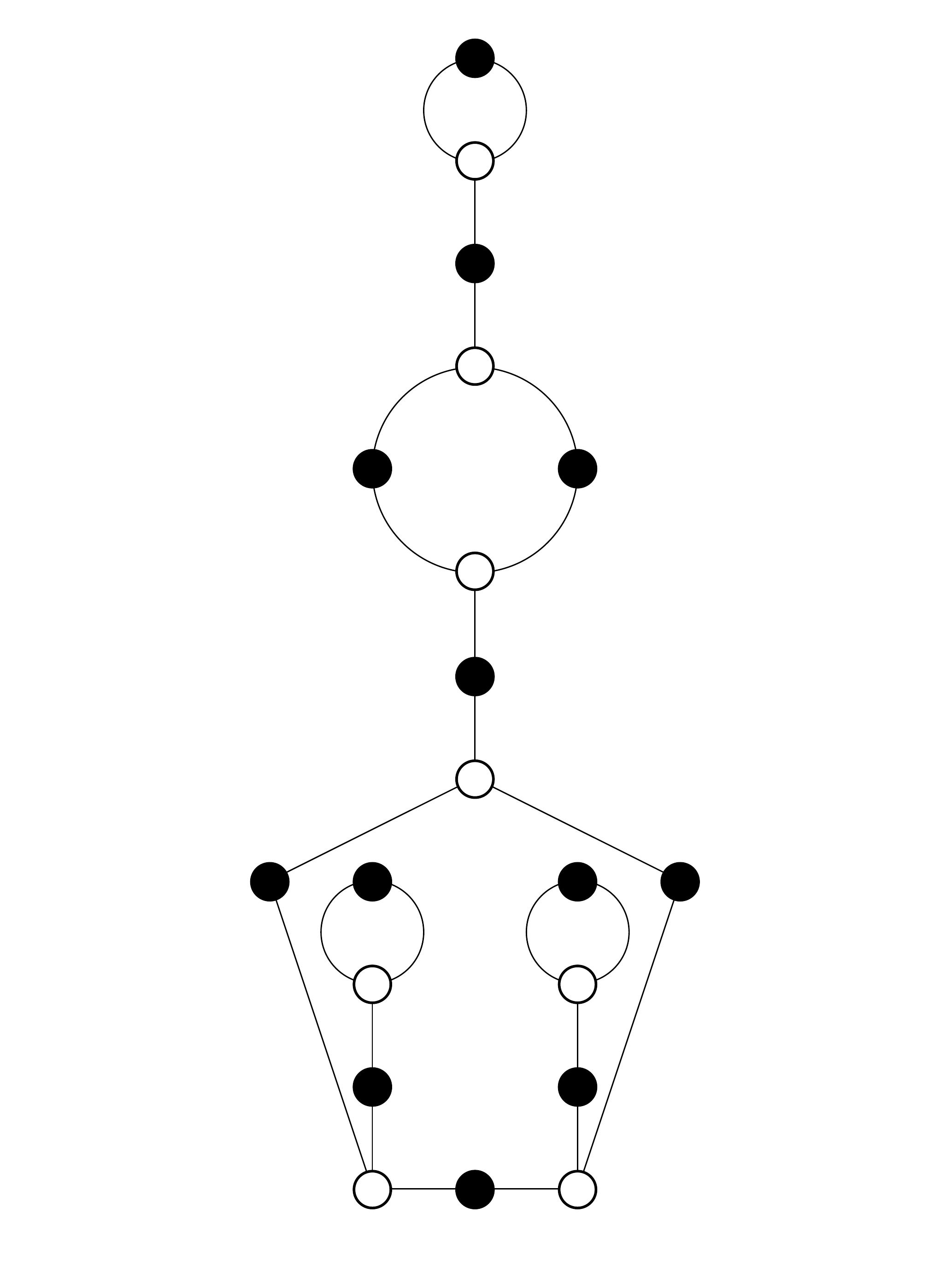}}
\par\end{center}{\scriptsize \par}

\begin{center}
{\scriptsize $10,9,2,1,1,1\;\left(\sqrt{5}\right)$}
\par\end{center}%
\end{minipage}{\scriptsize }%
\begin{minipage}[t]{0.33\textwidth}%
\begin{center}
{\scriptsize \includegraphics[scale=0.15]{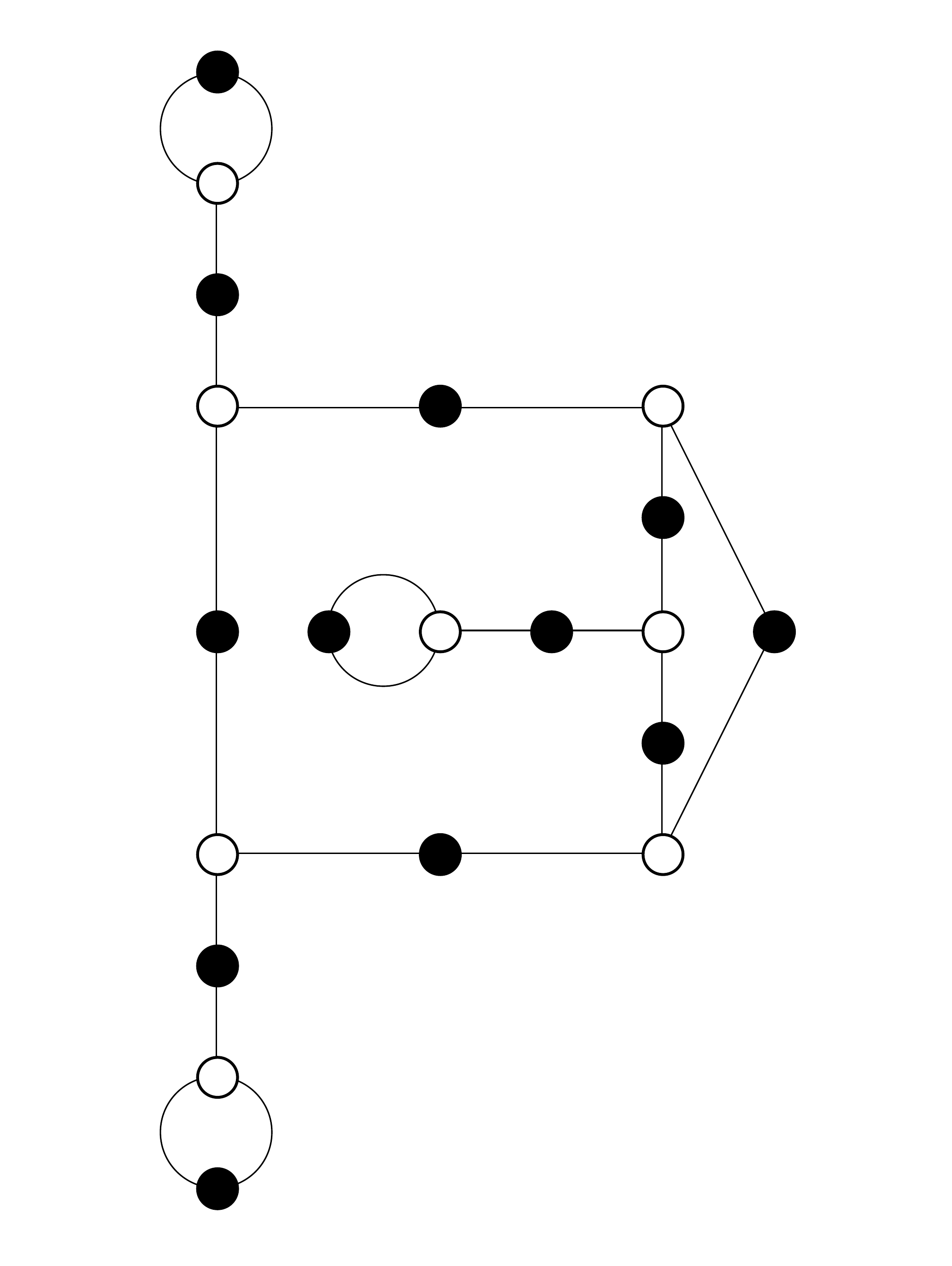}}
\par\end{center}{\scriptsize \par}

\begin{center}
{\scriptsize $10,8,3,1,1,1\;\left(\mathbb{Q}\right)$}
\par\end{center}%
\end{minipage}
\par\end{center}{\scriptsize \par}

\begin{center}
{\scriptsize }%
\begin{minipage}[t]{0.33\textwidth}%
\begin{center}
{\scriptsize \includegraphics[scale=0.15]{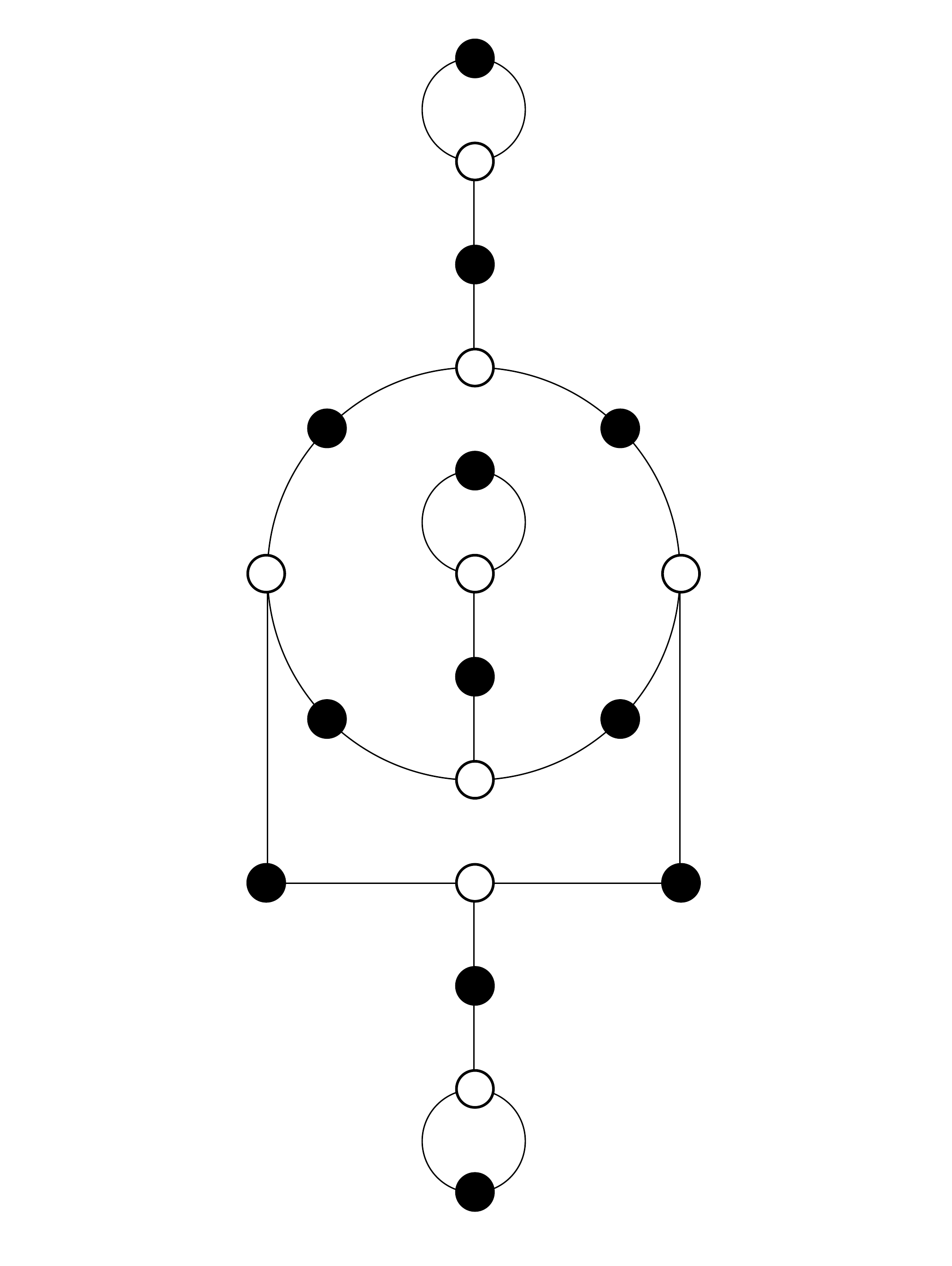}}
\par\end{center}{\scriptsize \par}

\begin{center}
{\scriptsize $10,7,4,1,1,1\;\left(\mathbb{Q}\right)$}
\par\end{center}%
\end{minipage}{\scriptsize }%
\begin{minipage}[t]{0.33\textwidth}%
\begin{center}
{\scriptsize \includegraphics[scale=0.15]{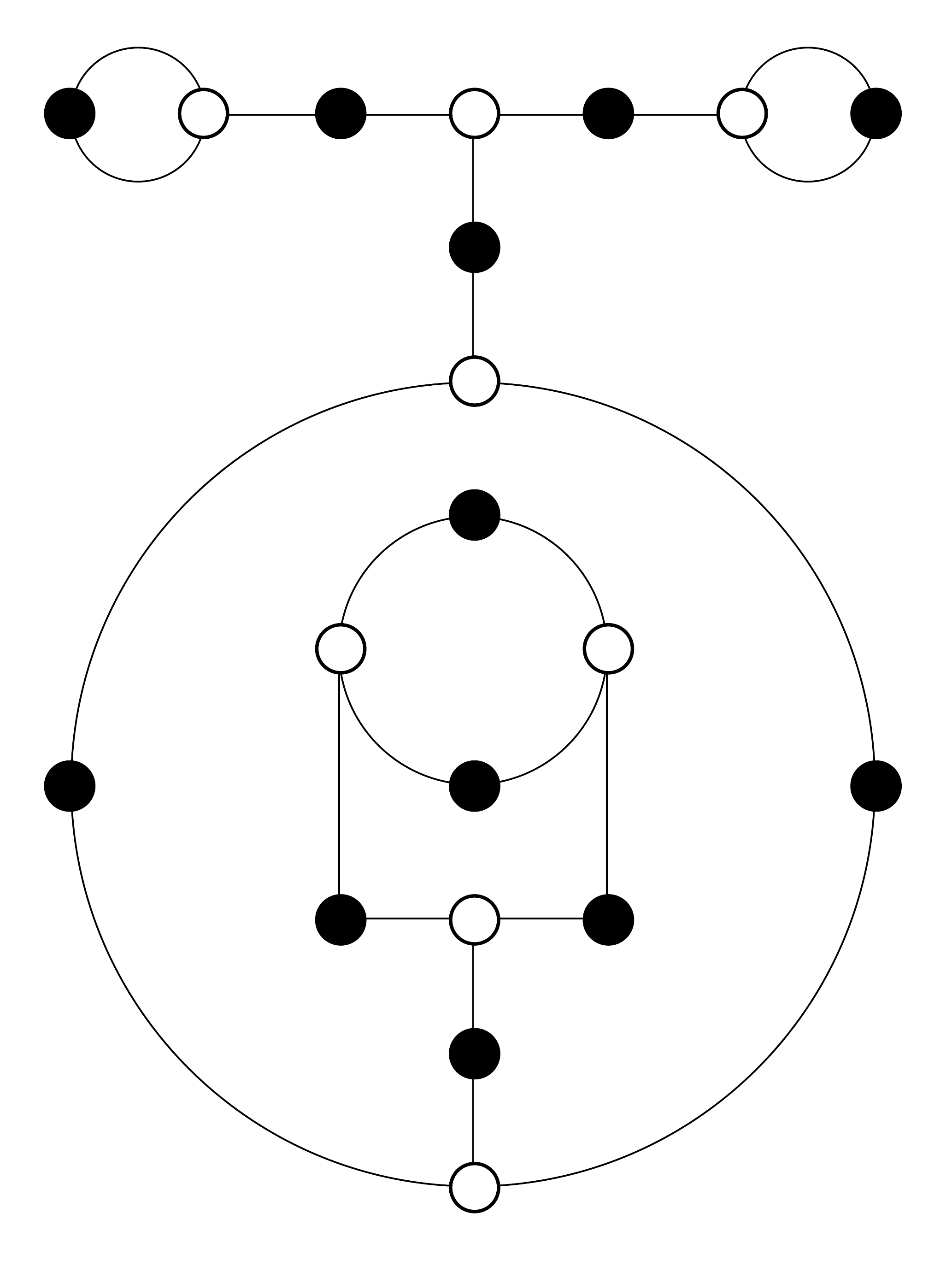}}
\par\end{center}{\scriptsize \par}

\begin{center}
{\scriptsize $10,7,3,2,1,1\;\left(\sqrt{21}\right)$}
\par\end{center}%
\end{minipage}{\scriptsize }%
\begin{minipage}[t]{0.33\textwidth}%
\begin{center}
{\scriptsize \includegraphics[scale=0.15]{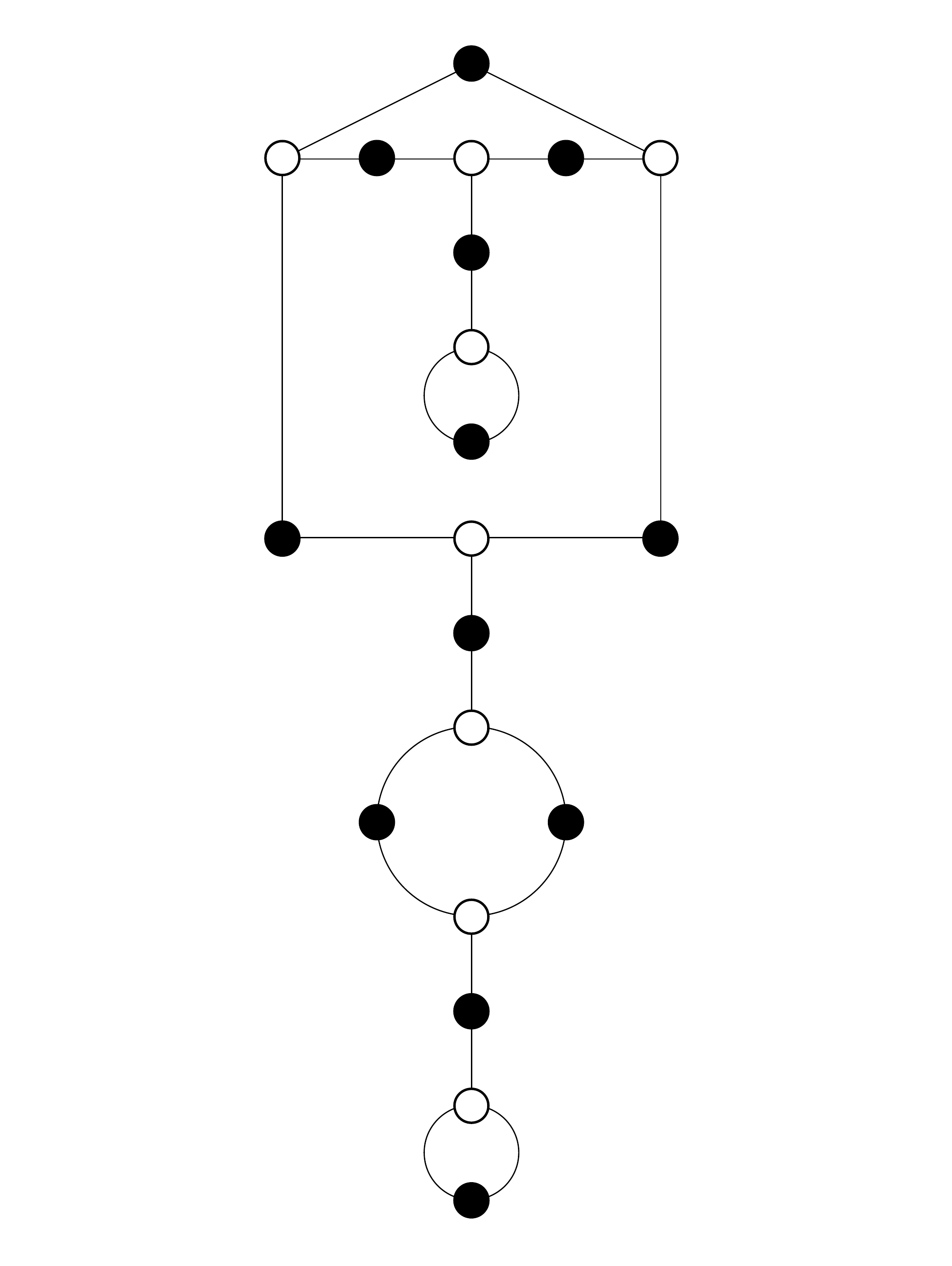}}
\par\end{center}{\scriptsize \par}

\begin{center}
{\scriptsize $10,7,3,2,1,1\;\left(\sqrt{21}\right)$}
\par\end{center}%
\end{minipage}
\par\end{center}{\scriptsize \par}

\begin{center}
{\scriptsize }%
\begin{minipage}[t]{0.33\textwidth}%
\begin{center}
{\scriptsize \includegraphics[scale=0.15]{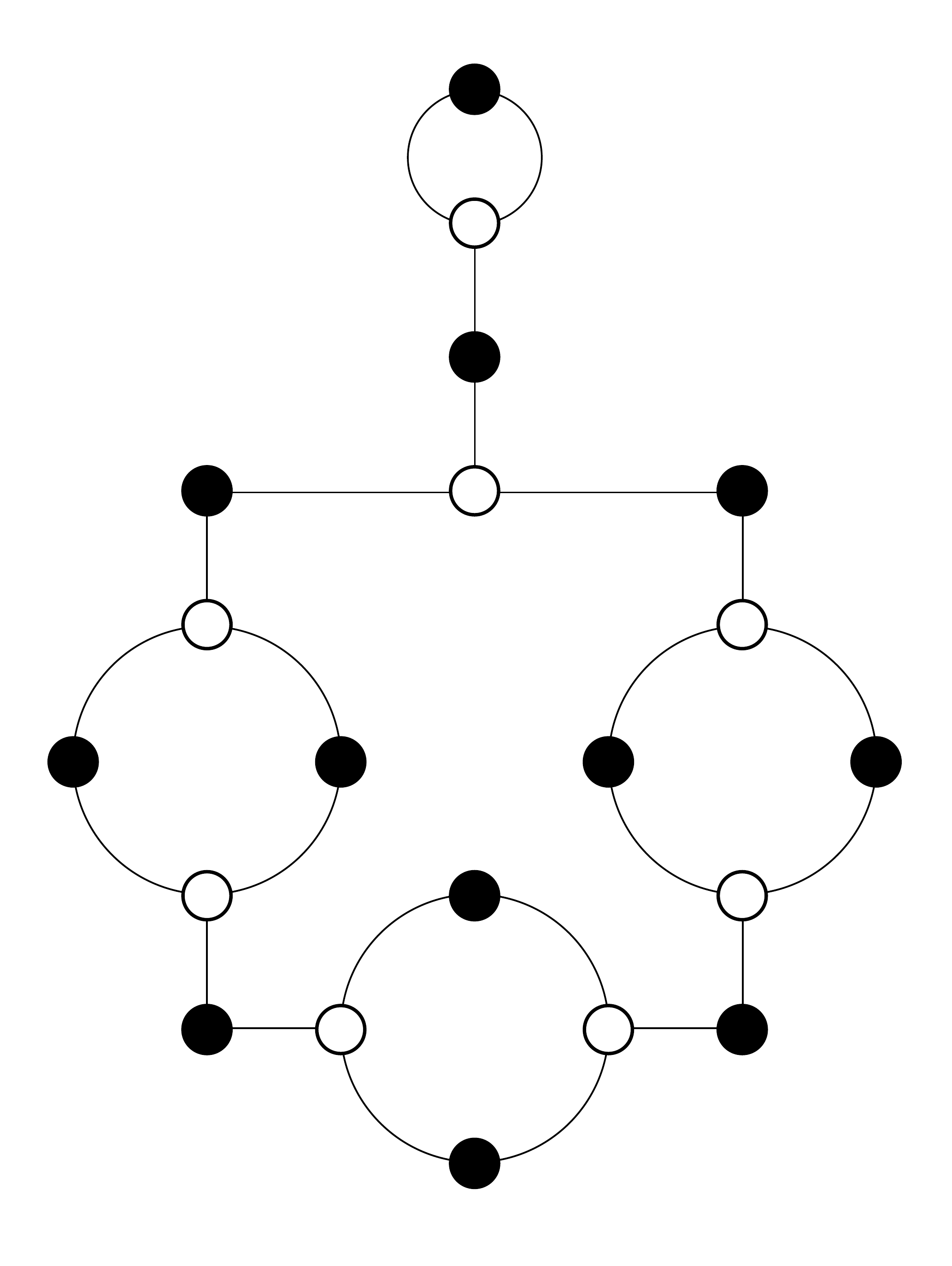}}
\par\end{center}{\scriptsize \par}

\begin{center}
{\scriptsize $10,7,2,2,2,1\;\left(\mathbb{Q}\right)$}
\par\end{center}%
\end{minipage}{\scriptsize }%
\begin{minipage}[t]{0.33\textwidth}%
\begin{center}
{\scriptsize \includegraphics[scale=0.15]{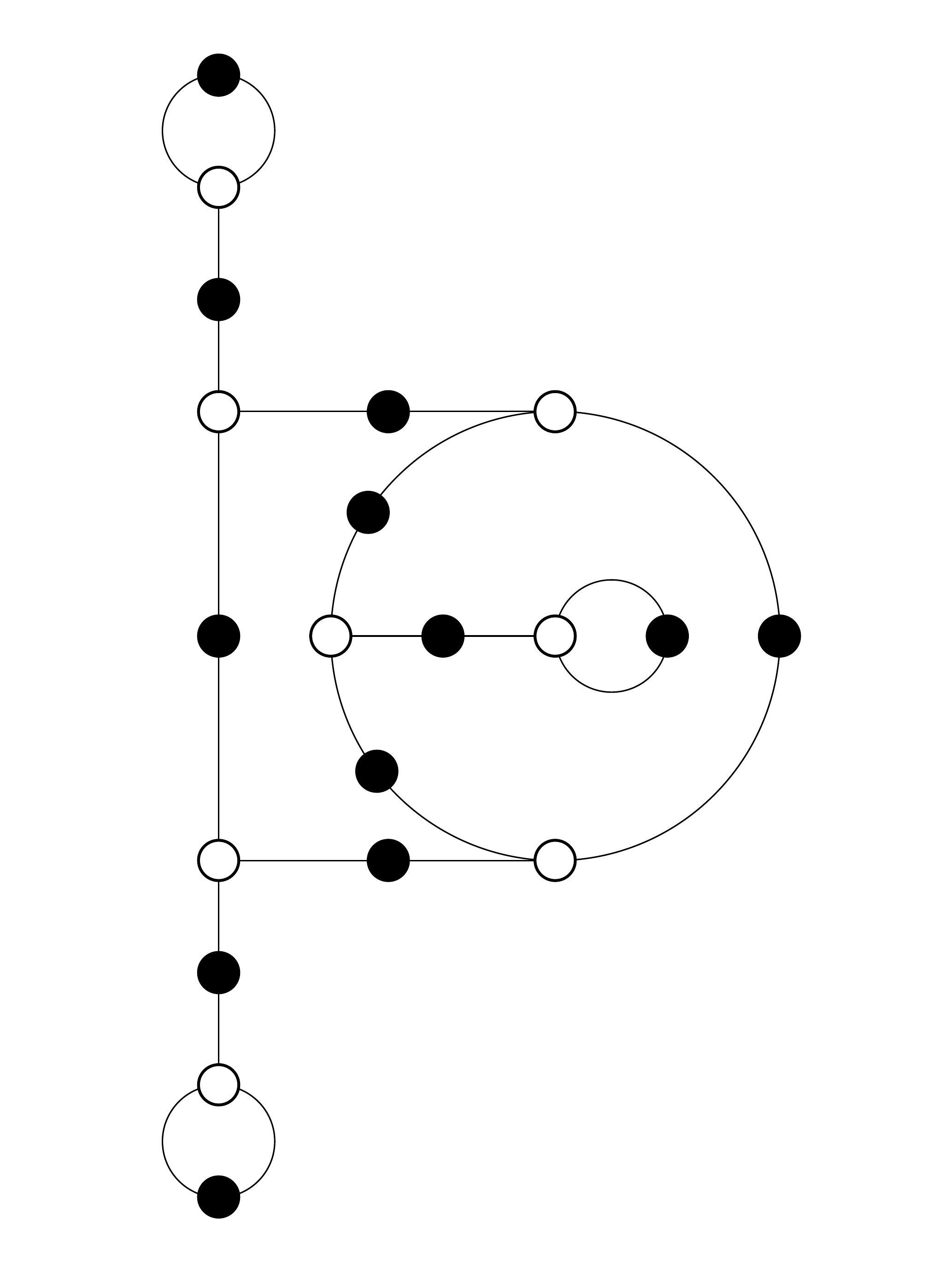}}
\par\end{center}{\scriptsize \par}

\begin{center}
{\scriptsize $10,6,5,1,1,1\;\left(\mathbb{Q}\right)$}
\par\end{center}%
\end{minipage}{\scriptsize }%
\begin{minipage}[t]{0.33\textwidth}%
\begin{center}
{\scriptsize \includegraphics[scale=0.15]{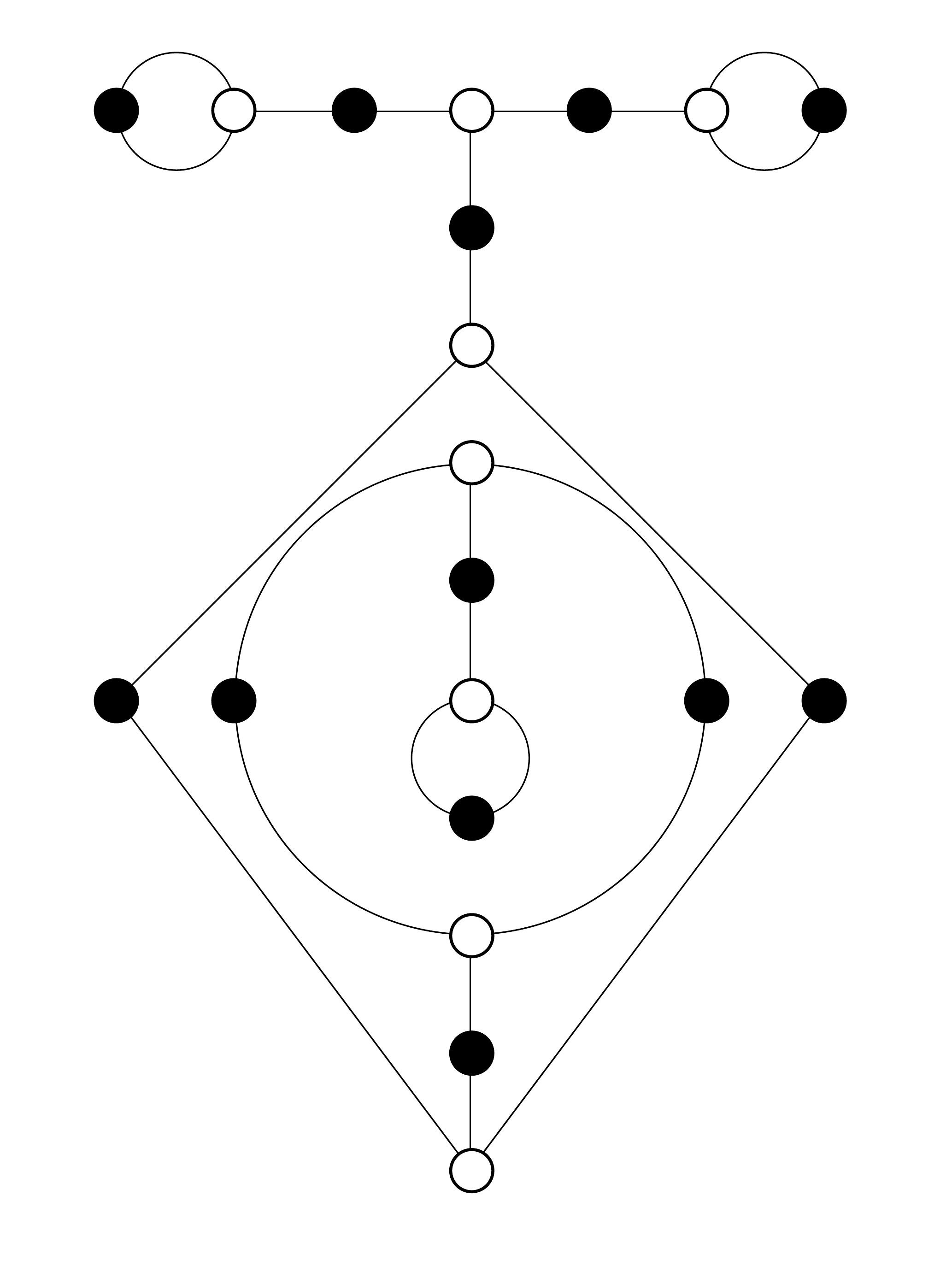}}
\par\end{center}{\scriptsize \par}

\begin{center}
{\scriptsize $10,6,5,1,1,1\;\left(\mathrm{cubic}\right)$}
\par\end{center}%
\end{minipage}
\par\end{center}{\scriptsize \par}

\begin{center}
{\scriptsize }%
\begin{minipage}[t]{0.33\textwidth}%
\begin{center}
{\scriptsize \includegraphics[scale=0.15]{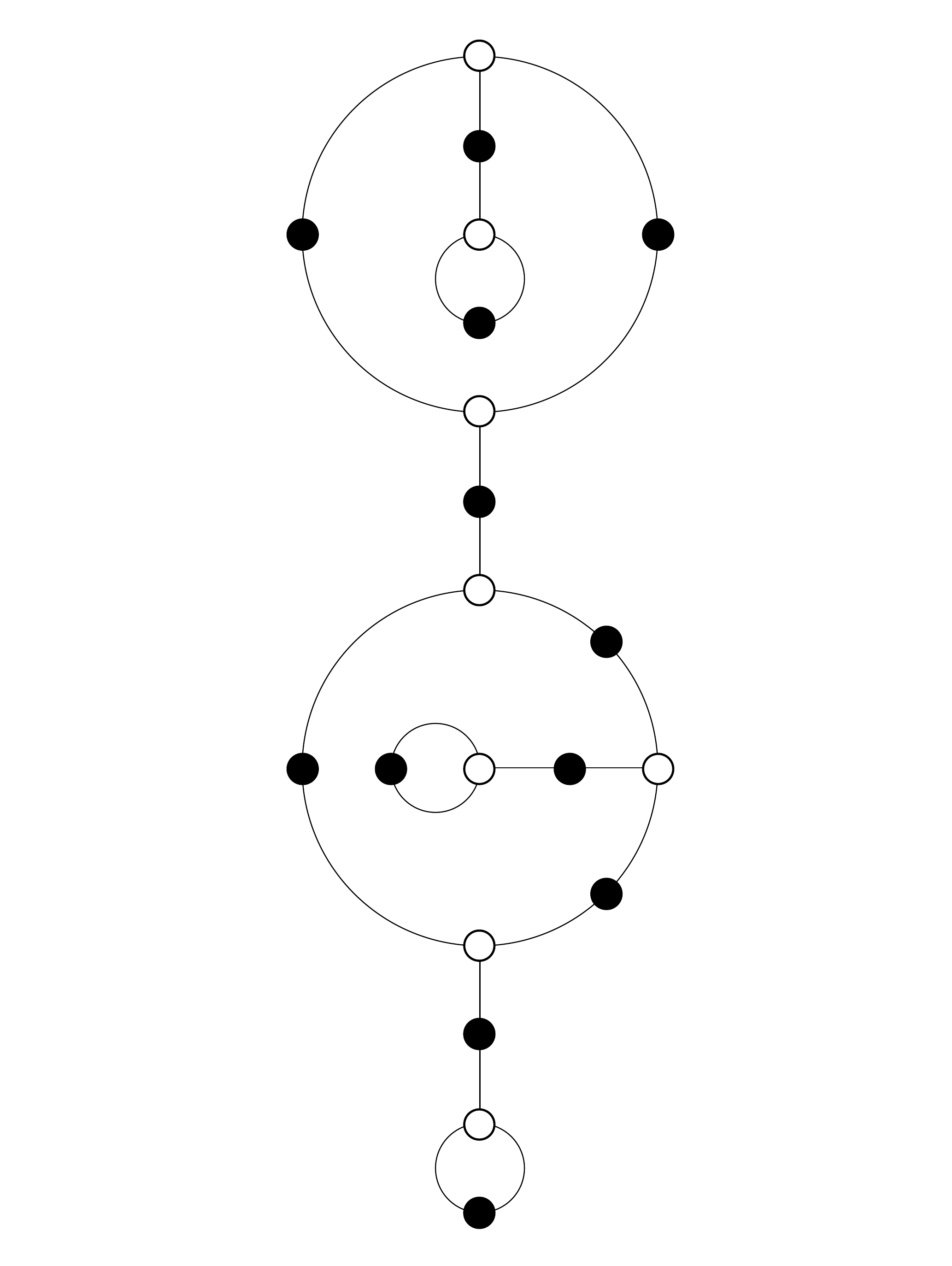}}
\par\end{center}{\scriptsize \par}

\begin{center}
{\scriptsize $10,6,5,1,1,1\;\left(\mathrm{cubic}\right)$}
\par\end{center}%
\end{minipage}{\scriptsize }%
\begin{minipage}[t]{0.33\textwidth}%
\begin{center}
{\scriptsize \includegraphics[scale=0.15]{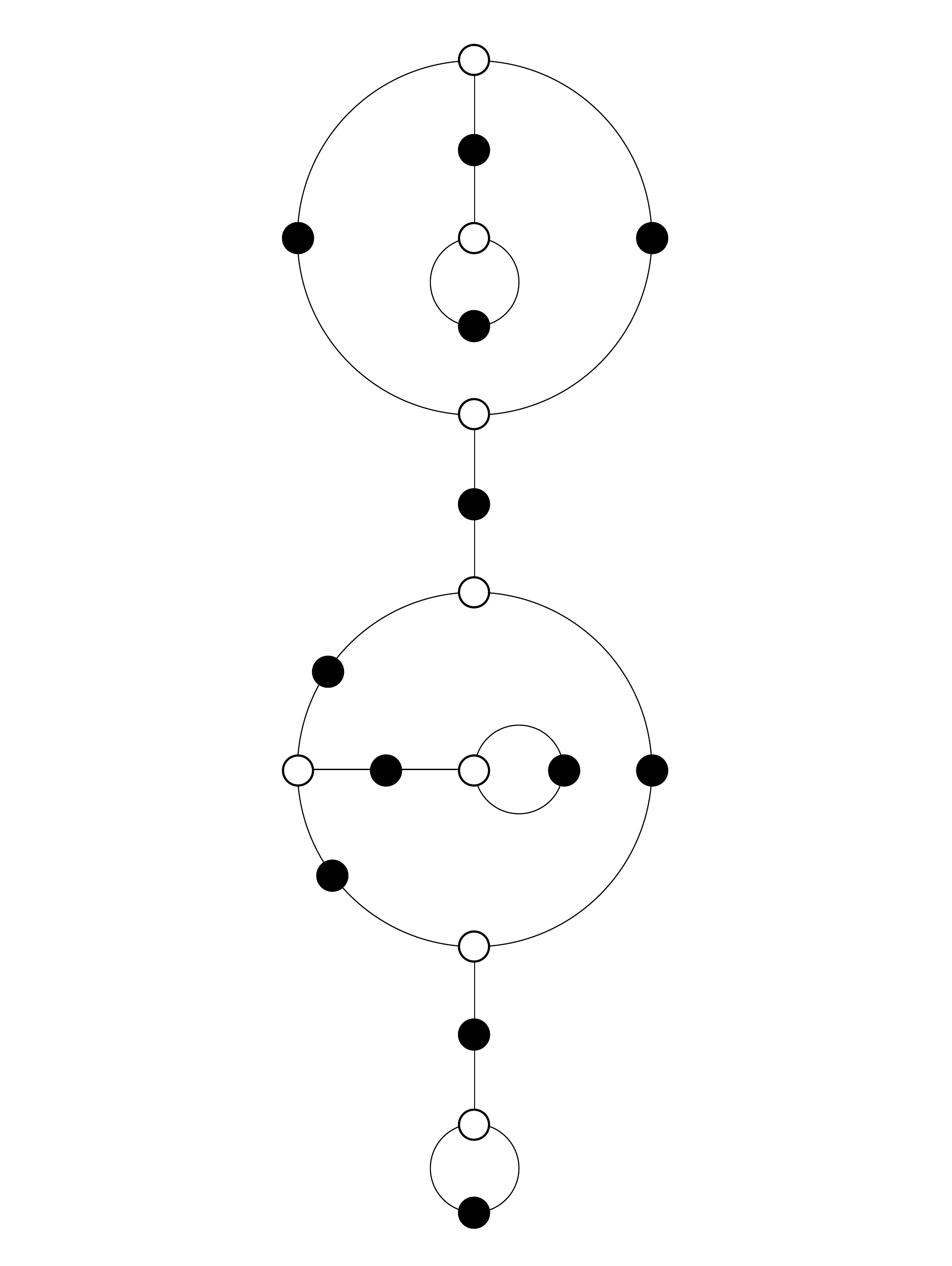}}
\par\end{center}{\scriptsize \par}

\begin{center}
{\scriptsize $10,6,5,1,1,1\;\left(\mathrm{cubic}\right)$}
\par\end{center}%
\end{minipage}{\scriptsize }%
\begin{minipage}[t]{0.33\textwidth}%
\begin{center}
{\scriptsize \includegraphics[scale=0.15]{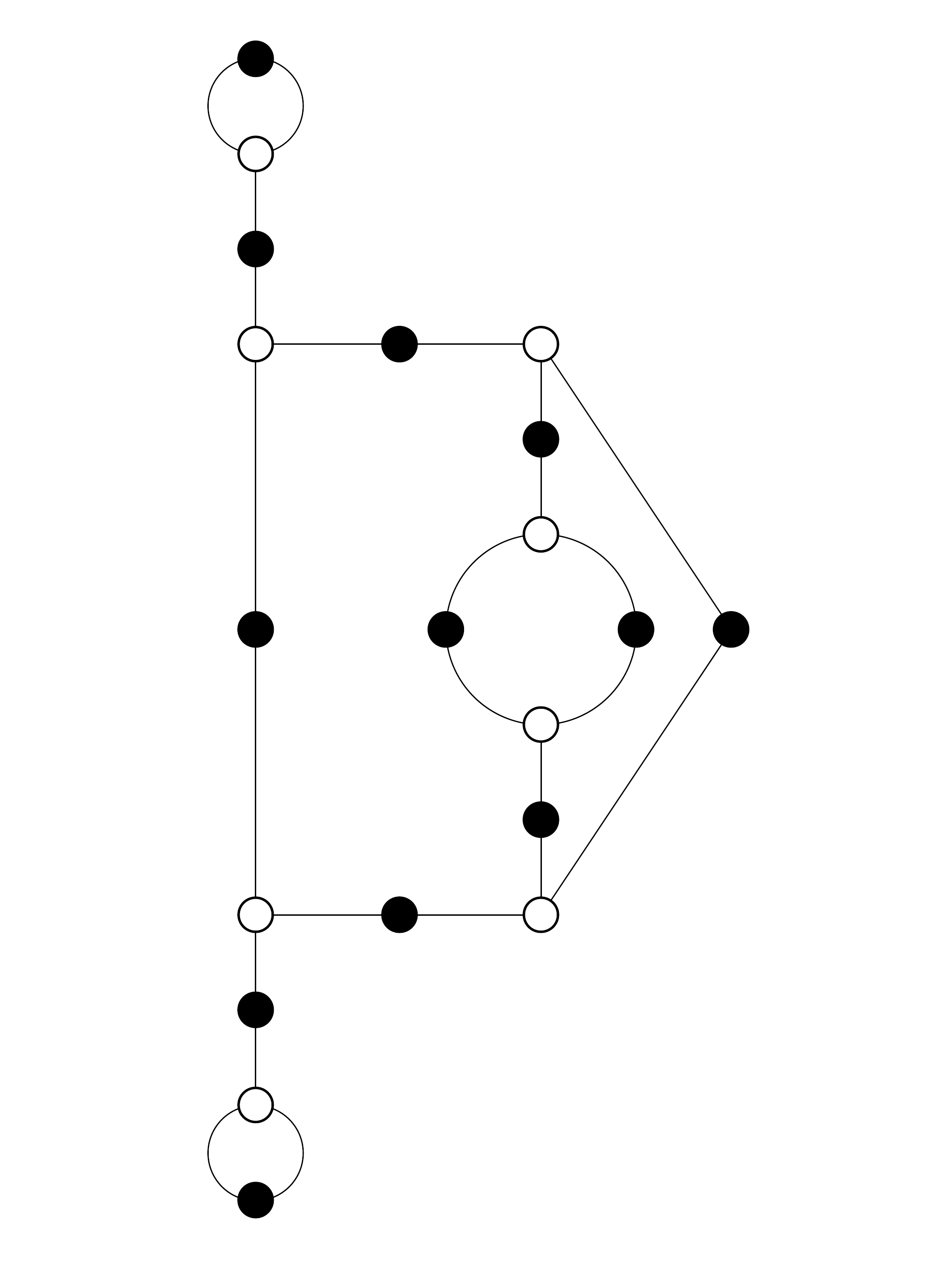}}
\par\end{center}{\scriptsize \par}

\begin{center}
{\scriptsize $10,6,4,2,1,1\;\left(\mathbb{Q}\right)$}
\par\end{center}%
\end{minipage}
\par\end{center}{\scriptsize \par}

\begin{center}
{\scriptsize }%
\begin{minipage}[t]{0.33\textwidth}%
\begin{center}
{\scriptsize \includegraphics[scale=0.15]{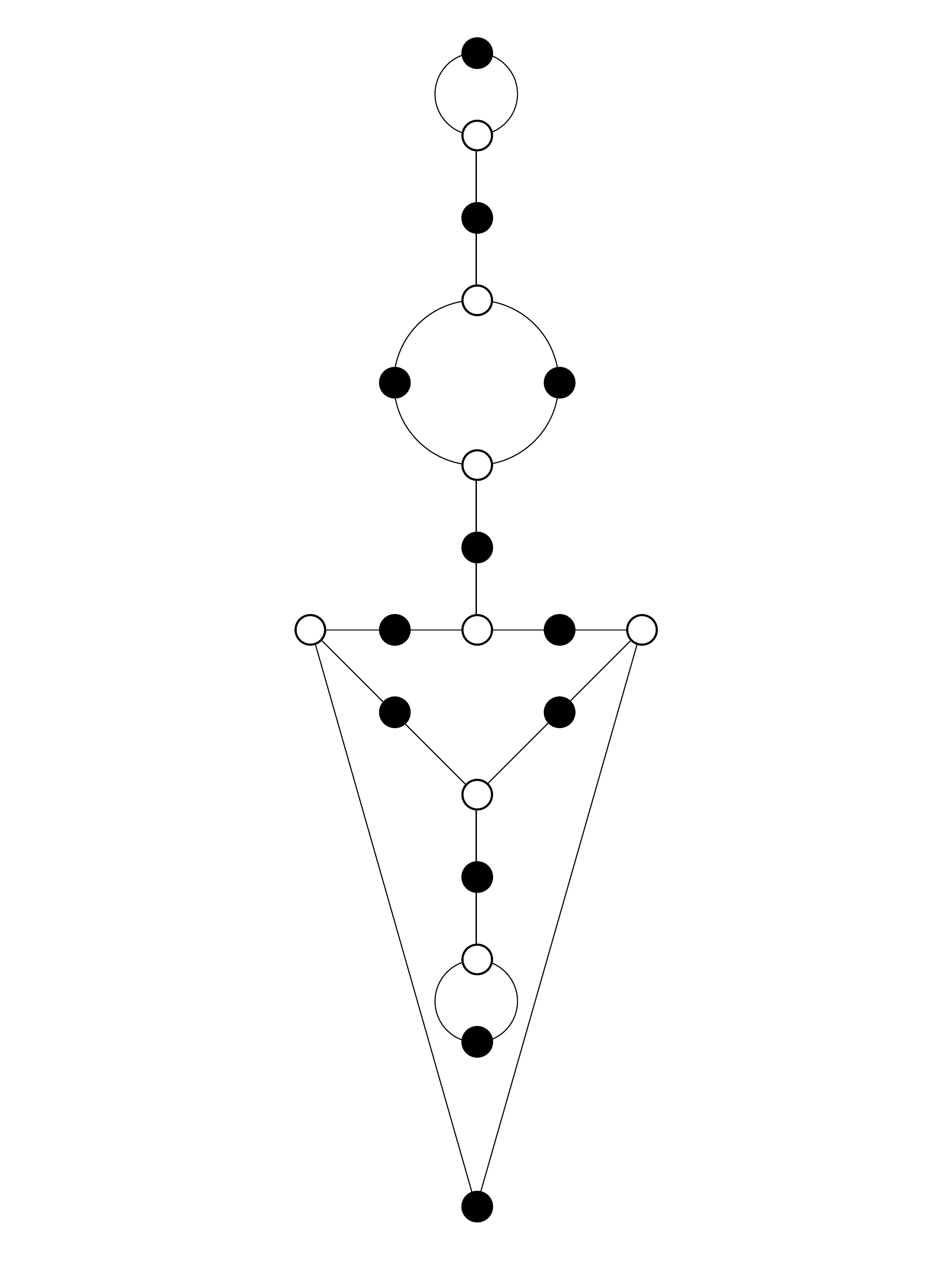}}
\par\end{center}{\scriptsize \par}

\begin{center}
{\scriptsize $10,6,4,2,1,1\;\left(\mathbb{Q}\right)$}
\par\end{center}%
\end{minipage}{\scriptsize }%
\begin{minipage}[t]{0.33\textwidth}%
\begin{center}
{\scriptsize \includegraphics[scale=0.15]{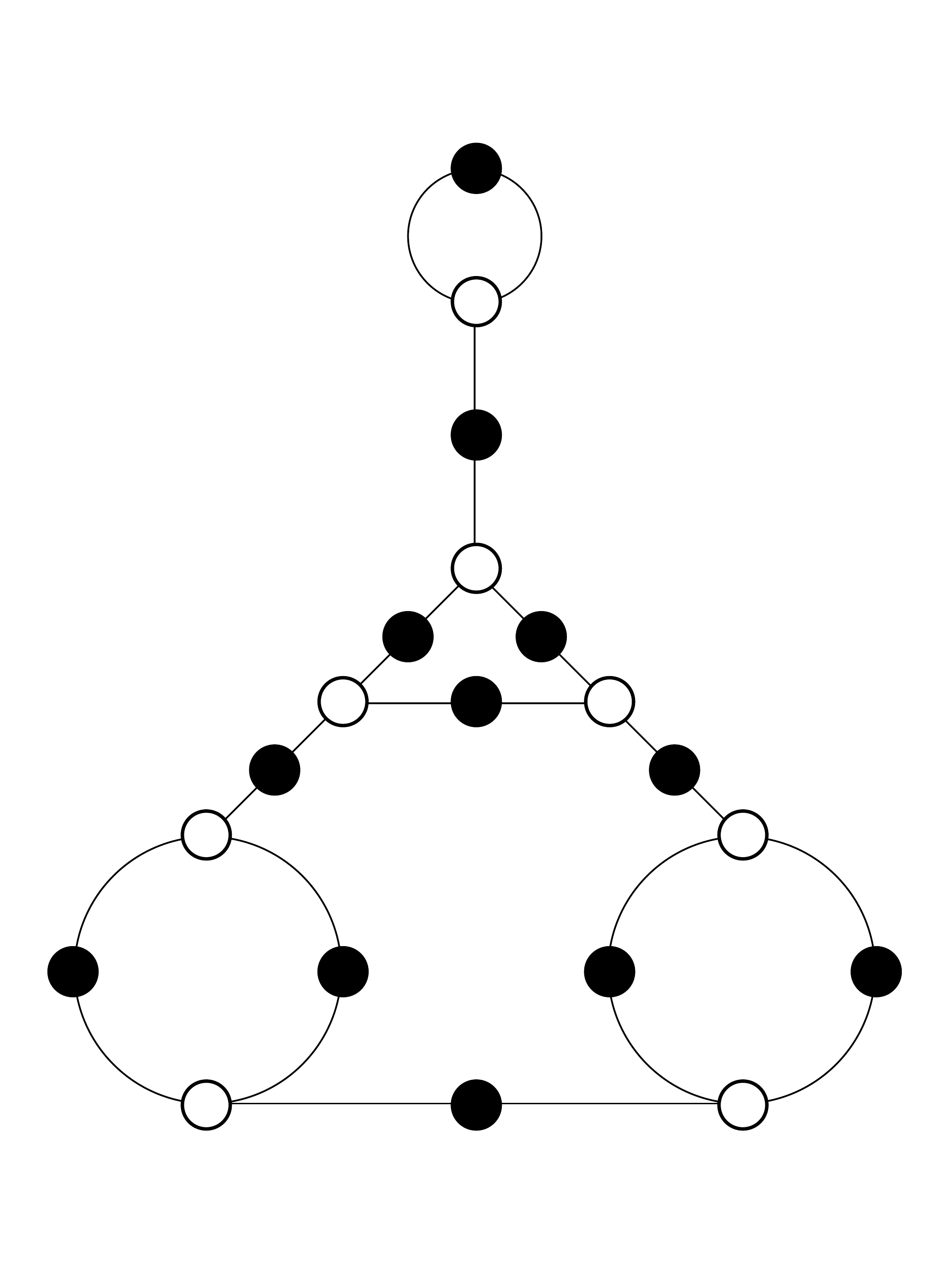}}
\par\end{center}{\scriptsize \par}

\begin{center}
{\scriptsize $10,6,3,2,2,1\;\left(\mathbb{Q}\right)$}
\par\end{center}%
\end{minipage}{\scriptsize }%
\begin{minipage}[t]{0.33\textwidth}%
\begin{center}
{\scriptsize \includegraphics[scale=0.15]{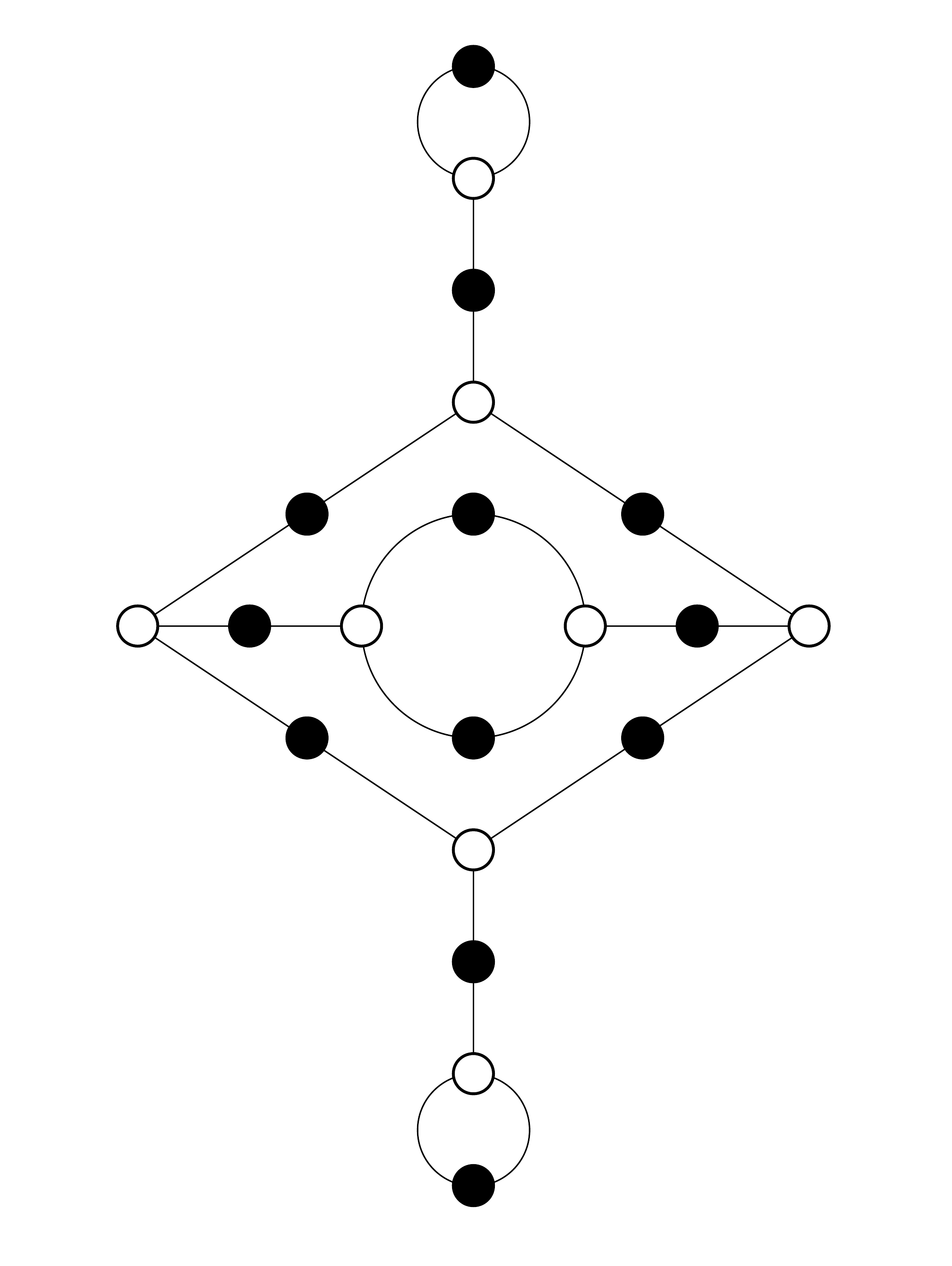}}
\par\end{center}{\scriptsize \par}

\begin{center}
{\scriptsize $10,5,5,2,1,1\;\left(\sqrt{5}\right)$}
\par\end{center}%
\end{minipage}
\par\end{center}{\scriptsize \par}

\begin{center}
{\scriptsize }%
\begin{minipage}[t]{0.33\textwidth}%
\begin{center}
{\scriptsize \includegraphics[scale=0.15]{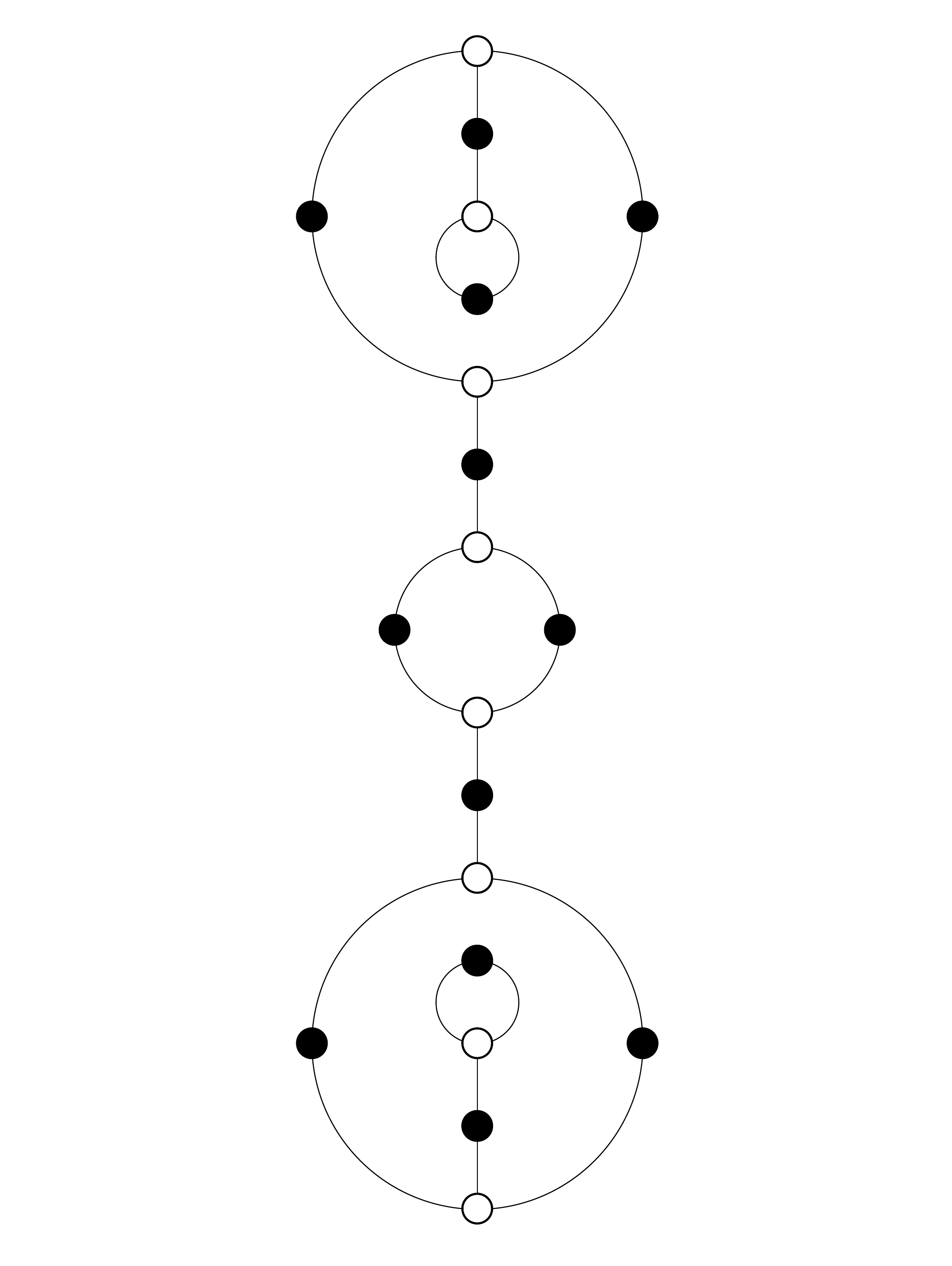}}
\par\end{center}{\scriptsize \par}

\begin{center}
{\scriptsize $10,5,5,2,1,1\;\left(\sqrt{5}\right)$}
\par\end{center}%
\end{minipage}{\scriptsize }%
\begin{minipage}[t]{0.33\textwidth}%
\begin{center}
{\scriptsize \includegraphics[scale=0.15]{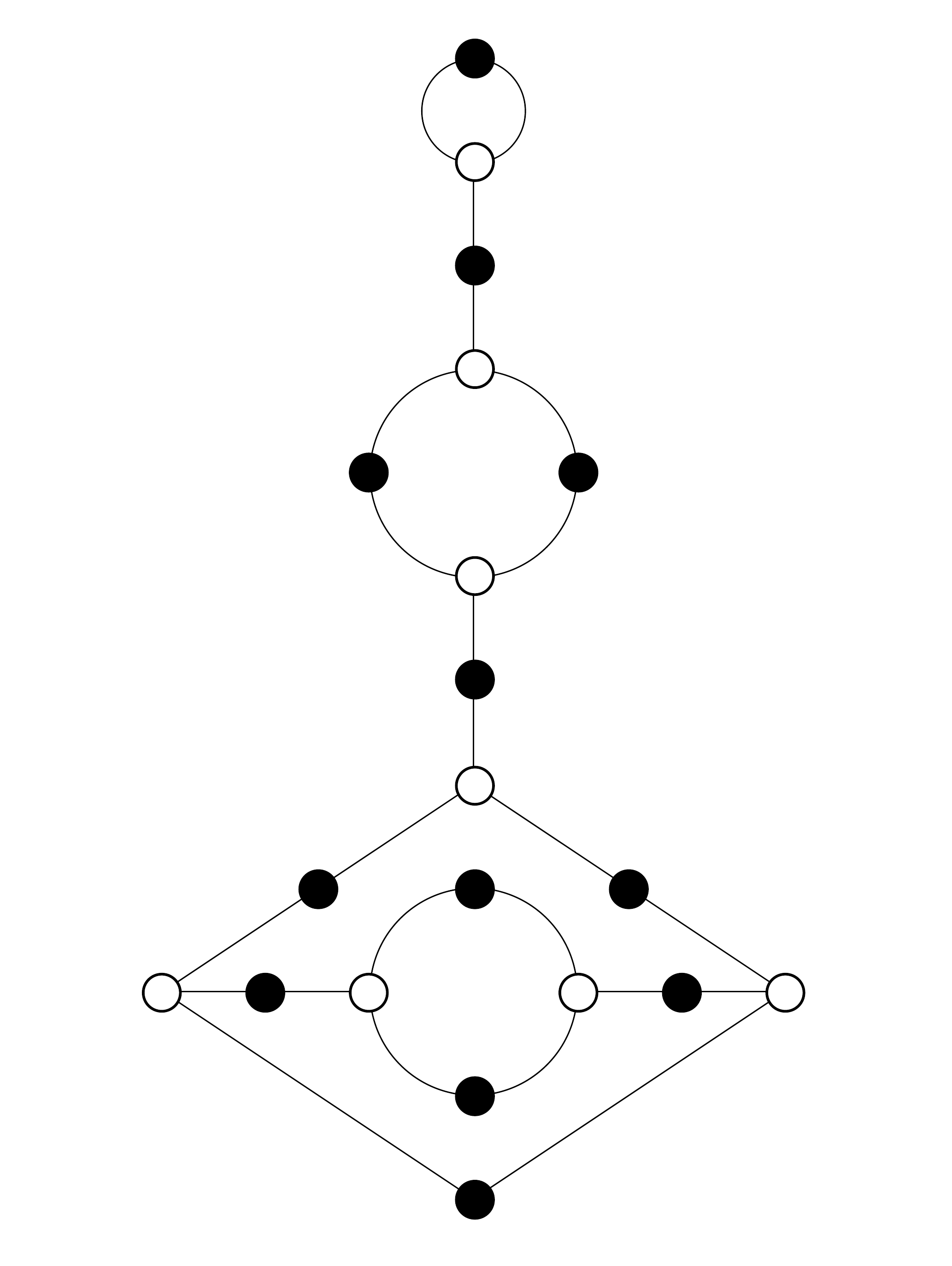}}
\par\end{center}{\scriptsize \par}

\begin{center}
{\scriptsize $10,5,4,2,2,1\;\left(\mathrm{cubic}\right)$}
\par\end{center}%
\end{minipage}{\scriptsize }%
\begin{minipage}[t]{0.33\textwidth}%
\begin{center}
{\scriptsize \includegraphics[scale=0.15]{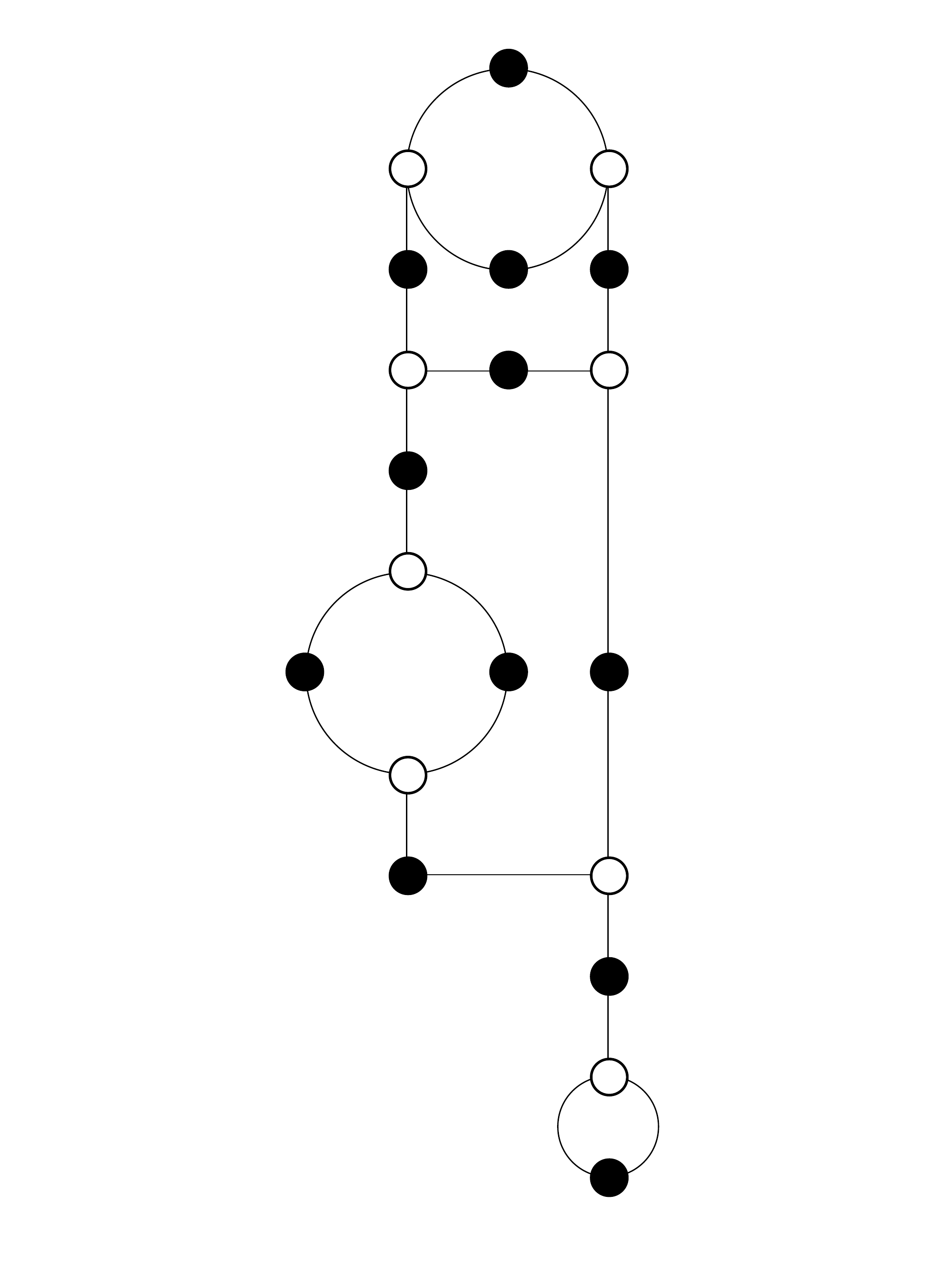}}
\par\end{center}{\scriptsize \par}

\begin{center}
{\scriptsize $10,5,4,2,2,1\;\left(\mathrm{cubic}\right)$}
\par\end{center}%
\end{minipage}
\par\end{center}{\scriptsize \par}

\begin{center}
{\scriptsize }%
\begin{minipage}[t]{0.33\textwidth}%
\begin{center}
{\scriptsize \includegraphics[scale=0.15]{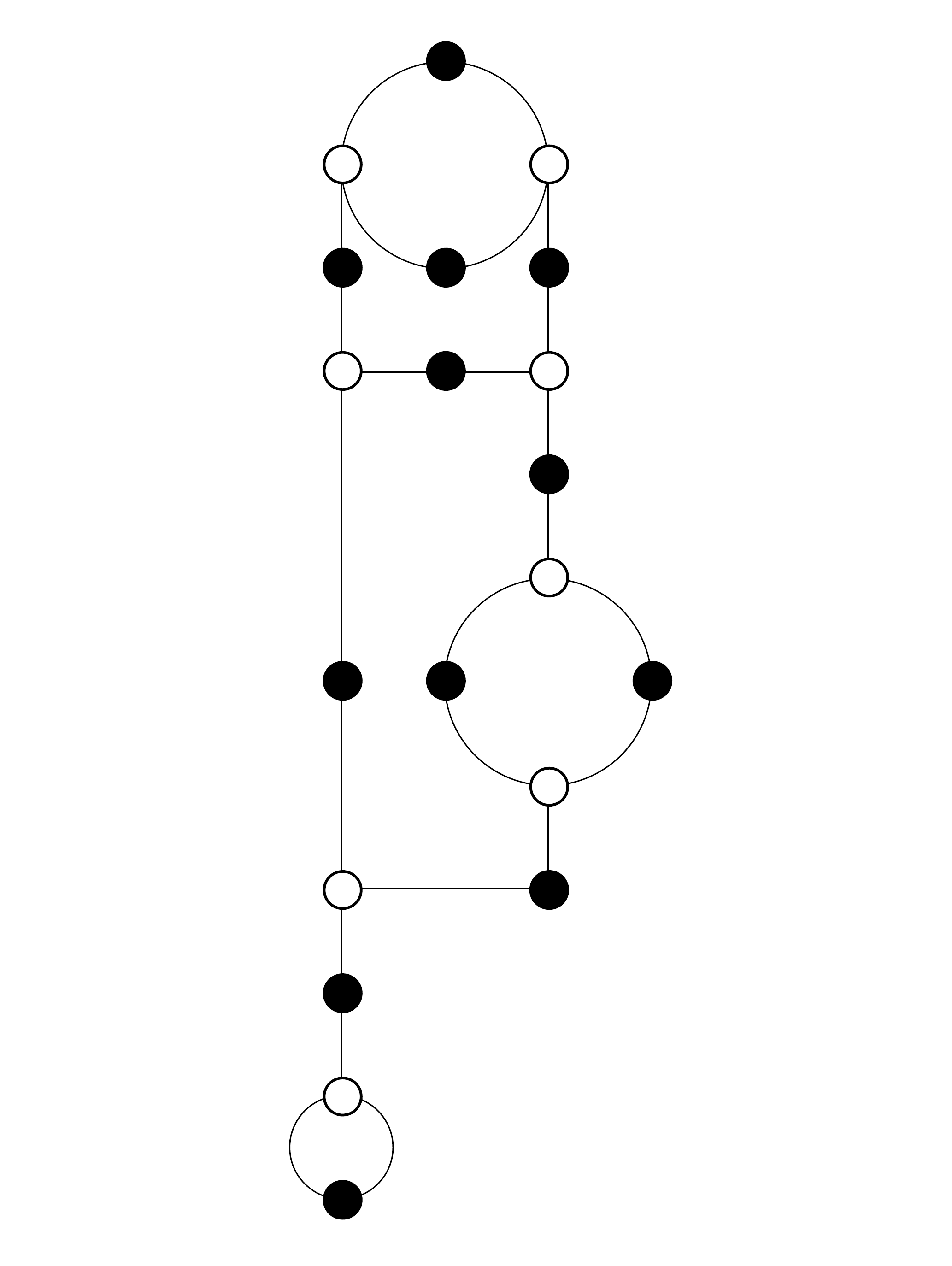}}
\par\end{center}{\scriptsize \par}

\begin{center}
{\scriptsize $10,5,4,2,2,1\;\left(\mathrm{cubic}\right)$}
\par\end{center}%
\end{minipage}{\scriptsize }%
\begin{minipage}[t]{0.33\textwidth}%
\begin{center}
{\scriptsize \includegraphics[scale=0.15]{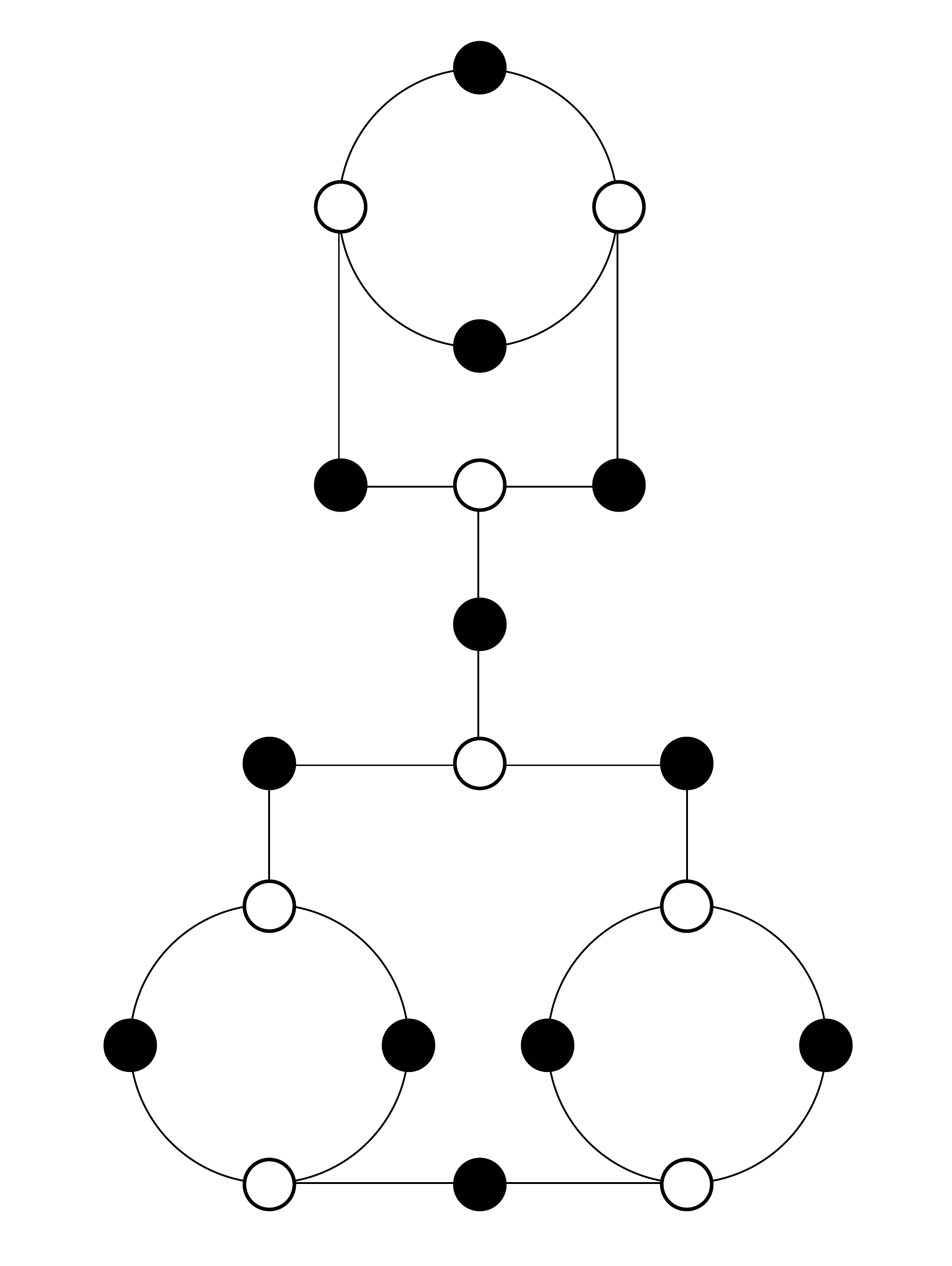}}
\par\end{center}{\scriptsize \par}

\begin{center}
{\scriptsize $10,5,3,2,2,2\;\left(\mathbb{Q}\right)$}
\par\end{center}%
\end{minipage}{\scriptsize }%
\begin{minipage}[t]{0.33\textwidth}%
\begin{center}
{\scriptsize \includegraphics[scale=0.15]{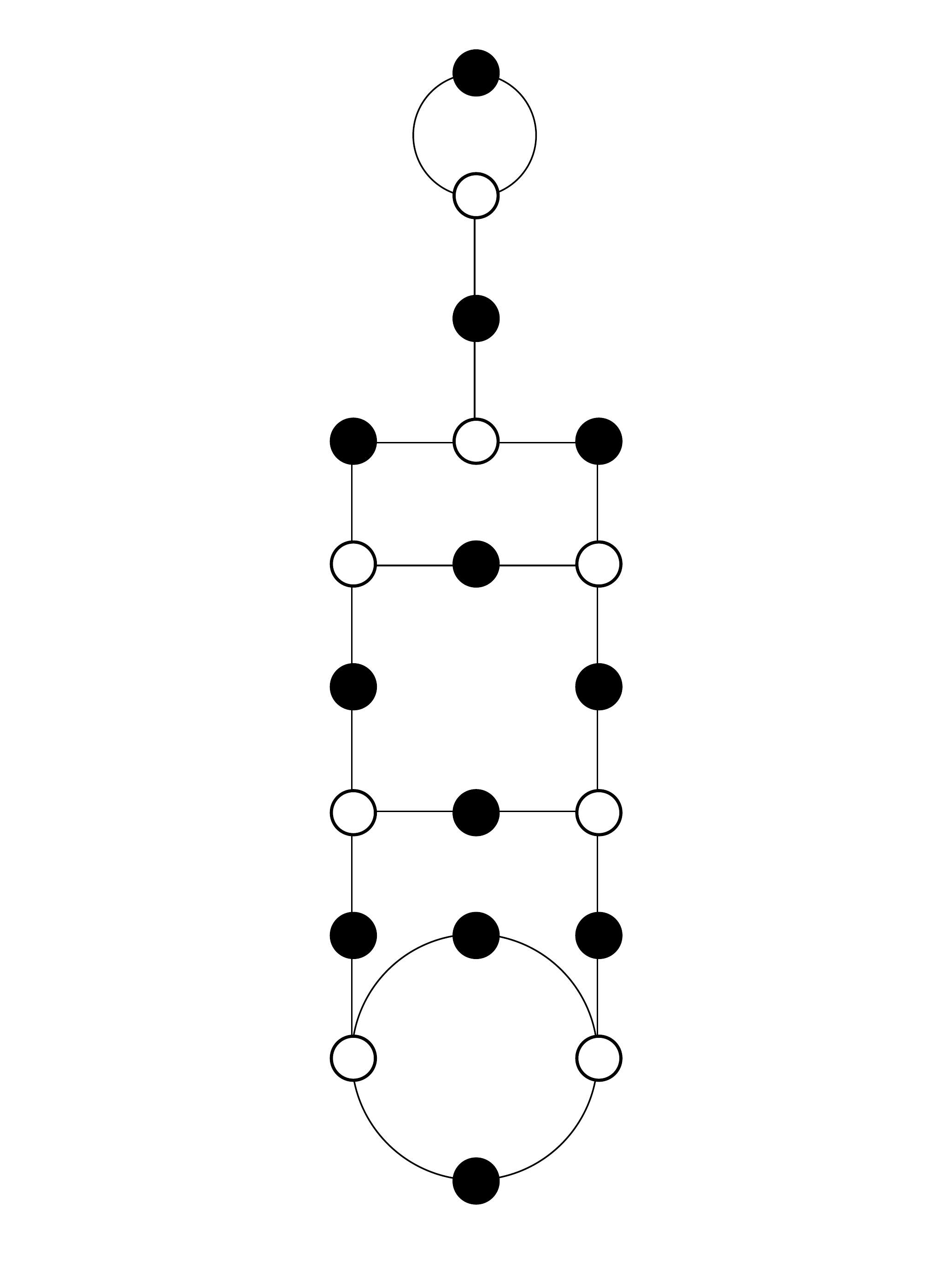}}
\par\end{center}{\scriptsize \par}

\begin{center}
{\scriptsize $10,4,4,3,2,1\;\left(\mathbb{Q}\right)$}
\par\end{center}%
\end{minipage}
\par\end{center}{\scriptsize \par}

\begin{center}
{\scriptsize }%
\begin{minipage}[t]{0.33\textwidth}%
\begin{center}
{\scriptsize \includegraphics[scale=0.15]{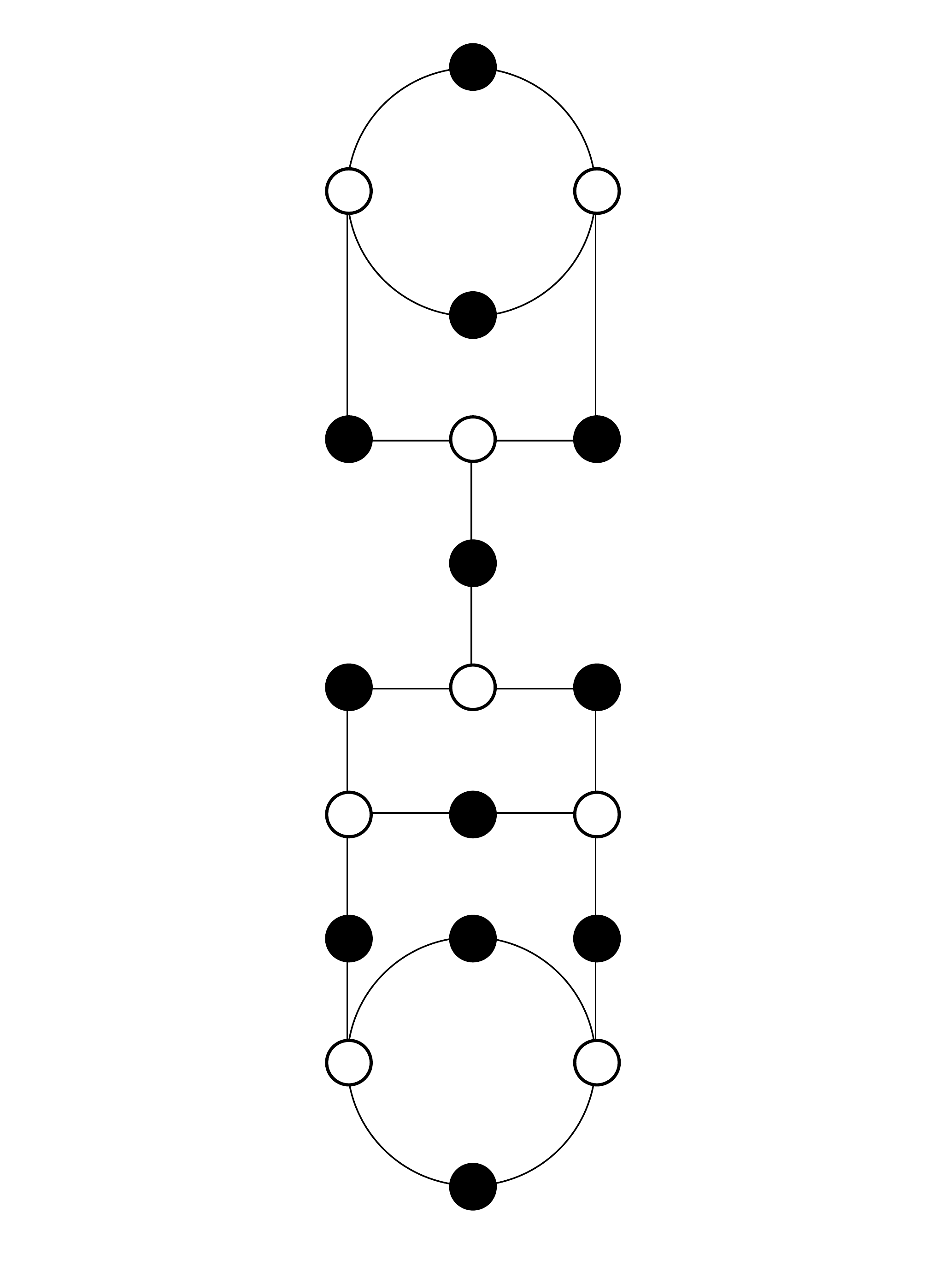}}
\par\end{center}{\scriptsize \par}

\begin{center}
{\scriptsize $10,4,3,3,2,2\;\left(\mathbb{Q}\right)$}
\par\end{center}%
\end{minipage}{\scriptsize }%
\begin{minipage}[t]{0.33\textwidth}%
\begin{center}
{\scriptsize \includegraphics[scale=0.15]{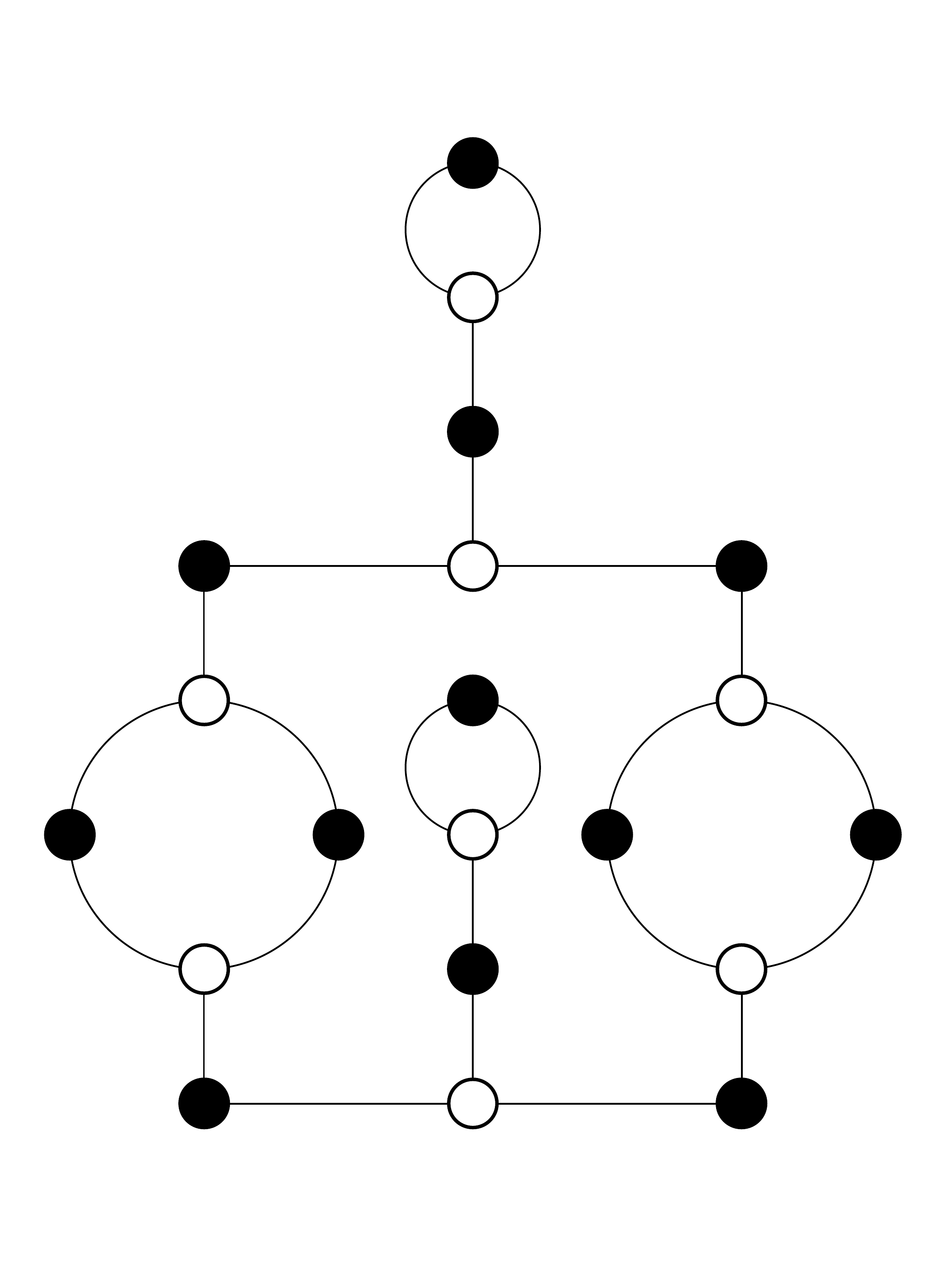}}
\par\end{center}{\scriptsize \par}

\begin{center}
{\scriptsize $9,9,2,2,1,1\;\left(\mathbb{Q}\right)$}
\par\end{center}%
\end{minipage}{\scriptsize }%
\begin{minipage}[t]{0.33\textwidth}%
\begin{center}
{\scriptsize \includegraphics[scale=0.15]{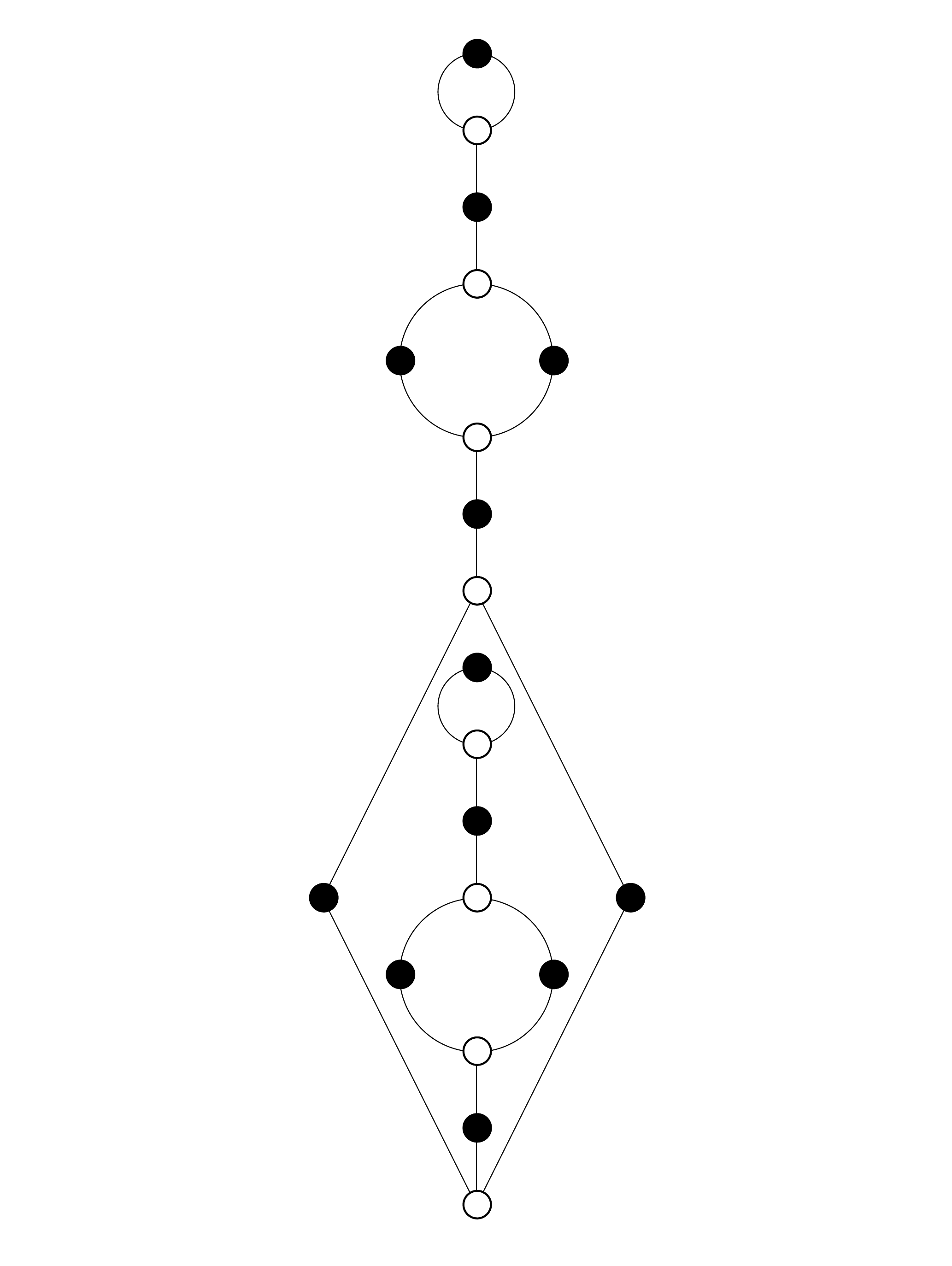}}
\par\end{center}{\scriptsize \par}

\begin{center}
{\scriptsize $9,9,2,2,1,1\;\left(\mathrm{cubic}\right)$}
\par\end{center}%
\end{minipage}
\par\end{center}{\scriptsize \par}

\begin{center}
{\scriptsize }%
\begin{minipage}[t]{0.33\textwidth}%
\begin{center}
{\scriptsize \includegraphics[scale=0.15]{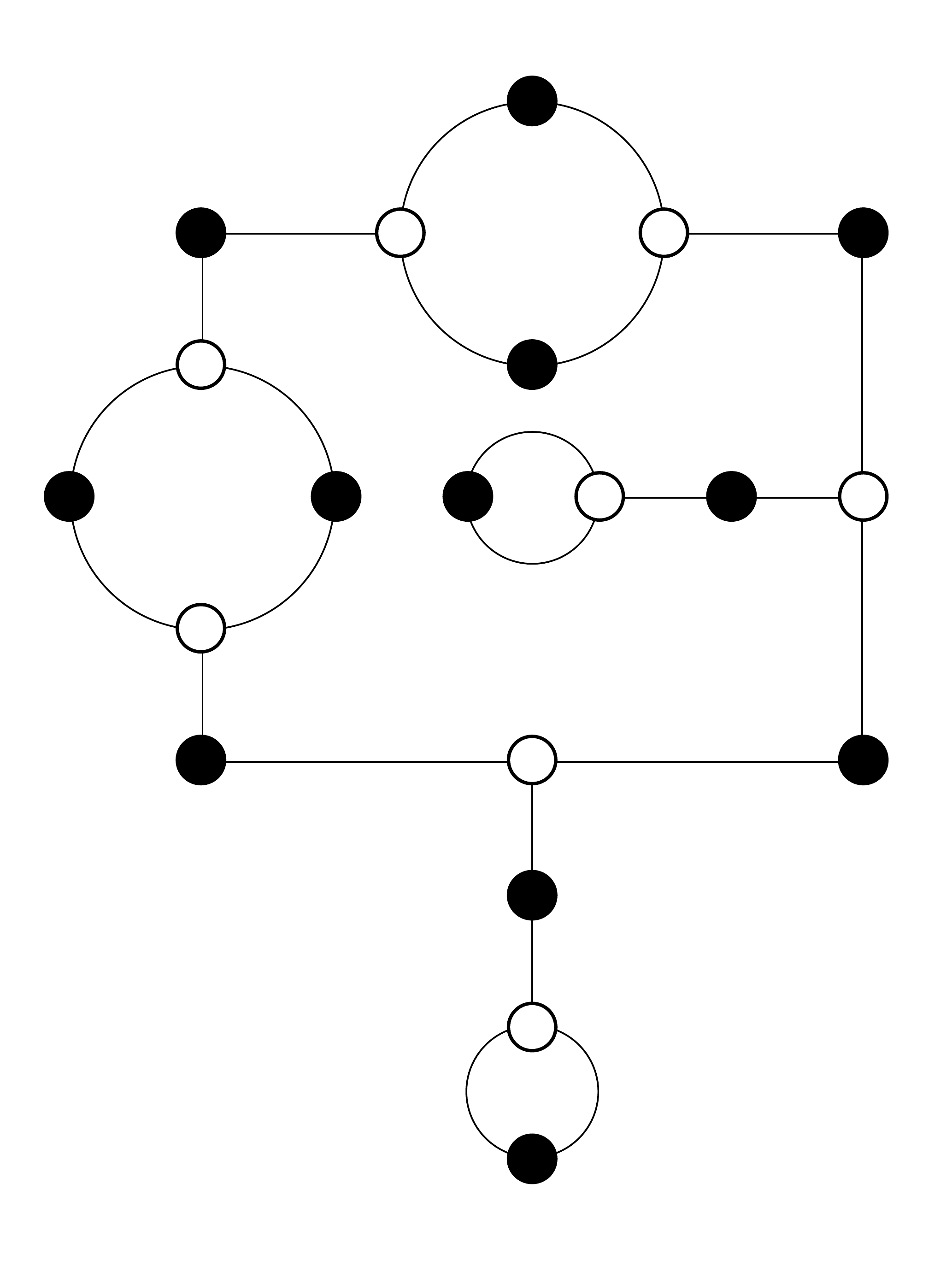}}
\par\end{center}{\scriptsize \par}

\begin{center}
{\scriptsize $9,9,2,2,1,1\;\left(\mathrm{cubic}\right)$}
\par\end{center}%
\end{minipage}{\scriptsize }%
\begin{minipage}[t]{0.33\textwidth}%
\begin{center}
{\scriptsize \includegraphics[scale=0.15]{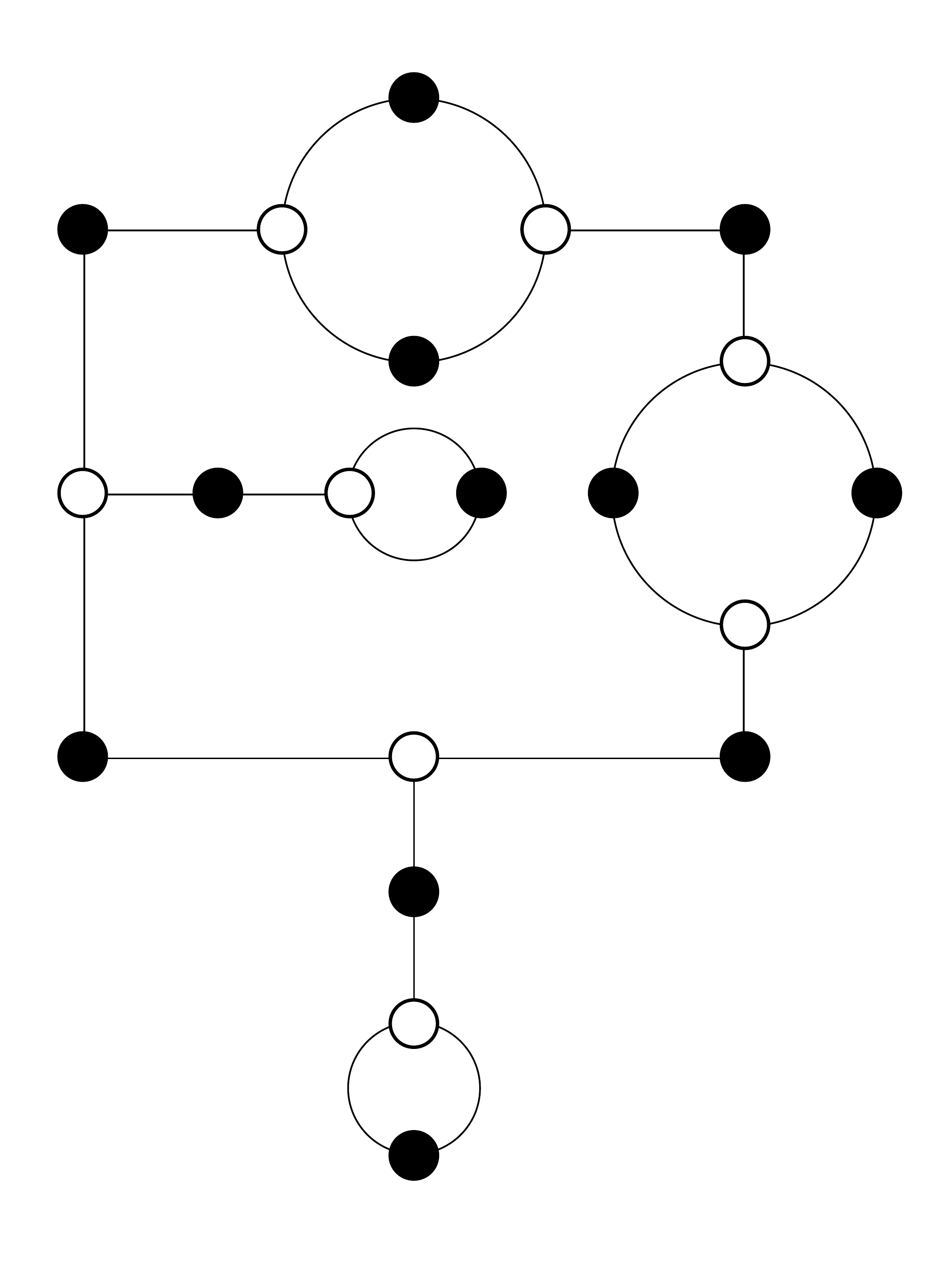}}
\par\end{center}{\scriptsize \par}

\begin{center}
{\scriptsize $9,9,2,2,1,1\;\left(\mathrm{cubic}\right)$}
\par\end{center}%
\end{minipage}{\scriptsize }%
\begin{minipage}[t]{0.33\textwidth}%
\begin{center}
{\scriptsize \includegraphics[scale=0.15]{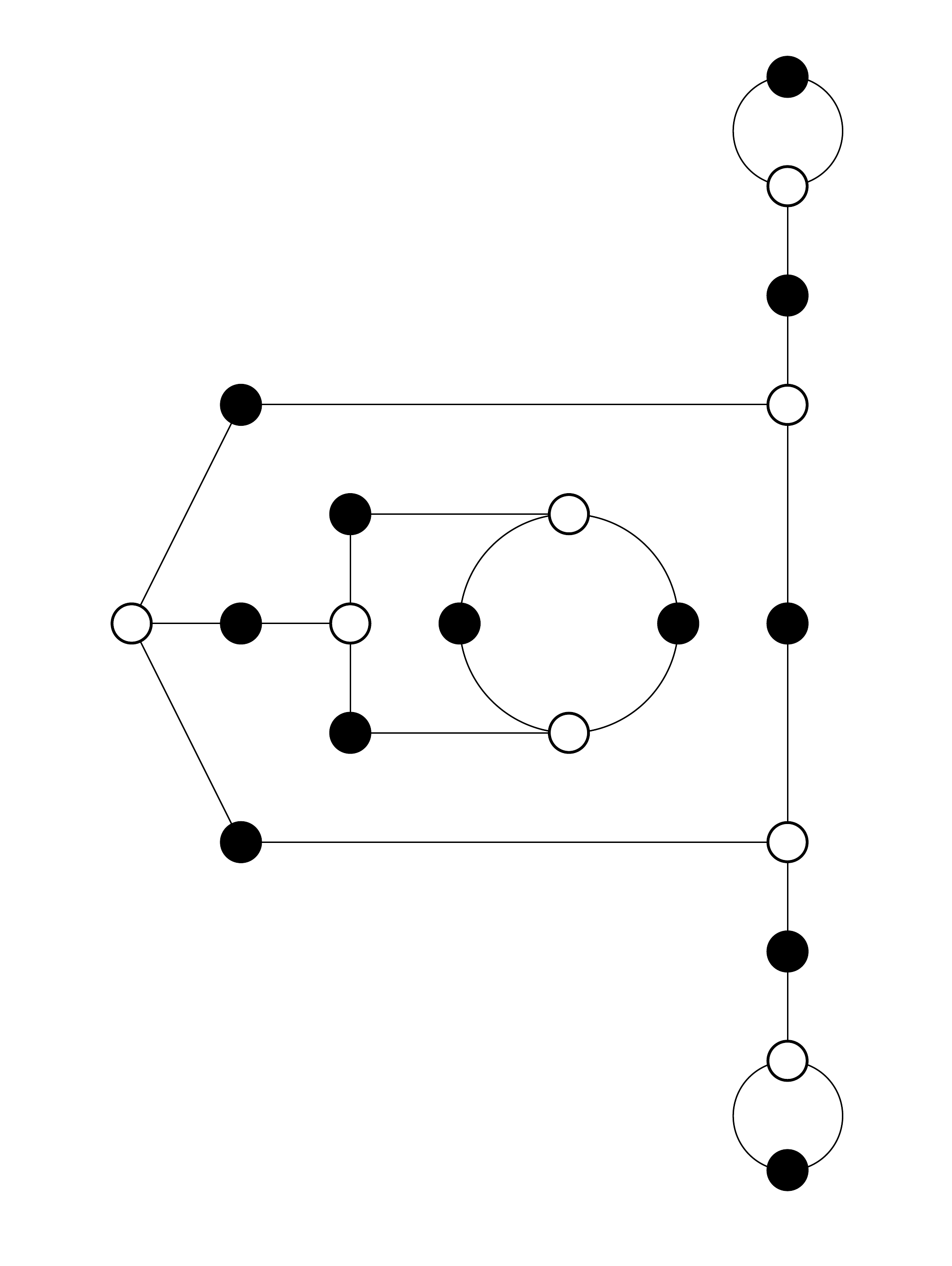}}
\par\end{center}{\scriptsize \par}

\begin{center}
{\scriptsize $9,8,3,2,1,1\;\left(\mathrm{cubic}\right)$}
\par\end{center}%
\end{minipage}
\par\end{center}{\scriptsize \par}

\begin{center}
{\scriptsize }%
\begin{minipage}[t]{0.33\textwidth}%
\begin{center}
{\scriptsize \includegraphics[scale=0.15]{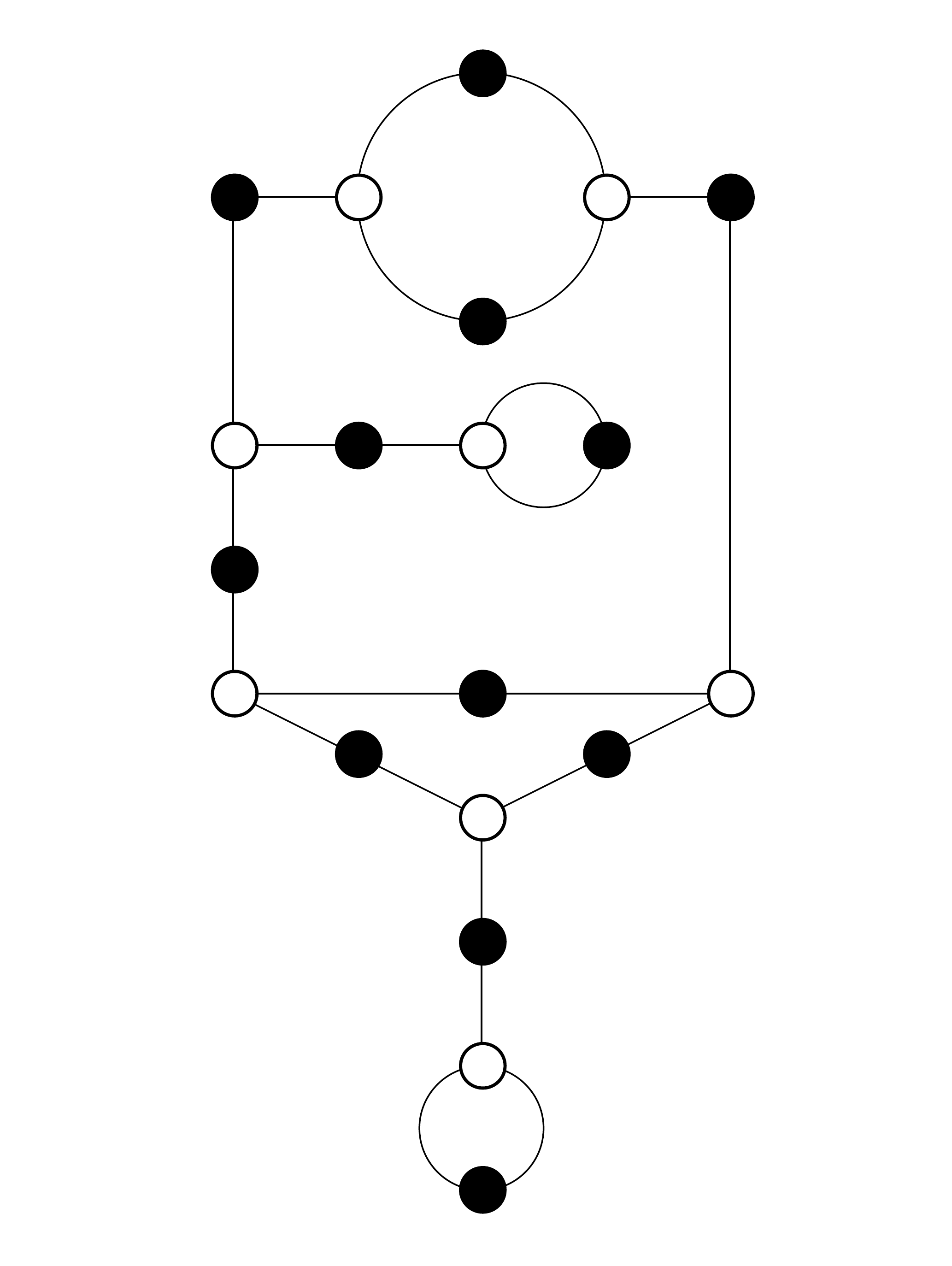}}
\par\end{center}{\scriptsize \par}

\begin{center}
{\scriptsize $9,8,3,2,1,1\;\left(\mathrm{cubic}\right)$}
\par\end{center}%
\end{minipage}{\scriptsize }%
\begin{minipage}[t]{0.33\textwidth}%
\begin{center}
{\scriptsize \includegraphics[scale=0.15]{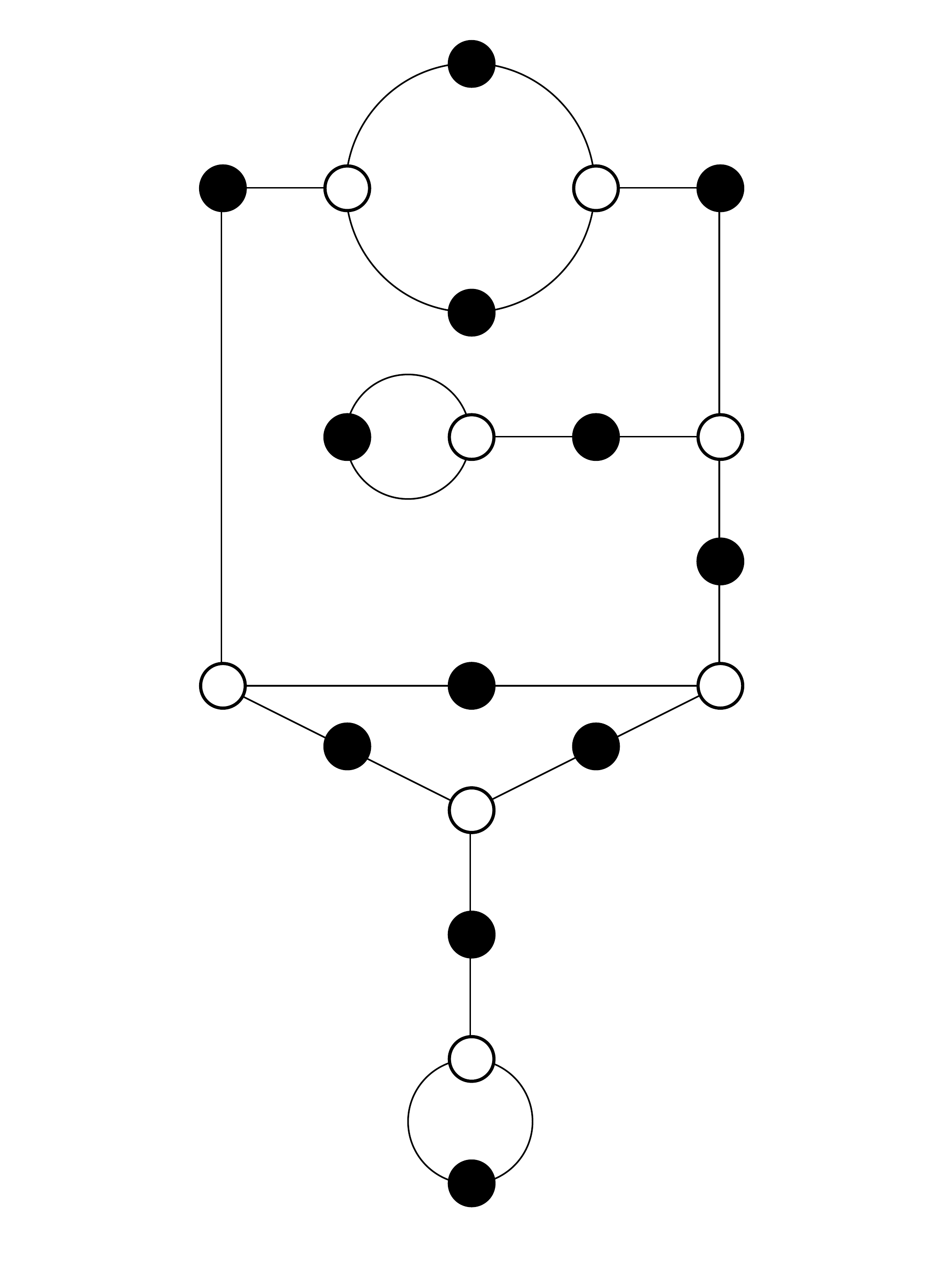}}
\par\end{center}{\scriptsize \par}

\begin{center}
{\scriptsize $9,8,3,2,1,1\;\left(\mathrm{cubic}\right)$}
\par\end{center}%
\end{minipage}{\scriptsize }%
\begin{minipage}[t]{0.33\textwidth}%
\begin{center}
{\scriptsize \includegraphics[scale=0.15]{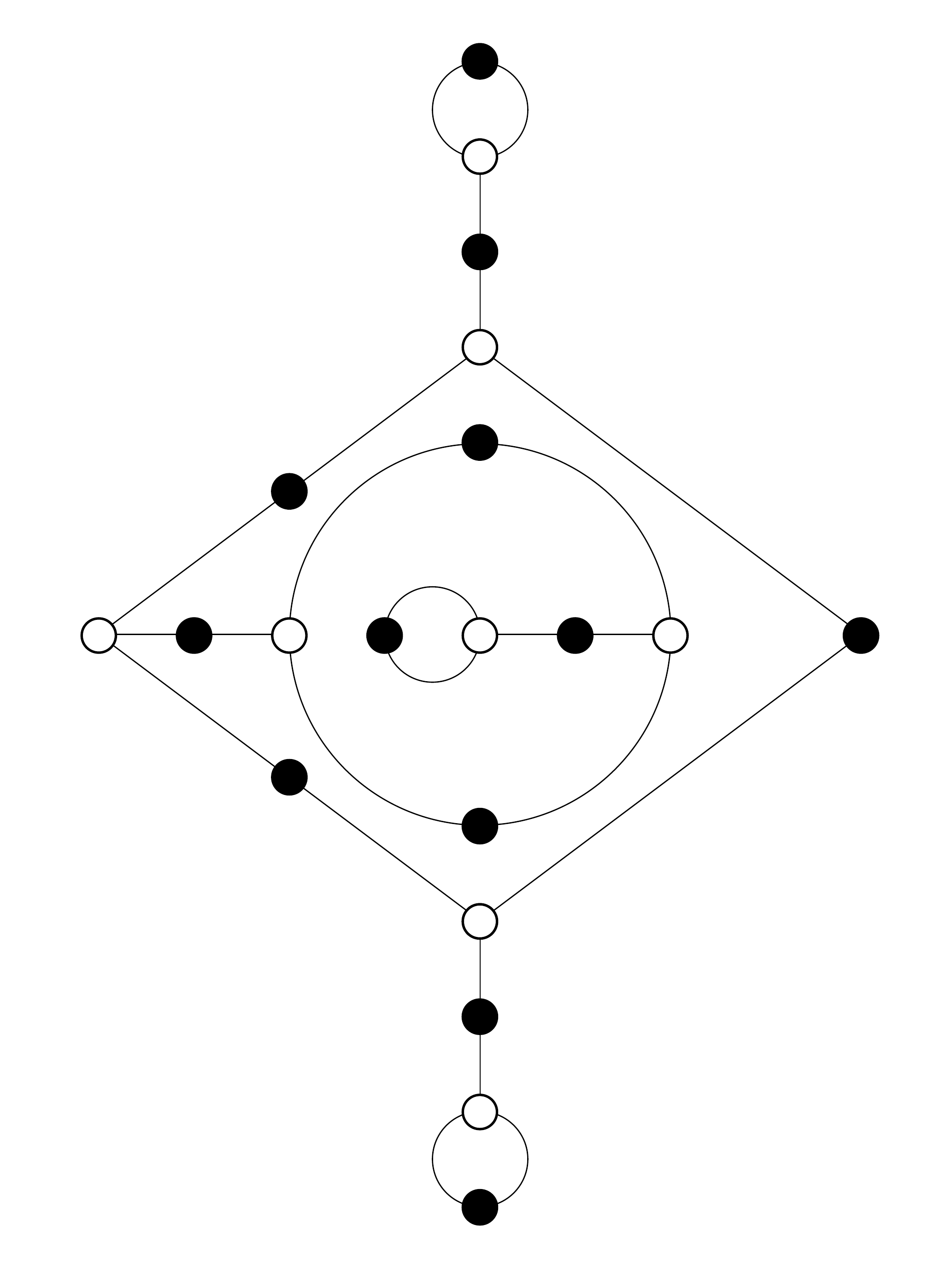}}
\par\end{center}{\scriptsize \par}

\begin{center}
{\scriptsize $9,7,5,1,1,1\;\left(\mathbb{Q}\right)$}
\par\end{center}%
\end{minipage}
\par\end{center}{\scriptsize \par}

\begin{center}
{\scriptsize }%
\begin{minipage}[t]{0.33\textwidth}%
\begin{center}
{\scriptsize \includegraphics[scale=0.15]{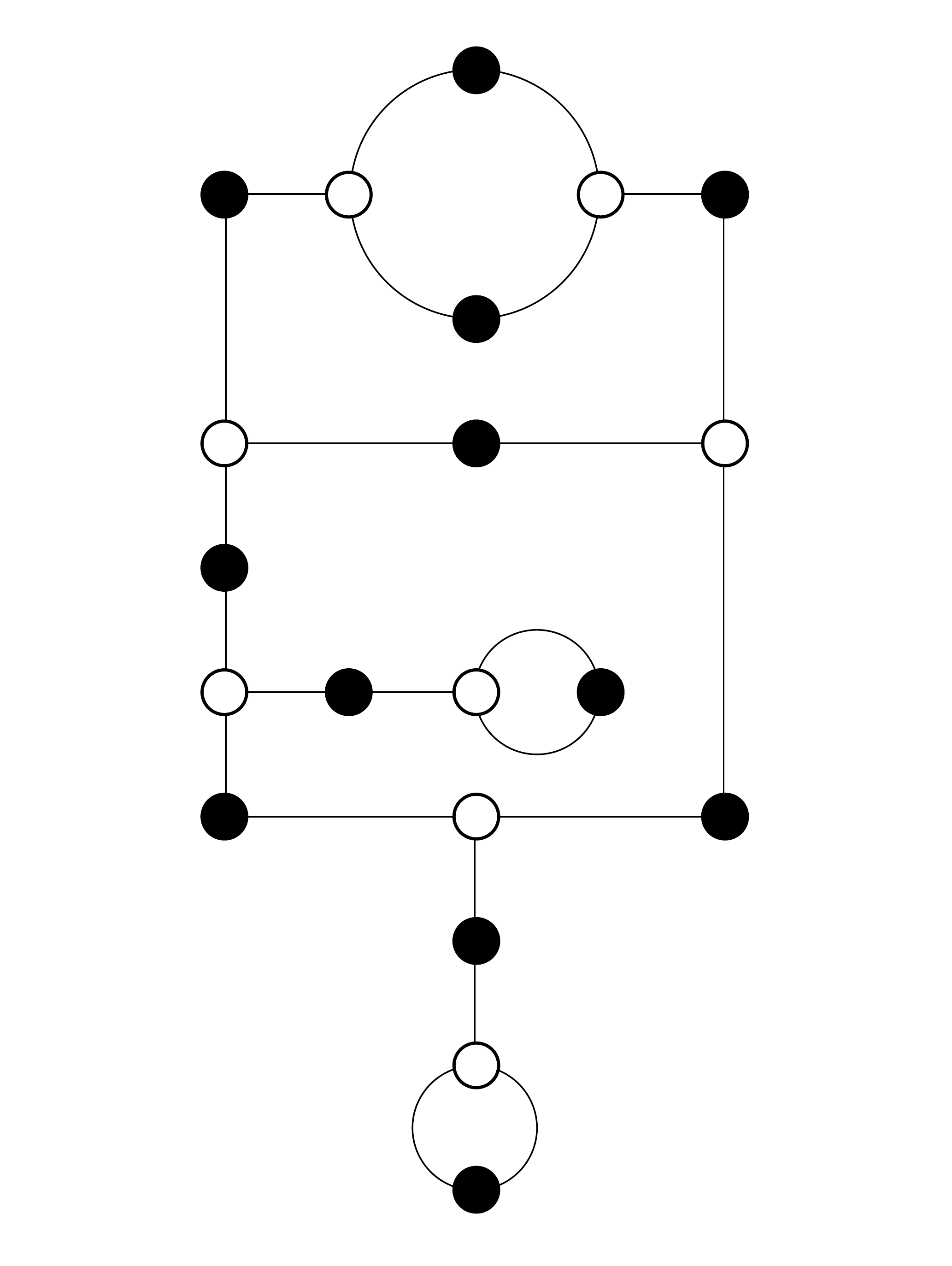}}
\par\end{center}{\scriptsize \par}

\begin{center}
{\scriptsize $9,7,4,2,1,1\;\left(\sqrt{-7}\right)$}
\par\end{center}%
\end{minipage}{\scriptsize }%
\begin{minipage}[t]{0.33\textwidth}%
\begin{center}
{\scriptsize \includegraphics[scale=0.15]{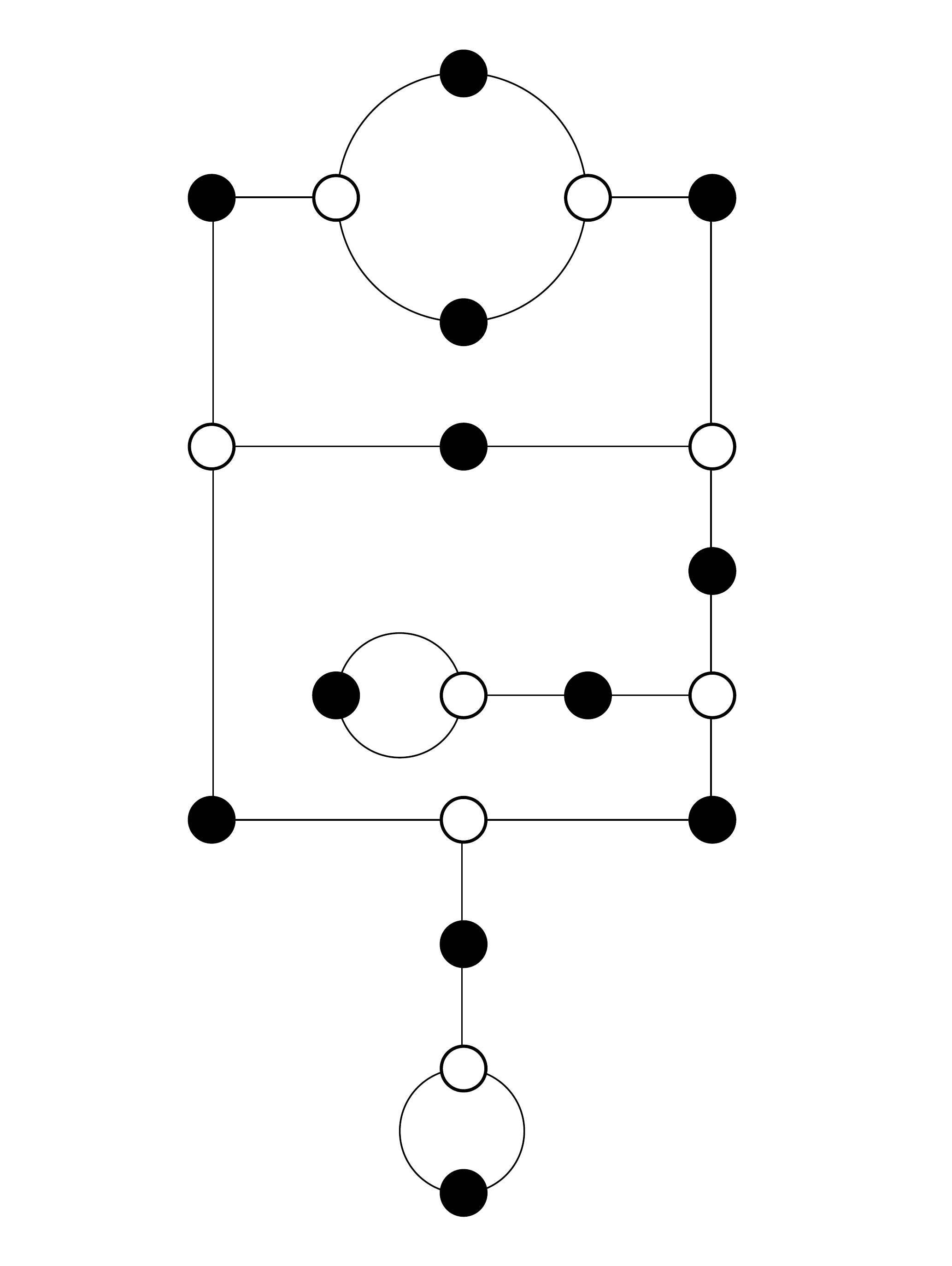}}
\par\end{center}{\scriptsize \par}

\begin{center}
{\scriptsize $9,7,4,2,1,1\;\left(\sqrt{-7}\right)$}
\par\end{center}%
\end{minipage}{\scriptsize }%
\begin{minipage}[t]{0.33\textwidth}%
\begin{center}
{\scriptsize \includegraphics[scale=0.15]{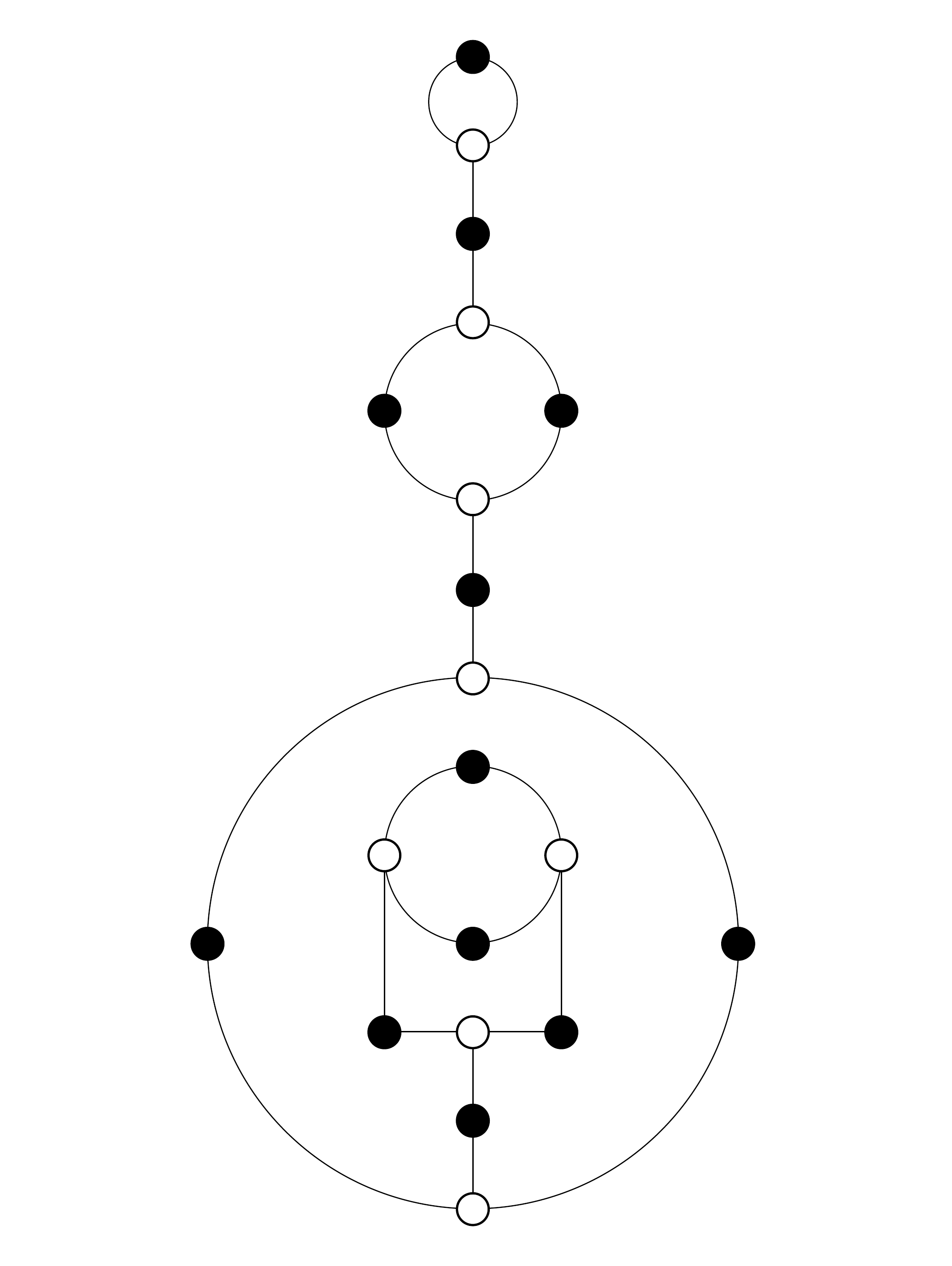}}
\par\end{center}{\scriptsize \par}

\begin{center}
{\scriptsize $9,7,3,2,2,1\;\left(\mathrm{cubic}\right)$}
\par\end{center}%
\end{minipage}
\par\end{center}{\scriptsize \par}

\begin{center}
{\scriptsize }%
\begin{minipage}[t]{0.33\textwidth}%
\begin{center}
{\scriptsize \includegraphics[scale=0.15]{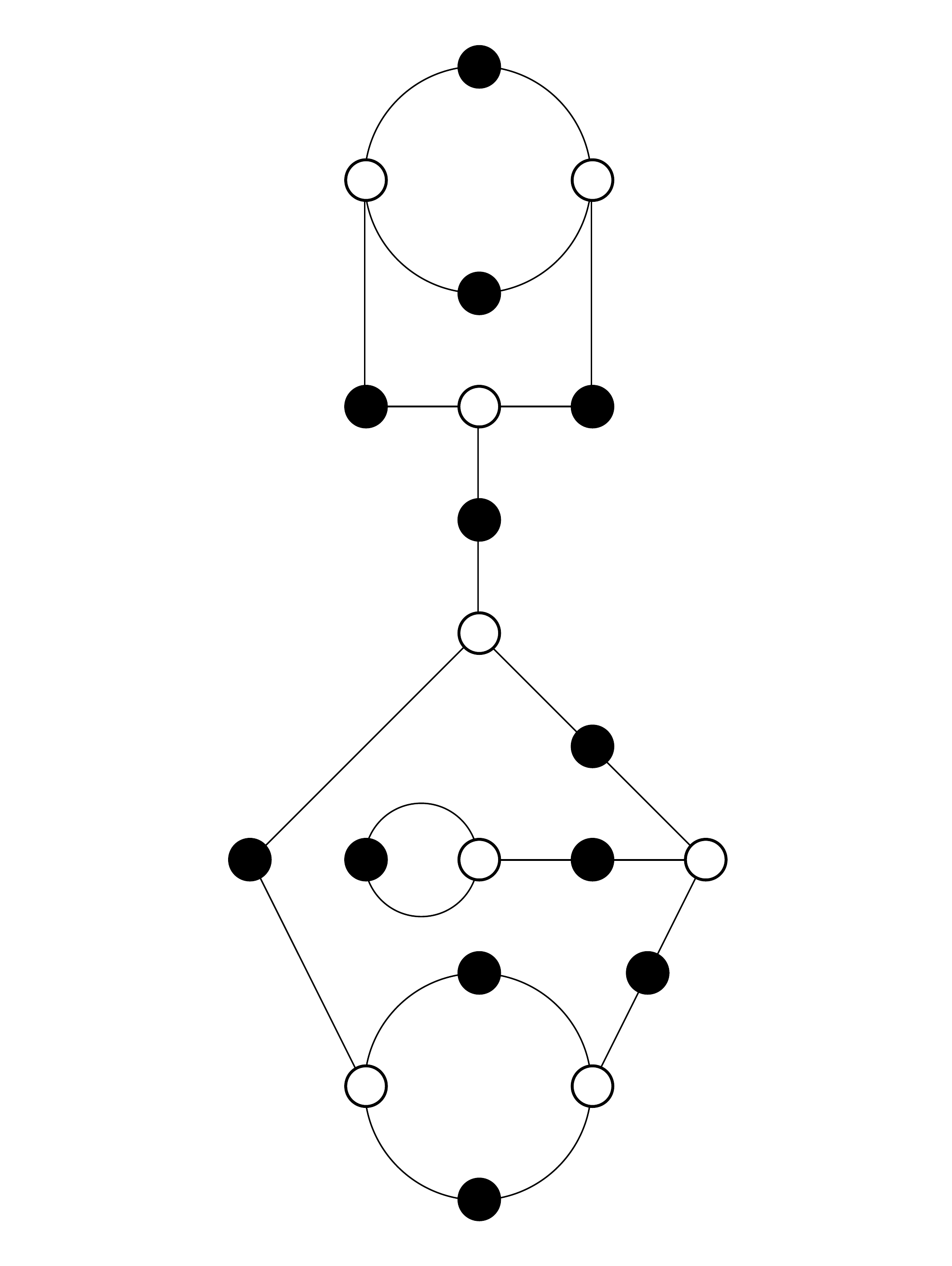}}
\par\end{center}{\scriptsize \par}

\begin{center}
{\scriptsize $9,7,3,2,2,1\;\left(\mathrm{cubic}\right)$}
\par\end{center}%
\end{minipage}{\scriptsize }%
\begin{minipage}[t]{0.33\textwidth}%
\begin{center}
{\scriptsize \includegraphics[scale=0.15]{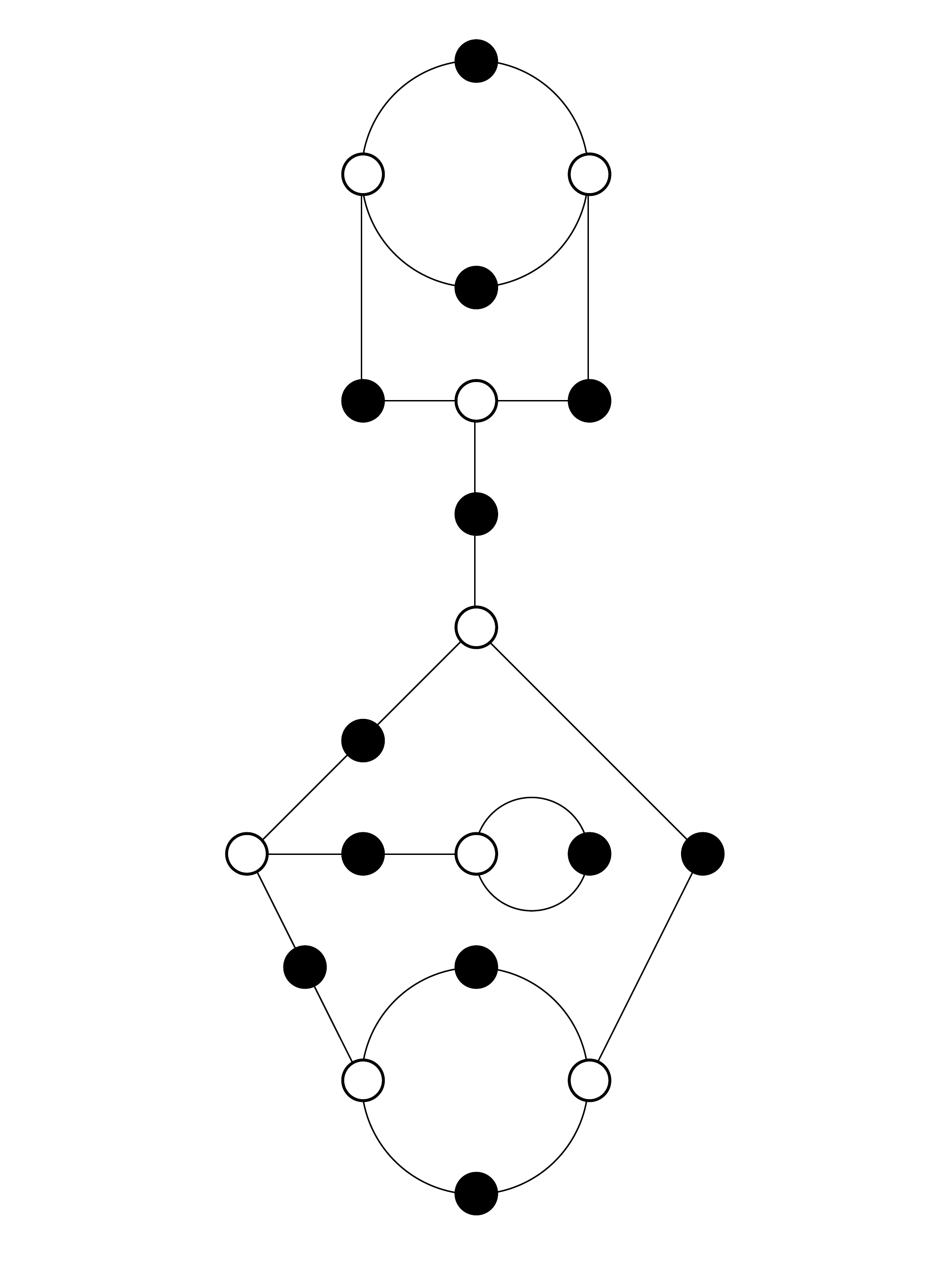}}
\par\end{center}{\scriptsize \par}

\begin{center}
{\scriptsize $9,7,3,2,2,1\;\left(\mathrm{cubic}\right)$}
\par\end{center}%
\end{minipage}{\scriptsize }%
\begin{minipage}[t]{0.33\textwidth}%
\begin{center}
{\scriptsize \includegraphics[scale=0.15]{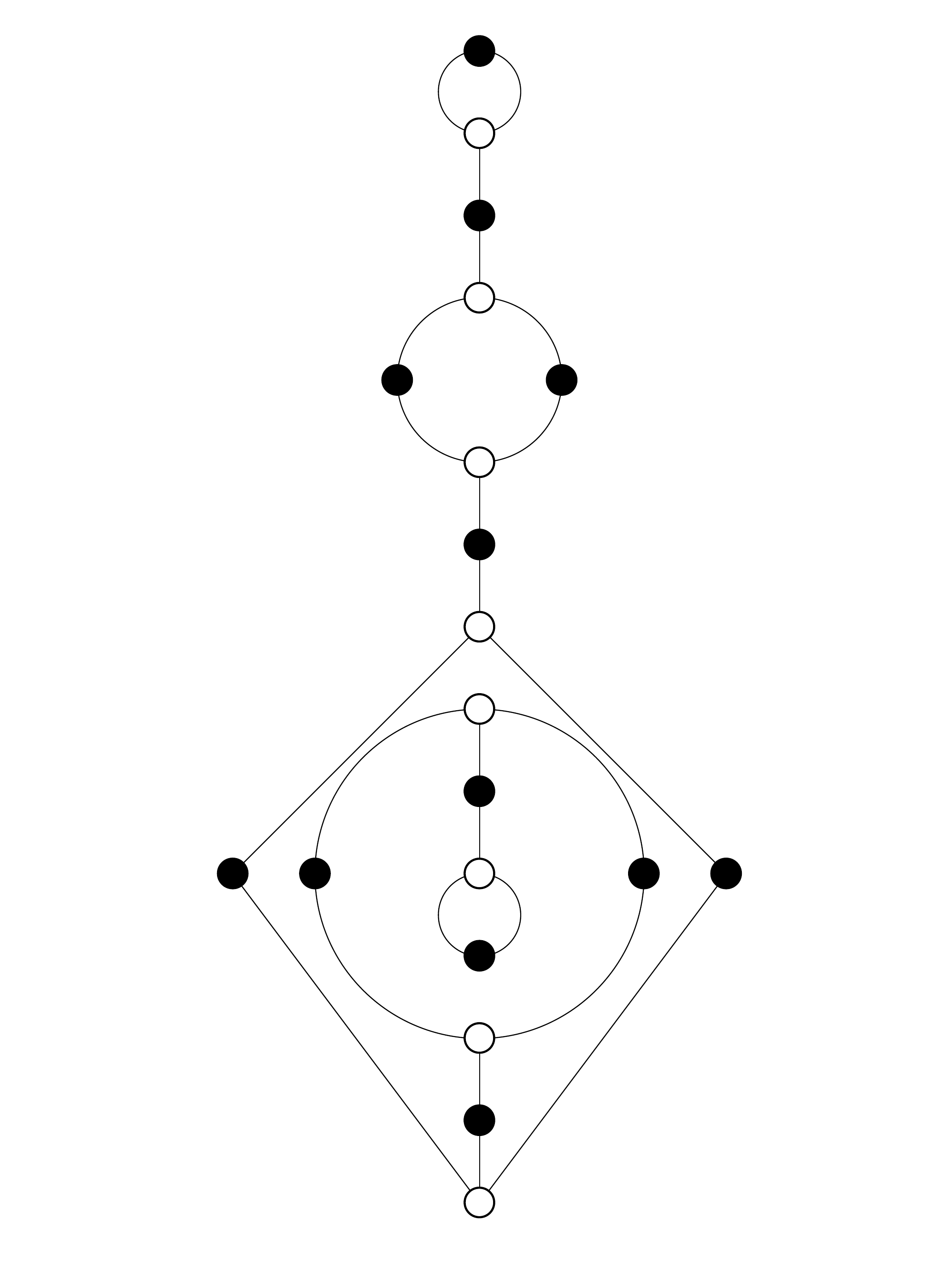}}
\par\end{center}{\scriptsize \par}

\begin{center}
{\scriptsize $9,6,5,2,1,1\;\left(\mathrm{cubic}\right)$}
\par\end{center}%
\end{minipage}
\par\end{center}{\scriptsize \par}

\begin{center}
{\scriptsize }%
\begin{minipage}[t]{0.33\textwidth}%
\begin{center}
{\scriptsize \includegraphics[scale=0.15]{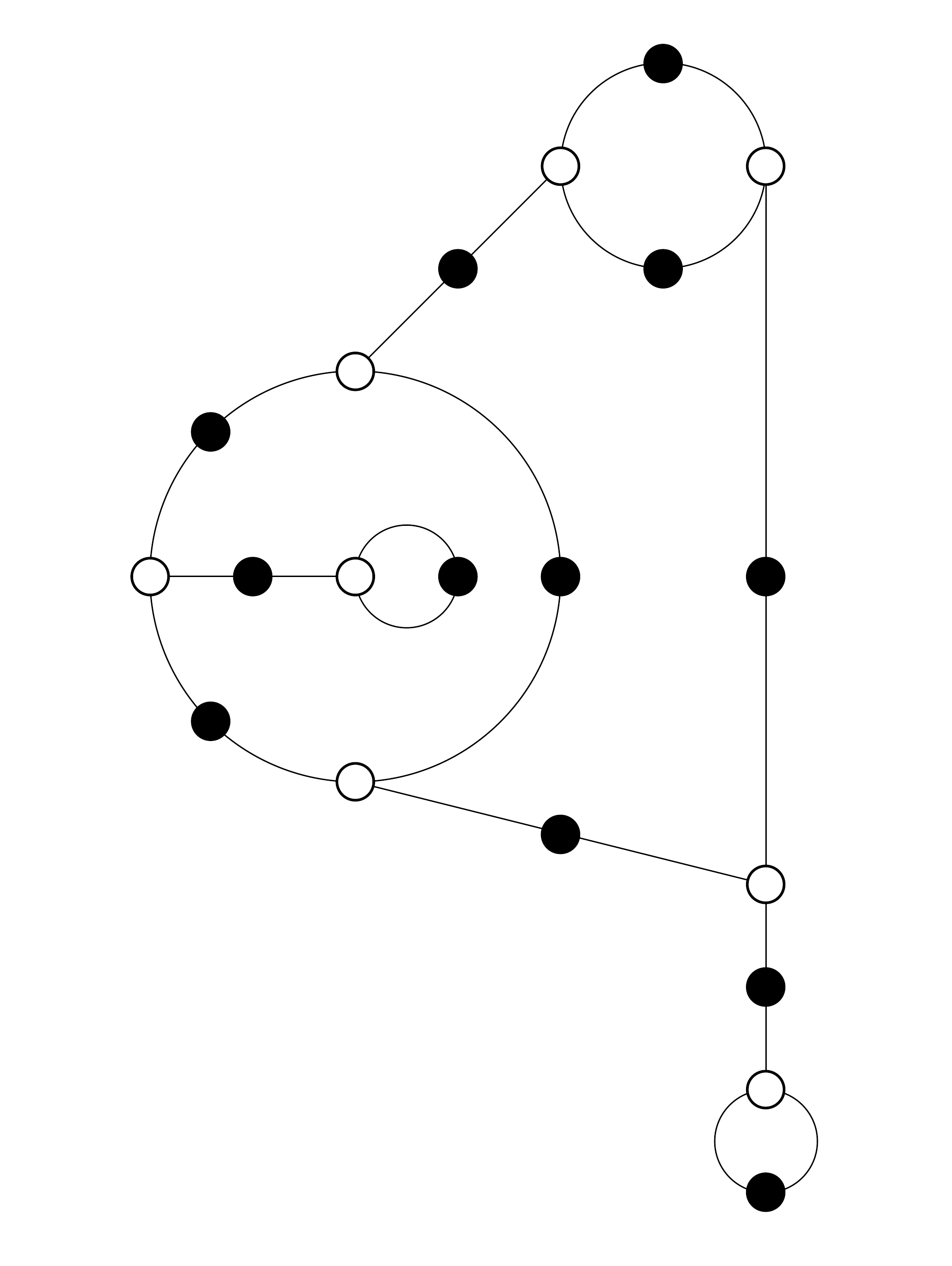}}
\par\end{center}{\scriptsize \par}

\begin{center}
{\scriptsize $9,6,5,2,1,1\;\left(\mathrm{cubic}\right)$}
\par\end{center}%
\end{minipage}{\scriptsize }%
\begin{minipage}[t]{0.33\textwidth}%
\begin{center}
{\scriptsize \includegraphics[scale=0.15]{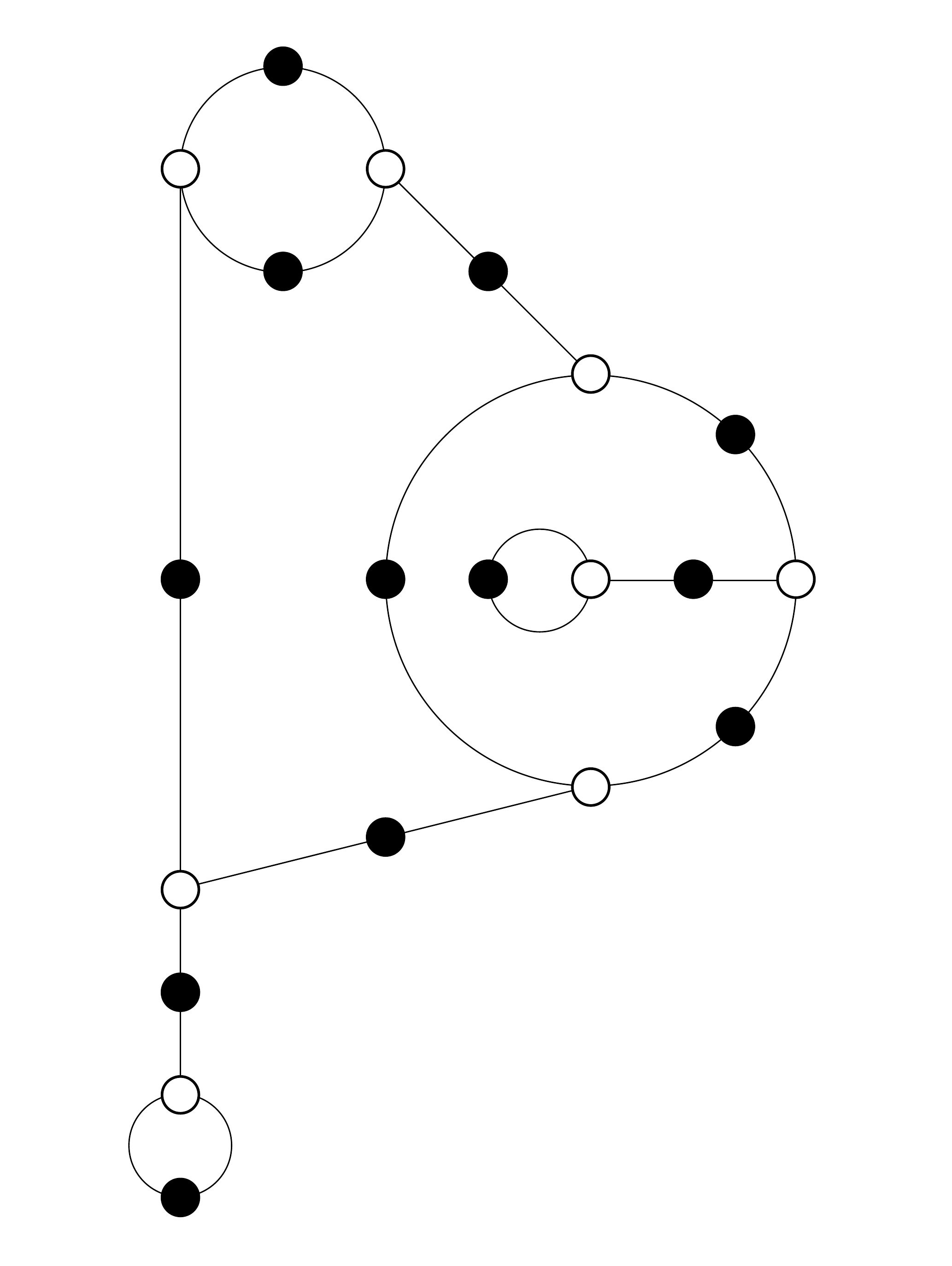}}
\par\end{center}{\scriptsize \par}

\begin{center}
{\scriptsize $9,6,5,2,1,1\;\left(\mathrm{cubic}\right)$}
\par\end{center}%
\end{minipage}{\scriptsize }%
\begin{minipage}[t]{0.33\textwidth}%
\begin{center}
{\scriptsize \includegraphics[scale=0.15]{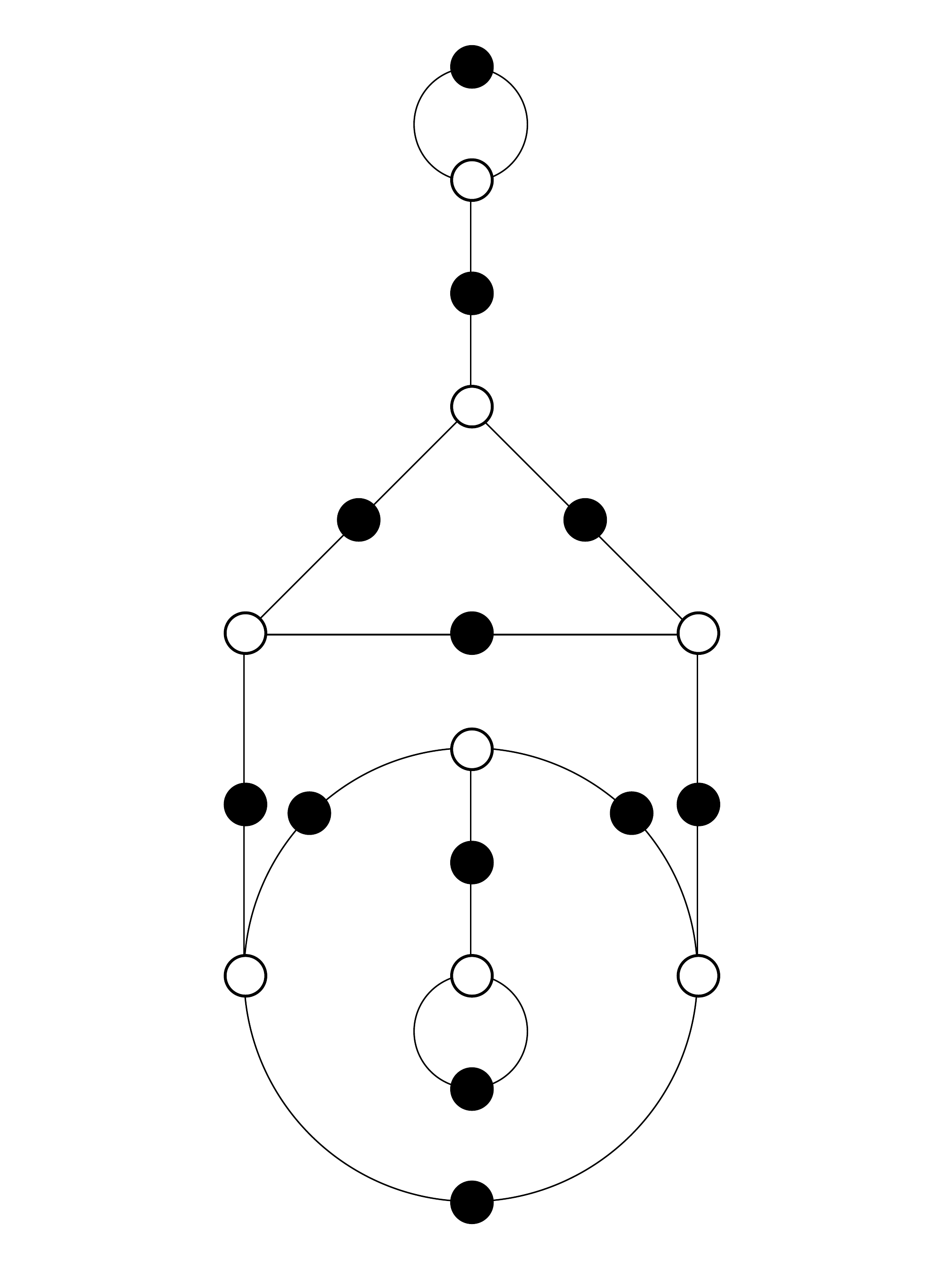}}
\par\end{center}{\scriptsize \par}

\begin{center}
{\scriptsize $9,6,4,3,1,1\;\left(\mathbb{Q}\right)$}
\par\end{center}%
\end{minipage}
\par\end{center}{\scriptsize \par}

\begin{center}
{\scriptsize }%
\begin{minipage}[t]{0.33\textwidth}%
\begin{center}
{\scriptsize \includegraphics[scale=0.15]{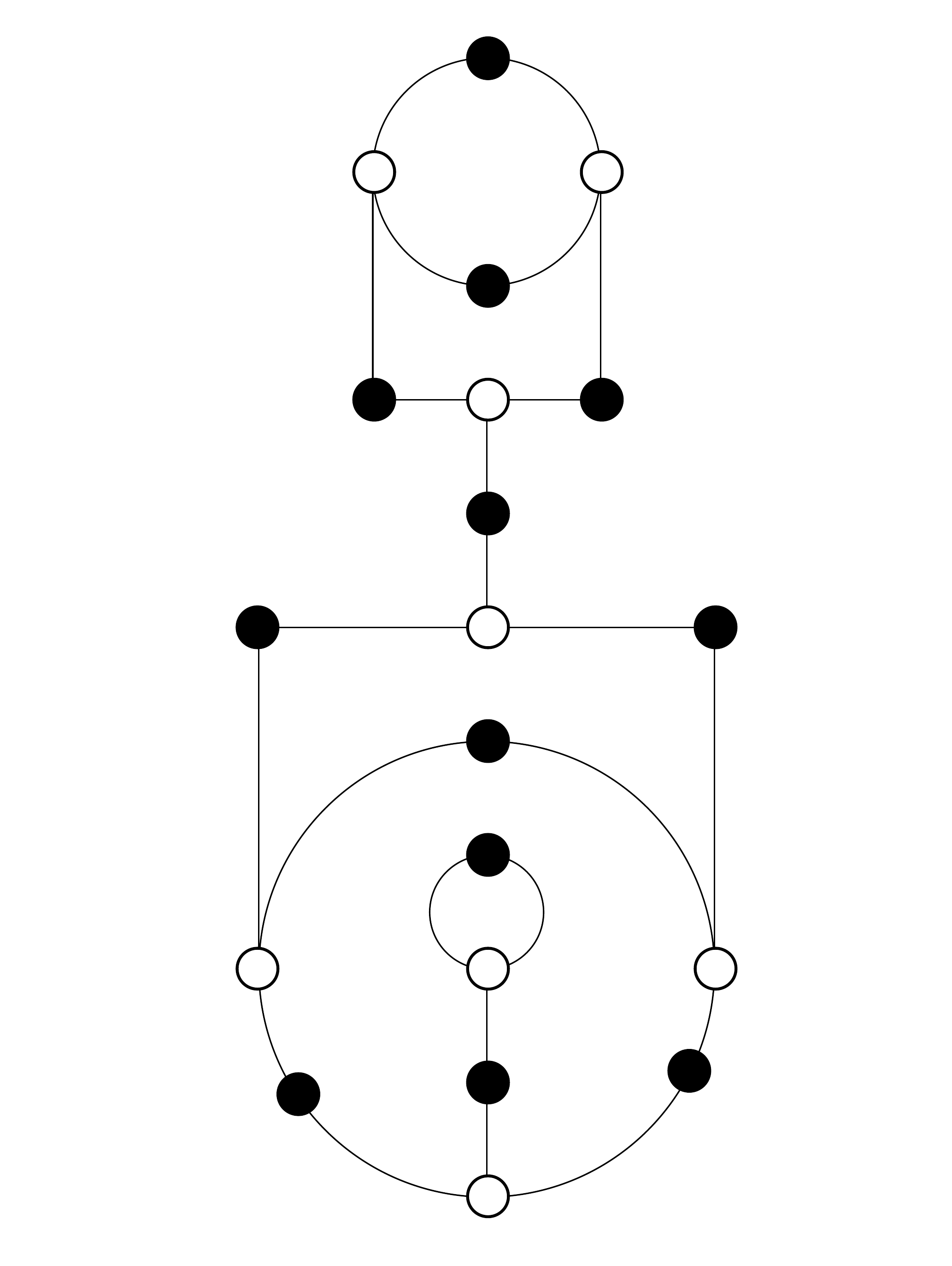}}
\par\end{center}{\scriptsize \par}

\begin{center}
{\scriptsize $9,6,3,3,2,1\;\left(\mathrm{cubic}\right)$}
\par\end{center}%
\end{minipage}{\scriptsize }%
\begin{minipage}[t]{0.33\textwidth}%
\begin{center}
{\scriptsize \includegraphics[scale=0.15]{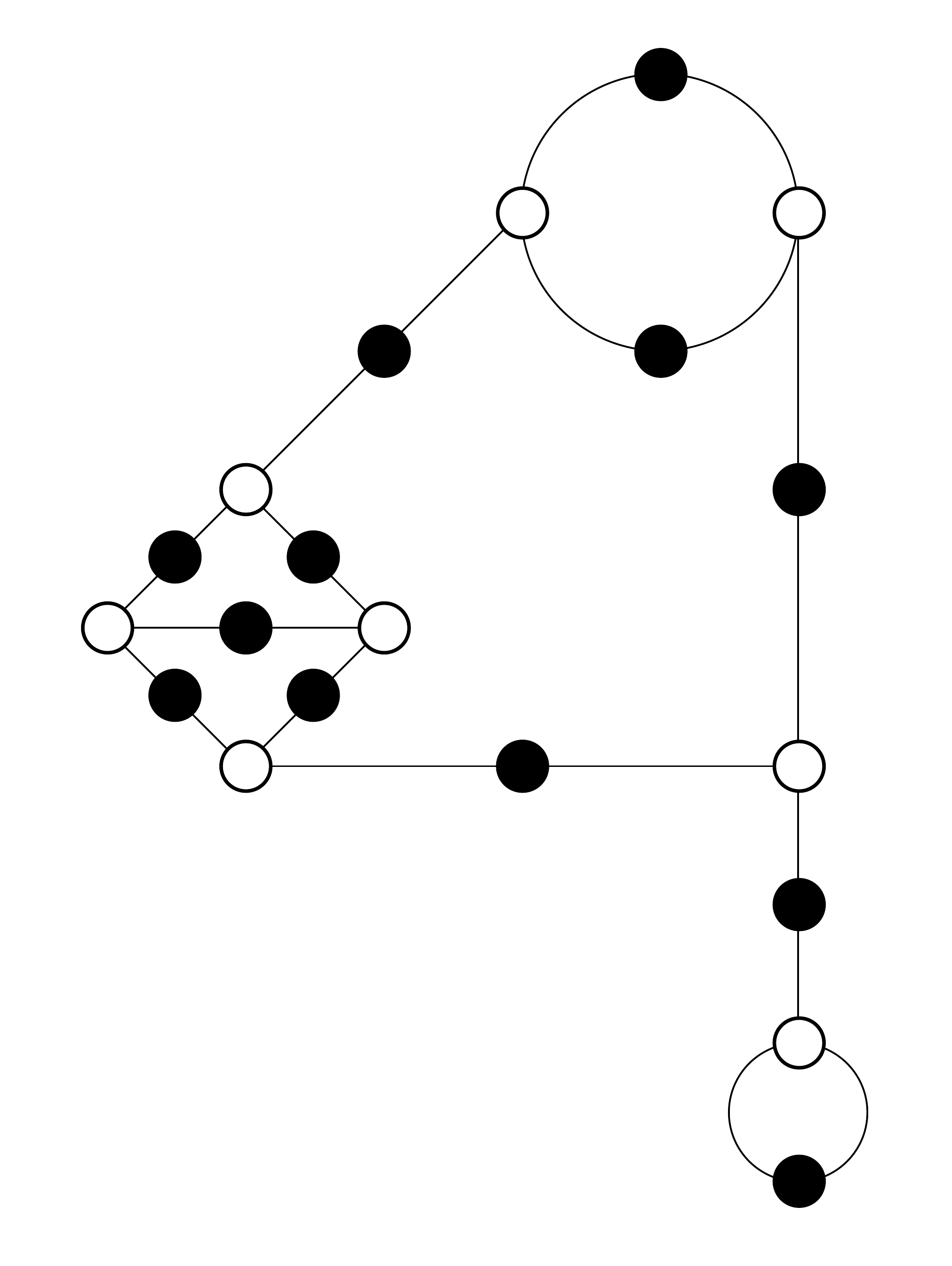}}
\par\end{center}{\scriptsize \par}

\begin{center}
{\scriptsize $9,6,3,3,2,1\;\left(\mathrm{cubic}\right)$}
\par\end{center}%
\end{minipage}{\scriptsize }%
\begin{minipage}[t]{0.33\textwidth}%
\begin{center}
{\scriptsize \includegraphics[scale=0.15]{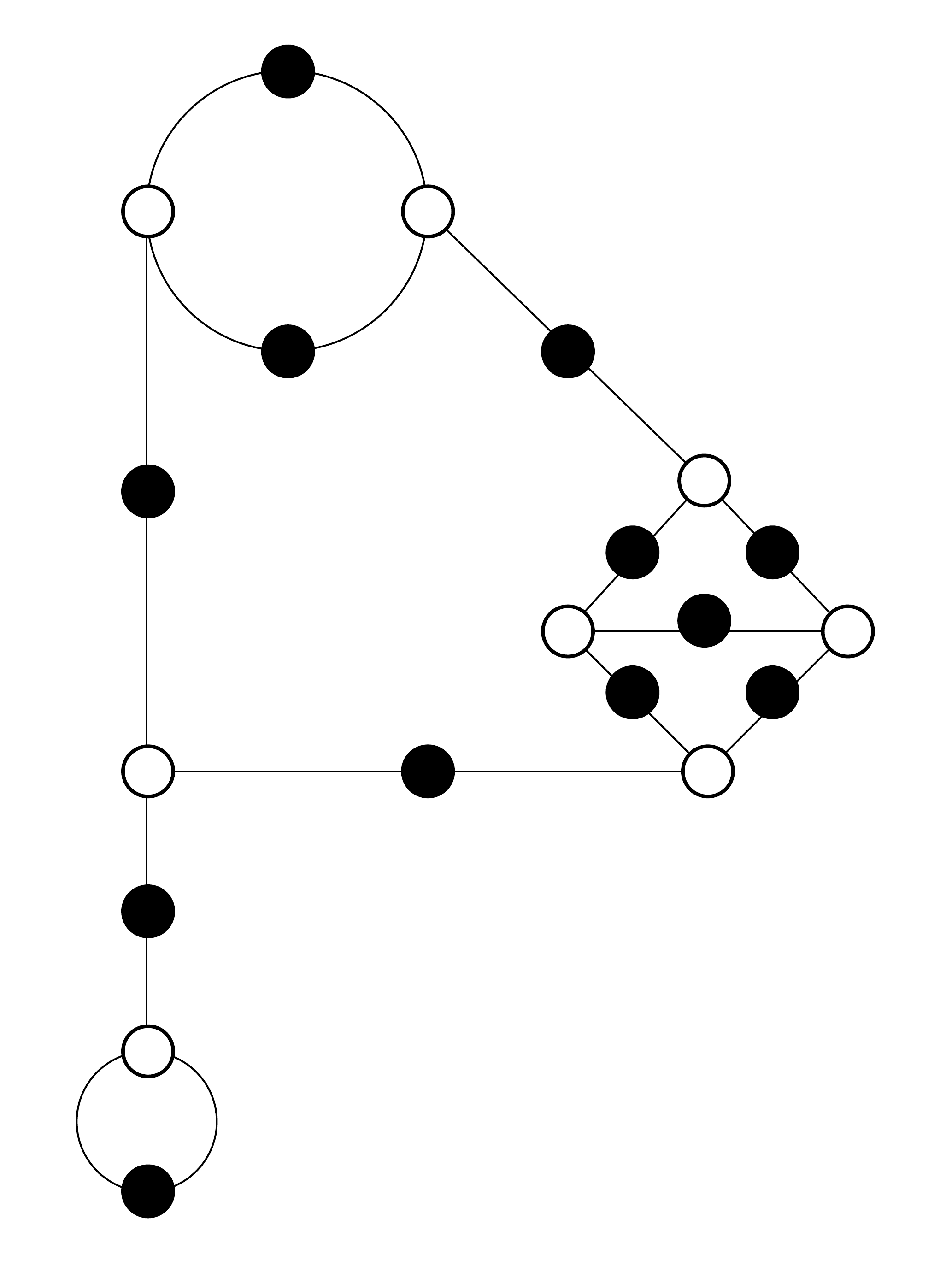}}
\par\end{center}{\scriptsize \par}

\begin{center}
{\scriptsize $9,6,3,3,2,1\;\left(\mathrm{cubic}\right)$}
\par\end{center}%
\end{minipage}
\par\end{center}{\scriptsize \par}

\begin{center}
{\scriptsize }%
\begin{minipage}[t]{0.33\textwidth}%
\begin{center}
{\scriptsize \includegraphics[scale=0.15]{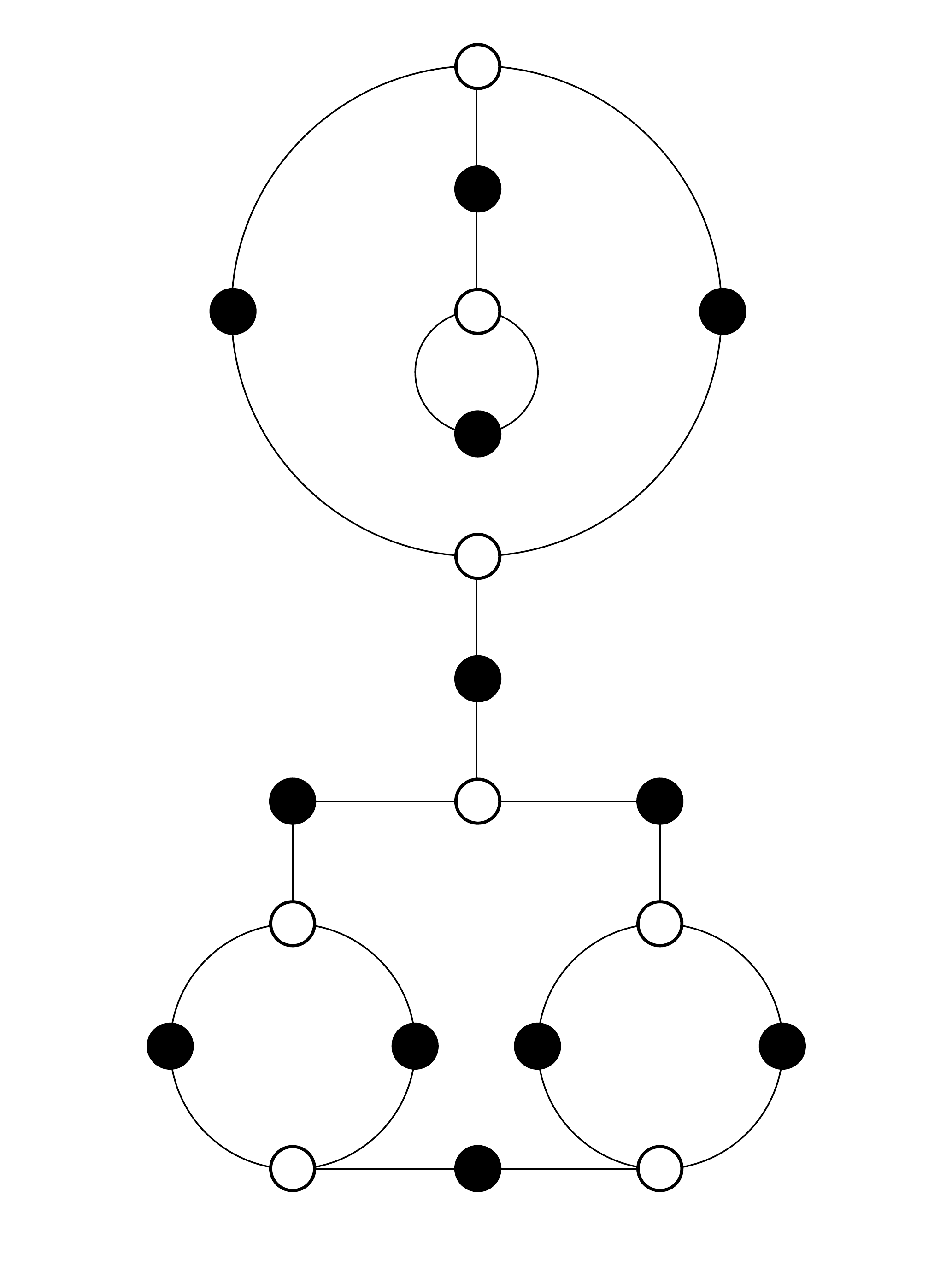}}
\par\end{center}{\scriptsize \par}

\begin{center}
{\scriptsize $9,5,5,2,2,1\;\left(\mathbb{Q}\right)$}
\par\end{center}%
\end{minipage}{\scriptsize }%
\begin{minipage}[t]{0.33\textwidth}%
\begin{center}
{\scriptsize \includegraphics[scale=0.15]{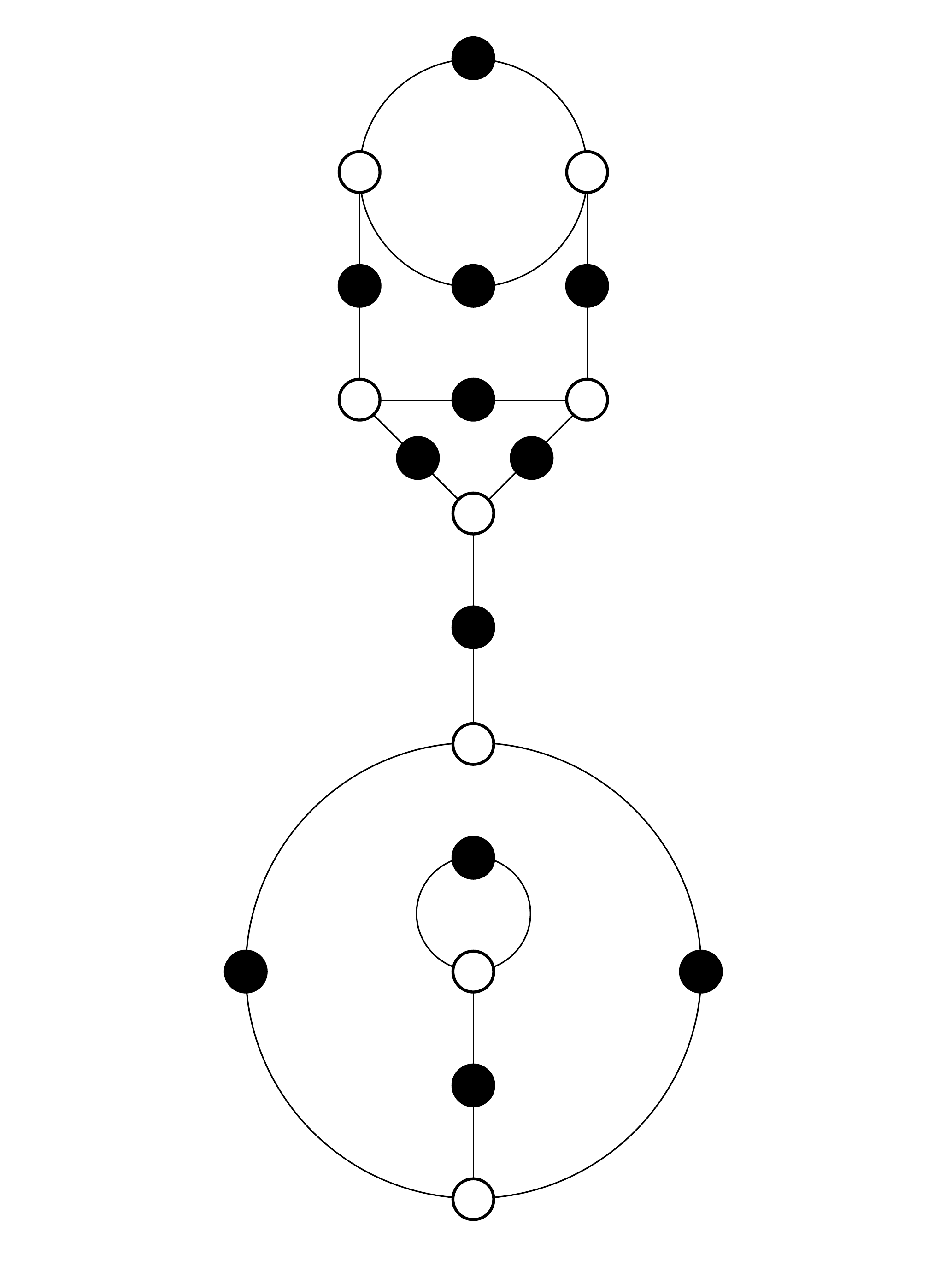}}
\par\end{center}{\scriptsize \par}

\begin{center}
{\scriptsize $9,5,4,3,2,1\;\left(\mathrm{cubic}\right)$}
\par\end{center}%
\end{minipage}{\scriptsize }%
\begin{minipage}[t]{0.33\textwidth}%
\begin{center}
{\scriptsize \includegraphics[scale=0.15]{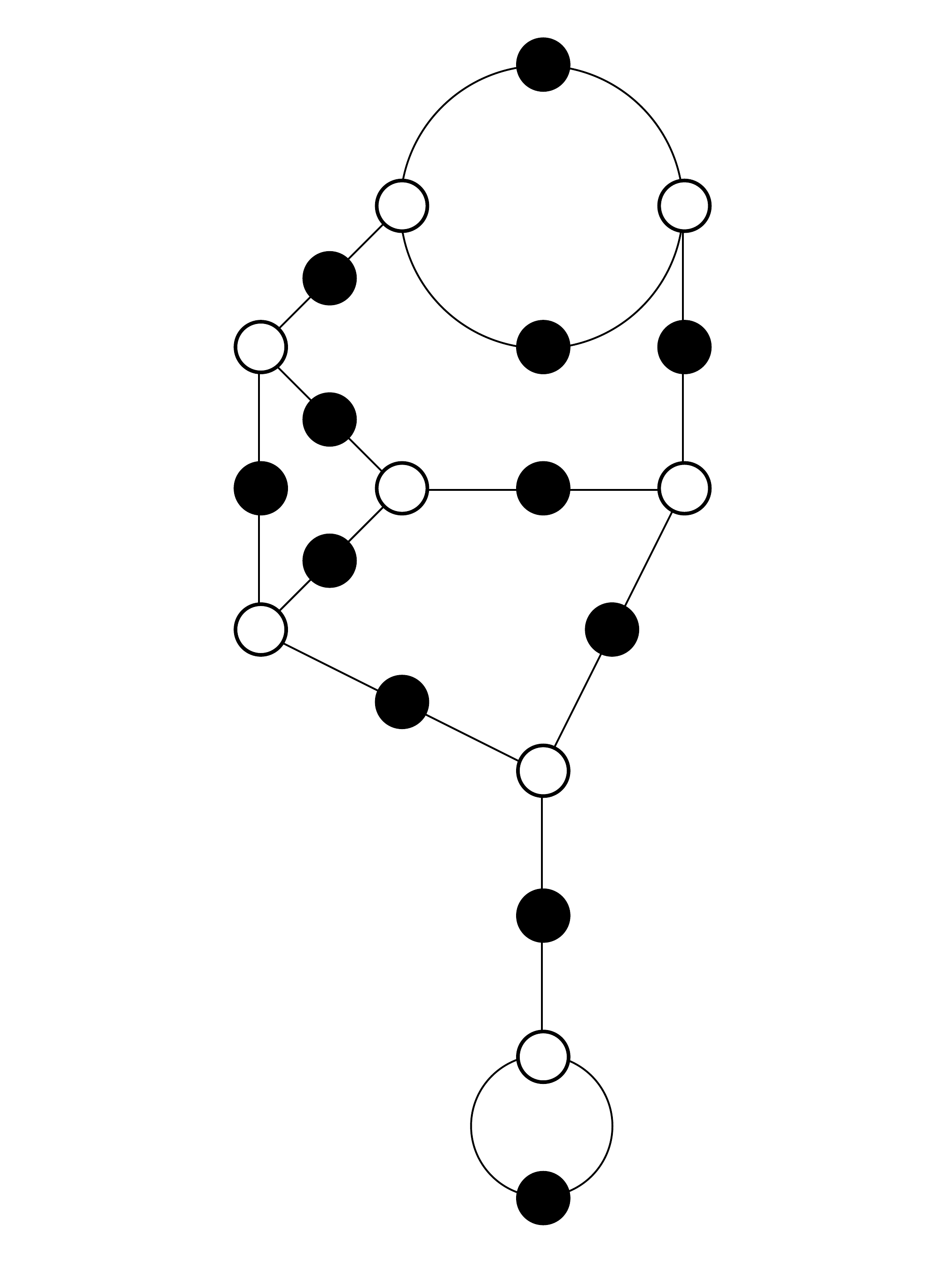}}
\par\end{center}{\scriptsize \par}

\begin{center}
{\scriptsize $9,5,4,3,2,1\;\left(\mathrm{cubic}\right)$}
\par\end{center}%
\end{minipage}
\par\end{center}{\scriptsize \par}

\begin{center}
{\scriptsize }%
\begin{minipage}[t]{0.33\textwidth}%
\begin{center}
{\scriptsize \includegraphics[scale=0.15]{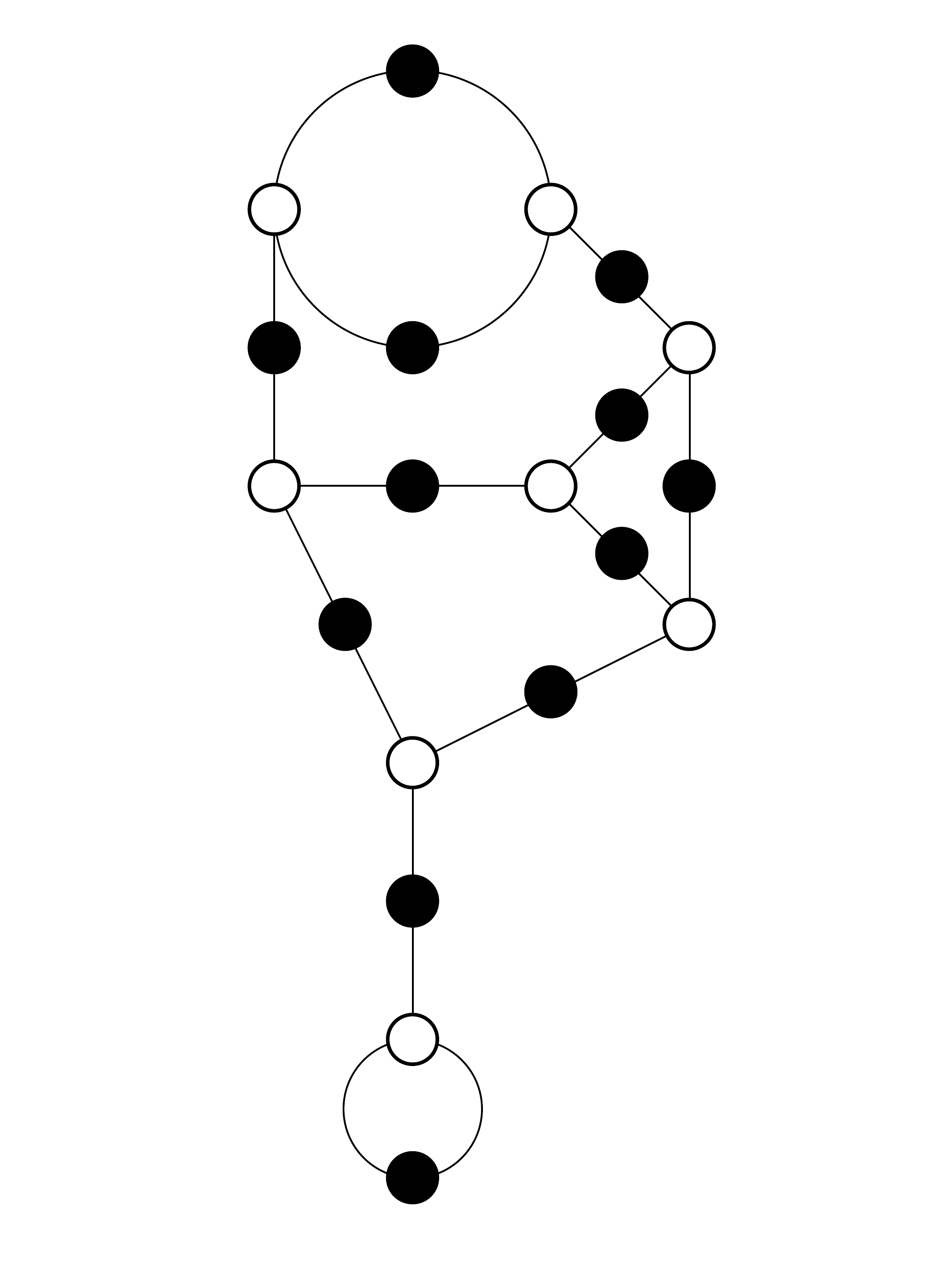}}
\par\end{center}{\scriptsize \par}

\begin{center}
{\scriptsize $9,5,4,3,2,1\;\left(\mathrm{cubic}\right)$}
\par\end{center}%
\end{minipage}{\scriptsize }%
\begin{minipage}[t]{0.33\textwidth}%
\begin{center}
{\scriptsize \includegraphics[scale=0.15]{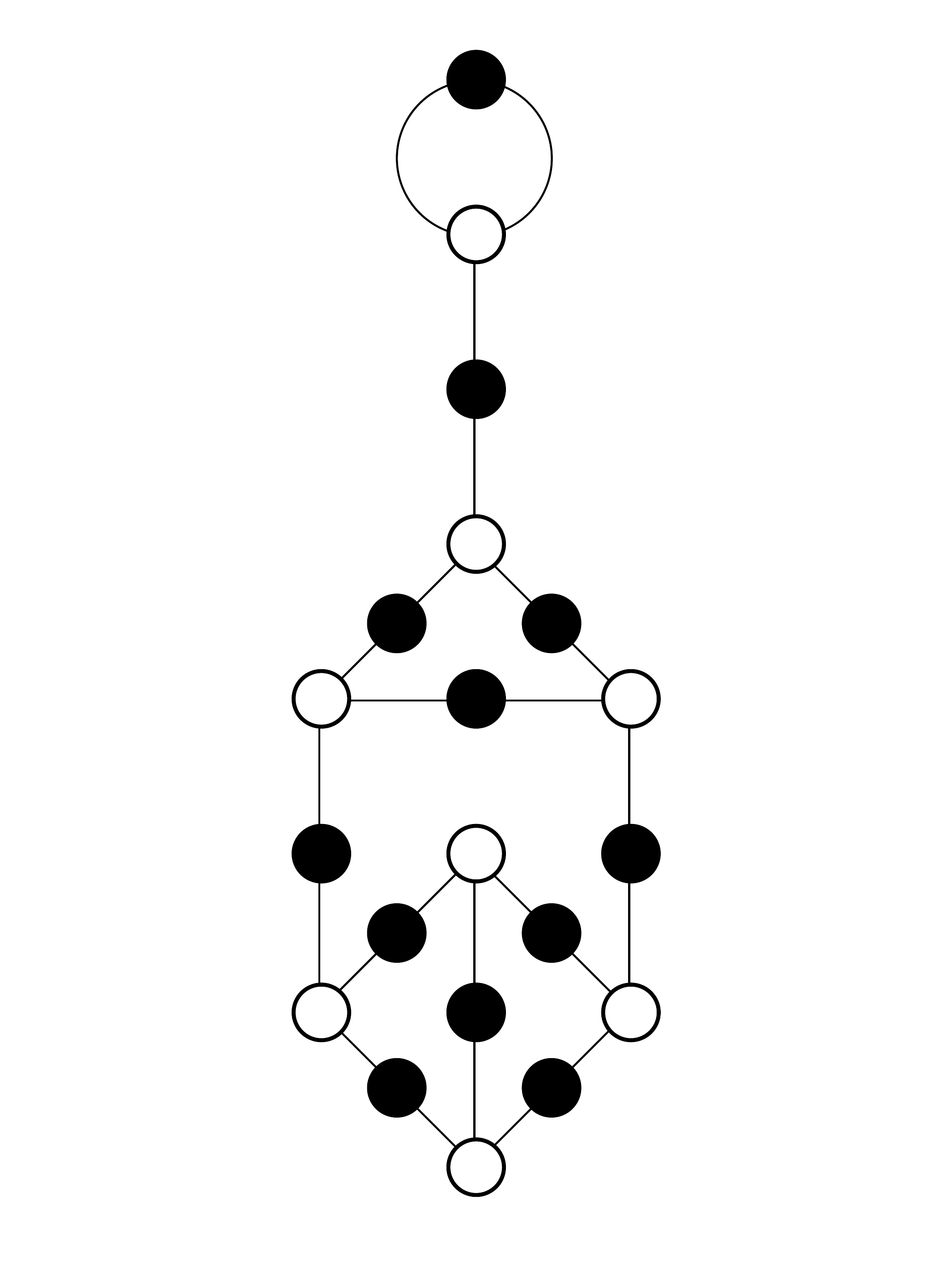}}
\par\end{center}{\scriptsize \par}

\begin{center}
{\scriptsize $9,5,3,3,3,1\;\left(\mathbb{Q}\right)$}
\par\end{center}%
\end{minipage}{\scriptsize }%
\begin{minipage}[t]{0.33\textwidth}%
\begin{center}
{\scriptsize \includegraphics[scale=0.15]{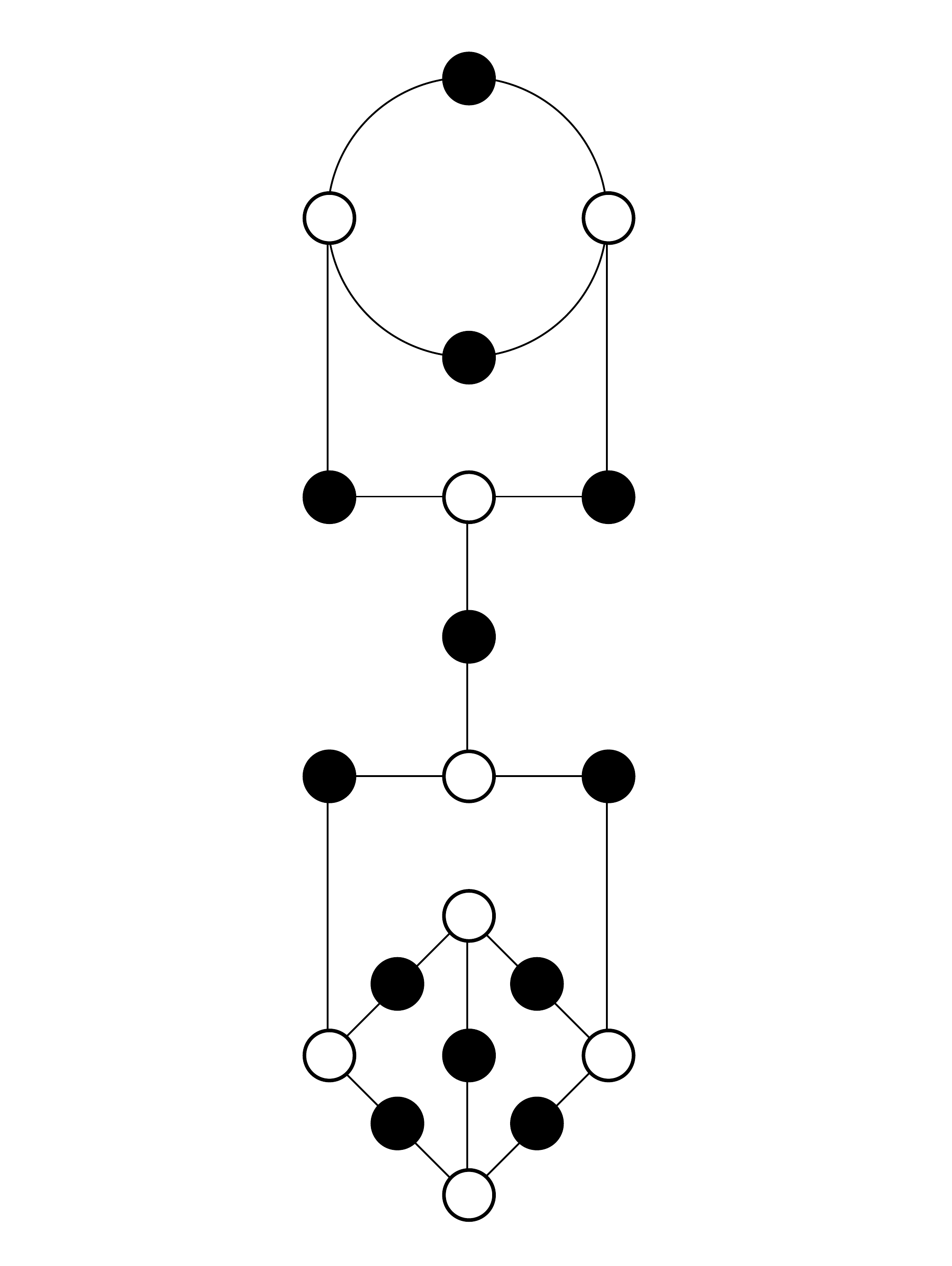}}
\par\end{center}{\scriptsize \par}

\begin{center}
{\scriptsize $9,4,3,3,3,2\;\left(\mathbb{Q}\right)$}
\par\end{center}%
\end{minipage}
\par\end{center}{\scriptsize \par}

\begin{center}
{\scriptsize }%
\begin{minipage}[t]{0.33\textwidth}%
\begin{center}
{\scriptsize \includegraphics[scale=0.15]{\string"PICT/8-8-4-2-1-1_Attempt_B\string".pdf}}
\par\end{center}{\scriptsize \par}

\begin{center}
{\scriptsize $8,8,4,2,1,1\;\left(\mathbb{Q}\right)$}
\par\end{center}%
\end{minipage}{\scriptsize }%
\begin{minipage}[t]{0.33\textwidth}%
\begin{center}
{\scriptsize \includegraphics[scale=0.15]{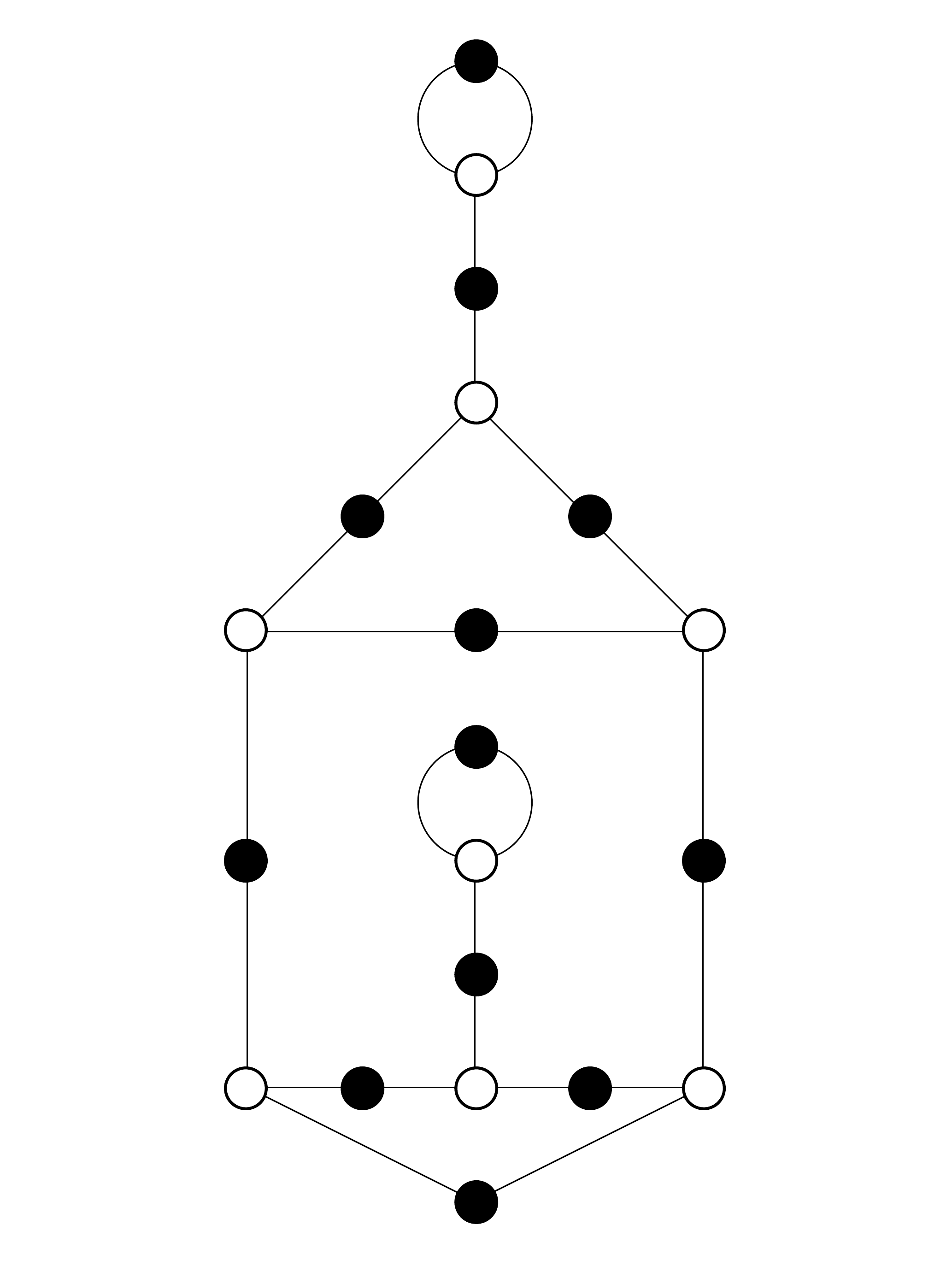}}
\par\end{center}{\scriptsize \par}

\begin{center}
{\scriptsize $8,8,3,3,1,1\;\left(\mathbb{Q}\right)$}
\par\end{center}%
\end{minipage}{\scriptsize }%
\begin{minipage}[t]{0.33\textwidth}%
\begin{center}
{\scriptsize \includegraphics[scale=0.15]{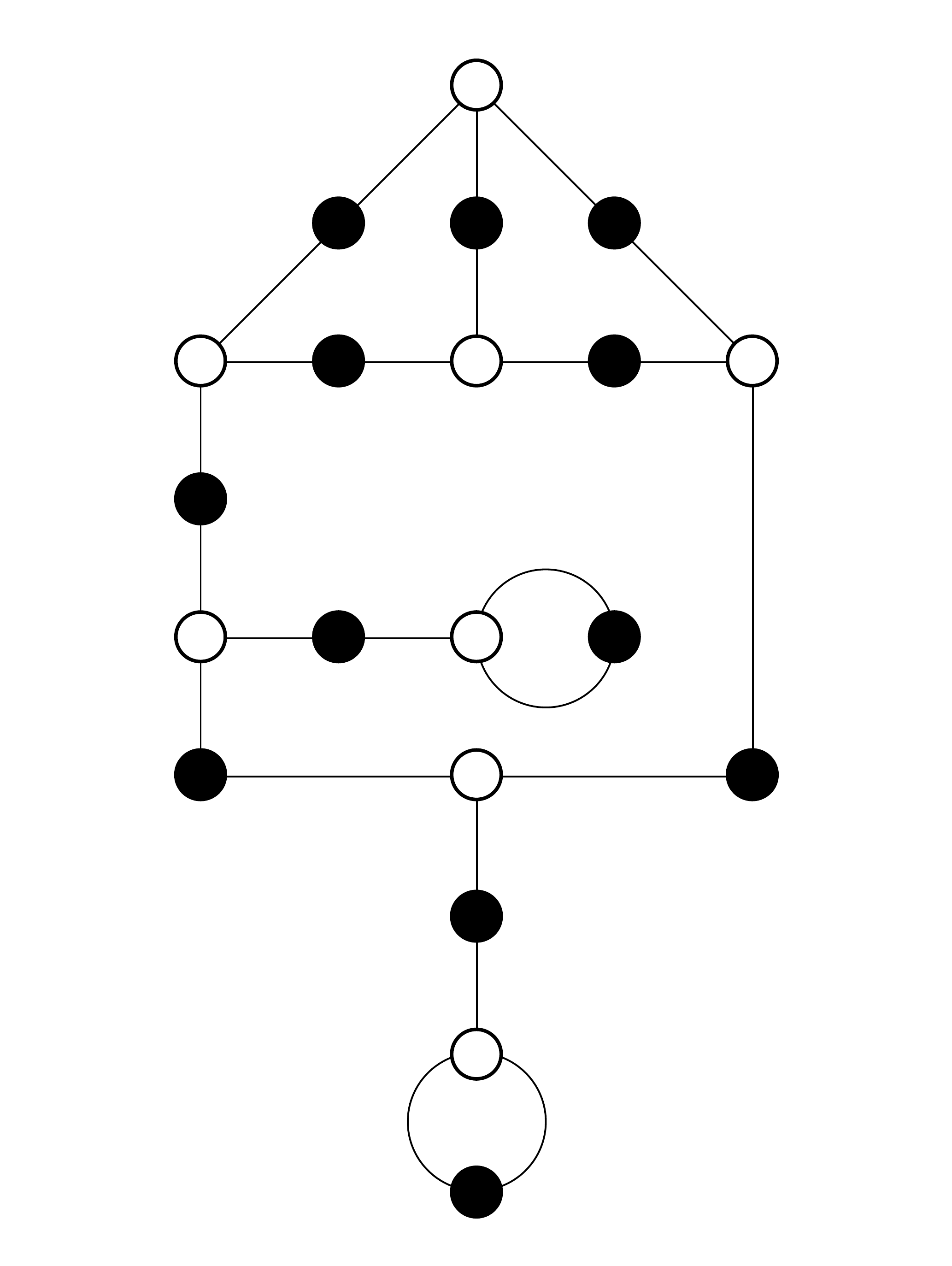}}
\par\end{center}{\scriptsize \par}

\begin{center}
{\scriptsize $8,8,3,3,1,1\;\left(\sqrt{-2}\right)$}
\par\end{center}%
\end{minipage}
\par\end{center}{\scriptsize \par}

\begin{center}
{\scriptsize }%
\begin{minipage}[t]{0.33\textwidth}%
\begin{center}
{\scriptsize \includegraphics[scale=0.15]{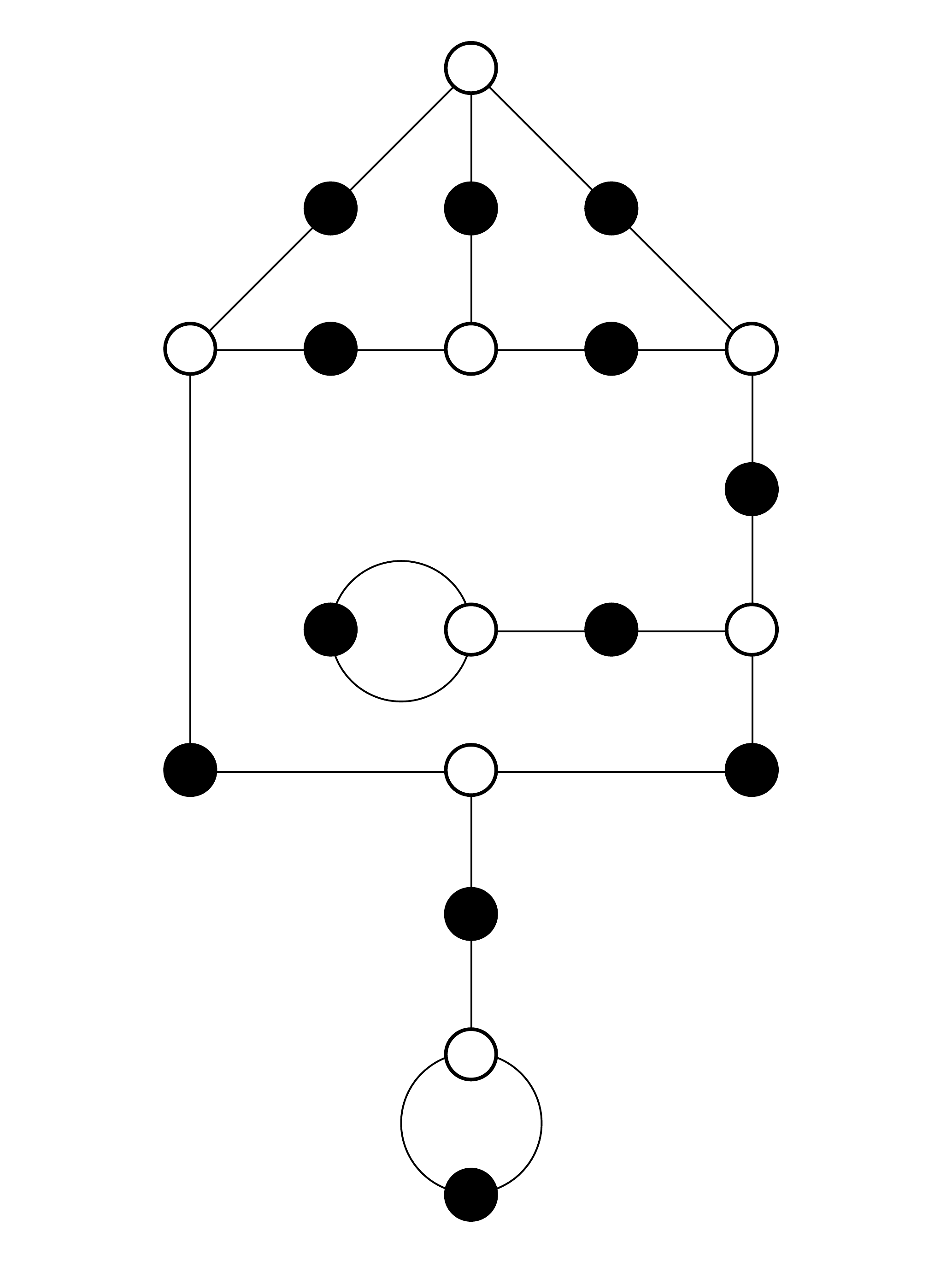}}
\par\end{center}{\scriptsize \par}

\begin{center}
{\scriptsize $8,8,3,3,1,1\;\left(\sqrt{-2}\right)$}
\par\end{center}%
\end{minipage}{\scriptsize }%
\begin{minipage}[t]{0.33\textwidth}%
\begin{center}
{\scriptsize \includegraphics[scale=0.15]{\string"PICT/8-8-2-2-2-2\string".pdf}}
\par\end{center}{\scriptsize \par}

\begin{center}
{\scriptsize $8,8,2,2,2,2\;\left(\mathbb{Q}\right)$}
\par\end{center}%
\end{minipage}{\scriptsize }%
\begin{minipage}[t]{0.33\textwidth}%
\begin{center}
{\scriptsize \includegraphics[scale=0.15]{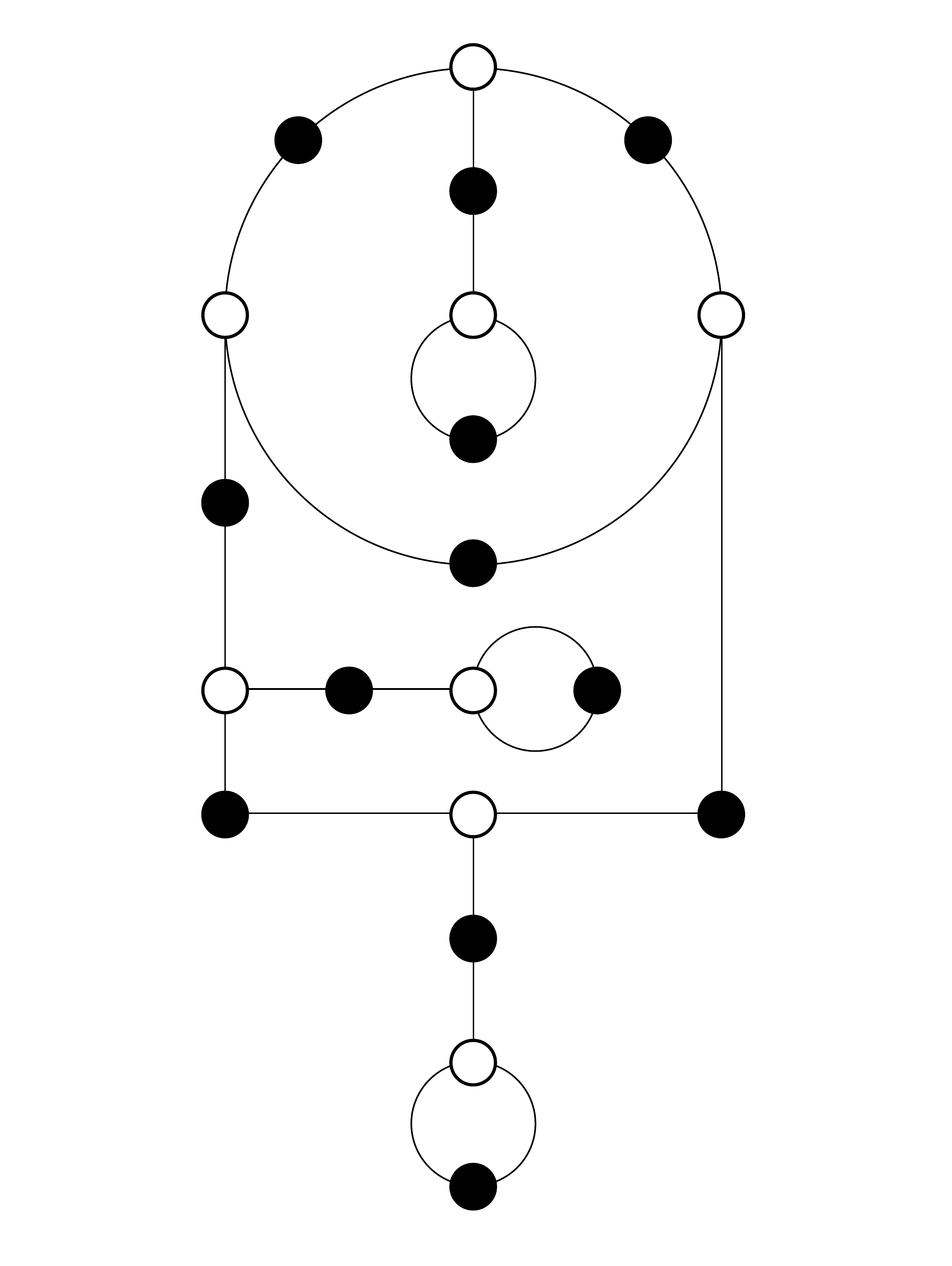}}
\par\end{center}{\scriptsize \par}

\begin{center}
{\scriptsize $8,7,6,1,1,1\;\left(\sqrt{-3}\right)$}
\par\end{center}%
\end{minipage}
\par\end{center}{\scriptsize \par}

\begin{center}
{\scriptsize }%
\begin{minipage}[t]{0.33\textwidth}%
\begin{center}
{\scriptsize \includegraphics[scale=0.15]{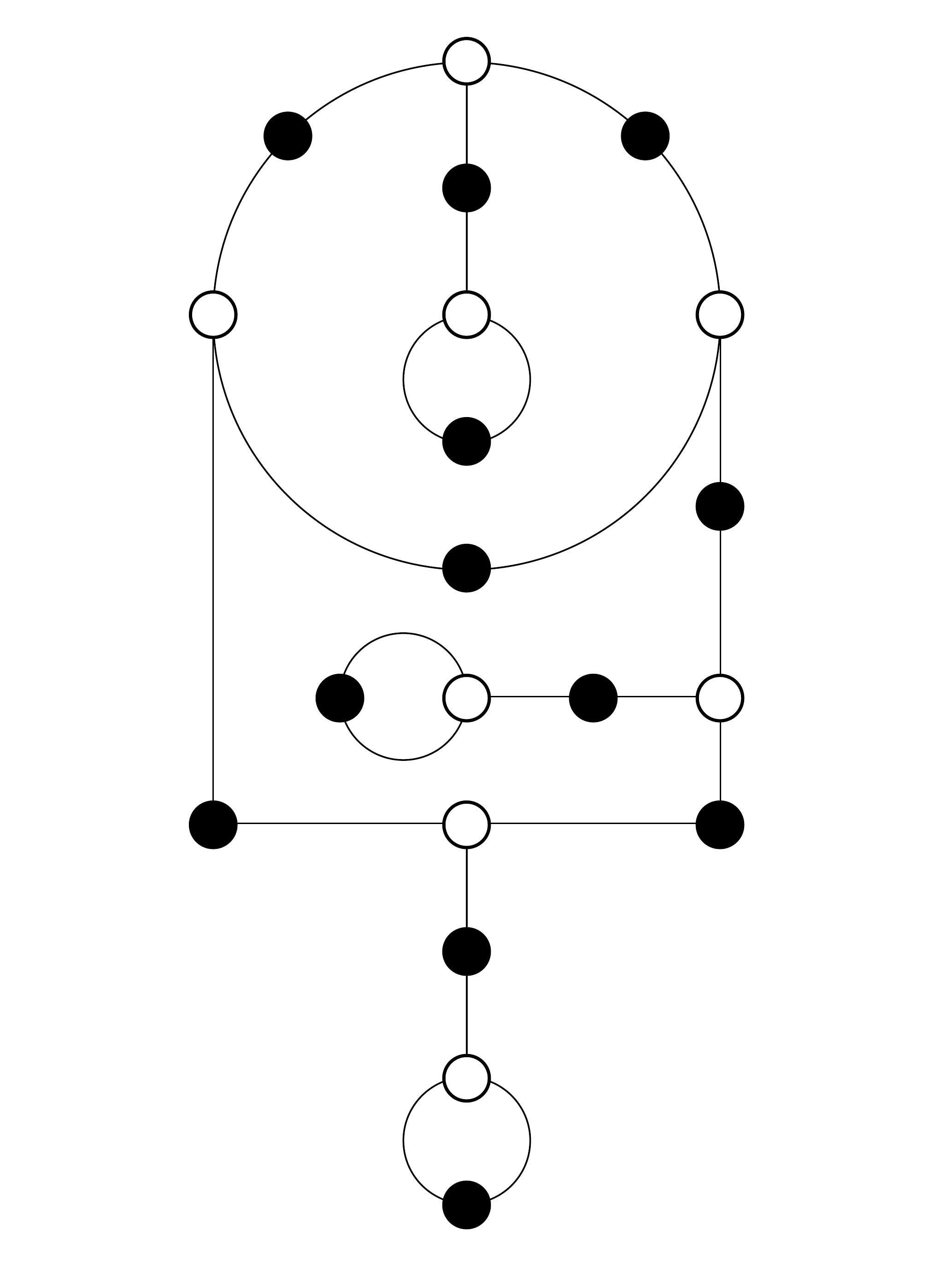}}
\par\end{center}{\scriptsize \par}

\begin{center}
{\scriptsize $8,7,6,1,1,1\;\left(\sqrt{-3}\right)$}
\par\end{center}%
\end{minipage}{\scriptsize }%
\begin{minipage}[t]{0.33\textwidth}%
\begin{center}
{\scriptsize \includegraphics[scale=0.15]{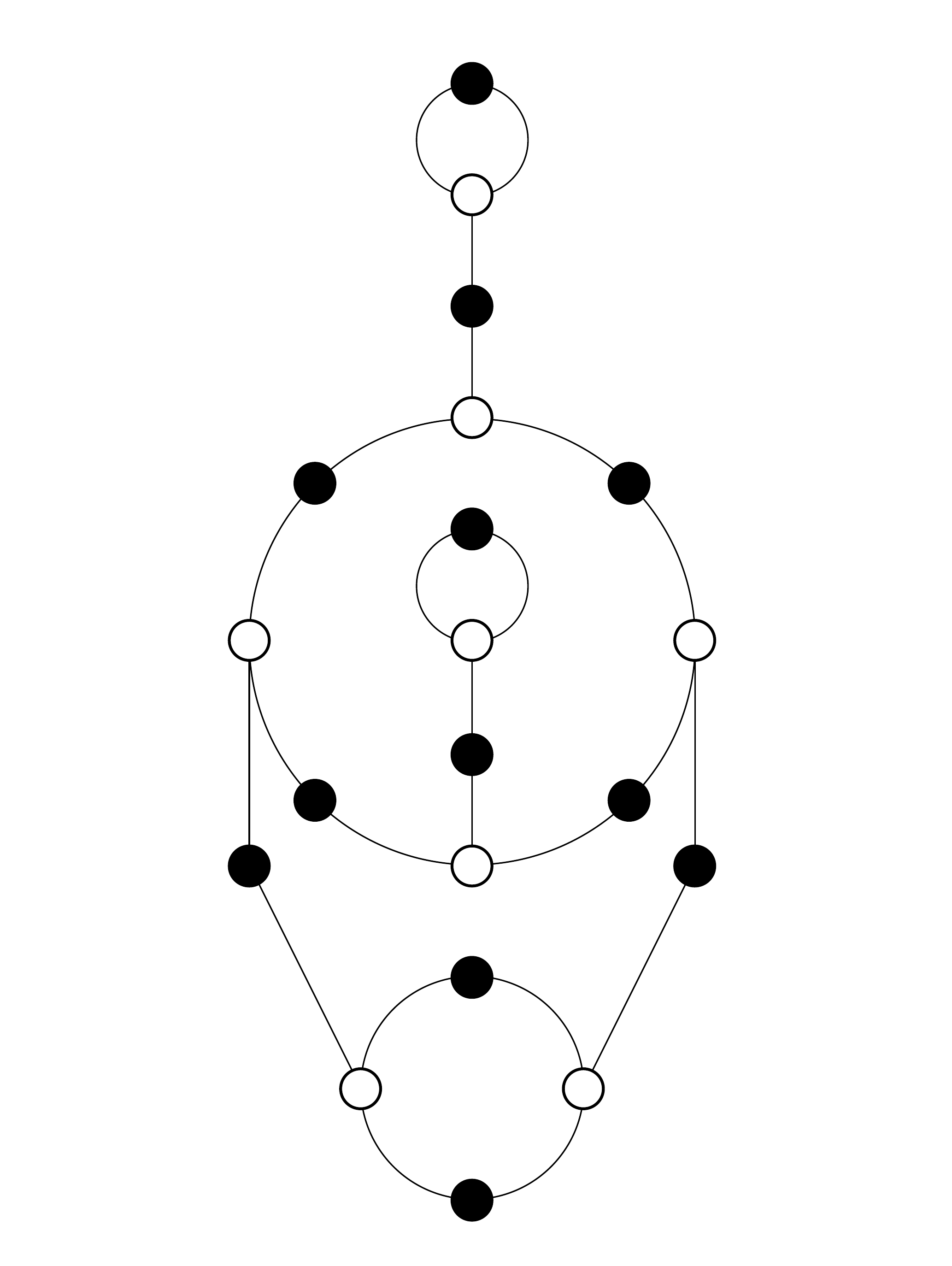}}
\par\end{center}{\scriptsize \par}

\begin{center}
{\scriptsize $8,7,5,2,1,1\;\left(\mathrm{cubic}\right)$}
\par\end{center}%
\end{minipage}{\scriptsize }%
\begin{minipage}[t]{0.33\textwidth}%
\begin{center}
{\scriptsize \includegraphics[scale=0.15]{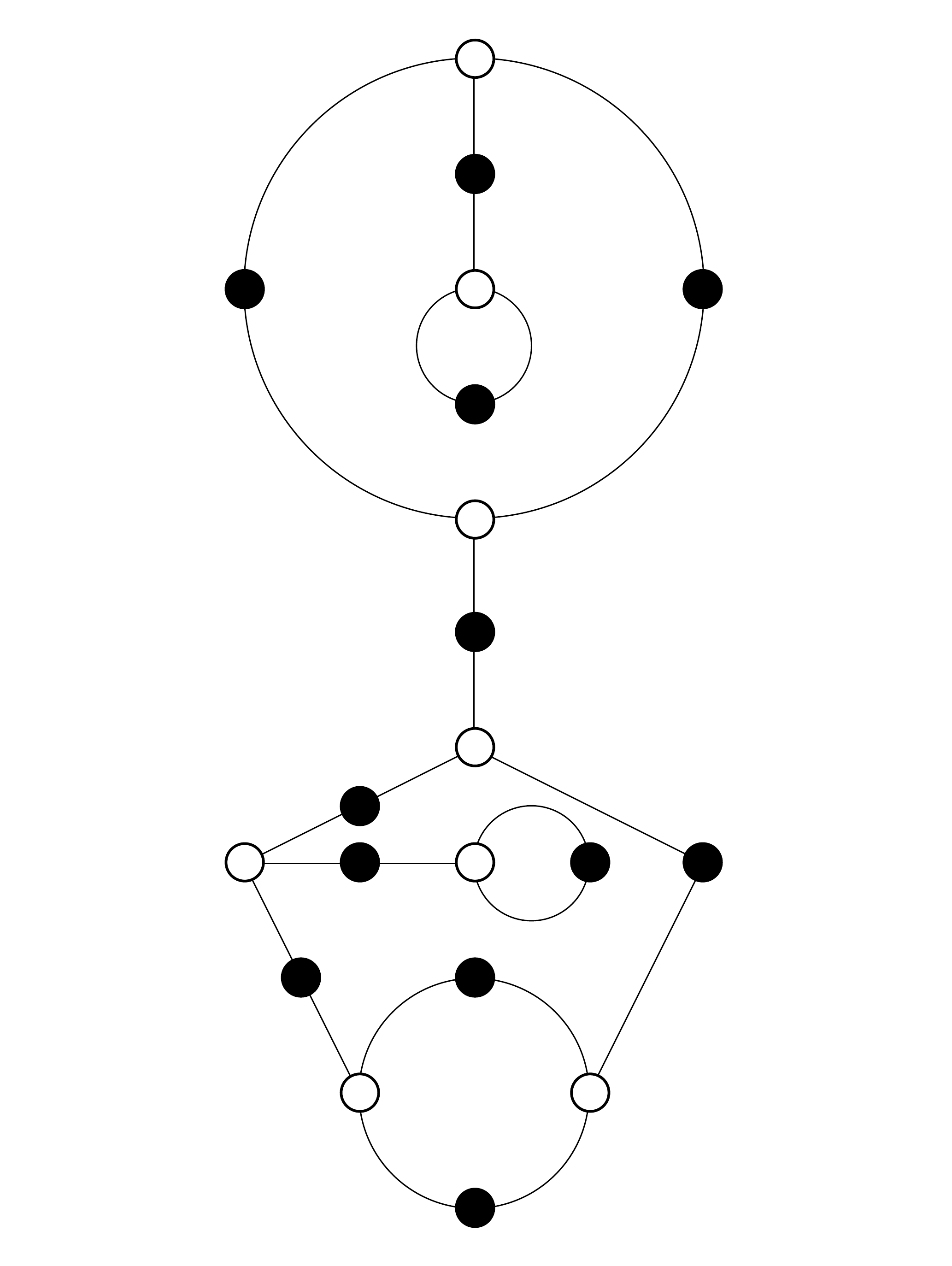}}
\par\end{center}{\scriptsize \par}

\begin{center}
{\scriptsize $8,7,5,2,1,1\;\left(\mathrm{cubic}\right)$}
\par\end{center}%
\end{minipage}
\par\end{center}{\scriptsize \par}

\begin{center}
{\scriptsize }%
\begin{minipage}[t]{0.33\textwidth}%
\begin{center}
{\scriptsize \includegraphics[scale=0.15]{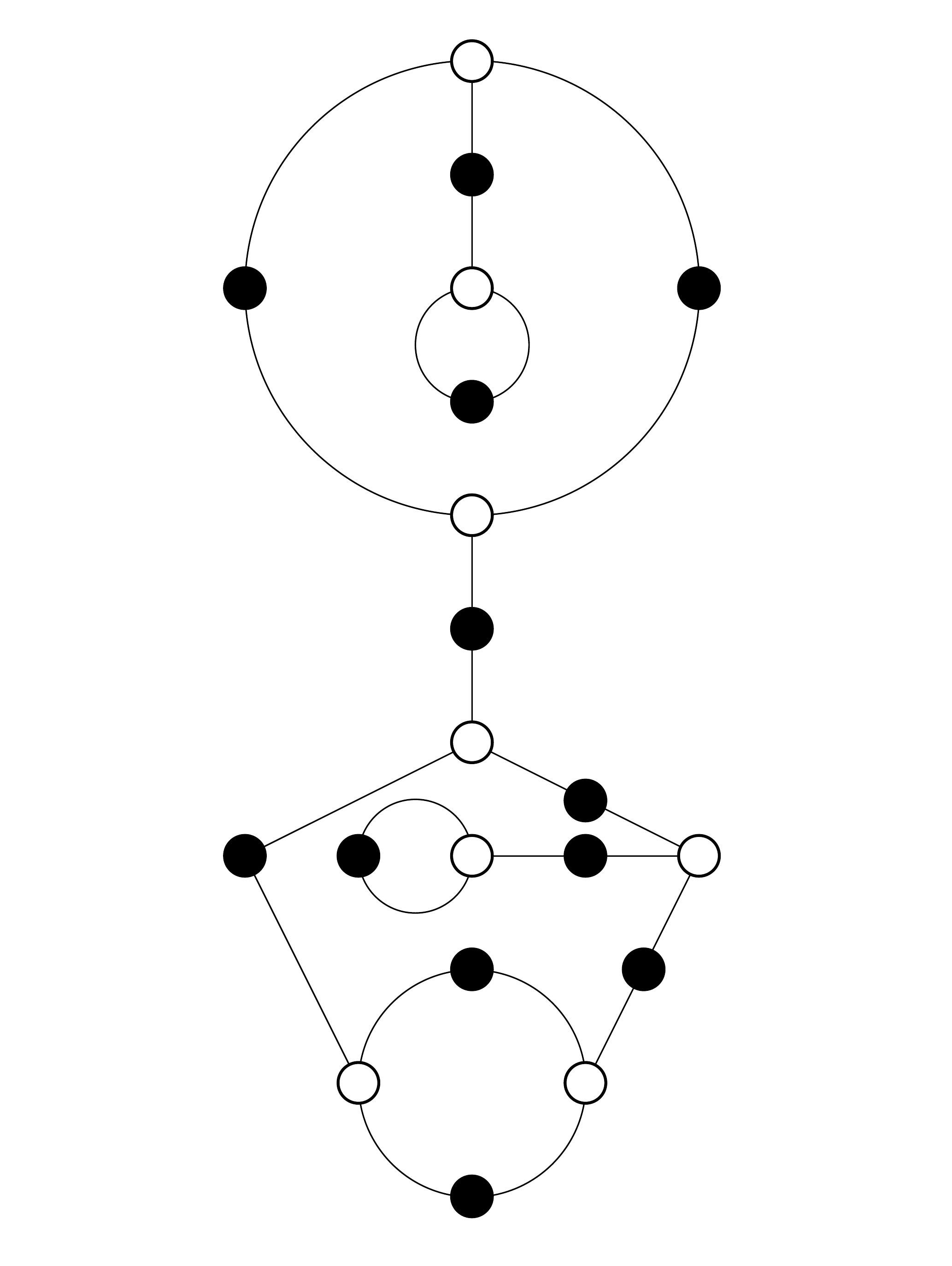}}
\par\end{center}{\scriptsize \par}

\begin{center}
{\scriptsize $8,7,5,2,1,1\;\left(\mathrm{cubic}\right)$}
\par\end{center}%
\end{minipage}{\scriptsize }%
\begin{minipage}[t]{0.33\textwidth}%
\begin{center}
{\scriptsize \includegraphics[scale=0.15]{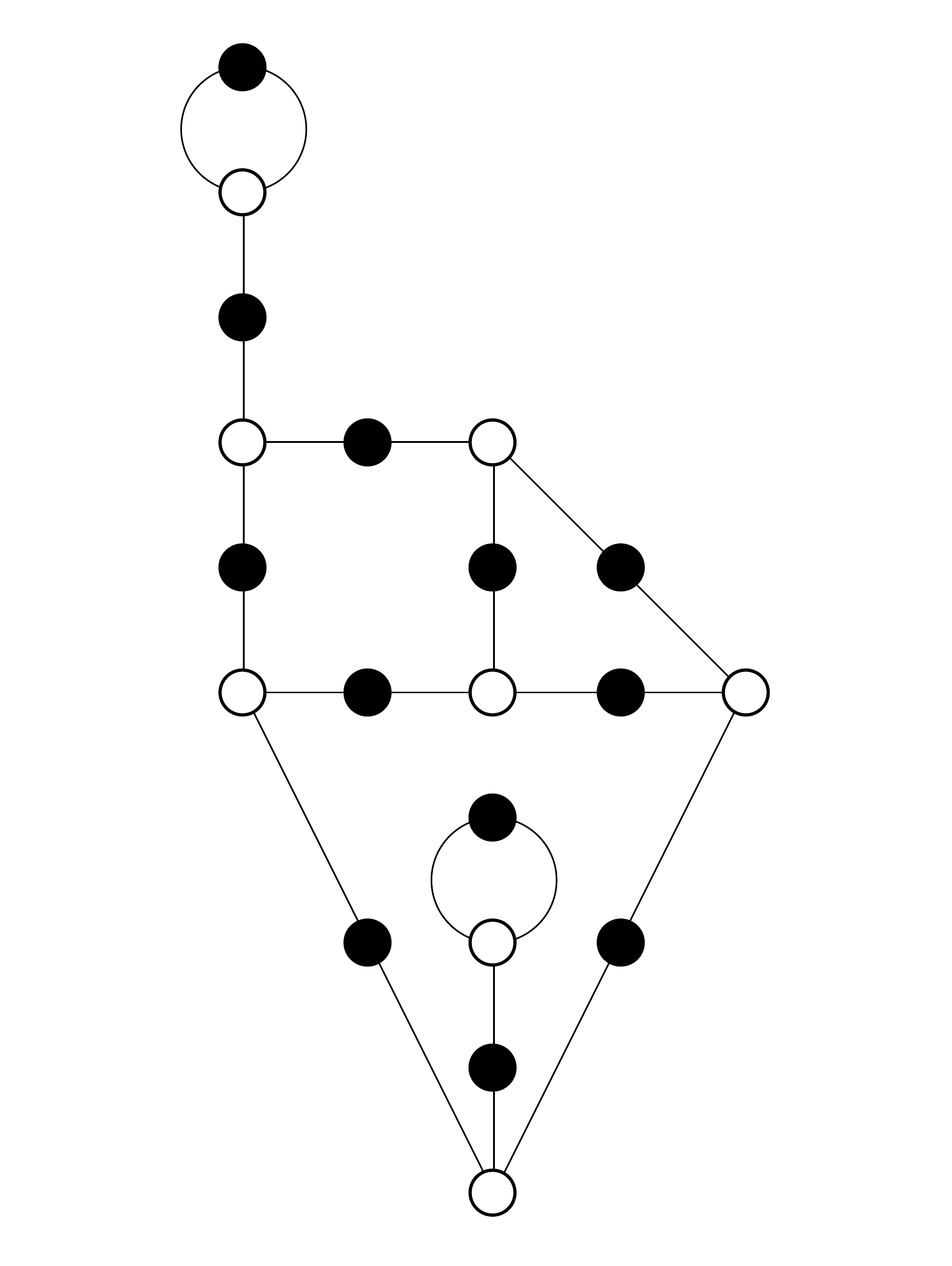}}
\par\end{center}{\scriptsize \par}

\begin{center}
{\scriptsize $8,7,4,3,1,1\;\left(\sqrt{-6}\right)$}
\par\end{center}%
\end{minipage}{\scriptsize }%
\begin{minipage}[t]{0.33\textwidth}%
\begin{center}
{\scriptsize \includegraphics[scale=0.15]{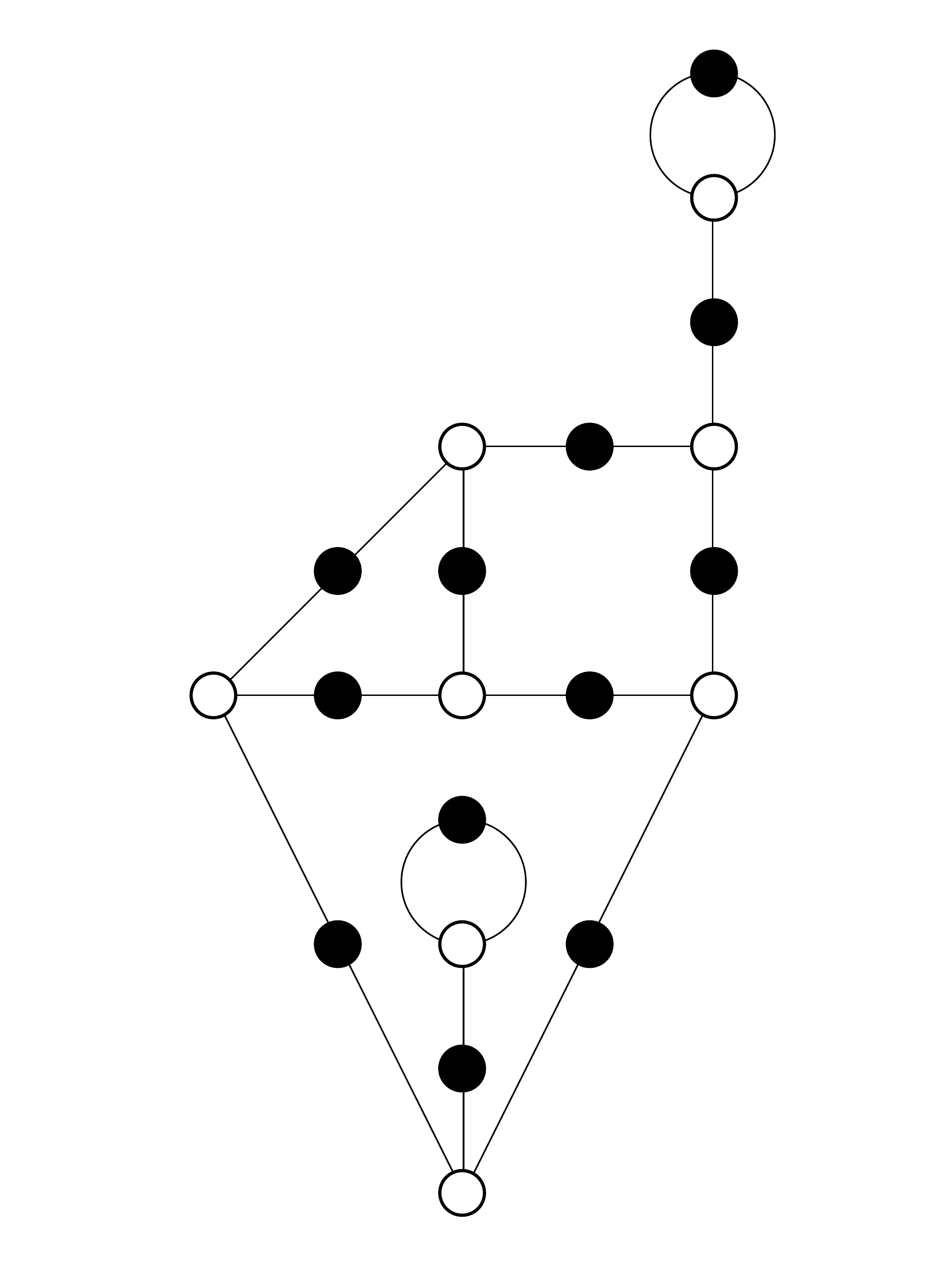}}
\par\end{center}{\scriptsize \par}

\begin{center}
{\scriptsize $8,7,4,3,1,1\;\left(\sqrt{-6}\right)$}
\par\end{center}%
\end{minipage}
\par\end{center}{\scriptsize \par}

\begin{center}
{\scriptsize }%
\begin{minipage}[t]{0.33\textwidth}%
\begin{center}
{\scriptsize \includegraphics[scale=0.15]{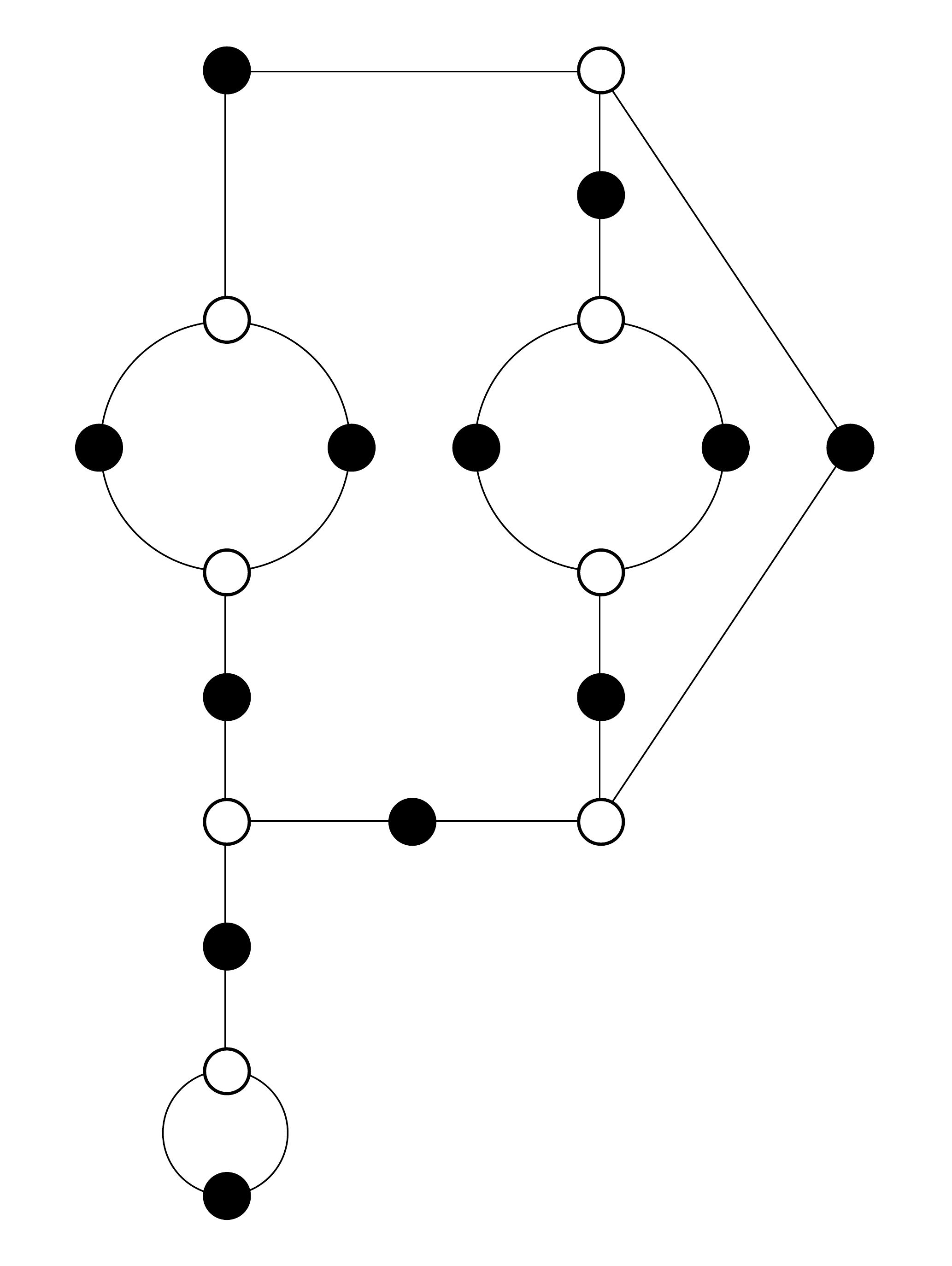}}
\par\end{center}{\scriptsize \par}

\begin{center}
{\scriptsize $8,7,4,2,2,1\;\left(\sqrt{-7}\right)$}
\par\end{center}%
\end{minipage}{\scriptsize }%
\begin{minipage}[t]{0.33\textwidth}%
\begin{center}
{\scriptsize \includegraphics[scale=0.15]{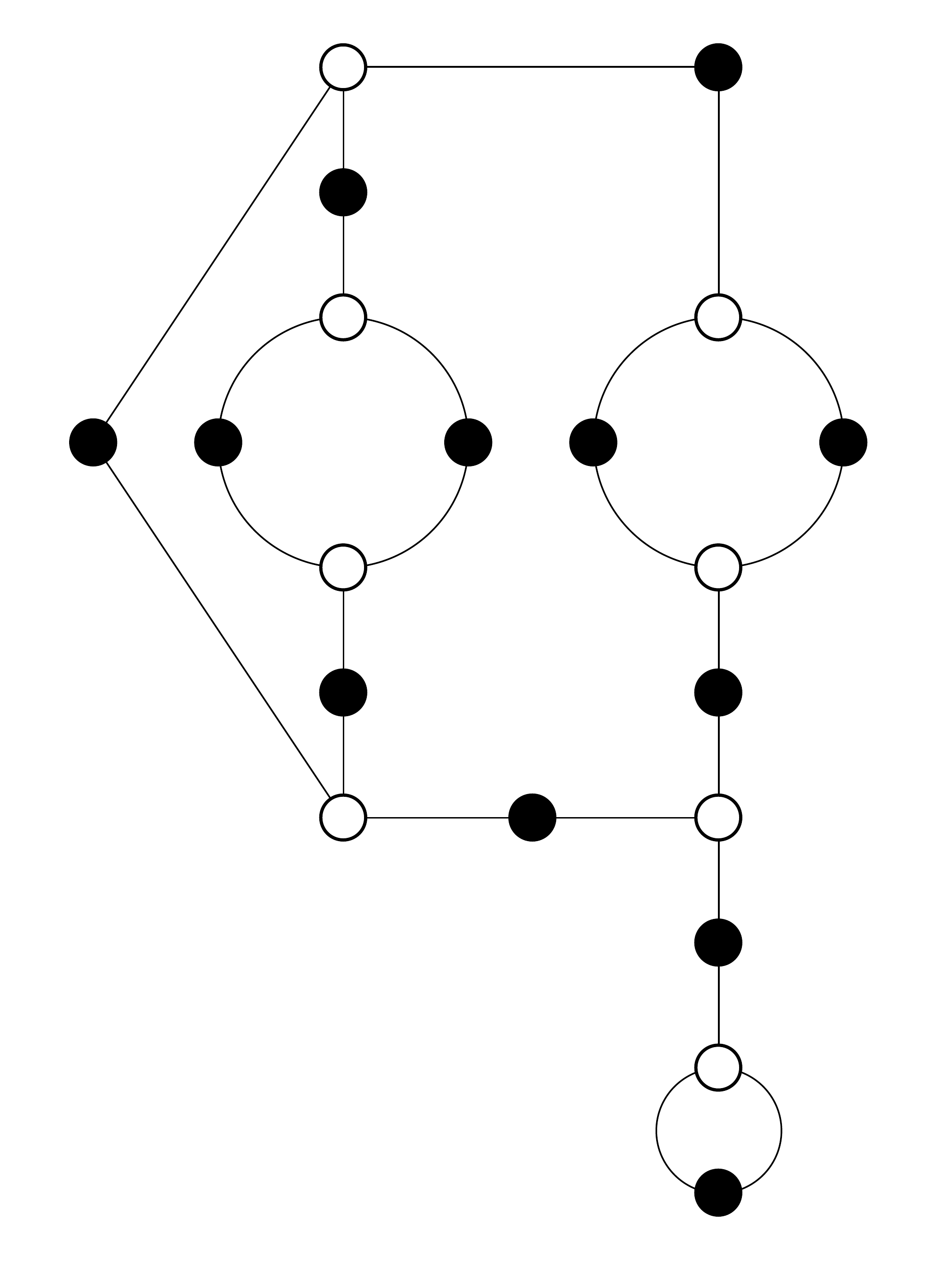}}
\par\end{center}{\scriptsize \par}

\begin{center}
{\scriptsize $8,7,4,2,2,1\;\left(\sqrt{-7}\right)$}
\par\end{center}%
\end{minipage}{\scriptsize }%
\begin{minipage}[t]{0.33\textwidth}%
\begin{center}
{\scriptsize \includegraphics[scale=0.15]{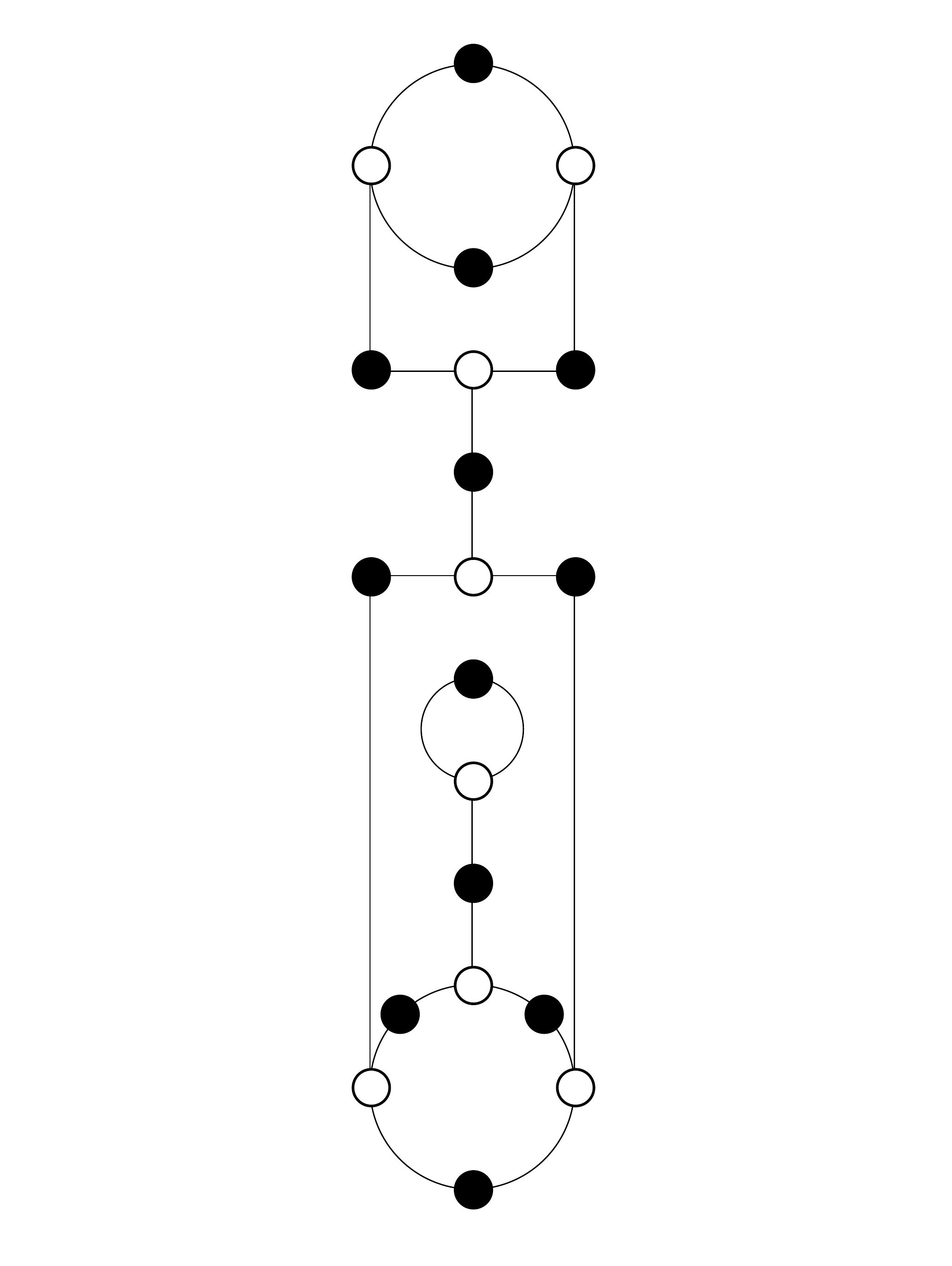}}
\par\end{center}{\scriptsize \par}

\begin{center}
{\scriptsize $8,7,3,3,2,1\;\left(\mathbb{Q}\right)$}
\par\end{center}%
\end{minipage}
\par\end{center}{\scriptsize \par}

\begin{center}
{\scriptsize }%
\begin{minipage}[t]{0.33\textwidth}%
\begin{center}
{\scriptsize \includegraphics[scale=0.15]{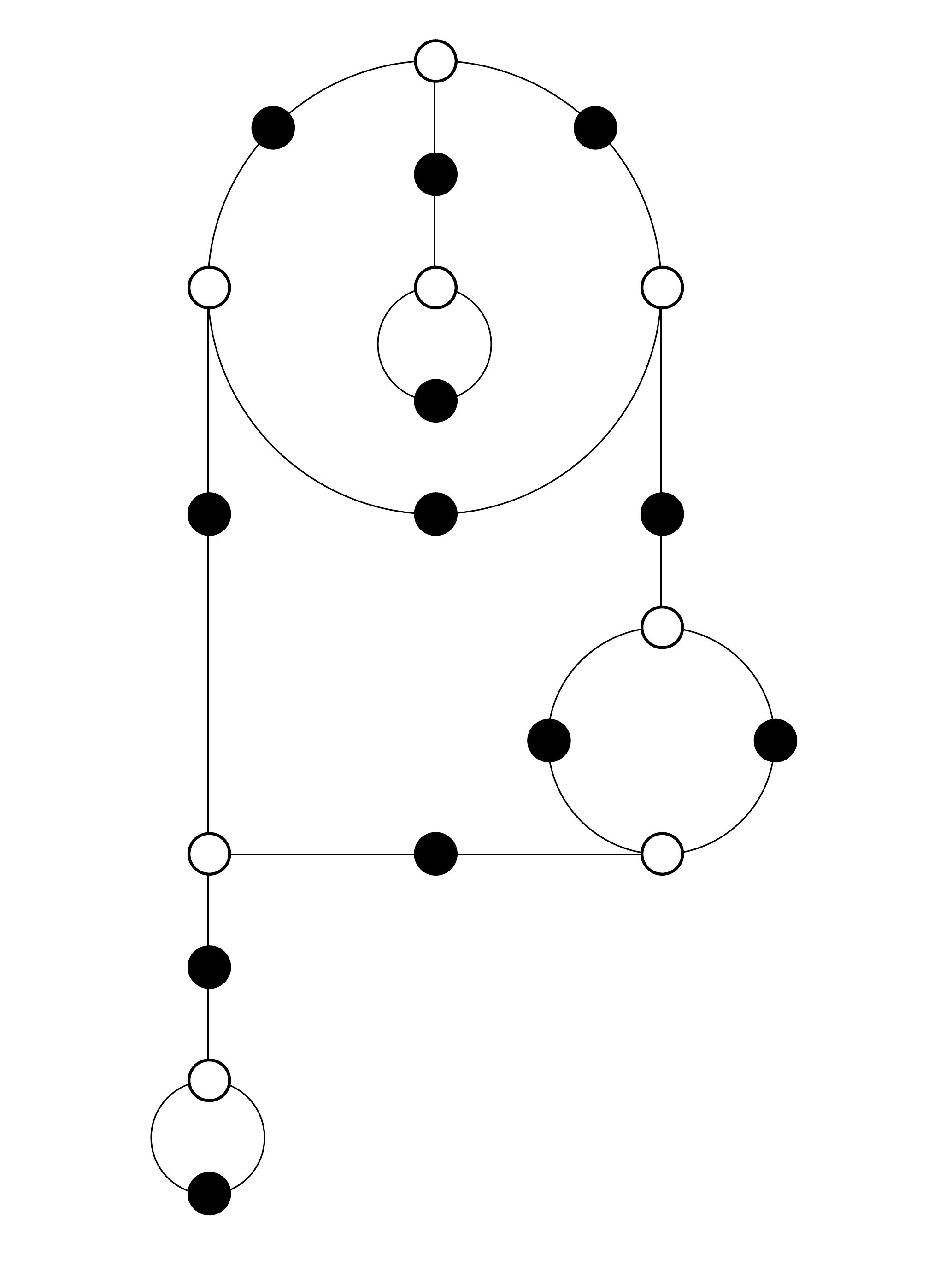}}
\par\end{center}{\scriptsize \par}

\begin{center}
{\scriptsize $8,6,6,2,1,1\;\left(\sqrt{-3}\right)$}
\par\end{center}%
\end{minipage}{\scriptsize }%
\begin{minipage}[t]{0.33\textwidth}%
\begin{center}
{\scriptsize \includegraphics[scale=0.15]{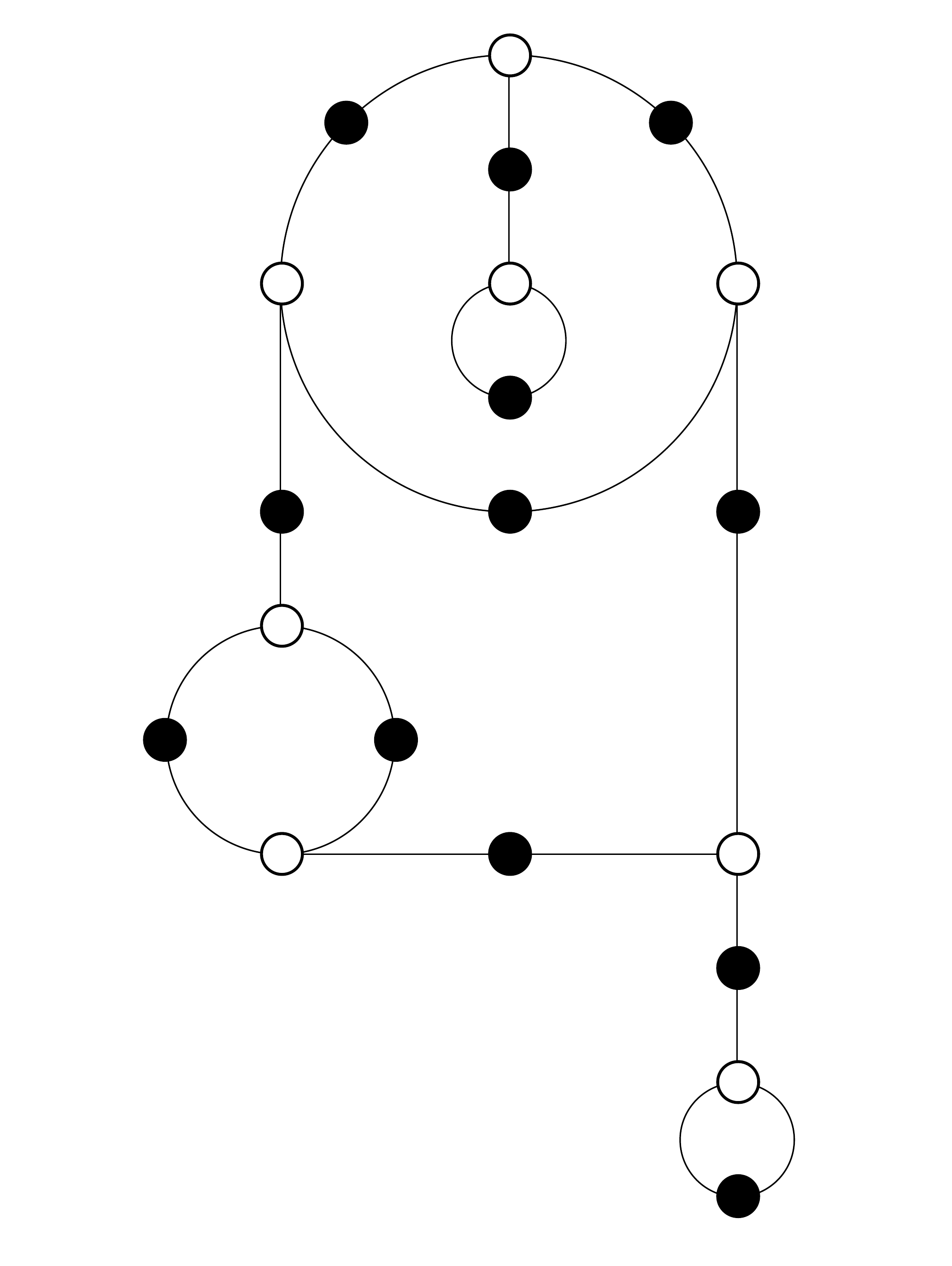}}
\par\end{center}{\scriptsize \par}

\begin{center}
{\scriptsize $8,6,6,2,1,1\;\left(\sqrt{-3}\right)$}
\par\end{center}%
\end{minipage}{\scriptsize }%
\begin{minipage}[t]{0.33\textwidth}%
\begin{center}
{\scriptsize \includegraphics[scale=0.15]{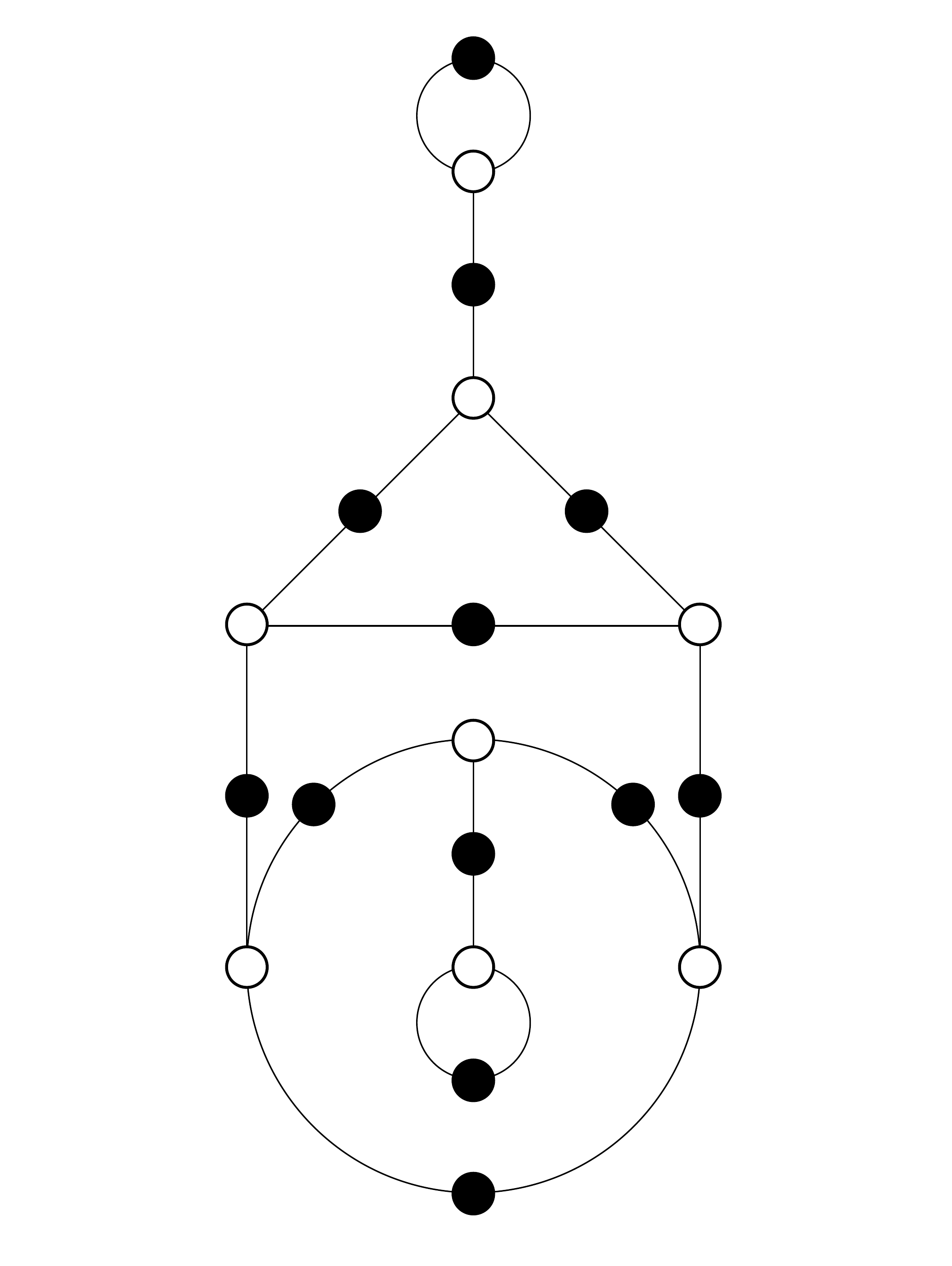}}
\par\end{center}{\scriptsize \par}

\begin{center}
{\scriptsize $8,6,5,3,1,1\;\left(\sqrt{5}\right)$}
\par\end{center}%
\end{minipage}
\par\end{center}{\scriptsize \par}

\begin{center}
{\scriptsize }%
\begin{minipage}[t]{0.33\textwidth}%
\begin{center}
{\scriptsize \includegraphics[scale=0.15]{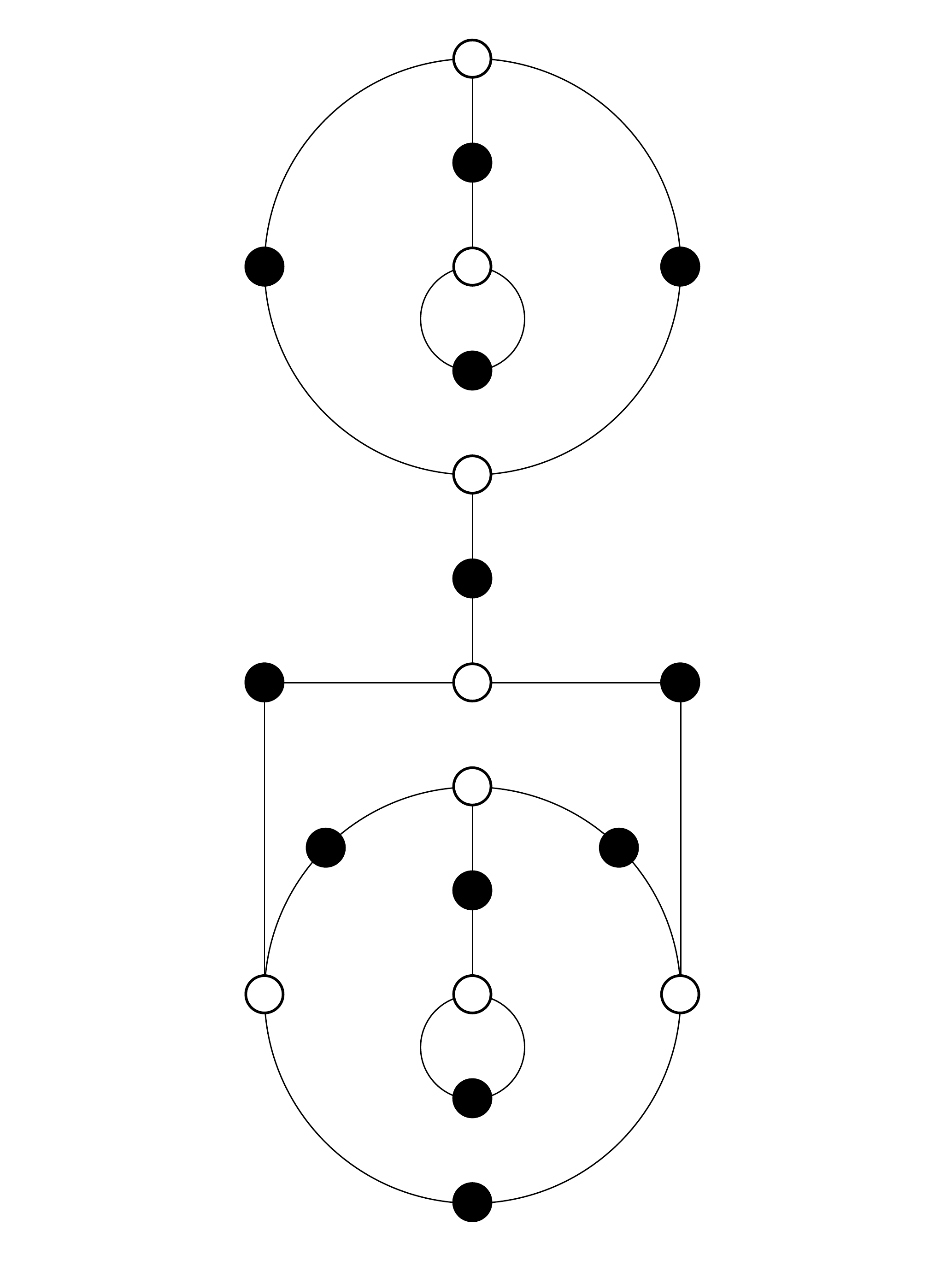}}
\par\end{center}{\scriptsize \par}

\begin{center}
{\scriptsize $8,6,5,3,1,1\;\left(\sqrt{5}\right)$}
\par\end{center}%
\end{minipage}{\scriptsize }%
\begin{minipage}[t]{0.33\textwidth}%
\begin{center}
{\scriptsize \includegraphics[scale=0.15]{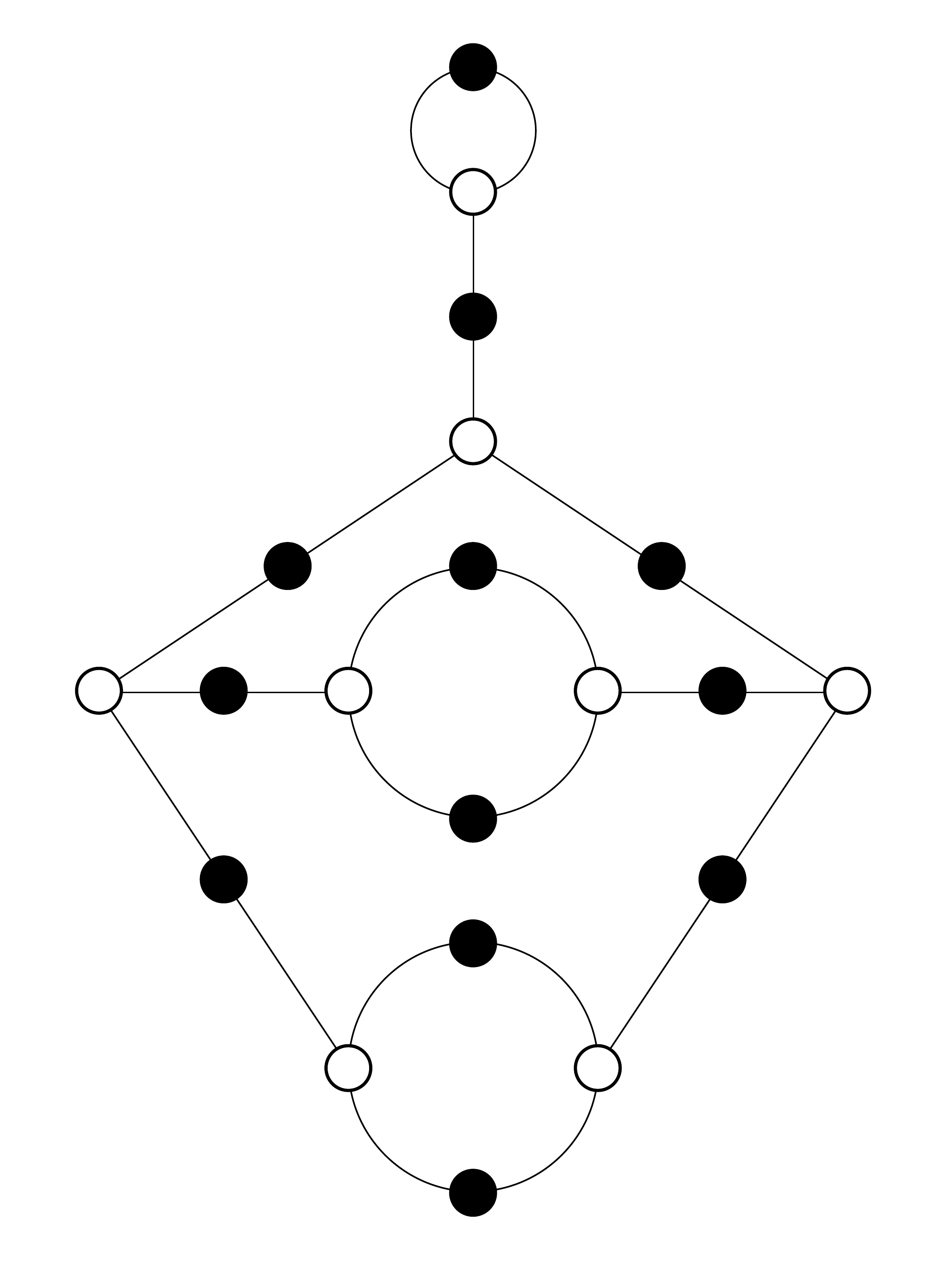}}
\par\end{center}{\scriptsize \par}

\begin{center}
{\scriptsize $8,6,5,2,2,1\;\left(\mathbb{Q}\right)$}
\par\end{center}%
\end{minipage}{\scriptsize }%
\begin{minipage}[t]{0.33\textwidth}%
\begin{center}
{\scriptsize \includegraphics[scale=0.15]{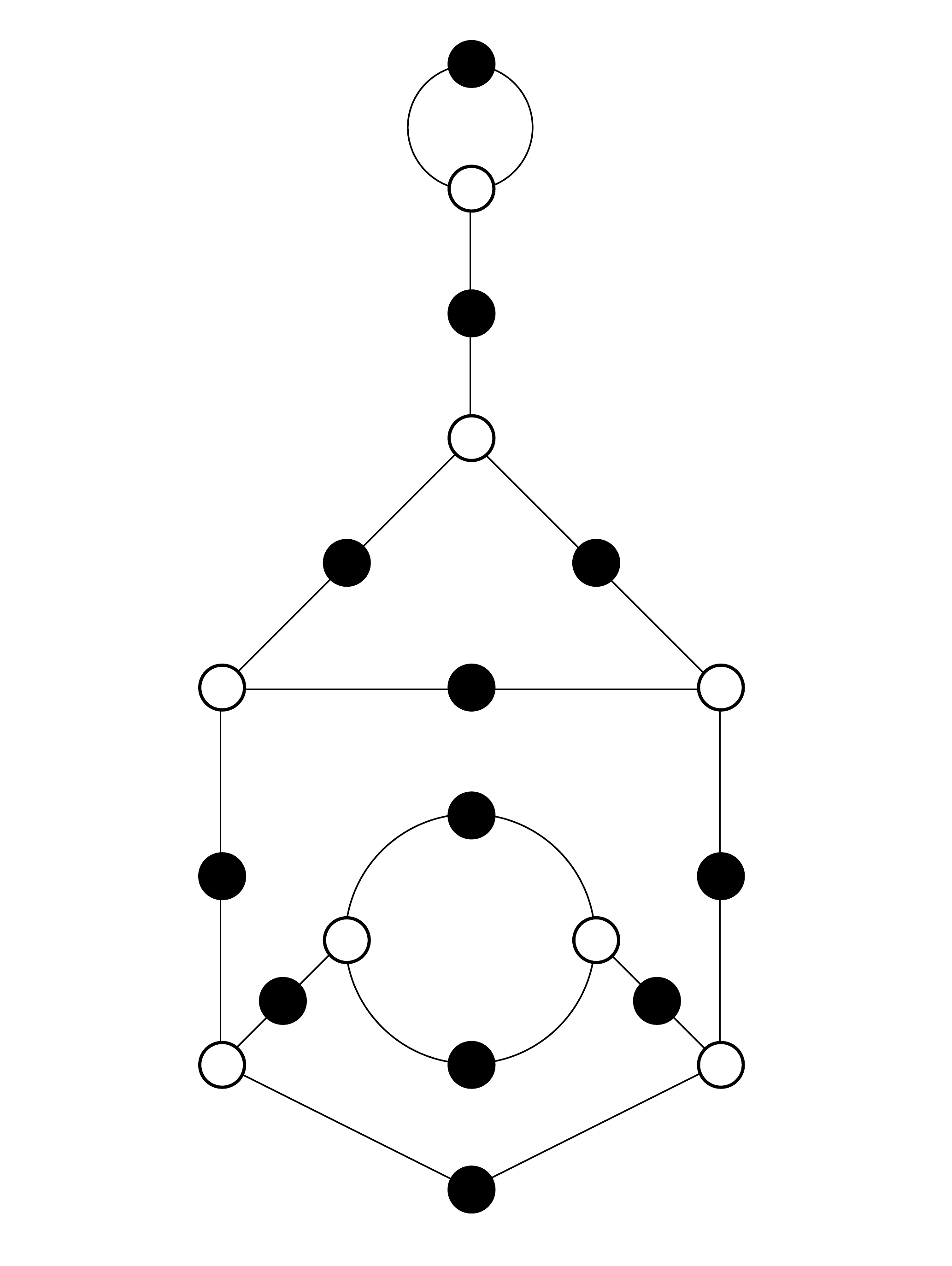}}
\par\end{center}{\scriptsize \par}

\begin{center}
{\scriptsize $8,6,4,3,2,1\;\left(\sqrt{2}\right)$}
\par\end{center}%
\end{minipage}
\par\end{center}{\scriptsize \par}

\begin{center}
{\scriptsize }%
\begin{minipage}[t]{0.33\textwidth}%
\begin{center}
{\scriptsize \includegraphics[scale=0.15]{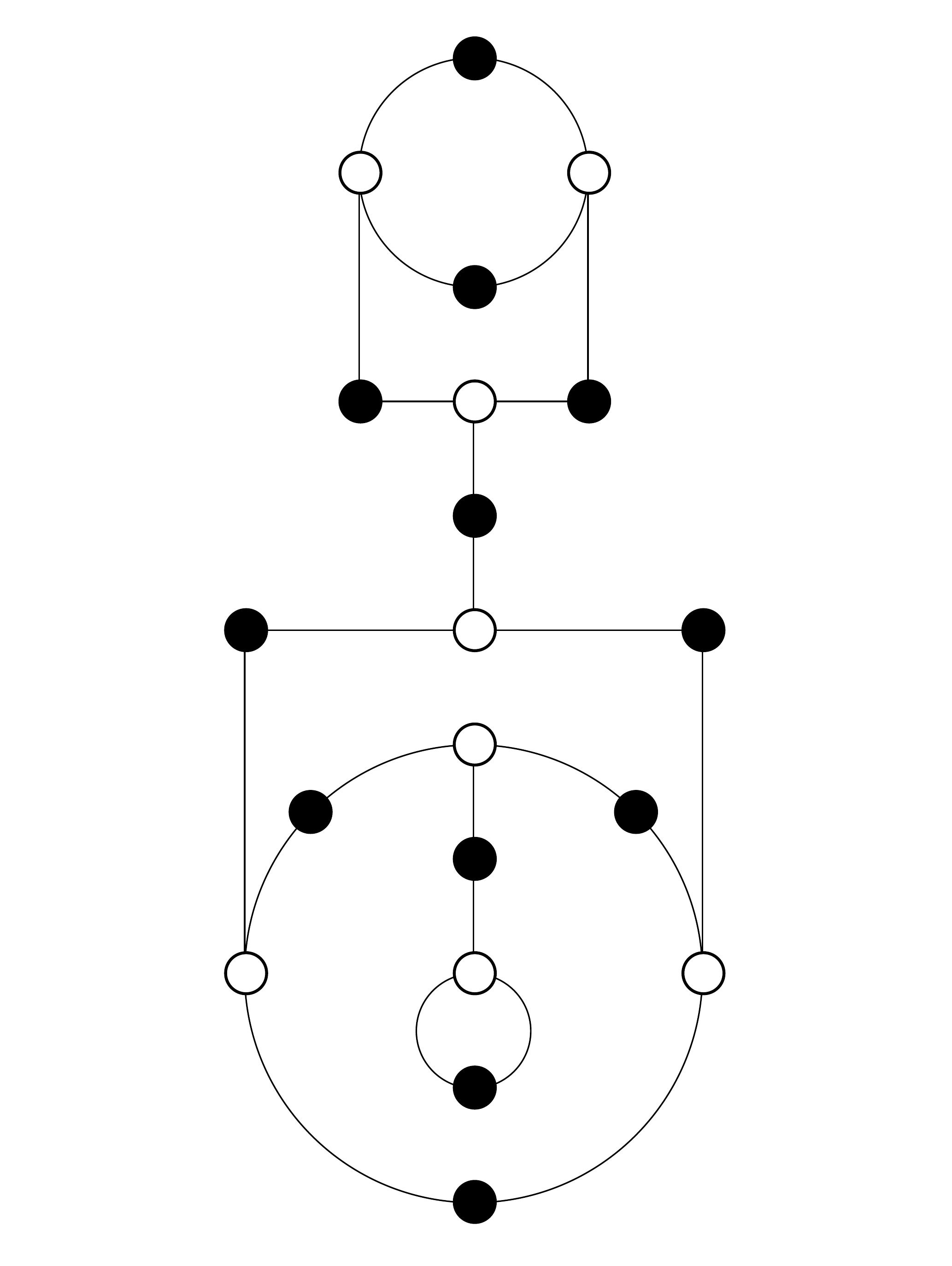}}
\par\end{center}{\scriptsize \par}

\begin{center}
{\scriptsize $8,6,4,3,2,1\;\left(\sqrt{2}\right)$}
\par\end{center}%
\end{minipage}{\scriptsize }%
\begin{minipage}[t]{0.33\textwidth}%
\begin{center}
{\scriptsize \includegraphics[scale=0.15]{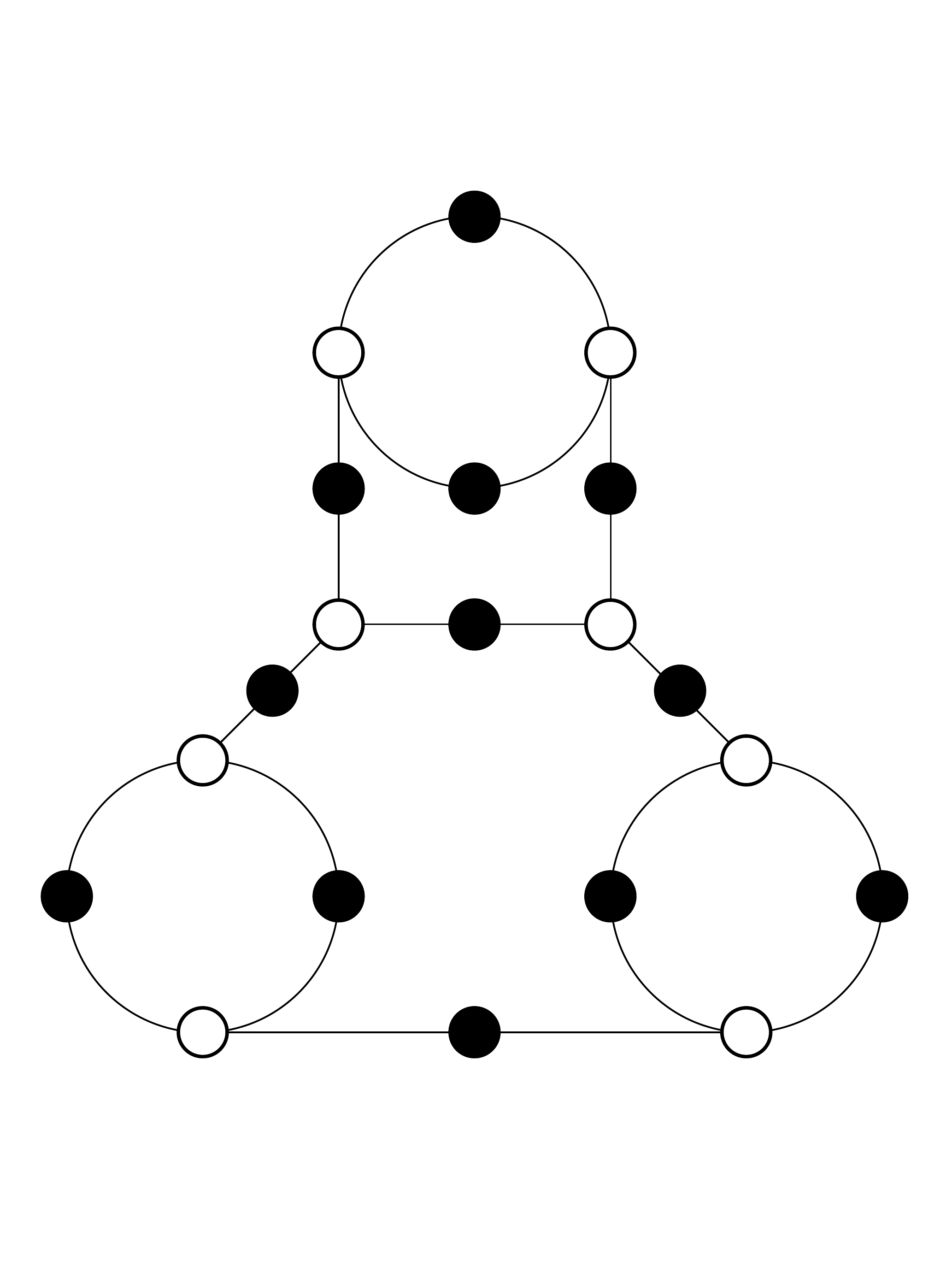}}
\par\end{center}{\scriptsize \par}

\begin{center}
{\scriptsize $8,6,4,2,2,2\;\left(\mathbb{Q}\right)$}
\par\end{center}%
\end{minipage}{\scriptsize }%
\begin{minipage}[t]{0.33\textwidth}%
\begin{center}
{\scriptsize \includegraphics[scale=0.15]{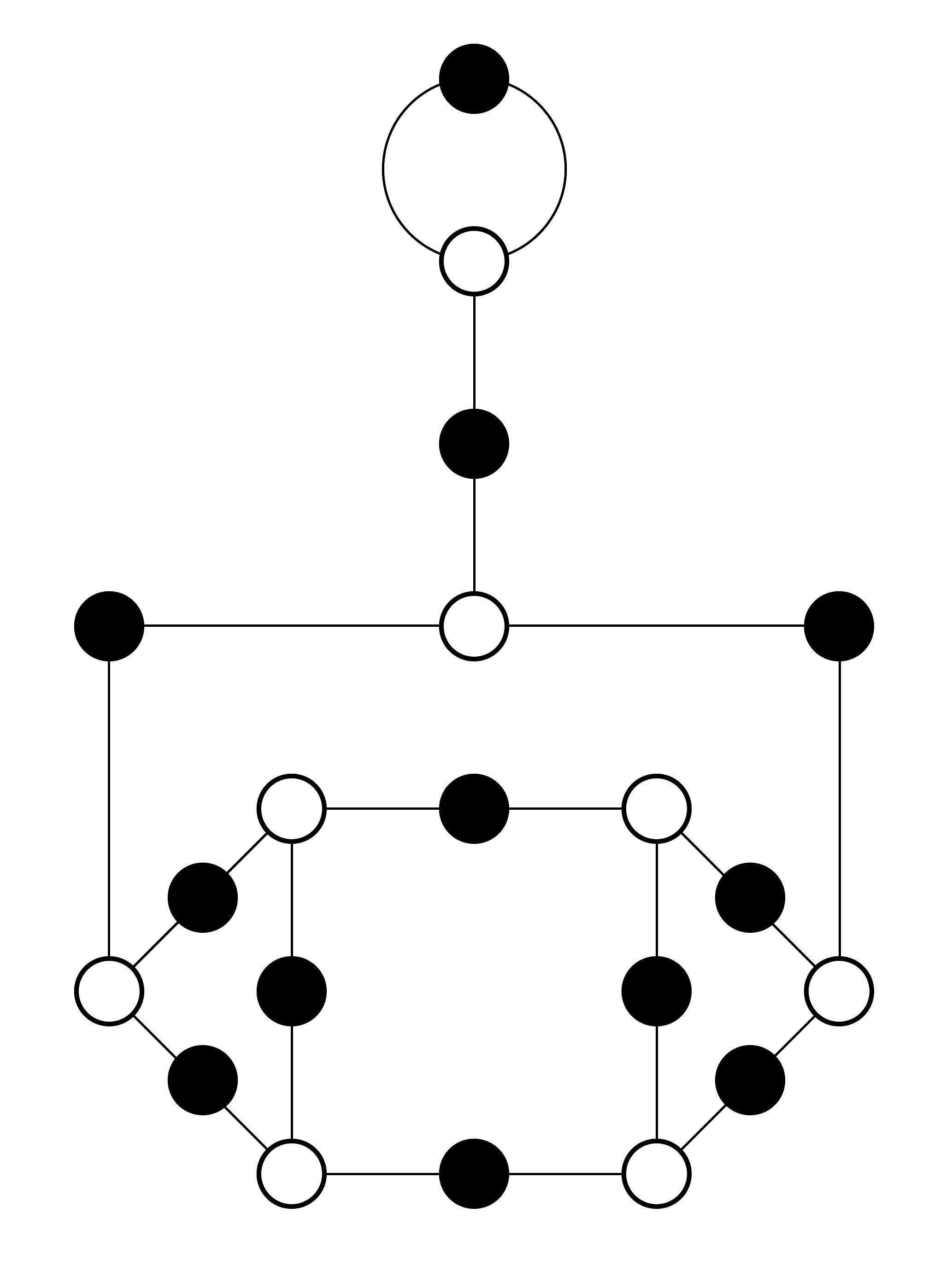}}
\par\end{center}{\scriptsize \par}

\begin{center}
{\scriptsize $8,5,4,3,3,1\;\left(\sqrt{10}\right)$}
\par\end{center}%
\end{minipage}
\par\end{center}{\scriptsize \par}

\begin{center}
{\scriptsize }%
\begin{minipage}[t]{0.33\textwidth}%
\begin{center}
{\scriptsize \includegraphics[scale=0.15]{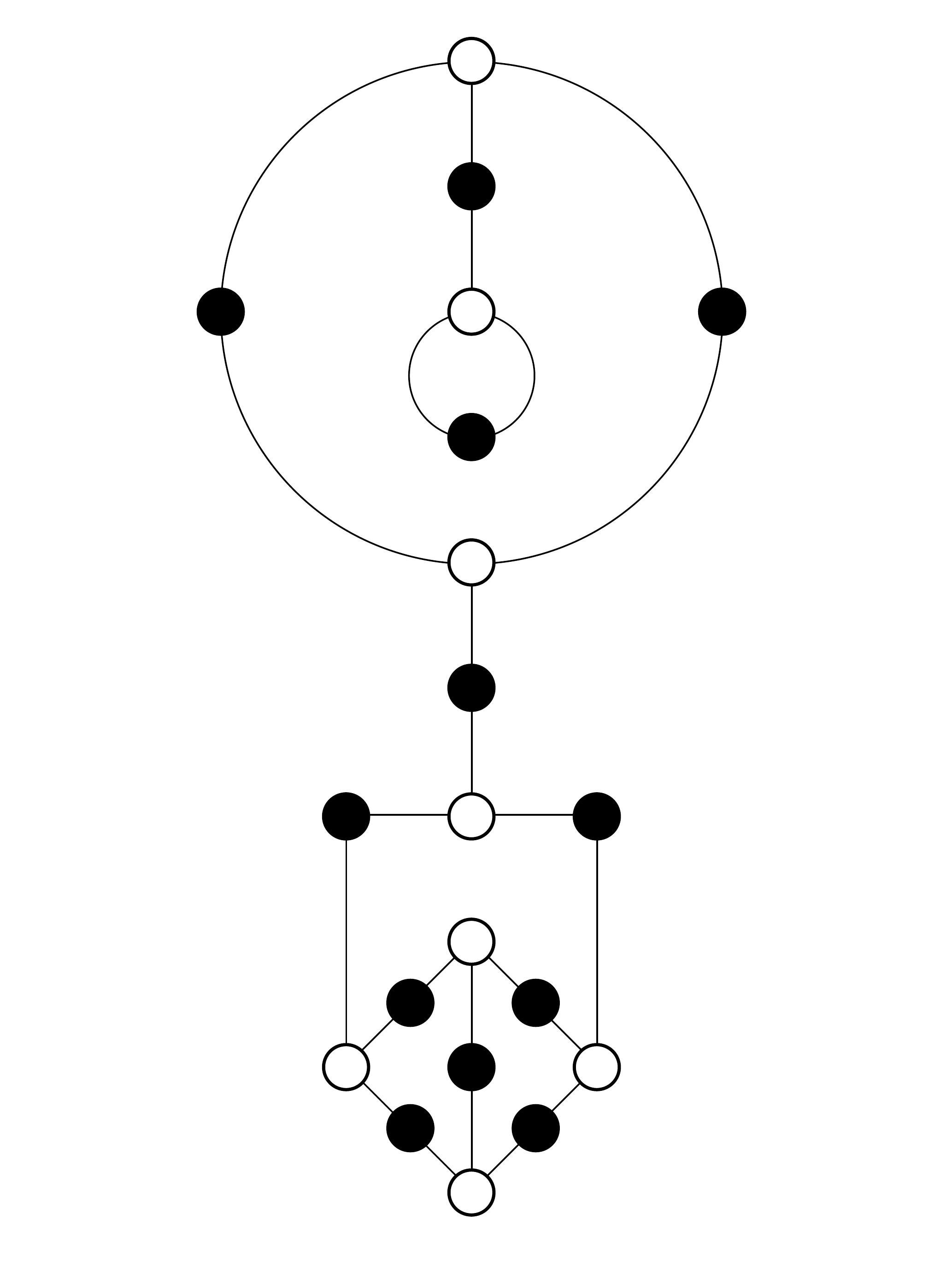}}
\par\end{center}{\scriptsize \par}

\begin{center}
{\scriptsize $8,5,4,3,3,1\;\left(\sqrt{10}\right)$}
\par\end{center}%
\end{minipage}{\scriptsize }%
\begin{minipage}[t]{0.33\textwidth}%
\begin{center}
{\scriptsize \includegraphics[scale=0.15]{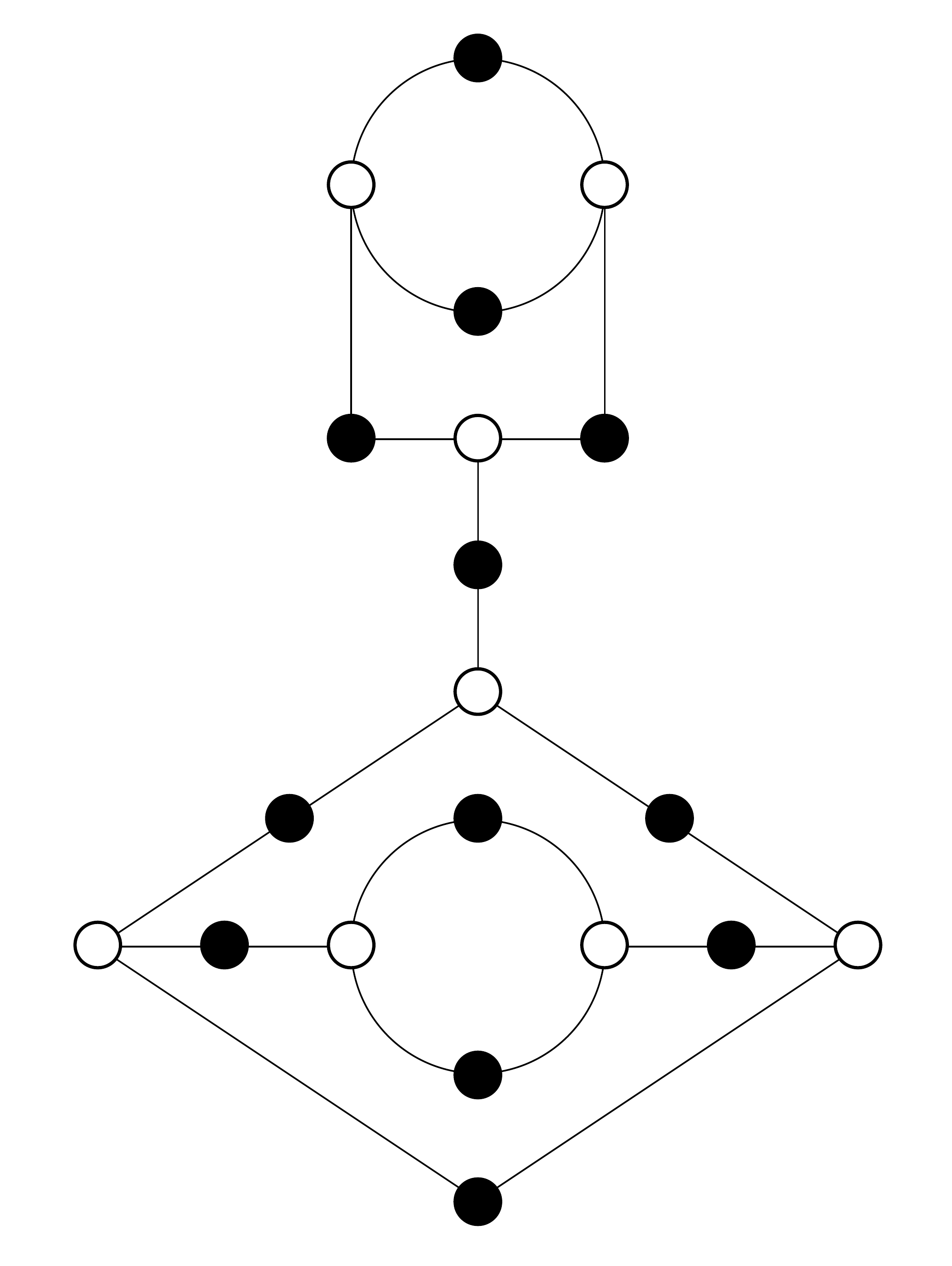}}
\par\end{center}{\scriptsize \par}

\begin{center}
{\scriptsize $8,5,4,3,2,2\;\left(\mathbb{Q}\right)$}
\par\end{center}%
\end{minipage}{\scriptsize }%
\begin{minipage}[t]{0.33\textwidth}%
\begin{center}
{\scriptsize \includegraphics[scale=0.15]{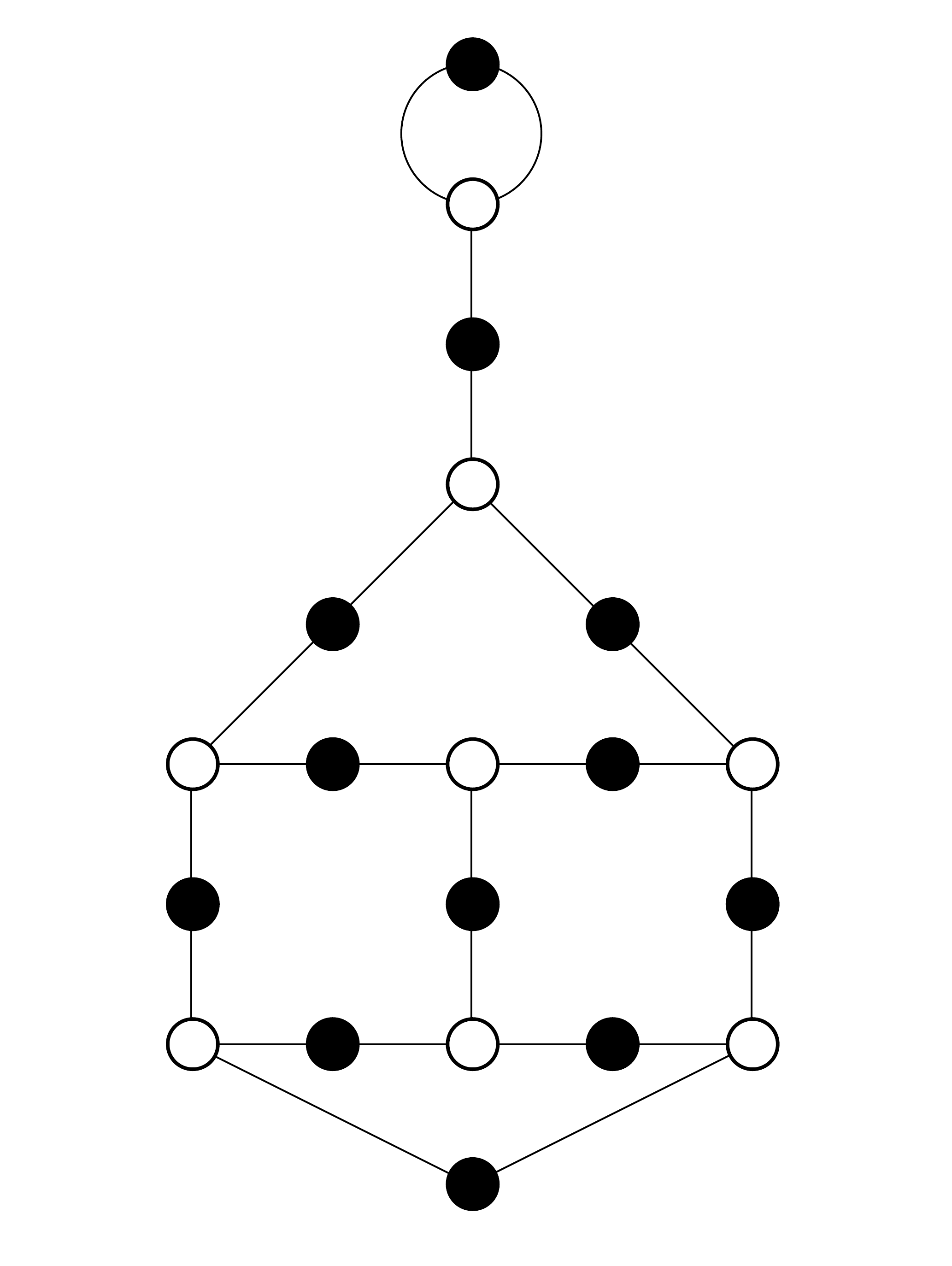}}
\par\end{center}{\scriptsize \par}

\begin{center}
{\scriptsize $8,4,4,4,3,1\;\left(\mathbb{Q}\right)$}
\par\end{center}%
\end{minipage}
\par\end{center}{\scriptsize \par}

\begin{center}
{\scriptsize }%
\begin{minipage}[t]{0.33\textwidth}%
\begin{center}
{\scriptsize \includegraphics[scale=0.15]{\string"PICT/8-4-4-4-2-2\string".pdf}}
\par\end{center}{\scriptsize \par}

\begin{center}
{\scriptsize $8,4,4,4,2,2\;\left(\mathbb{Q}\right)$}
\par\end{center}%
\end{minipage}{\scriptsize }%
\begin{minipage}[t]{0.33\textwidth}%
\begin{center}
{\scriptsize \includegraphics[scale=0.15]{\string"PICT/7-7-7-1-1-1\string".pdf}}
\par\end{center}{\scriptsize \par}

\begin{center}
{\scriptsize $7,7,7,1,1,1\;\left(\mathbb{Q}\right)$}
\par\end{center}%
\end{minipage}{\scriptsize }%
\begin{minipage}[t]{0.33\textwidth}%
\begin{center}
{\scriptsize \includegraphics[scale=0.15]{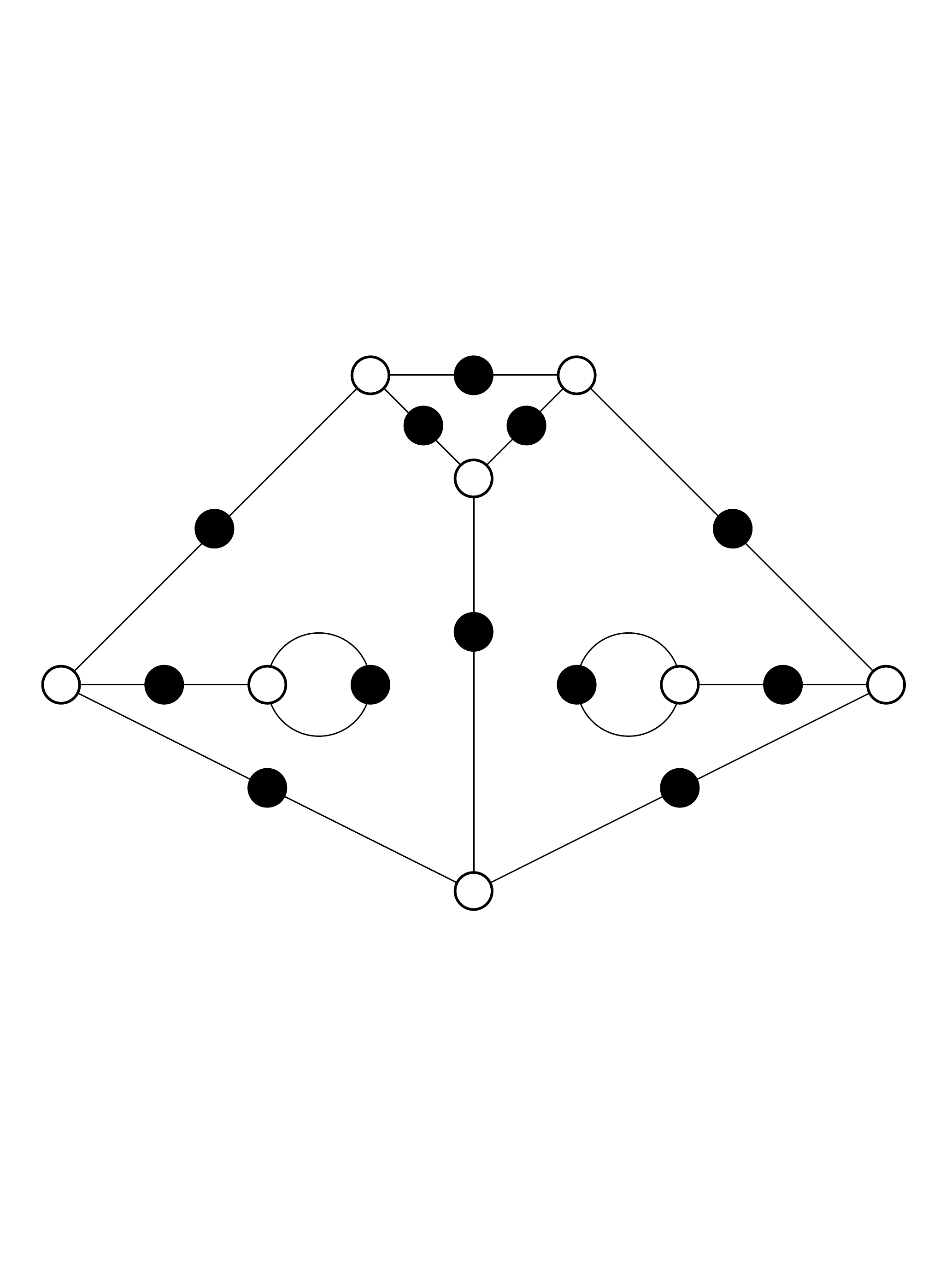}}
\par\end{center}{\scriptsize \par}

\begin{center}
{\scriptsize $7,7,5,3,1,1\;\left(\sqrt{21}\right)$}
\par\end{center}%
\end{minipage}
\par\end{center}{\scriptsize \par}

\begin{center}
{\scriptsize }%
\begin{minipage}[t]{0.33\textwidth}%
\begin{center}
{\scriptsize \includegraphics[scale=0.15]{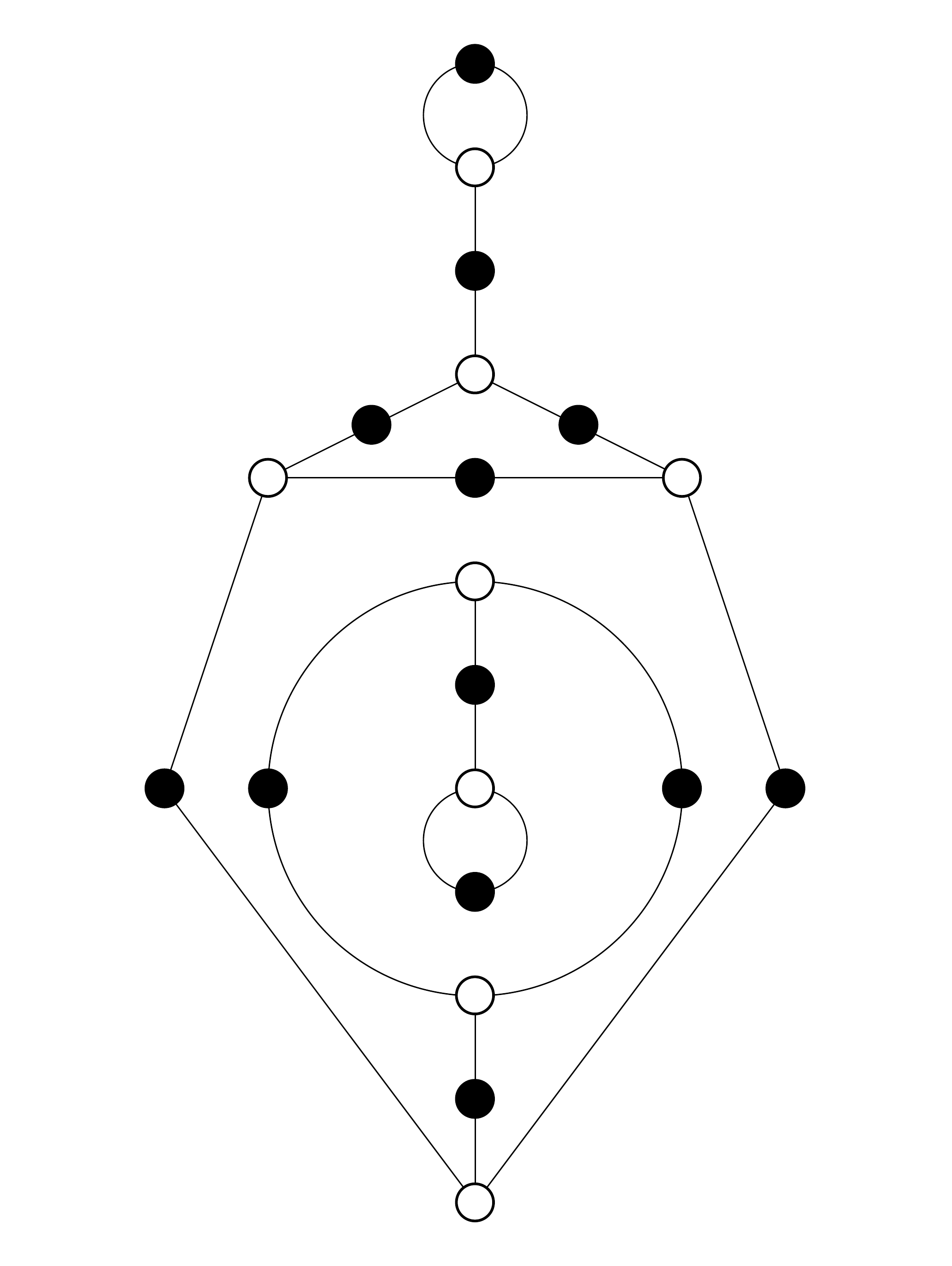}}
\par\end{center}{\scriptsize \par}

\begin{center}
{\scriptsize $7,7,5,3,1,1\;\left(\sqrt{21}\right)$}
\par\end{center}%
\end{minipage}{\scriptsize }%
\begin{minipage}[t]{0.33\textwidth}%
\begin{center}
{\scriptsize \includegraphics[scale=0.15]{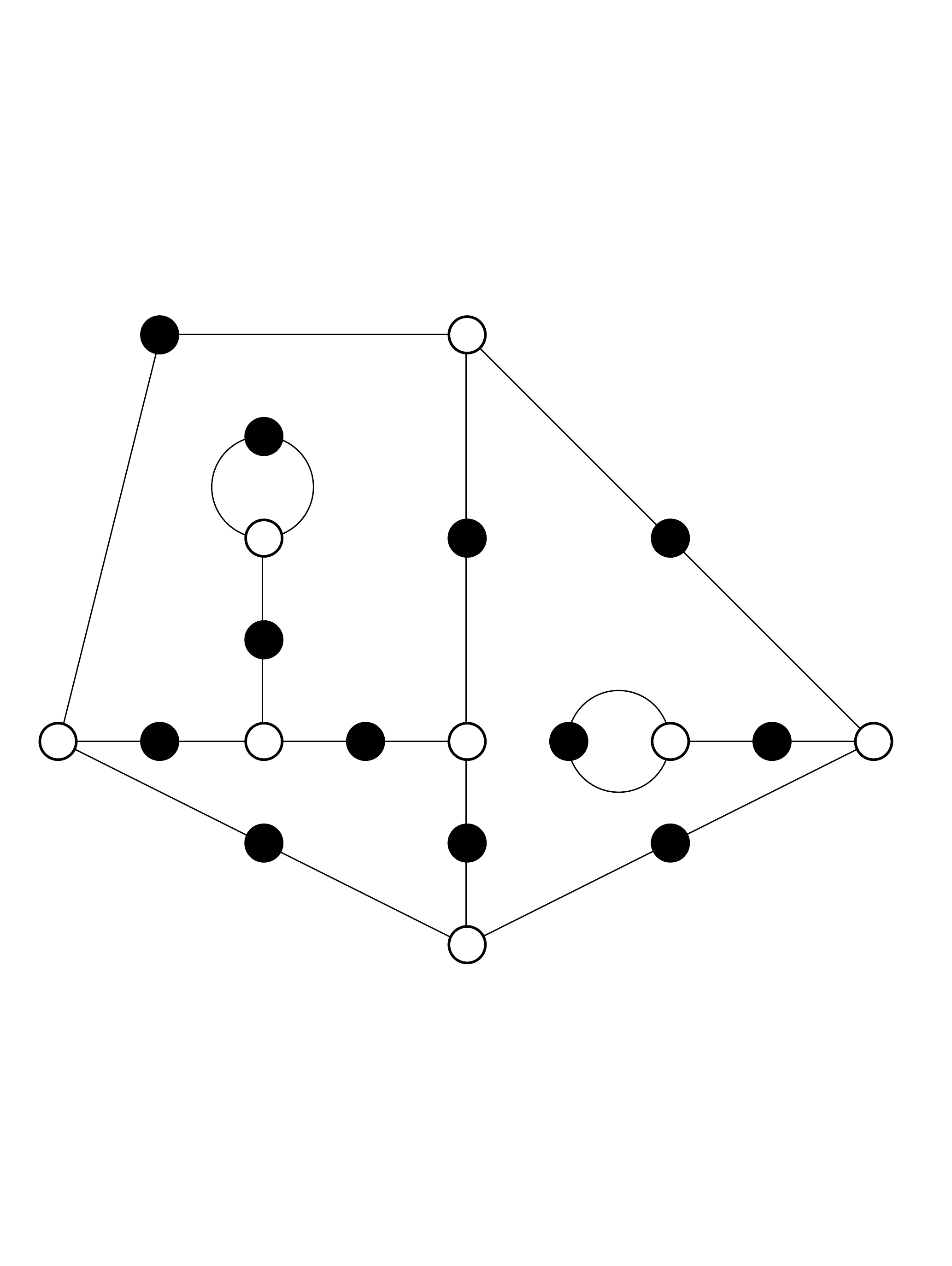}}
\par\end{center}{\scriptsize \par}

\begin{center}
{\scriptsize $7,7,4,4,1,1\;\left(\sqrt{-7}\right)$}
\par\end{center}%
\end{minipage}{\scriptsize }%
\begin{minipage}[t]{0.33\textwidth}%
\begin{center}
{\scriptsize \includegraphics[scale=0.15]{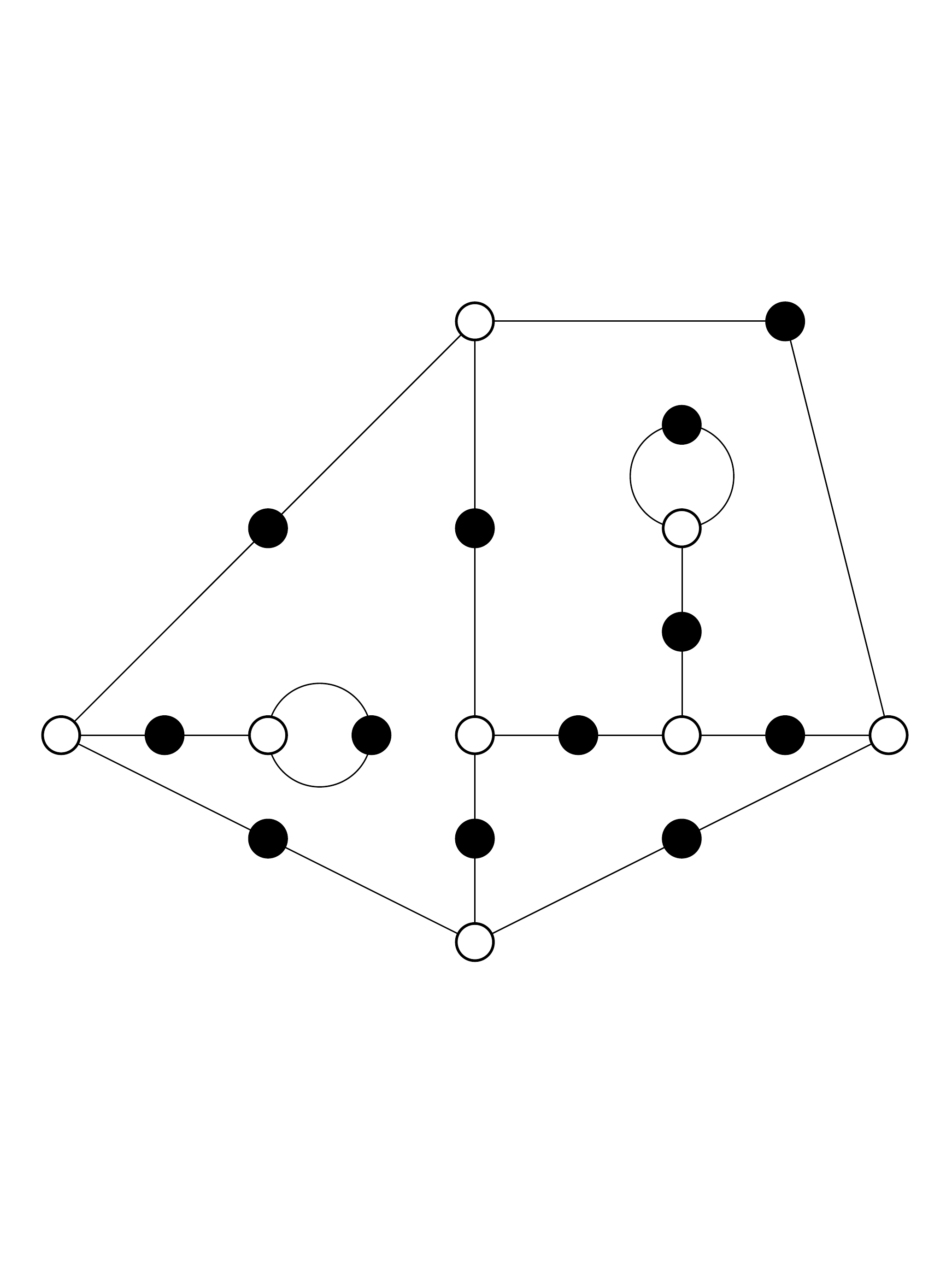}}
\par\end{center}{\scriptsize \par}

\begin{center}
{\scriptsize $7,7,4,4,1,1\;\left(\sqrt{-7}\right)$}
\par\end{center}%
\end{minipage}
\par\end{center}{\scriptsize \par}

\begin{center}
{\scriptsize }%
\begin{minipage}[t]{0.33\textwidth}%
\begin{center}
{\scriptsize \includegraphics[scale=0.15]{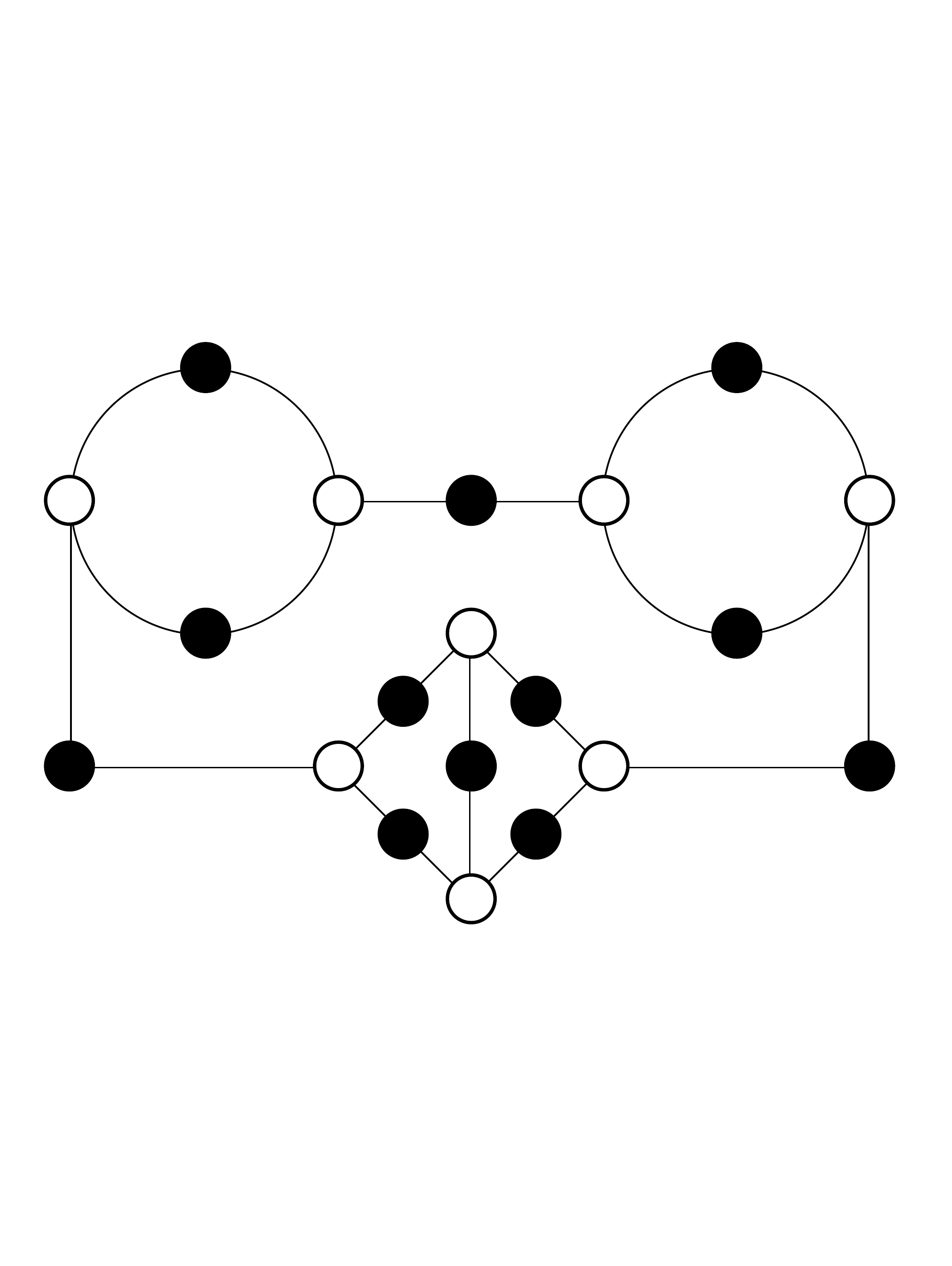}}
\par\end{center}{\scriptsize \par}

\begin{center}
{\scriptsize $7,7,3,3,2,2\;\left(\sqrt{7}\right)$}
\par\end{center}%
\end{minipage}{\scriptsize }%
\begin{minipage}[t]{0.33\textwidth}%
\begin{center}
{\scriptsize \includegraphics[scale=0.15]{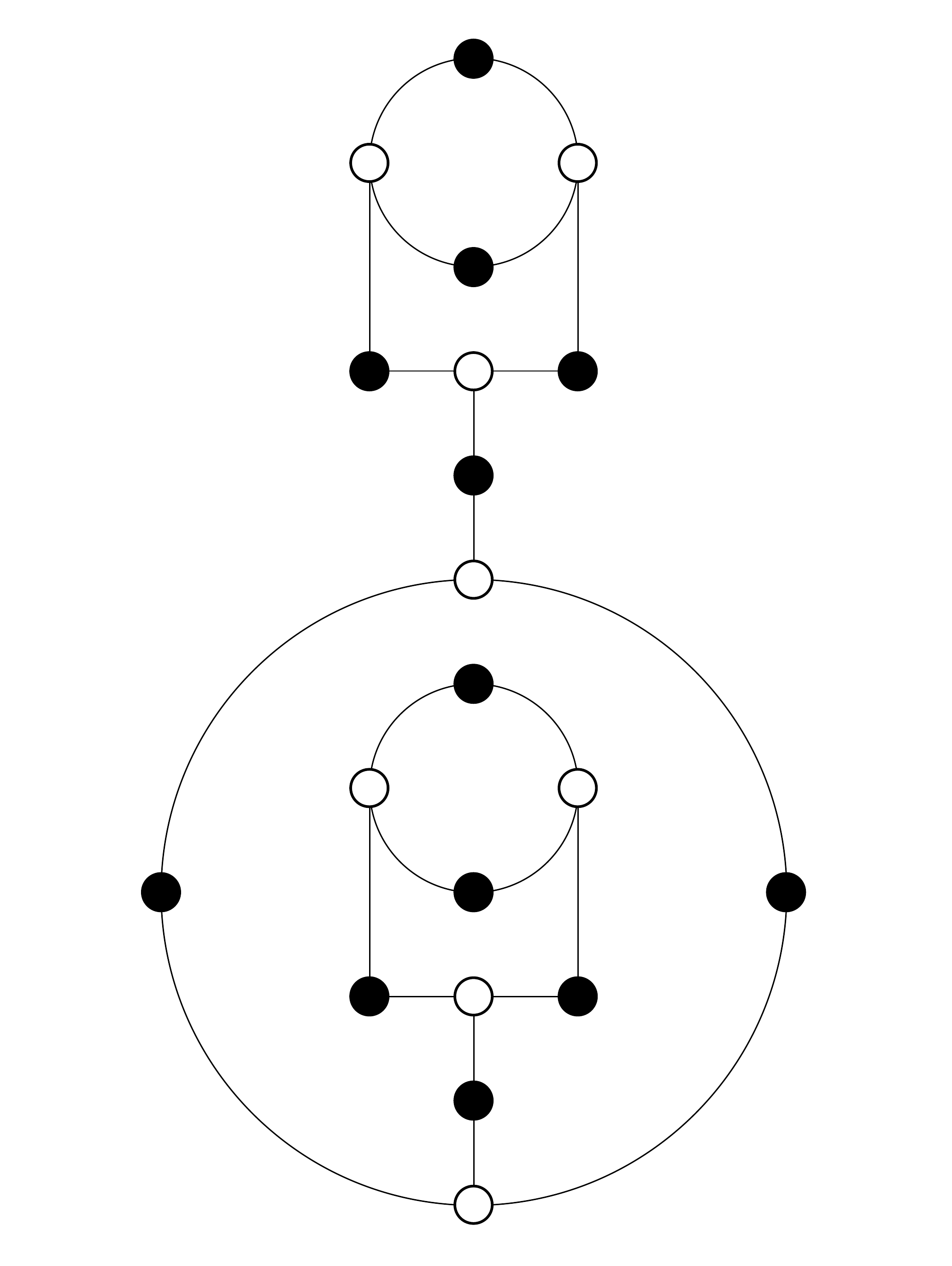}}
\par\end{center}{\scriptsize \par}

\begin{center}
{\scriptsize $7,7,3,3,2,2\;\left(\sqrt{7}\right)$}
\par\end{center}%
\end{minipage}{\scriptsize }%
\begin{minipage}[t]{0.33\textwidth}%
\begin{center}
{\scriptsize \includegraphics[scale=0.15]{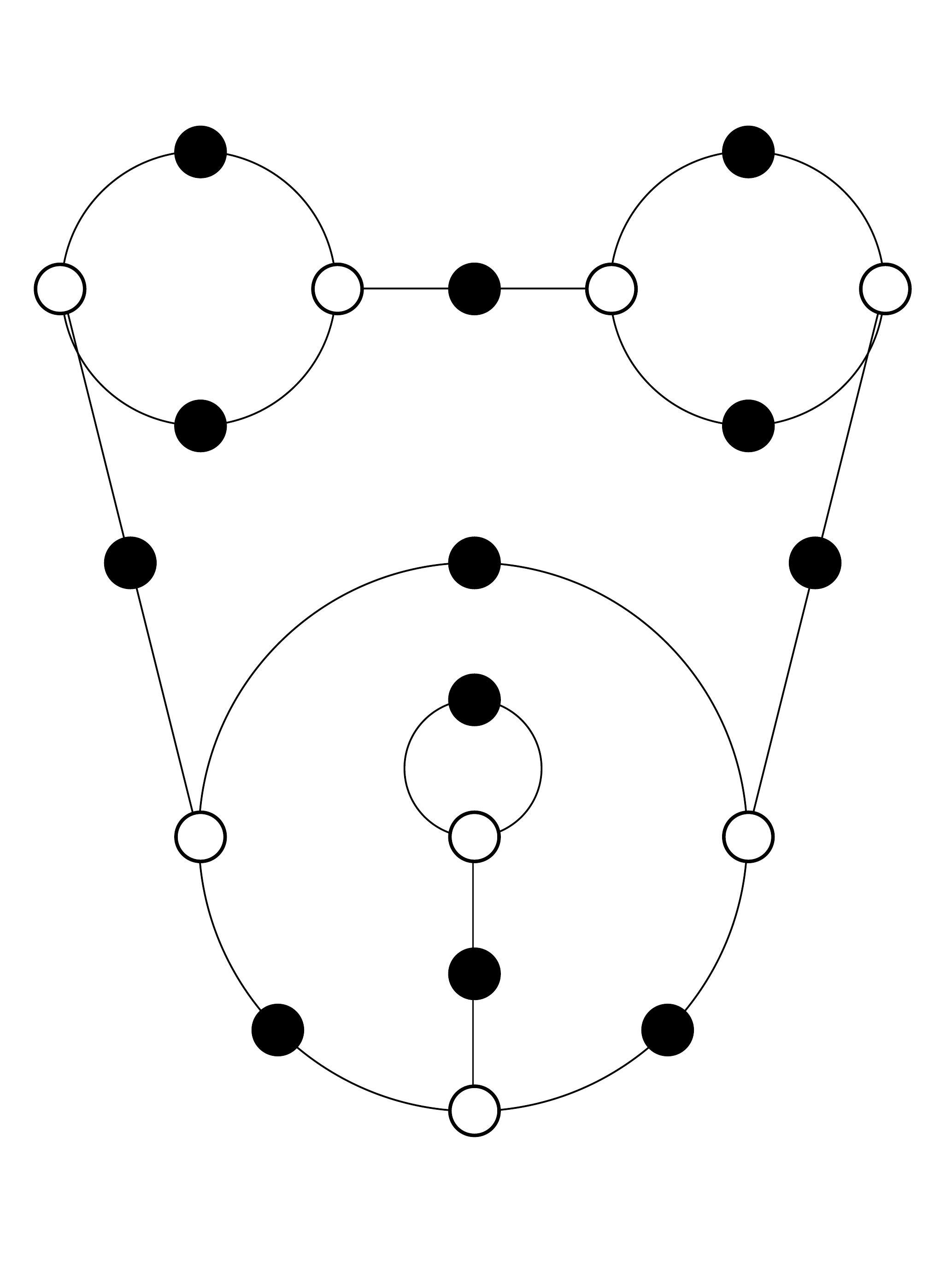}}
\par\end{center}{\scriptsize \par}

\begin{center}
{\scriptsize $7,6,6,2,2,1\;\left(\mathbb{Q}\right)$}
\par\end{center}%
\end{minipage}
\par\end{center}{\scriptsize \par}

\begin{center}
{\scriptsize }%
\begin{minipage}[t]{0.33\textwidth}%
\begin{center}
{\scriptsize \includegraphics[scale=0.15]{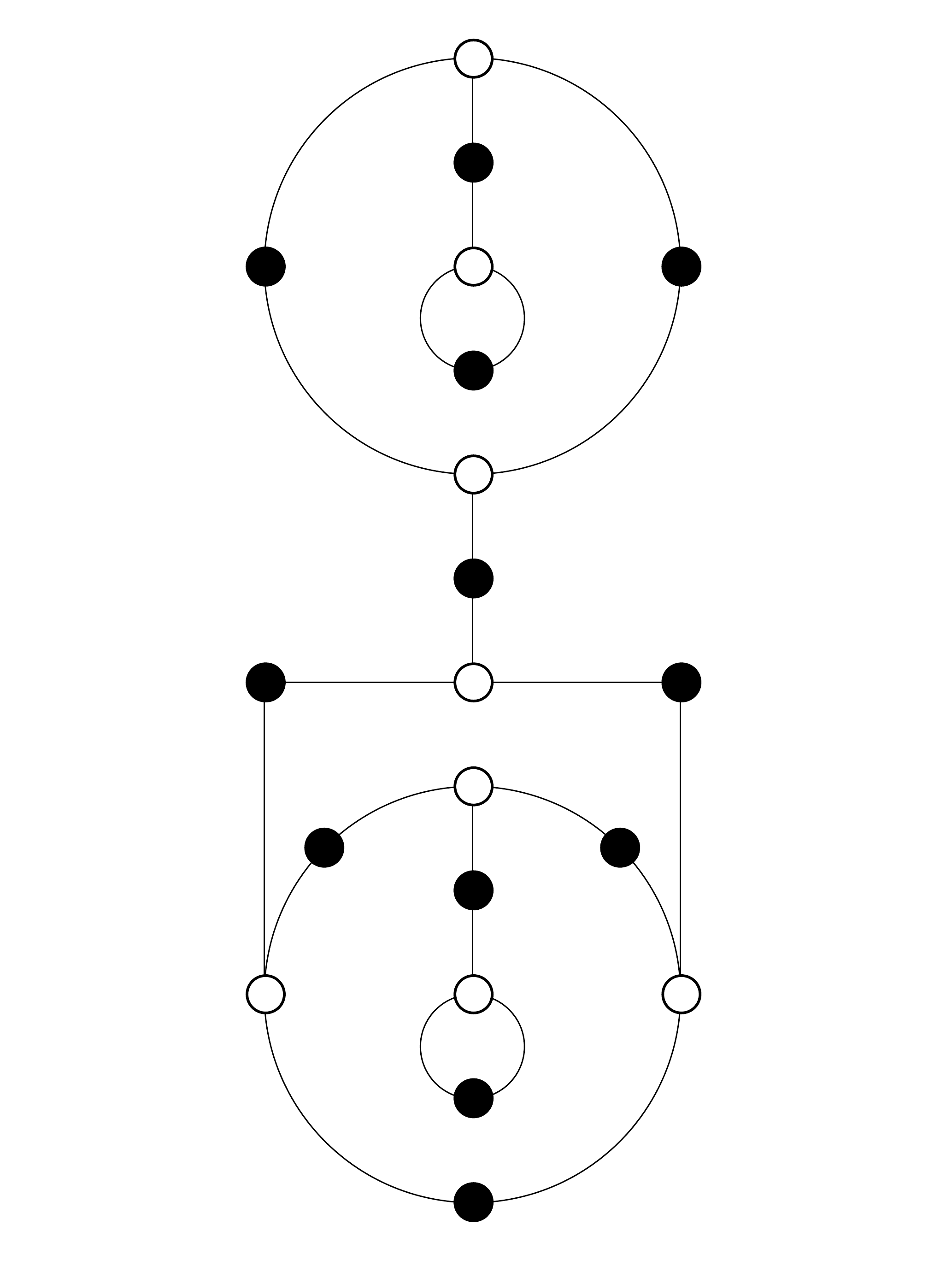}}
\par\end{center}{\scriptsize \par}

\begin{center}
{\scriptsize $7,6,5,4,1,1\;\left(\mathbb{Q}\right)$}
\par\end{center}%
\end{minipage}{\scriptsize }%
\begin{minipage}[t]{0.33\textwidth}%
\begin{center}
{\scriptsize \includegraphics[scale=0.15]{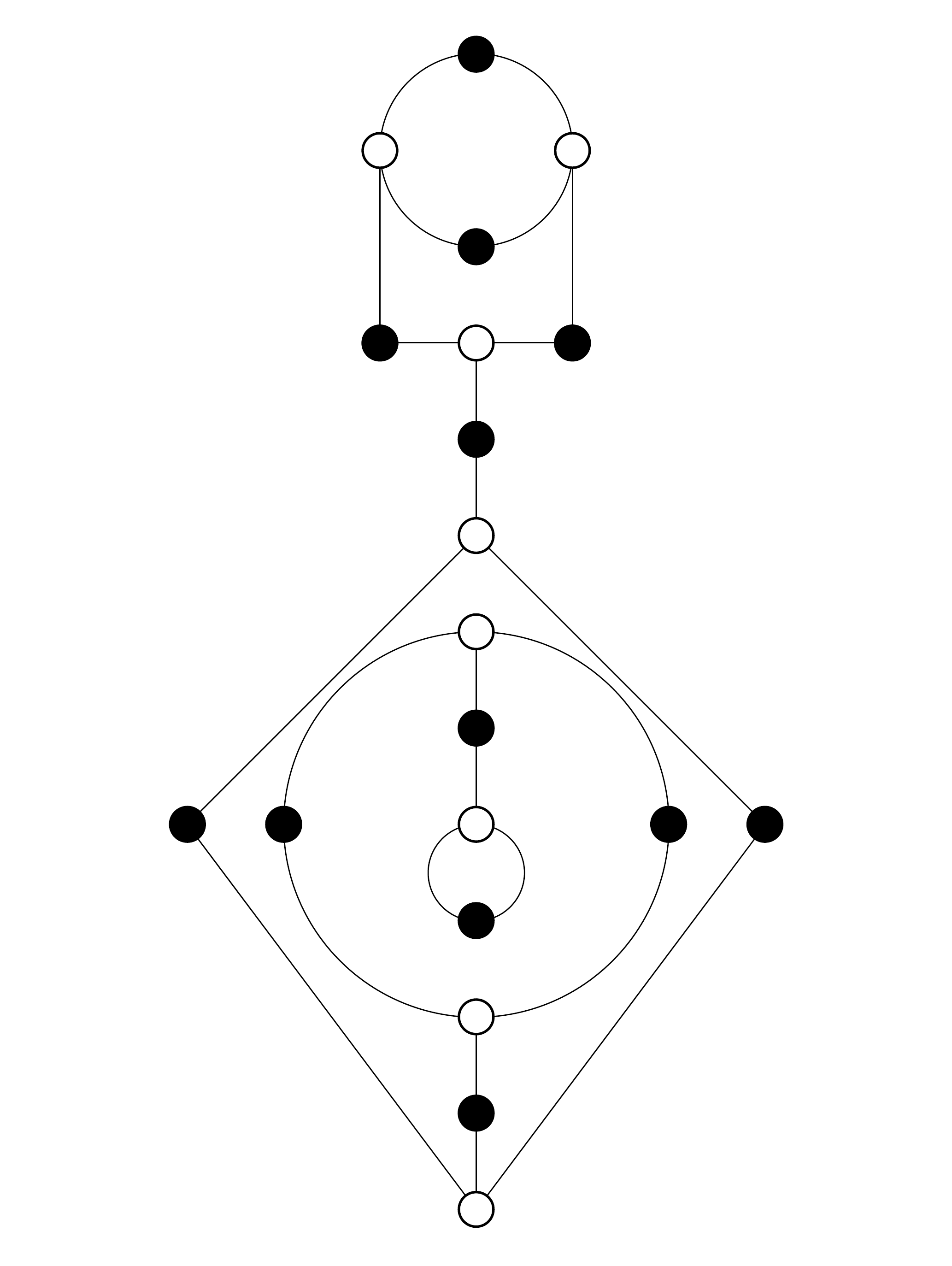}}
\par\end{center}{\scriptsize \par}

\begin{center}
{\scriptsize $7,6,5,3,2,1\;\left(\mathrm{cubic}\right)$}
\par\end{center}%
\end{minipage}{\scriptsize }%
\begin{minipage}[t]{0.33\textwidth}%
\begin{center}
{\scriptsize \includegraphics[scale=0.15]{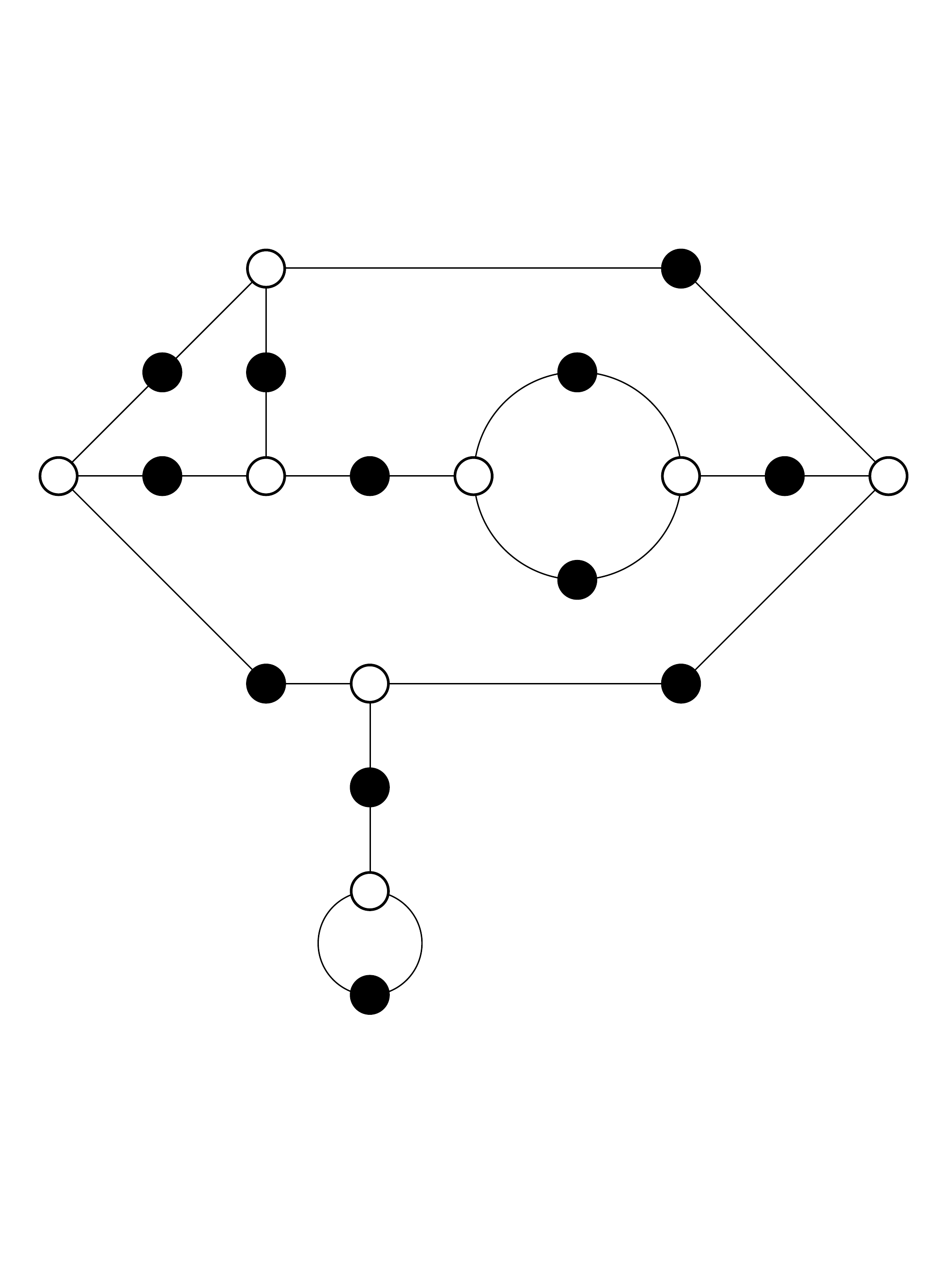}}
\par\end{center}{\scriptsize \par}

\begin{center}
{\scriptsize $7,6,5,3,2,1\;\left(\mathrm{cubic}\right)$}
\par\end{center}%
\end{minipage}
\par\end{center}{\scriptsize \par}

\begin{center}
{\scriptsize }%
\begin{minipage}[t]{0.33\textwidth}%
\begin{center}
{\scriptsize \includegraphics[scale=0.15]{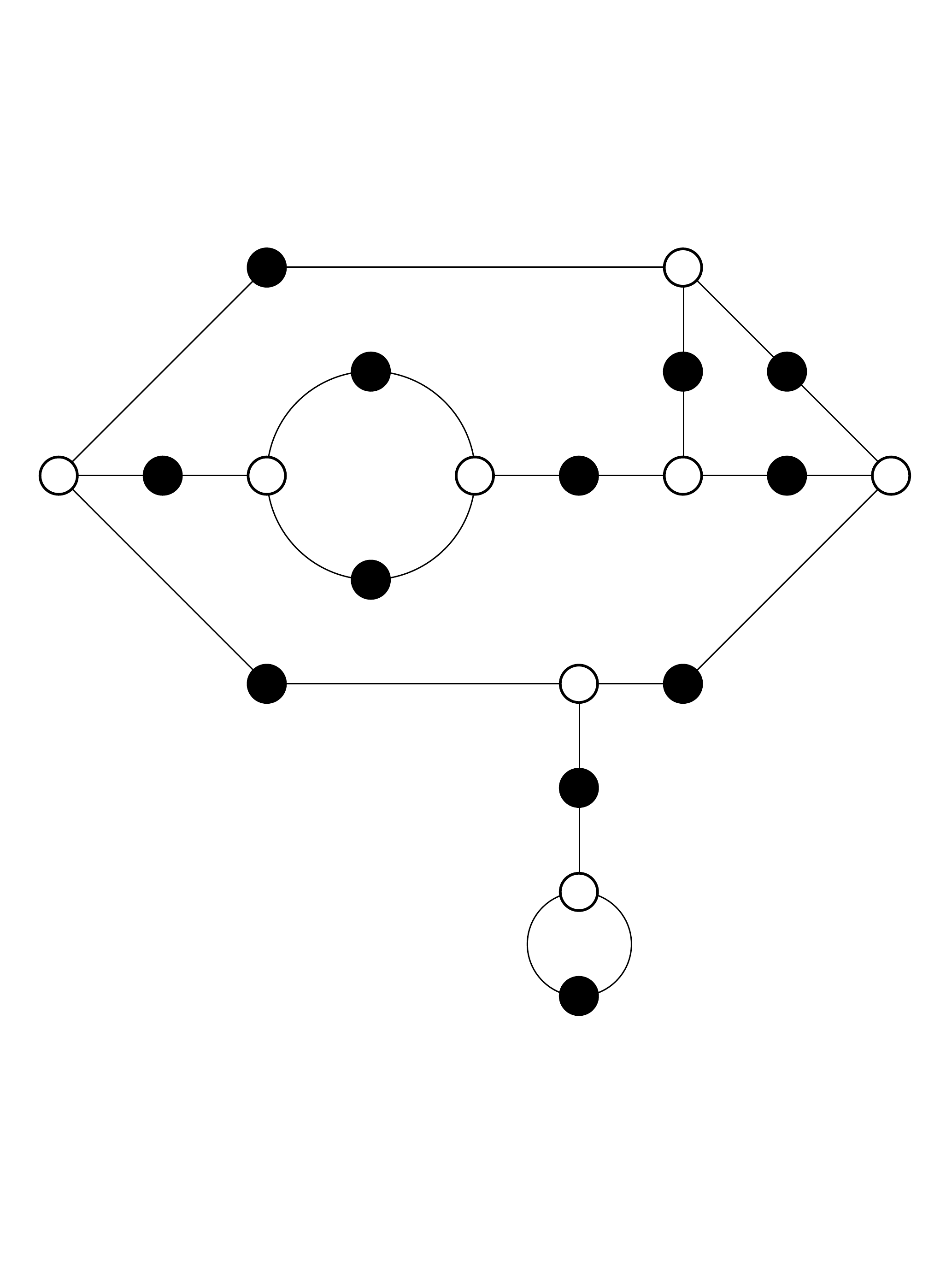}}
\par\end{center}{\scriptsize \par}

\begin{center}
{\scriptsize $7,6,5,3,2,1\;\left(\mathrm{cubic}\right)$}
\par\end{center}%
\end{minipage}{\scriptsize }%
\begin{minipage}[t]{0.33\textwidth}%
\begin{center}
{\scriptsize \includegraphics[scale=0.15]{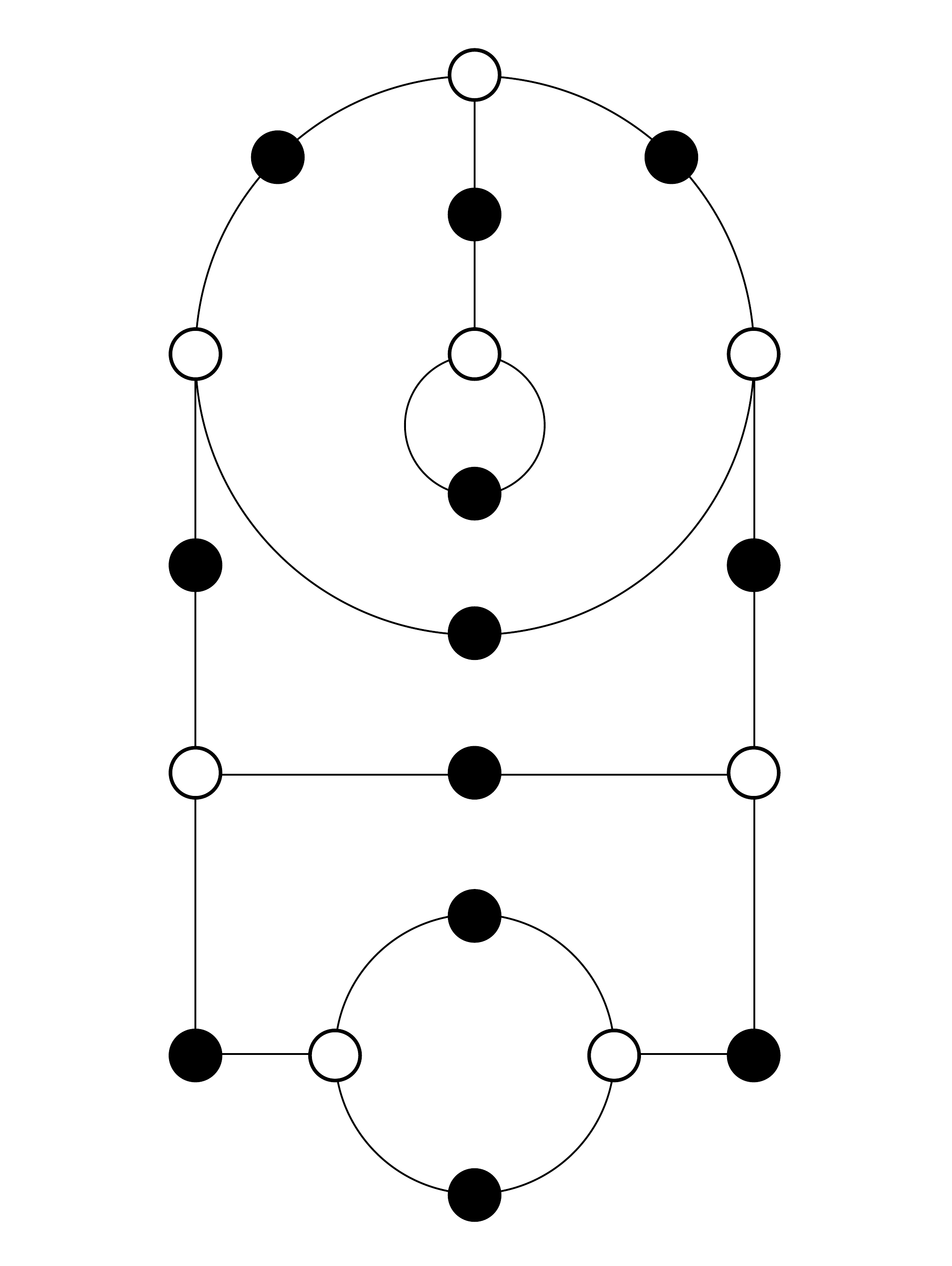}}
\par\end{center}{\scriptsize \par}

\begin{center}
{\scriptsize $7,6,4,4,2,1\;\left(\mathbb{Q}\right)$}
\par\end{center}%
\end{minipage}{\scriptsize }%
\begin{minipage}[t]{0.33\textwidth}%
\begin{center}
{\scriptsize \includegraphics[scale=0.15]{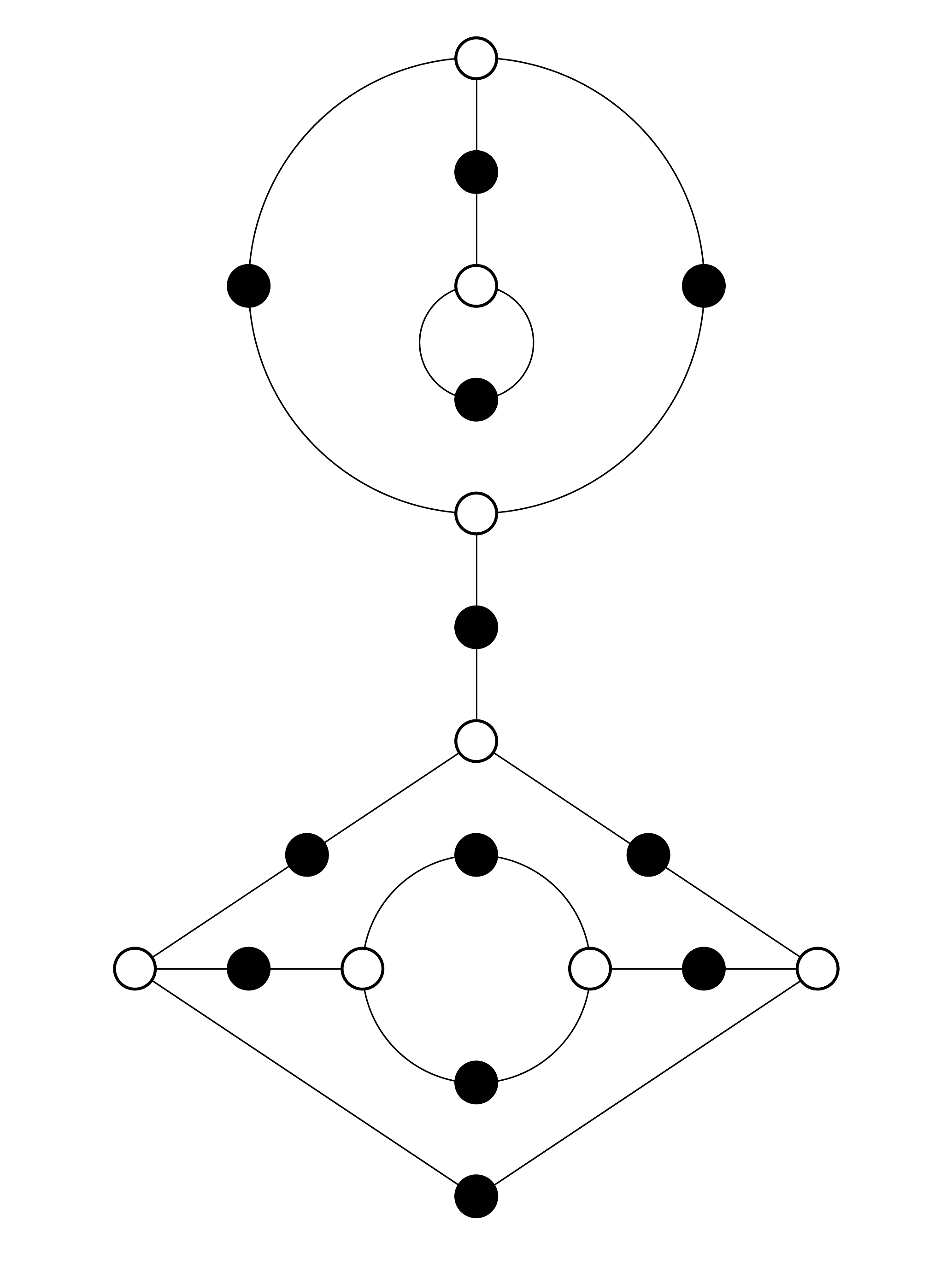}}
\par\end{center}{\scriptsize \par}

\begin{center}
{\scriptsize $7,5,5,4,2,1\;\left(\sqrt{2}\right)$}
\par\end{center}%
\end{minipage}
\par\end{center}{\scriptsize \par}

\begin{center}
{\scriptsize }%
\begin{minipage}[t]{0.33\textwidth}%
\begin{center}
{\scriptsize \includegraphics[scale=0.15]{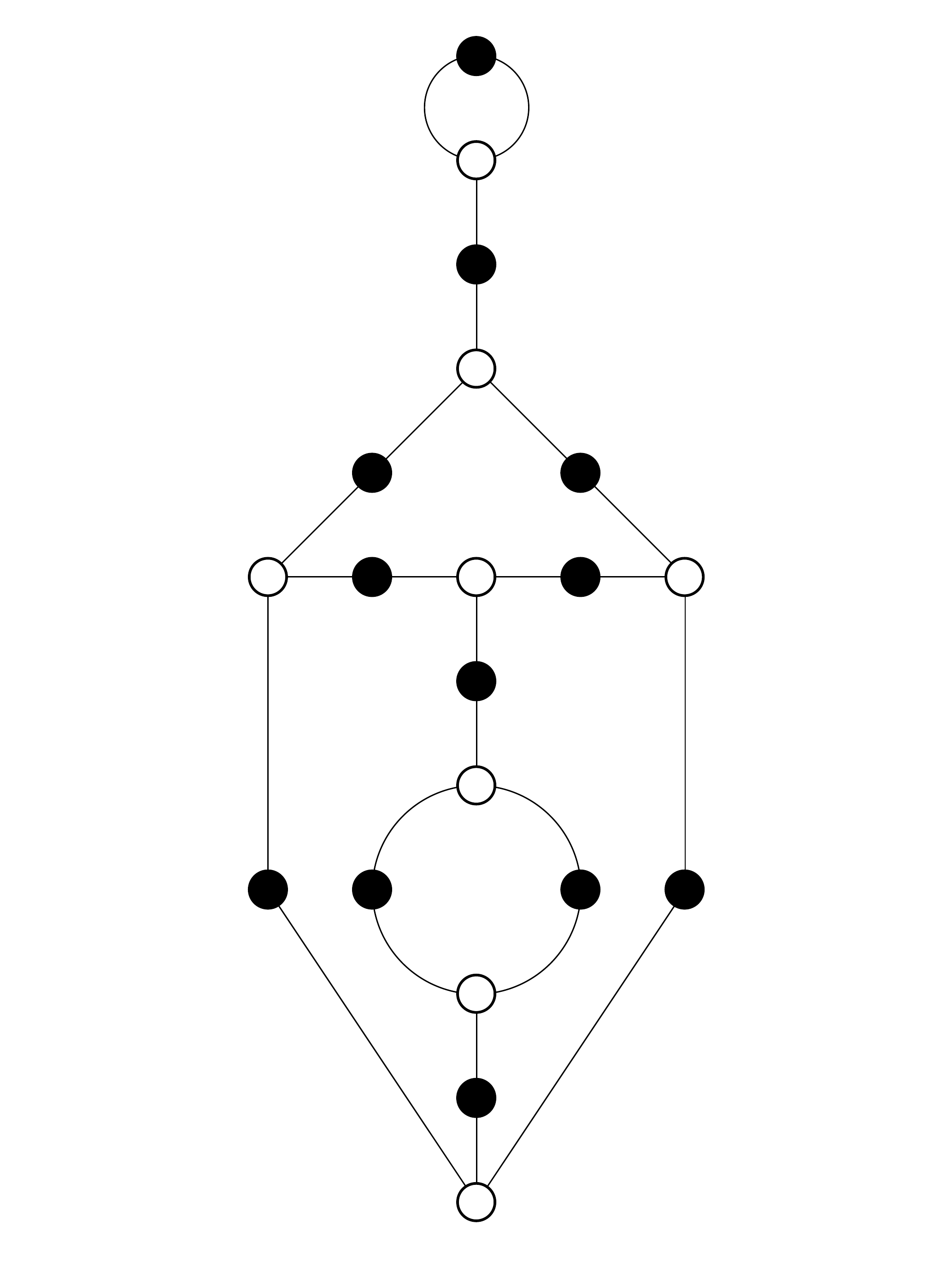}}
\par\end{center}{\scriptsize \par}

\begin{center}
{\scriptsize $7,5,5,4,2,1\;\left(\sqrt{2}\right)$}
\par\end{center}%
\end{minipage}{\scriptsize }%
\begin{minipage}[t]{0.33\textwidth}%
\begin{center}
{\scriptsize \includegraphics[scale=0.15]{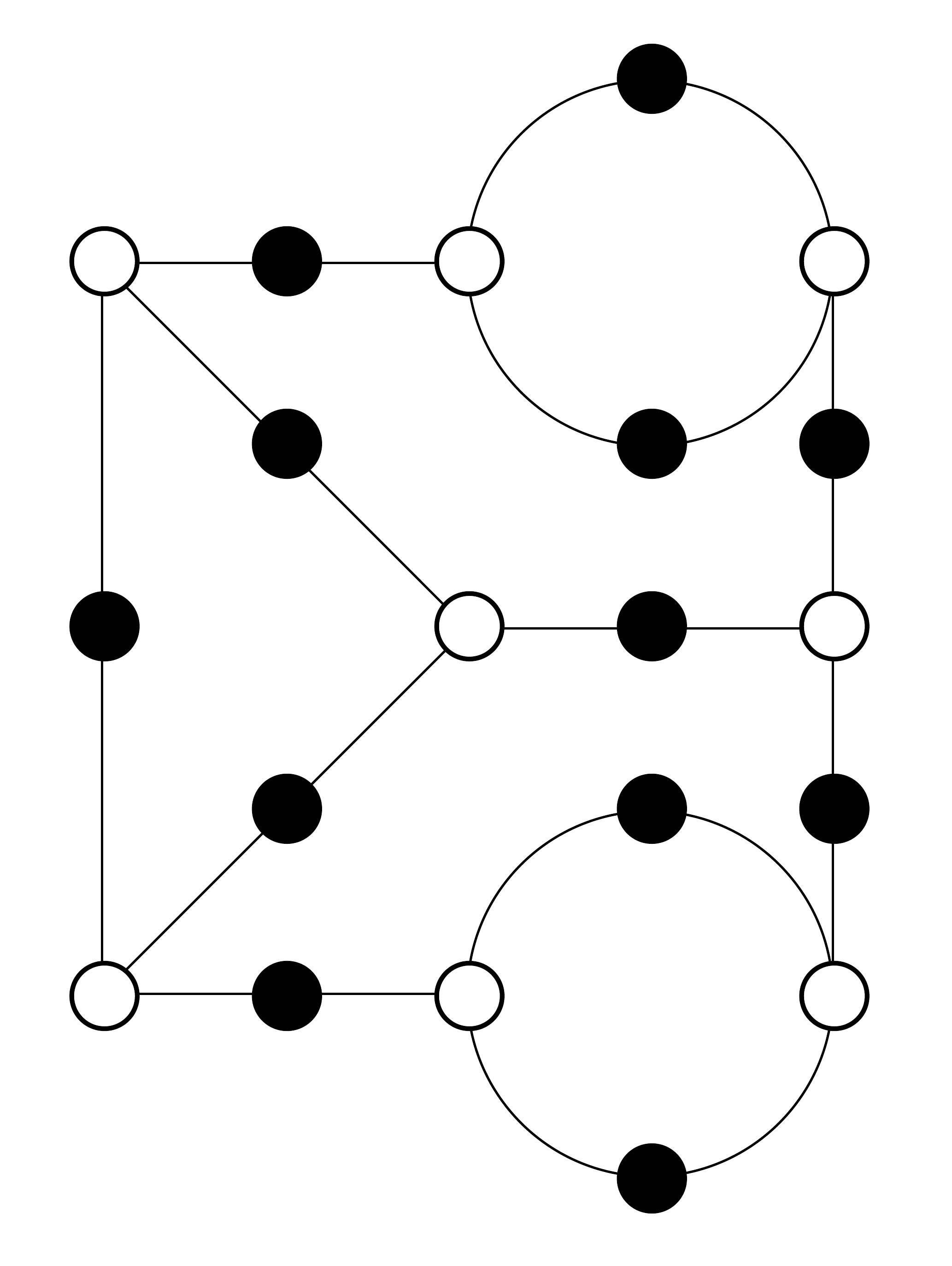}}
\par\end{center}{\scriptsize \par}

\begin{center}
{\scriptsize $7,5,5,3,2,2\;\left(\mathbb{Q}\right)$}
\par\end{center}%
\end{minipage}{\scriptsize }%
\begin{minipage}[t]{0.33\textwidth}%
\begin{center}
{\scriptsize \includegraphics[scale=0.15]{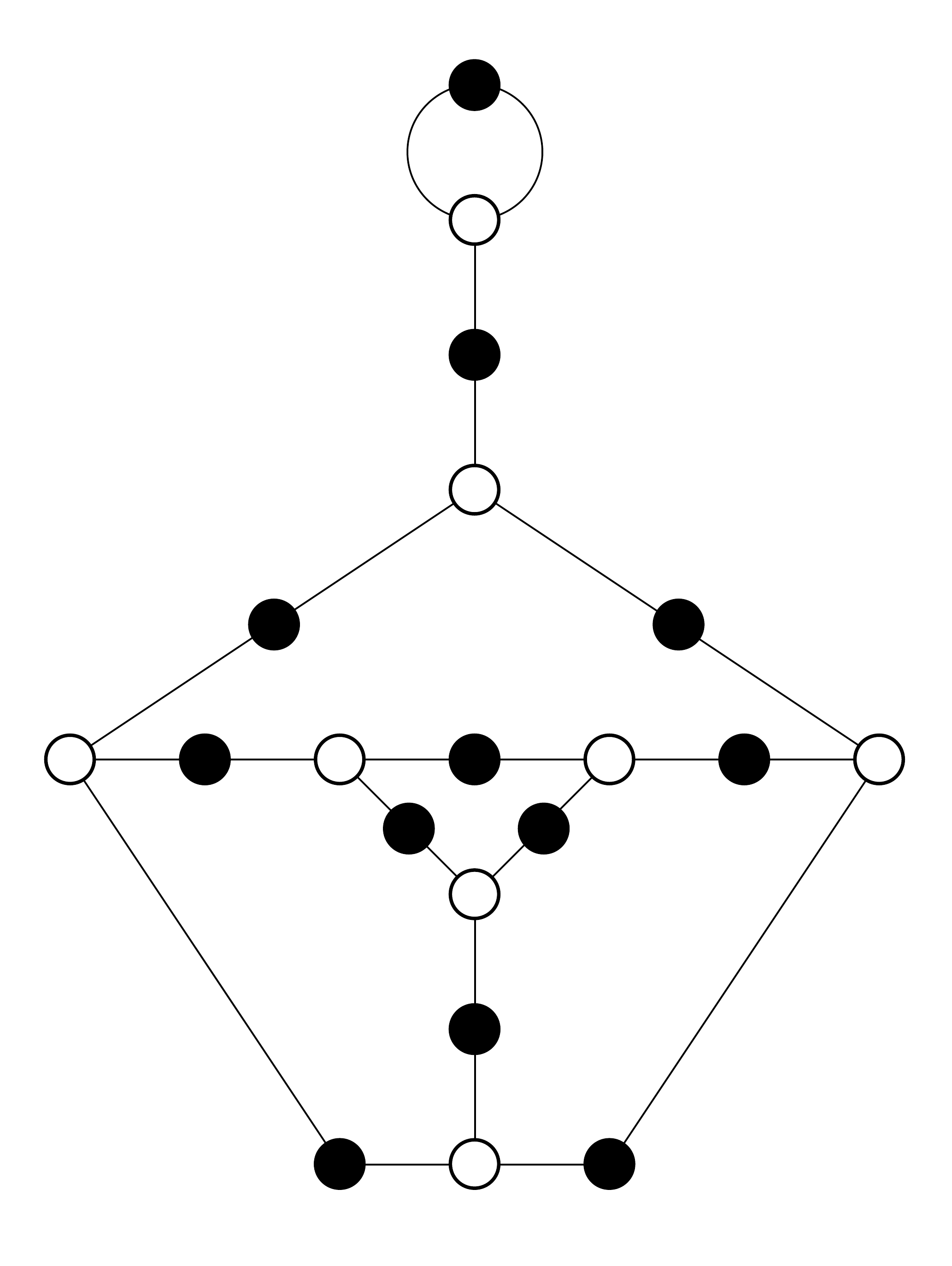}}
\par\end{center}{\scriptsize \par}

\begin{center}
{\scriptsize $7,5,4,4,3,1\;\left(\mathbb{Q}\right)$}
\par\end{center}%
\end{minipage}
\par\end{center}{\scriptsize \par}

\begin{center}
{\scriptsize }%
\begin{minipage}[t]{0.33\textwidth}%
\begin{center}
{\scriptsize \includegraphics[scale=0.15]{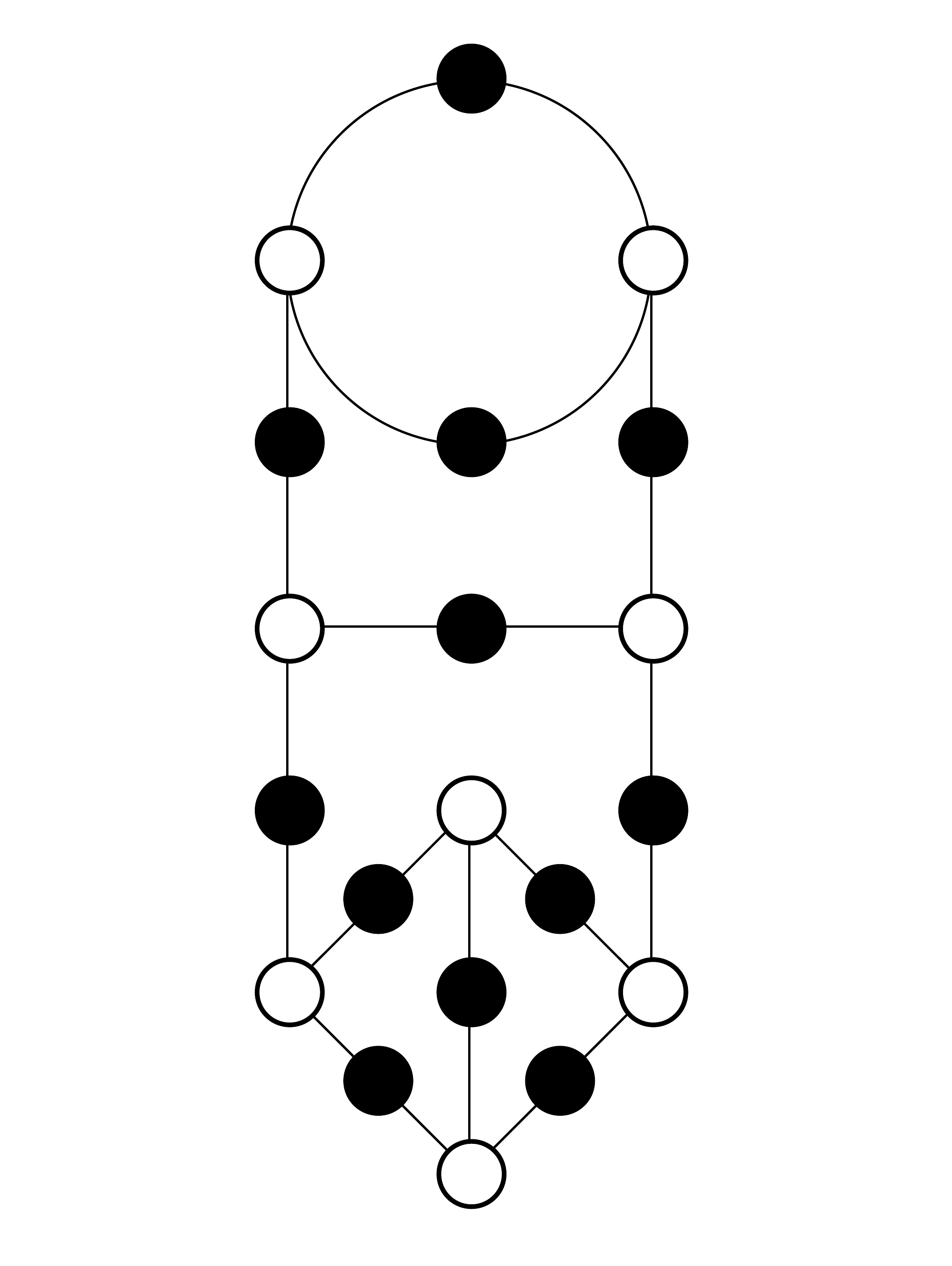}}
\par\end{center}{\scriptsize \par}

\begin{center}
{\scriptsize $7,5,4,3,3,2\;\left(\mathbb{Q}\right)$}
\par\end{center}%
\end{minipage}{\scriptsize }%
\begin{minipage}[t]{0.33\textwidth}%
\begin{center}
{\scriptsize \includegraphics[scale=0.15]{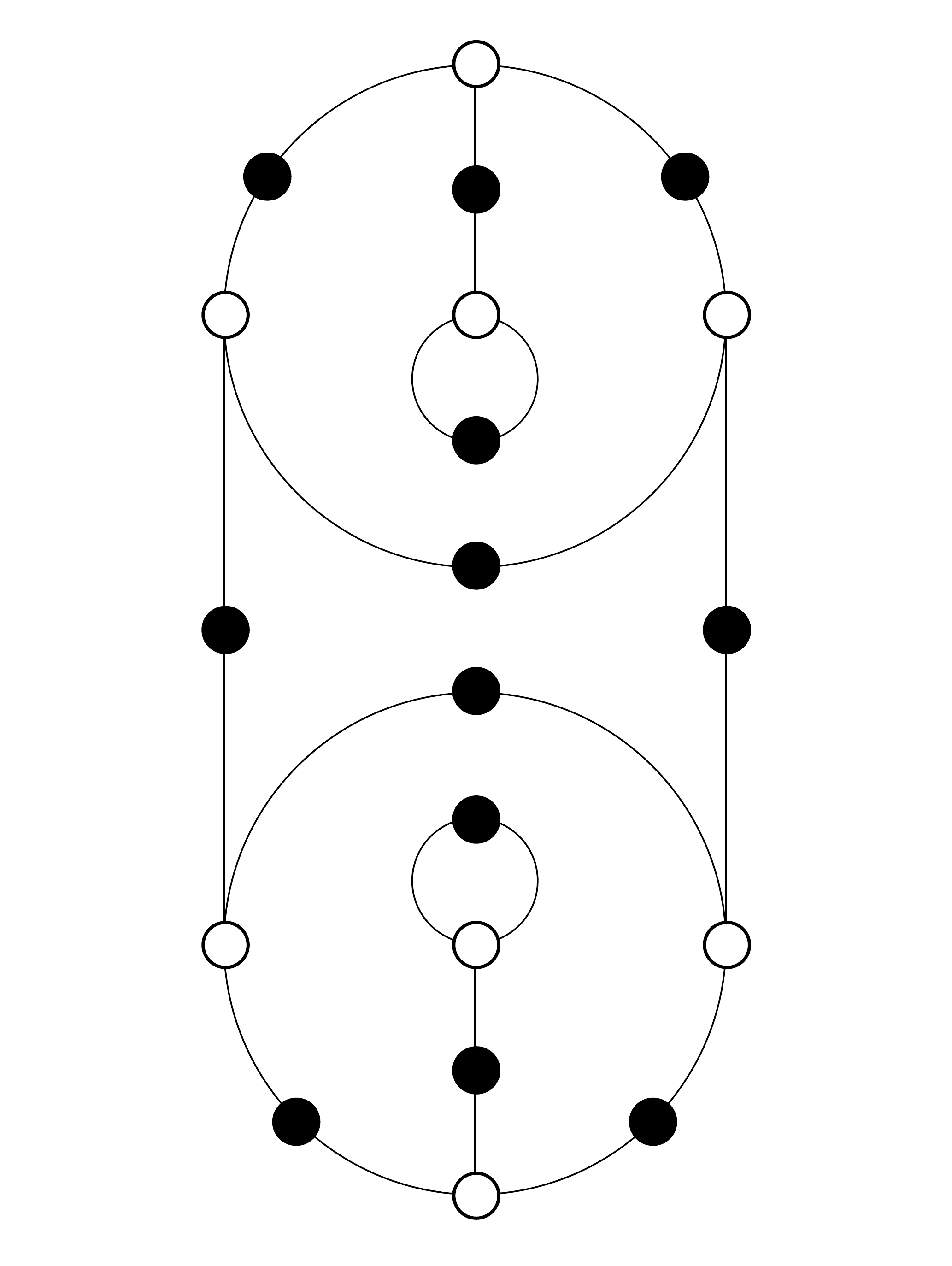}}
\par\end{center}{\scriptsize \par}

\begin{center}
{\scriptsize $6,6,6,4,1,1\;\left(\mathbb{Q}\right)$}
\par\end{center}%
\end{minipage}{\scriptsize }%
\begin{minipage}[t]{0.33\textwidth}%
\begin{center}
{\scriptsize \includegraphics[scale=0.15]{\string"PICT/6-6-6-2-2-2\string".pdf}}
\par\end{center}{\scriptsize \par}

\begin{center}
{\scriptsize $6,6,6,2,2,2\;\left(\mathbb{Q}\right)$}
\par\end{center}%
\end{minipage}
\par\end{center}{\scriptsize \par}

\begin{center}
{\scriptsize }%
\begin{minipage}[t]{0.33\textwidth}%
\begin{center}
{\scriptsize \includegraphics[scale=0.15]{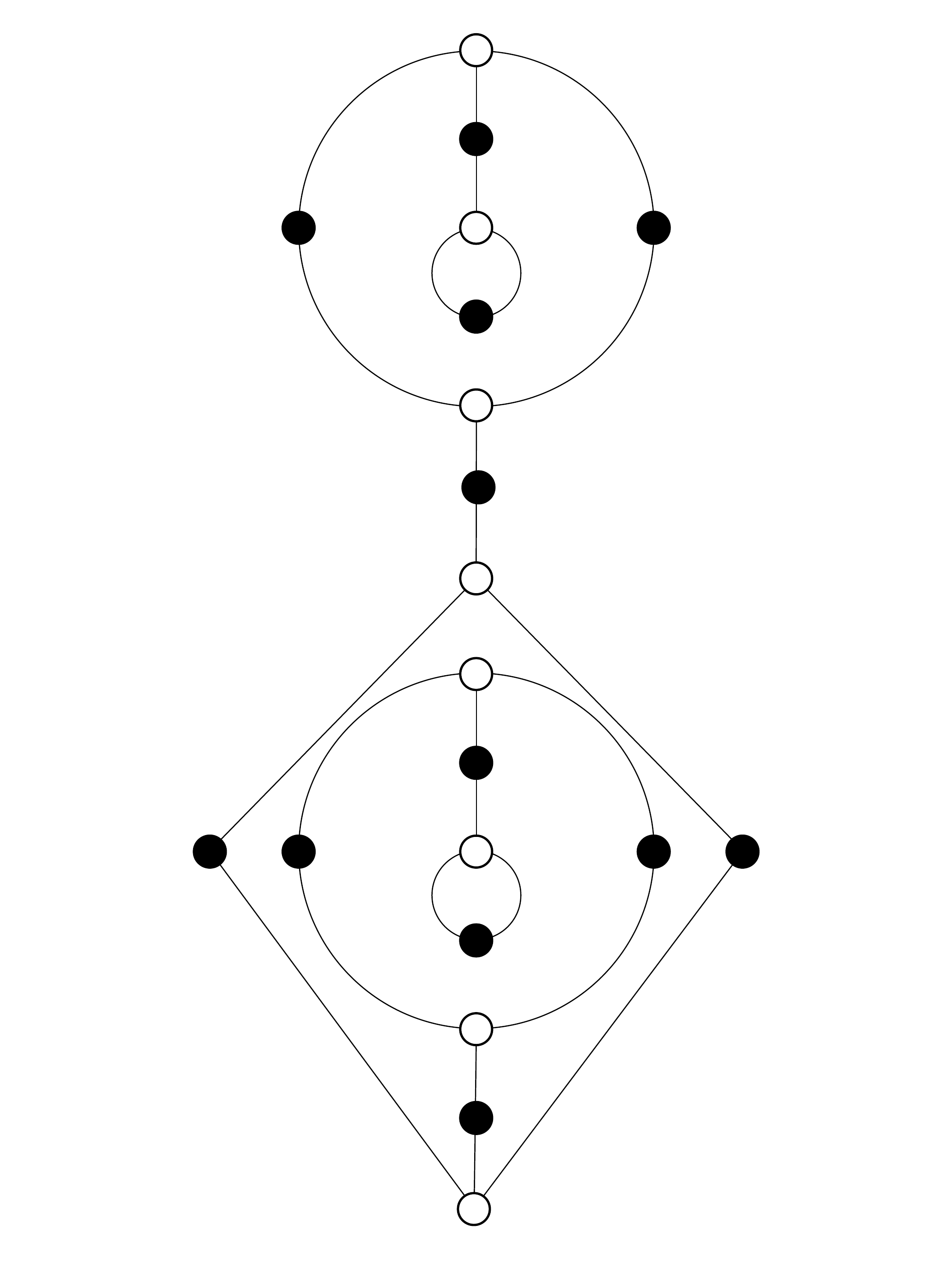}}
\par\end{center}{\scriptsize \par}

\begin{center}
{\scriptsize $6,6,5,5,1,1\;\left(\sqrt{3}\right)$}
\par\end{center}%
\end{minipage}{\scriptsize }%
\begin{minipage}[t]{0.33\textwidth}%
\begin{center}
{\scriptsize \includegraphics[scale=0.15]{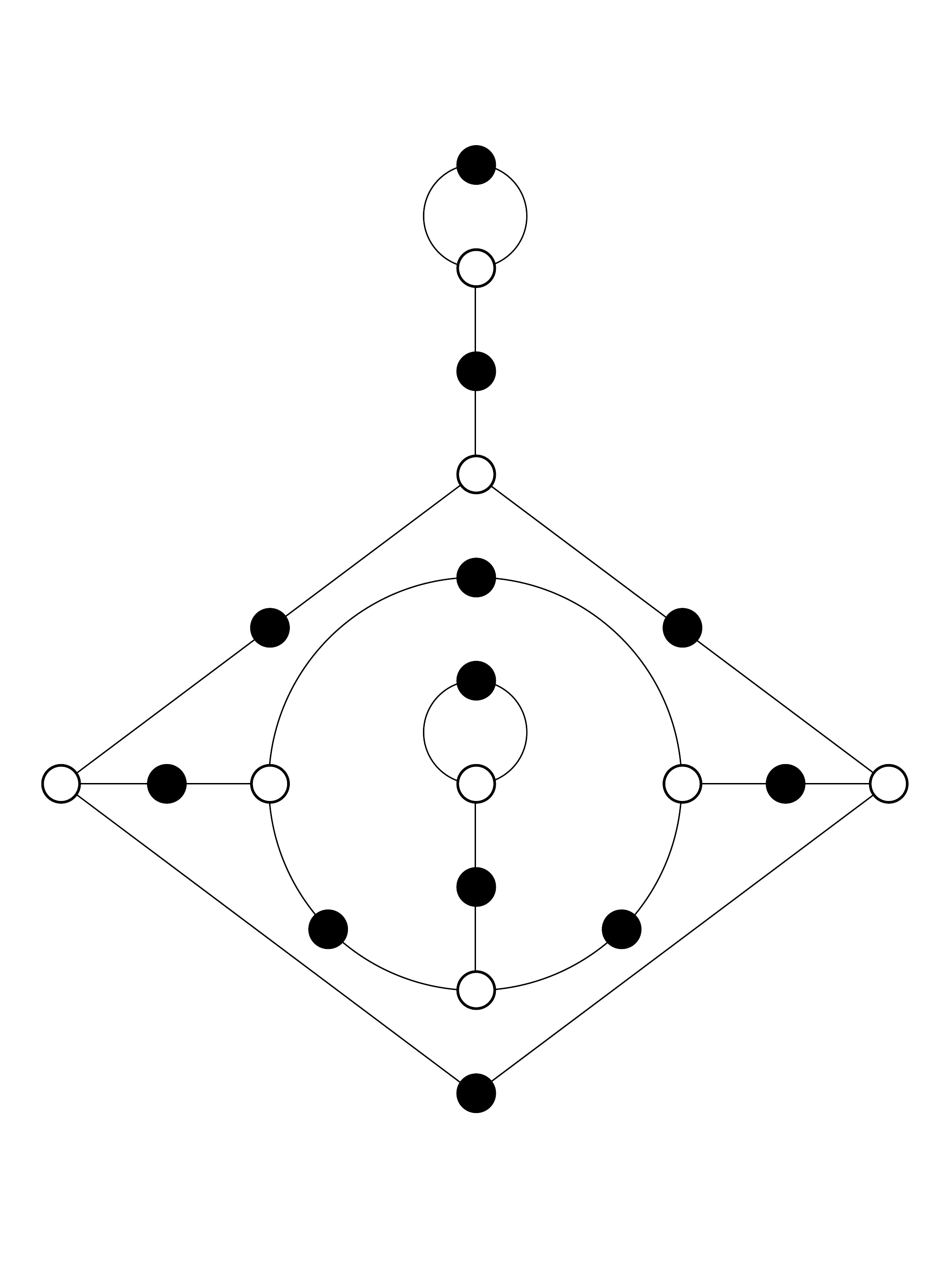}}
\par\end{center}{\scriptsize \par}

\begin{center}
{\scriptsize $6,6,5,5,1,1\;\left(\sqrt{3}\right)$}
\par\end{center}%
\end{minipage}{\scriptsize }%
\begin{minipage}[t]{0.33\textwidth}%
\begin{center}
{\scriptsize \includegraphics[scale=0.15]{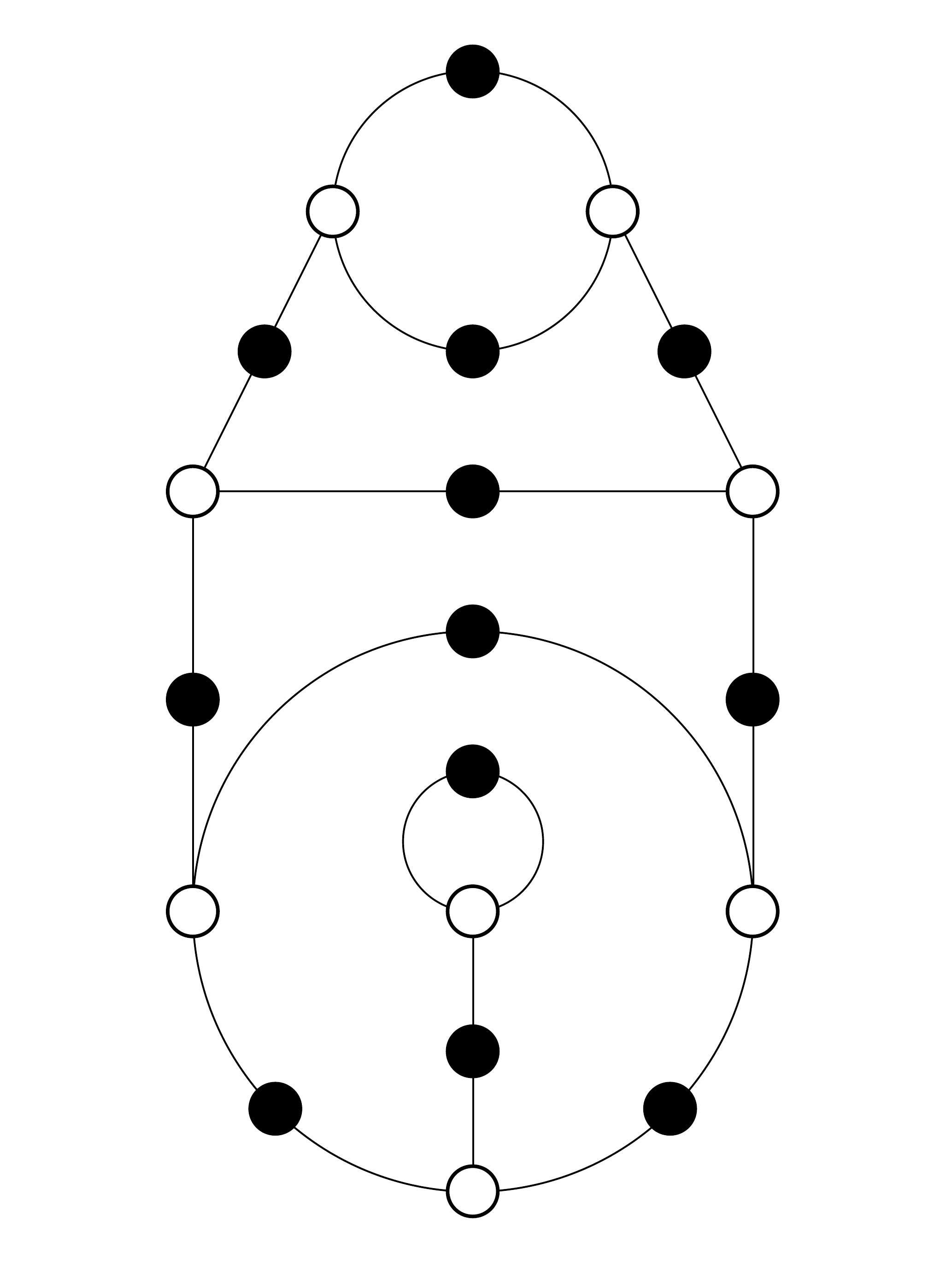}}
\par\end{center}{\scriptsize \par}

\begin{center}
{\scriptsize $6,6,5,4,2,1\;\left(\mathbb{Q}\right)$}
\par\end{center}%
\end{minipage}
\par\end{center}{\scriptsize \par}

\begin{center}
{\scriptsize }%
\begin{minipage}[t]{0.33\textwidth}%
\begin{center}
{\scriptsize \includegraphics[scale=0.15]{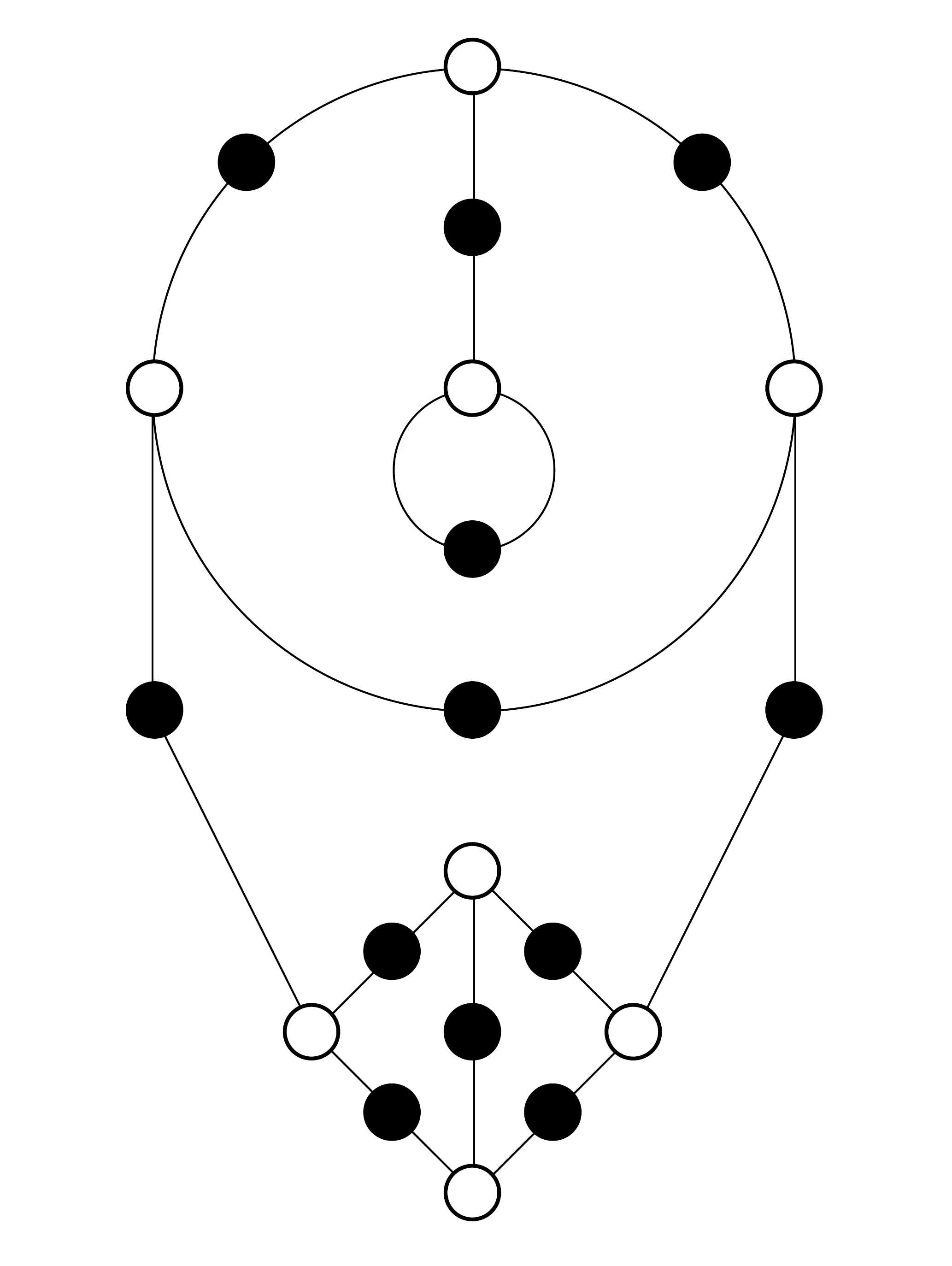}}
\par\end{center}{\scriptsize \par}

\begin{center}
{\scriptsize $6,6,5,3,3,1\;\left(\mathbb{Q}\right)$}
\par\end{center}%
\end{minipage}{\scriptsize }%
\begin{minipage}[t]{0.33\textwidth}%
\begin{center}
{\scriptsize \includegraphics[scale=0.15]{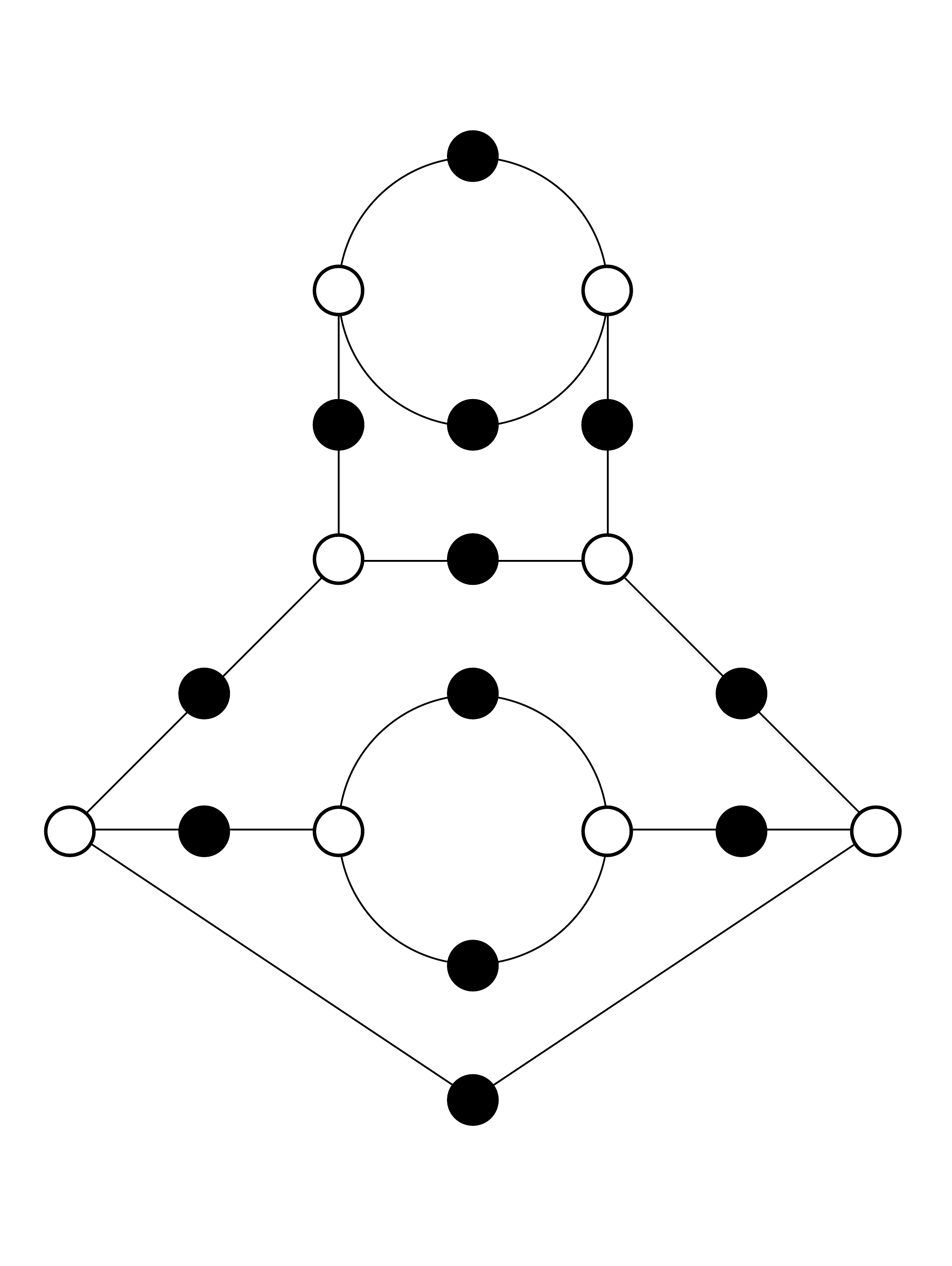}}
\par\end{center}{\scriptsize \par}

\begin{center}
{\scriptsize $6,6,4,4,2,2\;\left(\mathbb{Q}\right)$}
\par\end{center}%
\end{minipage}{\scriptsize }%
\begin{minipage}[t]{0.33\textwidth}%
\begin{center}
{\scriptsize \includegraphics[scale=0.15]{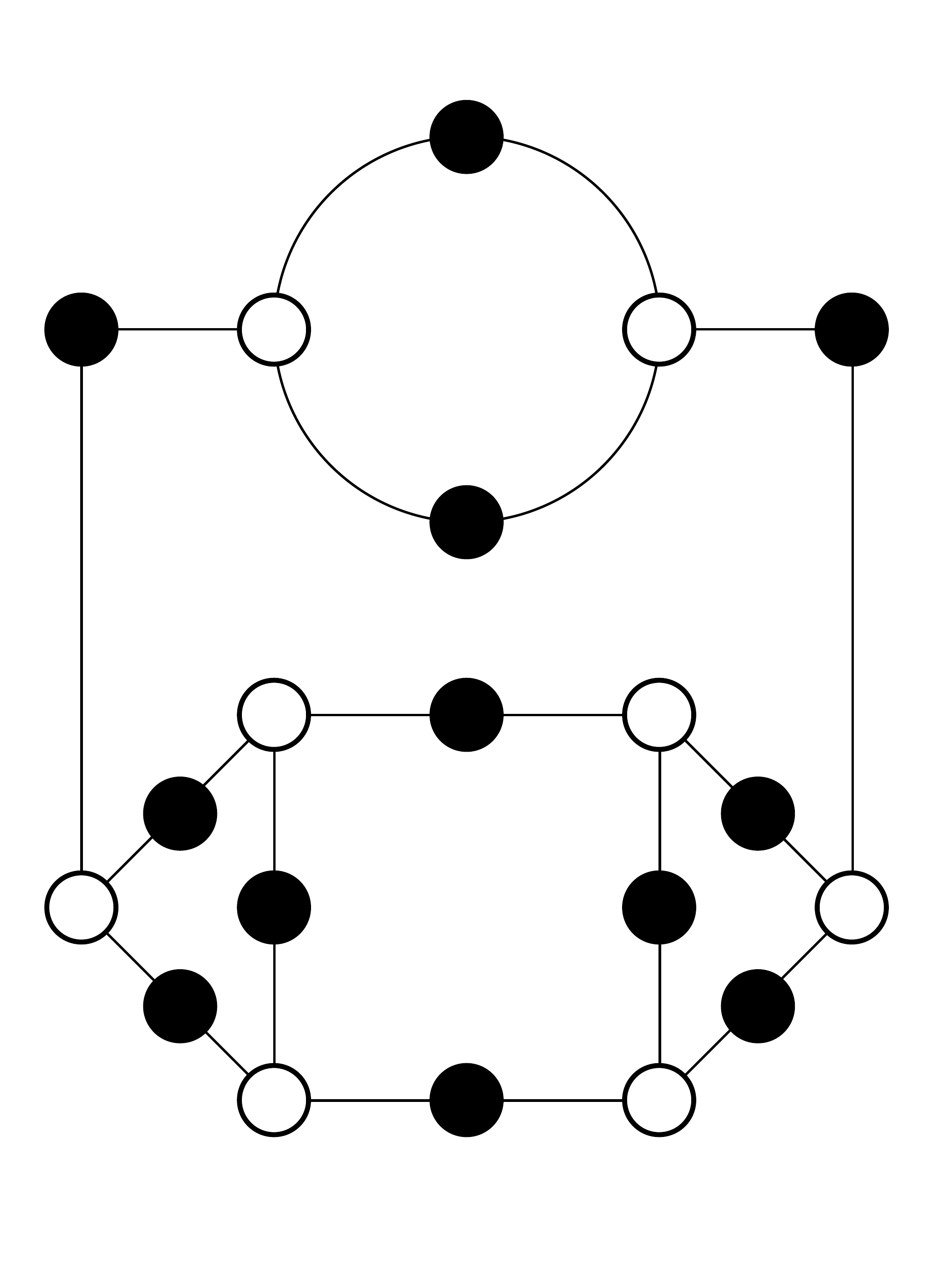}}
\par\end{center}{\scriptsize \par}

\begin{center}
{\scriptsize $6,6,4,3,3,2\;\left(\mathbb{Q}\right)$}
\par\end{center}%
\end{minipage}
\par\end{center}{\scriptsize \par}

\begin{center}
{\scriptsize }%
\begin{minipage}[t]{0.33\textwidth}%
\begin{center}
{\scriptsize \includegraphics[scale=0.15]{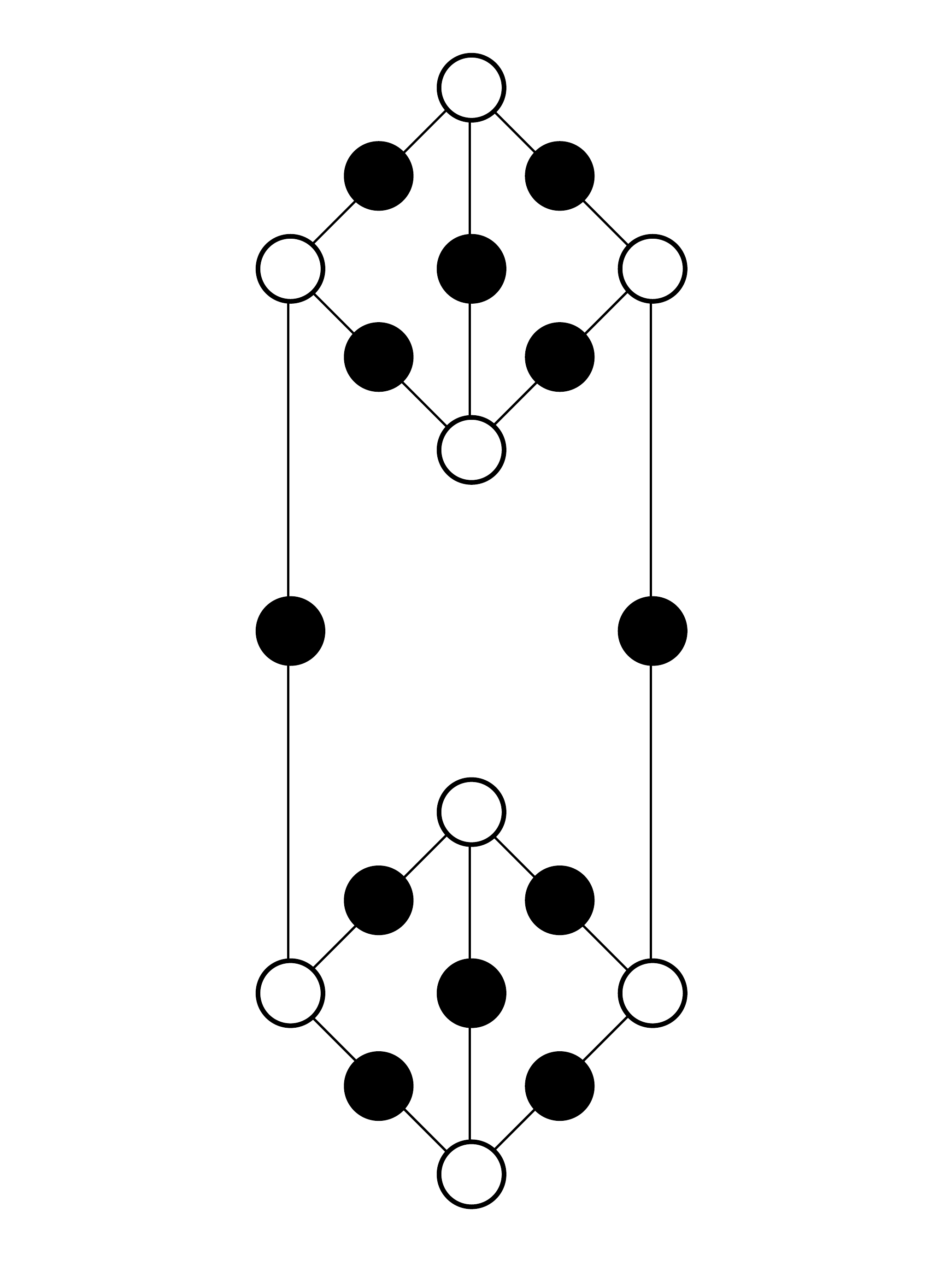}}
\par\end{center}{\scriptsize \par}

\begin{center}
{\scriptsize $6,6,3,3,3,3\;\left(\mathbb{Q}\right)$}
\par\end{center}%
\end{minipage}{\scriptsize }%
\begin{minipage}[t]{0.33\textwidth}%
\begin{center}
{\scriptsize \includegraphics[scale=0.15]{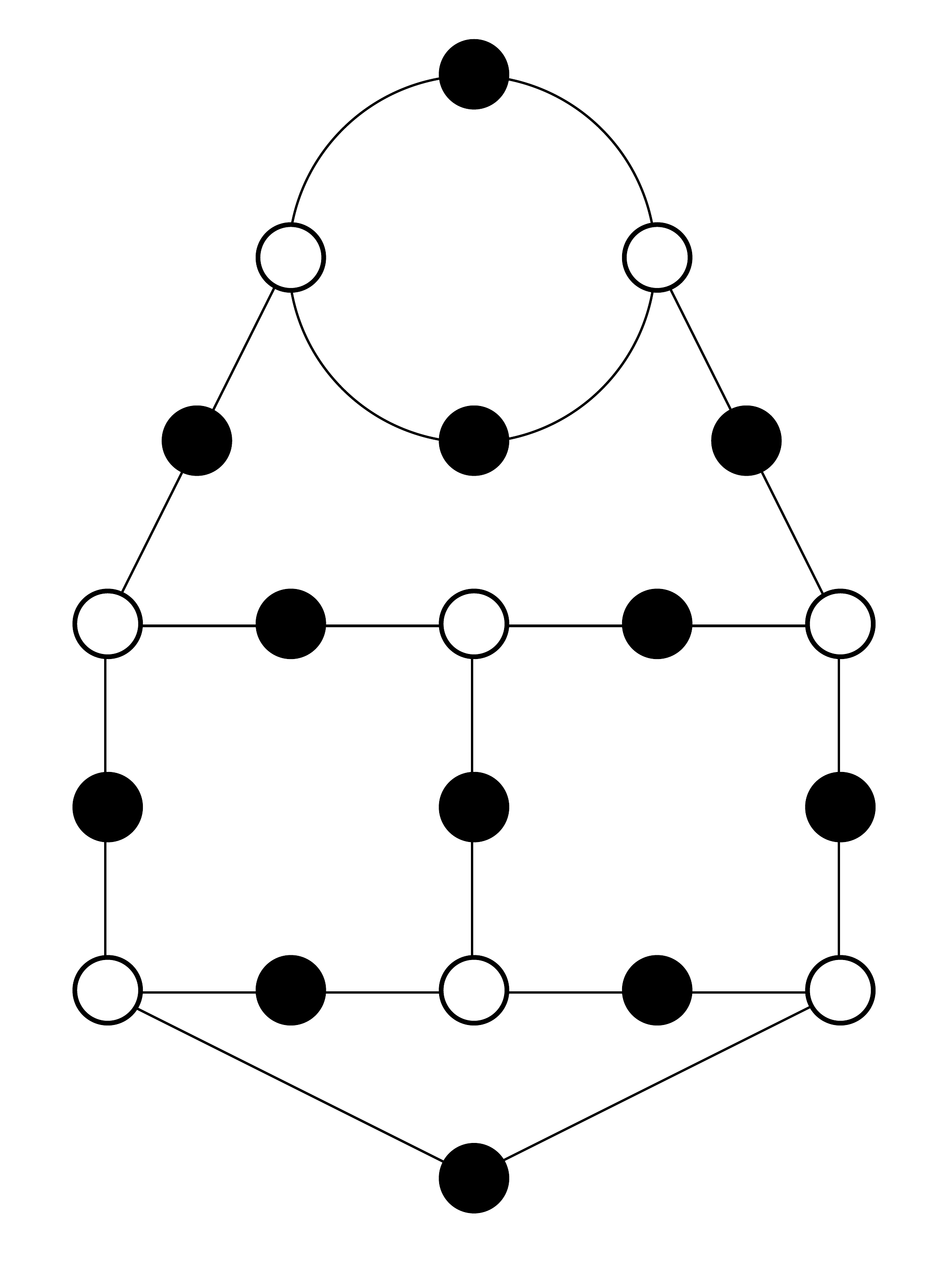}}
\par\end{center}{\scriptsize \par}

\begin{center}
{\scriptsize $6,5,4,4,3,2\;\left(\mathbb{Q}\right)$}
\par\end{center}%
\end{minipage}{\scriptsize }%
\begin{minipage}[t]{0.33\textwidth}%
\begin{center}
{\scriptsize \includegraphics[scale=0.15]{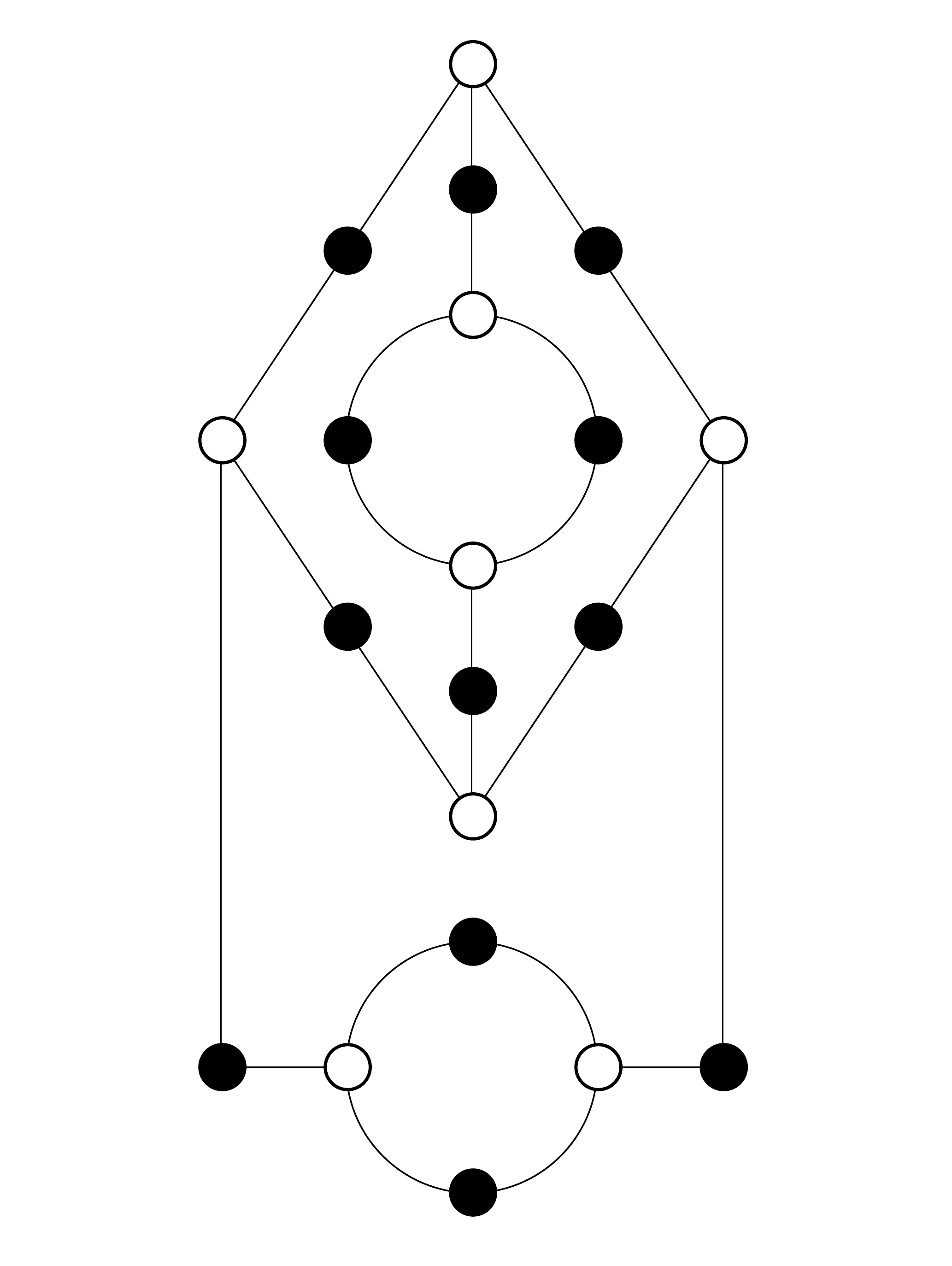}}
\par\end{center}{\scriptsize \par}

\begin{center}
{\scriptsize $5,5,5,5,2,2\;\left(\mathbb{Q}\right)$}
\par\end{center}%
\end{minipage}
\par\end{center}{\scriptsize \par}

\begin{center}
{\scriptsize }%
\begin{minipage}[t]{0.33\textwidth}%
\begin{center}
{\scriptsize \includegraphics[scale=0.15]{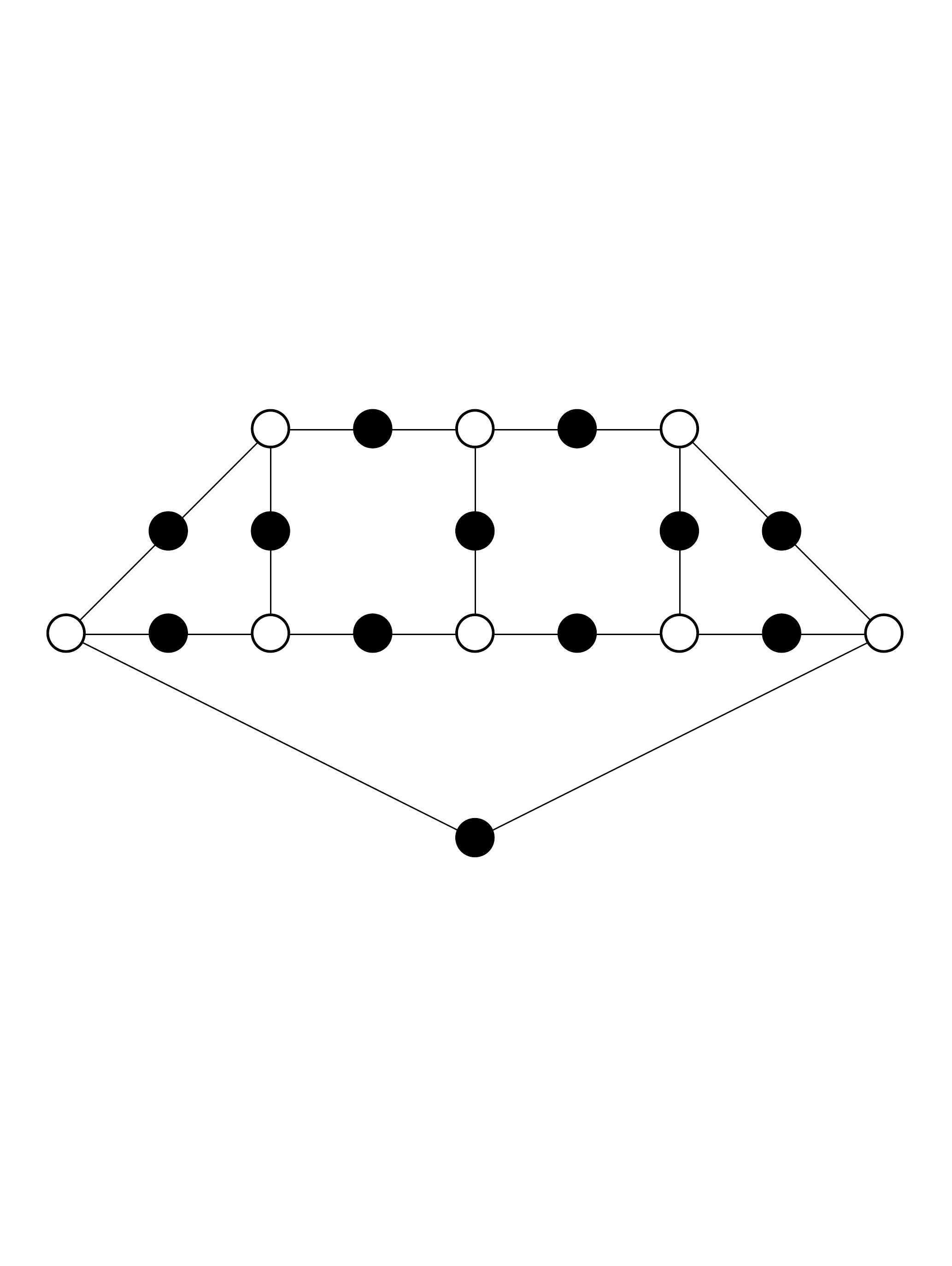}}
\par\end{center}{\scriptsize \par}

\begin{center}
{\scriptsize $5,5,4,4,3,3\;\left(\mathbb{Q}\right)$}
\par\end{center}%
\end{minipage}{\scriptsize }%
\begin{minipage}[t]{0.33\textwidth}%
\begin{center}
{\scriptsize \includegraphics[scale=0.15]{\string"PICT/4-4-4-4-4-4\string".pdf}}
\par\end{center}{\scriptsize \par}

\begin{center}
{\scriptsize $4,4,4,4,4,4\;\left(\mathbb{Q}\right)$}
\par\end{center}%
\end{minipage}
\par\end{center}{\scriptsize \par}

~\\
~\\
~\\


\bibliographystyle{plain}

\affiliationone{
  Yang-Hui He\\
   Department of Mathematics, City University\\
   Northampton Square, London EC1V 0HB, UK;\\
   School of Physics, NanKai University\\
   Tianjin, 300071, P.R. China;\\
   Merton College, University of Oxford\\
   OX1 4JD, UK\\
   \email{yang-hui.he@merton.ox.ac.uk}}
\affiliationtwo{
   John McKay\\
  Department of Mathematics and Statistics, Concordia University\\
  1455 de Maisonneuve Blvd. West\\
  Montreal, Quebec, H3G 1M8, Canada\\
   \email{mckay@encs.concordia.ca}}
\affiliationthree{%
   James Read\\
   Oriel College, University of Oxford\\
   OX1 4EW, UK
   \email{james.read@oriel.ox.ac.uk}}
\end{document}